\documentclass{article}
\usepackage{latexsym, amscd, amsfonts, eucal, mathrsfs, amsmath, amssymb, amsthm, xypic,xr, makeidx, stmaryrd}

\def\Z{\mathbf{Z}}

\DeclareMathOperator{\LPrest}{L\mathcal{P}r}

\DeclareMathOperator{\LPress}{\widehat{\Cat}_{\infty}^{\LPrest}}

\DeclareMathOperator{\Spectralheart}{\SSet_{\infty}^{\heartsuit}}
\DeclareMathOperator{\Mor}{Mor}

\DeclareMathOperator{\nounit}{nu}

\DeclareMathOperator{\Stabb}{\sigma}
\DeclareMathOperator{\qunit}{qu}

\DeclareMathOperator{\NUAdj}{Adj^{\nounit}}
\DeclareMathOperator{\UAdj}{Adj}
\DeclareMathOperator{\CatMod}{CatMod}

\DeclareMathOperator{\LinSeg}{{{\mathcal L}in}_{\ast}}

\DeclareMathOperator{\LPressStab}{\widehat{\Cat}_{\infty}^{\LPrest, \Stabb}}

\DeclareMathOperator{\SegMon}{Mon^{s}}
\DeclareMathOperator{\SegAlg}{Alg^{s}}

\DeclareMathOperator{\aperf}{aperf}

\DeclareMathOperator{\Spread}{\mathcal E}
\DeclareMathOperator{\im}{im}

\newcommand{\AInfty}{\mathfrak{A}_{\infty}}

\DeclareMathOperator{\Free}{Free}
\DeclareMathOperator{\Baar}{Bar}
\DeclareMathOperator{\perf}{perf}
\DeclareMathOperator{\disc}{disc}

\DeclareMathOperator{\conn}{conn}
\DeclareMathOperator{\wMon}{Lax}

\DeclareMathOperator{\AdjDiag}{ADat}
\DeclareMathOperator{\mSet}{\mathcal{S}et_{\Delta}^{+}}

\newcommand{\dplus}{-}
\newcommand{\dminus}{+}

\newcommand{\seg}[1]{{\langle #1 \rangle}_{\ast}}

\newcommand{\nostar}[1]{{\langle #1 \rangle}}

\newcommand{\Un}{Un}
\DeclareMathOperator{\connSpectra}{\SSet_{\infty}^{conn}}
\DeclareMathOperator{\RPres}{\mathcal{P}r^{R}}

\DeclareMathOperator{\calH}{\mathcal{H}}
\DeclareMathOperator{\Mon}{Mon}

\DeclareMathOperator{\Cart}{Cart}
\DeclareMathOperator{\Alg}{Alg}
\DeclareMathOperator{\CAlg}{CAlg}

\DeclareMathOperator{\LFun}{\Fun^{L}}
\DeclareMathOperator{\holim}{holim}
\DeclareMathOperator{\RFun}{\Fun^{R}}
\DeclareMathOperator{\rght}{R}

\newcommand{\MonLPres}{\widehat{\Cat}_{\infty}^{\pres,\otimes}}
\newcommand{\toposref}[1]{T.\ref{HTT-#1}}
\newcommand{\stableref}[1]{S.\ref{STA-#1}}

\DeclareMathOperator{\pres}{Pr}
\newcommand{\MonStab}{\widehat{\Cat}_{\infty}^{\Stabb, \otimes}}
\newcommand{\degree}{\text{o}}
\newcommand{\bfA}{{\mathbf A}}
\newcommand{\bfB}{{\mathbf B}}

\DeclareMathOperator{\cDelta}{{\bf \Delta}}
\DeclareMathOperator{\Set}{\mathcal{S}et}
\DeclareMathOperator{\sSet}{\mathcal{S}et_{\Delta}}
\DeclareMathOperator{\sCoNerve}{\mathfrak{C}}
\DeclareMathOperator{\Nerve}{N}
\DeclareMathOperator{\Cat}{\mathcal{C}at}
\newcommand{\h}[1]{\rm{h} \! #1}
\newcommand{\heart}[1]{{#1}^{\heartsuit}}
\newcommand{\Adjoint}[4]{\xymatrix@1{#2 \ar@<.4ex>[r]^-{#1} & #3 \ar@<.4ex>[l]^-{#4}}}

\DeclareMathOperator{\coker}{coker}

\newcommand{\et}{\'{e}t}

\newcommand{\bigdot}{\bullet}

\DeclareMathOperator{\Stab}{Stab}
\newcommand{\Spectraprime}{\SSet'_{\infty}}

\DeclareMathOperator{\bHom}{Map}

 \DeclareMathOperator{\End}{End}
 
\DeclareMathOperator{\Mod}{\mathcal{M}od}

\DeclareMathOperator{\calM}{\mathcal{M}}

 \DeclareMathOperator{\G}{\mathbf{G}}

\DeclareMathOperator{\Ext}{Ext} 
 
\DeclareMathOperator{\DerRing}{\mathcal{SCR}}

\DeclareMathOperator{\Ex}{Ex}

\DeclareMathOperator{\Tor}{Tor} 
\DeclareMathOperator{\colim}{colim}

\DeclareMathOperator{\calA}{\mathcal{A}}
\DeclareMathOperator{\bd}{\partial}

\newcommand{\Sphere}{S}

\DeclareMathOperator{\calE}{\mathcal{E}}

\DeclareMathOperator{\calO}{\mathcal{O}}

\DeclareMathOperator{\shift}{\eta}
\DeclareMathOperator{\Spectra}{\mathcal{S}_{\infty}}
\DeclareMathOperator{\FinSpace}{\mathcal{S}_{\ast}^{fin}}

\DeclareMathOperator{\FinSpectra}{\mathcal{S}_{\infty}^{fin}}

\DeclareMathOperator{\Spec}{Spec}

\DeclareMathOperator{\Hom}{Hom} 
\DeclareMathOperator{\HH}{H} 
\DeclareMathOperator{\id}{id} \DeclareMathOperator{\Fun}{Fun}
\DeclareMathOperator{\calC}{\mathcal{C}}
\DeclareMathOperator{\calI}{\mathcal{I}}
\DeclareMathOperator{\calN}{\mathcal{N}}
\DeclareMathOperator{\calJ}{\mathcal{J}}
\DeclareMathOperator{\calJJ}{\mathcal{J}}

\DeclareMathOperator{\SSet}{\mathcal{S}}

\DeclareMathOperator{\calX}{\mathcal{X}}
\DeclareMathOperator{\calY}{\mathcal{Y}}

\DeclareMathOperator{\calD}{\mathcal{D}}
\DeclareMathOperator{\Ind}{Ind} 
\DeclareMathOperator{\calP}{\mathcal{P}} \topmargin=0in

\newtheorem{theorem}{Theorem}[subsection]
\newtheorem{lemma}[theorem]{Lemma}
\newtheorem{proposition}[theorem]{Proposition}

\newtheorem{corollary}[theorem]{Corollary}

\theoremstyle{definition}
\newtheorem{definition}[theorem]{Definition}

\newtheorem{example}[theorem]{Example}

\newtheorem{notation}[theorem]{Notation}

\newtheorem{warning}[theorem]{Warning}
\newtheorem{remark}[theorem]{Remark}

\begin{document}

\title{Derived Algebraic Geometry II: Noncommutative Algebra}

\maketitle
\tableofcontents

\section*{Introduction}

Let $\calC$ be an abelian category, and let $X$ be an object of $\calC$. Then the set of endomorphisms $\Hom_{\calC}(X,X)$ is endowed with the structure of an associative ring: the addition on $\Hom_{\calC}(X,X)$ is determined by the additive structure on $\calC$, and the multiplication is
given by composition. In the higher categorical setting, one should consider not a {\em set} of maps
$\Hom_{\calC}(X,X)$, but instead a {\em space} of maps $\bHom_{\calC}(X,X)$. In this setting, the appropriate analogue of an abelian category is a {\it stable $\infty$-category} (see \cite{DAGStable}).
If $\calC$ is a stable $\infty$-category containing an object , then the space $\bHom_{\calC}(X,X)$ should itself be viewed as a kind of ring. There are several approaches to making this precise:

\begin{itemize}
\item[$(1)$] Let $\calH$ denote the homotopy category of topological spaces. If $\calC$ is an $\infty$-category, then the homotopy category $\h{\calC}$ is naturally enriched over $\calH$ (see \S \toposref{prereq1}). If $\calC$ is stable, then $\h{\calC}$ is additive; one can then use classical reasoning to deduce that for every $X \in \calC$, the space $E = \bHom_{\calC}(X,X)$ has the structure of a {\it ring object} of $\calH$: that is,  there exist addition and multiplication maps $E \times E \rightarrow E$ which satisfy the axioms for an associative ring, up to homotopy. Unfortunately, the theory of ring objects in the homotopy category $\calH$ is not a very useful one, because the category $\calH$ is very badly behaved: for example, $\calH$ does not admit finite limits or colimits, so it is impossible to carry out even very simple algebraic constructions in $\calH$.

\item[$(2)$] To avoid the difficulties inherent in approach $(1)$, it is tempting to work in a more rigid setting, such as the theory of {\it topological rings}. A topological ring is simply a ring $E$ equipped with a topology, for which the ring operations are continuous. In this case, one can develop a good theory of (homotopy) limits and colimits, and carry out a wide variety of constructions which mimic classical algebra. However, the setting of topological rings is very restrictive. Typically, if $\calC$ is a stable $\infty$-category containing an object $X$, the space $\bHom_{\calC}(X,X)$ is not homotopy equivalent to a topological abelian group, let alone a topological ring.
\end{itemize}

To obtain a good theory, it is necessary to find some middle ground between $(1)$ and $(2)$: we cannot generally arrange that the space $E = \bHom_{\calC}(X,X)$ satisfies the ring axioms ``on the nose'', but it is also not sufficient to merely assume that these axioms hold up to homotopy. Fortunately, there is a good theory which lies between these two extremes: namely, the theory of {\it associative ring spectra} or, as we will call them, {\it $A_{\infty}$-rings}. The theory of $A_{\infty}$-rings is a generalization of classical (noncommutative) ring theory, in the same sense that stable homotopy theory is a generalization of the classical theory of abelian groups. The theory of $A_{\infty}$-rings and their modules is a fundamental tool in the study of higher categorical mathematics, which can be applied in a great variety of situations. For example, a theorem of Schwede and Shipley asserts that many stable $\infty$-categories (and, by extension, many triangulated categories) can be realized as $\infty$-categories of modules over suitably chosen $A_{\infty}$-rings (Theorem \ref{schwedeshipley}). We also have a more specific reason to study $A_{\infty}$-rings: they are a natural stepping stone towards their commutative cousins, $E_{\infty}$-rings, which are essential to the foundations of derived algebraic geometry. 

Our goal in this paper is to introduce the definition of $A_{\infty}$-rings from an $\infty$-categorical point of view. An ordinary associative ring can be viewed as an algebra object of the category of abelian groups $\calA$, where $\calA$ is endowed with the structure of a monoidal category via the tensor product.
We wish to find an analogue of this definition, where the ordinary category $\calA$ is replaced by the $\infty$-category of spectra. For this, we will need a higher-categorical version of the theory of monoidal categories. The construction of this theory will occupy the bulk of this paper.

We will begin in \S \ref{monoid2.0} by giving the definition of a monoidal $\infty$-category. Roughly speaking, a monoidal $\infty$-category consists of an $\infty$-category $\calC$, equipped
with a distinguished object $1_{\calC} \in \calC$ (the unit) and a bifunctor
$\otimes: \calC \times \calC \rightarrow \calC$, which is unital and associative up to coherent homotopy. We will describe several sources of examples of monoidal $\infty$-categories, and various constructions for producing new monoidal $\infty$-categories from old. We will also
introduce the notion of an {\it algebra object} of a monoidal $\infty$-category $\calC$.
Roughly speaking, this is an object $A \in \calC$ equipped with a multiplication
$A \otimes A \rightarrow A$ and a unit map $1_{\calC} \rightarrow A$, which are again unital and associative up to coherent homotopy. The collection of all algebra objects of $\calC$ can itself be organized into an $\infty$-category, which we will denote by $\Alg(\calC)$. 
We will study the $\infty$-category $\Alg(\calC)$ in some detail. In particular, we will introduce criteria which imply the existence of a good supply of limits and colimits in $\Alg(\calC)$, and which guarantee that the forgetful functor $\Alg(\calC) \rightarrow \calC$ admits a left adjoint. 

Let $\calC$ be a monoidal $\infty$-category and $A$ an algebra object of $\calC$. In this case, we can speak of {\it $A$-modules} in $\calC$: that is, objects $M \in \calC$ equipped with a structure map $A \otimes M \rightarrow M$ which is unital and associative up to coherent homotopy. In \S \ref{hugr2.0}, we will study the theory of modules in a slightly more general setting, where
$M$ is taken to be an object of an $\infty$-category $\calM$ which is {\it tensored over} $\calC$ (Definition \ref{ulult}).

Let $\calM$ be a fixed $\infty$-category. Then there is a universal example of a monoidal $\infty$-category $\calC$ such that $\calM$ is tensored over $\calC$: namely, one can take
$\calC$ to be the $\infty$-category $\Fun(\calM, \calM)$ of {\it endofunctors} of $\calM$.
Algebra objects of $\Fun(\calM, \calM)$ are called {\it monads} on $\calM$. Given a monad
$T$ on $\calM$, one can define a new $\infty$-category $\Mod_{T}(\calM)$ of {\it $T$-modules}
in $\calM$. There is a forgetful functor $\Mod_{T}(\calM) \rightarrow \calM$, which admits a left
adjoint $\calM \rightarrow \Mod_{T}(\calM)$.
In general, given a pair of adjoint functors 
$$ \Adjoint{F}{\calM}{\calN}{G},$$
one can canonically associate a monad $T$ on $\calM$, such that the functor
$G: \calN \rightarrow \calM$ factors through some functor $G': \calN \rightarrow \Mod_{T}(\calM)$. In classical category theory, the {\it Barr-Beck theorem} provides necessary and sufficient conditions for
$G'$ to be an equivalence. We will prove an $\infty$-categorical version of the Barr-Beck theorem in \S \ref{barri} (Theorem \ref{barbeq}).

In \S \ref{monoid6.0} we will specialize to the main situation of interest: the case where $\calC$ is the $\infty$-category $\Spectra$ of {\em spectra}. We will see that $\Spectra$ admits an essentially unique monoidal structure with the property that the functor $\otimes: \Spectra \times \Spectra \rightarrow \Spectra$ preserves colimits in each variable (Corollary \ref{surcoi}). The bifunctor can be identified with the classical {\em smash product} operation on spectra. We then define an {\it $A_{\infty}$-ring} to be an algebra object of $\Spectra$. We will then proceed to show that a great deal of classical algebra can be carried out in the setting of $A_{\infty}$-rings and their modules: for example, we will introduce a theory of flat modules, and prove an analogue of Lazard's theorem:
an $A$-module $M$ is flat if and only if $M$ can be obtained as a filtered colimit of free $A$-modules (Theorem \ref{lazard}).

We should emphasize that the theory of $A_{\infty}$-rings and their modules is not new. There are various definitions available in the literature; see, for example, \cite{EKMM}. We have chosen to present the subject using the language of $\infty$-categories, which we feel is the natural home for these ideas. 

 \subsection*{Notation and Terminology}

For an introduction to the language of higher category theory (from the point of view taken in this paper), we refer the reader to \cite{topoi}. We will use the terminology and results of \cite{topoi}. We will also use \cite{DAGStable} as our reference for the theory of stable $\infty$-categories. References to \cite{topoi} will be indicated by use of the letter T, and references to \cite{DAGStable} will be indicated by use of the letter S. For example, Theorem \toposref{mainchar} refers to Theorem \ref{HTT-mainchar} of \cite{topoi}.

For each integer $n \geq -1$, we let $[n]$ denote the linearly ordered set
$\{ 0, \ldots, n \}$ (so that $[-1]$ denotes the empty set).
Throughout this paper, $\cDelta$ will denote the category of {\it combinatorial simplices}:
the objects of $\cDelta$ are the linearly ordered sets $\{ [n] \}_{ n \geq 0}$, and the morphisms
are (nonstrictly) increasing maps of linearly ordered sets. \index{ZZZcDelta@$\cDelta$}
The category $\cDelta$ is equivalent to the larger category of {\em all} nonempty linearly ordered finite sets. We will typically abuse notation by not distinguishing between $\cDelta$ and this larger category. In other words, if $J$ is a nonempty finite linearly ordered set, we will implicitly identify
$J$ with an object $[n] \in \cDelta$ via an isomorphism of linearly ordered sets $\alpha: J \simeq [n]$; we note that there is little risk in doing so, since the isomorphism $\alpha$ and the integer $n$ are uniquely determined.

If $p: X \rightarrow S$ is a map of simplicial sets and $s$ is a vertex of $S$, we will typically write
$X_{s}$ to denote the fiber $X \times_{S} \{s\}$. 
 
\section{Monoidal $\infty$-Categories}\label{monoid2.0}
  
Our main goal in this section is to introduce the definition of a monoidal $\infty$-category and to survey some of the basic examples. We will begin in \S \ref{monoid2.1} with a brief review of the classical theory of monoidal categories. We will then show how this theory can be reformulated in such a way that it admits a natural $\infty$-categorical generalization (Definition \ref{mainef}).


Many of the monoidal $\infty$-categories which we study in this paper will arise via one of the following constructions:
\begin{itemize}
\item[$(1)$] Let $\calC$ be an $\infty$-category which admits finite products. In this case,
we can endow $\calC$ with the {\it Cartesian monoidal structure}, in which the bifunctor
$\otimes: \calC \times \calC \rightarrow \calC$ is given by the Cartesian product. We will study this example in \S \ref{monoidcart}.

\item[$(2)$] Let $\calC$ be a monoidal $\infty$-category, and let $\calD \subseteq \calC$ be a full subcategory which contains the unit object and is stable under tensor products. Then $\calD$ inherits the structure of a monoidal $\infty$-category. In \S \ref{locol} we will study this construction together with a dual procedure, which produces monoidal structures on suitable {\em quotients} of $\calC$.

\item[$(3)$] Let $\calC$ be an ordinary monoidal category. Then the nerve $\Nerve(\calC)$ has the structure of a monoidal $\infty$-category. In \S \ref{monoidate} we will describe some variations on this observation. For example, we will show that if $\bfA$ is a (simplicial) monoidal model category, then the underlying $\infty$-category $\Nerve(\bfA^{\degree})$ of $\bfA$ inherits the structure of a monoidal $\infty$-category (Proposition \ref{hurgoven}). 
\end{itemize}
  
If $\calC$ is a monoidal $\infty$-category, then we can construct an $\infty$-category
$\Alg(\calC)$ of {\it algebra objects} of $\calC$; the definition will be given in \S \ref{monoid2.1}. Our second goal in this section is to analyze the $\infty$-category $\Alg(\calC)$. In \S \ref{limalg}, we will establish existence criteria for limits and colimits in $\Alg(\calC)$. The most difficult step is to prove that $\Alg(\calC)$ admits coproducts (given suitable assumptions on $\calC$). To prove this, we will need to know that the forgetful functor $\Alg(\calC) \rightarrow \calC$ admits a left adjoint; in other words, we need to be able to construct the {\it free algebra} generated by an object $C \in \calC$. 
We will describe the (rather technical) details of this construction in \S \ref{monoid5}. Our method requires a reformulation of the definition of monoidal $\infty$-categories, which we will discuss in \S \ref{segalapp}.

\subsection{Monoidal Structures and Algebra Objects}\label{monoid2.1}
  
Recall that a {\it monoid} is a set $M$ equipped with a multiplication $M \times M \rightarrow M$ and a unit object $1 \in M$ satisfying the identities
$$ 1x = x1 = x \quad x(yz)=(xy)z$$\index{monoid}\index{monoidal!category}\index{category!monoidal}
for all $x,y,z \in M$. Roughly speaking, a {\it monoidal category} is a category $\calC$ equipped with the same kind of structures: a unit object $1 \in \calC$, and a bifunctor $\otimes: \calC \times \calC$. However, in the categorical setting, it is unnatural to require the identities displayed above
to hold as equalities. In general, we do not expect $X \otimes (Y \otimes Z)$ to be {\em equal} to
$(X \otimes Y) \otimes Z$. Instead, the associative law should be formulated as the existence of
an isomorphism $\eta_{X,Y,Z}: (X \otimes Y) \otimes Z \simeq X \otimes (Y \otimes Z)$. 
Moreover, the isomorphisms $\eta_{X,Y,Z}$ are taken as additional data, and are required to satisfy further conditions (such as naturality in $X$, $Y$, and $Z$); we refer the reader to \S \toposref{monoidaldef} for a detailed definition. If we try to generalize this definition to higher categories, then the equations satisfied by the isomorphisms $\eta_{X,Y,Z}$ should themselves hold only up to isomorphism. It is possible to explicitly describe all of the relevant data (for example, using the theory of {\it Stasheff associahedra}), but the escalation in complexity is somewhat intimidating; it will be more convenient to proceed in another way.

We begin by considering an example of a monoidal category. Let $\calC$ be the category of complex vector spaces, with monoidal structure given by tensor products of vector spaces. Given a pair of vector spaces $U$ and $V$, the tensor product $U \otimes V$ is defined by the property that $\Hom_{\calC}(U \otimes V, W)$ can be identified with the set of {\em bilinear} maps $U \times V \rightarrow W$. In fact, this property really only determines $U \otimes V$ up to canonical isomorphism: in order to build an actual tensor product functor, we need to choose some particular
construction of $U \otimes V$. Because this requires making certain decision in an ad-hoc manner, it is unrealistic to expect an equality of vector spaces $U \otimes (V \otimes W) = (U \otimes V) \otimes W$. However, the existence of a canonical {\em isomorphism} between
$U \otimes (V \otimes W)$ and $(U \otimes V) \otimes W$ is easily explained: linear maps from either into a fourth vector space $X$ can be identified with trilinear maps $U \times V \times W \rightarrow X$.

The above example suggests that we might reformulate the definition of a monoidal category as follows. Rather than give a bifunctor $\otimes: \calC \times \calC \rightarrow \calC$, we instead specify, for each $n$-tuple $(C_1, \ldots, C_n)$ of objects of $\calC$ and each $D \in \calC$, the collection of morphisms $C_1 \otimes \ldots \otimes C_n \rightarrow D$.
Of course, we also need to specify how such morphisms are to be composed. The relevant data can be encoded in a new category $\calC^{\otimes}$:

\begin{definition}\label{converi}\index{ZZZCotimes@$\calC^{\otimes}$}
Let $(\calC, \otimes)$ be a monoidal category. We define a new category
$\calC^{\otimes}$ as follows:
\begin{itemize}
\item[$(i)$] An object of $\calC^{\otimes}$ is a finite (possibly empty) sequence
of objects of $\calC$, which we will denote by $[ C_1, \ldots, C_n ]$.
\item[$(ii)$] A morphism from $[C_1, \ldots, C_n]$ to $[C'_1, \ldots, C'_m]$ in $\calC^{\otimes}$ consists of a nonstrictly order-preserving map $f: [m] \rightarrow [n]$, and a collection of morphisms
$C_{ f(i-1)+1} \otimes \ldots \otimes C_{f(i)} \rightarrow C'_{i}$ for $1 \leq i \leq m$.
\item[$(iii)$] Composition in $\calC^{\otimes}$ is determined by composition of order preserving maps, composition in $\calC$, and the associativity and unit constraints of the monoidal structure on $\calC$.
\end{itemize}
\end{definition}

Let $(\calC, \otimes)$ be a monoidal category, and let $\calC^{\otimes}$ be as in Definition \ref{converi}. Then there is forgetful functor $p: \calC^{\otimes} \rightarrow \cDelta^{op}$, which carries an object $[ C_1, \ldots, C_n]$ to the linearly ordered set $[n]$. Moreover, $p$ has the following properties:

\begin{itemize}
\item[$(M1)$] The functor $p$ is an {\it op-fibration of categories}. In other words, for every object $[ C_1, \ldots, C_n] \in \calC^{\otimes}$ and every
morphism $f: [n] \rightarrow [m]$ in $\cDelta^{op}$, there exists a morphism
$\overline{f}: [C_1, \ldots, C_n] \rightarrow [C'_{1}, \ldots, C'_{m}]$ which covers $f$, and is
universal in the sense that composition with $\overline{f}$ induces a bijection
$$ \Hom_{ \calC^{\otimes} }( [C'_{1}, \ldots, C'_{m}], [ C''_{1}, \ldots, C''_{k}])
\rightarrow \Hom_{\calC^{\otimes}}( [C_{1}, \ldots, C_{n} ], [C''_{1}, \ldots, C''_{k}])
\times_{ \bHom_{\cDelta}( [k], [n]) } \bHom_{\cDelta}([k],[m])$$ 
for every object $[C''_{1}, \ldots, C''_{k}] \in \calC^{\otimes}$.
To achieve this, it suffices to choose $\overline{f}$ such that the maps $C_{ f(i-1)+1} \otimes \ldots \otimes C_{f(i)} \rightarrow C'_{i}$ are isomorphisms for $1 \leq i \leq m$.

\item[$(M2)$] Let $\calC^{\otimes}_{[n]}$ denote the fiber of $p$ over the object $[n] \in \cDelta^{op}$. 
Then $\calC^{\otimes}_{[1]}$ is equivalent to $\calC$. More generally, $\calC^{\otimes}_{[n]}$ is equivalent to an $n$-fold product of copies of $\calC$. The equivalence is induced by
functors associated to the inclusions $[1] \simeq \{ i-1, i \} \subseteq [n]$, $1 \leq i \leq n$.
\end{itemize}

We now observe that the monoidal structure on $\calC$ is determined, up to canonical equivalence, by $\calC^{\otimes}$. More precisely, suppose given a functor $p: \calD \rightarrow \cDelta^{op}$ satisfying conditions $(M1)$ and $(M2)$. Then:
\begin{itemize}
\item[$(a)$] Condition $(M2)$ implies that $\calD_{[0]}$ has a single object, up to equivalence. The projection $[1] \rightarrow [0]$ determines a functor $\calD_{[0]} \rightarrow \calD_{[1]} \simeq \calC$, which we can identify with an object $1 \in \calC$.

\item[$(b)$] The inclusion
$[1] \simeq \{ 0, 2\} \subseteq [2]$ determines a functor $\calC \times \calC \simeq \calD_{[2]}
\rightarrow \calD_{[1]} \simeq \calC$, which we may denote by $\otimes$. 

\item[$(c)$] The commutative diagram
$$ \xymatrix{ \{ 0, 3 \} \ar[r] \ar[d] & \{ 0, 1, 3 \} \ar[d] \\
\{ 0, 2, 3 \} \ar[r] & \{0, 1, 2, 3 \} } $$
in $\cDelta$ determines a
diagram of categories and functors (which commutes up to canonical isomorphism):
$$ \xymatrix{ \calD_{[1]} & \calD_{[2]} \ar[l] \\
\calD_{[2]} \ar[u] & \calD_{[3]}. \ar[u] \ar[l] }$$ 
Combining this with the equivalences $\calD_{[n]} \simeq \calC^{n}$, we obtain
a functorial isomorphism
$$ \eta_{A,B,C}: (A \otimes B) \otimes C \simeq A \otimes (B \otimes C).$$
A similar argument can be used to construct canonical isomorphisms
$$ 1 \otimes X \simeq X \simeq X \otimes 1.$$
\end{itemize}

It is not difficult to see that $(a)$, $(b)$, and $(c)$ endow $\calC$ with the structure of a monoidal category. For example, MacLane's pentagon axiom asserts that the diagram
$$ \xymatrix{ & ((A \otimes B) \otimes C) \otimes D \ar[dl]^{\eta_{A,B,C} \otimes \id_D}
\ar[dr]^{\eta_{ A \otimes B,C,D}} \\
(A \otimes (B \otimes C)) \otimes D \ar[d]^{ \eta_{A, B \otimes C, D}} & & (A \otimes B) \otimes (C \otimes D) \ar[d]^{\eta_{A,B, C \otimes D}} \\
A \otimes (( B \otimes C) \otimes D) \ar[rr]^{\id_A \otimes \eta_{B,C,D}} & & A \otimes (B \otimes (C \otimes D))}$$
is commutative. This follows from the fact that all five expressions can be canonically identified with the image of $(A,B,C,D)$ under the composite functor
$ \calC^{4} \simeq \calD_{[4]} \rightarrow \calD_{ \{0,4\} } \simeq \calC.$
In the case where $\calC$ is equipped with a monoidal structure and $\calD = \calC^{\otimes}$, then it is easy to see that the data provided by $(a)$, $(b)$ and $(c)$ recovers the original monoidal structure on $\calC$ (up to canonical isomorphism). Conversely, an arbitrary functor $\calD \rightarrow \cDelta^{op}$ satisfying $(M1)$ and $(M2)$ determines a monoidal structure on $\calC$ and an equivalence $\calD \simeq \calC^{\otimes}$. In other words, specifying a monoidal structure on $\calC$ is {\em equivalent} to specifying the functor $\calC^{\otimes} \rightarrow \cDelta^{op}$. However, the second approach has several advantages over the first:

\begin{itemize}
\item As we saw above in the case of vector spaces, it is sometimes easier to specify the
category $\calC^{\otimes}$ than to specify the bifunctor $\otimes$, in the sense that it requires fewer arbitrary choices. 

\item Axioms $(M1)$ and $(M2)$ concerning the functor $\calC^{\otimes} \rightarrow \cDelta^{op}$ are a bit simpler than the usual definition of a monoidal category. Complicated statements, such as the commutativity of MacLane's pentagon, are consequences of $(M1)$ and $(M2)$.
\end{itemize}

The significance of the latter point becomes more apparent in the $\infty$-categorical setting, where we expect the MacLane pentagon to be only the first step in a hierarchy of coherence conditions of ever-increasing complexity. Fortunately, the above discussion suggests an approach which does not require us to formulate these conditions explicitly.

\begin{definition}\label{mainef}\index{monoidal $\infty$-category}\index{$\infty$-category!monoidal}
A {\it monoidal $\infty$-category} is a coCartesian fibration of simplicial sets
$p: \calC^{\otimes} \rightarrow \Nerve(\cDelta)^{op}$ with the following property:
\begin{itemize}
\item[$(\ast)$] For each $n \geq 0$, the associated functors
$\calC^{\otimes}_{[n]} \rightarrow \calC^{\otimes}_{ \{i, i+1\} }$ determine an equivalence of $\infty$-categories
$$ \calC^{\otimes}_{[n]} \rightarrow \calC^{\otimes}_{ \{0,1\} } \times \ldots \times \calC^{\otimes}_{ \{n-1, n\} } \simeq (\calC^{\otimes}_{[1]})^{n}.$$
\end{itemize}
\end{definition}

\begin{remark}\index{ZZZotimes@$\otimes$}
Let $p: \calC^{\otimes} \rightarrow \Nerve(\cDelta)^{op}$ be a monoidal $\infty$-category. We will refer to the fiber $\calC = \calC^{\otimes}_{[1]}$ as the {\it underlying $\infty$-category} of $\calC^{\otimes}$; we will also say that $\calC^{\otimes}$ is a {\it monoidal structure} on $\calC$. 
Assertion $(\ast)$ implies that $\calC^{\otimes}_{[0]}$ is a contractible Kan complex. The projection $[1] \rightarrow [0]$ in $\cDelta$ determines a functor
$\calC^{\otimes}_{[0]} \rightarrow \calC^{\otimes}_{[1]}$, which we can identify with an object $1_{\calC} \in \calC^{\otimes}_{[1]}$, well-defined up to equivalence. We will refer to $1_{\calC} \in \calC^{\otimes}_{[1]}$ as the {\it unit object} of $\calC^{\otimes}_{[1]}$.\index{unit object}\index{unit}

The three embeddings of $[1]$ into $[2]$ determine a diagram
$$ \calC \times \calC \simeq \calC^{\otimes}_{ \{0,1\} } \times \calC^{\otimes}_{ \{1,2\} } \stackrel{\theta}{\leftarrow} \calC^{\otimes}_{[2]} \stackrel{\theta'}{\rightarrow} \calC^{\otimes}_{ \{0,2\}} \simeq \calC,$$
where $\theta$ is an equivalence. Composing $\theta'$ with a homotopy inverse to $\theta$, we obtain a functor $\otimes: \calC \times \calC \rightarrow \calC,$
which is again well-defined up to equivalence.
\end{remark}

\begin{remark}
Let $p: \calC^{\otimes} \rightarrow \Nerve(\cDelta)^{op}$ be a monoidal $\infty$-category in the sense of Definition \ref{mainef}, and let $f: \calC^{\otimes}_{[1]} \rightarrow \calC$ be an equivalence of $\infty$-categories. In this situation, we will generally abuse terminology by saying that $\calC$ is a {\it monoidal $\infty$-category}, or that $p: \calC^{\otimes} \rightarrow \Nerve(\cDelta)^{op}$ {\it exhibits $\calC$ as a monoidal $\infty$-category}.

Informally, we can think of $\calC^{\otimes}$ as encoding $\calC$ together with the
the tensor operation $\otimes: \calC \times \calC \rightarrow \calC$. Of course, Definition \ref{mainef} encodes a good deal more structure; this additional structure expresses the idea that the operation $\otimes$ is associative, up to {\em coherent} homotopy. 
\end{remark}

\begin{example}\label{itereat}
Let $\calC$ be a monoidal category, and let $\calC^{\otimes}$ be as in Definition \ref{converi}. Then
the induced map $\Nerve(\calC^{\otimes}) \rightarrow \Nerve(\cDelta)^{op}$ is a monoidal $\infty$-category. We will generalize this construction in \S \ref{monoidate}.
\end{example}

\begin{remark}\label{eatfine}
Let $\calC^{\otimes} \rightarrow \Nerve( \cDelta)^{op}$ be a monoidal structure on an $\infty$-category $\calC$. Then the induced functor $\h{ \calC^{\otimes} } \rightarrow \cDelta^{op}$ satisfies the axioms $(M1)$ and $(M2)$, and therefore determines a monoidal structure on the homotopy category $\h{\calC}$.
\end{remark}


In order to make effective use of the theory of monoidal $\infty$-categories, we will need an associated theory of {\it monoidal functors}.

\begin{definition}\label{defcondef}\index{convex}
A morphism $f: [m] \rightarrow [n]$ in $\cDelta$ is {\it convex} if $f$ is injective and the image
$\{ f(0), \ldots, f(m) \} \subseteq [n]$ is a convex subset of $[n]$.
\end{definition}

\begin{definition}\label{monfunc1}\index{monoidal!functor}\index{lax monoidal functor}\index{functor!monoidal}\index{functor!lax monoidal}
Let $\calC^{\otimes}$ and $\calD^{\otimes}$ be monoidal $\infty$-categories. We will say that a functor
$F: \calC^{\otimes} \rightarrow \calD^{\otimes}$ is {\it monoidal}
if the diagram
$$ \xymatrix{ \calC^{\otimes} \ar[rr]^{F} \ar[dr]^{p} & & \calD^{\otimes} \ar[dl]^{q} \\
& \Nerve(\cDelta)^{op} & }$$
commutes, and $F$ carries $p$-coCartesian morphisms to $q$-coCartesian morphisms.
We will say that $F$ is {\it lax monoidal} if
the above diagram commutes, and the following weaker condition is satisfied:
for every $p$-coCartesian morphism $\alpha$ in $\calC^{\otimes}$ such that $p(\alpha)$ is a convex morphism in $\cDelta$, the morphism $F(\alpha)$ is $q$-coCartesian.
We let $\Fun^{\wMon}(\calC^{\otimes}, \calD^{\otimes})$ denote the full subcategory
of $\bHom_{ \Nerve(\cDelta)^{op}}( \calC^{\otimes}, \calD^{\otimes})$ spanned by the lax monoidal functors, and $\Fun^{\Mon}(\calC^{\otimes}, \calD^{\otimes})$ the full subcategory
of $\Fun^{\wMon}(\calC^{\otimes}, \calD^{\otimes})$ spanned by the monoidal functors.
\end{definition}

\begin{remark}\label{labus}
Let $p: \calC^{\otimes} \rightarrow \Nerve(\cDelta)^{op}$ and $q: \calD^{\otimes} \rightarrow \Nerve(\cDelta)^{op}$ be monoidal $\infty$-categories, and let
$F \in \bHom_{ \Nerve(\cDelta)^{op} }( \calC^{\otimes}, \calD^{\otimes})$. Using the
equivalence $\calD^{\otimes}_{[n]} \simeq (\calD^{\otimes}_{[1]})^{n}$, we see that
$F$ is a (lax) monoidal functor if and only if, for every $p$-coCartesian morphism $f$ in $\calC^{\otimes}$ covering a (convex) map $[1] \rightarrow [n]$ in $\cDelta$, the image $F(f)$ is a $q$-coCartesian morphism in $\calD$. 
\end{remark}

\begin{remark}
Roughly speaking, a lax monoidal functor $F: \calC^{\otimes} \rightarrow \calD^{\otimes}$
consists of a functor $f$ between the underlying $\infty$-categories $\calC$ and $\calD$ equipped with a coherently associative collection of natural transformations 
$$f(C_1) \otimes \ldots \otimes f(C_n) \rightarrow f( C_1 \otimes \ldots \otimes C_{n} ).$$
The functor $F$ is monoidal if and only if each of these transformations is an equivalence.
\end{remark}

\begin{remark}
The notion of a {\it lax monoidal} functor as given by Definition \ref{monfunc1} is the $\infty$-categorical analogue of what we called a {\em right lax monoidal functor} in \S \toposref{monoidaldef}. There is also a dual notion which corresponds to {\em left} lax monoidal functors, but 
this notion is not so easily encoded using the formalism introduced above; see Remark \ref{selfop}.
\end{remark}

It follows immediately from the definition that the class of (lax) monoidal functors is stable under composition. Consequently, we may define simplicial categories
$\Cat_{\infty}^{\Delta, \Mon} \subseteq \Cat_{\infty}^{\Delta, \wMon}$ as follows:
\begin{itemize}
\item[$(1)$] The objects of $\Cat_{\infty}^{\Delta, \Mon}$ and $\Cat_{\infty}^{\Delta, \wMon}$
are monoidal $\infty$-categories $\calC^{\otimes} \rightarrow \Nerve(\cDelta)^{op}$.

\item[$(2)$] Given a pair of monoidal $\infty$-categories $\calC^{\otimes}$ and $\calD^{\otimes}$, we let
$$\bHom_{ \Cat_{\infty}^{\Delta,\Mon}}( \calC^{\otimes}, \calD^{\otimes}) \subseteq \Fun^{\Mon}( \calC^{\otimes}, \calD^{\otimes}) \quad
 \bHom_{ \Cat_{\infty}^{\Delta, \wMon}}( \calC^{\otimes}, \calD^{\otimes}) \subseteq \Fun^{\wMon}(\calC^{\otimes}, \calD^{\otimes})$$
be the largest Kan complexes contained in $\Fun^{\Mon}(\calC^{\otimes}, \calD^{\otimes})$ and
$\Fun^{\wMon}(\calC^{\otimes}, \calD^{\otimes})$, respectively. 
\end{itemize}

\begin{definition}\label{monfunc2}\index{ZZZCatinftyMon@$\Cat_{\infty}^{\Mon}$}
\index{ZZZCatinftywMon@$\Cat_{\infty}^{\wMon}$}
We let $\Cat_{\infty}^{\Mon}$ denote the simplicial nerve $\Nerve( \Cat_{\infty}^{\Delta, \Mon})$,
and $\Cat_{\infty}^{\wMon}$ the simplicial nerve $\Nerve( \Cat_{\infty}^{\Delta, \wMon})$. 
We will refer to $\Cat_{\infty}^{\Mon}$ as the {\it $\infty$-category of monoidal $\infty$-categories}.
\end{definition}

\begin{remark}\label{eggsal}
Let $F: \calC^{\otimes} \rightarrow \calD^{\otimes}$ be a monoidal functor between monoidal $\infty$-categories. Using Corollary \toposref{usefir}, we deduce that $F$ is an equivalence if and only if $F$ induces an equivalence of underlying ordinary categories $\calC^{\otimes}_{[1]} \rightarrow \calD^{\otimes}_{[1]}$. The analogous assertion for lax monoidal functors is false.
\end{remark}

Our next goal is to introduce the notion of an {\it algebra object} of a monoidal $\infty$-category $\calC^{\otimes} \rightarrow \Nerve(\cDelta)^{op}$. We begin by considering the classical case. Let $\calC$ be a monoidal category. An {\it algebra object} of $\calC$ is an object $A \in \calC$
equipped with maps
$$ 1 \rightarrow A \quad A \otimes A \rightarrow A$$
which satisfy the usual unit and associativity conditions.
In this case, we can define for each $n \geq 0$ a map $s_n: A^{ \otimes n} \rightarrow A$, given by iterated multiplication (or by the unit, in the case $n=0$); here $A^{\otimes n}$ denotes the $n$-fold tensor product of $A$ with itself (which is well-defined up to canonical isomorphism, in view of the associativity constraint on the tensor product $\otimes$). The associative law can then be reformulated as follows: given integers $n \geq 0$ and $k_1, \ldots, k_n \geq 0$, the composition
$$ A^{\otimes k_1 + \ldots + k_n} \stackrel{ s_{k_1} \otimes \ldots \otimes s_{k_n}}{\rightarrow}
A^{\otimes n} \stackrel{s_n}{\rightarrow} A$$
coincides with $s_{k_1 + \ldots + k_n}$. When the definition is phrased in this way, an algebra object of $\calC$ can be regarded as a {\em section} $\widetilde{A}$ of the projection $\calC^{\otimes} \rightarrow \cDelta^{op}$. More precisely, we associate to each $[n] \in \cDelta^{op}$ the lifting
$\widetilde{A}([n]) = [A, \ldots, A] \in \calC^{\otimes}$. To a map $f: [m] \rightarrow [n]$ in $\cDelta$, we associate the ``contraction'' map $\widetilde{A}([n]) \rightarrow \widetilde{A}([m])$, given by
$( s_{ f(1)-f(0)}, \ldots, s_{f(m)- f(m-1)})$. The functoriality of this construction encodes the associativity of the product on $A$. 

Of course, not {\em every} section of the projection
$\calC^{\otimes} \rightarrow \cDelta^{op}$ corresponds to an algebra object of $\calC$. If $\widetilde{A}: \cDelta^{op} \rightarrow \calC^{\otimes}$ is a general section, then we can view $A = \widetilde{A}([1]) \in \calC$ as our candidate for the ``underlying algebra''. For $n \geq 0$, we can view $\widetilde{A}([n])$ as an $n$-tuple of objects $( C_1, \ldots, C_n)$ in $\calC$. The inclusions $[1] \simeq \{ i-1, i \} \subseteq [n]$ induce
maps $\eta_{i}: C_{i} \rightarrow A$ for $1 \leq i \leq n$. It is not difficult to see that $\widetilde{A}$ is equivalent to the section arising from an algebra object (well-defined up to isomorphism) if and only if each $\eta_{i}$ is an isomorphism. This observation, together with Remark \ref{labus}, motivates the following definition:

\begin{definition}\label{suskin}\index{algebra object}\index{ZZZAlgC@$\Alg(\calC)$}
Let $\calC$ be a monoidal $\infty$-category. An {\it algebra object} of $\calC$ is a lax monoidal functor $\Nerve(\cDelta)^{op} \rightarrow \calC^{\otimes}$. We let
$\Alg(\calC) = \Fun^{\wMon}( \Nerve(\cDelta)^{op}, \calC^{\otimes})$ denote the $\infty$-category of algebra objects of $\calC$.
\end{definition}

More concretely, an algebra object of $p: \calC^{\otimes} \rightarrow \Nerve(\cDelta)^{op}$ consists of a section of $p$ which carries every convex morphism in $\cDelta$ to a $p$-coCartesian morphism in $\calC^{\otimes}$.

\begin{remark}
The notation of Definition \ref{suskin} is somewhat abusive: the $\infty$-category
$\Alg(\calC)$ depends not only on $\calC$, but also on its monoidal structure. There is little risk of confusion; the monoidal structure under consideration should always be clear from context.
\end{remark}

\begin{remark}
Let $\calC$ be a monoidal category, and let
$\calC^{\otimes}$ be as in Definition \ref{converi}. Then $\Nerve( \calC^{\otimes} )$ can be identified with a monoidal structure on the $\infty$-category $\Nerve(\calC)$. The above discussion shows that $\Alg( \Nerve(\calC) )$ is equivalent to the nerve of the ordinary category of algebra objects of $\calC$, interpreted in the classical sense.
\end{remark}

\begin{remark}
Let $\calC$ be an $\infty$-category equipped with a monoidal structure. Evaluation at $[1] \in \cDelta$ determines a functor $\Alg(\calC) \rightarrow \calC$, which we will refer to as the {\it forgetful functor}. 
\end{remark}

\begin{remark}\label{funtime}\index{monoidal $\infty$-category!of functors}
Let $\calC$ be a monoidal $\infty$-category and $K$ an arbitrary simplicial set. Then
$\Fun(K, \calC)$ inherits the structure of a monoidal $\infty$-category, where the tensor product is defined pointwise. More precisely, suppose that $\calC^{\otimes} \rightarrow \Nerve(\cDelta)^{op}$ exhibits $\calC$ as a monoidal $\infty$-category. Set
$$\Fun(K, \calC)^{\otimes} = \Fun(K, \calC^{\otimes}) \times_{ \Fun(K, \Nerve(\cDelta)^{op})} \Nerve(\cDelta)^{op}.$$
Then the projection $\Fun(K, \calC)^{\otimes} \rightarrow \Nerve(\cDelta)^{op}$ exhibits
$\Fun(K, \calC)^{\otimes}_{[1]} \simeq \Fun(K, \calC)$ as a monoidal $\infty$-category.
Moreover, we have a canonical isomorphism
$\Alg( \Fun(K, \calC) ) \rightarrow \Fun(K, \Alg(\calC) ).$
\end{remark}

\subsection{Cartesian Monoidal Structures}\label{monoidcart}

Let $\calC$ be an ordinary category which admits finite products. Then $\calC$ has the structure of a monoidal category, with the bifunctor $\otimes: \calC \times \calC \rightarrow \calC$ given by the Cartesian product. We will refer to this monoidal structure on $\calC$ as the {\it Cartesian monoidal structure}. Our goal in this section is to give an analogous construction in the $\infty$-categorical setting.

\begin{definition}\label{frieze}\index{monoidal structure!Cartesian}\index{Cartesian!monoidal structure}
Let $\calC$ be an $\infty$-category. We will say that a monoidal structure on $\calC$ is {\it Cartesian} if
the following conditions are satisfied:
\begin{itemize}
\item[$(1)$] The unit object $1_{\calC} \in \calC$ is final.
\item[$(2)$] For every pair of objects $C, D \in \calC$, the canonical maps
$$ C \simeq C \otimes 1_{\calC} \leftarrow C \otimes D \rightarrow 1_{\calC} \otimes D \simeq D$$
exhibit $C \otimes D$ as a product of $C$ and $D$ in the $\infty$-category $\calC$.
\end{itemize}
\end{definition}

If $\calC^{\otimes} \rightarrow \Nerve(\cDelta)^{op}$ is a Cartesian monoidal structure on
an $\infty$-category $\calC = \calC^{\otimes}_{[1]}$, then we can construct a functor $\pi: \calC^{\otimes} \rightarrow \calC$, which is given informally as follows.
To an object $C \in \calC^{\otimes}_{[n]}$, corresponding to an
$n$-tuple $(C_1, \ldots, C_n) \in \calC^n$, the functor $\pi$ associates the object
$\pi(C) = \prod_{1 \leq i \leq n} C_i$. We will give rigorous construction of $\pi$ below
(Proposition \ref{cardbade}); first, we axiomatize its properties.

\begin{definition}\label{puggybear}\index{Cartesian structure}\index{Cartesian structure!weak}\index{Cartesian structure!lax}\index{lax!Cartesian structure}
Let $p: \calC^{\otimes} \rightarrow \Nerve(\cDelta)^{op}$ be a monoidal $\infty$-category.
A {\it lax Cartesian structure} on $\calC^{\otimes}$ is a functor
$\pi: \calC^{\otimes} \rightarrow \calD$ satisfying the following condition:
\begin{itemize}
\item[$(\ast)$] Let $C$ be an object of $\calC^{\otimes}_{[n]}$. For each $1 \leq i \leq n$,
choose a $p$-coCartesian morphism $f_{i}: C \rightarrow C_{i}$ covering the inclusion
$[1] \simeq \{ i-1, i \} \subseteq [n]$ in $\cDelta$.
Then the morphisms
$\pi(f_i)$ exhibit $\pi(C)$ as a product
$\prod_{1 \leq i \leq n} \pi(C_i)$ in the $\infty$-category $\calD$.
\end{itemize}
We will say that $\pi$ is a {\it weak Cartesian structure} if it is a lax Cartesian structure, and
the following additional condition is satisfied:
\begin{itemize}
\item[$(\ast')$] Let $f: C \rightarrow C'$ be a $p$-coCartesian morphism covering the map $[1] \simeq \{0, n\} \rightarrow [n]$. Then $\pi(f)$ is an equivalence in $\calD$.
\end{itemize}
We will say that a weak Cartesian structure $\pi$ is a {\it Cartesian structure} if $\pi$ induces an equivalence $\calC^{\otimes}_{[1]} \rightarrow \calD$. 
\end{definition}

\begin{example}
Let $\calC$ be an ordinary category which admits finite products. Regard $\calC$ as endowed with the Cartesian monoidal structure, and let $\calC^{\otimes}$ be as in Definition \ref{converi}. Define $\theta: \calC^{\otimes} \rightarrow \calC$ by the formula
$\theta( [C_1, \ldots, C_n]) = C_1 \times \ldots C_n$. The nerve $\Nerve(\theta): \Nerve(\calC^{\otimes}) \rightarrow \Nerve(\calC)$ is a Cartesian structure on the monoidal
$\infty$-category $\Nerve(\calC^{\otimes})$.
\end{example}

It follows immediately from the definition that if $\calC$ is a monoidal $\infty$-category and there exists a Cartesian structure $\calC^{\otimes} \rightarrow \calD$, then the monoidal structure on $\calC$ is Cartesian. Our first result is a converse: if $\calC$ is a Cartesian monoidal $\infty$-category, then there exists an essentially unique Cartesian structure on $\calC$.

\begin{proposition}\label{cardbade}
Let $p: \calC^{\otimes} \rightarrow \Nerve(\cDelta)^{op}$ be a Cartesian monoidal structure on an
$\infty$-category $\calC$ and let $\calD$ be an $\infty$-category which admits finite products.
Let $\Fun^{\times}( \calC^{\otimes}, \calD)$ denote the full subcategory of $\Fun( \calC^{\otimes}, \calD)$ spanned by the weak Cartesian structures, and let $\Fun^{\times}(\calC, \calD)$ be the full subcategory of $\Fun(\calC, \calD)$ spanned by those functors which preserve finite products. The restriction map
$ \Fun^{\times}( \calC^{\otimes}, \calD) \rightarrow \Fun^{\times}(\calC, \calD)$ is an equivalence
of $\infty$-categories.
\end{proposition}

We will defer the proof until the end of this section.

Our next goal is to show that, if $\calC$ is an $\infty$-category which admits finite products, then there exists an (essentially unique) Cartesian monoidal structure on $\calC$. Our strategy is to give an explicit construction of this monoidal structure.

\begin{notation}\label{hugegrin}\index{ZZZcDeltatimes@$\cDelta^{\times}$}
The category $\cDelta^{\times}$ is defined as follows:
\begin{itemize}
\item[$(1)$] An object of $\cDelta^{\times}$ consists of an object $[n] \in \cDelta$ together
with a pair of integers $i$ and $j$ which satisfy $0 \leq i \leq j \leq n$.

\item[$(2)$] A morphism from $([n], i \leq j)$ to $([n'], i' \leq j')$ in $\cDelta^{\times}$
is a map $f: [n] \rightarrow [n']$ of linearly ordered sets, with the property that
$i' \leq f(i) \leq f(j) \leq j'$.
\end{itemize}
\end{notation}

The forgetful functor $\cDelta^{\times,op} \rightarrow \cDelta^{op}$ is a Grothendieck fibration, so that the induced map of $\infty$-categories $\Nerve(\cDelta^{\times})^{op} \rightarrow \Nerve(\cDelta)^{op}$ is a Cartesian fibration (Remark \toposref{gcart}).

\begin{remark}\label{hungerfor}
The forgetful functor $\cDelta^{\times} \rightarrow \cDelta$ admits a {\em unique} section $s$, defined by $s( [n] ) = ( [n], 0 \leq n )$. 
\end{remark}

\begin{notation}\index{ZZZcalCtimes@$\calC^{\times}$}
Let $\calC$ be an $\infty$-category. We define a simplicial set $\widetilde{\calC}^{\times}$ equipped with a map $\widetilde{\calC}^{\times} \rightarrow \Nerve(\cDelta)^{op}$ by the following universal property: for every map of simplicial sets $K \rightarrow \Nerve(\cDelta)^{op}$, we have a bijection
$$ \Hom_{ \Nerve(\cDelta)^{op} }(K, \widetilde{\calC}^{\times} )
\simeq \Hom_{ \sSet }( K \times_{ \Nerve(\cDelta)^{op} } \Nerve(\cDelta^{\times})^{op}, \calC).$$

Fix $n \geq 0$. We observe that the fiber $\widetilde{\calC}^{\times}_{[n]}$ can be identified with the
$\infty$-category of functors $f: \Nerve(P)^{op} \rightarrow \calC$, where $P$ is the partially ordered set of {\em intervals} in $[n]$: that is, the collection of all subsets of $[n]$ having the form
$\{ i, i+1, \ldots, j \} \subseteq [n]$. We let $\calC^{\times}$ be the full simplicial subset
of $\widetilde{\calC}^{\times}$ spanned by those vertices which correspond to those functors $f$
for which the maps $f( \{ i, i+1, \ldots, j \} ) \rightarrow f( \{k,k+1\} )$ exhibit $f( \{i, \ldots, j\} )$ as a product
$f( \{i,i+1\} ) \times \ldots \times f( \{j-1, j\})$. 
\end{notation}

\begin{proposition}\label{commonton}
Let $\calC$ be an $\infty$-category.
\begin{itemize}
\item[$(1)$] The projection $p: \widetilde{\calC}^{\times} \rightarrow \Nerve(\cDelta)^{op}$ is a coCartesian fibration.
\item[$(2)$] Let $\alpha: F \rightarrow F'$ be a morphism of $\widetilde{\calC}^{\times}$ whose image in $\Nerve(\cDelta)^{op}$ corresponds to a map $s: [n] \rightarrow [m]$. Then $\alpha$ is
$p$-coCartesian if and only if the induced map $F( \{ s(i), \ldots, s(j) \}) \rightarrow
F'( \{i, \ldots, j \})$ is an equivalence in $\calC$, for every $0 \leq i \leq j \leq m$.
\item[$(3)$] The projection $p$ restricts to a coCartesian fibration $\calC^{\times} \rightarrow \Nerve(\cDelta)^{op}$ $($with the same class of $p$-coCartesian morphisms$)$.
\item[$(4)$] The projection $\calC^{\times} \rightarrow \Nerve(\cDelta)^{op}$ is a monoidal $\infty$-category if and only if $\calC$ admits finite products.
\item[$(5)$] Suppose that $\calC$ admits finite products. 
Let $\pi: \calC^{\times} \rightarrow \calC$ be the map given by composition with
the section $\Nerve(\cDelta)^{op} \rightarrow \Nerve( \cDelta^{\times})^{op}$ $($see Remark \ref{hungerfor}$)$. Then $\pi$ is a Cartesian structure on $\calC^{\times}$. 
\end{itemize}
\end{proposition}

\begin{proof}
Assertions $(1)$ and $(2)$ follow immediately from Corollary \toposref{skinnysalad}, and $(3)$ follows from $(2)$ (since $\calC^{\times}$ is stable under the pushforward functors associated
to the coCartesian fibration $p$). 

We now prove $(4)$. If $\calC$ has no final object, then $\calC^{\times}_{[0]}$ is empty; consequently, we may assume without loss of generality that $\calC$ has a final object.
Then $\calC^{\times}_{[1]}$ is equivalent to the $\infty$-category of diagrams
$X_0 \leftarrow X \rightarrow X_1$ in $\calC$, where $X_0$ and $X_1$ are final. It follows that $\pi$ induces an equivalence $\calC^{\times}_{[1]} \simeq \calC$. Consequently,
$\calC^{\times}$ is a monoidal $\infty$-category if and only if, for each $n \geq 0$, the
natural map $\phi: \calC^{\times}_{[n]} \rightarrow \calC^{\times}_{ \{0,1\} } \times \ldots 
\times \calC^{\times}_{ \{n-1,n\} }$
is an equivalence of $\infty$-categories. Let $P$ denote the partially ordered set
of subintervals of $[n]$, and let $P_0 \subseteq P$ denote the subset consisting of the intervals $\{ k, k+1\}$, where $0 \leq k < n$. Then $\calC^{\times}_{[n]}$ can be identified with the
set of functors $F: \Nerve(P)^{op} \rightarrow \calC$ which are right Kan extensions of
$F| \Nerve(P_0)^{op}$, and $\phi$ coincides with the restriction map from $P$ to $P_0$. According to Proposition \toposref{lklk}, $\phi$ is fully faithful, and is essentially surjective if and only if
every functor $F_0: \Nerve(P_0)^{op} \rightarrow \calC$ admits a right Kan extension to
$\Nerve(P)^{op}$. Unwinding the definitions, we see that this is equivalent to the assertion that every finite collection of objects of $\calC$ admits a product in $\calC$. This completes the proof of $(4)$. Assertion $(5)$ follows immediately from the construction of $\calC^{\times}$ (and the description of the $p$-coCartesian morphisms supplied by Corollary \toposref{skinnysalad}).
\end{proof}
 
\begin{proposition}\label{sungto}
Let $p: \calC^{\otimes} \rightarrow \Nerve(\cDelta)^{op}$ be a monoidal $\infty$-category and let $\calD$ be an $\infty$-category which admits finite products. Let $\Fun^{\times, \wMon}( \calC^{\otimes}, \calD)$ denote the full subcategory of $\Fun(\calC^{\otimes}, \calD)$ spanned by the lax Cartesian structures and $\Fun^{\times}(\calC^{\otimes}, \calD) \subseteq \Fun^{\times, \wMon}(\calC^{\otimes},\calD)$ the full subcategory spanned by the weak Cartesian structures. 
Let $\pi: \calD^{\times} \rightarrow \calD$ be the Cartesian structure of Proposition \ref{commonton}. Then composition with $\pi$ induces trivial Kan fibrations
$$ \theta: \Fun^{\wMon}( \calC^{\otimes}, \calD^{\times}) \rightarrow \Fun^{\times,\wMon}( \calC^{\otimes}, \calD) \quad \theta_0: \Fun^{\Mon}(\calC^{\otimes}, \calD^{\times}) \rightarrow \Fun^{\times}( \calC^{\otimes}, \calD).$$
\end{proposition}

\begin{proof}
Unwinding the definitions, we can identify 
$\Fun^{\wMon}( \calC^{\otimes}, \calD^{\times} )$ with the full subcategory of
$$\bHom( \calC^{\otimes} \times_{ \Nerve(\cDelta)^{op} } \Nerve(\cDelta^{\times})^{op}, \calD)$$
spanned by those functors $F$ which satisfy the following conditions:
\begin{itemize}
\item[$(1)$] For every $0 \leq i \leq j \leq n$ and every $C_{[n]} \in \calC^{\otimes}_{[n]}$, $F$ induces an equivalence
$$F(C_{[n]}, i \leq j ) \rightarrow \prod_{0 \leq i < n} F(C_{[n]}, i \leq i+1)$$
in the $\infty$-category $\calD$.
\item[$(2)$] For every $p$-coCartesian morphism $C_{[n]} \rightarrow C_{[m]}$ covering a convex morphism $\alpha: [m] \rightarrow [n]$ in $\cDelta$, and every $0 \leq i \leq j \leq m$, the induced map
$F( C_{[n]}, \alpha(i) \leq \alpha(j) ) \rightarrow F( C_{[m]}, i \leq j)$ is an equivalence in $\calD$.
\end{itemize} 
The functor $F' = \pi \circ F$ can be described by the formula $F'( C_{[n]} ) = F( C_{[n]}, 0 \leq n )$. In other words, $F'$ can be identified with the restriction of $F$ to the full subcategory of
$\calC(0) \subseteq \calC^{\otimes} \times_{ \Nerve(\cDelta)^{op} } \Nerve(\cDelta^{\times})^{op}$ spanned by objects of the form $(C_{[n]}, 0 \leq n)$. 

We observe that for every object $X = ( C_{[n]}, i \leq j)$ of the fiber product $\calC^{\otimes} \times_{ \Nerve(\cDelta)^{op} } \Nerve(\cDelta^{\times})^{op}$, the $\infty$-category
$\calC(0)_{X/}$ has an initial object. More precisely, if we choose a $p$-coCartesian morphism
$\overline{\alpha}: C_{[n]} \rightarrow C_{\{ i, \ldots, j\} }$ lifting the inclusion $\alpha: \{i, \ldots, j\} \subseteq [n]$, then the induced map $\overline{\alpha}: (C_{[n]}, i \leq j) \rightarrow ( C_{\{ i, \ldots, j \}}, i \leq j )$
is an initial object of $\calC(0)_{X/}$. It follows that every functor $F': \calC(0) \rightarrow \calD$
admits a right Kan extension to $\calC^{\otimes} \times_{ \Nerve(\cDelta)^{op} } \Nerve(\cDelta^{\times})^{op}$, and that an arbitrary functor
$F: \calC^{\otimes} \times_{ \Nerve(\cDelta)^{op} } \Nerve(\cDelta^{\times})^{op} \rightarrow
\calD$ is a right Kan extension of $F| \calC(0)$ if and only if
$F( \overline{\alpha})$ is an equivalence, for every $\overline{\alpha}$ defined as above.

Let $\calE$ be the full subcategory of $\Fun( \calC^{\otimes} \times_{ \Nerve(\cDelta)^{op} } \Nerve(\cDelta^{\times})^{op}, \calD)$ spanned by those functors $F$ which satisfy the following conditions:
\begin{itemize}
\item[$(1')$] The restriction $F' = F | \calC(0)$ is a lax Cartesian structure on $\calC^{\otimes} \simeq \calC(0)$. 
\item[$(2')$] The functor $F$ is a right Kan extension of $F'$. 
\end{itemize}
Using Proposition \toposref{lklk}, we conclude that the restriction map
$\calE \rightarrow \Fun^{\times, \wMon}( \calC^{\otimes}, \calD)$ is a trivial fibration of simplicial sets.
To prove that $\theta$ is a trivial Kan fibration, it will suffice to show that conditions $(1)$ and $(2)$ are equivalent to conditions $(1')$ and $(2')$. 

Suppose first that $(1')$ and $(2')$ are satisfied by a functor $F$. Condition then $(1)$ follows easily. Choose a map $C_{[n]} \rightarrow C_{[m]}$, covering a convex morphism $[m] \rightarrow [n]$ in $\cDelta$, and let $0 \leq i \leq j \leq m$ be as in the statement of $(2)$. We have a homotopy commutative diagram
$$ \xymatrix{ F( C_{[n]}, \alpha(i) \leq \alpha(j) ) \ar[r] \ar[d] & F( C_{[m]}, i \leq j) \ar[d] \\
F( C_{\{ \alpha(i), \ldots ,  \alpha(j)\} }, \alpha(i) \leq \alpha(j) ) \ar@{=}[r] & F( C_{\{ i, \ldots, j \} }, i \leq j ).}$$
Condition $(2')$ implies that the vertical maps are equivalences, and the lower horizontal map is the identity. It follows that the upper horizontal map is an equivalence, which proves $(2)$.

Now suppose that $(1)$ and $(2)$ are satisfied. The implication $(2) \Rightarrow (2')$ is obvious; it will therefore suffice to verify $(1')$. Let $C_{[n]}$ be an object of $\calC^{\otimes}_{[n]}$, and choose $p$-coCartesian morphisms $g_i: C_{[n]} \rightarrow C_{ \{i, i+1\} }$.
We wish to show that the induced map
$$ F( C_{[n]}, 0 \leq n ) \rightarrow
F( C_{ \{0,1\}}, 0 \leq 1) \times \ldots \times F( C_{ \{n-1, n\}}, n-1 \leq n )$$
is an equivalence, which follows immediately from $(1)$ and $(2')$. This completes the proof that $\theta$ is a trivial Kan fibration.

To prove that $\theta_0$ is a trivial Kan fibration, it will suffice to prove that $\theta_0$ is a pullback of $\theta$. In other words, it will suffice to show that if 
$F: \calC^{\otimes} \times_{ \Nerve(\cDelta)^{op} } \Nerve(\cDelta^{\times})^{op} \rightarrow
\calD$ is a functor satisfying conditions $(1)$ and $(2)$, then $F| \calC(0)$ is a weak Cartesian structure on $\calC^{\otimes}$ if and only if $F$ determines a monoidal functor from
$\calC^{\otimes}$ into $\calD^{\times}$. Let $q: \calD^{\times} \rightarrow \Nerve(\cDelta)^{op}$ denote the projection map. Using the description of the $q$-coCartesian morphisms provided by Proposition \ref{commonton}, we see that the latter condition is equivalent to

\begin{itemize}
\item[$(2_{+})$] For every $p$-coCartesian morphism $\overline{\alpha}: C_{[n]} \rightarrow C_{[m]}$ covering a map $\alpha: [m] \rightarrow [n]$ in $\cDelta$, and every $0 \leq i \leq j \leq m$, the induced map
$F( C_{[n]}, f(i) \leq f(j) ) \rightarrow F( C_{[m]}, i \leq j)$ is an equivalence in $\calD$.
\end{itemize}

Moreover, $F| \calC(0)$ is a weak Cartesian structure if and only if $F$ satisfies the following:

\begin{itemize}
\item[$(3')$] For every $n \geq 0$ and every $p$-coCartesian morphism
$\overline{\alpha}: C_{[n]} \rightarrow C_{[1]}$ in $\calC^{\otimes}$ lifting the map
$[1] \simeq \{0, n\} \subseteq [n]$, the induced map
$F( C_{[n]}, 0 \leq n) \rightarrow F( C_[1], 0 \leq 1)$ is an equivalence in $\calD$.
\end{itemize}

It is clear that $(2_{+})$ implies $(3')$. Conversely, suppose that $(3')$ is satisfied, and let
$\overline{\alpha}$ and $0 \leq i \leq j \leq m$ be as in the statement of $(2_+)$. Consider the diagram
$$ \xymatrix{ F( C_{[n]}, \alpha(i) \leq \alpha(j) ) \ar[r] \ar[d] & F( C_{[m]}, i \leq j) \ar[d] \\
F( C_{\{ \alpha(i), \ldots ,  \alpha(j)\} }, \alpha(i) \leq \alpha(j) ) \ar[r] & F( C_{\{ i, \ldots, j \} }, i \leq j ).}$$
Condition $(2')$ implies that the vertical maps are equivalences, and two applications
of $(3')$ implies that the lower horizontal map is an equivalence as well. It follows that the upper horizontal map is also an equivalence, as desired.
\end{proof}

We are now in a position to establish the uniqueness of Cartesian monoidal structures. Let
$\Cat_{\infty}^{\Mon, \times}$ denote the full subcategory of $\Cat_{\infty}^{\Mon}$ spanned by the Cartesian monoidal $\infty$-categories. Let $\Cat_{\infty}^{\Cart}$ denote the subcategory
of $\Cat_{\infty}$ whose objects are $\infty$-categories which admit finite products, and whose morphisms are functors which preserve finite products. 

\begin{corollary}
The forgetful functor
$\theta: \Cat_{\infty}^{\Mon, \times} \rightarrow \Cat_{\infty}^{\Cart}$ is an equivalence of $\infty$-categories.
\end{corollary}

\begin{proof}
We first observe that if $\calC$ is an $\infty$-category which admits finite products, then
Proposition \ref{commonton} implies that the monoidal $\infty$-category $\calC^{\times}$ is a preimage of $\calC$ under the forgetful functor $\theta$. It follows that $\theta$ is essentially surjective.
Moreover, if $\calC^{\otimes}$ is {\em any} Cartesian monoidal structure on $\calC$, then
Proposition \ref{cardbade} guarantees the existence of a Cartesian structure
$\pi: \calC^{\otimes} \rightarrow \calC$. Applying Proposition \ref{sungto}, we can lift
$\pi$ to a monoidal functor $F: \calC^{\otimes} \rightarrow \calC^{\times}$. Remark \ref{eggsal} implies that $F$ is an equivalence. 

We now show that $\theta$ is fully faithful. Let $\calC^{\otimes}$ and $\calD^{\otimes}$ be Cartesian monoidal structures on $\infty$-categories $\calC$ and $\calD$. We wish to show that the restriction map
$$ \bHom_{ \Cat_{\infty}^{\Mon, \times} }( \calC^{\otimes}, \calD^{\otimes} )
\rightarrow \bHom_{ \Cat_{\infty}^{\Cart} }( \calC, \calD )$$
is a homotopy equivalence. We will prove a slightly stronger assertion: namely, that the restrictio map
$\psi: \Fun^{\Mon}( \calC^{\otimes}, \calD^{\otimes} ) \rightarrow \Fun^{\times}(\calC, \calD)$
is a categorical equivalence, where $\Fun^{\times}(\calC, \calD)$ denotes the full subcategory of
$\Fun(\calC, \calD)$ spanned by those functors which preserve finite products. In view of the above remarks, it suffices to prove this in the case where $\calD^{\otimes} = \calD^{\times}$. 
In this case, the map $\psi$ factors as a composition
$$ \Fun^{\Mon}(\calC^{\otimes}, \calD^{\otimes} )
\stackrel{\psi'}{\rightarrow} \Fun^{\times}( \calC^{\otimes}, \calD)
\stackrel{\psi''}{\rightarrow} \Fun^{\times}(\calC, \calD).$$
Proposition \ref{sungto} implies that $\psi'$ is a categorical equivalence, and 
Proposition \ref{cardbade} implies that $\psi''$ is a categorical equivalence.
\end{proof}

Our next goal is to study the algebra objects in an $\infty$-category equipped with a Cartesian monoidal structure. We begin by reviewing a bit of classical category theory.
Let $\calC$ be an ordinary category which admits finite products. A {\it monoid object} of
$\calC$ is an object $M \in \calC$, equipped with maps\index{monoid object!of an ordinary category}
$$ \ast \rightarrow M,  \quad M \times M \rightarrow M$$
which satisfy the usual associativity and unit conditions; here $\ast$ denotes a final object of $\calC$. Equivalently, we can define a monoid object of $\calC$ to be a contravariant functor
from $\calC$ to the category of monoids, such that the underlying functor $\calC^{op} \rightarrow \Set$ is representable by an object $M \in \calC$.

\begin{example}\label{yikyik}
If $\calC$ is the category of sets, then a monoid object of $\calC$ is simply a monoid $M$. We can identify the monoid $M$ with a category $\calD_M$, having only a single object $E$ with
$\Hom_{\calD_M}(E,E) = M$. The nerve $\Nerve( \calD_M )$ is a simplicial set, which is typically denoted by $BM$ and called the {\it classifying space} of $M$. Concretely, the set of $n$-simplices of $BM$ can be identified with an $n$-fold product of $M$ with itself, and the face and degeneracy operations on $BM$ encode the multiplication and unit operations on $M$. The functor
$M \mapsto BM$ is a fully faithful embedding of the category of monoids into the category of simplicial sets. Moreover, a simplicial set $X$ is isomorphic to the classifying space of a monoid if and only if, for each $n \geq 0$, the natural map
$X( [n]) \rightarrow X( \{ 0,1 \}) \times \ldots \times X( \{n-1,n\})$
is a bijection. In this case, the underlying monoid is given by $X([1])$, with unit determined by the degeneracy map $\ast \simeq X([0]) \rightarrow X([1])$ and multiplication by the face map $X([1]) \times X([1]) \simeq X([2]) \stackrel{d_1}{\rightarrow} X([1])$. 
\end{example}

It follows from Example \ref{yikyik} that we can identify monoids in an arbitrary category $\calC$ with certain simplicial objects of $\calC$. This observation allows us to generalize the notion of a monoid to higher category theory.

\begin{definition}\label{monob}\index{monoid object!of an $\infty$-category}\index{ZZZMonC@$\Mon(\calC)$}
Let $\calC$ be an $\infty$-category. A {\it monoid object} of $\calC$ is a simplicial object
$X: \Nerve(\cDelta)^{op} \rightarrow \calC$ with the property that, for each $n \geq 0$, the collection of maps $X( [n] ) \rightarrow X( \{i, i+1\} )$ exhibits $X([n])$ as a product
$X( \{ 0,1 \}) \times \ldots \times X( \{n-1,n\})$. We let $\Mon(\calC)$ denote the full subcategory of
$\Fun( \Nerve(\cDelta)^{op}, \calC)$ spanned by the monoid objects of $\calC$.
\end{definition}

\begin{example}
Let $\calC$ be an $\infty$-category. Every group object of $\calC$ (see Definition \toposref{gropab}) is a monoid object of $\calC$. 
\end{example}

\begin{proposition}\label{ungbat}
Let $\calC^{\otimes} \rightarrow \Nerve(\cDelta)^{op}$ be a monoidal structure on an $\infty$-category $\calC = \calC^{\otimes}_{[1]}$, and let $q: \calC^{\otimes} \rightarrow \calD$
be a Cartesian structure. Then composition with $q$ induces an equivalence
$\Alg(\calC) \rightarrow \Mon(\calD).$
\end{proposition}

\begin{proof}
In view of Proposition \ref{sungto} and Remark \ref{eggsal}, we may assume without loss of generality that
$\calC^{\otimes} = \calD^{\times}$. We now apply Proposition \ref{sungto} again to deduce that the map $$\Alg(\calC) = \Fun^{\wMon}( \Nerve(\cDelta)^{op}, \calC^{\otimes}) \rightarrow
\Fun^{\times, \wMon}( \Nerve(\cDelta)^{op}, \calD) = \Mon(\calD)$$
is a trivial Kan fibration.
\end{proof}

\begin{remark}\label{otherlander}\index{monoidal $\infty$-category}\index{$\infty$-category!monoidal}
Definition \ref{monob} allows us to reformulate the notion of a monoidal $\infty$-category. More precisely, we claim that a monoidal $\infty$-category is essentially the same thing as a monoid object in the $\infty$-category $\Cat_{\infty}$. To see this, we let $(\mSet)_{/ \Nerve(\cDelta)^{op} }$ denote the category of marked simplicial sets
equipped with a map to $\Nerve(\cDelta)^{op}$ (see \S \toposref{twuf}), endowed with the {\em opposite} of the model structure defined in \S \toposref{markmodel} (so that the fibrant objects
of $(\mSet)_{/ \Nerve(\cDelta)^{op} }$ can be identified with coCartesian fibrations
$X \rightarrow \Nerve(\cDelta)^{op}$. Then $\Cat_{\infty}^{\Mon}$ can be identified with a full subcategory of the underlying $\infty$-category $\Nerve( (\mSet)_{/ \Nerve(\cDelta)^{op}}^{\degree})$.

Theorem \toposref{straightthm} and Proposition \toposref{gumby4} furnish equivalences
of $\infty$-categories
$$ \Fun( \Nerve(\cDelta)^{op}, \Cat_{\infty}) \leftarrow \Nerve( \Fun( \cDelta^{op}, \mSet)^{\degree}) \rightarrow \Nerve( (\mSet)_{/ \Nerve(\cDelta)^{op}}^{\degree}).$$
Under the composite equivalence, the full subcategory
$\Cat_{\infty}^{\Mon} \subseteq \Nerve( (\mSet)_{/ \Nerve(\cDelta)^{op}}^{\degree})$
can be identified with $\Mon( \Cat_{\infty}) \subseteq \Fun( \Nerve(\cDelta)^{op}, \Cat_{\infty})$. 
In other words, we may identify monoidal $\infty$-categories with monoid objects in the
$\infty$-category $\Cat_{\infty}$. 

The above argument suggests yet another possible definition: we can identify monoidal $\infty$-categories with certain functors $F$ from $\cDelta^{op}$ into the {\em ordinary} category of $\infty$-categories, which have the property that the induced map
$F([n]) \rightarrow F( \{0,1\} ) \times \ldots \times F( \{n-1, n\})$
is a categorical equivalence for each $n \geq 0$. Such functors can be identified with
{\em bisimplicial sets} satisfying appropriate extension conditions; we leave the details to the reader. We will later obtain an even more concrete model for the $\infty$-category $\Mon(\Cat_{\infty})$: namely, monoidal $\infty$-categories can be identified with $\infty$-categories $\calC$ equipped with an object $1 \in \calC$ and a multiplication $\otimes: \calC \times \calC \rightarrow \calC$ which
is strictly associative (Example \ref{exalcun}).

Although these alternative approaches are perhaps more concrete, they are more difficult to use in practice. For example, the theory of algebra objects in a monoidal $\infty$-category is most easily formulated in terms of Definition \ref{mainef}.
\end{remark}

\begin{remark}\label{selfop}\index{monoidal $\infty$-category!opposite of}
The definition of a monoidal $\infty$-category is not manifestly self-dual. However, it is neverthless true that any monoidal structure on an $\infty$-category $\calC$ determines a monoidal structure on $\calC^{op}$, which is unique up to contractible ambiguity. Roughly speaking, we can use 
Remark \ref{otherlander} to identify a monoidal $\infty$-category $\calC$ with a monoid object
$\Nerve(\cDelta)^{op} \rightarrow \Cat_{\infty}$. We would then like to obtain a new monoid object by composing with an involution of $\Cat_{\infty}$ which carries each $\infty$-category to its opposite.

To carry out the details in practice, it is convenient to replace $\Cat_{\infty}$ by an equivalent
$\infty$-category with a slightly more elaborate definition. Recall that $\Cat_{\infty}$ is defined to be the simplicial nerve of a simplicial category $\Cat_{\infty}^{\Delta}$, whose objects are $\infty$-categories, where $\bHom_{\Cat_{\infty}^{\Delta}}(X,Y)$ is the largest Kan complex contained in
$\Fun(X,Y)$. The construction $X \mapsto X^{op}$ does not induce a simplicial functor from
$\Cat_{\infty}^{\Delta}$ to itself; instead we have a canonical isomorphism
$\bHom_{\Cat_{\infty}^{\Delta}}( X^{op}, Y^{op} )
\simeq \bHom_{\Cat_{\infty}^{\Delta}}(X,Y)^{op}.$
However, if we let $\Cat_{\infty}^{\top}$ denote the topological category obtained by geometrically realizing the morphism spaces in 
$\Cat_{\infty}^{\Delta}$, then $i$ induces an autoequivalence of $\Cat_{\infty}^{\top}$ as a topological category (via the natural homeomorphisms $| X | \simeq |X^{op}|$, which is defined for every simplicial set $X$). We now define $\Cat'_{\infty}$ to be the topological nerve of $\Cat_{\infty}^{\top}$ (see 
Definition \toposref{topnerve}). Then $\Cat'_{\infty}$ is an $\infty$-category equipped with a canonical equivalence $\Cat_{\infty} \rightarrow \Cat'_{\infty}$, and the involution $i$ induces
an involution of $\Cat'_{\infty}$, which carries each object $\calC \in \Cat'_{\infty}$ to 
the opposite $\infty$-category $\calC^{op}$.

Remark \ref{otherlander} shows that the theory of monoidal $\infty$-categories
is equivalent to the theory of monoid objects of $\Cat_{\infty}$, which is in turn equivalent
to the theory of monoid objects of $\Cat'_{\infty}$. Composing with the involution $i$ allows us pass between monoidal structures on an $\infty$-category $\calC$ and monoidal structures on the opposite $\infty$-category $\calC^{op}$.
\end{remark}

\begin{remark}\index{coCartesian monoidal structure}
According to Remark \ref{selfop}, any monoidal structure on an $\infty$-category $\calC^{op}$ determines a monoidal structure on $\calC$, up to contractible ambiguity. In particular, if $\calC$ admits finite coproducts, then the Cartesian monoidal structure on $\calC^{op}$ determines a monoidal structure on $\calC$, which we will call the {\it coCartesian monoidal structure}. It is characterized up to equivalence by the following properties:
\begin{itemize}
\item[$(1)$] The unit object $1_{\calC} \in \calC$ is initial.
\item[$(2)$] For every pair of objects $C,D \in \calC$, the canonical maps
$$ C \simeq C \otimes 1_{\calC} \rightarrow C \otimes D
\leftarrow 1_{\calC} \otimes D \simeq D$$
exhibit $C \otimes D$ as a coproduct of $C$ and $D$ in the $\infty$-category $\calC$.
\end{itemize}
There is an analogue of Proposition \ref{ungbat}, which applies in the situation where $\calC$ is equipped with a coCartesian monoidal structure: in this case, the forgetful functor $\theta: \Alg(\calC) \rightarrow \calC$ is a trivial Kan fibration. We will prove an analogue of this statement (for symmetric monoidal $\infty$-categories) in \cite{symmetric}.
\end{remark}

We close this section with the proof of Proposition \ref{cardbade}.

\begin{proof}[Proof of Proposition \ref{cardbade}]
We define a subcategory $\calI \subseteq \cDelta^{op} \times [1]$ as follows:

\begin{itemize}
\item[$(a)$] An object $([n], i) \subseteq \cDelta^{op} \times [1]$ belongs to $\calI$ if and only if
either $i = 1$ or $n > 0$.
\item[$(b)$] A morphism $( [m] , i ) \rightarrow ( [n], j)$ in $[1] \times \cDelta^{op}$ belongs
to $\calI$ if and only if either $j=1$ or the induced map $[n] \rightarrow [m]$ preserves endpoints.
\end{itemize}

Let $\calC'$ denote the fiber product $\calC^{\otimes} \times_{ \Nerve(\cDelta)^{op} } \Nerve(\calI)$, which we regard as a subcategory of $\calC^{\otimes} \times \Delta^1$, and let
$p': \calC' \rightarrow \Nerve(\calI)$ denote the projection.
Let $\calC'_0$ and $\calC'_{1}$ denote the intersections of $\calC'$ with
$\calC^{\otimes} \times \{0\}$ and $\calC^{\otimes} \times \{1\}$, respectively.
We note that there is a canonical isomorphism $\calC'_{1} \simeq \calC^{\otimes}$.

Let $\calE$ denote the full subcategory of $\Fun( \calC', \calD)$ spanned by those functors $F$ which satisfy the following conditions:
\begin{itemize}
\item[$(i)$] For every object $C \in \calC^{\otimes}_{[n]}$, where $n > 0$, the induced map
$F(C,0) \rightarrow F(C,1)$ is an equivalence in $\calD$.
\item[$(ii)$] The restriction $F | \calC'_{1}$ is a weak Cartesian structure on $\calC^{\otimes}$.
\end{itemize}

It is clear that if $(i)$ and $(ii)$ are satisfied, then the restriction 
$F_0 = F | \calC'_0$ satisfies the following additional conditions:

\begin{itemize}
\item[$(iii)$] The restriction $F_0 | \calC^{\otimes}_{[1]} \times \{0\}$ is a functor from $\calC$ to $\calD$ which preserves finite products.

\item[$(iv)$] For every $p'$-coCartesian morphism $\alpha$ in $\calC'_{0}$, the induced map
$F_0(\alpha)$ is an equivalence in $\calD$. 
\end{itemize}

Moreover, $(i)$ is equivalent to the assertion that $F$ is a right Kan extension of $F| \calC'_{1}$. Proposition \toposref{lklk} implies that the restriction map $r: \calE \rightarrow \Fun^{\times}( \calC^{\otimes}, \calD)$ induces a trivial Kan fibration onto its essential image. The map $r$ has a section $s$, given by composition with the projection map
$\calC' \rightarrow \calC^{\otimes}$. The restriction map $\Fun^{\times}( \calC^{\otimes}, \calD) \rightarrow \Fun^{\times}(\calC, \calD)$
factors as a composition
$$ \Fun^{\times}( \calC^{\otimes}, \calD) \stackrel{s}{\rightarrow} \calE \stackrel{e}{\rightarrow} \Fun^{\times}(\calC, \calD),$$
where $e$ is induced by composition with the inclusion
$\calC \subseteq \calC'_0 \subseteq \calC'$. Consequently, it will suffice to prove that $e$ is an equivalence of $\infty$-categories.

Let $\calE_0 \subseteq \Fun( \calC'_0, \calD)$ be the full subcategory spanned by those functors which satisfy conditions $(iii)$ and $(iv)$. The map $e$ factors as a composition 
$$ \calE \stackrel{e'}{\rightarrow} \calE_0 \stackrel{e''}{\rightarrow} \Fun^{\times}(\calC, \calD).$$
Consequently, it will suffice to show that $e'$ and $e''$ are trivial Kan fibrations.

Let $f: \calC'_0 \rightarrow \calD$ be an arbitrary functor, and let
$C \in \calC^{\otimes}_{[n]} \subseteq \calC'_{0}$, where $n > 0$. There exists a unique map
$\alpha: ( [n], 0) \rightarrow ([1], 0)$ in $\calI$; choose a $p'$-coCartesian morphism
$\overline{\alpha}: C \rightarrow C'$ lifting $\alpha$. We observe that $C'$ is
an initial object of $\calC \times ( \calC'_0 )_{/C'} \times_{ \calC'_0} \calC$. Consequently,
$f$ is a right Kan extension of $f | \calC$ at $C$ if and only if $f( \overline{\alpha})$ is an equivalence.
It follows that $f$ satisfies $(iv)$ if and only if $f$ is a right Kan extension of $f | \calC$. The same argument (and Lemma \toposref{kan0}) shows that every functor $f_0: \calC \rightarrow \calD$ admits a right Kan extension to $\calC'_0$. Applying Proposition \toposref{lklk}, we deduce that $e''$ is a trivial Kan fibration.

It remains to show that $e'$ is a trivial Kan fibration. In view of Proposition \toposref{lklk}, it will suffice to prove the following pair of assertions, for every functor $f \in \calE_0$:

\begin{itemize}
\item[$(1)$] There exist a functor $F: \calC' \rightarrow \calD$ which is a left Kan extension of
$f = F | \calC'_0$.
\item[$(2)$] An arbitrary functor $F: \calC' \rightarrow \calD$ which extends $f$ is a
left Kan extension of $f$ if and only if $F$ belongs to $\calE$.
\end{itemize}

Let $(C,1) \in \calC^{\otimes}_{[n]} \times \{1\} \subseteq \calC'$. Since there exists a final object
$1_{\calC} \in \calC$, the $\infty$-category $\calC'_0 \times_{ \calC' } \calC'_{/C}$ also has a final object,
given by the map $\alpha: (C',0) \rightarrow (C,1)$, where $C' \in \calC^{\otimes}_{[0] \star [n] \star [0]}$
corresponds, under the equivalence
$$ \calC^{\otimes}_{ [0] \star [n] \star [0] } \simeq \calC \times \calC^{\otimes}_{[n]} \times \calC,$$
to the triple $(1_{\calC}, C, 1_{\calC})$. We now apply Lemma \toposref{kan2} to deduce $(1)$, together with the following analogue of $(2)$:

\begin{itemize}
\item[$(2')$] An arbitrary functor $F: \calC' \rightarrow \calD$ which extends $f$ is a left
Kan extension of $f$ if and only if, for every morphism $\alpha: (C',0) \rightarrow (C,1)$ as above,
the induced map $F(C',0) \rightarrow F(C,1)$ is an equivalence in $\calD$.
\end{itemize}

To complete the proof, it will suffice to show that $F$ satisfies the conditions stated in $(2')$ if and only if $F \in \calE$. We first prove the ``if'' direction. Let $\alpha: (C', 0) \rightarrow (C,1)$ be as above; we wish to prove that $F(\alpha): F(C',0) \rightarrow F(C,1)$ is an equivalence in $\calD$. If $n=0$, this is clear, since both sides are final objects of $\calD$. Let us therefore assume that $n > 0$. 
The map $\alpha$ factors as a composition
$$ (C',0 ) \stackrel{\alpha'}{\rightarrow} (C', 1) \stackrel{\alpha''}{\rightarrow} (C,1).$$
Condition $(i)$ guarantees that $F(\alpha')$ is an equivalence. Condition $(ii)$ guarantees that
$F( C', 1)$ is equivalent to a product $F(1_{\calC}, 1) \times F(C,1) \times F(1_{\calC},1)$, and that
$F(\alpha'')$ can be identified with the projection onto the middle factor. Moreover, since
$1_{\calC}$ is a final object of $\calC$, condition $(ii)$ also guarantees that
$F(1_{\calC}, 1)$ is a final object of $\calD$. It follows that $F(\alpha'')$ is an equivalence, so that
$F(\alpha)$ is an equivalence as desired.

Now let us suppose that $F$ satisfies the condition stated in $(2')$. We wish to prove that
$F \in \calE$. Here we must invoke our assumption that the monoidal structure on $\calC$ is Cartesian.
We begin by verifying condition $(i)$.
Let $C \in \calC^{\otimes}_{[n]}$, for $n > 0$, and let $\alpha: ( C', 0) \rightarrow (C, 1)$ be
defined as above. There is a unique map $\beta: ( [n], 0 ) \rightarrow 
([0] \star [n] \star [0], 0)$ in $\calI$, which is the identity on $[n] \subseteq [0] \star [n] \star [0]$. Choose
a $p'$-coCartesian morphism $\overline{\beta}: (C,0) \rightarrow (C'', 0)$ lifting $\beta$. Since the final object
$1_{\calC} \in \calC$ is also the unit object of $\calC$, we can identify $C''$ with $C'$. The composition
$ (C, 0) \stackrel{\overline{\beta}}{\rightarrow} (C', 1) \stackrel{\alpha}{\rightarrow} (C,1)$
is homotopic to the canonical map $\gamma: (C,0) \rightarrow (C,1)$ appearing in the statement of $(i)$. Condition $(iv)$ guarantees that $F( \overline{\beta} )$ is an equivalence, and
$(2')$ guarantees that $F(\alpha)$ is an equivalence. Using the two-out-of-three property, we deduce that $F(\gamma)$ is an equivalence, so that $F$ satisfies $(i)$.

To prove that $F$ satisfies $(ii)$, we must verify two conditions:
\begin{itemize}
\item[$(ii_0)$] If $\beta: (C, 1) \rightarrow (D,1)$ is a $p'$-coCartesian morphism
in $\calC'$, and the underlying morphism in $\cDelta$ is endpoint-preserving, then
$F(\beta)$ is an equivalence.
\item[$(ii_1)$] Let $C \in \calC^{\otimes}_{[n]}$, and choose $p$-coCartesian morphisms
$\gamma_i: C \rightarrow C_{i}$ covering the inclusions $\{i-1, i \} \subseteq [n]$, for $1 \leq i \leq n$. Then the maps $\gamma_i$ exhibit $F(C,1)$ as a product
$F( C_1, 1) \times \ldots \times F( C_{n}, 1)$ in $\calD$.
\end{itemize}

Condition $(ii_0)$ follows immediately from $(i)$ and $(iv)$, provided that
$C \in \calC^{\otimes}_{[n]}$ for $n > 0$. If $n = 0$, then the source and target of
$F(\beta)$ are both final objects of $\calD$, so $F(\beta)$ is automatically an equivalence.
To prove $(ii_1)$, we consider the maps $\alpha: (C',0) \rightarrow (C,1)$ and
$\alpha_i: (C'_i, 0) \rightarrow (C_i,1)$ which appear in the statement of $(2')$. We
have a collection of commutative diagrams
$$ \xymatrix{ (C',0) \ar[r]^{\alpha} \ar[d]^{\gamma'_i} & (C,1) \ar[d]^{\gamma_i} \\
(C'_{i}, 0) \ar[r]^{\alpha_i} & (C_{i}, 1). }$$
Condition $(2')$ guarantees that the maps $F(\alpha)$ and 
$F(\alpha_i)$ are equivalences in $\calD$. Consequently, it will suffice to show that
the maps $f(\gamma'_i)$ exhibit $f(C',0)$ as a product of $f( C'_{i}, 0)$ in $\calD$.
Let $f_0 = f | \calC$. Using condition $(iv)$, we obtain canonical equivalences
$$ f(C',0) \simeq f_0( 1_{\calC} \otimes C_1 \otimes \ldots \otimes C_{n} \otimes 1_{\calC} )$$
$$ f(C'_i, 0) \simeq f_0 ( 1_{\calC} \otimes C_i \otimes 1_{\calC} )$$
Since condition $(iii)$ guarantees that $f_0$ preserves products, it will suffice to show that
the canonical map
$$ 1_{\calC} \otimes C_1 \otimes \ldots \otimes C_{n} \otimes 1_{\calC}
\simeq \prod_{1 \leq i \leq n} 1_{\calC} \otimes C_i \otimes 1_{\calC}$$
is an equivalence in the $\infty$-category $\calC$. This follows easily by induction on $n$ (the case $n=2$ reduces to our assumption that the monoidal structure on $\calC$ is Cartesian.)
\end{proof}

\subsection{Subcategories of Monoidal $\infty$-Categories}\label{locol}

Let $\calC$ be an $\infty$-category. Our goal in this section is to show that a monoidal structure on $\calC$ determines a monoidal structure on suitable (full) subcategeroies $\calD \subseteq \calC$. 
Suppose that $\calD$ is a full subcategory of $\calC$, which is stable under equivalence in $\calC$: that is, if $D \in \calD$ and there is an equivalence $C \simeq D$ in $\calC$, then $C \in \calD$. Let $\calC^{\otimes} \rightarrow \cDelta^{op}$ be a monoidal structure on $\calC= \calC^{\otimes}_{[1]}$. We define a full subcategory $\calD^{\otimes} \subseteq \calC^{\otimes}$ as follows: an object $C \in \calC^{\otimes}_{[n]}$ belongs to
$\calD^{\otimes}$ if and only if the image of $C$ under the equivalence
$\calC^{\otimes}_{[n]} \simeq \calC^{n}$ belongs to $\calD^{n}$. Our goal in this section is to study circumstances under which $\calD^{\otimes}$ inherits the structure of a monoidal $\infty$-category.
The most obvious case is that in which $\calD$ is stable under tensor products in $\calC$:

\begin{proposition}\label{yukyik}\index{monoidal $\infty$-category!colocalization of}
\index{monoidal $\infty$-category!subcategory of}
Let $\calC^{\otimes} \rightarrow \Nerve(\cDelta)^{op}$ be a monoidal $\infty$-category, and let $\calD \subseteq \calC$ be a full subcategory which is stable under equivalence, contains the unit object of $\calC$, and is stable under tensor products. 
\begin{itemize}
\item[$(1)$] The restricted map $\calD^{\otimes} \rightarrow \Nerve(\cDelta)^{op}$ is a monoidal $\infty$-category.
\item[$(2)$] The inclusion $\calD^{\otimes} \subseteq \calC^{\otimes}$ is a monoidal functor.
\item[$(3)$] Suppose that the inclusion $\calD \subseteq \calC$ admits a right adjoint $L$ $($so that $\calD$ is a {\it colocalization} of $\calC${}$)$. Then
the inclusion $\calD^{\otimes} \subseteq \calC^{\otimes}$ admits a right adjoint $L^{\otimes}$.
\item[$(4)$] Under the hypothesis of $(3)$, the functor $L^{\otimes}$ is lax monoidal.
\end{itemize}
\end{proposition}

\begin{proof}
Assertions $(1)$ and $(2)$ follow immediately from the definitions. Suppose that $\calD$ is a colocalization of $\calC$. Let us say that a morphism
$\alpha: X' \rightarrow X$ in $\calC^{\otimes}$ is a {\it colocalization} if $X \in \calD^{\otimes}$, and composition with $\alpha$ induces a homotopy equivalence $\bHom_{\calC^{\otimes}}( Y, X') \rightarrow \bHom_{\calC^{\otimes}}(Y,X)$ for all $Y \in \calC^{\otimes}$ for all $Y \in \calD^{\otimes}$. According to Proposition \toposref{testreflect}, $(3)$ is equivalent to the following assertion: 
\begin{itemize}
\item[$(3')$] For every $X \in \calC^{\otimes}$, there exists a colocalization $\alpha: X' \rightarrow X$.
\end{itemize}
Assuming this for the moment, $(4)$ is equivalent to the following pair of assertions:
\begin{itemize}
\item[$(4')$] For every colocalization $\alpha: X' \rightarrow X$, the image of $\alpha$
in $\Nerve(\cDelta)^{op}$ is an equivalence (and therefore degenerate).
\item[$(4'')$] If $\alpha: X' \rightarrow X$ is a colocalization in $\calC^{\otimes}_{[n]}$, and 
$[m] \rightarrow [n]$ is a convex morphism in $\cDelta$, then the image of $\alpha$ in
$\calC^{\otimes}_{[m]}$ is a colocalization.
\end{itemize}

Suppose that $X \in \calC^{\otimes}_{[n]}$. Choose a morphism $\alpha: X' \rightarrow X$ in $\calC^{\otimes}_{[n]}$ which corresponds, under the equivalence
$\calC^{\otimes}_{[n]} \simeq \calC^{n}$, to the canonical map $(LC_1, \ldots, LC_{n}) \rightarrow (C_1, \ldots, C_n)$. We will show that $\alpha$ is a colocalization map in $\calC^{\otimes}$. 
This will prove $(3')$. Assertions $(4')$ and $(4'')$ will likewise follow, since colocalization maps are unique up to equivalence.

Let $Y \in \calD^{\otimes}_{[m]}$ be any object. We wish to show that composition with $\alpha$ induces a homotopy equivalence
$\bHom_{\calC^{\otimes}}( Y, X') \rightarrow \bHom_{\calC^{\otimes}}(Y,X).$
Both sides are given a disjoint union, taken over the set of morphisms $\beta: [n] \rightarrow [m]$ in $\cDelta$. It will therefore suffice to prove the corresponding result for each summand. Using the product structure on $\calC^{\otimes}_{[n]}$, we can reduce to the case where $n=1$ and $\beta$ is endpoint-preserving. We can identify $Y$ with an object
$( D_1, \ldots, D_m ) \in \calD^{m}$, and must show that the canonical map
$$\bHom_{\calC}( D_1 \otimes \ldots \otimes D_m, X') \rightarrow \bHom_{\calC}( D_1 \otimes \ldots \otimes D_m, X)$$
is an equivalence. This follows, since $X' \rightarrow X$ is a colocalization in $\calC$, and
$D_1 \otimes \ldots \otimes D_m$ belongs to $\calD$.
\end{proof}

\begin{example}\label{satline}\index{invertible object}\index{object!invertible}\index{monoidal $\infty$-category!invertible object of}
Let $\calC$ monoidal $\infty$-category. We say that an object $C \in \calC$ is
{\it invertible} if it is an invertible object of the homotopy category $\h{\calC}$: that is, if there exists
an object $D \in \calC$ and equivalences
$$ C \otimes D \simeq 1_{\calC} \simeq D \otimes C.$$
Let $\calC^{0} \subseteq \calC$ denote the full subcategory spanned by the invertible objects. It is easy to see that the hypotheses of Proposition \ref{yukyik} are satisfied, so that $\calC^{0}$ inherits the structure of a monoidal $\infty$-category.
\end{example}

\begin{remark}\label{yuka}
Let $\calD \subseteq \calC$ be as in the statement of Proposition \ref{yukyik}, and suppose that the inclusion $\calD \subseteq \calC$ admits a right adjoint $L$.
The inclusion $i: \calD^{\otimes} \subseteq \calC^{\otimes}$ is a monoidal functor, which induces
a fully faithful embedding $\Alg(\calD) \subseteq \Alg( \calC)$. Let $L^{\otimes}$ be a right adjoint to $i$, so that $L^{\otimes}$ induces a functor $f: \Alg(\calC) \rightarrow \Alg(\calD)$. It is easy to see that $f$ is a right adjoint to the inclusion $\Alg(\calD) \subseteq \Alg(\calC)$. Moreover, if $\theta: \Alg(\calC) \rightarrow \calC$ denotes the forgetful functor, then for each $A \in \Alg(\calC)$, applying $\theta$ to the canonical map $f(A) \rightarrow A$ exhibits $\theta( f(A) )$ as a colocalization of $\theta(A)$.
\end{remark}

Proposition \ref{yukyik} has an evident converse: if $\calD^{\otimes}$ is a monoidal $\infty$-category and the inclusion $i: \calD^{\otimes} \subseteq \calC^{\otimes}$ is a monoidal functor, then
$\calD$ is stable under tensor products in $\calC$. However, it is possible for this converse to fail if $i$ is only assumed to be weakly monoidal. We now discuss a general class of examples where
$\calD$ is {\em not} stable under tensor products, yet $\calD$ nonetheless inherits a monoidal structure from $\calC$.

\begin{definition}\label{compatmon}\index{monoidal $\infty$-category!compatible localization of}
Let $\calC$ be a monoidal $\infty$-category, and let 
$L: \calC \rightarrow \calC$ be a localization functor. We will say that $L$ is {\it compatible with the monoidal structure on $\calC$} if the following condition is satisfied: 
\begin{itemize}
\item[$(\ast)$] Let $f: X \rightarrow Y$ be a morphism in $\calC$, $Z \in \calC$ another object, and
$f': X \otimes Z \rightarrow Y \otimes Z$, $f'': Z \otimes X \rightarrow Z \otimes Y$ the induced morphisms. If $Lf$ is an equivalence, then $Lf'$ and $Lf''$ are equivalences.
\end{itemize}
\end{definition}

\begin{remark}\label{suferr}
In the situation of Definition \ref{compatmon}, it suffices to verify $(\ast)$ in the case where $f: X \rightarrow Y$ is equivalent to the localization map $X \rightarrow LX$.
\end{remark}

\begin{remark}\label{intus}
Suppose that $\calC$ is a {\em presentable} $\infty$-category equipped with a monoidal structure, and that the tensor product $\otimes: \calC \times \calC \rightarrow \calC$ preserves small colimits separately in each variable. It is easy to see that the collection of all maps $f$ which satisfy $(\ast)$ is strongly saturated (see Definition \toposref{saturated2}). Consequently, if $S$ is a collection of morphisms of $\calC$ such that $L \calC = S^{-1} \calC$, then it suffices to verify $(\ast)$ when
$f \in S$.
\end{remark}

\begin{example}\label{jumper}
Let $\calC$ be a presentable $\infty$-category equipped with a monoidal structure, and suppose that the tensor product $\otimes: \calC \times \calC \rightarrow \calC$ preserves small colimits separately in each variable. Let $n \geq -2$ be an integer, and let $\tau_{\leq n} \calC$ be the collection of $n$-truncated objects of $\calC$ (see \S \toposref{truncintro}). Then the localization
functor $\tau_{\leq n}$ is compatible with the monoidal structure on $\calC$. 

To see this, 
let $F: \calC \rightarrow \calC$ be a colimit of the constant map $\bd \Delta^{n+2} \rightarrow \{\id_{\calC} \} \subseteq
\Fun(\calC, \calC)$, let $\alpha: F \rightarrow \id_{\calC}$ be the canonical transformation, and let
$S$ be the collection of all morphisms in $\calC$ equivalent to those of the form
$\alpha(C): F(C) \rightarrow C$. The proof of Proposition \toposref{maketrunc} implies that
$\tau_{\leq n} \calC = S^{-1} \calC$. 

Since the tensor product preserves colimits, we have canonical equivalences
$$ F(C) \otimes D \simeq F( C \otimes D) \simeq C \otimes F(D).$$
It follows that if $f \in S$ and $C \in \calC$, then $f \otimes \id_{C}$ and $\id_{C} \otimes f$ also belong to $S$. Now apply Remark \ref{intus}.
\end{example}

Our next goal is to prove that if $\calC$ is a monoidal $\infty$-category, and $L: \calC \rightarrow \calC$ is a localization functor which is compatible with the monoidal structure on $\calC$, then
the full subcategory $L \calC \subseteq \calC$ inherits the structure of a monoidal category (Proposition \ref{localjerk}). First, we need a lemma.

\begin{lemma}\label{surine}
Let $p: \calC \rightarrow \calD$ be a coCartesian fibration of $\infty$-categories. Let
$L: \calC \rightarrow \calC$ and $L': \calD \rightarrow \calD$ be localization functors, with essential images $L\calC$ and $L'\calD$. Suppose that $L$ and $L'$ are compatible in the following sense:
\begin{itemize}
\item[$(i)$] The functor $p$ restricts to a functor $p': L \calC \rightarrow L' \calD$.
\item[$(ii)$] If $f$ is a morphism in $\calC$ such that $Lf$ is an equivalence, then
$L' p(f)$ is an equivalence in $\calD$.
\end{itemize}
Then:
\begin{itemize}
\item[$(1)$] The functor $L$ carries $p$-coCartesian morphisms of $\calC$ to
$p'$-coCartesian morphisms of $L \calC$.
\item[$(2)$] The functor $p'$ is a coCartesian fibration.
\end{itemize}
\end{lemma}

\begin{proof}
Let $f: X \rightarrow Y$ be a $p$-coCartesian morphism of $\calC$. We wish to prove that $Lf$ is $p'$-coCartesian. According to Proposition \toposref{charCart}, it will suffice to show that for every $Z' \in L\calC$, the diagram of Kan complexes
$$ \xymatrix{ \calC_{Lf/} \times_{\calC} \{Z'\} \ar[r] \ar[d] & \calC_{LX/} \times_{\calC} \{Z'\} \ar[d] \\
\calD_{p(Lf)/} \times_{\calD} \{ p(Z') \} \ar[r] & \calD_{p(LX)/} \times_{\calD} \{p(Z') \} }$$
is homotopy Cartesian. Let $Z' = LZ$. Since $L^2 \simeq L$, we can assume without loss of generality that $Z \in L\calC$. 

Since $f$ is $p$-coCartesian, Proposition \toposref{charCart} implies that the diagram
$$ \xymatrix{ \calC_{f/} \times_{\calC} \{Z\} \ar[r] \ar[d] & \calC_{X/} \times_{\calC} \{Z\} \ar[d] \\
\calD_{p(f)/} \times_{\calD} \{ p(Z) \} \ar[r] & \calD_{p(X)/} \times_{\calD} \{p(Z) \} }$$
is homotopy Cartesian. Choose a natural transformation $\alpha: \id_{\calC} \rightarrow L$ which exhibits $L$ as a localization functor. Then $\alpha$ induces a natural transformation between the above diagrams. It will therefore suffice to show that each of the induced maps
$$ \calC_{f/} \times_{\calC} \{Z\} \rightarrow \calC_{Lf/} \times_{\calC} \{LZ\} \quad  \calC_{X/} \times_{\calC} \{Z\} \rightarrow \calC_{LX/} \times_{\calC} \{LZ\} $$
$$ \calD_{p(f)/} \times_{\calD} \{p(Z)\} \rightarrow \calD_{p(Lf)/} \times_{\calD} \{ p(LZ) \} \quad \calD_{p(LX)/} \times_{\calD} \{ p(Z) \} \rightarrow \calD_{p(LX)/} \times_{\calD} \{ p(LZ) \}$$
is a homotopy equivalence. For the first pair of maps, this follows from the fact that $Z \in L\calC$.
For the second pair, we observe that $(i)$ and $(ii)$ imply that for every $C \in \calC$, the map
$p( \alpha(C) ): p(C) \rightarrow p(LC)$ is equivalent to the $L'$-localization
$p(C) \rightarrow L'p(C)$, and $p(Z) \in L' \calD$. This completes the proof of $(1)$.

To prove $(2)$, choose any object $C \in L\calC$ and a morphism $f: p(C) \rightarrow D$ in
$L' \calD$. Choose a $p$-coCartesian morphism $\overline{f}: C \rightarrow \overline{D}$ in $\calC$. According to $(1)$, the morphism $L(\overline{f}): LC \rightarrow L\overline{D}$ is $p'$-coCartesian. We now use the fact that $p'$ is a categorical fibration to lift
the equivalence $p( \alpha( \overline{f}) )$ to an equivalence $L(\overline{f}) \simeq \overline{f}'$, where $\overline{f}': C \rightarrow \overline{D}'$ is a $p'$-coCartesian morphism lifting $f$.
\end{proof}

\begin{proposition}\label{localjerk}
Let $p: \calC^{\otimes} \rightarrow \Nerve(\cDelta)^{op}$ be a monoidal structure on the $\infty$-category $\calC = \calC^{\otimes}_{[1]}$, and let $L: \calC \rightarrow \calC$ be a compatible localization functor with essential image $\calD \subseteq \calC$. Then:

\begin{itemize}
\item[$(1)$] The inclusion $D^{\otimes} \subseteq \calC^{\otimes}$ has a left adjoint $L^{\otimes}$.
\item[$(2)$] The restriction $p| \calD^{\otimes}: \calD^{\otimes} \rightarrow \Nerve(\cDelta)^{op}$ is a monoidal structure on the $\infty$-category $\calD$.
\item[$(3)$] The inclusion functor $\calD^{\otimes} \subseteq \calC^{\otimes}$ is lax monoidal, and its left adjoint $L^{\otimes}: \calC^{\otimes} \rightarrow \calD^{\otimes}$ is monoidal.
\end{itemize}
\end{proposition}

\begin{proof}
We first prove $(1)$. According to Proposition \toposref{testreflect}, it will suffice to show that
for every $X \in \calC^{\otimes}$, there exists a localization map $X \rightarrow X'$ where
$X' \in \calD^{\otimes}$.
Let $X$ lie over $[n] \in \cDelta$, corresponding to an $n$-tuple $(C_1, \ldots, C_n) \in \calC^{n}$ under the equivalence $\calC^{\otimes}_{[n]} \simeq \calC^{n}$.
Let $\alpha: X \rightarrow X'$ be a morphism in $\calC^{\otimes}_{[n]}$, which corresponds to 
a localization map $(C_1, \ldots, C_n) \rightarrow (LC_1, \ldots, LC_n)$. We must show that for every object $Y \in \calD^{\otimes}_{[m]}$, composition with $\alpha$ induces a homotopy equivalence
$\bHom_{\calC^{\otimes}}( X', Y) \rightarrow \bHom_{\calC^{\otimes}}( X, Y).$
Both sides can be expressed as a disjoint union, taken over the set of all maps $\beta: [m] \rightarrow [n]$. Using the product structure on $\calC^{\otimes}_{[m]}$, we can reduce to the case $m=1$ and where $\beta$ is endpoint-preserving. In this case, we can identify $Y$ with an object $D \in \calD$. Unwinding the definitions, we are reduced to proving that the natural map
$ \bHom_{\calC}( LC_1 \otimes \ldots \otimes LC_n, D) \rightarrow
\bHom_{\calC}(C_1 \otimes \ldots \otimes C_n, D)$
is a homotopy equivalence. Since $D \in \calD$, it will suffice to show that the map
$ L( C_1 \otimes \ldots \otimes C_n) \rightarrow L( LC_1 \otimes \ldots \otimes LC_n)$ is an equivalence, which follows easily from the compatibility of $L$ with the monoidal structure on $\calC$. This completes the proof of $(1)$. Moreover, the proof of $(1)$ shows that $L^{\otimes}$ fits into a commutative diagram
$$ \xymatrix{ \calC^{\otimes} \ar[dr]^{p} \ar[rr]^{L'} & & \calC^{\otimes} \ar[dl]^{p} \\
& \Nerve(\cDelta)^{op}. & }$$ Applying the second part of Lemma \ref{surine}, we conclude that 
$p| \calD^{\otimes}$ is a coCartesian fibration. The equivalence $\calD^{\otimes}_{[n]} \simeq
(\calD)^n$ follows by construction. This proves $(2)$. The fact that the inclusion $L \calC^{\otimes} \subseteq \calC^{\otimes}$ is lax monoidal follows by inspection. The assertion that $L^{\otimes}$ is a monoidal functor follows from the first part of Lemma \ref{surine}.
\end{proof}

We now discuss a slightly more concrete situation, where the localization functors in question are given by truncation with respect to a t-structure on a stable $\infty$-category.

\begin{definition}\label{tuppa}\index{t-structure!compatible with a monoidal structure}
Let $\calC$ be a stable $\infty$-category equipped with a t-structure, and let
$\calC^{\otimes} \rightarrow \Nerve(\cDelta)^{op}$ be a monoidal structure on $\calC$. We will say that the monoidal structure on $\calC$ is {\it compatible} with the t-structure on $\calC$ if the following conditions are satisfied:
\begin{itemize}
\item[$(1)$] For each $C \in \calC$, the functors $\bigdot \otimes C$ and $C \otimes \bigdot$ are
exact functors from $\calC$ to itself.
\item[$(2)$] The full subcategory $\calC_{\geq 0}$ contains the unit object of $\calC$ and is closed under tensor products.
\end{itemize}
\end{definition}

\begin{remark}\label{sunderfone}
Let $\calC$ be as in Definition \ref{tuppa}. Then, for every pair of integers $m,n \in \Z$, the
tensor product functor $\otimes$ carries
$\calC_{\geq m} \times \calC_{\geq n}$ into $\calC_{ \geq n+m}$. This follows immediately from the exactness of the tensor product in each variable, and the assumption that the desired result holds in the case $m=n=0$.
\end{remark}

\begin{proposition}\label{jumperr}
Let $\calC$ be a stable $\infty$-category equipped with a t-structure, and let
$\calC^{\otimes} \rightarrow \Nerve(\cDelta)^{op}$ be a compatible monoidal structure on $\calC$.
Then, for every $n \geq 0$, the localization functor
$\tau_{\leq n}: \calC_{\geq 0} \rightarrow \calC_{\geq 0}$
is compatible $($in the sense of Definition \ref{compatmon}$)$ with the induced monoidal structure on $\calC_{\geq 0}$ $($Proposition \ref{yukyik}$)$. 
\end{proposition}

\begin{proof}
Let $X,Y \in \calC_{\geq 0}$, and let $f: X \rightarrow \tau_{\leq n} X$ be the canonical map. In view of Remark \ref{suferr}, we need only show that the maps $f \otimes \id_{Y}$ and $\id_{Y} \otimes f$ become equivalences after applying the functor $\tau_{\leq n}$. By symmetry, it will suffice to treat the case of $f \otimes \id_{Y}$. Since the functor $\bigdot \otimes Y$ is exact, we have a long exact sequence
$$ \ldots \rightarrow \pi_{k} ( \tau_{\geq n+1} X \otimes Y) \rightarrow \pi_{k} (X \otimes Y)
\stackrel{\beta_{k}}{\rightarrow} \pi_{k} (\tau_{\leq n} X \otimes Y) \rightarrow 
\pi_{k-1} (\tau_{\geq n+1} X \otimes Y) \rightarrow \ldots $$
in the abelian category $\heart{\calC}$. We wish to show that $\beta_{k}$ is an isomorphism for
$k \leq n$. In view of the long exact sequence, it will suffice to show that $\pi_{k} ( \tau_{\geq n+1} X \otimes Y)$ vanishes for $k \leq n$. In other words, we must show that $\tau_{\geq n+1} X \otimes Y$ belongs to $\calC_{\geq n+1}$, which follows immediately from Remark \ref{sunderfone}.
\end{proof}

\begin{remark}
In the case where $\calC$ is presentable, the tensor product $\otimes$ preserves colimits separately in each variable, and the t-structure on $\calC$ is accessible, Proposition \ref{jumperr} is a special case of Example \ref{jumper}.
\end{remark}

\subsection{Free Algebras}\label{monoid5}

Let $(\calC,\otimes)$ be a monoidal category, and let $\Alg(\calC)$ be the category of algebra objects of $\calC$. Suppose that $\calC$ admits countable coproducts, and that the bifunctor
$\otimes: \calC \times \calC \rightarrow \calC$ preserves countable coproducts separately in each variable. In this case, the forgetful functor $\Alg(\calC) \rightarrow \calC$ admits a left adjoint.
This left adjoint associates to each object $C \in \calC$ the {\it free algebra}
$$\Free_{\calC}(C) = \coprod_{n \geq 0} C^{\otimes n}.$$
Our goal in this section is to prove an $\infty$-categorical analogue of this statement.

\begin{definition}\index{algebra object!free}\index{free!algebra object}
Let $\calC$ be a monoidal $\infty$-category, and let
$\theta: \Alg(\calC) \rightarrow \calC$ be the forgetful functor. Let $C \in \calC$ be an object.
A {\it free algebra generated by $C$} is an object $A \in \Alg(\calC)$ equipped with a map
$\phi: C \rightarrow \theta(A)$, such that composition with $\phi$ induces a homotopy equivalence
$$ \bHom_{\Alg(\calC)}(A,B) \rightarrow \bHom_{\calC}(C, \theta(B) )$$
for every $B \in \Alg(\calC)$.
\end{definition}


The main result of this section is the following:

\begin{theorem}\label{hutmunn}
Let $\calC$ be a monoidal $\infty$-category. Assume that
$\calC$ admits countable coproducts, and that the bifunctor
$\otimes: \calC \times \calC \rightarrow \calC$ preserves countable coproducts separately in each variable. 

\begin{itemize}
\item[$(1)$] The forgetful functor $\Alg(\calC) \rightarrow \calC$ admits a left adjoint.

\item[$(2)$] For every object $C \in \calC$, there exists a map $C \rightarrow A( [1]) $ which
exhibits $A \in \Alg(\calC)$ as the free algebra generated by $\calC$.

\item[$(3)$] For every object $C \in \calC$ and every algebra object $A \in \Alg(\calC)$, a
map $C \rightarrow A([1])$ exhibits $A$ as a free algebra generated by $C$ in $\calC$ if and only if the induced map
$ \coprod_{n \geq 0} C^{\otimes n} \rightarrow A([1])$ is an equivalence in $\calC$ $($that is, if and only if $A$ is a free algebra generated by $C$ in the homotopy category $\h{\calC}${}$)$.
\end{itemize}
\end{theorem}

The proof of Theorem \ref{hutmunn} is long and complicated, and will occupy our attention throughout most of this section. Before embarking on the proof, we will treat a special case which can be established using much easier arguments (and which requires fewer hypotheses):

\begin{proposition}\label{gurgle}\index{algebra object!initial}
Let $\calC$ be a monoidal $\infty$-category. Then:
\begin{itemize}
\item[$(1)$] The $\infty$-category $\Alg(\calC)$ has an initial object.
\item[$(2)$] An object $A \in \Alg(\calC)$ is initial if and only if the unit map
$1_{\calC} \rightarrow A([1])$ is an equivalence in $\calC$.
\end{itemize}
\end{proposition}

\begin{proof}
Let $q: \calC^{\otimes} \rightarrow \Nerve(\cDelta)^{op}$ be a monoidal structure on
$\calC = \calC^{\otimes}_{[1]}$. The fiber $\calC^{\otimes}_{[0]}$ is contractible, so it contains an object $C_{[0]}$. Using Lemma \toposref{kan2} we deduce that there exists a map $A: \Nerve(\cDelta)^{op} \rightarrow \calC^{\otimes}$ which is a section of $q$, where $A([0]) = C_{[0]}$ and $A$ is a $q$-left Kan extension of $A | \{ [0] \}$. It is easy to see that $A$ carries each morphism in $\cDelta$ to a $q$-coCartesian edge of $\calC^{\otimes}$; in particular we conclude that
$A \in \Alg(\calC)$. Using Proposition \toposref{leftkanadj}, we deduce that $A$ is an initial object
of $\Alg(\calC)$ (in fact, it is an initial object in the larger $\infty$-category of {\em all} sections
of $q$). This proves $(1)$. 

We observe that the unit map $1_{\calC} \rightarrow A([1])$ is an equivalence. The ``only if'' direction of $(2)$ now follows from the fact that initial objects of $\Alg(\calC)$ are unique up to equivalence. To prove the converse, we consider an arbitrary algebra object $B \in \Alg(\calC)$. Since $A$ is an initial object of $\Alg(\calC)$, there exists a map $f: A \rightarrow B$ which induces a commutative triangle
$$ \xymatrix{ & 1_{\calC} \ar[dl]^{u_A} \ar[dr]^{u_B} & \\
A([1]) \ar[rr]^{f([1])} & & B([1]) }$$
in $\calC$. Here $u_A$ and $u_B$ denote the unit maps of $A$ and $B$, respectively. If $u_B$ is an equivalence, then we conclude that $f([1])$ is an equivalence, so that $f$ is itself an equivalence (Corollary \ref{jumunj}). This proves the ``if'' direction of $(2)$.
\end{proof}

We now return to the proof of Theorem \ref{hutmunn}. The argument is long and technical, and can be safely skipped; the details will not be important elsewhere in this paper. We begin with an outline of our strategy:

\begin{itemize}
\item[$(a)$] Using the results of \S \ref{segalapp}, we can reduce to proving an analogous assertion for {\em Segal} monoidal $\infty$-categories. Let us therefore assume that 
we are given a map $q: \calC^{\otimes} \rightarrow \Nerve( \LinSeg)$ that exhibits
$\calC = \calC^{\otimes}_{\seg{1}}$ as a Segal monoidal $\infty$-category (we will also assume that the reader is familiar with the notation and terminology from \S \ref{segalapp}).

\item[$(b)$] We will define a sequence of categories $\calJ_0 \subseteq \calJ_1 \subseteq \calJ_2 \subseteq \calJ_3$
equipped with (compatible) functors $\psi_i: \calJ_i \rightarrow \LinSeg$.
This will allow us to factor the forgetful functor $\theta: \SegAlg(\calC) \rightarrow \calC$ as a composition
\begin{eqnarray*}
\SegAlg(\calC) & \rightarrow & \bHom_{ \Nerve(\LinSeg)}( \Nerve(\calJ_3), \calC^{\otimes}) \\
& \rightarrow & \bHom_{ \Nerve(\LinSeg)}( \Nerve(\calJ_2), \calC^{\otimes})
 \\
& \rightarrow & \bHom_{ \Nerve(\LinSeg)}( \Nerve(\calJ_1), \calC^{\otimes}) \\
& \rightarrow & \bHom_{ \Nerve(\LinSeg)}( \Nerve(\calJ_0), \calC^{\otimes}) \\
& \rightarrow & \calC. 
\end{eqnarray*}

\item[$(c)$] We will refine the above picture by selecting appropriate subcategories
$ \Spread_i(\calC) \subseteq \bHom_{ \Nerve(\LinSeg)}( \Nerve(\calJ_i), \calC^{\otimes})$
so that the forgetful functor $\theta$ factors as a composition
$$ \Alg(\calC) \stackrel{\theta_4}{\rightarrow} \Spread_3(\calC)
\stackrel{\theta_3}{\rightarrow} \Spread_2(\calC) \stackrel{\theta_2}{\rightarrow}
\Spread_1(\calC) \stackrel{\theta_1}{\rightarrow} \Spread_0(\calC) \stackrel{\theta_0}{\rightarrow} \calC.$$

\item[$(d)$] We will show that the functor $\theta_0$ is a trivial Kan fibration (Lemma \ref{smallbird}).

\item[$(e)$] We will show that the restriction functor $\theta_1$ admits a left adjoint $\phi_1$, given by $q$-left Kan extension along the inclusion $\Nerve(\calJ_0) \subseteq \Nerve(\calJ_1)$.

\item[$(f)$] We will show that the restriction functor $\theta_2$ admits a left adjoint $\phi_2$. This is a slightly tricky part of the argument, since the left adjoint to $\theta_2$ is {\em not} given by a $q$-left Kan extension (the relevant Kan extension typically does not exist).

\item[$(g)$] The functor $\theta_3$ does not generally admit a left adjoint. Nevertheless, we will prove that there exists a partially-defined left adjoint to $\theta_3$, whose domain includes the essential image of $\phi_2 \circ \phi_1$. This adjoint is again given by a $q$-left Kan extension along the inclusion $\Nerve(\calJ_2) \subseteq \Nerve(\calJ_3)$.

\item[$(h)$] Finally, it will follow immediately from the construction that the functor $\theta_4$ has a left adjoint, given by composition with a section to the functor $\psi_3: \calJ_3 \rightarrow \LinSeg$.
\end{itemize}

The definitions of the categories $\calJ_i$ become more complicated as the index $i$ increases. We therefore begin with the simplest case, where $i = 0$.

\begin{notation}
The category $\calJ_0$ is defined to be the following subcategory of $\LinSeg$:
\begin{itemize}
\item[$(J0)$] Every object of $\LinSeg$ belongs to $\calJ_0$.
\item[$(J0')$] A morphism $\alpha: J_{\ast} \rightarrow J'_{\ast}$ of $\LinSeg$ belongs to 
$\calJ_0$ if and only if every element of $J'$ has exactly one preimage under $\alpha$.
\end{itemize}
We let $\psi_0: \calJ_0 \rightarrow \LinSeg$ denote the inclusion.
\end{notation}

\begin{notation}\label{pinspread}
Let $q: \calC^{\otimes} \rightarrow \Nerve(\LinSeg)$ be a Segal monoidal $\infty$-category.
We let $\Spread_0(\calC)$ denote the full subcategory of
$\bHom_{\Nerve(\LinSeg)}( \Nerve(\calJ_0), \calC^{\otimes})$ spanned by those functors
$f \in \bHom_{\Nerve(\LinSeg)}( \Nerve(\calJ_0), \calC^{\otimes})$ which carry
every morphism in $\calJ_0$ to a $q$-coCartesian morphism in $\calC^{\otimes}$.
\end{notation}

\begin{example}
Let $q: \calC^{\otimes} \rightarrow \Nerve(\LinSeg)$ be a Segal monoidal $\infty$-category and let
$A$ be a section of $q$. We observe that $A$ is a Segal algebra object of $\calC$ if and only if $A| \Nerve(\calJ_0) \in \Spread_0(\calC)$. In particular, restriction determines a functor
$\SegAlg(\calC) \rightarrow \Spread_0(\calC)$.
\end{example}

More informally, a functor $f: \Nerve(\calJ_0) \rightarrow \calC^{\otimes}$ belongs to $\Spread_0(\calC)$ if there exists an object $C = f( \seg{1} ) \in \calC$, such that for all $n \geq 0$, the object $f( \seg{n} ) \in \calC^{\otimes}_{ \seg{n} }$ corresponds to $(C,C, \ldots, C)$ under the equivalence
$\calC^{\otimes}_{ \seg{n} } \simeq \calC^{n}$. The functor $F$ is then determined by $C$, up to canonical equivalence. Our next result makes this statement more precise:

\begin{lemma}\label{smallbird}
Let $q: \calC^{\otimes} \rightarrow \Nerve(\LinSeg)$ be a Segal monoidal structure on an $\infty$-category $\calC = \calC^{\otimes}_{ \seg{1} }$. Then evaluation at $\seg{1} \in \calJ_0$ induces
a trivial Kan fibration $\theta_0: \Spread_0(\calC) \rightarrow \calC$.
\end{lemma}

\begin{proof}
Fix $n \geq 0$, and regard $\seg{n}$ as an object of $\calJ_0$. Let $\calJ_{-1}$ denote the full subcategory of  $\calJ_0$ spanned by the object $\seg{1}$. We observe that a functor
$f \in \bHom_{ \Nerve(\LinSeg)} ( (\Nerve(\calJ_{-1})_{\seg{n}/})^{\triangleleft}, \calC^{\otimes})$
is a $q$-limit diagram if and only if it carries each morphism of
$((\Nerve( \calJ_{-1})^{op})_{\seg{n}/})^{\triangleleft}$ to a $q$-coCartesian morphism in $\calC^{\otimes}$. Using Lemma \toposref{kan2}, we deduce that every functor
$g_0 \in \bHom_{ \Nerve(\LinSeg)}( \Nerve( \calJ_{-1} ), \calC^{\otimes})$ admits a $q$-right Kan extension
to a functor $g \in \bHom_{ \Nerve(\LinSeg)}( \Nerve(\calJ_0), \calC^{\otimes})$. Let
$\calD \subseteq \bHom_{ \Nerve(\LinSeg)}( \Nerve(\calJ_0), \calC^{\otimes})$ be the full subcategory spanned by those functors $g$ which are $q$-right Kan extensions of $g | \Nerve(\calJ_{-1})$. Proposition \toposref{lklk} implies that restriction determines a trivial Kan fibration
$ \calD \rightarrow  \bHom_{ \Nerve(\LinSeg)}( \Nerve( \calJ_{-1}) , \calC^{\otimes}) \simeq \calC.$
To complete the proof, it will suffice to show that $\calD = \Spread_0(\calC)$. Using the description of the $q$-limit diagrams given above, we deduce that a functor $g \in \bHom_{ \Nerve(\LinSeg)}( \Nerve(\calJ_0), \calC^{\otimes})$ belongs to $\calD$ if and only if the following condition is satisfied:
\begin{itemize}
\item[$(\ast)$] For every morphism $\alpha: \seg{n} \rightarrow \seg{1}$ in $\calJ_0$, the morphism
$g(\alpha)$ is $q$-coCartesian.
\end{itemize}
From this characterization it follows immediately that $\Spread_0(\calC) \subseteq \calD$. 
To prove the reverse inclusion, let us suppose that $g \in  \bHom_{ \Nerve(\LinSeg)}( \Nerve(\calJ_0), \calC^{\otimes})$ satisfies $(\ast)$. We wish to show that $g$ is a coCartesian section of $q$. Let $\alpha: \seg{n} \rightarrow \seg{m}$ be a morphism in $\calJ_0$, let $\overline{\alpha} = g(\alpha)$,
and let $\widetilde{\alpha}: g( \seg{n} ) \rightarrow C$ be a $q$-coCartesian morphism of $\calC^{\otimes}$ which projects to $\alpha$. There exists an (essentially unique) $2$-simplex
$$ \xymatrix{ & C \ar[dr]^{\widetilde{\beta}} & \\
g( \seg{n} ) \ar[ur]^{\widetilde{\alpha}} \ar[rr]^{\overline{\alpha}} & & g( \seg{m} ) }$$
in $\calC^{\otimes}$. We wish to show that $\beta$ is an equivalence. For this, it will suffice to show that, for every morphism $\gamma: \seg{m} \rightarrow \seg{1}$ in $\calJ_0$, the image of
$\widetilde{\beta}$ under the associated functor $\gamma_!: \calC^{\otimes}_{\seg{m}} \rightarrow \calC^{\otimes}_{\seg{1}}$ is an equivalence. We have a commutative diagram
$$ \xymatrix{ & C \ar[dr]^{\widetilde{\beta}} \ar[rr]^{\widetilde{\gamma}} & & \gamma_! C \ar[dr]^{\gamma_! \widetilde{\beta}} & \\
g( \seg{n} ) \ar[ur]^{\widetilde{\alpha}} \ar[rr]^{\overline{\alpha}} & & g( \seg{m} ) \ar[rr]^{\overline{\gamma}} & & \gamma_! g( \seg{m} ). }$$
Consequently, it will suffice to show that the morphisms $\overline{\gamma} \circ \overline{\alpha}$
and $\widetilde{\gamma} \circ \widetilde{\alpha}$ are $q$-coCartesian. In the second case, this follows from the fact that $q$-coCartesian morphisms are stable under composition. In the first case, we can use $(\ast)$ to assume, without loss of generality, that $\overline{\gamma} = g(\gamma)$. Then $\overline{\gamma} \circ \overline{\alpha} \simeq g( \gamma \circ \alpha)$ is $q$-coCartesian (invoking $(\ast)$ once more), as desired.
\end{proof}

\begin{notation}
The category $\calJ_1$ is defined as follows:
\begin{itemize}
\item[$(J1)$] Objects of $\calJ_1$ are given by morphisms $\alpha: \seg{n} \rightarrow \seg{m}$ in
$\LinSeg$.
\item[$(J1')$] Given a pair of objects $\alpha, \alpha' \in \calJ_1$, a morphism from $\alpha$ to
$\alpha'$ is a commutative diagram
$$ \xymatrix{ \seg{n} \ar[r]^{\alpha} \ar[d]^{\gamma} & \seg{m} \ar[d] \\
\seg{n'} \ar[r]^{\alpha'} & \seg{m'} }$$
in $\LinSeg$, such that $\gamma$ is a morphism in $\calJ_0$.
\end{itemize}
We will identify $\calJ_0$ with the full subcategory of $\calJ_1$ spanned by those morphisms $\alpha$ which are isomorphisms in $\LinSeg$.
The inclusion $\psi_0: \calJ_0 \rightarrow \LinSeg$ extends to a functor
$\psi_1: \calJ_1 \rightarrow \LinSeg$, given by the formula
$$ ( \alpha: \seg{n} \rightarrow \seg{m} ) \mapsto \seg{m}.$$
\end{notation}

\begin{notation}
Let $\calC$ be a Segal monoidal $\infty$-category. We define $\Spread_1(\calC) \subseteq \bHom_{ \Nerve(\LinSeg)}( \Nerve(\calJ_1), \calC^{\otimes})$ to be the full subcategory spanned by those functors $f$ such that $f | \Nerve(\calJ_0)$ belongs to $\Spread_0(\calC)$.
\end{notation}

\begin{remark}
Let $\calC$ be a monoidal $\infty$-category. Then composition with $\psi_1: \calJ_1 \rightarrow \LinSeg$ induces a functor $\SegAlg(\calC) \rightarrow \Spread_1(\calC)$.
\end{remark}

We now analyze the relationship between $\Spread_0(\calC)$ and $\Spread_1(\calC)$.

\begin{lemma}\label{firstpart}
Let $q: \calC^{\otimes} \rightarrow \Nerve(\LinSeg)$ be a Segal monoidal $\infty$-category. Then:
\begin{itemize}
\item[$(1)$] A functor $f \in \bHom_{\Nerve(\LinSeg)}( \Nerve(\calJ_1), \calC^{\otimes})$ is a $q$-left Kan extension of $f_0 = f | \Nerve(\calJ_0)$ if and only if, for every
object $\alpha: \seg{n} \rightarrow \seg{m}$ in $\calJ_1$, if 
$e_{\alpha}: \id_{[n]} \rightarrow \alpha$ denotes the associated morphism of $\calJ_1$, then
$f(e_{\alpha})$ is a $q$-coCartesian morphism in $\calC^{\otimes}$.

\item[$(2)$] Every functor $f_0 \in \bHom_{ \Nerve(\LinSeg)}( \Nerve(\calJ_0), \calC^{\otimes})$ admits an extension $f \in \bHom_{ \Nerve(\LinSeg)}( \Nerve(\calJ_1), \calC^{\otimes})$ which is a $q$-left Kan extension of $f$.

\item[$(3)$] Let $\bHom'_{ \Nerve(\LinSeg)}( \Nerve(\calJ_1), \calC^{\otimes})$ denote the full subcategory of $\bHom_{ \Nerve(\LinSeg)}( \Nerve(\calJ_1), \calC^{\otimes})$ spanned by those functors $f$ which are $q$-left Kan extensions of $f_0 = f| \Nerve(\calJ_0)$. Then
the restriction map
$$\bHom'_{ \Nerve(\LinSeg)}( \Nerve(\calJ_1), \calC^{\otimes}) \rightarrow
\bHom_{ \Nerve(\LinSeg)}( \Nerve(\calJ_0), \calC^{\otimes})$$
is a trivial Kan fibration.

\item[$(4)$] The restriction map $\Spread_1(\calC) \rightarrow \Spread_0(\calC)$ admits a section $\phi_1$ which is simultaneously a left adjoint. The essential image of $\phi_1$ consists of those functors $f \in \Spread_1(\calC)$ which satisfy the condition stated in $(1)$.
\end{itemize}
\end{lemma}

\begin{proof}
We observe that if $\alpha \in \calJ_1$, then the morphism $e_{\alpha}$ is a final object of the
category $(\calJ_1)_{/ \alpha} \times_{ \calJ_1} \calJ_0$. Assertion $(1)$ now follows immediately from Proposition \toposref{relexists}, and assertion $(2)$ follows from Lemma \toposref{kan2} and Proposition \toposref{relexist}, since $q$ is a coCartesian fibration. Assertion $(3)$ follows from Proposition \toposref{lklk}, and $(4)$ from Proposition \toposref{leftkanadj}.
\end{proof}

\begin{notation}
The category $\calJ_2$ is defined as follows.
\begin{itemize}
\item[$(J2)$] An object of $\calJ_2$ consists of a quadruple
$( J, J_0, \seg{n}, \alpha)$, where $J$ is a finite linearly ordered set,
$J_0$ is a subset of $J$, and $\alpha: \seg{n} \rightarrow (J_0)_{\ast}$ is a morphism in
$\LinSeg$.

\item[$(J2')$] Given a pair of objects $( J, J_0, \seg{n}, \alpha), (J', J'_0, \seg{n'} , \alpha') \in \calJ_2,$
a morphism from $( J, J_0, \seg{n}, \alpha)$ to $( J', J'_0, \seg{n'}, \alpha')$ in $\calJ_2$ consists of a pair of morphisms $\beta: J_{\ast} \rightarrow J'_{\ast}$, $\gamma: \seg{n} \rightarrow \seg{n'}$ in $\LinSeg$, satisfying the following conditions:

\begin{itemize}
\item The morphism $\gamma$ belongs to $\calJ_0$.

\item Suppose $j \in J'$ and $\beta^{-1} \{j\} \subseteq J_0$. Then $j \in J'_0$, and $\gamma$ induces a bijection $(\alpha \circ \beta)^{-1} \{j\} \rightarrow (\alpha')^{-1} \{j\}$. 
\end{itemize}
\end{itemize}

We will identify $\calJ_1$ with the full subcategory of $\calJ_2$ spanned by those objects
$( J, J_0, \seg{n}, \alpha)$ such that $J_0 = J$.
We note that the functor $\psi_1: \calJ_1 \rightarrow \LinSeg$ admits an extension
$\psi_2: \calJ_2 \rightarrow \LinSeg$, given by the formula
$$( J, J_0, \seg{n}, \alpha) \mapsto J_{\ast}.$$
\end{notation}

\begin{notation}\label{minspread}
Let $q: \calC^{\otimes} \rightarrow \Nerve(\LinSeg)$ be a Segal monoidal $\infty$-category. We let $\Spread_2(\calC)$ denote the full subcategory
of $\bHom_{\Nerve(\LinSeg)}( \Nerve(\calJ_2), \calC^{\otimes})$ spanned by those functors $f$ which possess the following properties:
\begin{itemize}
\item[$(i)$] The restriction $f_0 = f | \Nerve(\calJ_0)$ belongs to $\Spread_0(\calC)$.
\item[$(ii)$] Let $(J, J_0, \seg{n}, \alpha)$ be an object of $\calJ_2$, and let
$$(\beta, \id): (J, J_0, \seg{n}, \alpha) \rightarrow (J_0, J_0, \seg{n}, \alpha)$$
be given by the formula
$$\beta(j) = \begin{cases} j & \text{if } j \in J_0 \\
\ast & \text{if } j \notin J_0. \end{cases}$$
Then $f( \beta, \id)$ is a $q$-coCartesian morphism in $\calC^{\otimes}$.
\end{itemize}
\end{notation}

\begin{remark}\label{umpertine}
Composition with $\psi_2: \calJ_2 \rightarrow \LinSeg$ induces a functor
$\SegAlg(\calC) \rightarrow \Spread_2(\calC)$.
\end{remark}

\begin{lemma}\label{2ndpart}
Let $q: \calC^{\otimes} \rightarrow \Nerve(\LinSeg)$ be a Segal monoidal $\infty$-category. Assume that the underlying $\infty$-category $\calC$ has an initial object, and a tensor product
$C_1 \otimes \ldots \otimes C_k$ in $\calC$ is initial whenever any of the factors $C_i$ is initial.
Then:
\begin{itemize}
\item[$(1)$] The restriction functor $\theta_2: \Spread_2(\calC) \rightarrow \Spread_1(\calC)$ admits a section $\phi_2$, which is simultaneously a left adjoint to $\theta_2$.
\item[$(2)$] A functor $f \in \Spread_2(\calC)$ belongs to the essential image of $\phi_2$ if and only if the following condition is satisfied:
\begin{itemize}
\item[$(\ast)$] For every object $( J, J_0, \seg{n}, \alpha) \in \calJ_2$ and every
element $j \in J - J_0$, if $\beta: J_{\ast} \rightarrow \{j\} \cup \{\ast\}$ is the morphism of
$\LinSeg$ such that $\beta^{-1} \{j\} = \{j\}$, then the induced functor
$\beta_{!}: \calC^{\otimes}_{ J_{\ast} } \rightarrow \calC^{\otimes}_{ \{j\}_{\ast}} = \calC$
carries $f(J, J_0, \seg{n}, \alpha)$ to an initial object of $\calC$.
\end{itemize}
\end{itemize}
\end{lemma}

\begin{proof}
We observe that the inclusion $\calJ_1 \subseteq \calJ_2$ has a left adjoint $U$, given by
$$ (J, J_0, \seg{n}, \alpha) \mapsto (J_0, J_0, \seg{n}, \alpha).$$
Let $T: \Delta^1 \times \Nerve(\calJ_2) \rightarrow \Nerve(\calJ_2)$ be a natural transformation from $\id$ to $U$ which is a unit for the adjunction. Let $T'$ be the composition of $T$ with
$\psi_2$, regarded as a functor from $\Nerve(\calJ_2)$ to $\Fun( \Delta^1, \Nerve(\LinSeg))$.

Let $\calD$ denote the full subcategory of $\Fun( \Delta^1, \calC^{\otimes})$ spanned by the $q$-coCartesian morphisms, and let $\calM$ be the fiber product
$ \Nerve(\calJ_2) \times_{ \Fun( \Delta^1, \Nerve(\LinSeg) )} \calD.$ Similarly, we let
$\calM_0$ denote the fiber product $\Nerve(\calJ_2) \times_{ \Nerve(\LinSeg)} \calC^{\otimes}$ where $\Nerve(\calJ_2)$ maps to $\Nerve(\LinSeg)$ via $\psi_2$, and
$\calM_1$ the fiber product $\Nerve(\calJ_2) \times_{ \Nerve(\LinSeg)} \calC^{\otimes}$ where
$\Nerve(\calJ_2)$ maps to $\Nerve(\LinSeg)$ via $\psi_1 \circ U$. We have canonical maps
$ \calM_0 \stackrel{e_0}{\leftarrow} \calM \stackrel{e_1}{\rightarrow} \calM_1,$
where $e_0$ is a trivial Kan fibration.

Let $\overline{\Spread}_2(\calC) \subseteq \bHom_{ \Nerve(\calJ_2)}( \Nerve(\calJ_2), \calM)$ denote the inverse image of $\Spread_2(\calC) \subseteq \bHom_{ \Nerve(\calJ_2)}( \Nerve(\calJ_2), \calM_0)$, so that $e_0$ induces a trivial Kan fibration $\overline{\Spread}_2(\calC)$.
Since the restriction of $T'$ to $\Nerve(\calJ_1)$ is an equivalence, we have a canonical map
$\overline{\Spread}_2(\calC) \times \Delta^1 \rightarrow \Spread_1(\calC)$ which induces an
equivalence from the composition
$$ \overline{\Spread}_2(\calC) \stackrel{e_0 \circ}{\rightarrow} \Spread_{2}(\calC) \stackrel{| \Nerve(\calJ_1)}{\rightarrow} \Spread_{1}(\calC)$$
to the composition
$$ \overline{\Spread}_2(\calC) \stackrel{e_1 \circ}{\rightarrow} \bHom_{ \Nerve(\calJ_2)}( \Nerve(\calJ_2), \calM_1) \rightarrow \bHom_{ \Nerve(\LinSeg)}( \Nerve(\calJ_1), \calC^{\otimes}).$$
Unwinding the definitions, we see that there is a pullback diagram
$$ \xymatrix{ \overline{\Spread}_2(\calC) \ar@{^{(}->}[r] \ar[d] & \bHom_{ \Nerve(\calJ_2)}( \Nerve(\calJ_2), \calM) \ar[d]^{g} \\
\bHom'_{ \Nerve(\calJ_2)}( \Nerve(\calJ_2), \calM_1) \ar@{^{(}->}[r] & 
\bHom_{\Nerve(\calJ_2)}( \Nerve(\calJ_2), \calM_1), }$$
where $\bHom'_{ \Nerve(\calJ_2)}( \Nerve(\calJ_2), \calM_1)$ denotes the full subcategory of
$\bHom_{ \Nerve(\calJ_2)}( \Nerve(\calJ_2), \calM_1)$ spanned by those functors $f$ for
which $f_0 = f | \Nerve(\calJ_1)$ belongs to $\Spread_1(\calC)$, and $f$ is a $q$-right Kan extension of $f_0$. Using Proposition \toposref{lklk}, we deduce that the projection
$ \bHom'_{ \Nerve(\calJ_2)}( \Nerve(\calJ_2), \calM_1) \rightarrow \Spread'(\calC)$
is a trivial Kan fibration. Consequently, we may reduce to proving the following analogues of statements $(1)$ and $(2)$ for the map $g$:

\begin{itemize}
\item[$(1')$] The projection
$g: \bHom_{ \Nerve(\calJ_2)}( \Nerve(\calJ_2), \calM) \rightarrow 
\bHom_{ \Nerve(\calJ_2)}( \Nerve(\calJ_2, \calM_1)$
admits a section $s$, which is simultaneously a left adjoint to $g$.

\item[$(2')$] A functor $F \in \bHom_{ \Nerve(\calJ_2)}( \Nerve(\calJ_2), \calM)$ belongs to the essential image of $s$ if and only if the following condition is satisfied:
\begin{itemize}
\item[$(\ast')$] For every object $( J, J_0, \seg{n}, \alpha) \in \calJ_2$ and every
element $j \in J - J_0$, if $\alpha^{j,J}: J_{\ast} \rightarrow \{j\}_{\ast}$ is defined as in Notation \ref{ajack}, then the induced functor
$\alpha^{j,J}_{!}: \calC^{\otimes}_{ J_{\ast} } \rightarrow \calC^{\otimes}_{ \{j\}_{\ast}} = \calC$
carries $F(J, J_0, \seg{n}, \alpha)$ to an initial object of $\calC$.
\end{itemize}
\end{itemize}

We observe that the projection $e_1: \calM \rightarrow \calM_1$ is a coCartesian fibration.
Moreover, each fiber of $e_1$ is equivalent to a product of copies of $\calC$, and therefore admits an initial object. Let $\calM'$ denote the collection of objects $M \in \calM$ such that
$M$ is an initial object of the fiber $\calM \times_{ \calM_1} \{ e_1(M) \}$. We observe that
$(2')$ is equivalent to the assertion that the essential image of $s$ is
the $\infty$-category $\bHom_{ \Nerve(\calJ_2)}( \Nerve(\calJ_2), \calM')$. 
Using our assumption that a tensor product $C \otimes D \in \calC$ is initial provided that either $C$ or $D$ is initial, and the definition of the category $\calJ_2$, we deduce that the collection of objects of $\calM'$ is stable under the collection of functors $\calM_{x} \rightarrow \calM_{y}$ associated to morphisms $x \rightarrow y$ in $\calM_1$. 
Proposition \toposref{relcolfibtest} implies that every object of $\calM'$ is an $e_1$-initial object of $\calM$. It now follows from Proposition \toposref{lklk} that $g$ restricts to a trivial Kan fibration
$$g': \bHom_{ \Nerve(\calJ_2)}( \Nerve(\calJ_2), \calM') \rightarrow 
\bHom_{ \Nerve(\calJ_2)}( \Nerve(\calJ_2, \calM_1). $$
Let $s$ be a section to $g'$. To complete the proof, it will suffice to show that $s$ is a left adjoint to $g$. In other words, it suffices to show that if $X \in \bHom_{ \Nerve(\calJ_2)}( \Nerve(\calJ_2), \calM_1)$ and $Y \in  \bHom_{ \Nerve(\calJ_2)}( \Nerve(\calJ_2), \calM)$, then $g$ induces a homotopy equivalence of Kan complexes
$$ \Hom^{\rght}_{  \bHom_{ \Nerve(\calJ_2)}( \Nerve(\calJ_2), \calM) }( s(X), Y) 
\rightarrow \Hom^{\rght}_{  \bHom_{ \Nerve(\calJ_2)}( \Nerve(\calJ_2), \calM_1) }( X, g(Y) ).$$
But this map is a trivial Kan fibration, since $s(X)$ is a $g$-initial object of
$ \bHom_{ \Nerve(\calJ_2)}( \Nerve(\calJ_2), \calM)$ by Lemma \ref{surtybove}.
\end{proof}

\begin{notation}
The category $\calJ_3$ is defined to be the categorical mapping cylinder of 
the functor $\psi_2: \calJ_2 \rightarrow \LinSeg$. More precisely, this category may be described as follows:
\begin{itemize}
\item[$(J3)$] An object of $\calJ_3$ is either an object of $\LinSeg$ or an object of $\calJ_2$.
\item[$(J3')$] Given a pair of objects $\seg{m}, \seg{n} \in \LinSeg$, we have
$\Hom_{\calJ_3}( \seg{m}, \seg{n} ) = \Hom_{\LinSeg}( \seg{m}, \seg{n} )$.
\item[$(J3'')$] Given a pair of objects $X, Y \in \calJ_2$, we have
$ \Hom_{\calJ_3}(X,Y) = \Hom_{\calJ_2}(X,Y)$.
\item[$(J3''')$] Given a pair of objects $X= ( J, J_0, \seg{m}, \alpha) \in \calJ_2$, $\seg{n} \in \LinSeg$, we have 
$$\Hom_{\calJ_3}( X, \seg{n}) = \Hom_{\LinSeg}(J_{\ast}, \seg{n} ) \quad \Hom_{\calJ_3}( \seg{n},X ) = \emptyset.$$
\end{itemize}

By construction, we may regard $\calJ_2$ as a full subcategory of $\calJ_3$.
The functor $\psi_2: \calJ_2 \rightarrow \LinSeg$ extends canonically to a retraction
$\psi_3: \calJ_3 \rightarrow \LinSeg$.
\end{notation}

\begin{notation}\label{notebard}
Let $q: \calC^{\otimes} \rightarrow \Nerve(\LinSeg)$ be a Segal monoidal $\infty$-category. We let $\Spread'_2(\calC) \subseteq \Spread_2(\calC)$ denote the essential image of the functor
$\phi_2 \circ \phi_1: \Spread_0(\calC) \rightarrow \Spread_2(\calC)$.
In other words, $\Spread'_2(\calC)$ is the full subcategory of $\bHom_{ \Nerve(\LinSeg)}(
\Nerve(\calJ_2), \calC^{\otimes})$ spanned by
those functors $f \in \Nerve( \calJ_2) \rightarrow \calC^{\otimes}$ which satisfy the following conditions: 

\begin{itemize}
\item[$(A0)$] The composition $q \circ f$ coincides with the map
$\Nerve(\calJ_2) \rightarrow \Nerve( \LinSeg)$ induced by $\psi_2$.

\item[$(A1)$] For every morphism $\alpha: J_{\ast} \rightarrow J'_{\ast}$ in
$\calJ_0$, $f$ carries the associated map
$$ (\alpha, \alpha): ( J, J, J_{\ast}, \id_{ J_{\ast} }) \rightarrow ( J', J', J'_{\ast}, \id_{J'_{\ast}})$$
to a $q$-coCartesian morphism in $\calC^{\otimes}$ (Notation \ref{pinspread}).

\item[$(A2)$] Let $(J, J_0, \seg{n}, \alpha)$ be an object of $\calJ_2$, and let
$(\beta, \id): (J, J_0, \seg{n}, \alpha) \rightarrow (J_0, J_0, \seg{n}, \alpha)$
be given by the formula
$$\beta(j) = \begin{cases} j & \text{if } j \in J_0 \\
\ast & \text{if } j \notin J_0. \end{cases}$$
Then $f( \beta, \id)$ is a $q$-coCartesian morphism in $\calC^{\otimes}$ (Notation \ref{minspread}).

\item[$(A3)$] For every object $( J, J_0, \seg{n}, \alpha) \in \calJ_2$ and every
element $j \in J - J_0$, if $\alpha^{j,J}: J_{\ast} \rightarrow \{j\}_{\ast}$ is 
defined as in Notation \ref{ajack}, then the induced functor
$\beta_{!}: \calC^{\otimes}_{ J_{\ast} } \rightarrow \calC^{\otimes}_{ \{j\}_{\ast}} = \calC$
carries $f(J, J_0, \seg{n}, \alpha)$ to an initial object of $\calC$ (Lemma \ref{2ndpart}). 

\item[$(A4)$] For every morphism $\alpha: J_{\ast} \rightarrow J'_{\ast}$ in
$\LinSeg$, $f$ carries the associated map
$$ (\alpha, \id_{J_{\ast}}): ( J, J, J_{\ast}, \id_{ J_{\ast} }) \rightarrow ( J', J', J_{\ast}, \alpha)$$
to a $q$-coCartesian morphism in $\calC^{\otimes}$ (Lemma \ref{firstpart}).
\end{itemize}
\end{notation}

\begin{lemma}\label{upperness}
Let $A$ and $B$ be marked simpicial sets (see \S \toposref{twuf}), and let
$j: A \rightarrow B$ be a marked equivalence. Let $\calC$ be an $\infty$-category, and $p: B \rightarrow \calC^{\natural}$ a map of marked simplicial sets. Then
the induced map $q: \calC_{/p} \rightarrow \calC_{/pj}$ is a categorical equivalence.
\end{lemma}

\begin{proof}
Choose a marked anodyne map $A \rightarrow A'$, where $A'$ is a fibrant object
of $\mSet$. Now choose a marked anodyne map $B \coprod_{A} A' \rightarrow B'$, where $B'$ is fibrant. Let $j': A' \rightarrow B'$ denote the induced map; since $(\mSet)$ is left proper, $j'$ is a marked equivalence. The inclusion $B \rightarrow B'$ is marked anodyne, so the map
$p$ extends to a map $p': B' \rightarrow \calC^{\natural}$. We have a commutative diagram
$$ \xymatrix{ \calC_{/p'} \ar[r]^{q'} \ar[d] & \calC_{/p'j'} \ar[d] \\
\calC_{/p} \ar[r]^{q} & \calC_{/pj}. }$$
Proposition \toposref{eggwhite} implies that the vertical maps are trivial Kan fibrations. It will therefore suffice to show that $q'$ is an equivalence of $\infty$-categories. We now observe that
Theorem \toposref{bigdiag} implies that $j'$ is a categorical equivalence of the underlying simplicial sets, so that $q'$ is a categorical equivalence by Proposition \toposref{gorban3}.
\end{proof}

\begin{lemma}\label{gottamake}
Let $p: \calD \rightarrow \calC$ be a functor between $\infty$-categories. Let
$\calD_0 \subseteq \calD$ be the essential image of a localization functor $L: \calD \rightarrow \calD$, and suppose that for every $D \in \calD$, $p$ carries the localization map
$D \rightarrow LD$ to an equivalence in $\calC$. 
Let $p_0 = p | \calD_0$. Then
the induced map $q: \calC_{/p} \rightarrow \calC_{/p_0}$ is a trivial fibration of simplicial sets.
\end{lemma}

\begin{proof}
Since the map $q$ is a right fibration, it will suffice to show that $q$ is a categorical equivalence. 
Let $\alpha: \id_{\calD} \rightarrow L$ be a unit for the adjunction between $L$ and the inclusion
$\calD_0 \subseteq \calD$.
Let $B$ be the marked simplicial set $(\calD, S)$, where $S$ is the collection of morphisms $f$ of $\calD$ such that $L(f)$ is an equivalence, and let $A = \calD_0^{\natural} \subseteq B$.
Then $\alpha$ determines a homotopy $( \Delta^1)^{\sharp} \times B \rightarrow B$
from $\id_{B}$ to the map $L: B \rightarrow A$, and restricts to a homotopy
from $\id_{A}$ to the induced map $L|A: A \rightarrow A$.
It follows that $L$ is a homotopy inverse to the inclusion $j: A \subseteq B$, so that
$j$ is a marked equivalence. The desired result now follows from Lemma \ref{upperness}.
\end{proof}

\begin{lemma}\label{sumotin}
Let $q: \calC^{\otimes} \rightarrow \Nerve(\LinSeg)$ be a Segal monoidal $\infty$-category. Assume that the underlying $\infty$-category $\calC = \calC^{\otimes}_{\seg{1}}$ admits
countable coproducts, and that the bifunctor $\otimes: \calC \times \calC \rightarrow \calC$ preserves countable coproducts separately in each variable. Let
$f_0 \in \Spread'_2( \calC)$. Then:
\begin{itemize}
\item[$(1)$] There exists a commutative diagram
$$ \xymatrix{ \Nerve(\calJ_2) \ar[r]^{f_0} \ar@{^{(}->}[d] & \calC^{\otimes} \ar[d]^{q} \\
\Nerve(\calJ_3) \ar[r]^{\psi_3} \ar[ur]^{f} & \Nerve(\LinSeg), }$$
where $f$ is a $q$-left Kan extension of $f_0$.
\item[$(2)$] The restriction $A = f | \Nerve(\LinSeg)$ is a Segal algebra object of $\calC$.
\item[$(3)$] Let $C = f_0( \nostar{1}, \nostar{1}, \seg{1}, \id_{\seg{1}} ) \in \calC$. Then the canonical map $C \rightarrow A( \seg{1} )$ and the algebra structure on $A$ determine an equivalence
$$ \coprod_{n \geq 0} C^{\otimes n} \rightarrow A( \seg{1} ).$$
\end{itemize}
\end{lemma}

\begin{remark}\label{plaque}
The hypotheses of Lemma \ref{sumotin} have the following consequence, which we will use repeatedly: given a collection of objects $C_1, \ldots, C_k \in \calC$, if any $C_i$ is an initial object of $\calC$, then the tensor product $C_1 \otimes \ldots \otimes C_k$ is initial.
\end{remark}

\begin{proof}
We begin by formulating a more precise version of $(1)$, which will be used to establish $(2)$ and $(3)$. First, we need to introduce a bit of notation. Let $J$ be a finite linearly ordered set, and let $j$ be an element of $J$. Let $\alpha^{j,J}: J_{\ast} \rightarrow \{j\}_{\ast}$ be defined as in Notation \ref{ajack}, and let $\alpha^{j,J}_{!}: \calC^{\otimes}_{J_{\ast}} \rightarrow \calC^{\otimes}_{\{j\}_{\ast}} = \calC$
denote the associated functor. For each $n \geq 0$, we let
$e^{j,J}_{n}$ denote the object of $\calJ_2$ given by
$(J, \{j\}, \seg{n}, \beta)$, where the map $\beta: \seg{n} \rightarrow \{ j\}_{\ast}$ satisfies
$\beta^{-1}(\ast) = \ast$. We note that $\psi_2( e^{j, J}_{n} ) = J_{\ast}$, so there is a canonical map
$e^{j,J}_{n} \rightarrow J_{\ast}$ in $\calJ_3$.
In addition to $(1)$, we will need the following characterization of the left Kan extensions of $f_0$:
\begin{itemize}
\item[$(1')$] Suppose that $f \in \bHom_{ \Nerve(\LinSeg)}( \Nerve(\calJ_3), \calC^{\otimes})$
is such that $f_0 = f | \Nerve(\calJ_2)$. Then $f$ is a $q$-left Kan extension of $f_0$ if and only if,
for every finite linearly ordered set $J$ and every element $j$, $f$ exhibits
$\alpha^{j,J}_{!} f(J_{\ast})$ as a coproduct of the objects $\alpha^{j,J}_{!} f( e^{j,J}_{n} )$ in $\calC$. Here $n$ ranges over all nonnegative integers.
\end{itemize}

Fix an object $J_{\ast} \in \LinSeg$, let
$\calI$ denote the category $\calJ_2 \times_{ \LinSeg} ( \LinSeg)_{/J_{\ast}}$, and let
$f_1$ denote the composition
$$ \Nerve( \calI) \rightarrow \Nerve(\calJ_2) \stackrel{f_0}{\rightarrow} \calC^{\otimes}.$$
According to Lemma \toposref{kan2}, assertion $(1)$ is equivalent to the assertion that, for any choice of $J_{\ast}$, the lifting problem
$$ \xymatrix{ \Nerve(\calI) \ar[r]^{f_1} \ar@{^{(}->}[d] & \calC^{\otimes} \ar[d]^{q} \\
\Nerve( \calI)^{\triangleright} \ar@{-->}[ur]^{\overline{f}_1} \ar[r] & \Nerve(\LinSeg) }$$
admits a solution $\overline{f}_1$ which is a $q$-colimit diagram. Choose a $q$-coCartesian transformation
$f_1 \rightarrow f_2$ covering the canonical natural transformation from $q \circ f_1$ to the constant functor $\Nerve( \calI) \rightarrow \Nerve(\LinSeg)$ taking the value $J_{\ast}$. 
Proposition \toposref{chocolatelast} shows that we can identify $q$-colimits
of $f_1$ with $q$-colimits of $f_2$. Combining this observation with Proposition \toposref{relcolfibtest}, we are reduced to proving the following pair of assertions:

\begin{itemize}
\item[$(a)$] The diagram $f_2: \Nerve(\calI) \rightarrow \calC^{\otimes}_{J_{\ast}}$ admits
a colimit $\overline{f}_2: \Nerve(\calI)^{\triangleright} \rightarrow \calC^{\otimes}_{J_{\ast}}$.
\item[$(b)$] Let $\delta: J_{\ast} \rightarrow J'_{\ast}$ be a morphism in $\LinSeg$ and
$\delta_!: \calC^{\otimes}_{J_{\ast}} \rightarrow \calC^{\otimes}_{J'_{\ast}}$ the associated functor. Then $\delta_! \circ \overline{f}_2: \Nerve(\calI)^{\triangleright} \rightarrow \calC^{\otimes}_{J'_{\ast}}$ is a colimit diagram.
\end{itemize}

Moreover, $(1')$ translates into the following condition:

\begin{itemize}
\item[$(c)$] Let $\overline{f}_2: \Nerve( \calI)^{\triangleright} \rightarrow \calC^{\otimes}_{J_{\ast}}$
be an arbitrary extension of $f_2$. Then $\overline{f}_2$ is a colimit diagram if and only if, for each $j \in J$, $\alpha^{j,J}_{!} \circ \overline{f}_2$ exhibits the image of the cone point as a coproduct
of the objects $(\alpha^{j,J}_{!} \circ f_2)(e^{j,J}_{n})$ in $\calC$, where $n$ ranges over all nonnegative integers.
\end{itemize}

We will begin with the proofs of $(a)$ and $(c)$, and return later to the proof of $(b)$.
Using the fact that $q: \calC^{\otimes} \rightarrow \Nerve(\LinSeg)$ is a Segal monoidal $\infty$-category, we see that $(a)$ and $(c)$ can be reformulated as follows:

\begin{itemize}
\item[$(a')$] For each $j \in J$, the map $f_3 = \alpha^{j,J}_{!} \circ f_2$
admits a colimit in $\calC$.
\item[$(c')$] An arbitrary map $\overline{f}_3: \Nerve( \calI)^{\triangleright} \rightarrow \calC$
extending $f_3$ is a colimit diagram if and only if it exhibits the image of the cone point as a coproduct of the objects $f_3( e^{j,J}_{n}) \in \calC$; here $n$ ranges over all nonnegative integers.
\end{itemize}

Let us now fix $j \in J$ and prove the assertions $(a')$ and $(c')$.
An object of $\calI$ can be identified with a quintuple
$( I, I_0, \seg{n}, \alpha, \alpha')$, where
$I_0 \subseteq I$ are linearly ordered sets, $\seg{n}$ is an object of $\LinSeg$, and
$\alpha: \seg{n} \rightarrow (I_0)_{\ast}$ and $\alpha': I_{\ast} \rightarrow J_{\ast}$ are morphisms
in $\LinSeg$. Let $\calI_0$ denote the full subcategory of $\calI$ spanned by those objects for which $\alpha'$ is an isomorphism. We observe that the inclusion $\calI_0 \subseteq \calI$ has a left adjoint, given by
$$ (I, I_0, \seg{n}, \alpha, \alpha') \mapsto (J, J_0, \seg{n}, \beta \circ \alpha, \id),$$
where $J_0 = \{ j \in J: (\alpha')^{-1} \{j\} \subseteq I_0 \}$ and
$\beta: (I_0)_{\ast} \rightarrow (J_0)_{\ast}$ has the property that
$\beta^{-1} \{j\} = (\alpha')^{-1} \{j\}$ for all $j \in J_0$. It follows that the inclusion
$\Nerve(\calI_0) \subseteq \Nerve(\calI)$ is cofinal. We therefore obtain the following reformulations of $(a')$ and $(c')$:

\begin{itemize}
\item[$(a'')$] The map $f_4 = f_3 | \Nerve(\calI_0)$ admits a colimit in $\calC$.
\item[$(c'')$] An arbitrary extension $\overline{f}_4: \Nerve(\calI_0)^{\triangleright} \rightarrow \calC$ of $f_4$ is a colimit diagram if and only if it exhibits the image of the cone point as a coproduct of the objects $f_4( e^{j,J}_{n}) \in \calC$; here $n$ ranges over all nonnegative integers.
\end{itemize}

We may identify objects of $\calI_0$ with triples $(J_0, \seg{n}, \alpha)$, where
$J_0$ is a subset of $J$, $\seg{n}$ is an object of $\LinSeg$, and 
$\alpha: \seg{n} \rightarrow (J_0)_{\ast}$ is a morphism in $\LinSeg$. Let
$\calI_1$ be the full subcategory of $\calI_0$ spanned by those objects for which
$j \in J_0$. Using condition $(A3)$ of Notation \ref{notebard}, we conclude that
if $X \in \calI_0$ and $X \notin \calI_1$, then $f_4(X)$ is an initial object of $\calC$.
It follows that $f_4$ is a left Kan extension of $f_5 = f_4 | \Nerve(\calI_1)$. 
Lemma \toposref{kan0} shows that $(a'')$ and $(c'')$ are equivalent to the following conditions:

\begin{itemize}
\item[$(a''')$] The map $f_5$ admits a colimit in $\calC$.
\item[$(c''')$] An arbitrary extension $\overline{f}_5: \Nerve(\calI_1)^{\triangleright} \rightarrow \calC$ of $f_5$ is a colimit diagram if and only if it exhibits the image of the cone point as a coproduct of the objects $f_5( e^{j,J}_{n}) \in \calC$; here $n$ ranges over all nonnegative integers.
\end{itemize}

Let $\calI_2$ denote the full subcategory of $\calI_1$ spanned by those objects
$( J_0, \seg{n}, \alpha)$ where $J_0 = \{j\}$. We observe that the inclusion
$ \calI_2 \subseteq \calI_1$ has a right adjoint, given by
$$ (J_0, \seg{n}, \alpha) \mapsto ( \{j\}, \seg{n}, \alpha^{j,J_0} \circ \alpha).$$
Using condition $(A2)$ of Notation \ref{notebard}, we observe that $f_5$ carries each
counit map $( \{j\}, \seg{n}, \alpha^{j, J_0} \circ \alpha) \rightarrow
(J_0, \seg{n}, \alpha)$ to an equivalence in $\calC$. Applying (the dual of) Lemma \ref{gottamake}, we obtain the following reformulations of $(a''')$ and $(c''')$: 

\begin{itemize}
\item[$(a'''')$] The map $f_6 = f_5 | \Nerve(\calI_2)$ admits a colimit in $\calC$.
\item[$(c'''')$] An arbitrary extension $\overline{f}_6: \Nerve(\calI_2)^{\triangleright} \rightarrow \calC$ of $f_6$ is a colimit diagram if and only if it exhibits the image of the cone point as a coproduct of the objects $f_6( e^{j,J}_{n}) \in \calC$; here $n$ ranges over all nonnegative integers.
\end{itemize}

Let $\calI_3$ denote the full subcategory of $\calI_2$ spanned by triples of the form
$( \{j\}, \seg{n}, \alpha)$, where $\alpha: \seg{n} \rightarrow \{j\}$ has the property that
$\alpha^{-1} \{\ast\} = \{ \ast \}$. The inclusion $\calI_3 \subseteq \calI_2$ has a left adjoint, so the induced inclusion $\Nerve(\calI_3) \subseteq \Nerve(\calI_2)$ is cofinal. We therefore obtain the following reformulations of $(a'''')$ and $(c'''')$: 

\begin{itemize}
\item[$(a''''')$] The map $f_7 = f_6 | \Nerve(\calI_3)$ admits a colimit in $\calC$.
\item[$(c''''')$] An arbitrary extension $\overline{f}_7: \Nerve(\calI_3)^{\triangleright} \rightarrow \calC$ of $f_7$ is a colimit diagram if and only if it exhibits the image of the cone point as a coproduct of the objects $f_7( e^{j,J}_{n}) \in \calC$; here $n$ ranges over all nonnegative integers.
\end{itemize}

We now observe that the category $\calI_3$ is discrete, and its objects can be identified
with nonnegative integers $n$ via the bijection
$ n \mapsto ( e^{j,J}_{n} \rightarrow J_{\ast} ).$
Consequently, assertion $(c''''')$ is a tautology, and $(a''''')$ follows from our assumption that
the $\infty$-category $\calC$ admits countable coproducts. This completes the proofs of $(a)$ and $(c)$.

We now return to the proof of $(b)$. Using the fact that $q: \calC^{\otimes} \rightarrow \Nerve(\LinSeg)$ is a Segal monoidal $\infty$-category, we see that $(b)$ is equivalent to the following slightly weaker statement:

\begin{itemize}
\item[$(b')$] Let $\beta: J_{\ast} \rightarrow \seg{1}$ be a morphism in $\LinSeg$ and
$\beta_!: \calC^{\otimes}_{J_{\ast}} \rightarrow \calC^{\otimes}_{\seg{1}} = \calC$ the associated functor. Then
$\overline{g}_3 = \delta_! \circ \overline{f}_2: \Nerve(\calI)^{\triangleright} \rightarrow \calC$ is a colimit diagram.
\end{itemize}

Noting once again that the inclusion $\Nerve(\calI_0) \subseteq \Nerve(\calI)$ is cofinal, we 
can reformulate $(b')$ as follows:

\begin{itemize}
\item[$(b'')$] The map $\overline{g}_4 = \overline{g}_3 | \Nerve(\calI_0)^{\triangleright} \rightarrow \calC$ is a colimit of $g_4 = \overline{g}_4 | \Nerve(\calI_0)$. 
\end{itemize}

The morphism $\beta: J_{\ast} \rightarrow \seg{1}$ is uniquely determined by the the subset
$J_1 = J - \beta^{-1} \{ \ast \} \subseteq J$. 
Let $\calI'_1$ denote the full subcategory of $\calI_0$ spanned by those objects
$(J_0, \seg{n}, \alpha)$ such that $J_1 \subseteq J_0$. Using Condition $(A3)$ of Notation \ref{notebard} and Remark \ref{plaque}, we deduce that if $(J_0, \seg{n}, \alpha) \in \calI_0$ does
not belong to $\calI'_1$, then $g_4( J_0, \seg{n}, \alpha)$ is an initial object of $\calC$. It follows that $g_4$ is a left Kan extension of $g_5 = g_4 | \Nerve(\calI'_1)$. Invoking Lemma \toposref{kan0}, we see that $(b'')$ is equivalent to:

\begin{itemize}
\item[$(b''')$] The map $\overline{g}_5 = \overline{g}_4 | \Nerve(\calI'_1)^{\triangleright} \rightarrow \calC$ is a colimit diagram.
\end{itemize}

Let $\calI'_2$ denote the full subcategory of $\calI'_1$ spanned by those objects
$(J_0, \seg{n}, \alpha)$ for which $J_0 = J_1$. The inclusion $\calI'_2 \subseteq \calI'_1$ admits a right adjoint $V$. Moreover, condition $(A3)$ of Notation \ref{notebard} guarantees that
$g_5$ carries each of the unit maps $V(J_1, \seg{n}, \alpha) \rightarrow (J_0, \seg{n}, \alpha)$ to an equivalence in $\calC$. Applying (the dual of) Lemma \ref{gottamake}, we deduce that
$(b''')$ is equivalent to:

\begin{itemize}
\item[$(b'''')$] The map $\overline{g}_6 = \overline{g}_5 | \Nerve(\calI'_2)^{\triangleright} \rightarrow \calC$ is a colimit diagram.
\end{itemize}

Finally, let $\calI'_3$ denote the full subcategory of $\calI'_2$ spanned by those objects
$(J_1, \seg{n}, \alpha)$ for which $\alpha^{-1} \{\ast\} = \{ \ast \}$. The inclusion
$\calI'_3 \subseteq \calI'_2$ admits a left adjoint, so the induced inclusion
$\Nerve( \calI'_3) \subseteq \Nerve(\calI'_2)$ is cofinal. We therefore obtain the following final reformulation of $(b'''')$: 

\begin{itemize}
\item[$(b''''')$] The map $\overline{g}_7 = \overline{g}_6 | \Nerve(\calI'_3)^{\triangleright} \rightarrow \calC$ is a colimit diagram.
\end{itemize}

We now observe that the category $\calI'_3$ is discrete. If we let
$J_1 = \{ j_0, \ldots, j_k \}$, then the objects of $\calI'_3$ are in bijection with finite sequences
of nonnegative integers $( n_0, \ldots, n_k)$. The bijection is given by the formula
$$ (n_0, \ldots, n_k) \mapsto ( J_1, \{ n_0 + \ldots + n_k \}, \alpha ),$$
where $\alpha^{-1} \{ j_i \} = \{ m | n_0 + \ldots + n_{i-1} \leq m < n_0 + \ldots + n_i \}$. 
Using conditions $(A1)$ and $(A4)$ of Notation \ref{notebard}, we conclude that
$\overline{g}_7( n_0, \ldots, n_k)$ is canonically isomorphic to the tensor product
$$ f_1( e^{j_0, J}_{n_0}) \otimes \ldots \otimes f_1( e^{ j_k, J}_{n_k} )$$
in the homotopy category $\h{\calC}$. On the other hand, the same argument and condition $(c)$
imply that the image of the cone point under $\overline{g}_7$ is equivalent to the tensor product
$$ ( \coprod_{n_0 \geq 0} f_1( e^{j_0, J}_{n_0} ) )
\otimes \ldots \otimes ( \coprod_{ n_k \geq 0 } f_1( e^{j_k, J}_{n_k} ) ).$$
Condition $(b''''')$ now follows from our assumption that the tensor product on $\calC$ preserves countable coproducts in each variable. This completes our proof of $(b)$, and therefore also our proofs of $(1)$ and $(1')$.

For the remainder of the proof, we will assume that $f \in \bHom_{\Nerve(\LinSeg) }( \Nerve(\calJ_3), \calC^{\otimes})$ is a $q$-left Kan extension of $f_0$. Let us prove $(2)$.
We wish to show that, if $\alpha: J_{\ast} \rightarrow J'_{\ast}$ is a morphism in
$\LinSeg$ which induces a bijection from $\alpha^{-1} J'$ onto $J'$, then $f(\alpha)$ is a $q$-coCartesian morphism in
$\calC^{\otimes}$. Let $j'$ be an element of $J'$, and $j \in J$ its preimage under $\alpha$.
In view of our assumption that
$q: \calC^{\otimes} \rightarrow \Nerve( \LinSeg)$ is a Segal monoidal $\infty$-category, it will suffice to show that $\alpha$ induces (for every choice of $j' \in J'$) an equivalence
$ \alpha^{j,J}_{!} f( J_{\ast} ) \rightarrow \alpha^{j', J'}_{!} f(J'_{\ast}).$
Invoking $(1)$, we see that $f$ induces equivalences
$$ \alpha^{j,J}_{!} f( J_{\ast} ) \simeq \coprod_{ n \geq 0} \alpha^{j,J}_{!}( f_0( e^{j,J}_{n} ) ) \quad \alpha^{j',J'}_{!} f( J'_{\ast} ) \simeq \coprod_{ n \geq 0} \alpha^{j',J'}_{!}( f_0( e^{j',J'}_{n} ) ).$$
It will therefore suffice to show that each of the maps
$ \alpha^{j,J}_{!} f_0( e^{j,J}_{n}) \rightarrow \alpha^{j',J'}_{!} f_0(e^{j',J'}_{n})$
is an equivalence; this follows from $(A2)$ of Notation \ref{notebard}.

We now prove $(3)$. Invoking $(1')$, we deduce that $f$ exhibits $A( \seg{1} )$ as a coproduct
of the objects $f_0(\nostar{1}, \nostar{1}, \seg{n}, \alpha(n) )$, where $\alpha(n)$ denotes the 
unique map $\seg{n} \rightarrow \seg{1}$ such that $\alpha^{-1} \{ \ast \} = \{ \ast \}$. 
The functor $f$ induces a commutative diagram
$$ \xymatrix{ f_0(\nostar{n}, \nostar{n}, \seg{n}, \id) \ar[r] \ar[d] & f_0( \nostar{1}, \nostar{1}, \seg{n}, \alpha(n) ) \ar[d] \\
A( \seg{n} ) \ar[r] & A( \seg{1} ) }$$
in $\calC^{\ast}$. Condition $(A1)$ of Notation \ref{notebard} guarantees that the left vertical map corresponds, under the equivalence $\calC^{\otimes}_{\seg{n}} \simeq \calC^{n}$, to the $n$th power of the map $C \rightarrow A( \seg{1} )$. 
Let $\alpha(n)_!: \calC^{\otimes}_{\seg{n}} \rightarrow \calC^{\otimes}_{\seg{1}} = \calC$ be the functor induced by $\alpha(n)$. Then the above diagram
induces a commutative square
$$ \xymatrix{ \alpha(n)_! f_0(\nostar{n}, \nostar{n}, \seg{n}, \id ) \ar[r]^{u} \ar[d] & f_0( \nostar{1}, \nostar{1}, \seg{1}, \alpha(n) ) \ar[d] \\
\alpha(n)_{!} A( \seg{n} ) \ar[r] & A( \seg{1} ) }$$
in $\calC$. Condition $(A4)$ of Notation \ref{notebard} guarantees that the map
$u$ is an equivalence. Consequently, the map $f_0(\nostar{1}, \nostar{1}, \seg{n}, \alpha(n) ) \rightarrow A( \seg{1} )$ can
be identified with the composition
$$ C^{\otimes n} \rightarrow A( \seg{1} )^{\otimes n} \simeq \alpha(n)_{!} A( \seg{n} ) \rightarrow A( \seg{1} ),$$
so that $(3)$ follows from $(1')$.
\end{proof}

\begin{notation}
Let $\calC$ be a monoidal $\infty$-category. We define $\Spread_3(\calC)$ to be the full subcategory of $$\bHom_{ \Nerve(\LinSeg) }( \Nerve(\calJ_3), \calC^{\otimes})$$ spanned
by those functors $f \in \bHom_{ \Nerve(\LinSeg) }( \Nerve(\calJ_3), \calC^{\otimes})$
with the following properties:
\begin{itemize}
\item[$(i)$] The restriction $f| \Nerve(\calJ_2)$ belongs to $\Spread_2(\calC)$.
\item[$(ii)$] The restriction $f| \Nerve(\LinSeg)$ belongs to $\SegAlg(\calC)$.
\end{itemize}
\end{notation}

\begin{lemma}\label{abletime}
Let $q: \calC^{\otimes} \rightarrow \Nerve(\LinSeg)$ be a Segal monoidal $\infty$-category. Then:
\begin{itemize}
\item[$(1)$] Composition with
$\psi_3: \calJ_3 \rightarrow \LinSeg$ and the inclusion $\LinSeg \rightarrow \calJ_3$ defines a 
pair of adjoint functors
$$ \Adjoint{ \phi'_4}{ \bHom_{ \Nerve(\LinSeg)}( \Nerve(\calJ_3), \calC^{\otimes})}{ \bHom_{ \Nerve(\LinSeg)}( \Nerve(\LinSeg), \calC^{\otimes})}{ \theta'_4}.$$
\item[$(2)$] The functors $\phi'_4$ and $\theta'_4$ restrict to a pair of adjoint functors
$\Adjoint{ \phi_4}{ \Spread_3(\calC) }{ \SegAlg(\calC) }{\theta_4}.$
\end{itemize}
\end{lemma}

\begin{proof}
In view of Proposition \toposref{leftkanadj}, to prove $(1)$ it will suffice to show that
$\theta'_4$ is a $q$-right Kan extension functor. In other words, it will suffice to show that
if $A \in \bHom_{ \Nerve(\LinSeg)}( \Nerve(\LinSeg), \calC^{\otimes})$, then
$A' = A \circ \psi_3$ is a $q$-right Kan extension of $A$. Fix an object $(J,J_0, \seg{n}, \alpha)  \in \calJ_2$. We observe that
$\LinSeg \times_{ \calJ_3} ( \calJ_3)_{/J}$ has a final object, given by $J_{\ast} \in \LinSeg$.
Consequently, $A'$ is a $q$-right Kan extension of $A$ at
$J$ if and only if the induced map $A'(J, J_0, \seg{n}, \alpha) \rightarrow A'(J_{\ast})$ is $q$-Cartesian. We now observe that this map is an equivalence by construction.

To prove $(2)$, it will suffice to show that
$\phi'_4( \Spread_3(\calC) ) \subseteq \SegAlg(\calC)$ and
$\theta'_4( \SegAlg(\calC) ) \subseteq \Spread_3(\calC)$. The first inclusion is obvious, and the second follows from Remark \ref{umpertine}.
\end{proof}

\begin{proof}[Proof of Theorem \ref{hutmunn}]
The equivalence of $(1)$ and $(2)$ follows from Proposition \toposref{simpex}. We will prove $(2)$.
Using Remark \ref{segmonex}, we may assume without loss of generality that the monoidal structure on $\calC$ is the restriction of a {\em Segal} monoidal category $q: \calC^{\otimes} \rightarrow \Nerve(\LinSeg)$. Let $\theta: \SegAlg(\calC) \rightarrow \calC$ denote the forgetful functor. Using Proposition \ref{algcompare}, we are reduced to proving the following analogue of $(2)$:

\begin{itemize}
\item[$(2')$] For every $C \in \calC$, there exists a Segal algebra object $A' \in \SegAlg(\calC)$
and a map $\eta: C \rightarrow \theta(A')$ such that, for every $B' \in \SegAlg(\calC)$, composition with $\eta$ induces a homotopy equivalence
$ \bHom_{\SegAlg(\calC)}( A', B') \rightarrow \bHom_{\calC}( C, \theta(B') ).$
\end{itemize}

Fix an object $C \in \calC$. Using Lemmas \ref{smallbird}, \ref{firstpart}, and \ref{2ndpart}, we can choose $\overline{C} \in \Spread'_{2}(\calC)$ such that 
$\overline{C}(\nostar{1}, \nostar{1}, \seg{1}, \id_{ \seg{1} }) = C$. Using Lemma \ref{sumotin}, we can choose $\widetilde{C} \in \Spread_3(\calC)$ which is a $q$-left Kan extension of $\overline{C}$. Let 
$A' = \widetilde{C} | \Nerve(\LinSeg)$, so that $\widetilde{C}$ determines a morphism $C \rightarrow \theta(A')$ in $\calC$. We claim that this morphism has the desired property. Let $B'$ be an arbitrary object of $\SegAlg(\calC)$. The canonical map
$g: \bHom_{\SegAlg(\calC)}(A', B') \rightarrow \bHom_{\calC}( C, \theta(B') )$ factors as a composition
\begin{eqnarray*}
\bHom_{\SegAlg(\calC)}(A', B') & \stackrel{g_4}{\rightarrow} & \bHom_{\Spread_3(\calC)}( \widetilde{C}, B' \circ \psi_3) \\
& \stackrel{g_3}{\rightarrow} & \bHom_{\Spread_2(\calC)}( \overline{C}, B' \circ \psi_2) \\
& \stackrel{g_2}{\rightarrow} & \bHom_{\Spread_1(\calC)}( \overline{C} | \Nerve(\calJ_1), B' \circ \psi_1) \\
&  \stackrel{g_1}{\rightarrow} & \bHom_{\Spread_0(\calC)}( \overline{C} | \Nerve(\calJ_0), B' \circ \psi_0) \\
& \stackrel{g_0}{\rightarrow} & \bHom_{\calC}( C, \theta(B') ).
\end{eqnarray*}
The maps $g_0$, $g_1$, $g_2$, $g_3$, and $g_4$ are all homotopy equivalences
(Lemmas \ref{smallbird}, \ref{firstpart}, \ref{2ndpart}, \toposref{kan1}, and \ref{abletime}, respectively). It follows that $g$ is a homotopy equivalence, as desired. This completes the proof of $(2')$. Taking $A \in \Alg(\calC)$ to be the composition of $A'$ with the map
$\Nerve(\cDelta)^{op} \rightarrow \Nerve(\LinSeg)$, we deduce $(2)$.

The third assertion of Lemma \ref{sumotin} implies that the canonical map
$\coprod_{n} C^{\otimes n} \rightarrow A([1]) \simeq \theta(A')$ is an equivalence. This (together with observation that free algebras generated by $C$ are uniquely determined up to equivalence) proves the ``only if'' direction of $(3)$. For the converse, suppose given an algebra object $B \in \Alg(\calC)$ and a map $C \rightarrow B([1])$ which induces an equivalence $\coprod_{n} C^{\otimes n} \rightarrow B([1])$. Since $A$ is freely generated by $C$, there exists a map
$f: A \rightarrow B$ in $\Alg(\calC)$ and a commutative triangle
$$ \xymatrix{ & A([1]) \ar[dr]^{f'} & \\
C \ar[ur] \ar[rr] & & B([1])  }$$
in $\calC$. Consider the induced diagram
$$ \xymatrix{ & A([1]) \ar[dr]^{f'} & \\
\coprod_{n \geq 0} C^{\otimes n} \ar[rr] \ar[ur] & & B([1]). }$$
Using the two-out-of-three property, we deduce that $f'$ is an equivalence.
Corollary \ref{jumunj} now implies that $f$ is an equivalence, so that $B$ is freely generated by $C$ as desired.
\end{proof}

\begin{remark}
The proof of Theorem \ref{hutmunn} can be adapted to prove the following analogue of
Theorem \ref{hutmunn}:
\begin{itemize}
\item[$(\ast)$]  Suppose that $\calC$ is a monoidal $\infty$-category which admits countable coproducts, and that the tensor product $\otimes: \calC \times \calC \rightarrow \calC$ preserves countable coproducts separately in each variable. Then the forgetful functor
$\theta: \Alg^{\nounit}(\calC) \rightarrow \calC$ admits a left adjoint (here
$\Alg^{\nounit}(\calC)$ denotes the $\infty$-category of nonunital algebra objects of $\calC$;
see Definition \ref{nounitalg}). Moreover, the composition of 
$\theta$ with this left adjoint is (canonically) identified with the functor
$C \mapsto \coprod_{ n > 0} C^{\otimes n}.$
\end{itemize}
In fact, assertion $(\ast)$ is slightly easier to prove that Theorem \ref{hutmunn}, since it is possible to avoid the formalism of Segal monoidal categories. We leave the details to reader.
\end{remark}

\subsection{Limits and Colimits of Algebras}\label{limalg}

Let $\calC$ be a monoidal $\infty$-category. Our goal in this section is to prove the existence of limits and colimits in $\Alg(\calC)$, given suitable assumptions on $\calC$. The case of limits is straightforward: our main result is Corollary \ref{slimycomp}, which asserts that limits in $\Alg(\calC)$ can be computed at the level of the underlying objects of $\calC$. 

Colimits in $\Alg(\calC)$ are more complicated. One might naively guess that one can form a colimit of a diagram $K \rightarrow \Alg(\calC)$ by first forming a colimit of the underlying diagram $K \rightarrow \calC$, and then endowing it with an algebra structure. However, this obviously fails in some simple cases. For example, if $\calC$ is the (nerve of the) category of sets, equipped with its Cartesian monoidal structure then $\Alg(\calC)$ is (the nerve of) the category of associative monoids. In this case, both $\calC$ and $\Alg(\calC)$ admit coproducts, but the forgetful functor $\theta: \Alg(\calC) \rightarrow \calC$ does not preserve coproducts: if $M$ and $N$ are monoids, then the disjoint union $M \coprod N$ does not inherit the structure of a monoid. 

Nevertheless, there are {\em some} colimits in $\Alg(\calC)$ which can (often) be computed at the level of the underlying objects of $\calC$. Namely, we will see that this is true when $K$ is a {\em sifted} simplicial set (Definition \stableref{siftdef}); see Proposition \ref{fillfem} for a precise statement.

The construction of general colimits in $\Alg(\calC)$ is more difficult. Using general arguments, we can reduce to problem of constructing sifted colimits (which is addressed by Proposition \ref{fillfem}) and the problem of constructing coproducts. To handle the latter, we will not proceed directly. 
Instead, we will use the results of \S \ref{monoid7} and \S \ref{monoid5} to resolve arbitrary objects of $\Alg(\calC)$ by free algebras. We can then reduce to the problem of constructing coproducts of free algebras, which are easily computed in terms of coproducts in $\calC$. 

We now begin our analysis by considering limits in $\Alg(\calC)$.

\begin{lemma}\label{surinor}
Let $\calC^{\otimes}$ and $\calD^{\otimes}$ be monoidal structures on $\infty$-categories
$\calC$ and $\calD$, respectively, and let $\alpha: F \rightarrow F'$ be a morphism in
$\Fun^{\wMon}(\calC^{\otimes}, \calD^{\otimes})$. The following are equivalent:
\begin{itemize}
\item[$(1)$] The transformation $\alpha$ is an equivalence in $\Fun^{\wMon}(\calC^{\otimes}, \calD^{\otimes})$. 
\item[$(2)$] For every object $C \in \calC$, the morphism $\alpha(C)$ is an equivalence in $\calD$.
\end{itemize}
\end{lemma}

\begin{proof}
The implication $(1) \Rightarrow (2)$ is obvious. Conversely, suppose that $(2)$ is satisfied. Let
$X \in \calC^{\otimes}_{[n]}$; we wish to prove that $\alpha(X)$ is an equivalence in
$\calD^{\otimes}_{[n]}$. Since $\calD^{\otimes}$ is a monoidal $\infty$-category, it suffices to show that for $0 \leq i < n$, the image of $\alpha(X)$ under the map 
$\calD^{\otimes}_{[n]} \rightarrow \calD^{\otimes}_{[1]} \simeq \calD$ induced by the inclusion
$[1] \simeq \{i, i+1\} \hookrightarrow [n]$ is an equivalence in $\calD$. Since $F$ and $F'$
are lax monoidal functors, this morphism can be identified with $\alpha(X_i)$, where
$X_i$ is the image of $X$ under the corresponding map $\calC^{\otimes}_{[n]} \rightarrow \calC^{\otimes}_{[1]} \simeq \calC$. The desired result now follows immediately from $(2)$.
\end{proof}

\begin{proposition}\label{limycomp}
Let $\calC^{\otimes}$ and $\calD^{\otimes}$ be monoidal structures on $\infty$-categories
$\calC$ and $\calD$, respectively, and let $q: K \rightarrow \Fun^{\wMon}( \calC^{\otimes}, \calD^{\otimes})$. Suppose that, for every object $C \in \calC$, the induced diagram
$K \rightarrow \calD$ admits a limit in $\calD$. Then:

\begin{itemize}
\item[$(1)$] The diagram $q$ has a limit in $\Fun^{\wMon}(\calC^{\otimes}, \calD^{\otimes})$. 

\item[$(2)$] An arbitrary extension $\overline{q}: K^{\triangleleft} \rightarrow \Fun^{\wMon}(\calC^{\otimes}, \calD^{\otimes})$ of $q$ is a limit diagram if and only if, for every object $C \in \calC$, the induced map $K^{\triangleleft} \rightarrow \calD$ is a limit diagram.
\end{itemize}
\end{proposition}

\begin{proof}
We will prove that there exists a limit diagram $\overline{q}: K^{\triangleleft} \rightarrow \Fun^{\wMon}(\calC^{\otimes}, \calD^{\otimes})$ satisfying the condition of $(2)$. This will prove $(1)$. Moreover, the ``only if'' direction of $(2)$ will follow from the uniqueness of limit diagrams up to equivalence. To prove the ``if'' direction of $(2)$, suppose that $\overline{q}': K^{\triangleleft} \rightarrow \Fun^{\wMon}(\calC^{\otimes}, \calD^{\otimes})$ is an arbitrary extension of $q$ which satisfies the condition of $(2)$. Since $\overline{q}$ is a limit diagram, we obtain a natural transformation $\alpha: \overline{q}' \rightarrow \overline{q}$ in $\Fun^{\wMon}(\calC^{\otimes}, \calD^{\otimes})_{/q}$. It follows from Lemma \ref{surinor} that $\alpha$ is an equivalence.

For every object $X \in \calC^{\otimes}_{[n]}$, let $q_{X}: K \rightarrow \calD^{\otimes}_{[n]}$ be the diagram induced by $q$. We observe that if $X$ corresponds to a sequence of objects
$(C_1, \ldots, C_n)$ under the equivalence $\calC^{\otimes}_{[n]} \simeq \calC^{n}$, then $q_X$ has a limit, which projects under the equivalence $\calD^{\otimes}_{[n]} \simeq \calD^n$ to a collection of limits for the diagrams $\{ q_{C_i} \}_{1 \leq i \leq n}$. Applying Proposition \toposref{limiteval}, we conclude that $q$ admits a limit in the $\infty$-category $\bHom_{\Nerve(\cDelta)^{op}}( \calC^{\otimes}, \calD^{\otimes})$, and that the limit of $q$ belongs to 
$\Fun^{\wMon}( \calC^{\otimes}, \calD^{\otimes})$. It is readily verified that this limit has the desired properties.
\end{proof}

\begin{corollary}\label{slimycomp}\index{algebra object!limit of}
Let $\calD$ be an $\infty$-category equipped with a monoidal structure, let $\theta: \Alg(\calD) \rightarrow \calD$ be the forgetful functor, and let $q: K \rightarrow \Alg(\calD)$ be a diagram. Suppose that $\theta \circ q$ has a limit in $\calD$. Then:
\begin{itemize}
\item[$(1)$] The diagram $\theta$ has a limit in $\Alg(\calD)$. 

\item[$(2)$] An arbitrary extension $\overline{q}: K^{\triangleleft} \rightarrow \Alg(\calD)$ is a limit diagram if and only if $\theta \circ \overline{q}$ is a limit diagram.
\end{itemize}
\end{corollary}

\begin{proof}
Apply Proposition \ref{limycomp} in the case where $\calC^{\otimes} = \Nerve(\cDelta)^{op}$.
\end{proof}

\begin{corollary}\label{jumunj}
Let $\calC$ be an $\infty$-category equipped with a monoidal structure. Then the forgetful functor
$\theta: \Alg(\calC) \rightarrow \calC$ is conservative. In other words, a morphism
$f: A \rightarrow A'$ of algebra objects of $\calC$ is an equivalence in $\Alg(\calC)$ if and only if
$\theta(f)$ is an equivalence in $\calC$.
\end{corollary}

\begin{corollary}\label{firebaugh}\index{algebra object!final}
Let $\calC$ be a monoidal $\infty$-category, and suppose that $\calC$ admits a final object. 
Then:
\begin{itemize}
\item[$(1)$] The $\infty$-category $\Alg(\calC)$ admits a final object.
\item[$(2)$] An object $A \in \Alg(\calC)$ is final if and only if its image in $\calC$ is final.
\end{itemize}
\end{corollary}

In other words, if $\calC$ is a monoidal $\infty$-category which admits a final object $U$, then $U$ can be promoted (in an essentially unique way) to an algebra object of $\calC$. 

\begin{remark}\label{firebaugh2}
Let $\calC$ be a monoidal $\infty$-category which admits a final object. Corollary \ref{firebaugh} has the following analogue for nonunital algebras, which can be proven using a similar argument (see \S \ref{digunit} for an explanation of the notation):
\begin{itemize}
\item[$(1)$] The $\infty$-category $\Alg^{\nounit}(\calC)$ admits a final object.
\item[$(2)$] An object $A \in \Alg^{\nounit}(\calC)$ is final if and only if its image in $\calC$ is final.
\end{itemize}
\end{remark}

\begin{definition}\label{indexco}\index{colimit!compatible with monoidal structure}\index{monoidal $\infty$-category!and compatible colimits}
Let $K$ be a simplicial set, and let $\calC$ be an $\infty$-category which admits $K$-indexed colimits. We will say that a monoidal structure on $\calC$ is {\it compatible with $K$-indexed colimits} if, for every object $C \in \calC$, the functors 
$$ \bigdot \otimes C: \calC \rightarrow \calC \quad C \otimes \bigdot: \calC \rightarrow \calC$$
preserve $K$-indexed colimits.
\end{definition}

In the case where $K$ is {\em sifted} (see Definition \stableref{siftdef}), Definition \ref{indexco} admits a convenient reformulation.

\begin{lemma}\label{hungerin}
Let $K$ be a sifted simplicial set, and let $\calC^{\otimes} \rightarrow \Nerve(\cDelta)^{op}$ be a monoidal $\infty$-category which is compatible with $K$-indexed colimits. For each morphism $f: [m] \rightarrow [n]$ in $\Nerve(\cDelta)^{op}$, the associated functor $\theta: \calC^{\otimes}_{[n]} \rightarrow \calC^{\otimes}_{[m]}$ preserves $K$-indexed colimits.
\end{lemma}

\begin{proof}
Let $k$ be the largest element of the set $\{ f(i+1)-f(i) \}_{0 \leq i < m}$; we will work by induction on $k$. Let $\calC = \calC^{\otimes}_{[1]}$. Using the equivalence $\calC^{\otimes}_{[m]} \simeq \calC^{m}$, we can reduce to the case $m=1$. Similarly, we may reduce to the case where the map $f: [1] \rightarrow [n]$ preserves initial and final objects. If $n > 2$, then we can factor $f$ as a composition $[1] \stackrel{f'}{\rightarrow} [2] \stackrel{f''}{\rightarrow} [n]$, and the desired result holds for $f'$ and $f''$ by the inductive hypothesis. We may therefore suppose that $n \leq 2$. If $n =1$, then $\theta$ is an equivalence and there is nothing to prove. If $n = 0$, then $\theta$ is equivalent to a constant map; the desired result then follows from Corollary \toposref{silt}, since $K$ is weakly contractible. When $n=2$, we apply Proposition \stableref{urbil}. 
\end{proof}


\begin{lemma}\label{hunkerin}
Let $p: X \rightarrow S$ be a coCartesian fibration of simplicial sets, and let $K$ be an arbitrary simplicial set. Suppose that:
\begin{itemize}
\item[$(i)$] For each vertex $s$ of $S$, the fiber $X_{s}$ admits $K$-indexed colimits.
\item[$(ii)$] For each edge $s \rightarrow s'$ of $S$, the associated functor
$X_{s} \rightarrow X_{s'}$ preserves $K$-indexed colimits. 
\end{itemize}
Then:
\begin{itemize}
\item[$(1)$] Every diagram $K \rightarrow \bHom_{S}(S,X)$ admits a colimit.
\item[$(2)$] An arbitrary map $K^{\triangleright} \rightarrow \bHom_{S}(S,X)$ is a colimit diagram
if and only, for each $s \in S$, the associated map $K^{\triangleright} \rightarrow X_{s}$ is a colimit diagram.
\end{itemize}
\end{lemma}

\begin{proof}
Combine Proposition \toposref{relcolfibtest} with (the dual of) Lemma \ref{surtybove}.
\end{proof}

\begin{proposition}\label{fillfem}\index{algebra object!sifted colimit of}
Let $\calC^{\otimes}$ and $\calD^{\otimes}$ be monoidal structures on $\infty$-categories
$\calC$ and $\calD$, respectively. Let $K$ be a sifted simplicial set. Suppose that
$\calD$ admits $K$-indexed colimits, and that the monoidal structure on $\calD$ is compatible with $K$-indexed colimits. Then:
\begin{itemize}
\item[$(1)$] The $\infty$-category $\bHom_{ \Nerve(\cDelta)^{op} }( \calC^{\otimes}, \calD^{\otimes})$ admits $K$-indexed colimits.
\item[$(2)$] The full subcategories 
$$\Fun^{\Mon}(\calC^{\otimes}, \calD^{\otimes}) \subseteq 
\Fun^{\wMon}( \calC^{\otimes}, \calD^{\otimes}) \subseteq 
\bHom_{ \Nerve(\cDelta)^{op} }( \calC^{\otimes}, \calD^{\otimes})$$ are stable under $K$-indexed colimits.
\item[$(3)$] The forgetful functor $\theta: \Fun^{\wMon}(\calC^{\otimes}, \calD^{\otimes}) \rightarrow \Fun(\calC, \calD)$ detects $K$-indexed colimits. More precisely, a map $\overline{q}: K^{\triangleright} \rightarrow \Fun^{\wMon}(\calC^{\otimes}, \calD^{\otimes})$ is a colimit diagram if and only if $\theta \circ \overline{q}$ is a colimit diagram.
\end{itemize}
\end{proposition}

\begin{proof}
Assertions $(1)$, $(2)$, and the ``only if'' part of $(3)$ follow immediately from Lemmas \ref{hungerin} and \ref{hunkerin}. To prove the ``if'' direction of $(3)$, let us suppose that
$\overline{q}: K^{\triangleright} \rightarrow \Fun^{\wMon}( \calC^{\otimes}, \calD^{\otimes})$ is
is such that $\theta \circ \overline{q}$ is a colimit diagram. We wish to prove that $\overline{q}$ is a colimit diagram. In view of $(2)$ and Lemma \ref{hunkerin}, it will suffice to show that for each object
$X \in \calC^{\otimes}_{[n]}$, the 
induced map $\overline{q}_{X}: K^{\triangleright} \rightarrow \calD^{\otimes}_{[n]}$ is a colimit diagram. In view of the equivalence
$\calD^{\otimes}_{[n]} \rightarrow \calD^{n}$, it will suffice to show that each of the composite maps
$$K^{\triangleright} \stackrel{ \overline{q}_{X} }{\rightarrow} \calD^{\otimes}_{[n]} \simeq \calD^n \rightarrow \calD$$
is a colimit diagram. Since the values assumed by $\overline{q}$ are lax monoidal functors, these compositions can be identified with $K^{\triangleright} \stackrel{q_{X_i}}{\rightarrow}
\calD$, where $X$ corresponds to $(X_1, \ldots, X_n) \in \calC^{n}$ under the equivalence $\calC^{\otimes}_{[n]} \simeq \calC^{n}$.
\end{proof}

\begin{corollary}\label{filtfem}
Let $K$ be a sifted simplicial set, let $\calD$ be an $\infty$-category which admits $K$-indexed colimits, and let $\calD^{\otimes} \rightarrow \Nerve(\cDelta)^{op}$ be a monoidal structure on
$\calD$ which is compatible with $K$-indexed colimits. Then:
\begin{itemize}
\item[$(1)$] The $\infty$-category $\Alg(\calD)$ admits $K$-indexed colimits.

\item[$(2)$] The forgetful functor $\theta: \Alg(\calD) \rightarrow \calD$ detects $K$-indexed colimits. More precisely, a map $\overline{q}: K^{\triangleright} \rightarrow \Alg(\calD)$ is a colimit diagram if and only if $\theta \circ \overline{q}$ is a colimit diagram.
\end{itemize}
\end{corollary}

\begin{proof}
Apply Proposition \ref{fillfem} in the case $\calC^{\otimes} = \Nerve(\cDelta)^{op}$.
\end{proof}

Let us now turn to the problem of constructing general colimits in $\Alg(\calC)$. In view of Corollary \toposref{uterrr} and Lemma \stableref{simpenough}, arbitrary colimits in $\Alg(\calC)$ can be built out of filtered colimits, geometric realizations, and finite coproducts. The cases of filtered colimits and geometric realizations are addressed by Corollary \ref{filtfem} (in view of Examples \stableref{bin1} and \stableref{bin2}). It will therefore sufice construct finite coproducts in $\Alg(\calC)$. 

\begin{proposition}\label{coprodmake}\index{algebra object!arbitrary colimit of}
Let $\kappa$ be an uncountable regular cardinal.
Let $\calC$ be an $\infty$-category which admits $\kappa$-small colimits, endowed with a monoidal structure which is compatible with $\kappa$-small colimits. Then:
\begin{itemize}
\item[$(1)$] The forgetful functor $\theta: \Alg(\calC) \rightarrow \calC$ has a left adjoint $\psi$.
\item[$(2)$] For every object $A \in \Alg(\calC)$, there exists a simplicial object $A_{\bigdot}$
of $\Alg(\calC)$ having $A$ as a colimit, such that each $A_{n}$ belongs to the essential image of $\psi$.
\item[$(3)$] The $\infty$-category $\Alg(\calC)$ admits $\kappa$-small colimits.
\end{itemize}
\end{proposition}

\begin{proof}
Part $(1)$ was established as Proposition \ref{hutmunn}, and $(2)$ follows from Corollary \ref{filtfem} and Proposition \ref{littlebeck}. We now prove $(3)$.
Since $\Alg(\calC)$ admits $\kappa$-small sifted colimits (Corollary \ref{filtfem}), it will suffice to show that $\Alg(\calC)$ admits $\kappa$-small coproducts (in fact, it suffices to treat the case of {\em finite} coproducts, but we will not need this). 
Choose a $\kappa$-small collection $\{ A^{\beta} \}_{\beta \in B}$ of objects of $\Alg(\calC)$. We wish to show that there exists a coproduct for this collection in $\Alg(\calC)$. Let $\psi: \calC \rightarrow \Alg(\calC)$ be a left adjoint to the inclusion functor. Let us say that an object of $\Alg(\calC)$ is {\it free} if it belongs to the essential image of $\psi$. Since $\psi$ preserves all colimits which exist in $\calC$ (Proposition \toposref{adjointcol}), the coproduct of the collection 
$\{ A^{\beta} \}_{\beta \in B}$ exists whenever each $A^{\beta}$ is free.

We now treat the general case. According to $(2)$, each $A^{\beta}$ can be obtained as the geometric realization of a simplicial object $A^{\beta}_{\bigdot}$ of $\Alg(\calC)$, where
each $A^{\beta}_{n}$ is free. Let $p: \coprod_{\beta \in B} \Nerve(\cDelta)^{op} \rightarrow \Alg(\calC)$ be the result of amalgamating all of these simplicial objects. We can now apply the methods of \S \toposref{quasilimit1} to decompose the diagram $p$ in two different ways:
\begin{itemize}
\item[$(i)$] Since each summand of $p$ has the object $A^{\beta}$ as a colimit, we can identify colimits of $p$ with coproducts of the family $\{ A^{\beta} \}_{\beta \in B}$. 
\item[$(ii)$] For each object $[n] \in \cDelta$, let $p_{[n]}$ denote the restriction of $p$
to $\coprod_{ \beta \in B} (\Nerve(\cDelta)^{op})_{/[n]}$, and let $p'_{[n]}$ be the restriction of
$p$ to $\coprod_{ \beta \in B} \{ [n] \}$. Since the inclusion 
$$ \coprod_{ \beta \in B } \{[n]\} \subseteq \coprod_{ \beta \in B } (\Nerve(\cDelta)^{op})_{/[n]}$$
is cofinal, we can identify colimits of $p_{[n]}$ with colimits of $p'_{[n]}$, which are coproducts
of the family $\{ A^{\beta}_{n} \}_{\beta \in B}$. The coproducts $A_{n} = \coprod_{\beta \in B} A^{\beta}_{n}$ exist since each $A^{\beta}_n$ is free. Using the methods of \S \toposref{quasilimit1}, we can organize the coproducts $A_{n}$ into a simplicial object $A_{\bigdot}$, such that colimits of $p$ can be identified with colimits of $A_{\bigdot}$. 
\end{itemize}
We now observe that the simplicial object $A_{\bigdot}$ of $\Alg(\calC)$ has a colimit in $\Alg(\calC)$, by virtue of Corollary \ref{filtfem}.
\end{proof}

\begin{corollary}\label{algcol}
Let $\calC$ be an $\infty$-category which admits small colimits, and suppose that
$\calC$ is endowed with a monoidal structure which is compatible with small colimits. Then
$\Alg(\calC)$ admits small colimits.
\end{corollary}

We now treat the question of whether or not $\Alg(\calC)$ is a presentable $\infty$-category.

\begin{lemma}\label{sillyfun}
Let $\calC^{\otimes} \rightarrow \Nerve(\cDelta)^{op}$ be a monoidal structure on an $\infty$-category $\calC = \calC^{\otimes}_{[1]}$. The following conditions are equivalent:
\begin{itemize}
\item[$(1)$] The $\infty$-category $\calC$ is accessible, and the monoidal structure on $\calC$ is compatible with $\kappa$-filtered colimits for all sufficiently large $\kappa$.
\item[$(2)$] Each fiber of $p$ is accessible, and for every morphism $[m] \rightarrow [n]$ in $\cDelta$, the associated functor $\calC^{\otimes}_{[n]} \rightarrow \calC^{\otimes}_{[m]}$ is accessible.
\end{itemize}
\end{lemma}

\begin{proof}
Suppose first that $(2)$ is satisfied. We deduce immediately that $\calC$ is accessible, and that
the tensor product $\otimes: \calC \times \calC \simeq \calC^{\otimes}_{[2]} \rightarrow \calC$
is an accessible functor. It now suffices to observe that for each $C \in \calC$, the inclusion
$\{C \} \times \calC \subseteq \calC \times \calC$ is also an accessible functor.

Now suppose that $(1)$ holds. Each fiber $\calC^{\otimes}_{[n]}$ is equivalent to an $n$-fold product of copies of $\calC$, and is therefore accessible. Now choose a morphism $[m] \rightarrow [n]$ in $\cDelta$; we wish to show that the associated functor $\calC^{\otimes}_{[n]} \rightarrow \calC^{\otimes}_{[m]}$ is accessible. Using the equivalences $\calC^{\otimes}_{[m]} \simeq \calC^{m}$ and $\calC^{\otimes}_{[n]} \simeq \calC^{n}$, we reduce easily to the case where $m=1$ and 
the map $[m] \rightarrow [n]$ preserves initial and final elements. The functor $\calC^{\otimes}_{[n]} \rightarrow \calC^{\otimes}_{[1]}$ can be identified with an iterated tensor product
$$ \otimes: \calC \times \ldots \times \calC \rightarrow \calC.$$
If $n=0$ then this functor is constant, if $n=1$ it is equivalent to the identity, and if $n=2$ it is accessible in view of our assumption. For $n > 2$, we can use the associativity of $\otimes$ to reduce to the case $n=2$.
\end{proof}

\begin{proposition}\label{algprec}\index{algebra object!of a presentable $\infty$-category}
Let $\calC$ be a monoidal $\infty$-category. Then:
\begin{itemize}
\item[$(1)$] If $\calC$ is accessible and the monoidal structure on $\calC$ is compatible with $\kappa$-filtered colimits, for $\kappa$ sufficiently large, then $\Alg(\calC)$ is an accessible $\infty$-category.
\item[$(2)$] If $\calC$ is a presentable $\infty$-category and the monoidal structure on $\calC$ is compatible with small colimits, then $\Alg(\calC)$ is a presentable $\infty$-category.
\end{itemize}
\end{proposition}

\begin{proof}
Assertion $(1)$ follows from Lemma \ref{sillyfun} and Proposition \toposref{prestorkus}. To prove $(2)$, we combine $(1)$ with Corollary \ref{algcol}.
\end{proof}

\subsection{Monoidal Model Categories}\label{monoidate}

Let $\calC$ be an $\infty$-category. In Example \ref{itereat}, we saw that if $\calC$ is the nerve of a monoidal category, then $\calC$ inherits the structure of a monoidal $\infty$-category. Our goal in this section is to obtain some generalizations of this result. We will be principally interested in the case where $\calC = \Nerve( \bfA^{\degree})$ is the $\infty$-category underlying a simplicial model category $\bfA$. We will show that if $\bfA$ admits a monoidal structure which is suitably compatible with its simplicial and model structures (Definition \ref{hyrt0}), then $\calC$ again inherits the structure of a monoidal $\infty$-category. Moreover, there is a close relationship between the $\infty$-category of algebra objects $\Alg(\calC)$ and the ordinary category of
(strictly associative) algebras in $\bfA$: see Theorem \ref{beckify}.

\begin{definition}\label{hyrt}\index{monoidal category!with compatible simplicial structure}
Let $\calC$ be simplicial category. We will say that a monoidal structure on $\calC$ is {\it weakly compatible} with the simplicial structure on $\calC$ provided that
the operation $\otimes: \calC \times \calC \rightarrow \calC$ is endowed with the structure of a simplicial functor, which is compatible with associativity and unit transformations.

We will say that a closed monoidal structure on $\calC$ is {\it compatible} with the simplicial structure on $\calC$ if it is weakly compatible and for every triple of objects $A,B,C \in \calC$, the natural maps
$$ \bHom( B,  {}^A\!C) \rightarrow \bHom(A \otimes B, C) \leftarrow \bHom(A, C^B)$$
are isomorphisms of simplicial sets (see \S \toposref{monoidaldef}). 
\end{definition}

\begin{remark}
The compatibility of Definition \ref{hyrt} is not merely a condition, but additional data which must be supplied. However, in practice the simplicial structure on the functor $\otimes$ tends to be clear in practice.
\end{remark}

Let $\calC$ be a simplicial category equipped with a weakly compatible monoidal structure, and let $\calC^{\otimes}$ the category described in Definition \ref{converi}. Then $\calC^{\otimes}$ is naturally endowed with the structure of a simplicial category, where we define the 
simplicial mapping space $\bHom_{\calC^{\otimes}}( [C_1, \ldots, C_n], [C'_1, \ldots, C'_m] )$
to be the disjoint union 
$$\coprod_{ f: [m] \rightarrow [n] } \prod_{1 \leq i \leq m} \bHom_{\calC}( C_{f(i-1)+1} \otimes \ldots \otimes C_{f(i)}, C'_i )$$

\begin{proposition}\label{hurgove}
Let $\calC$ be a fibrant simplicial category with a weakly compatible monoidal structure, and let
$\calC^{\otimes}$ be the simplicial category constructed above. Then the induced map
$p: \Nerve( \calC^{\otimes} ) \rightarrow \Nerve(\cDelta)^{op}$
is a monoidal structure on the $\infty$-category $\Nerve(\calC)$. 
\end{proposition}

\begin{proof}
We first show that $p$ is a coCartesian fibration. It is obviously an inner fibration, since
$\Nerve( \calC^{\otimes})$ is an $\infty$-category and $\Nerve( \cDelta)^{op}$ is the nerve of an ordinary category. Choose an object
$[C_1, \ldots, C_n] \in \calC^{\otimes}$ and a morphism $f: [m] \rightarrow [n]$ in $\cDelta$.
Then there exists a morphism $\overline{f}: [C_1, \ldots, C_n] \rightarrow [C'_{1}, \ldots, C'_{m}]$ in
$\Nerve( \calC^{\otimes})$ which covers $f$, where $C'_{i} \simeq C_{ f(i-1) + 1} \otimes \ldots \otimes C_{f(i)}$. It follows immediately from Proposition \toposref{trainedg} that $\overline{f}$ is $p$-coCartesian. We conclude the proof by observing that the fiber $\Nerve(\calC^{\otimes})_{[n]}$ is isomorphic to $\Nerve(\calC)^{n}$, and the projection onto the $i$th factor can be identified with
the functor induced by the inclusion $[1] \simeq \{ i-1, i\} \subseteq [n]$. 
\end{proof}

In practice, the hypotheses of Proposition \ref{hurgove} are often too strong. 
Suppose, for example, that $\bfA$ is a simplicial model category, and let $\bfA^{\degree}$ be the full (simplicial) subcategory of $\bfA$ spanned by the fibrant-cofibrant objects. Then $\bfA^{\degree}$ is a fibrant simplicial category, and we refer to the simplicial nerve $\Nerve( \bfA^{\degree})$ as the underlying $\infty$-category of $\bfA$. 
A monoidal structure on $\bfA$ need not restrict to a monoidal structure on $\bfA^{\degree}$, even if it is compatible with the model structure of $\bfA$ (see Definition \ref{hyrt0} below). Nevertheless, we can use a variant of Proposition \ref{hurgove} to endow $\Nerve( \bfA^{\degree})$ with the structure of a monoidal $\infty$-category.

\begin{definition}\label{hyrt0}\index{monoidal model category}\index{model category!monoidal}
Let $\bfA$ be a simplicial model category. We will say that a monoidal structure $\otimes$ on $\bfA$ is {\it compatible} with the simplicial model structure if it is compatible with the simplicial structure
on $\bfA$ (Definition \ref{hyrt}) and satisfies the following conditions:
\begin{itemize}
\item[$(1)$] For every pair of cofibrations $i: A \rightarrow A'$, $j: B \rightarrow B'$ in $\bfA$, the induced map
$$k: (A \otimes B') \coprod_{ A \otimes B} (A' \otimes B) \rightarrow A' \otimes B'$$
is a cofibration. Moreover, if either $i$ or $j$ is a weak equivalence, then $k$ is a weak equivalence.

\item[$(2)$] The unit object $1$ of $\bfA$ is cofibrant.

\item[$(3)$] The monoidal structure on $\bfA$ is closed. 
\end{itemize}
\end{definition}

\begin{proposition}\label{hurgoven}
Let $\bfA$ be a simplicial model category with a compatible monoidal structure, let
$\bfA^{\otimes}$ be defined as in Definition \ref{converi}, endowed with the simplicial structure considered in Proposition \ref{hurgove}. Let $\bfA^{\degree}$ be the full subcategory of $\bfA$ spanned by the fibrant-cofibrant objects, and let $\bfA^{\otimes,\degree}$ be the full subcategory
of $\bfA^{\otimes}$ spanned by those objects $[C_1, \ldots, C_n]$ such that each
$C_i$ belongs to $\bfA^{\degree}$. Then the natural map
$p: \Nerve( \bfA^{\otimes,\degree} ) \rightarrow \Nerve( \cDelta)^{op}$
determines a monoidal structure on the $\infty$-category $\Nerve( \bfA^{\degree} )$.
\end{proposition}

\begin{proof}
Our first step is to prove that $\Nerve( \bfA^{\otimes, \degree})$ is an $\infty$-category.
To prove this, it will suffice to show that $\bfA^{\otimes, \degree}$ is a fibrant simplicial category. 
Let $[C_1, \ldots, C_n]$ and $[C'_1, \ldots, C'_{m}]$ be objects of $\bfA^{\otimes, \degree}$. Then
$$\bHom_{ \bfA^{\otimes, \degree} }( [ C_1, \ldots, C_n ], [C'_1, \ldots, C'_m ])$$
is a disjoint union of products of simplicial sets of the form
$\bHom_{\bfA}( C_i \otimes \ldots \otimes C_j, C'_k)$. Each of these simplicial sets
is a Kan complex, since $C_i \otimes \ldots \otimes C_j$ is cofibrant and $C'_k$ is fibrant.

Since $p$ is a map from an $\infty$-category to the nerve of an ordinary category, it is automatically an inner fibration. We next claim that $p$ is a coCartesian fibration. Let
$[ C_1, \ldots, C_n]$ be an object of $\bfA^{\otimes, \degree}$, and let 
$f: [m] \rightarrow [n]$ be a map in $\cDelta$. For each $1 \leq i \leq m$, choose a trivial
cofibration $$\eta_i: C_{ f(i-1)+1 } \otimes \ldots \otimes C_{f(i)} \rightarrow C'_{i},$$ where
$C'_{i}$ is a fibrant object of $\bfA$. Together these determine a map $\overline{f}: [C_1, \ldots, C_n] \rightarrow [C'_1, \ldots, C'_m]$ in $\Nerve( \bfA^{\otimes, \degree})$. We claim that $\overline{f}$ is $p$-coCartesian. In view of Proposition \toposref{trainedg}, it will suffice to show that for every morphism $g: [k] \rightarrow [m]$ in $\cDelta$ and every 
$[ C''_1, \ldots, C''_{k} ] \in \bfA^{\otimes, \degree}$, the induced map
$$\bHom_{ \bfA^{\otimes} }( [C'_1, \ldots, C'_m] , [C''_1, \ldots, C''_k] )
\times_{ \Hom_{\cDelta}( [k], [m] ) } \{g\} 
\rightarrow \bHom_{ \bfA^{\otimes} }( [C_1, \ldots, C_n], [C''_1, \ldots, C''_k])$$
determines a homotopy equivalence onto the summand of
$\bHom_{ \bfA^{\otimes} }( [C_1, \ldots, C_n], [C''_{1}, \ldots, C''_k] )$ spanned
by those morphisms which cover the map $h = f \circ g: [k] \rightarrow [n]$.
Unwinding the definitions, it will suffice to prove that for $1 \leq i \leq k$, the induced map
$$ \bHom_{\bfA}( C'_{ g(i-1)+1} \otimes \ldots \otimes C'_{g(i)}, C''_{i} )
\rightarrow \bHom_{\bfA}( C_{h(i-1)+1} \otimes \ldots \otimes C_{h(i)}, C''_{i} )$$
is a homotopy equivalence of Kan complexes. For this, we need only show that the map 
$$\eta: C_{h(i-1)+1} \otimes \ldots \otimes C_{h(i)} \rightarrow 
C'_{g(i-1) + 1} \otimes \ldots \otimes C'_{g(i)}$$
is a weak equivalence. This follows from the observation that $\eta$ can be identified
with the tensor product of the maps $\{ \eta_{j} \}_{g(i-1) < j \leq g(i) }$, each of which is a weak equivalence between cofibrant objects.

We now conclude the proof by observing that the fiber $\Nerve( \bfA^{\otimes, \degree} )_{[n]}$
is isomorphic to an $n$-fold product of $\Nerve( \bfA^{\degree})$ with itself, and that the projection onto the $j$th factor can be identified with the functor associated to the inclusion
$[1] \simeq \{j-1, j\} \subseteq [n]$.
\end{proof}

\begin{example}\label{exalcu}
Let $\bfA$ be a simplicial model category. Suppose that the Cartesian monoidal structure on $\bfA$
is compatible with the model structure (in other words, that the final object of $\bfA$ is cofibrant, and that for any pair of cofibrations $i: A \rightarrow A'$, $j: B \rightarrow B'$, the induced map
$i \wedge j: (A \times B') \coprod_{A \times B} (A' \times B) \rightarrow A' \times B'$ is a cofibration, trivial if either $i$ or $j$ is trivial). Then there is a canonical map of simplicial categories
$\theta: \bfA^{\otimes} \rightarrow \bfA$, given on objects by the formula
$ \theta( [A_1, \ldots, A_n]) = A_1 \times \ldots \times A_n.$
Since the collection of fibrant-cofibrant objects of $\bfA$ is stable under finite products,
$\theta$ induces a map $\bfA^{\otimes, \degree} \rightarrow \bfA^{\degree}$. It follows that
$\Nerve( \bfA^{\otimes, \degree}) \rightarrow \Nerve(\cDelta)^{op}$ induces the Cartesian monoidal structure on the $\infty$-category $\Nerve( \bfA^{\degree} )$. 
\end{example}

\begin{definition}\label{sumpy}\index{left closed}\index{right closed}\index{closed}\index{monoidal $\infty$-category!closed}
Let $\calC$ be a monoidal $\infty$-category. We will say that $\calC$ is {\it left closed} if, for each $C \in \calC$, the functor $D \mapsto C \otimes D$ admits a right adjoint. Similarly, we will say that $\calC$ is right closed if, for each $C \in \calC$, the functor
$D \mapsto D \otimes C$ admits a right adjoint. We will say that $\calC$ is {\it closed} if it is both left closed and right closed.
\end{definition}

\begin{remark}
In view of Proposition \toposref{sumpytump}, the condition that a monoidal $\infty$-category $\calC$ be closed can be checked at the level of the ($\calH$-enriched) homotopy category of $\calC$, with its induced monoidal structure. More precisely, $\calC$ is {\em right closed} if and only if, for
every pair of objects $C, D \in \calC$, there exists another object $D^{C}$ and a map
$D^{C} \otimes C \rightarrow D$ with the following universal property: for every $E \in \calC$, 
the induced map
$\bHom_{\calC}( E, D^{C} ) \rightarrow \bHom(E \otimes C, D)$
is a homotopy equivalence. In this case, the construction $D \mapsto D^{C}$ determines a right adjoint to the functor $E \mapsto E \otimes C$. 
\end{remark}

\begin{remark}\label{hurto}
Let $\bfA$ be a monoidal model category equipped with a compatible simplicial structure, and let
$C$ be a cofibrant object of $\bfA$. Then the construction $D \mapsto C \otimes D$ determines a left Quillen functor $f_{C}: \bfA \rightarrow \bfA$. Suppose that $C$ is also fibrant, so that $C$ can be identified with an object of the underlying $\infty$-category $\Nerve( \bfA^{\degree})$. Then tensor product with $C$ also induces a functor $Lf_{C}$ from $\Nerve( \bfA^{\degree} )$ to itself. The proof of Proposition \ref{hurgoven} shows that $Lf_{C}$ can be identified with a left derived functor of $f_{C}$. In particular, $Lf_{C}$ admits a right adjoint. It follows that the monoidal $\infty$-category
$\Nerve( \bfA^{\degree} )$ is left closed in the sense of  Definition \ref{sumpy} (the same argument shows that $\Nerve( \bfA^{\degree})$ is also right closed).  
In particular, the functor $Lf_{C}$ preserves all colimits which exist in $\Nerve( \bfA^{\degree} )$. 
\end{remark}

Let $\calC$ be a monoidal category. An {\it algebra object} of $\calC$ is an object $A \in \calC$ equipped with maps
$$ 1 \rightarrow A, \quad  A \otimes A \rightarrow A$$
which satisfy the usual unit and associativity axioms. The collection of algebra objects of $\calC$ can be organized into a category which we will denote by $\Alg(\calC)$.\index{algebra object}\index{ZZZAlgC@$\Alg(\calC)$}

Suppose now that $\calC$ is a monoidal category equipped with a compatible simplicial structure. Then $\Alg(\calC)$ inherits the structure of a simplicial category, where we let $\bHom_{\Alg(\calC)}(A,B) \subseteq \bHom_{\calC}(A,B)$ be the simplicial subset described by the following property: a map $K \rightarrow \bHom_{\calC}(A,B)$ factors through $\bHom_{\Alg(\calC)}(A,B)$ if and only if the diagrams
$$ \xymatrix{ \Delta^0 \times K \ar[d] \ar[r] & \bHom_{\calC}(1,A) \times \bHom_{\calC}(A,B) \ar[d] \\
\Delta^0 \ar[r] & \bHom_{\calC}(1,B) }$$
$$ \xymatrix{ K \ar[r] \ar[dr] & K \times K \ar[r] & \bHom_{\calC}(A,B) \times \bHom_{\calC}(A,B) \ar[r] &
\bHom_{\calC}(A \otimes A, B \otimes B) \ar[d] \\
& \Delta^0 \times K \ar[r] & \bHom_{\calC}(A \otimes A,A) \times \bHom_{\calC}(A,B) \ar[r] & \bHom_{\calC}(A \otimes A,B)}$$
commute. The simplicial nerve of $\Alg(\calC)$ is typically not equivalent to the category of algebra objects of $\Nerve(\calC)$, where $\Nerve(\calC)$ is endowed with the monoidal structure of Proposition \ref{hurgove}. In fact, $\Alg(\calC)$ need not be a fibrant simplicial category, even when $\calC$ is itself fibrant. However, in the case where $\calC$ is a monoidal {\em model} category we can often remedy the situation by passing to a suitable subcategory of $\Alg(\calC)$. 

In the arguments which follow, we will need to invoke the following hypothesis (formulated originally by Schwede and Shipley; see \cite{monmod}):

\begin{definition}[Monoid Axiom]\label{monax}\index{monoid axiom}
Let $\bfA$ be a combinatorial symmetric monoidal model category. Let $U$ be the collection of all morphisms of $\bfA$ having the form
$$ X \otimes Y \stackrel{ \id_X \otimes f}{\rightarrow} X \otimes Y',$$
where $f$ is a trivial cofibration, and let $\overline{U}$ denote the saturated class of morphisms generated by $U$ (Definition \toposref{saturated}). We will say that
{\it $\bfA$ satsifies the monoid axiom} if every morphism of $\overline{U}$ is a weak equivalence in $\bfA$.
\end{definition}

\begin{remark}
Let $\bfA$ be a combinatorial symmetric monoidal model category in which every object is cofibrant, and let $U$ and $\overline{U}$ be as in Definition \ref{monax}. Then every morphism belonging to $U$ is a trivial cofibration. Since the collection of trivial cofibrations in $\bfA$ is saturated, we conclude that $\bfA$ satisfies the monoid axiom.
\end{remark}

\begin{notation}
Let $f: X \rightarrow X'$ and $g: Y \rightarrow Y'$ be morphisms in a monoidal category $\bfA$ which admits pushouts. We define the {\it pushout product} of $i$ and $j$ to be the induced map
$$ f \smash g: (X \otimes Y') \coprod_{ X \otimes Y} (X' \otimes Y) \rightarrow X' \otimes Y'.$$
The operation $\wedge$ endows the category $\Fun( [1], \bfA)$ with a monoidal structure, which is symmetric if the monoidal structure on $\bfA$ is symmetric.
\end{notation}

\begin{lemma}\label{suitus}
Let $\bfA$ be a combinatorial symmetric monoidal model category which satisfies the monoid axiom, and let $\overline{U}$ be as in Definition \ref{monax}. Then:
\begin{itemize}
\item[$(1)$] If $f: X \rightarrow X'$ belongs to $\overline{U}$ and $Y$ is an object
of $\bfA$, then $f \otimes \id_{Y}: X \otimes Y \rightarrow X' \otimes Y$ belongs to $\overline{U}$.

\item[$(2)$] If $f,g \in \overline{U}$, then $f \wedge g \in \overline{U}$.
\end{itemize}
\end{lemma}

\begin{proof}
To prove $(1)$, let $S$ denote the collection of all morphisms $f$ in $\bfA$ such that
$f \otimes \id_{Y}$ belongs to $\overline{U}$. It is easy to see that $S$ is saturated. It will therefore suffice to show that $U \subseteq S$, which is obvious.

To prove $(2)$, we use the same argument. Fix $g$, and let $S'$ be the set of all morphisms
$f \in \bfA$ such that $f \wedge g$ belongs to $\overline{U}$. We wish to prove that
$\overline{U} \subseteq S'$. Since $S'$ is saturated, it will suffice to show that $U \subseteq S'$.
In other words, we may assume that $f$ is of the form $f_0 \otimes \id_{A}$, where
$f_0$ is a trivial cofibration in $\bfA$. Similarly, we may assume that $g = g_0 \otimes \id_{B}$.
Then $f \wedge g = (f_0 \wedge g_0) \otimes (\id_{A \otimes B})$, which belongs to $U$ since
$f_0 \wedge g_0$ is a trivial cofibration in $\bfA$.
\end{proof}

\begin{proposition}\label{comboline}\index{model category!of algebras}[Schwede-Shipley]
Let $\bfA$ be a combinatorial monoidal model category. Assume that either
every object of $\bfA$ is cofibrant, or that $\bfA$ is a symmetric monoidal model category which satisfies the monoid axiom. Then:

\begin{itemize}
\item[$(1)$] The category $\Alg(\bfA)$ admits a combinatorial model structure, where:
\begin{itemize}
\item[$(W)$] A morphism $f: A \rightarrow B$ of algebra objects of $\bfA$ is a weak equivalence if it is a weak equivalence when regarded as a morphism in $\bfA$.
\item[$(F)$] A morphism $f: A \rightarrow B$ of algebra objects of $\bfA$ is a fibration if it is a fibration when regarded as a morphism in $\bfA$.
\end{itemize}
\item[$(2)$] The forgetful functor $\theta: \Alg(\bfA) \rightarrow \bfA$ is a right Quillen functor.
\item[$(3)$] If $\bfA$ is equipped with a compatible simplicial structure, then $\Alg(\bfA)$ inherits the structure of a simplicial model category.
\end{itemize}
\end{proposition}

\begin{proof}
We first observe that the category $\Alg(\bfA)$ is presentable (this follows, for example, from Proposition \ref{algprec}). Recall that a collection $S$ of morphisms in a presentable category $\calC$ is {\it saturated} if it is stable under pushouts, retracts, and transfinite composition (see Definition \toposref{saturated}); we will say that $S$ is {\it generated by} a subset $S_0 \subset S$ if $S$ is the smallest saturated collection of morphisms containing $S_0$.

Since $\calC$ is combinatorial, there exists a (small) collection of morphisms $I = \{ i_{\alpha}: C \rightarrow C' \}$ which generates the class of cofibrations in $\bfA$, and a (small) collection of morphisms $J = \{ j_{\alpha}: D \rightarrow D' \}$ which generates the class of trivial cofibrations in $\bfA$.

Let $F: \bfA \rightarrow \Alg(\bfA)$ be a left adjoint to the forgetful functor. Let $\overline{F(I)}$ be the saturated class of morphisms in $\Alg(\bfA)$ generated by $\{ F(i): i \in I \}$, and let $\overline{F(J)}$ be defined similarly. Unwinding the definitions, we see that a morphism in $\Alg(\bfA)$ is a trivial fibration if and only if it has the right lifting property with respect to $F(i)$, for every $i \in I$. Invoking the small object argument, we deduce that every morphism $f: A \rightarrow C$ in $\Alg(\bfA)$ admits a factorization
$ A \stackrel{f'}{\rightarrow} B \stackrel{f''}{\rightarrow} C$ where $f' \in \overline{F(I)}$ and
$f''$ is a trivial fibration. Similarly, we can find an analogous factorization where
$f' \in \overline{F(J)}$ and $f''$ is a fibration.

Using a standard argument, we may reduce the proof of $(1)$ to the problem of showing that every morphism belonging to $\overline{ F(J) }$ is a weak equivalence in $\Alg(\bfA)$. If
$\bfA$ is symmetric and satisfies the monoid axiom, then let $\overline{U}$ be as in Definition \ref{monax}; otherwise let $\overline{U}$ be the collection of all trivial cofibrations in $\bfA$.
Let $S$ be the collection of all morphisms in $\Alg(\bfA)$ such that the induced map in
$\bfA$ belongs to $\overline{U}$. Then $S$ is saturated, and every element of $S$ is a weak equivalence. To complete the proof, it will sufice to show that $F(J) \subseteq S$. In other words, we must prove:

\begin{itemize}
\item[$(\ast)$] Let 
$$ \xymatrix{ F(C) \ar[r]^{F(i)} \ar[d] & F(C') \ar[d] \\
A \ar[r]^{f} & A' }$$ be a pushout diagram in $\Alg(\bfA)$. If $i$ is a trivial cofibration in
$\bfA$, then $f \in S$.
\end{itemize}

Let $\emptyset$ be an initial object of $\bfA$, and let $j: \emptyset \rightarrow A$ be the unique morphism.
We now observe that $A'$ can be obtained as the direct limit of a sequence
$$ A = A^{(0)} \stackrel{f_1}{\rightarrow} A^{(1)} \stackrel{f_2}{\rightarrow} \ldots$$
of objects of $\bfA$, where each $f_n$ is a pushout of $j \wedge i \wedge j \wedge \ldots \wedge i \wedge j$; here
the factor $i$ appears $n$ times. If every object of $\bfA$ is cofibrant, then we conclude
that $f_{n}$ is a trivial cofibration using the definition of a monoidal model category.
If the monoidal structure on $\bfA$ is symmetric and satisfies the monoid axiom, then
repeated application of Lemma \ref{suitus} shows that $f_{n} \in \overline{U}$. 
Since $\overline{U}$ is saturated, it follows that $f \in S$ as desired. This completes the proof of $(1)$.

Assertion $(2)$ is obvious. To prove $(3)$, we observe both $\bfA$ and $\Alg(\bfA)$ are cotensored  over simplicial sets, and that we have canonical isomorphisms $\theta( A^K ) \simeq \theta(A)^K$
for $A \in \Alg(\bfA)$, $K \in \sSet$. To prove that $\Alg(\bfA)$ is a simplicial model category, it will suffice to show that $\Alg(\bfA)$ is tensored over simplicial sets, and that given a fibration
$i: A \rightarrow A'$ in $\Alg(\bfA)$ and a cofibration $j: K \rightarrow K'$ in $\sSet$, the induced map $A^{K'} \rightarrow A^{K} \times_{ {A'}^{K}} {A'}^{K'}$ is a fibration, trivial if either $i$ or
$j$ is a fibration. The second claim follows from the fact that $\theta$ detects fibrations and trivial fibrations. For the first, it suffices to prove that for $K \in \sSet$, the functor $A \mapsto A^{K}$
has a left adjoint; this follows from the adjoint functor theorem.
\end{proof}

We now come to the main result of this section, which asserts that if $\bfA$ is a simplicial model category with a compatible monoidal structure, then every algebra object of the $\infty$-category $\Nerve( \bfA^{\degree} )$ is equivalent to a strictly associative algebra object in $\bfA$. This is a kind of ``straightening theorem'' which we will use in \S \ref{monoid6.0} to compare our approach to the theory of $A_{\infty}$-rings with more classical definitions.

\begin{lemma}\label{GGG}
Let $\bfA$ be a combinatorial monoidal model category with a compatible simplicial structure, and let $\calC$ be a small category such that $\Nerve(\calC)$ is sifted (Definition \stableref{siftdef}). 
Assume either that every object of $\bfA$ is cofibrant, or that $\bfA$ satisfies the following
pair of conditions:
\begin{itemize}
\item[$(A)$] The monoidal structure on $\bfA$ is symmetric, and $\bfA$ satisfies the monoid axiom.
\item[$(B)$] The class of cofibrations in $\bfA$ is generated by cofibrations between cofibrant objects (this is automatic if every object of $\bfA$ is cofibrant).
\end{itemize}
Then the forgetful functor $\Nerve( \Alg( \bfA)^{\degree} ) \rightarrow \Nerve( \bfA^{\degree} )$
preserves $\Nerve(\calC)$-indexed colimits.
\end{lemma}

\begin{proof}
In view of Theorem \toposref{colimcompare} and Proposition \toposref{gumby4}, it will suffice to prove that the forgetful functor $\theta: \Alg(\bfA) \rightarrow \bfA$ preserves homotopy colimits indexed by $\calC$. Let us regard $\Alg(\bfA)^{\calC}$ and $\bfA^{\calC}$ as endowed with
the projective model structure (see \S \toposref{quasilimit3}). Let $F: \bfA^{\calC} \rightarrow \bfA$
and $F_{\Alg}: \Alg(\bfA)^{\calC} \rightarrow \Alg(\bfA)$ be colimit functors, and let
$\theta^{\calC}: \Alg(\bfA)^{\calC} \rightarrow \bfA^{\calC}$ be given by composition with
$\theta$. Since $\Nerve(\calC)$ is sifted, there is a canonical isomorphism of functors
$\alpha: F \circ \theta^{\calC} \simeq \theta \circ F_{\Alg}$. We wish to prove that this isomorphism persists after deriving all of the relevant functors. Since $\theta$ and $\theta^{\calC}$ preserve weak equivalences, they can be identified with their right derived functors. Let $LF$ and $LF_{\Alg}$
be the left derived functors of $F$ and $F_{\Alg}$, respectively. Then $\alpha$ induces a natural transformation $\overline{\alpha}: LF \circ \theta^{\calC} \rightarrow \theta \circ LF_{\Alg}$; we wish to show that $\overline{\alpha}$ is an isomorphism. Let $A: \calC \rightarrow \Alg(\bfA)$ be a strongly cofibrant object of $\Alg(\bfA)^{\calC}$; we must show that the natural map
$$LF( \theta^{\calC}(A)) \rightarrow \theta( LF_{\Alg}(A)) \simeq  \theta( F_{\Alg}(A))
\simeq F( \theta^{\calC}(A))$$ is a weak equivalence in $\bfA$. 

Let us say that an object $X \in \bfA^{\calC}$ is {\it good} if each of the objects
$X(C) \in \bfA$ is cofibrant, the object $F(X) \in \bfA$ is cofibrant, and the natural map $LF(X) \rightarrow F(X)$ is a weak equivalence in $\bfA$: in other words, if the colimit of $X$ is also a homotopy colimit of $X$.
To complete the proof, it will suffice to show that $\theta^{\calC}(A)$ is good, whenever $A$
is a strongly cofibrant object of $\Alg(\bfA)^{\calC}$. This is not obvious, since $\theta^{\calC}$ is a right Quillen functor and does not preserve strongly cofibrant objects in general (note that we have not yet used the full strength of our assumption that $\Nerve(\calC)$ is sifted). To continue the proof, we will need a relative version of the preceding condition. We will say that a morphism
$f: X \rightarrow Y$ in $\bfA^{\calC}$ is {\it good} if the following conditions are satisfied:
\begin{itemize}
\item[$(i)$] The objects $X,Y \in \bfA^{\calC}$ are good.
\item[$(ii)$] For each $C \in \calC$, the induced map $X(C) \rightarrow Y(C)$ is a cofibration in $\bfA$.
\item[$(iii)$] The map $F(X) \rightarrow F(Y)$ is a cofibration in $\bfA$.
\end{itemize}

We now make the following observations:
\begin{itemize}
\item[$(1)$] The collection of good morphisms is stable under transfinite composition.
More precisely, suppose given an ordinal $\alpha$ and a direct system of objects
$\{ X^{\beta} \}_{\beta < \alpha}$ of $\bfA^{\calC}$. Suppose further that for every $0 < \beta < \alpha$, the map $\varinjlim \{ X^{\gamma} \}_{ \gamma < \beta } \rightarrow X^{\beta}$ is
good. Then the induced map $X^0 \rightarrow \varinjlim \{ X^{\beta} \}_{\beta < \alpha} $ is good.
The only nontrivial point is to verify that the object $X = \varinjlim \{ X^{\beta} \}_{\beta < \alpha}$ is good. For this, we observe $X$ is a homotopy colimit of the system $\{ X^{\beta} \}$ (in virtue of $(ii)$), while
$F(X)$ is a homotopy colimit of the system $\{ F(X^{\beta}) \}$ (in virtue of $(iii)$), and that the collection of homotopy colimit diagrams is stable under homotopy colimits.

\item[$(2)$] Suppose given a pushout diagram
$$ \xymatrix{ X \ar[r]^{f} \ar[d] & Y \ar[d] \\
X' \ar[r]^{f'} & Y' }$$
in $\bfA^{\calC}$. If $f$ is good and $X'$ is good, then $f'$ is good. Once again, the only nontrivial point is to show that $Y'$ is good. To see this, we observe that our hypotheses imply that $Y'$ is homotopy pushout of $Y$ with $X'$ over $X$. Similarly, $F(Y')$ is a homotopy pushout of $F(Y)$ with $F(X')$ over $F(X)$. We now invoke once again the fact that the class of homotopy colimit diagrams is stable under homotopy colimits.

\item[$(3)$] Let $F: \calC \rightarrow \bfA$ be a constant functor whose value is a cofibrant object of $\bfA$. Then $F$ is good. This follows from the fact that $\Nerve(\calC)$ is weakly contractible (Proposition \stableref{siftcont}). 

\item[$(4)$] Every strongly cofibrant object of $\bfA^{\calC}$ is good. Every strong
cofibration between strongly cofibrant objects of $\bfA^{\calC}$ is good.

\item[$(5)$] If $X$ and $Y$ are good objects of $\bfA^{\calC}$, then $X \otimes Y$ is good. To prove this, we first observe that the collection of cofibrant objects of $\bfA$ is stable under tensor products. Because $\Nerve(\calC)$ is sifted, we have a chain of isomorphisms
in $\h{\bfA}$: $$LF(X \otimes Y) \simeq LF(X) \otimes LF(Y) \simeq F(X) \otimes F(Y) \simeq F(X \otimes Y).$$

\item[$(6)$] Let $f: X \rightarrow X'$ be a good morphism in $\bfA^{\calC}$, and let
$Y$ be a good object of $\bfA^{\calC}$. Then the morphism
$f \otimes \id_{Y}$ is good. Condition $(i)$ follows from $(5)$, condition $(ii)$ follows from the fact that tensoring with each $Y(C)$ preserves cofibrations (since $Y(C)$ is cofibrant), and condition $(iii)$ follows by applying the same argument to $F(Y)$ (and invoking the fact that $F$ commutes with tensor products). 

\item[$(7)$] Let $f: X \rightarrow X'$ and $g: Y \rightarrow Y'$ be good morphisms in $\bfA^{\calC}$. Then $$f \wedge g: (X \otimes Y') \coprod_{ X \otimes Y} (X' \otimes Y) \rightarrow X' \otimes Y'$$
is good. Condition $(ii)$ follows immediately from the fact that $\bfA$ is a monoidal model category.
Condition $(iii)$ follows from the same argument, together with the observation that $F$ commutes with pushouts and tensor products. Condition $(i)$ follows by combining $(5)$, $(6)$, and $(2)$.
\end{itemize}

We observe that our assumption $(B)$ (which is a consequence of the assumption that every object of $\bfA$ is cofibrant) implies an analogous result for $\bfA^{\calC}$:
\begin{itemize}
\item[$(B')$] The collection of all strong cofibrations in $\bfA^{\calC}$ is generated by
strong cofibrations between strongly cofibrant objects.
\end{itemize}

Let $T: \bfA^{\calC} \rightarrow \Alg(\bfA)^{\calC}$ be a left adjoint to $\theta^{\calC}$. 
Using the small object argument and $(B')$, we conclude that for every strongly cofibrant object $A \in \Alg(\bfA)^{\calC}$ there exists a transfinite sequence $\{ A^{\beta} \}_{\beta \leq \alpha}$ in $\Alg(\bfA)^{\calC}$ with the following properties:
\begin{itemize}
\item[$(a)$] The object $A^{0}$ is initial in $\Alg(\bfA)^{\calC}$.
\item[$(b)$] The object $A$ is a retract of $A^{\alpha}$.
\item[$(c)$] If $\lambda \leq \alpha$ is a limit ordinal, then
$A^{\lambda} \simeq \colim \{ A^{\beta} \}_{ \beta < \lambda}$.  
\item[$(d)$] For each $\beta < \alpha$, there is a pushout diagram
$$ \xymatrix{ T(X') \ar[r]^{T(f)} \ar[d] & T(X) \ar[d] \\
A^{\beta} \ar[r] & A^{\beta+1} }$$
where $f$ is a strong cofibration between strongly cofibrant objects of $\bfA^{\calC}$. 
\end{itemize}

We wish to prove that $\theta^{\calC}(A)$ is good. In view of $(b)$, it will suffice to show
that $\theta^{\calC}(A^{\alpha})$ is good. We will prove a more general assertion: for
every $\gamma \leq \beta \leq \alpha$, the induced morphism $u_{\gamma, \beta}: 
\theta^{\calC}(A^{\gamma}) \rightarrow \theta^{\calC}( A^{\beta})$ is good. The proof is by induction on $\beta$. If $\beta = 0$, then we are reduced to proving that $\theta^{\calC}(A^0)$ is good. This follows from $(a)$ and $(3)$. If $\beta$ is a nonzero limit ordinal, then the desired result follows from $(c)$ and $(1)$. It therefore suffices to treat the case where $\beta = \beta' + 1$ is a successor ordinal. Moreover, we may suppose that
$\gamma = \beta'$: if $\gamma < \beta'$, then we observe that
$u_{\gamma, \beta} = u_{\beta', \beta} \circ u_{ \gamma, \beta'}$ and invoke $(1)$, while
if $\gamma > \beta'$, then $\gamma = \beta$ and we are reduced to proving that
$\theta^{\calC}(A^{\beta})$ is good, which follows from the assertion that
$u_{\beta', \beta}$ is good. We are now reduced to proving the following:

\begin{itemize}
\item[$(\ast)$] Let $$ \xymatrix{ T(X') \ar[r]^{T(f)} \ar[d] & T(X) \ar[d] \\
B' \ar[r]^{v} & B }$$
be a pushout diagram in $\Alg(\bfA)^{\calC}$, where $f: X' \rightarrow X$ is a strong cofibration between strongly cofibrant objects of $\bfA^{\calC}$. If $\theta^{\calC}(B')$ is good, then $\theta^{\calC}(v)$ is good.
\end{itemize}

To prove $(\ast)$, we set $Y = \theta^{\calC}(B) \in \bfA^{\calC}$, $Y' = \theta^{\calC}(B') \in \bfA^{\calC}$. Let $g: \emptyset \rightarrow Y'$ the unique morphism, where $\emptyset$
denotes an initial object of $\bfA^{\calC}$. As in the proof of Proposition \ref{comboline}, $Y$
can be identified with the colimit of a sequence
$$ Y^{(0)} \stackrel{w_1}{\rightarrow} Y^{(1)} \stackrel{w_2}{\rightarrow} \ldots $$
where $Y^{(0)} = Y'$, and $w_{k}$ is a pushout of the morphism
$f^{(k)} = g \wedge f \wedge g \wedge \ldots \wedge f \wedge g$, where the factor $f$ appears $k$ times. In view of $(1)$ and $(2)$, it will suffice to prove that each $f^{(k)}$ is a good morphism.
Since $Y'$ is good, we conclude immediately that $g$ is good. It follows from
$(4)$ that $f$ is good. Repeated application of $(7)$ allows us to deduce that $f^{(k)}$ is good, and to conclude the proof.
\end{proof}

We now come to our main result:

\begin{theorem}\label{beckify}
Let $\bfA$ be a combinatorial monoidal model category equipped with a compatible simplicial structure. Assume either:
\begin{itemize}
\item[$(A)$] Every object of $\bfA$ is cofibrant.
\item[$(B)$] The class of cofibrations in $\bfA$ is generated by cofibrations between cofibrant objects, the monoidal structure on $\bfA$ is symmetric, and $\bfA$ satisfies the monoid axiom.
\end{itemize}
Then the canonical map
$$ \Nerve(\Alg(\bfA)^{\degree}) \rightarrow \Alg( \Nerve(\bfA^{\degree}) )$$
is an equivalence of $\infty$-categories.
\end{theorem}

\begin{example}\label{fielder}
Let $k$ be a field, and let $\bfA$ denote the category of complexes of $k$-vector spaces
$$ \ldots \rightarrow M_{n+1} \rightarrow M_{n} \rightarrow M_{n-1} \rightarrow \ldots, $$
with monoidal structure given by the formation of tensor products of complexes.
The category $\bfA$ admits a model structure, where:
\begin{itemize}
\item[$(C)$] A map of complexes $f: M_{\bigdot} \rightarrow N_{\bigdot}$ is a {\it cofibration} if it
induces an injection $M_{n} \rightarrow N_{n}$ of $k$-vector spaces, for each
$n \in \Z$.
\item[$(F)$] A map of complexes $f: M_{\bigdot} \rightarrow N_{\bigdot}$ is a {\it fibration} if it
induces an surjection $M_{n} \rightarrow N_{n}$ of $k$-vector spaces, for each
$n \in \Z$.
\item[$(W)$] A map of complexes $f: M_{\bigdot} \rightarrow N_{\bigdot}$ is a {\it weak equivalence} if it is a quasi-isomorphism; that is, if $f$ induces an isomorphism on homology groups
$H_{n}( M_{\bigdot}) \rightarrow H_{n}( N_{\bigdot})$ for each $n \in \Z$.
\end{itemize}
The category $\bfA$ also admits a simplicial structure, compatible with its monoidal model structure; see the discussion in \S \stableref{stable10}. Moreover, every object of $\bfA$ is cofibrant.

The $\infty$-category $\Nerve( \bfA^{\degree})$ can be identified with the (unbounded) {\it derived $\infty$-category of $k$-vector spaces} (see Definition \stableref{smucky}); let us denote this $\infty$-category by $\calD(k)$. According to Theorem \ref{beckify}, the $\infty$-category of 
algebra object $\Alg( \calD(k) )$ is equivalent to $\Nerve( \Alg( \bfA)^{\degree} )$: that is, to a suitable $\infty$-category of {\it differential graded algebras} over $k$.

An analogous result holds if we replace the field $k$ by an arbitrary commutative ring, and endow
$\bfA$ with the {\em projective} model structure. In this case, not every object of $\bfA$ is cofibrant, but $\bfA$ still satisfies hypothesis $(B)$ of Theorem \ref{beckify} (see
\cite{monmod}).
\end{example}

\begin{example}\label{bulwork}\index{symmetric spectra}
Let $\bfA$ be the category of {\it symmetric spectra}, as defined in \cite{symmetricspectra}. Then
$\bfA$ admits several model structures which satisfy assumption $(B)$ of Theorem \ref{beckify}.
Using Corollary \ref{surcoi}, we deduce that $\Nerve(\bfA^{\degree})$ is equivalent, as a monoidal $\infty$-category, to the $\infty$-category $\Spectra$ of spectra (endowed with the smash product monoidal structure; see \S \ref{hummingburg}). Using Theorem \ref{beckify}, we deduce that our $\infty$-category $\AInfty = \Alg( \Spectra)$ of $A_{\infty}$-rings is equivalent to the $\infty$-category underlying the category $\Alg(\bfA)$ of {\it symmetric ring spectra}.
\end{example}

\begin{example}\label{exalcun}[$\infty$-Categorical MacLane Coherence Theorem]\index{MacLane coherence theorem}
Let $\bfA$ be the category of marked simplicial sets (see \S \toposref{twuf}). Then $\bfA$ is a simplicial model category, which satisfies the hypotheses of Example \ref{exalcu}. The underlying $\infty$-category $\Nerve(\bfA^{\degree})$ can be identified with $\Cat_{\infty}$, the $\infty$-category of $\infty$-categories. Proposition \ref{ungbat} implies that composition with the
Cartesian structure $\Nerve( \bfA^{\otimes, \degree}) \rightarrow \Nerve( \bfA^{\degree} )$ induces an equivalence of $\infty$-categories $\Alg( \Nerve(\bfA^{\degree}) ) \rightarrow \Mon( \Cat_{\infty})$. Combining this observation with Theorem \ref{beckify}, we conclude that the $\infty$-category of monoid objects of $\Cat_{\infty}$ is equivalent to the $\infty$-category underlying the category of {\em strictly} associative monoids in $\bfA$. In other words, every monoidal $\infty$-category $\calC$ is equivalent (as a monoidal $\infty$-category) to an $\infty$-category
$\calC'$ equipped with a strictly associative monoid structure $\calC' \times \calC' \rightarrow \calC'$ (which determines a monoidal structure on $\calC'$ via Proposition \ref{sulken}). We regard this assertion as an $\infty$-categorical analogue of {\it MacLane's coherence theorem}, which asserts that every monoidal category is equivalent to a strict monoidal category (that is, a monoidal category in which the tensor product operation $\otimes$ is associative up to equality, and the
associativity isomorphisms are simply the identity maps).
\end{example}

\begin{proof}[Proof of Theorem \ref{beckify}]
Consider the diagram
$$ \xymatrix{ \Nerve( \Alg(\bfA)^{\degree}) \ar[rr] \ar[dr]^{G} & & \Alg( \Nerve( \bfA^{\degree}) )\ar[dl]^{G'} \\
& \Nerve(\bfA^{\degree}). & }$$
It will suffice to show that this diagram satisfies the hypotheses of Corollary \ref{littlerbeck}:
\begin{itemize}
\item[$(a)$] The $\infty$-categories $\Nerve( \Alg(\bfA)^{\degree} )$ and
$\Alg( \Nerve(\bfA^{\degree}))$ admit geometric realizations of simplicial objects. In fact, both of these $\infty$-categories are presentable. For $\Nerve( \Alg(\bfA)^{\degree})$, this follows from
Propositions \toposref{notthereyet} and \ref{comboline}. For $\Alg( \Nerve(\bfA)^{\degree})$, we first observe that $\Nerve(\bfA)^{\degree}$ is presentable (Proposition \toposref{notthereyet}) and that
the tensor product preserves colimits separately in each variable, and apply Proposition \ref{algprec}.

\item[$(b)$] The functors $G$ and $G'$ admit left adjoints $F$ and $F'$. The existence of a left adjoint to $G$ follows from the fact that $G$ is determined by a right Quillen functor. The existence of a left adjoint to $G'$ follows from Proposition \ref{hutmunn}.

\item[$(c)$] The functor $G'$ is conservative and preserves geometric realizations of simplicial objects. This follows from Corollaries \ref{filtfem} and \ref{jumunj}.

\item[$(d)$] The functor $G$ is conservative and preserves geometric realizations of simplicial objects. The first assertion is immediate from the definition of the weak equivalences in
$\Alg(\bfA)$, and the second follows from Lemma \ref{GGG}.

\item[$(e)$] The canonical map $G' \circ F' \rightarrow G \circ F$ is an equivalence of functors. This follows from the observation that both sides induce, on the level of homotopy categories, the free algebra functor $C \mapsto \coprod_{ n \geq 0} C^{\otimes n}$ (Proposition \ref{hutmunn}).
\end{itemize}
\end{proof}

\begin{remark}
Proposition \ref{hurgoven} admits a converse. Suppose that $\calC$ is a presentable $\infty$-category endowed with a monoidal structure, and that the
associated bifunctor $\otimes: \calC \times \calC \rightarrow \calC$ preserves small colimits separately in each variable. Then $\calC$ is equivalent (as a monoidal $\infty$-category) to $\Nerve(\bfA^{\degree})$, where
$\bfA$ is a combinatorial simplicial model category, endowed with a compatible monoidal structure. Since we will not need this fact, we will only give a sketch of proof.

First, we apply Example \ref{exalcun} to reduce to the case where $\calC$ is a {\em strict} monoidal $\infty$-category; that is, $\calC$ is a simplicial monoid. Now choose a regular 
cardinal $\kappa$ such that $\calC$ is $\kappa$-accessible. Enlarging $\kappa$ if necessary, we may suppose that the full subcateogry $\calC^{\kappa} \subseteq \calC$ spanned by the $\kappa$-compact objects contains the unit object of $\calC$ and is stable under tensor products.

The $\infty$-category $\calC^{\kappa}$ is essentially small. We define a sequence of simplicial subsets $$ \calD(0) \subseteq \calD(1) \subseteq \ldots \subseteq \calC^{\kappa}$$ as follows. Let 
$\calD(0) = \emptyset$, and for $i \geq 0$ let $\calD(i+1)$ be a small simplicial subset of
$\calC^{\kappa}$ which is categorically equivalent to $\calC^{\kappa}$ and contains the submonoid of $\calC^{\kappa}$ generated by $\calD(i)$. Let $\calD = \bigcup \calD(i)$, so that $\calD$ is a small simplicial submonoid of $\calC^{\kappa}$ such that the inclusion $\calD \subseteq \calC^{\kappa}$ is a categorical equivalence.

The proof of Theorem \toposref{pretop} shows that $\calC$ can be identified with an accessible
localization of $\calP(\calD)$. According to Proposition \toposref{othermod}, we can identify
$\calP(\calD)$ with $\Nerve( \bfA^{\degree})$, where $\bfA$ denotes the category $(\sSet)_{/\calD}$ endowed with the contravariant model structure (see \S \toposref{contrasec}). Let $L: \calP(\calD) \rightarrow \calC$ be a localization functor, and let $\bfB$ be the category $(\sSet)_{/\calD}$ endowed with the following {\em localized} model structure:
\begin{itemize}
\item[$(C)$] A morphism $\alpha: X \rightarrow Y$ in $(\sSet)_{/\calD}$ is a cofibration in
$\bfB$ if and only if $\alpha$ is a monomorphism of simplicial sets.
\item[$(W)$] A morphism $\alpha: X \rightarrow Y$ in $(\sSet)_{/\calD}$ is a weak equivalence in $\bfB$ if and only if the $L(\beta)$ is an isomorphism in the homotopy category $\h{\calC}$, where $\beta$ denotes the corresponding morphism in $\h{ \calP(\calD)} \simeq \h{\bfA}$. 
\item[$(F)$] A morphism $\alpha: X \rightarrow Y$ in $(\sSet)_{/\calD}$ is a fibration in $\bfB$ if and only if it has the right lifting property with respect to every morphism which is simultaneously a cofibration and a weak equivalence in $\bfB$.
\end{itemize}
Proposition \toposref{suritu} implies that $\bfB$ is a (combinatorial) simplicial model category, and that the underlying $\infty$-category $\Nerve(\bfB^{\degree})$ is equivalent to $\calC$.

The category $(\sSet)_{/\calD}$ is endowed with a monoidal structure,
which may be described as follows: given a finite collection of objects $X_1, \ldots, X_n \in \SSet_{/\calD}$, we let
$X_1 \otimes \ldots \otimes X_n$ denote the product $X_1 \times \ldots \times X_n$ of the underlying simplicial sets, mapping to $\calD$ via the composition
$X_1 \times \ldots \times X_n \rightarrow \calD^{n} \rightarrow \calD$,
where the second map is given by the monoid structure on $\calD$. It is not difficult to verify that this monoidal structure is compatible with the model structure on $\bfB$. Applying Proposition \ref{hurgoven}, we deduce that $\Nerve( \bfB^{\degree, \otimes} )$ determines a monoidal structure on
$\Nerve(\bfB^{\degree}) \simeq \calC$. One can show that this monoidal structure coincides (up to equivalence) with the structure determined by the associative multiplication on $\calC$.
\end{remark}

\subsection{Digression: Segal Monoidal $\infty$-Categories}\label{segalapp}

In this section, we will introduce the definition of a {\it Segal monoidal $\infty$-category}. The theory of a Segal monoidal $\infty$-categories is equivalent to the theory of monoidal $\infty$-categories which we use throughout this paper (see Remark \ref{segmonex}). However, the formalism of Segal monoidal $\infty$-categories is better suited to the technical arguments we use in \S \ref{monoid5} to construct free algebras. The results of this section will be used {\em only} in \S \ref{monoid5}, and may be skipped without loss of continuity.

\begin{notation}\index{ZZZJast@$J_{\ast}$}\index{ZZZnostar@$\nostar{n}$}
\index{ZZZnostarr@$\seg{n}$}\label{labusit}
If $J$ is a linearly ordered set, we let $J_{\ast}$ denote the set $J \cup \{\ast \}$, where
$\ast$ is a new distinct element. We do not regard this new element as in any way related to the ordering on $J$. For every nonnegative integer $n$, we let $\nostar{n}$ denote the linearly ordered set
$\{1 < 2 < \ldots < n \}$.
\end{notation}

\begin{definition}\index{ZZZLinSeg@$\LinSeg$}
The category $\LinSeg$ is defined as follows.
\begin{itemize}
\item[$(1)$] The objects of $\LinSeg$ are the sets $\seg{n}$, where $n$ is a nonnegative integer.
\item[$(2)$] Given a pair of objects $\seg{m}, \seg{n} \in \LinSeg$, we let
$\Hom_{\LinSeg}( \seg{m}, \seg{n})$ denote the collection of all maps
$\alpha: \seg{m} \rightarrow \seg{n}$ with the following properties:
\begin{itemize}
\item[$(i)$] The map $\alpha$ satisfies $\alpha(\ast) = \ast$.
\item[$(ii)$] The restriction $\alpha | ( \seg{m} - \alpha^{-1} \{ \ast \} )$ is a (nonstrictly) order-preserving map.
\end{itemize}
\end{itemize}
\end{definition}

\begin{remark}
We note that every finite linearly ordered set $J$ is uniquely isomorphic to
$\nostar{n}$, where $n$ is the cardinality of $J$. We will typically abuse notation
by referring to $J_{\ast}$ as an object of $\LinSeg$; in this case, we are implicitly identifying
$J_{\ast}$ with $\seg{n}$.
\end{remark}

\begin{notation}\label{defps}
We define a functor $\psi: \cDelta^{op} \rightarrow \LinSeg$ as follows:
\begin{itemize}
\item[$(1)$] For each $n \geq 0$, we have $\psi( [n] ) = \seg{n}$.
\item[$(2)$] Given a morphism $\alpha: [n] \rightarrow [m]$ in $\cDelta$, the associated morphism
$\psi(\alpha): \seg{m} \rightarrow \seg{n}$ is given by the formula
$$\psi(\alpha)(i) = \begin{cases} j & \text{if } (\exists j) [\alpha(j-1) < i \leq \alpha(j)] \\
\ast & \text{otherwise.} \end{cases}$$
\end{itemize}
More informally, $\psi$ assigns to a nonempty linearly ordered set $[n]$ the set of
all ``gaps'' between adjacent elements of $[n]$.
\end{notation}

\begin{notation}\label{alps}\index{ZZZalphajJ@$\alpha^{j,J}$}\label{ajack}
If $J$ is a finite linearly ordered set containing an element $j$, we let
$\alpha^{j,J}: J_{\ast} \rightarrow \{j\}_{\ast}$
denote the map defined by the formula
$$\alpha^{j,J}(i) = \begin{cases} j & \text{if } i= j \\
\ast & \text{otherwise}. \end{cases}$$
\end{notation}

\begin{definition}\index{monoid object!Segal}\index{Segal!monoid object}
Let $\calC$ be an $\infty$-category. A {\it Segal monoid object} of $\calC$ is a functor
$f: \Nerve( \LinSeg ) \rightarrow \calC$ with the following property:
\begin{itemize}
\item[$(\ast)$] For every finite linearly ordered set $J$, the collection of maps
$\{ f( \alpha^{j,J} ) \}$ exhibits $f( J_{\ast})$ as a product of the objects
$\{ f( \{j\}_{\ast} ) \}_{j \in J}$. 
\end{itemize}
We let $\SegMon( \calC)$ denote the full subcategory of
$\Fun(\Nerve(\LinSeg), \calC)$ spanned by the Segal monoid objects.\index{ZZZSegMon@$\SegMon(\calC)$}
\end{definition}

\begin{proposition}\label{scabby}
Let $\calC$ be an $\infty$-category. Composition with the functor
$\psi: \cDelta^{op} \rightarrow \LinSeg$ of Notation \ref{defps} induces an
equivalence of $\infty$-categories $q: \SegMon(\calC) \rightarrow \Mon(\calC)$.
\end{proposition}

We will give the proof of Proposition \ref{scabby} at the end of this section.

\begin{definition}\index{monoidal $\infty$-category!Segal}\index{Segal!monoidal $\infty$-category}
A {\it Segal monoidal category} is a coCartesian fibration
$\calC^{\otimes} \rightarrow \Nerve( \LinSeg)$ with the following property:
for every finite linearly ordered set $J$, the associated functors $\alpha^{j,J}_{!}: \calC^{\otimes}_{J_{\ast}} \rightarrow \calC^{\otimes}_{ \{j\}_{\ast} }$ induce an equivalence
$\calC^{\otimes}_{J_{\ast}} \simeq \prod_{j \in J} \calC^{\otimes}_{ \{j\}_{\ast} }.$
\end{definition}

\begin{remark}
Let $q: \calC^{\otimes} \rightarrow \Nerve(\LinSeg)$ be a Segal monoidal $\infty$-category.
We will write $\calC$ to indicate the fiber $\calC^{\otimes} \times_{ \Nerve(\LinSeg)} \{ \seg{1} \}$. We will refer to $\calC$ as the {\it underlying $\infty$-category of $\calC^{\otimes}$}, and say that
$q: \calC^{\otimes} \rightarrow \Nerve(\LinSeg)$ is a {\it Segal monoidal structure} on $\calC$.
\end{remark}

\begin{remark}\label{utilbeltil}
Let $q: \calC^{\otimes} \rightarrow \Nerve(\LinSeg)$ be a Segal monoidal $\infty$-category, and let
$\psi: \cDelta^{op} \rightarrow \LinSeg$ be defined as in Notation \ref{defps}. The fiber product
$\calC^{\otimes} \times_{ \Nerve( \LinSeg)} \Nerve(\cDelta)^{op}$ is a monoidal $\infty$-category in the sense of Definition \ref{mainef}.
\end{remark}

\begin{remark}\label{segmonex}
Let $q: \calC^{\otimes} \rightarrow \Nerve(\LinSeg)$ be a coCartesian fibration, classified by a functor $f: \Nerve(\LinSeg) \rightarrow \Cat_{\infty}$. Then $\calC^{\otimes}$ is a Segal monoidal $\infty$-category if and only if $f$ is a Segal monoid object of $\Cat_{\infty}$. Proposition \ref{scabby} implies that the pullback functor $\SegMon( \Cat_{\infty} ) \rightarrow \Mon(\Cat_{\infty})$ is an equivalence of $\infty$-categories. Consequently, up to equivalence, every monoidal $\infty$-category arises from a Segal monoidal $\infty$-category via the construction of Remark \ref{utilbeltil}.
\end{remark}

Our final objective in this section is to study the analogue of the theory of algebras in the setting of Segal monoidal $\infty$-categories.

\begin{definition}\index{Segal!algebra object}\index{algebra object!Segal}\index{ZZZSegAlgC@$\SegAlg(\calC)$}
Let $q: \calC^{\otimes} \rightarrow \Nerve( \LinSeg)$ be a Segal monoidal structure on an $\infty$-category $\calC = \calC^{\otimes}_{ \seg{1} }$. A {\it Segal algebra object} of $\calC$ is a
section $A: \Nerve(\LinSeg) \rightarrow \calC^{\otimes}$ of $q$ with the following property:
for every morphism $\alpha: J_{\ast} \rightarrow J'_{\ast}$ in $\LinSeg$ which induces a bijection
$\alpha^{-1}( J') \rightarrow J'$, the morphism $A(\alpha)$ is $q$-coCartesian. We let $\SegAlg(\calC)$ denote the full subcategory of $\bHom_{ \Nerve(\LinSeg)}( \Nerve(\LinSeg), \calC^{\otimes})$ spanned by the Segal algebra objects.
\end{definition}

\begin{remark}\label{splurk}
Let $q: \calC^{\otimes} \rightarrow \Nerve( \LinSeg)$ be a Segal monoidal $\infty$-category and
let $A: \Nerve(\LinSeg) \rightarrow \calC^{\otimes}$ be a section of $q$. Then $A$ is a
Segal algebra object if and only if, for every finite linearly ordered set $J$ and every element
$j \in J$, the morphism $A( \alpha^{j,J})$ is $q$-coCartesian (see Notation \ref{alps}).
\end{remark}

\begin{remark}
Let $q: \calC^{\otimes} \rightarrow \Nerve( \LinSeg)$ be a Segal monoidal structure on an $\infty$-category $\calC = \calC^{\otimes}_{\seg{1}}$ and
let $A: \Nerve(\LinSeg) \rightarrow \calC^{\otimes}$ be a section of $q$. Using Remark \ref{splurk}, we deduce that $A$ is a Segal algebra object of $\calC$ if and only if
$A \circ \psi: \Nerve( \cDelta)^{op} \rightarrow \Nerve(\cDelta)^{op} \times_{ \Nerve(\LinSeg)} \calC^{\otimes}$
is an algebra object of $\calC$, in the sense of Definition \ref{suskin}.
\end{remark}

The following comparison result will be needed in \S \ref{monoid5}:

\begin{proposition}\label{algcompare}
Let $q: \calC^{\otimes} \rightarrow \Nerve(\LinSeg)$ be a Segal monoidal structure on an $\infty$-category $\calC$. Composition with the functor $\psi: \cDelta^{op} \rightarrow \LinSeg$ of Notation \ref{defps} determines an equivalence of $\infty$-categories
$\theta: \SegAlg(\calC) \rightarrow \Alg(\calC)$.
\end{proposition}

Proposition \ref{scabby} should be regarded as a special case of Proposition \ref{algcompare}, where we take $\calC$ to be the $\infty$-category $\Cat_{\infty}$, equipped with the Cartesian Segal monoidal structure. Since we do not wish to press any further into the theory, we will not make this precise; instead, we will give two separate (but nearly identical) proofs for
Propositions \ref{algcompare} and \ref{scabby}.

\begin{proof}[Proof of Proposition \ref{scabby}]
It follows immediately from the definition that composition with the functor $\psi$ carries $\SegMon(\calC)$ into
$\Mon(\calC)$. We define a category $\calI$ as follows:
\begin{itemize}
\item[$(1)$] An object of $\calI$ is either an object of $\cDelta^{op}$ or an object of $\LinSeg$.
\item[$(2)$] Morphisms in $\calI$ are give by the formulas
$$\Hom_{\calI}( [m], [n] ) = \Hom_{ \cDelta^{op} }( [m], [n]) \quad \Hom_{\calI}( \seg{m}, \seg{n} ) = \Hom_{\LinSeg}( \seg{m}, \seg{n} )$$
$$ \bHom_{ \calI}( \seg{m}, [n] ) = \bHom_{ \LinSeg}( \seg{m}, \psi( [n] ) ) \quad
 \bHom_{ \calI}( [n], \seg{m} ) = \emptyset.$$ 
\end{itemize}

The map $\psi$ extends to a retraction $\overline{\psi}: \calI \rightarrow \LinSeg$.
Let $\overline{\Mon}(\calC)$ denote the full subcategory of $\Fun( \Nerve(\calI), \calC)$ consisting of those functors $f: \Nerve(\calI) \rightarrow \calC$ which possess the following properties:
\begin{itemize}
\item[$(i)$] For each $n \geq 0$, $f$ carries the canonical map $\seg{n} \rightarrow [n]$ in
$\calI$ to an equivalence in $\calC$.

\item[$(ii)$] The restriction $f | \Nerve(\cDelta)^{op}$ is a monoid object of $\calC$.
\item[$(ii')$] The restriction $f | \Nerve(\LinSeg)$ is a Segal monoid object of $\calC$.
\end{itemize}

We observe that if $(i)$ is satisfied, then $(ii)$ and $(ii')$ are equivalent to one another. Moreover, $(i)$ is equivalent to the assertion that $f$ is a left Kan extension of $f| \Nerve(\LinSeg)$. Since
every functor $f_0: \Nerve(\LinSeg) \rightarrow \calC$ admits a left Kan extension (given, for example, by $f_0 \circ \overline{\psi}$), Proposition \toposref{lklk} implies that the restriction map
$p: \overline{\Mon}(\calC) \rightarrow \SegMon(\calC)$ is a trivial Kan fibration. The map
$q$ is the composition of a section to $p$ (given by composition with $\overline{\psi}$) 
and the restriction map $p': \overline{\Mon}(\calC) \rightarrow \Mon(\calC)$. It will therefore suffice to show that $p'$ is a trivial fibration. In view of Proposition \toposref{lklk}, this will follow from the following pair of assertions:

\begin{itemize}
\item[$(a)$] Every $f_0 \in \Mon(\calC)$ admits a right Kan extension
$f: \Nerve(\calI) \rightarrow \calC$.
\item[$(b)$] Given $f: \Nerve(\calI) \rightarrow \calC$ such that $f_0 = f | \Nerve(\cDelta)^{op}$ is a monoid object of $\calC$, $f$ is a right Kan extension of $f_0$ if and only if $f$ satisfies condition $(i)$ above.
\end{itemize}

To prove $(a)$, we fix an object $J_{\ast} \in \LinSeg$. Let $\calJ$ be the category
$ ( \cDelta)^{op} \times_{ \LinSeg} ( \LinSeg )_{J_{\ast}/ },$
and let $g$ denote the composition
$\Nerve(\calJ) \rightarrow \Nerve(\cDelta)^{op} \stackrel{f_0}{\rightarrow} \calC.$
According to Lemma \toposref{kan2}, it will suffice to show that $g$ admits a limit in $\calC$ (for every choice of $J$). The objects of $\calJ$ can be identified with
morphisms $\alpha: J_{\ast} \rightarrow I_{\ast}$ in $\LinSeg$. Let $\calJ_0 \subseteq \calJ$
denote the full subcategory spanned by those objects for which every element of
$I$ has a unique preimage under $\alpha$. The inclusion $\calJ_0 \subseteq \calJ$
has a right adjoint, so that $\Nerve(\calJ_0)^{op} \rightarrow \Nerve(\calJ)^{op}$ is cofinal. Consequently, it will suffice to show that $g_0 = g | \Nerve(\calJ_0)$ admits a limit in $\calC$.

Let $\calJ_1$ denote the full subcategory of $\calJ_0$ spanned by the morphisms
$\alpha^{j,J}: J_{\ast} \rightarrow \{j\}_{\ast}$. Using our assumption that
$f_0$ is a monoid object of $\calC$, we deduce that $g_0$ is a right Kan extension of
$g_1 = g_0 | \Nerve(\calJ_1)$. In view of Lemma \toposref{kan0}, it will suffice to show that
the map $g_1$ has a limit in $\calC$. But this is clear; our assumption that $f_0$ is a monoid object of $\calC$ guarantees that $f$ exhibits $f( [n])$ as a limit of $g_1$. This proves $(a)$. Moreover, the proof shows that $f$ is a right Kan extension of $f_0$ at $\seg{n}$ if and only if $f$ induces an
equivalence $f( \seg{n} ) \rightarrow f( [n] )$; this immediately implies $(b)$ as well. 
\end{proof}

\begin{proof}[Proof of Proposition \ref{algcompare}]
We will essentially repeat the proof of Proposition \ref{scabby}.
Let $\calI$ be the category introduced in the proof of Proposition \ref{scabby}, and observe again that the functor $\psi: \cDelta^{op} \rightarrow \LinSeg$ extends canonically to a retraction $\overline{\psi}: \calI \rightarrow \LinSeg$.
Let $\overline{\Alg}(\calC)$ denote the full subcategory of $\bHom_{ \Nerve( \LinSeg)}( \Nerve(\calI), \calC^{\otimes})$ consisting of those functors $f: \Nerve(\calI) \rightarrow \calC^{\otimes}$ 
such that $q \circ f = \overline{\psi}$ and the following additional conditions are satisfied:
\begin{itemize}
\item[$(i)$] For each $n \geq 0$, $f$ carries the canonical map $\seg{n} \rightarrow [n]$ in
$\calI$ to an equivalence in $\calC^{\otimes}$.

\item[$(ii)$] The restriction $f | \Nerve(\cDelta)^{op}$ is an algebra object of $\calC$.
\item[$(ii')$] The restriction $f | \Nerve(\LinSeg)$ is a Segal algebra object of $\calC$.
\end{itemize}

We observe that if $(i)$ is satisfied, then $(ii)$ and $(ii')$ are equivalent to one another. Moreover, $(i)$ is equivalent to the assertion that $f$ is a $q$-Kan extension of $f| \Nerve(\LinSeg)$. Since
every functor $f_0: \Nerve(\LinSeg) \rightarrow \calC$ admits a $q$-left Kan extension (given, for example, by $f_0 \circ \overline{\psi}$), Proposition \toposref{lklk} implies that the restriction map
$p: \overline{\Alg}(\calC) \rightarrow \SegAlg(\calC)$ is a trivial Kan fibration. The map
$\theta$ is the composition of a section to $p$ (given by composition with $\overline{\psi}$) 
and the restriction map $p': \overline{\Alg}(\calC) \rightarrow \Alg(\calC)$. It will therefore suffice to show that $p'$ is a trivial fibration. In view of Proposition \toposref{lklk}, this will follow from the following pair of assertions:

\begin{itemize}
\item[$(a)$] Every $f_0 \in \Alg(\calC)$ admits a $q$-right Kan extension
$f \in \bHom_{ \Nerve(\LinSeg)}( \Nerve(\calI), \calC^{\otimes})$.
\item[$(b)$] Given $f \in \bHom_{ \Nerve(\LinSeg)}( \Nerve(\calI), \calC^{\otimes})$
such that $f_0 = f | \Nerve(\cDelta)^{op}$ is an algebra object of $\calC$, $f$ is a $q$-right Kan extension of $f_0$ if and only if $f$ satisfies condition $(i)$ above.
\end{itemize}

To prove $(a)$, we fix an object $J_{\ast} \in \LinSeg$. Let $\calJ$ denote the category
$( \cDelta)^{op} \times_{ \LinSeg} ( \LinSeg )_{J_{\ast}/ },$
and let $g$ denote the composition
$ \Nerve(\calJ) \rightarrow \Nerve(\cDelta)^{op} \stackrel{f_0}{\rightarrow} \calC.$
According to Lemma \toposref{kan2}, it will suffice to show that $g$ admits a $q$-limit in $\calC$ (for every choice of $J$). The objects of $\calJ$ can be identified with
morphisms $\alpha: J_{\ast} \rightarrow I_{\ast}$ in $\LinSeg$. Let $\calJ_0 \subseteq \calJ$
denote the full subcategory spanned by those objects for which every element of
$I$ has a unique preimage under $\alpha$. The inclusion $\calJ_0 \subseteq \calJ$
has a right adjoint, so that $\Nerve(\calJ_0)^{op} \rightarrow \Nerve(\calJ)^{op}$ is cofinal. Consequently, it will suffice to show that $g_0 = g | \Nerve(\calJ_0)$ admits a $q$-limit in $\calC$.

Let $\calJ_1$ denote the full subcategory of $\calJ_0$ spanned by the morphisms
$\alpha^{j,J}: J_{\ast} \rightarrow \{j\}_{\ast}$. Using our assumption that
$f_0$ is a monoid object of $\calC$, we deduce that $g_0$ is a $q$-right Kan extension of
$g_1 = g_0 | \Nerve(\calJ_1)$. In view of Lemma \toposref{kan0}, it will suffice to show that
the map $g_1$ has a $q$-limit in $\calC$. But this is clear; our assumption that $f_0$ is an algebra object of $\calC$ guarantees that $f$ exhibits $f( [n])$ as a limit of $g_1$. This proves $(a)$. Moreover, the proof shows that $f$ is a $q$-right Kan extension of $f_0$ at $\seg{n}$ if and only if $f$ induces an equivalence $f( \seg{n} ) \rightarrow f( [n] )$; this immediately implies $(b)$ as well.
\end{proof}

\section{Modules}\label{hugr2.0}

Let $\calC$ be a monoidal category, and let $A$ be an algebra object of $\calC$. Then
we can consider {\it left $A$-modules} in $\calC$: that is, objects $M \in \calC$
equipped with a multiplication $A \otimes M \rightarrow M$ satisfying suitable unit and associativity conditions. The purpose of this section is to introduce an analogue of these ideas in the case where $\calC$ is replaced by a monoidal $\infty$-category. With an eye towards later applications, we will work in a somewhat more general setting, where $M$ is allowed to be an object of an $\infty$-category which is {\em left-tensored} over $\calC$. We will introduce the formalism of tensored $\infty$-categories in \S \ref{hugr}. In particular, we will see that for every $\infty$-category $\calM$ which is left-tensored over $\calC$, and every algebra object $A \in \Alg(\calC)$, we can introduce an $\infty$-category of (left) $A$-modules $\Mod_{A}(\calM)$. Moreover, every monoidal $\infty$-category $\calC$ can be regarded as (left) tensored over itself in a canonical way (Example \ref{sumai}).

In \S \ref{modcolim} we will study the $\infty$-categories $\Mod_{A}(\calM)$ in more detail. In particular, we will explain how to compute limits and colimits in $\Mod_{A}(\calM)$. Our main results, Corollaries \ref{goop} and \ref{gloop}, assert that (in favorable cases) limits and colimits in $\Mod_{A}(\calM)$ are detected by the forgetful functor $\theta: \Mod_{A}(\calM) \rightarrow \calM$.
Our technique for computing colimits is based on a  general result, Proposition \ref{poststorkus}, which can be used to show that $\Mod_{A}(\calM)$ inherits many pleasant features of $\calM$. 

The forgetful functor $\theta$ is generally not an equivalence, but it always admits a left adjoint. In \S \ref{freetea} we will construct this adjoint, and show that it can be described by
the formula $M \mapsto A \otimes M$.

In \S \ref{hutman}, we will study the theory of module objects in the underlying $\infty$-category of a monoidal model category $\bfA$. In this case, we will show that our $\infty$-categorical theory of module objects of $\Nerve( \bfA^{\degree})$ is closely related to the classical theory of (strictly associative) modules in $\bfA$. 

The theory of  $\infty$-categories left-tensored over a given monoidal $\infty$-category $\calC$ is really a special case of the general theory of module objects. More precisely, we have already seen that $\calC$ can be identified with an algebra object of $\Cat_{\infty}$ (Remark \ref{otherlander}). In \S \ref{hugree} we will extend this result by showing $\infty$-categories left-tensored over $\calC$ can be identified with $\calC$-module objects of $\Cat_{\infty}$. These identifications allow us to apply our study of module objects to deduce results about tensored $\infty$-categories. For example, we can use this method to show that
$\calC$ is freely generated (as an $\infty$-category tensored over $\calC$) by the unit object
$1_{\calC} \in \calC$ (Corollary \ref{specialtime}). 

In \S \ref{cupper} we will shift our focus somewhat; rather than studying the $\infty$-category $\Mod_{A}(\calM)$ over a fixed algebra object $A \in \Alg(\calC)$, we will instead fix an object $M \in \calM$ and study algebras which act on $M$. We will show that algebra objects of $\calC$ which act on $M$ can be identified with algebra objects in a larger $\infty$-category $\calC[M]$ (Proposition \ref{poofer}), and that in some cases one can construct a universal example of such an algebra, which we will denote by $\End(M)$.

For certain applications, it is convenient to work with {\em nonunital} algebras in a monoidal $\infty$-category $\calC$: that is, objects $A \in \calC$ equipped with a multiplication
$A \otimes A \rightarrow A$ which is coherently associative, but no unit map $1_{\calC} \rightarrow A$. In \S \ref{digunit} we will give a precise definition for the category $\Alg^{\nounit}(\calC)$ of nonunital algebra objects of $\calC$, and introduce an accompanying theory of nonunital modules.
In \S \ref{giddug}, we will show that the $\infty$-categories $\Alg^{\nounit}(\calC)$ and
$\Alg(\calC)$ are closely related: namely, $\Alg(\calC)$ is equivalent to a subcategory of
$\Alg^{\nounit}(\calC)$ (Theorem \ref{uniqueunit}). In other words, if a nonunital algebra object of $\calC$ admits a unit, then that unit is essentially unique. 

\subsection{Tensored $\infty$-Categories}\label{hugr}

Let $\calC$ be a monoidal $\infty$-category. Our first goal is to introduce
the notion of an $\infty$-category $\calM$ {\it left-tensored} over $\calC$. Roughly speaking, this means we have an action
$ \otimes: \calC \times \calM \rightarrow \calM$
of $\calC$ on $\calM$, which is unital and associative up to coherent homotopy.

\begin{definition}\label{ulult}\index{$\infty$-category!tensored}\index{tensored $\infty$-category}
Let $p: \calC^{\otimes} \rightarrow \Nerve( \cDelta)^{op}$ be a monoidal $\infty$-category.
An {\it $\infty$-category left-tensored over $\calC^{\otimes}$} is a categorical fibration
$q: \calM^{\otimes} \rightarrow \calC^{\otimes}$ with the following properties:
\begin{itemize}
\item[$(1)$] The composition $(p \circ q): \calM^{\otimes} \rightarrow \Nerve(\cDelta)^{op}$ is a coCartesian fibration.
\item[$(2)$] The map $q$ carries $(p \circ q)$-coCartesian edges of $\calM^{\otimes}$ to $p$-coCartesian edges of $\calC^{\otimes}$.
\item[$(3)$] For each $n \geq 0$, the inclusion $\{ n \} \subseteq [n]$ induces an
equivalence of $\infty$-categories $\calM^{\otimes}_{[n]} \rightarrow \calC^{\otimes}_{[n]} \times  \calM^{\otimes}_{ \{n\} }$.
\end{itemize}
\end{definition}

\begin{remark}
Let $q: \calM^{\otimes} \rightarrow \calC^{\otimes}$ be as in Definition \ref{ulult}. We will refer to the
fiber $\calM = \calM^{\otimes}_{[0]}$ as the {\it underlying $\infty$-category} of $\calM^{\otimes}$. Let $\calC = \calC^{\otimes}_{[1]}$ be the underlying $\infty$-category of $\calC^{\otimes}$. Then
$\calM^{\otimes}_{[n]}$ is equivalent to the product $\calC^{n} \times \calM$. 
The coCartesian fibration $q$ induces functors
$$\calM^{\otimes}_{ \{0\} } \leftarrow \calM^{\otimes}_{[1]} \rightarrow \calC^{\otimes}_{[1]} \times \calM^{\otimes}_{ \{1\} },$$
where the right map is an equivalence. We therefore obtain a bifunctor
$\otimes: \calC \times \calM \rightarrow \calM$, well-defined up to homotopy. Moreover, the structure associated to $\calM^{\otimes}_{[n]}$ for $n > 1$ ensures that the bifunctor $\otimes: \calC \times \calM \rightarrow \calM$ is coherently associative; in particular, the homotopy category
$\h{\calC}$ is tensored over $\h{\calM}$ in the sense of classical category theory.
We will generally abuse terminology by saying simply that $\calM$ is left-tensored over $\calC$.
\end{remark}\index{ZZZotimes@$\otimes$}

\begin{example}\label{sumai}\index{monoidal $\infty$-category!tensored over itself}
Let $\calC$ be a monoidal $\infty$-category. Then the tensor product
$\otimes: \calC \times \calC \rightarrow \calC$
exhibits $\calC$ as left-tensored over itself.
More precisely, suppose that the monoidal structure on $\calC$ is given by a projection $p: \calC^{\otimes} \rightarrow \Nerve(\cDelta)^{op}$. 
Let $e: \cDelta \rightarrow \cDelta$ denote the ``shift'' functor $[n] \mapsto [n] \star [0],$
where $\star$ denotes the operation of concatenating linearly ordered sets. Then $e$ induces a functor $\Nerve(\cDelta)^{op} \rightarrow \Nerve(\cDelta)^{op}$, which we will also denote by $e$. The inclusion $[n] \subseteq [n] \star [0]$ determines a natural transformation $\alpha: e \rightarrow \id$, which we can view as a functor from $\Delta^1 \times \Nerve(\cDelta)^{op}$ to $\Nerve(\cDelta)^{op}$. We now define a simplicial set $\calC^{\otimes,L}$ equipped with a map
$\calC^{\otimes,L} \rightarrow \Nerve( \cDelta)^{op}$ via the formula
$$ \Hom_{ \Nerve(\cDelta)^{op}}(K, \calC^{\otimes,L})
= \Hom'_{ \Nerve( \cDelta)^{op}}(\Delta^1 \times K, \calC^{\otimes})$$
where $\Delta^1 \times K$ maps to $\Nerve( \cDelta)^{op}$ via the composition
$ \Delta^1 \times K \rightarrow \Delta^1 \times \Nerve(\cDelta)^{op} \stackrel{\alpha}{\rightarrow} \Nerve(\cDelta)^{op},$
and $\Hom'_{\Nerve(\cDelta)^{op}}( \Delta^1 \times K, \calC^{\otimes} ) \subseteq
\Hom_{\Nerve(\cDelta)^{op}}( \Delta^1 \times K, \calC^{\otimes})$ is the subset consisting of those maps $\Delta^1 \times K \rightarrow \calC^{\otimes}$ which carry each edge
$\Delta^1 \times \{k\}$ to a $p$-coCartesian edge of $\calC^{\otimes}$.
Using Proposition \toposref{doog}, we deduce that the projection $\calC^{\otimes,L} \rightarrow \Nerve(\cDelta)^{op}$ is a coCartesian fibration, and that the projections
$ \calC^{\otimes}_{L} \stackrel{\psi}{\leftarrow} \calC^{\otimes,L} \stackrel{q}{\rightarrow} \calC^{\otimes}$
preserve coCartesian edges. Here $\calC^{\otimes}_{L}$ denotes the pullback of
$p$ along the map $e: \Nerve(\cDelta)^{op} \rightarrow \Nerve(\cDelta)^{op}$.
It is easy to see that $q$ is a categorical fibration (since the inclusion $\{1\} \times K \subseteq \Delta^1 \times K$ is a cofibration of simplicial sets). Since $\psi$ is an equivalence, we deduce
that $\calC^{\otimes,L}_{[n]}$ is canonically equivalent to the $\infty$-category $\calC^{n+1}$. Using this equivalence, we can easily verify that $q$ satisfies condition $(4)$ of Definition \ref{ulult}. Consequently, $q: \calC^{\otimes, L} \rightarrow \calC^{\otimes}$ can be viewed as an $\infty$-category left-tensored over $\calC$. The map $\psi$ induces an equivalence
$\calC^{\otimes,L}_{[0]} \rightarrow \calC^{\otimes}_{[1]} \simeq \calC$. We will summarize the situation by saying that $\calC^{\otimes, L}$ {\it exhibits $\calC$ as left-tensored over itself}.
\end{example}\index{ZZZCotimesL@$\calC^{\otimes,L}$}

\begin{definition}\label{defmond}\index{module object}
Let $p: \calC^{\otimes} \rightarrow \Nerve(\cDelta)^{op}$ be a monoidal $\infty$-category, and let
$q: \calM^{\otimes} \rightarrow \calC^{\otimes}$ be an $\infty$-category equipped with a left action of $\calC^{\otimes}$. A {\it module object} of $\calM^{\otimes}$ is a functor
$F: \Nerve(\cDelta)^{op} \rightarrow \calM^{\otimes}$ with the following properties:
\begin{itemize}
\item[$(1)$] The composition $q \circ F$ is an algebra object of $\calC^{\otimes}$. In particular,
$p \circ q \circ F$ is the identity on $\Nerve(\cDelta)^{op}$.
\item[$(2)$] Let $\alpha: [m] \rightarrow [n]$ be a convex map in $\cDelta$ such that
$\alpha(m) = n$. Then $F(\alpha)$ is a $p \circ q$-coCartesian morphism of $\calM^{\otimes}$.
\end{itemize}
In this case, we will also abuse terminology by saying that $F$ is a {\it module object} of the underlying $\infty$-category $\calM = \calM^{\otimes}_{[0]}$. We let $\Mod(\calM)$ denote the full subcategory of $\bHom_{\Nerve(\cDelta)^{op} }(\Nerve(\cDelta)^{op}, \calM^{\otimes})$ spanned by the module objects of $\calM$. 
\end{definition}\index{ZZZModM@$\Mod(\calM)$}

\begin{remark}\label{hugrr}
There is an analogous notion of an {\it $\infty$-category {\em right-tensored over $\calC^{\otimes}$}}, obtained by replacing the inclusion $\{n\} \subseteq [n]$ with the inclusion $\{0\} \subseteq [n]$ in condition $(\ast)$ of Definition \ref{ulult}. Similarly, we have a dual notion of {\em module object} in such a category. When necessary, we will distinguish these two notions by
referring to {\em left modules} and {\em right modules}, respectively.
\end{remark}

\begin{remark}
Let $\calC$ be a monoidal $\infty$-category, and let $\calC^{\otimes, L}$ be defined as in Example \ref{sumai}. We will abuse terminology by referring to module objects of $\calC^{\otimes, L}$ as {\it module objects of $\calC$}. The $\infty$-category of module objects of $\calC^{\otimes, L}$ will be denoted by $\Mod(\calC)$, or
$\Mod^{L}(\calC)$ when it is necessary to emphasize that we are working with {\em left} modules. Note that this notation is not quite compatible with the notation introduced in Definition \ref{defmond}, since $\calC^{\otimes,L}_{[0]}$ is equivalent but not isomorphic to $\calC$. 
\end{remark}

\begin{example}\label{algitself}
Let $p: \calC^{\otimes} \rightarrow \Nerve(\cDelta)^{op}$ be a monoidal structure on an $\infty$-category $\calC = \calC^{\otimes}_{[1]}$. Let $\alpha: \Delta^1 \times \Nerve(\cDelta)^{op} \rightarrow \Nerve(\cDelta)^{op}$ be defined as in Example \ref{sumai}. Composition with
$\alpha$ determines a section $s$ of the forgetful functor
$\theta: \Mod(\calC) \rightarrow \Alg(\calC)$. This section carries each $A \in \Alg(\calC)$
to a left $A$-module $s(A)$. We can view $s(A)$ as $A$, equipped with the canonical left action on itself. For this reason, we will generally not distinguish in our notation between $A$ and $s(A)$.
\end{example}

\begin{remark}
Let $q: \calM^{\otimes} \rightarrow \calC^{\otimes}$ be as in Definition \ref{defmond}. Then
composition with $q$ induces a functor $\theta: \Mod(\calM) \rightarrow \Alg(\calC)$, which we will refer to as the {\it forgetful functor}. If $A$ is an algebra object of $\calC$, we let
$\Mod_{A}(\calM)$ denote the fiber $\theta^{-1} \{A\}$. It is easy to see that $\theta$ is a categorical fibration of simplicial sets, so that each $\Mod_{A}(\calM)$ is an $\infty$-category, which we will call
the {\it $\infty$-category of $($left$)$ $A$-modules in $\calM$}. We will soon see (Corollary \ref{thetacart}) that $\theta$ is a Cartesian fibration, so that $\Mod_{A}(\calM)$ is an $\infty$-category which is contravariantly functorial in $A$.
In the case where $\calC$ is acting on itself (as in Example \ref{sumai}), we will also denote this $\infty$-category by $\Mod_{A}(\calC)$ or simply $\Mod_{A}$, and refer to it as the {\it $\infty$-category of $($left$)$ $A$-modules}.\index{ZZZModAM@$\Mod_{A}(\calM)$}
\end{remark}

In classical category theory, one can define a similar notion of a category (left) tensored over a fixed monoidal category $\calC$. However, this notion plays a secondary role to the theory of {\em enriched} categories (see, for example, \S \toposref{enrichcat}). For example, every category $\calC$ can be regarded as enriched over the category $\Set$ of sets. However, $\calC$ is enriched and tensored over $\Set$ if and only if $\calC$ admits arbitrary coproducts. 

It is possible to modify Definition \ref{ulult} to obtain the definition of an $\infty$-category {\em enriched over} a given monoidal $\infty$-category $\calC$. This notion is useful in a variety of situations. For example, suppose we wanted to develop the theory of {\it $(\infty, n)$-categories}: that is, higher categories in which we allow noninvertible morphisms of order $\leq n$. One reasonable definition is inductive: one defines an $(\infty, n)$-category to be an $\infty$-category which is enriched over $\Cat_{(\infty, n-1)}$, the $\infty$-category of $(\infty, n-1)$-categories (endowed with the Cartesian monoidal structure). We will not pursue the theory of enriched $\infty$-categories in this paper, since the requisite modifications to Definition \ref{ulult} are somewhat cumbersome to work with in practice. However, it requires very little additional effort to formulate the following {\em strengthening} of Definition \ref{ulult}:

\begin{definition}\label{supiner}\index{$\infty$-category!enriched}\index{$\infty$-category!cotensored}
Let $\calC$ be a monoidal $\infty$-category, and let $\calM$ be an $\infty$-category which is left-tensored over $\calC$. 

\begin{itemize}
\item[$(1)$] Let $M$ and $N$ be objects of $\calM$. A {\it morphism object} for
$M$ and $N$ is an object $\Mor_{\calM}(M,N) \in \calC$ together with a map
$\alpha: \Mor_{\calM}(M, N) \otimes M \rightarrow N$
with the following universal property: for every object $C \in \calC$, composition with $\alpha$ induces a homotopy equivalence
$\bHom_{\calC}( C, \Mor_{\calM}(M,N) ) \rightarrow
\bHom_{\calM}( C \otimes M, N).$

\item[$(2)$] Let $M$ be an object of $\calM$ and $C$ an object of $\calC$. An
{\it exponential object} is an object $^{C}\!M \in \calM$ together with a map
$\beta: C \otimes ^{C}\!M \rightarrow M$ with the following universal property: for every object $N \in \calM$, composition with $\beta$ induces a homotopy equivalence
$\bHom_{ \calM} ( N, ^{C}\!M) \rightarrow \bHom_{\calM}( C \otimes N, M).$
\end{itemize}\index{morphism object}\index{exponential object}

We will say that $\calM$ is {\it tensored and  enriched over $\calC$} if, for every pair of objects $M, N \in \calM$, there exists a morphism object $\Mor(M,N) \in \calC$. We will say that $\calM$ is
{\it tensored and cotensored over $\calC$} if, for every $M \in \calM$ and $C \in \calC$, there exists an exponential object $^{C}\!M \in \calM$. 
\end{definition}

\begin{remark}
Let $\calC$ and $\calM$ be as in Definition \ref{supiner}. Suppose that
$\calM$ is enriched over $\calC$. Then the morphism object
$\Mor_{\calM}(M, N) \in \calC$ can be constructed as functor of
$M$ and $N$. To see this, we observe that we have a trifunctor
$(\calC \times \calM)^{op}  \times \calM \rightarrow \SSet,$
given by composing the tensor product
$\otimes: \calC \times \calM \rightarrow \calM$, the Yoneda embedding
$j: \calM^{op} \rightarrow \Fun( \calM, \SSet)$, and the
evaluation map $\Fun(\calM, \SSet)  \times \calM \rightarrow \SSet$. 
We may identify this map with a bifunctor
$e: \calM^{op} \times \calM \rightarrow \Fun( \calC^{op}, \SSet ).$
If $\calM$ is enriched over $\calC$, then the image of $e$ is contained in full
subcategory $\Fun'( \calC^{op}, \SSet ) \subseteq \Fun(\calC^{op}, \SSet)$ spanned
by the essential image of the Yoneda embedding $j': \calC \rightarrow \Fun'( \calC^{op}, \SSet)$. 
Composing $e$ with a homotopy inverse to $j'$, we obtain the desired functor
$\Mor_{\calM}: \calM^{op} \times \calM \rightarrow \calC.$
Similarly, if $\calM$ is cotensored over $\calC$, then the exponential
$^{C}\!M$ can be regarded as a functor
$\calC^{op} \times \calM \rightarrow \calM.$
\end{remark}

\begin{example}
Let $\calC$ be a monoidal $\infty$-category, and regard $\calC$ as left-tensored over itself (Example \ref{sumai}). Then $\calC$ is enriched over itself if and only if it is {\em right closed}, and cotensored over itself if and only if it is {\em left closed} (Definition \ref{sumpy}). 
\end{example}

The following criterion provides a large supply of examples of enriched and cotensored $\infty$-categories:

\begin{proposition}\label{enterich}
Let $\calC$ be a monoidal $\infty$-category, and let $\calM$ be an $\infty$-category which is left-tensored over $\calC$. Suppose further that $\calC$ and $\calM$ are presentable.
\begin{itemize}
\item[$(1)$] If, for every $C \in \calC$, the functor
$C \otimes \bigdot: \calM \rightarrow \calM$ preserves small colimits, then
$\calM$ is cotensored over $\calC$.
\item[$(2)$] If, for every $M \in \calM$, the functor 
$\bigdot \otimes M: \calC \rightarrow \calM$ preserves small colimits, then
$\calM$ is enriched over $\calC$.
\end{itemize}
\end{proposition}

\begin{proof}
This follows immediately from the representability criterion of Proposition \toposref{representable}.
\end{proof}

We conclude this section with a few technical observations which will be needed later, when we need study tensored $\infty$-categories in detail. Let 
$ \calM^{\otimes} \stackrel{q}{\rightarrow} \calC^{\otimes} \stackrel{p}{\rightarrow} \Nerve(\cDelta)^{op}$ be as in Definition \ref{ulult}. Then $p$ and $(p \circ q)$ are coCartesian fibrations, but $q$ is usually not a coCartesian fibration. Nevertheless, $q$ shares some of the pleasant features of coCartesian fibrations.

\begin{lemma}\label{excel}
Let $\calC$ be a monoidal $\infty$-category, and let $\calM$ be an $\infty$-category which is left tensored over $\calC$. Then the associated functor $p: \calM^{\otimes} \rightarrow \calC^{\otimes}$ is a locally coCartesian fibration.
\end{lemma}

\begin{proof}
For each $n \geq 0$, the map $p_{[n]}: \calM^{\otimes}_{[n]} \rightarrow \calC^{\otimes}_{[n]}$
is equivalent to the projection $\calC^{\otimes}_{[n]} \times \calM \rightarrow \calC^{\otimes}_{[n]}$, and therefore a coCartesian fibration. The desired result now follows from Proposition \toposref{fibertest}.
\end{proof}

\begin{remark}\label{marbits}
In the situation of Lemma \ref{excel}, suppose we are given a morphism
$\overline{\alpha}: C \rightarrow D$ in $\calC^{\otimes}$, covering a map $\alpha: [m] \rightarrow [n]$ in $\cDelta$. We can identify $C$ with an $n$-tuple of objects
$(C_1, \ldots, C_n) \in \calC^n$, $D$ with an $m$-tuple $(D_1, \ldots, D_m) \in \calC^{m}$, and
$\overline{\alpha}$ with a collection of morphisms
$C_{\alpha(i-1)+1} \otimes \ldots \otimes C_{\alpha(i)} \rightarrow D_i.$
The fibers of the map $p: \calM^{\otimes} \rightarrow \calC^{\otimes}$ over the objects
$C$ and $D$ are both canonically equivalent to $\calM$, and the induced functor
$\overline{\alpha}_{!}: \calM \rightarrow \calM$ is given (up to equivalence) by the formula
$ M \mapsto C_{ \alpha(m)+1} \otimes \ldots \otimes C_{n} \otimes M.$
\end{remark}

\begin{remark}\label{greft}
Let $p: \calM^{\otimes} \rightarrow \calC^{\otimes}$ be as in Lemma \ref{excel}, and suppose given a commutative triangle
$$ \xymatrix{ & M' \ar[dr]^{g} & \\
M \ar[ur]^{f} \ar[rr]^{h} & & M'' }$$
in $\calM^{\otimes}$, covering a triangle
$$ \xymatrix{ & [m] \ar[dl]^{\alpha} & \\
[n] & & [k] \ar[ll]^{\gamma} \ar[ul]^{\beta} }$$
in the category $\cDelta$. Suppose furthermore that $f$ and $g$ are locally $p$-coCartesian, and
that $\alpha$ induces a bijection
$$ \{ i \in [m] : \beta(k) < i \} \simeq \{ j \in [n]: \gamma(k) < j \leq \alpha(m) \}$$
(so that, in particular, $\beta(k) < i \leq m$ implies $\gamma(k) < \alpha(i)$). Then the description given in Remark \ref{marbits} implies that $h$ is also locally $p$-coCartesian.
\end{remark}

\begin{lemma}\label{smarties}
Let $p: \calM^{\otimes} \rightarrow \calC^{\otimes}$ be as in Lemma \ref{excel}, and let $\overline{\alpha}$ be a morphism in $\calM^{\otimes}$ which covers a map $\alpha: [m] \rightarrow [n]$ in $\cDelta$. Then:
\begin{itemize}
\item[$(1)$] Suppose that $\alpha(m) = n$ and $\overline{\alpha}$ is locally $p$-Cartesian. Then
$\overline{\alpha}$ is $p$-Cartesian.
\item[$(2)$] Suppose that $m \leq n$, that $\alpha: [m] \rightarrow [n]$ is the canonical inclusion, and that $\overline{\alpha}$ is locally $p$-coCartesian. Then $\overline{\alpha}$ is $p$-coCartesian.
\end{itemize}
\end{lemma}

\begin{proof}
Assertion $(1)$ follows from Remark \ref{greft} and Lemma \toposref{gruft}, while
$(2)$ follows from Remark \ref{greft} and Lemma \toposref{charloccart}.
\end{proof}

\subsection{Nonunital Algebras and Modules}\label{digunit}

Let $A$ be an associative ring. A {\it nonunital $A$-module} is an abelian group $M$ equipped with a bilinear multiplication map $A \times M \rightarrow M$
which satisfies the associativity formula
$a(bm) = (ab)m.$
This is slightly weaker than the condition that $M$ be an $A$-module, since we do not require the unit object $1 \in A$ to act by the identity on $M$. For example, there is a trivial nonunital $A$-module
structure on any abelian group $M$, given by the zero map
$ A \times M \stackrel{0}{\rightarrow} M.$
Of course, the category of $A$-modules can be identified with a full subcategory of the
category of nonunital $A$-modules. Our goal in this section is to formulate and prove an $\infty$-categorical analogue of this statement (Proposition \ref{uniquemo}). This result will play an essential role in our study of the Barr-Beck theorem in \S \ref{barri}.

We first observe that the notion of a nonunital $A$-module makes no reference to
the identity element of $A$; it therefore makes sense in the case where $A$ is allowed to be a nonunital algebra. Our first goal in this section is to set up a theory of nonunital algebras in the $\infty$-categorical context.


\begin{definition}\label{nounitalg}
We define a subcategory $\cDelta^{\nounit} \subseteq \cDelta$ as follows:
\begin{itemize}
\item[$(1)$] The objects of $\cDelta^{\nounit}$ are the objects of $\cDelta$: that is, they
are linearly ordered sets $[n]$ where $n \geq 0$.
\item[$(2)$] A morphism $\alpha: [m] \rightarrow [n]$ of $\cDelta$ belongs to $\cDelta^{\nounit}$
if and only if $\alpha$ is an {\em injective} map of linearly ordered sets.
\end{itemize}\index{ZZZcDeltanounit@$\cDelta^{\nounit}$}

Let $q: \calC^{\otimes} \rightarrow \Nerve(\cDelta)^{op}$ be a monoidal structure on an
$\infty$-category $\calC = \calC^{\otimes}_{[1]}$. A {\it nonunital algebra object of $\calC$}
is a functor $A: \Nerve( \cDelta^{\nounit})^{op} \rightarrow \calC^{\otimes}$ such that
$q \circ A$ coincides with the inclusion $\Nerve(\cDelta^{\nounit})^{op} \subseteq \Nerve(\cDelta)^{op}$, and $A$ carries convex morphisms in $\cDelta^{\nounit}$ to $q$-coCartesian morphisms
in $\calC^{\otimes}$. We let $\Alg^{\nounit}(\calC)$ denote the full subcategory of
$\bHom_{ \Nerve(\cDelta)^{op} }( \Nerve( \cDelta^{\nounit})^{op}, \calC^{\otimes})$ spanned by the nonunital algebra objects of $\calC$.\index{algebra object!nonunital}\index{nonunital!algebra object}\index{ZZZAlgnounitC@$\Alg^{\nounit}(\calC)$}
\end{definition}

\begin{warning}
In \cite{topoi}, we denoted the category $\Delta^{\nounit}$ by $\Delta_{s}$.
\end{warning}

\begin{remark}
The definition of a nonunital algebra object of a monoidal $\infty$-category $\calC$ makes sense in the more general setting of {\em nonunital monoidal $\infty$-categories}. We will have no need for this additional level of generality.
\end{remark}

\begin{remark}
Let $\calC$ be a monoidal $\infty$-category. Evaluation on the object $[1] \in \cDelta^{\nounit}$ defines a forgetful functor $\Alg^{\nounit}(\calC) \rightarrow \calC$. We will generally abuse notation by not distinguishing between a nonunital algebra $A \in \Alg^{\nounit}(\calC)$ and its image in $\calC$.
\end{remark}

\begin{remark}\label{sugarbear}
Let $q:\calC^{\otimes} \rightarrow \Nerve(\cDelta)^{op}$ be a monoidal structure on an $\infty$-category $\calC = \calC^{\otimes}_{[1]}$, and let $A: \Nerve(\cDelta)^{op} \rightarrow
\calC^{\otimes}$ be a section of $q$. Then $A$ is an algebra object of $\calC$ if and only if
$A | \Nerve( \cDelta^{\nounit})^{op}$ is a nonunital algebra object of $\calC$. Consequently, restriction to the subcategory $\cDelta^{\nounit} \subseteq \cDelta$ defines a forgetful functor
$\theta: \Alg(\calC) \rightarrow \Alg^{\nounit}(\calC).$
We will later show that the functor $\theta$ induces an equivalence of
$\Alg(\calC)$ with a subcategory of $\Alg^{\nounit}(\calC)$ (Theorem \ref{uniqueunit}).
\end{remark}

\begin{definition}\label{slipper}\index{module object!nonunital}\index{nonunital!module object}
Let $\calM^{\otimes} \stackrel{q}{\rightarrow} \calC^{\otimes} \stackrel{p}{\rightarrow} \Nerve(\cDelta)^{op}$ exhibit $\calM = \calM^{\otimes}_{[0]}$ as left-tensored over the monoidal $\infty$-category
$\calC = \calC^{\otimes}_{[1]}$. A {\it nonunital module object} of $\calM$ is a functor
$M: \Nerve( \cDelta^{\nounit})^{op} \rightarrow \calM^{\otimes}$ with the following properties:
\begin{itemize}
\item[$(1)$] The composition $q \circ M$ is a nonunital algebra object of $\calC$.
\item[$(2)$] Let $\alpha: [m] \rightarrow [n]$ be a convex map in $\cDelta^{\nounit}$ such that
$\alpha(m) = n$. Then $M(\alpha)$ is a $(p \circ q)$-coCartesian morphism of $\calM^{\otimes}$.
\end{itemize}
We let $\Mod^{\nounit}(\calM)$ denote the full subcategory of
$\bHom_{ \Nerve(\cDelta)^{op} }( \Nerve( \cDelta^{\nounit})^{op}, \calM^{\otimes})$ spanned by the nonunital module objects. If $A \in \Alg^{\nounit}(\calC)$, then we let
$\Mod_{A}^{\nounit}(\calM)$ denote the fiber
$\Mod^{\nounit}(\calM) \times_{ \Alg^{\nounit}( \calM) } \{ A\}$. If
$A \in \Alg(\calC)$, then we let $\Mod_{A}^{\nounit}(\calM) = \Mod_{\theta(A)}^{\nounit}(\calM)$, where $\theta: \Alg(\calC) \rightarrow \Alg^{\nounit}(\calC)$ denotes the forgetful functor (see Remark \ref{sugarbear}). \index{ZZZModnounitM@$\Mod^{\nounit}(\calM)$}\index{ZZZModnounitAM@$\Mod^{\nounit}_{A}(\calM)$}
\end{definition}

\begin{remark}
The terminology of Definition \ref{slipper} is perhaps slightly confusing: remember that a module object of $\calM$ always refers to a pair $(A, M)$, where $A$ is an algebra object of $\calC$ and $M$ is a module over $A$. To obtain the notion of a nonunital module object, we allow $A$ to be a nonunital algebra (and drop the requirement that the unit of $A$, if it exists, acts by the identity on $M$).
\end{remark}

\begin{remark}
In the situation of Definition \ref{slipper}, evaluation at $[0] \in \cDelta^{\nounit}$ induces a forgetful functor $\Mod^{\nounit}(\calM) \rightarrow \calM$. We will generally abuse notation by not distinguishing between a nonunital module object $M \in \Mod^{\nounit}(\calM)$ and its image in $\calM$.
\end{remark}

\begin{remark}\label{sugarbear2}
Let $\calM^{\otimes} \stackrel{q}{\rightarrow} \calC^{\otimes} \stackrel{p}{\rightarrow} \Nerve(\cDelta)^{op}$ exhibit $\calM = \calM^{\otimes}_{[0]}$ as left-tensored over the monoidal $\infty$-category
$\calC = \calC^{\otimes}_{[1]}$, and let $M: \Nerve(\cDelta)^{op} \rightarrow \calM^{\otimes}$ be a section of $p \circ q$. Then $M$ is a module object of $\calM$ if and only if
$M | \Nerve(\cDelta^{\nounit})^{op}$ is a nonunital module object of $\calM$. In particular, restriction induces a functor $\Mod(\calM) \rightarrow \Mod^{\nounit}(\calM)$.
\end{remark}

\begin{definition}\label{quasiunitalg}\index{left unit}\index{unit}\index{right unit}\index{quasi-unit}
Let $\calC$ be a monoidal $\infty$-category, let $A$ be a nonunital algebra object of $\calC$, and let $u: 1_{\calC} \rightarrow A$ be an arbitrary map. We say that $u$ is a {\it left unit} if
the composition
$$ A \simeq A \otimes 1_{\calC} \stackrel{u}{\rightarrow} A \otimes A \rightarrow A$$
is homotopic to the identity in $\calC$. Similarly, we will say that $u$ is a 
{\em right unit} if the composition
$$ A \simeq 1_{\calC} \otimes A \stackrel{u}{\rightarrow} A \otimes A \rightarrow A$$
is homotopic to the identity in $\calC$. We will say that $u$ is a {\it quasi-unit} if it is both
a left unit and a right unit. We will say that $A$ is {\it quasi-unital} if there exists a quasi-unit
$u: 1_{\calC} \rightarrow A$.
Let $f: A \rightarrow A'$ be a morphism between
nonunital algebra objects of $\calC$. We will say that $f$ is {\it quasi-unital} if
there exists a quasi-unit $u: 1_{\calC} \rightarrow A$, and the composition
$f \circ u: 1_{\calC} \rightarrow A'$ is a quasi-unit for $A'$.
We let $\Alg^{\qunit}(\calC)$ denote the subcategory of $\Alg^{\nounit}(\calC)$ whose objects
are quasi-unital algebra objects of $\calC$, and whose morphisms are quasi-unital maps of algebras.
\index{algebra object!quasi-unital}\index{ZZZAlgqunitC@$\Alg^{\qunit}(\calC)$}
\end{definition}

\begin{remark}\label{pounit}
Let $A$ be a nonunital algebra object of a monoidal $\infty$-category $\calC$. 
Suppose that $u: 1_{\calC} \rightarrow A$ is a left unit and
$v: 1_{\calC} \rightarrow A$ is a right unit. Then the composition
$$ 1_{\calC} \simeq 1_{\calC} \otimes 1_{\calC} \stackrel{u \otimes v}{\rightarrow}
A \otimes A \rightarrow A$$
is homotopic to both $u$ and $v$, so that $u$ and $v$ are homotopic to one another.
It follows that $A$ is quasi-unital if and only if $A$ admits both a left and right unit. In this case, a quasi-unit $u: 1_{\calC} \rightarrow A$ is uniquely determined (up to homotopy).
\end{remark}

\begin{remark}
We will see later that there is essentially no difference between algebras and quasi-unital algebras. More precisely, we will show that for every monoidal $\infty$-category $\calC$, the forgetful functor
$\Alg(\calC) \rightarrow \Alg^{\qunit}(\calC)$ is a trivial Kan fibration (Theorem \ref{uniqueunit}).
\end{remark}

\begin{definition}\label{spreck}
Let $\calC$ be a monoidal $\infty$-category, let $\calM$ be an $\infty$-category which is left-tensored over $\calC$, and let $A$ be a quasi-unital algebra object of $\calC$. We will say that
an object $M \in \Mod^{\nounit}_{A}(\calM)$ is {\it quasi-unital} if the composition
$$ \phi: M \simeq 1_{\calC} \otimes M \stackrel{u}{\rightarrow} A \otimes M \rightarrow M$$
is homotopic to the identity in $\calM$, where $u$ is a quasi-unit for $A$.
We let $\Mod^{\qunit}_{A}(\calM)$ denote the full subcategory of
$\Mod^{\nounit}_{A}(\calM)$ spanned by the quasi-unital $A$-modules.
\index{module object!quasi-unital}
\index{ZZZModqunitAM@$\Mod^{\qunit}_{A}(\calM)$}
\end{definition}

\begin{remark}
In view of Remark \ref{pounit}, the condition of Definition \ref{spreck} does not depend on the choice of a quasi-unit $u: 1_{\calC} \rightarrow A$.
\end{remark}

\begin{remark}
In the situation of Definition \ref{spreck}, the condition that $\psi$ be homotopic to the identity is equivalent to the (apparently weaker) condition that $\psi$ be an equivalence. For suppose that
$\psi$ is an equivalence. Since the composition
$$ 1_{\calC} \simeq 1_{\calC} \otimes 1_{\calC} \stackrel{u \otimes u}{\rightarrow} A \otimes A \rightarrow A$$
is equivalent to $u$, we conclude that $\psi^2$ is homotopic to $\psi$ (that is, $\psi^2$ and
$\psi$ belong to the same connected component of $\bHom_{\calM}(M,M)$ ). If $\psi$
is invertible in the homotopy category $\h{\calM}$, this forces $\psi$ to be homotopic to the identity.
\end{remark}

Our goal in this section is to prove the following result:

\begin{proposition}\label{uniquemo}
Let $\calM$ be an $\infty$-category which is left-tensored over a monoidal $\infty$-category 
$\calC$, and let $A$ be an algebra object of $\calC$. The restriction map $\Mod_{A}(\calM) \rightarrow \Mod^{\nounit}_{A}(\calM)$ induces a trivial Kan fibration $\theta: \Mod_{A}(\calM) \rightarrow \Mod^{\qunit}_{A}(\calM)$.
\end{proposition}

In other words, the theory of $A$-modules is equivalent to the theory of nonunital $A$-modules in which the unit of $A$ happens to act by the identity.

\begin{proof}
It is clear that $\theta$ is a categorical fibration. It will therefore suffice to show that $\theta$ is a categorical equivalence.

We define a subcategory $\calI \subseteq [1] \times \cDelta^{op}$ as follows:
\begin{itemize}
\item[$(1)$] Every object of $[1] \times \cDelta^{op}$ belongs to $\calI$.
\item[$(2)$] A morphism $\alpha: (i, [m]) \rightarrow (j, [n])$ in $[1] \times \cDelta^{op}$ belongs to $\calI$ if and only if either $i = 0$ or the map $[n] \rightarrow [m]$ is injective.
\end{itemize}
For $i = 0, 1$, we let $\calI_{i}$ denote the full subcategory of $\calI$ spanned by the objects
$\{ ( i, [n] ) \}_{n \geq 0}$. 

Form a pullback diagram
$$ \xymatrix{ \calN \ar[r] \ar[d]^{p} & \calM^{\otimes} \ar[d] \\
\Nerve(\cDelta)^{op} \ar[r]^{A} & \calC^{\otimes}. }$$
Let $\calD$ denote the full subcategory of $\bHom_{ \Nerve(\cDelta)^{op} }( \Nerve(\calI), \calN)$
spanned by those functors $f: \Nerve(\calI) \rightarrow \calN$ with the following properties:
\begin{itemize}
\item[$(i)$] The restriction $f| \Nerve(\calI_0)$ belongs to $\Mod_{A}(\calM)$.
\item[$(ii)$] The restriction $f| \Nerve(\calI_1)$ belongs to $\Mod^{\qunit}_{A}(\calM)$.
\item[$(iii)$] For each $n \geq 0$, the induced map
$f( 0, [n]) \rightarrow f(1, [n])$ is an equivalence in $\calN_{[n]} \simeq \calM$.
\end{itemize}

We observe that, if $(iii)$ is satisfied, then $(i)$ and $(ii)$ are equivalent. For each $n \geq 0$, the
category $\calI_0 \times_{ \calI } \calI_{/ (1, [n])}$ has a final object, given by the map
$(0, [n]) \rightarrow (1, [n])$. It follows that a functor $f \in \bHom_{ \Nerve(\cDelta)^{op} }( \Nerve(\calI), \calN)$ satisfies $(iii)$ if and only if $f$ is a $p$-left Kan extension of
$f| \Nerve(\calI_0)$. Invoking Proposition \toposref{lklk}, we deduce that
the restriction map $r_0: \calD \rightarrow \Mod_{A}(\calM)$ is a trivial Kan fibration onto its essential image. Moreover, $r_0$ has a section $s$, given by composition with the projection $\calI \rightarrow \cDelta^{op}$. Consequently, $r_0$ is a trivial Kan fibration, so $s$ is a categorical equivalence. The restriction map $\theta$ factors as a composition
$ \Mod_{A}(\calM) \stackrel{s}{\rightarrow} \calD \stackrel{r_1}{\rightarrow} \Mod_{A}^{\qunit}(\calM).$
Consequently, it will suffice to show that $r_1$ is a categorical equivalence. In view of Proposition \toposref{lklk}, it will suffice to prove the following:

\begin{itemize}
\item[$(a)$] Every $f_1 \in \Mod^{\qunit}_{A}(\calM)$ admits a $p$-right Kan extension
$f \in \bHom_{ \Nerve(\cDelta)^{op} }( \Nerve(\calI), \calN)$.
\item[$(b)$] Let $f \in \bHom_{ \Nerve(\cDelta)^{op}}( \Nerve(\calI), \calN)$ be an arbitrary extension of $f_1 = f | \Nerve(\calI_1) \in \Mod^{\qunit}_{A}(\calM)$. Then $f$ is a $p$-right Kan extension of
$f_1$ if and only if $f$ satisfies condition $(iii)$.
\end{itemize}

In what follows, let us fix a quasi-unital module $f_1 \in \Mod_{A}^{\qunit}(\calM)$ and an
integer $n \geq 0$. Let $\calJ = \calI_1 \times_{ \calI} \calI_{(0, [n])/}$, and let
$g: \Nerve(\calJ) \rightarrow \calN$ be the composition of $f_1$ with the projection
$\calJ \rightarrow \Nerve(\calI_1)$. Let $v \in \calJ$ denote the object corresponding to the morphism $(0, [n]) \rightarrow (1, [n])$. In view of Lemma \toposref{kan2}, it will suffice to prove the following assertions (for every choice of $f_1$, $n \geq 0$):

\begin{itemize}
\item[$(a')$] There exists a $p$-limit diagram $\overline{g}: \Nerve(\calJ)^{\triangleright} \rightarrow \calN$ rendering the following diagram commutative:
$$ \xymatrix{ \Nerve(\calJ) \ar[r]^{g} \ar@{^{(}->}[d] & \calN \ar[d]^{p} \\
\Nerve(\calJ)^{\triangleright} \ar[r] \ar@{-->}[ur]^{\overline{g}} & \Nerve(\cDelta)^{op}. }$$

\item[$(b')$] Given an arbitrary map $\overline{g}$ which renders the above diagram commutative, $\overline{g}$ is a $p$-limit diagram if and only if $\overline{g}$ carries $\{v\}^{\triangleright} \subseteq \Nerve(\calJ)^{\triangleright}$ to an equivalence in $\calN_{[n]} \simeq \calM$.
\end{itemize}

Let $\calJ_0$ denote the full subcategory of $\calJ$ spanned by those maps
$\alpha: ( 0, [n]) \rightarrow (1, [m])$ for which the image of the induced map
$[m] \rightarrow [n]$ contains $n$. We claim that the inclusion
$\Nerve(\calJ_0)^{op} \subseteq \Nerve(\calJ)^{op}$ is cofinal. In view of
Theorem \toposref{hollowtt}, it will suffice to show that, for every object
$\alpha \in \calJ$, the category $\calJ_0 \times_{ \calJ } \calJ_{/ \alpha}$ has weakly contractible nerve. If $\alpha \in \calJ_0$, this is obvious. Suppose that $\alpha$ classifies a map $\gamma: [m] \rightarrow [n]$. Then $\calJ_0 \times_{ \calJ} \calJ_{/ \alpha}$ can be identified with
a product of categories $\{ \calE^{op}_{i} \}_{0 \leq i \leq n}$, where 
$$\calE_{i} \simeq \begin{cases} (\cDelta^{\nounit}_{+})_{\gamma^{-1}(i)/} & 
\text{if } i<n \\
\cDelta^{\nounit} & \text{if } i=n.\end{cases}$$ 
The categories $\calE_i$ have initial objects for $i < n$, and $\calE_n$ has weakly contractible nerve (because the inclusion $\Nerve(\cDelta^{\nounit})^{op} \subseteq \Nerve(\cDelta)^{op}$ is cofinal (Lemma \toposref{bball3}), cofinal maps are weak homotopy equivalences (Proposition \toposref{cofbasic}), and $\Nerve(\cDelta)^{op}$ is weakly contractible (Example \stableref{bin2} and Proposition \stableref{siftcont})), so $\prod_{0 \leq i \leq n} \Nerve(\calE_i)^{op}$ is likewise weakly contractible.

Let $g_0 = g | \Nerve(\calJ_0)$. In view of the cofinality statement above, $(a')$ and $(b')$ are equivalent to the following assertions:

\begin{itemize}
\item[$(a'')$] There exists a $p$-limit diagram $\overline{g}_0: \Nerve(\calJ_0)^{\triangleright} \rightarrow \calN$ rendering the following diagram commutative:
$$ \xymatrix{ \Nerve(\calJ_0) \ar[r]^{g_0} \ar@{^{(}->}[d] & \calN \ar[d]^{p} \\
\Nerve(\calJ_0)^{\triangleright} \ar[r] \ar@{-->}[ur]^{\overline{g}_0} & \Nerve(\cDelta)^{op}. }$$

\item[$(b'')$] Given an arbitrary map $\overline{g}_0$ which renders the above diagram commutative, $\overline{g}_0$ is a $p$-limit diagram if and only if $\overline{g}_0$ carries $\{v\}^{\triangleright} \subseteq \Nerve(\calJ_0)^{\triangleright}$ to an equivalence in $\calN_{[n]} \simeq \calM$.
\end{itemize}

We now observe that, for every morphism $\alpha: [m] \rightarrow [n]$ in $\cDelta$ for which
$\alpha(m) = n$ classifying a map $\Delta^1 \rightarrow \Nerve(\cDelta)^{op}$
the pullback $\calN \times_{ \Nerve(\cDelta)^{op} } \Delta^1$ is equivalent to a product
$\calM \times \Delta^1$. It follows that for every object $N \in \calN_{[m]}$, there exists
a locally $p$-Cartesian morphism $\overline{\alpha}: N' \rightarrow N$ in $\calN$ covering $\alpha$. Lemma \ref{smarties} implies that $\overline{\alpha}$ is $p$-Cartesian.

Let $h_1: \Nerve(\calJ_0) \rightarrow \Nerve(\cDelta)^{op}$ denote the
composition
$\Nerve(\calJ_0) \rightarrow \Nerve(\calI) \rightarrow \Nerve(\cDelta)^{op},$
so that we have a natural transformation $h: \Delta^1 \times \Nerve(\calJ_0) \rightarrow \Nerve(\cDelta)^{op}$ from $h_0 = h | \{0\} \times \Nerve(\calJ_0)$ to
$h_1$, where $h_0$ is the constant functor taking the value $[n]$. For each object
$x \in \Nerve(\calJ_0)^{\triangleright}$, the restriction of $h$ to $\Delta^1 \times \{x\}$ classifies a morphism $\alpha: [m] \rightarrow [n]$ satisfying $\alpha(m) = n$. It follows that we can lift $h$ to a $p$-Cartesian transformation $\widetilde{h}: \Delta^1 \times \Nerve(\calJ_0) \rightarrow \calN$
with $\widetilde{h} | \{1\} \times \Nerve(\calJ_0) = g_0$. Let $g'_0 = \widetilde{h}| \{0\} \times \Nerve(\calJ_0)$. Using Proposition \toposref{chocolatelast}, we obtain the following reformulations
of $(a'')$ and $(b'')$: 

\begin{itemize}
\item[$(a''')$] There exists a $p$-limit diagram $\overline{g}'_0: \Nerve(\calJ_0)^{\triangleright} \rightarrow \calN$ rendering the following diagram commutative:
$$ \xymatrix{ \Nerve(\calJ_0) \ar[r]^{g'_0} \ar@{^{(}->}[d] & \calN \ar[d]^{p} \\
\Nerve(\calJ_0)^{\triangleright} \ar[r]^{h_0} \ar@{-->}[ur]^{\overline{g}'_0} & \Nerve(\cDelta)^{op}. }$$

\item[$(b''')$] Given an arbitrary map $\overline{g}_0$ which renders the above diagram commutative, $\overline{g}'_0$ is a $p$-limit diagram if and only if $\overline{g}'_0$ carries $\{v\}^{\triangleright} \subseteq \Nerve(\calJ_0)^{\triangleright}$ to an equivalence in $\calN_{[n]} \simeq \calM$.
\end{itemize}

We now prove $(a''')$. Let us first regard $g'_0$ as a functor from
$\Nerve(\calJ_0)$ to $\calN_{[n]} \simeq \calM$. Let $M = f_1( [0]) \in \calM$. Unwinding the definitions, we see that the values assumed by $g'_0$ can be identified with $M$, and the morphisms between these values are given by iterated multiplication by the unit $1_{\calC} \rightarrow A$. Since $f_1$ is a quasi-unital module, it follows that $g'_0$ carries every morphism in $\calJ_0$ to an equivalence in $\calN_{[n]}$. The simplicial set $\Nerve(\calJ_0)$ is weakly contractible, since it is isomorphic to the product
$(\prod_{ 0 \leq i < n} \Nerve( \cDelta^{\nounit}_{+} )^{op}) \times \Nerve(\cDelta^{\nounit})^{op}.$
Applying Corollary \toposref{silt}, we deduce that $g'_0$ admits a colimit
$\overline{g}'_0: \Nerve(\calJ_0)^{\triangleright} \rightarrow \calN_{[n]}$, and that
$\overline{g}'_0$ carries $\{v\}^{\triangleright}$ to an equivalence in $\calN_{[n]}$. 
Since the $\infty$-category $\calM^{\otimes}_{[n]}$ is equivalent to the product
$\calC^{\otimes}_{[n]} \times \calM$, we conclude that $\overline{g}'_0$ is a 
$q_{[n]}$-limit diagram, where $q_{[n]}: \calM^{\otimes}_{[n]} \rightarrow \calC^{\otimes}_{[n]}$ denotes the production. Applying Corollary \toposref{pannaheave} to the diagram
$$ \xymatrix{ \calM^{\otimes} \ar[rr]^{q} \ar[dr] & & \calC^{\otimes} \ar[dl] \\
& \Nerve(\cDelta)^{op}, & }$$
we conclude that $\overline{g}'_0$ is a $q$-limit diagram. Since $p$ is a pullback of $q$, we 
deduce that $\overline{g}'_0$ is a $p$-limit diagram. This proves $(a''')$. 
Moreover, the uniqueness of $p$-limit diagrams implies the ``only if'' direction of $(b''')$.

It remains only to prove the ``if'' direction of $(b''')$. Let $\widetilde{g}'_0: \Nerve(\calJ_0)^{\triangleleft} \rightarrow \calN_{[n]}$ be an arbitrary extension of $g'_0$ such that
$\widetilde{g}'_0$ carries $\{v\}^{\triangleleft}$ to an equivalence. Since
$\overline{g}'_0$ is a limit of $g'_0$, there exists a natural transformation
$\gamma: \widetilde{g}'_0 \rightarrow \overline{g}'_0$ which is an equivalence when
restricted to $\Nerve(\calJ_0)$. Since both $\widetilde{g}'_0$ and
$\overline{g}'_0$ carry $\{v\}^{\triangleleft}$ to an equivalence in
$\calN_{[n]}$, we conclude that $\gamma$ is also an equivalence at the cone point
of $\Nerve(\calJ_0)^{\triangleleft}$. It follows that $\widetilde{g}'_0$ is equivalent to
$\overline{g}'_0$, so that $\widetilde{g}'_0$ is also a $p$-limit diagram as desired.
\end{proof}

\subsection{Limits and Colimits of Modules}\label{modcolim}

Let $\calC$ be a monoidal $\infty$-category, $\calM$ an $\infty$-category which is left-tensored over $\calC$, and $A$ an algebra object of $\calC$. Our goal in this section is to construct limits and colimits in the $\infty$-category $\Mod_{A}(\calM)$. 

We begin with the study of limits in $\Mod_{A}(\calM)$. In fact, we will work a little bit more generally, and consider {\em relative} limits with respect to the forgetful functor $\Mod(\calM) \rightarrow \Alg(\calC)$ (we refer the reader to \S \toposref{relacoim} for a discussion of relative limits in general).
Our basic result is Proposition \ref{eggur}, which assert that $\theta$ admits relative limits provided that the corresponding limits exist in $\calM$. We begin with a somewhat technical lemma, whose proof we defer until the end of this section.

\begin{lemma}\label{surtybove}
Let $p: \calC \rightarrow \calD$ be an inner fibration of $\infty$-categories, let $\calE$ and $K$ be
simplicial sets, and suppose given a diagram
$$ \xymatrix{ K \times \calE \ar@{^{(}->}[d] \ar[r]^{f} & \calC \ar[d]^{p} \\
K^{\triangleleft} \times \calE \ar@{-->}[ur]^{\overline{f}} \ar[r] & \calD. }$$
Suppose further that for each vertex $E$ of $\calE$, there exists an extension
$\overline{f}_{E}: K^{\triangleleft} \rightarrow \calC$ of $f_{E}$ which is compatible with the above diagram and is a $p$-limit. Then:
\begin{itemize}
\item[$(1)$] There exists a map $\overline{f}: K^{\triangleleft} \times \calE \rightarrow \calC$
rendering the above diagram commutative, with the property that for each vertex $E$ of $\calE$, the induced map $\overline{f}_{E}: K^{\triangleleft} \rightarrow \calC$ is a $p$-limit diagram.

\item[$(2)$] Let $\overline{f}: K^{\triangleleft} \times \calE \rightarrow \calC$
be an arbitrary map which renders the above diagram commutative. Then $\overline{f}$ satisfies the condition of $(1)$ if and only if the adjoint map $K^{\triangleleft} \rightarrow \Fun(\calE,\calC)$ is a $p^{\calE}$-limit diagram, where $p^{\calE}: \Fun(\calE, \calC) \rightarrow \Fun(\calE, \calD)$ is given by composition with $p$. 
\end{itemize}
\end{lemma}

\begin{proposition}\label{eggur}\index{module object!limit of}
Let $\calC$ be a monoidal $\infty$-category, and let $\calM$ be an $\infty$-category which is left-tensored over $\calC$. Let $K$ be a simplicial set such that $\calM$ admits $K$-indexed limits, and let $\theta: \Mod(\calM) \rightarrow \Alg(\calC)$ be the forgetful functor. Then:

\begin{itemize}
\item[$(1)$] For every diagram
$$ \xymatrix{ K \ar@{^{(}->}[d] \ar[r] & \Mod(\calM) \ar[d]^{\theta} \\
K^{\triangleleft} \ar[r] \ar@{-->}[ur] & \Alg(\calC) }$$
there exists a dotted arrow as indicated, which is a $\theta$-limit diagram.

\item[$(2)$] An arbitrary map $\overline{g}: K^{\triangleleft} \rightarrow \Mod(\calM)$ is a $\theta$-limit diagram if and only if the induced map $K^{\triangleleft} \rightarrow \calM$ is a limit diagram.
\end{itemize}
\end{proposition}

\begin{proof}
Let $n \geq 0$. We first observe that the equivalence $\calM^{\otimes}_{[n]} \simeq \calC^{\otimes}_{[n]} \times \calM^{\otimes}_{ \{n\} }$ implies the following:
\begin{itemize}
\item[$(1')$] For every diagram 
$$ \xymatrix{ K \ar@{^{(}->}[d] \ar[r] & \calM^{\otimes}_{[n]} \ar[d]^{q_{[n]}} \\
K^{\triangleleft} \ar[r] \ar@{-->}[ur] & \calC^{\otimes}_{[n]} }$$
there exists a dotted arrow as indicated, which is a $q_{[n]}$-limit diagram.
\item[$(2')$] An arbitrary diagram $K^{\triangleleft} \rightarrow \calM^{\otimes}_{[n]}$
is a $q_{[n]}$-limit diagram if and only the composition
$$K^{\triangleleft} \rightarrow \calM^{\otimes}_{[n]} \rightarrow \calM^{\otimes}_{ \{n\} }
\simeq \calM$$ is a limit diagram.
\end{itemize}
Combining this observation with Corollary \toposref{pannaheave}, we deduce:
\begin{itemize}
\item[$(1'')$] For every diagram 
$$ \xymatrix{ K \ar@{^{(}->}[d] \ar[r] & \calM^{\otimes}_{[n]} \ar[d]^{q_{[n]}} \\
K^{\triangleleft} \ar[r] \ar@{-->}[ur] & \calC^{\otimes}_{[n]} }$$
there exists a dotted arrow as indicated, which is a $q$-limit diagram.
\item[$(2'')$] An arbitrary diagram $K^{\triangleleft} \rightarrow \calM^{\otimes}_{[n]}$
is a $q$-limit diagram if and only the composition
$$K^{\triangleleft} \rightarrow \calM^{\otimes}_{[n]} \rightarrow \calM^{\otimes}_{ \{n\} }
\simeq \calM$$ is a limit diagram.
\end{itemize}

Let $\Mod'(\calM)$ be the full subcategory of $\bHom_{ \Nerve(\cDelta)^{op} }( \Nerve(\cDelta)^{op}, \calM^{\otimes})$ spanned by those objects which project to algebra objects of $\calC$.
Combining $(1'')$, $(2'')$, and Lemma \ref{surtybove}, we deduce:

\begin{itemize}
\item[$(1''')$] For every diagram
$$ \xymatrix{ K \ar@{^{(}->}[d] \ar[r] & \Mod'(\calM) \ar[d]^{\theta'} \\
K^{\triangleleft} \ar[r] \ar@{-->}[ur] & \Alg(\calC) }$$
there exists a dotted arrow as indicated, which is a $\theta'$-limit diagram.

\item[$(2''')$] An arbitrary map $\overline{p}: K^{\triangleleft} \rightarrow \Mod'(\calM)$ is a $\theta'$-limit diagram if and only if, for every $n \geq 0$, the composition
$K^{\triangleleft} \rightarrow \calM^{\otimes}_{[n]} \rightarrow
\calM^{\otimes}_{ \{n\} } \simeq \calM$ is a limit diagram.
\end{itemize}

To deduce $(1)$ from these assertions, it will suffice to show that if
if $\overline{g}: K^{\triangleright} \rightarrow \Mod'(\calM)$ satisfies $(2''')$, and $g = \overline{g}| K$ factors through
$\Mod(\calM)$, then $\overline{g}$ factors through $\Mod(\calM)$. Let $f: [m] \rightarrow [n]$ be a convex morphism in $\cDelta$ such that $f(m) = n$. Then $f$ induces a natural transformation
$\overline{g}_{[n]} \rightarrow \overline{g}_{[m]}$ of functors $K^{\triangleleft} \rightarrow \calM^{\otimes}$. We wish to show that this natural transformation is $(p \circ q)$-coCartesian.
Since $f(m) = n$, we have a homotopy commutative diagram
$$ \xymatrix{ \calM^{\otimes}_{[n]} \ar[rr] \ar[dr]^{\alpha} & & \calM^{\otimes}_{[m]} \ar[dl]^{\beta} \\
& \calM, & }$$
and it suffices to show that the associated transformation $\overline{t}: \alpha \circ \overline{g}_{[n]} 
\rightarrow \beta \circ \overline{g}_{[m]}$ is an equivalence. Our hypothesis implies
that $\overline{t}$ restricts to an equivalence $t: \alpha \circ g_{[n]} \rightarrow \beta \circ g_{[m]}$.
Since $\overline{g}$ satisfies $(2'')$, the maps $\alpha \circ \overline{g}_{[n]}$ and
$\beta \circ \overline{g}_{[m]}$ are both limit diagrams in $\calM$. It follows that
$\overline{t}$ is an equivalence as well, as desired.

We now complete the proof by observing that if $\overline{g}: K^{\triangleleft} \rightarrow \Mod'(\calM)$ factors through $\Mod(\calM)$, then the criteria of $(2)$ and $(2''')$ are equivalent.
\end{proof}

\begin{corollary}\label{thetacart}
Let $\calC$ be a monoidal $\infty$-category, $\calM$ an $\infty$-category which is left-tensored over $\calC$, and $\theta: \Mod(\calM) \rightarrow \Alg(\calC)$ the forgetful functor. Then $\theta$ is a Cartesian fibration. Moreover, a morphism $f$ in $\Mod(\calM)$ is $\theta$-Cartesian if and only if the image of $f$ in $\calM$ is an equivalence.
\end{corollary}

\begin{proof}
Apply Proposition \ref{eggur} in the case $K = \Delta^0$.
\end{proof}

\begin{remark}\label{thetacart2}
Proposition \ref{eggur} has an analogue for nonunital modules, with exactly the same proof.
In particular, we obtain the following nonunital version of Corollary \ref{thetacart}:
\begin{itemize}
\item[$(\ast)$] The forgetful functor $\theta: \Mod^{\nounit}(\calM) \rightarrow \Alg^{\nounit}(\calC)$ 
is a Cartesian fibration. Moreover, a morphism $f: M \rightarrow N$ of nonunital module objects is
$\theta$-Cartesian if and only if its image in $\calM$ is an equivalence.
\end{itemize}
\end{remark}

\begin{corollary}\label{goop}
Let $\calC$ be a monoidal $\infty$-category, $\calM$ an $\infty$-category which is left-tensored over $\calC$, and $\theta: \Mod(\calM) \rightarrow \Alg(\calC)$ the forgetful functor. Let $A$ be an algebra object of $\calC$. Suppose that $\calM$ admits $K$-indexed limits. Then:
\begin{itemize}
\item[$(1)$] The $\infty$-category $\Mod_{A}(\calM)$ admits $K$-indexed limits.
\item[$(2)$] A map $\overline{p}: K^{\triangleleft} \rightarrow \Mod_{A}(\calM)$ is a limit diagram if and only if the induced map $K^{\triangleleft} \rightarrow \calM$ is a limit diagram.
\item[$(3)$] Given a morphism $\phi: B \rightarrow A$ of algebra objects of $\calA$, the induced functor $\Mod_{A}(\calM) \rightarrow \Mod_{B}(\calM)$ preserves $K$-indexed limits.
\end{itemize}
\end{corollary}

We now turn to the problem of constructing colimits in $\infty$-categories of modules.
We begin with the following very general principle:

\begin{proposition}\label{poststorkus}
Let $\calA \subseteq \widehat{\Cat}_{\infty}$ be a subcategory of the $\infty$-category of
$($not necessarily small$)$ $\infty$-categories. Assume that $\calA$ has the following properties:
\begin{itemize}
\item[$(a)$] The $\infty$-category $\calA$ admits small limits, and the inclusion
$\calA \subseteq \widehat{\Cat}_{\infty}$ preserves small limits.
\item[$(b)$] If $X$ belongs to $\calA$, then $\Fun(\Delta^1, X)$ belongs to $\calA$.
\item[$(c)$] If $X$ and $Y$ belong to $\calC$, then a functor $X \rightarrow \Fun(\Delta^1,Y)$ is a morphism in $\calA$ if and only if, for every vertex $v$ of $\Delta^1$, the composite functor
$X \rightarrow \Fun(\Delta^1,Y) \rightarrow \Fun( \{v\}, Y) \simeq Y$ is a morphism of $\calA$.
\end{itemize}
Let $\calC$ be a monoidal $\infty$-category, $A$ an algebra object of $\calC$, $\calM$ an $\infty$-category which is left tensored over $\calC$. Suppose $\calM$ is an object of $\calA$, and that the functor $A \otimes \bigdot: \calM \rightarrow \calM$ is a morphism of $\calA$. Then:
\begin{itemize}
\item[$(1)$] The $\infty$-category $\Mod_{A}(\calM)$ is an object of $\calA$.
\item[$(2)$] For every $\infty$-category $\calN$ belonging to $\calA$, a functor
$\calN \rightarrow \Mod_{A}(\calM)$ is a morphism in $\calA$ if and only if the composite functor
$\calN \rightarrow \Mod_{A}(\calM) \rightarrow \calM$ is a morphism in $\calA$.
\end{itemize}
In particular, the forgetful functor $\Mod_{A}(\calM) \rightarrow \calM$ is a morphism in
$\calA$.
\end{proposition}

\begin{proof}
Let $p: \calM^{\otimes} \rightarrow \calC^{\otimes}$ exhibit $\calM$ as left-tensored over $\calC$. Form a pullback diagram
$$ \xymatrix{ X \ar[d]^{p'} \ar[r] & \calM^{\otimes} \ar[d]^{p} \\
\Nerve(\cDelta)^{op} \ar[r]^{A} & \calC^{\otimes}. }$$
We observe that $p'$ is a locally coCartesian fibration (Lemma \ref{excel}), each fiber of $p'$ is equivalent to $\calM$, and each of the
associated functors can be identified with an iterate of the functor $A \otimes \bigdot: \calM \rightarrow \calM$. Assertion $(1)$ is now an immediate consequence of Proposition \toposref{prestorkus}.

Now suppose that $f: \calN \rightarrow \Mod_{A}(\calM)$ is as in $(2)$. Proposition \toposref{prestorkus} implies that $f$ is a morphism of $\calA$ if and only if, for every $n \geq 0$, the composite map $\calN \rightarrow \Mod_{A}(\calM) \rightarrow X_{[n]}$ belongs to $\calA$.
We complete the proof by observing that each of the functors $\Mod_{A}(\calM) \rightarrow X_{[n]}$ is equivalent to the forgetful functor $\Mod_{A}(\calM) \rightarrow \calM$.
\end{proof}

\begin{corollary}\label{gloop}\index{module object!colimit of}
Let $\calC$ be a monoidal $\infty$-category, let $\calM$ be an $\infty$-category which is left-tensored over $\calC$. Let $K$ be a simplicial set such that $\calM$ admits $K$-indexed colimits.
Let $A \in \Alg(\calC)$, and suppose that the functor $A \otimes \bigdot: \calM \rightarrow \calM$ preserves $K$-indexed colimits. Then:
\begin{itemize}
\item[$(1)$] Every diagram $f: K \rightarrow \Mod_{A}(\calM)$ has a colimit.
\item[$(2)$] An arbitrary diagram $\overline{f}: K^{\triangleright} \rightarrow \Mod_{A}(\calM)$ is a colimit diagram if and only if it induces a colimit diagram $K^{\triangleright} \rightarrow \calM$.
\end{itemize}
\end{corollary}

\begin{proof}
Let $\calA \subseteq \widehat{\Cat}_{\infty}$ be the subcategory whose objects are $\infty$-categories which admit $K$-indexed colimits, and whose morphisms are functors which preserve $K$-indexed colimits. Applying Proposition \ref{poststorkus}, we conclude that
$\Mod_{A}(\calM)$ belongs to $\calA$, and that the forgetful functor
$\Mod_{A}(\calM) \rightarrow \calM$ is a morphism of $\calA$. This proves $(1)$, and the
``only if'' direction of $(2)$. The ``if'' direction then follows formally, since the forgetful functor
$\Mod_{A}(\calM) \rightarrow \calM$ is conservative.

\end{proof}

\begin{corollary}\label{underwhear}
Let $\calC$ be an $\infty$-category equipped with a monoidal structure and let $\calM$ be left-tensored over $\calC$. Suppose that $\calM$ is presentable and that, for each
$C \in \calC$, the functor $C \otimes \bigdot: \calM \rightarrow \calM$ preserves small colimits. Then:
\begin{itemize}
\item[$(1)$] For every $A \in \Alg(\calC)$, the $\infty$-category $\Mod_{A}(\calM)$ is presentable.
\item[$(2)$] For every morphism $A \rightarrow B$ of algebra objects of $\calC$, the associated functor $\Mod_{B}(\calM) \rightarrow \Mod_{A}(\calM)$ preserves small limits and colimits.
\item[$(3)$] The forgetful functor $\theta: \Mod(\calM) \rightarrow \Alg(\calC)$ is a presentable fibration $($Definition \toposref{urtus}$)$. 
\end{itemize}
\end{corollary}

\begin{proof}
Assertion $(1)$ is a special case of Proposition \ref{poststorkus}.
Corollary \ref{thetacart} implies that the diagram
$$ \xymatrix{ \Mod_{B}(\calM) \ar[dr] \ar[rr] & & \Mod_{A}(\calM) \ar[dl] \\
& \calM & }$$
commutes up to homotopy. Assertion $(2)$ follows immediately from Corollaries \ref{goop} and \ref{gloop}. Assertion $(3)$ follows from $(1)$, Corollaries \ref{thetacart} and \ref{goop}, and Proposition \toposref{surtog}.
\end{proof}

\begin{remark}\label{sillybear}
Under the hypotheses of Corollary \ref{underwhear}, if $A \rightarrow B$ is a morphism of
algebra objects of $\calC$. then the forgetful functor $\psi: \Mod_{B}(\calM) \rightarrow \Mod_{A}(\calM)$ admits both left and right adjoints (Corollary \toposref{adjointfunctor}). In \S \ref{balpair} we will prove the existence of a left adjoint to $\psi$ under much weaker assumptions (Lemma \ref{peacestick}).
\end{remark}

We conclude this section by giving a proof of Lemma \ref{surtybove}.

\begin{lemma}\label{niki}
Let $p: \calD \rightarrow \calC$ be a right fibration of $\infty$-categories, let
$\calC^{0} \subseteq \calC$ be a full subcategory, and let $\calD^{0} = \calC^{0} \times_{ \calC} \calD$.
Let $q: X \rightarrow S$ be an inner fibration of simplicial sets, and let
$F: \calC \rightarrow X$ be a map which is a $q$-left Kan extension of $F | \calC^{0}$. 
Then $F \circ p$ is a $q$-left Kan extension of $F \circ p | \calD^{0}$.
\end{lemma}

\begin{proof}
Let $D$ be an object of $\calD$, $C = p(D)$, and define
$$\calC^{0}_{/C} = \calC^{0} \times_{\calC} \calC_{/C} \quad
\calD^{0}_{/D} = \calD^{0} \times_{\calD} \calD_{/D}.$$
We wish to show that the composition
$$ (\calD^{0}_{/D})^{\triangleright} 
( \calC^{0}_{/C})^{\triangleright} \rightarrow \calC \stackrel{F}{\rightarrow} X$$
is a $q$-colimit diagram. Since $F$ is a $q$-left Kan extension of $F | \calC^{0}$, it will suffice to show that the map $\phi_0: \calD^{0}_{/D} \rightarrow \calC^{0}_{/C}$ is a trivial Kan fibration. The map
$\phi_0$ is a pullback of the map $\phi: \calD_{/D} \rightarrow \calC_{/C}$, which is a trivial Kan fibration by Proposition \toposref{sharpen2}.
\end{proof}

\begin{lemma}\label{kanexp}
Let $p: X \rightarrow S$ and
$q: Y \rightarrow Z$ be maps of simplicial sets. Assume that $q$ is a categorical fibration, and that
$p$ is either a Cartesian fibration or a coCartesian fibration.

Define new simplicial sets
$Y'$ and $Z'$ equipped with maps $Y' \rightarrow S$, $Z' \rightarrow S$ via the formulas
$$ \Hom_{S}(K, Y') \simeq \Hom( X \times_{S} K, Y)$$
$$ \Hom_{S}(K,Z') \simeq \Hom( X \times_{S} K, Z).$$
Let $\calC'$ be an $\infty$-category equipped with a functor
$f: \calC' \rightarrow Y'$, and let $\calC$ be a full subcategory of $\calC'$. 

Then:
\begin{itemize}
\item[$(1)$] Composition with $q$ determines a categorical fibration
$q': Y' \rightarrow Z'$.
\item[$(2)$] Let $F: X \times_{S} \calC' \rightarrow Y$ be the map classified by $f$, and suppose that $F$ is a $q$-left Kan extension of $F | X \times_{S} \calC$. Then
$f$ is a $q'$-left Kan extension of $f | \calC$.
\end{itemize}
\end{lemma}

\begin{proof}
We first prove $(1)$. We wish to show that $q'$ has the right lifting property with respect to every
inclusion $i: A \rightarrow B$ of simplicial sets which is a categorical equivalence. For this, it suffices to show that $q$ has the right lifting property with respect to every inclusion of the form
$i': X \times_{S} A \rightarrow X \times_{S} B$. Since $q$ is a categorical fibration, it suffices to prove that $i'$ is a categorical equivalence. This follows from Proposition \toposref{basechangefunky}.

Let $C$ be an object of $\calC'$, and let $\calC_{/C}$ denote the fiber product
$\calC \times_{ \calC'} \calC'_{/C}$. We wish to show that the composition
$$ \calC_{/C}^{\triangleright} \rightarrow \calC' \stackrel{f}{\rightarrow} Y'$$
is a $q'$-colimit diagram. Replacing $\calC \subseteq \calC'$ by the inclusion
$\calC_{/C} \subseteq \calC^{\triangleright}_{/C}$ (and applying Lemma \ref{niki}), we can reduce to the case $\calC' = \calC^{\triangleright}$. Replacing $\calC$ by a minimal model if necessary (Proposition \toposref{minimod}), we may assume that $\calC$ is minimal. Similarly, we may assume that $p$ is minimal as an inner fibration.

Let $n > 0$, and suppose we are given a diagram
$$ \xymatrix{ \calC \star \bd \Delta^n \ar[r]^{f'} \ar@{^{(}->} & Y' \ar[d]^{q'} \\
\calC \star \Delta^n \ar[r]^{g} \ar@{-->}[ur] & Z', }$$
where $f' | \calC \star \{ 0\}$ coincides with $f$. We wish to show that there exists a dotted arrow, as indicated in the diagram. Composing $g$ with the map $Z' \rightarrow S$, we obtain a map
$\calC \star \Delta^n \rightarrow S$. Let $\calD$ denote the fiber product
$X \times_{S} ( \calC \star \Delta^n)$, and let $\calD^{0} = X \times_{S} ( \calC \star \bd \Delta^n)$. 
Unwinding the definitions, we are reduced to solving a lifting problem depicted in the diagram
$$ \xymatrix{ \calD^{0} \ar@{^{(}->}[d] \ar[r]^{F'} & Y \ar[d]^{q} \\
\calD \ar[r]^{G} & Z. }$$

Let $\calD'$ denote the inverse image of $\Delta^n$ in $\calD$. For every simplicial subset
$W \subseteq \calD'$, we let $\calD(W)$ denote the collection of all simplices of
$\calD$ whose intersection with $\calD'$ belongs to $W$. Let $P$ denote the partially ordered
collection of pairs $(W, F_{W})$, where $W \subseteq \calD'$ contains the inverse image of
$\bd \Delta^n$, and $F_{W}: \calD(W) \rightarrow Y$ satisfies $F_{W} | \calD^{0} = F'$, and
$F_W$ fits into a commutative diagram
$$ \xymatrix{ \calD(W) \ar@{^{(}->}[d] \ar[r]^{F_W} & Y \ar[d]^{q} \\
\calD \ar[r]^{G} & Z. }$$

Applying Zorn's lemma, we deduce that $P$ has a maximal element. To complete the proof, it will suffice to show that $W = \calD'$. Assume otherwise, and let $\sigma$ be a simplex of $\calD'$ of minimal dimension which does not belong to $W$. We note that $\sigma$ surjects onto $\Delta^n$, so the dimension of $\sigma$ is necessarily positive. Let $\calE$ denote the fiber product
$\calD_{/ \sigma} \times_{ \calC \star \Delta^n} \calC$. Using the minimality of $p$, 
$\calC \star \Delta^n$, and Proposition \toposref{minstrict}, we deduce the existence of a pushout diagram
$$ \xymatrix{ \calE \star \bd \sigma \ar[r] \ar@{^{(}->}[d] & \calD(W) \ar[d] \\
\calE \star \sigma \ar[r] & \calD(W'). }$$ 
In view of the maximality of $(W,F_W)$, it follows that $F_{W}$ does not admit an
extension to $\calD(W')$. In other words, there is no solution to the associated lifting problem
$$ \xymatrix{ \calE \star \bd \sigma \ar@{^{(}->}[d] \ar[r]^{h} & Y \ar[d]^{q} \\
\calE \star \sigma \ar[r] \ar@{-->}[ur] & Z. }$$
Let $h'$ denote the restriction of $h$ to $\calE \star \{v\}$, where $v$ is the initial vertex of
$\sigma$. Then $h'$ is not a $q$-colimit diagram. However, $h'$ can be written as a composition
$$ \calE \star \{v \} \rightarrow (\calD_{/v} \times_{\calC \star \Delta^n} \calC)^{\triangleright}
\stackrel{h''}{\rightarrow} Y,$$
where $h''$ is a $q$-colimit diagram. The map
$\calE \rightarrow \calD_{/v} \times_{ \calC \star \Delta^n} \calC$ is a pullback of the
trivial fibration $\calD_{/\sigma} \rightarrow \calD_{/v}$, and therefore cofinal; we therefore obtain a contradiction and the proof is complete.
\end{proof}

\begin{proof}[Proof of Lemma \ref{surtybove}]
Without loss of generality, we may suppose that $\calE$ and $K$ are $\infty$-categories.
Let $\infty$ denote the cone point of $K^{\triangleleft}$. For each object $E \in \calE$, the inclusion
$K \times \{ \id_E \} \subseteq K \times \calE_{E/} \simeq (K \times \calE)_{(\infty, E)/}$ is left anodyne. Consequently, $\overline{f}$ satisfies $(1)$ if and only if $\overline{f}$ is a $p$-right Kan extension of $f'$. The existence of $\overline{f}$ follows from Lemma \toposref{kan2}.

The ``only if'' direction of $(2)$ follows immediately from Lemma \ref{kanexp}. The converse follows from the uniqueness of limits (up to equivalence).
\end{proof}

\subsection{Free Modules}\label{freetea}

Let $\calC$ be a monoidal $\infty$-category, and let $A$ be an algebra object of $\calC$.
Our goal in this section is to construct a left adjoint to the forgetful functor $\Mod_{A}(\calC) \rightarrow \calC$. In fact, we will treat a slightly more general problem, by allowing modules taking values in an arbitrary $\infty$-category left tensored over $\calC$.

\begin{definition}\label{squirch}\index{free!module object}\index{module object!free}
Let $\calC$ be a monoidal $\infty$-category, let $\calM$ be an $\infty$-category which is left-tensored over $\calC$. Let $A \in \Alg(\calC)$, let $\theta: \Mod_{A}(\calM) \rightarrow \calM$ denote the forgetful functor, and let $M_0 \in \calM$. A {\it free $A$-module generated by $M_0$} is an object $M \in \Mod_{A}(\calM)$ together with a map $\alpha: M_0 \rightarrow \theta(M)$ with the following universal property: for every $N \in \Mod_{A}(\calM)$, composition with $\alpha$ induces a homotopy equivalence
$\bHom_{\Mod_{A}(\calM)}(M,N) \rightarrow \bHom_{\calM}( M_0, \theta(N) ).$
\end{definition}

The goal of this section is to prove the following result:

\begin{proposition}\label{pretara}
Let $\calC$ be a monoidal $\infty$-category, and let $\calM$ be an $\infty$-category which is left-tensored over $\calC$. Let $A \in \Alg(\calC)$. Then:
\begin{itemize}
\item[$(1)$] The forgetful functor $\theta: \Mod_{A}(\calM) \rightarrow \calM$ admits a left adjoint.
\item[$(2)$] For every $M_0 \in \calM$, there exists a free $A$-module $M$ generated by $M_0$.
\item[$(3)$] An arbitrary map $M_0 \rightarrow \theta(N)$ exhibits $N$ as a free $A$-module generated by $M_0$ if and only if the composition 
$$A \otimes M_0 \rightarrow A \otimes \theta(N) \rightarrow \theta(N)$$ is an equivalence in $\calM$.
\end{itemize}
\end{proposition}

\begin{corollary}\label{ingurtime}
Let $\calC$ be a monoidal $\infty$-category, $1_{\calC}$ the unit object of $\calC$, $A$ an algebra object of $\calC$, $\theta: \Mod_{A}(\calC) \rightarrow \calC$ the forgetful functor, and $M \in \Mod_{A}(\calC)$. Then composition with the unit map of $A$ induces a homotopy equivalence
$\bHom_{ \Mod_{A}(\calC) }(A, M) \rightarrow \bHom_{\calC}( 1_{\calC}, \theta(M) ).$
Here we identify $A$ with the corresponding left $A$-module $($Example \ref{algitself}$)$. 
\end{corollary}

\begin{proof}
Proposition \ref{pretara} implies that $A$ is freely generated by the unit map $1_{\calC} \rightarrow A$ as a left $A$-module.
\end{proof}

\begin{corollary}\label{puterry}
Let $\calC$ be an $\infty$-category equipped with a monoidal structure, and let $\calM$ be an $\infty$-category which is left-tensored over $\calC$. Let $A$ be an algebra object of $\Alg(\calC)$ such that the unit map $1_{\calC} \rightarrow A$ is an equivalence in $\calC$. Then the forgetful functor $\theta: \Mod_{A}(\calM) \rightarrow \calM$ is an equivalence of $\infty$-categories.
\end{corollary}

\begin{proof}
Let $F$ be the left adjoint to $\theta$ supplied by Proposition \ref{pretara}, and let
$$ u: \id_{\calM} \rightarrow \theta \circ F, \quad  v: F \circ \theta \rightarrow \id_{\calM}$$
be a compatible unit and counit for the adjunction. We wish to prove that $u$ and $v$
are equivalences of functors.

We first consider the functor $u$. Proposition \ref{pretara} implies that the composition
$\theta \circ F$ can be identified with the functor $M \mapsto A \otimes M$. The unit map
$u$ is given by tensor product with the unit map $u_0: 1_{\calC} \rightarrow A$ of the algebra $A$.
By hypothesis, $u_0$ is an equivalence in $\calC$, so that $u$ is an equivalence in
$\Fun( \calM, \calM)$. 

We now show that $v$ is an equivalence. Since the functor $\theta$ detects equivalences (Corollary \ref{thetacart}), it will suffice to show that the induced transformation
$\alpha: \theta \circ F \circ \theta \rightarrow \theta$ is an equivalence of functors. We now observe that $u$ provides a right inverse to this $\alpha$. Since $u$ is an equivalence, we conclude also that $\alpha$ is an equivalence.
\end{proof}

\begin{proof}[Proof of Proposition \ref{pretara}]
Let $\calM^{\otimes} \stackrel{p}{\rightarrow} \calC^{\otimes} \stackrel{q}{\rightarrow} \Nerve(\cDelta)^{op}$ realize $\calM$ as left-tensored over $\calC$.
Let $\calI$ be the subcategory of $[1] \times \cDelta^{op}$ defined as follows:
\begin{itemize}
\item[$(a)$] Every object of $[1] \times \cDelta^{op}$ belongs to $\calI$.
\item[$(b)$] Given $0 \leq i \leq j \leq 1$ and a map $\alpha: [m] \rightarrow [n]$ in $\cDelta$, the corresponding morphism $(i, [n]) \rightarrow (j, [m])$ belongs to $\calI$ if and only if
either $j=1$ or $\alpha(m) = n$.
\end{itemize}
We let $\calI_0 \subseteq \cDelta^{op}$ and $\calI_1 \simeq \cDelta^{op}$ denote the inverse images in $\cal$ of the objects $0,1 \in [1]$. We regard $\Nerve(\calI)$ as an object of $(\sSet)_{/ \calC^{\otimes}}$ via the composition
$\Nerve(\calI) \subseteq \Nerve([1] \times \cDelta^{op}) \rightarrow \Nerve(\cDelta)^{op}
\stackrel{A}{\rightarrow} \calC^{\otimes}.$ Let $\calN$
denote the fiber product $ \Nerve(\calI) \times_{ \calC^{\otimes} } \calM^{\otimes}$, and let
$p': \calN \rightarrow \Nerve(\calI)$ be the projection. Lemma \ref{excel} implies that $p'$ is a locally coCartesian fibration. Let $\calD$ denote the full subcategory of $\bHom_{\Nerve(\calI)}( \Nerve(\calI), \calN)$
spanned by those functors $f: \Nerve(\calI) \rightarrow \calM^{\otimes}$ with the following properties:
\begin{itemize}
\item[$(i)$] The restriction $f| \Nerve(\calI_0)$ carries each morphism in $\calI_0$ to a
$p'$-coCartesian edge of $\calN$.
\item[$(ii)$] The restriction $f| \Nerve(\calI_1)$ is a module object of $\calM^{\otimes}$.
\end{itemize}
Let $\calD_0$ be the full subcategory of $\calD$ spanned by those functors $f$ which satisfy the following additional condition:
\begin{itemize}
\item[$(iii_0)$] The functor $f$ is a $p'$-left Kan extension of $f | \Nerve(\calI_0)$.
\end{itemize}
Let $\calD_1$ be the full subcategory of $\calD$ spanned by those functors which satisfy the following additional condition:
\begin{itemize}
\item[$(iii_1)$] The functor $f$ is a $p'$-right Kan extension of $f| \Nerve(\calI_1)$. 
\end{itemize}
We observe that every for every $n \geq 0$, the $\infty$-category $\Nerve(\calI_1)_{(0, [n])/}$ has a final object $(1, [n])$, and the $\infty$-category $\Nerve(\calI_0)_{/(1,[n])}$ has an initial object
$(0, [n] \star [0] )$. Using Lemma \toposref{kan2}, we deduce:
\begin{itemize}
\item[$(\ast)$] Every functor $f_0 \in \bHom_{ \Nerve(\calI)}( \Nerve(\calI_0), \calN)$ admits a $p'$-left Kan extension $f \in \bHom_{ \Nerve(\calI)}( \Nerve(\calI), \calN)$. Moreover, an arbitrary extension $f'$ is a $p'$-left Kan extension if and only if $p'$ carries
each morphism $( 0, [n] \star [0] ) \rightarrow (1, [n])$ to a locally $p'$-coCartesian edge of $\calN$ (any such edge is automatically $p'$-coCartesian, in view of Lemma \ref{smarties}).

\item[$(\ast')$] Every functor $f_1 \in \bHom_{ \Nerve(\calI) }( \Nerve(\calI_1), \calN)$ admits a $p'$-right Kan extension $f \in \bHom_{ \Nerve(\calI) }( \Nerve(\calI), \calN)$. Moreover, an arbitrary extension $f'$ is a $p'$-right Kan extension if and only if $q$ carries
each morphism $( 0, [n] ) \rightarrow (1, [n])$ to a $p'$-Cartesian edge of $\calN$ (these are the morphisms which project to equivalences in $\calM^{\otimes}$).
\end{itemize} 

Using $(\ast)$, we deduce that $(i)$ and $(iii_0)$ imply $(ii)$. Similarly, using $(\ast')$, we deduce that $(ii)$ and $(iii_1)$ imply $(i)$. Let $\calD'_0$ be the full subcategory of
$\bHom_{ \Nerve(\calI)}( \Nerve(\calI_0), \calN)$ spanned by the coCartesian sections, and let
$\calD'_1$ be the full subcategory of $\bHom_{ \Nerve(\calI)}( \Nerve(\calI_1), \calN)$ spanned by those functors which classify $A$-modules. Proposition \toposref{lklk} implies that the restriction maps
$$ \phi_0: \calD_0 \rightarrow \calD'_0, \quad \phi_1: \calD_1 \rightarrow \calD'_1$$
are trivial Kan fibrations. The functor $\phi_1$ has a canonical section $s$, given by composition with the projection $\calI \rightarrow \cDelta^{op} \simeq \calI_1$.

Let $\calN_0 = \calN \times_{ \Nerve(\calI)} \Nerve(\calI_0)$. Using Proposition \toposref{gotta}, we deduce that the projection $p'_0: \calN_0 \rightarrow \Nerve(\calI_0)$ is a coCartesian fibration. Each fiber of $p'_0$ is equivalent to $\calM$, and each of the associated functors between the fibers is an equivalence. Since $\Nerve(\calI_0)$ is weakly contractible (the category $\calI_0$ has a zero object $(0, [0])$), we deduce that $\calN_0$ is equivalent to a product $\Nerve(\calI_0) \times \calM$. It follows also that the $\infty$-category $\calD'_0$ is equivalent, via evaluation at
the zero object of $\calI_0$, to $\calM$. 

The forgetful functor $\Mod_{A}(\calM) \rightarrow \calM$ is equivalent to the composition
$$ \theta: \Mod_{A}(\calM) \simeq \calD'_1 \stackrel{s}{\rightarrow} \calD_1 \rightarrow
\calD'_0 \rightarrow \calM,$$
where the final pair of maps are given by restriction and evaluation at the zero object of $\calI_0$.
Using Proposition \toposref{leftkanadj}, we deduce that $\theta$ has a left adjoint, given by the composition
$$ \calM \stackrel{s'}{\rightarrow} \calD'_0 \stackrel{s''}{\rightarrow} \calD_0
\subseteq \calD \rightarrow \calD'_1 \simeq \Mod_{A}(\calM)$$
where the penultimate map is given by restriction, and the maps $s'$ and $s''$ are sections of the trivial fibrations $\calD_0 \rightarrow \calD'_0 \rightarrow \calM$. This proves $(1)$. 

The equivalence $(1) \Leftrightarrow (2)$ follows from Proposition \toposref{simpex}. 
More precisely, given $M_0 \in \calM$, a free $A$-module $M$ generated by $M_0$ can be constructed as follows. Choose a morphism $\alpha: (s'' \circ s')(M_0) \rightarrow D$ in $\calD$, where $D \in \calD_1$. Then $D$ restricts to an object $M \in \Mod_{A}(\calM)$, and the natural transformation $\alpha$ determines a map 
$$\psi: M_0 = ((s'' \circ s')(M_0))(0, [0]) \rightarrow ((s'' \circ s')(M_0))(1, [0])
\simeq D(1, [0]) = \theta(M).$$
Since $(s'' \circ s')(M_0)$ is a $q$-left Kan extension of $D|\calI_0$, $(\ast)$ implies that
the edge $$((s'' \circ s')(M_0))(0, [1]) \rightarrow ((s'' \circ s')(M_0))(1,[0])$$ is $q$-coCartesian.
It follows that the canonical map $A \otimes M_0 \rightarrow \theta(M)$ is an equivalence.
This proves the ``only if'' direction of $(3)$, since a free $A$-module generated by $M_0$ is uniquely determined up to equivalence. 

We now prove the ``if'' direction of $(3)$. Suppose given an $A$-module $N$ and a map
$\psi': M_0 \rightarrow \theta(N)$ which induces an equivalence $A \otimes M_0 \simeq N$.
Let $M \in \Mod_{A}(\calM)$ be defined as above. Then $\psi'$ is equivalent to a composition
$M_0 \stackrel{\psi}{\rightarrow} \theta(M) \stackrel{\theta(g)}{\rightarrow} \theta(N)$
for some map $g: M \rightarrow N$. We have a commutative diagram
$$ \xymatrix{ & \theta(M) \ar[dr]^{\theta(g)} & \\
A \otimes M_0 \ar[ur] \ar[rr] & & \theta(N). }$$
It follows from the two-out-of-three property that $\theta(g)$ is an equivalence in $\calM$.
Corollary \ref{goop} now implies that $g$ is an equivalence, so that $\psi'$ exhibits $N$ as an $A$-module freely generated by $M_0$.
\end{proof}

\subsection{Modules in a Monoidal Model Category}\label{hutman}

Let $\bfA$ be a simplicial model category, equipped with a compatible monoidal structure. In \S \ref{monoidate}, we saw that the underlying $\infty$-category $\Nerve(\bfA^{\degree})$ inherits the structure of a monoidal $\infty$-category. Moreover, in good cases we can show that every
algebra object $A_0 \in \Alg( \Nerve(\bfA^{\degree}))$ arises (up to equivalence) from a strictly associative algebra object $A \in \Alg(\bfA)$ (Theorem \ref{beckify}). Our goal in this section is to compare the $\infty$-categorical theory of $A_0$-modules in $\Nerve(\bfA^{\degree})$ with the classical theory of $A$-modules in $\bfA$. 

We begin by observing that if $\bfA$ is a monoidal category equipped with a compatible simplicial structure, and $A$ is a (strictly associative) algebra object of $\bfA$, then $\Mod_{A}(\bfA)$ inherits a simplicial structure in a canonical way. Moreover:

\begin{proposition}\label{marchingduck}
Let $\bfA$ be a combinatorial monoidal model category, and let $A$ be an algebra object of $\bfA$ which is cofibrant as an object of $\bfA$. Then
$\Mod_{A}(\bfA)$ has the structure of a combinatorial model category, where:
\begin{itemize}
\item[$(W)$] A morphism $f: M \rightarrow N$ is a weak equivalence in $\Mod_{A}(\bfA)$ if and only if it is a weak equivalence in $\bfA$.
\item[$(F)$] A morphism $f: M \rightarrow N$ is a fibration in $\Mod_{A}(\bfA)$ if and only if it is a fibration in $\bfA$.
\end{itemize}
The forgetful functor $\Mod_{A}(\bfA) \rightarrow \bfA$ is both a left Quillen functor and a right Quillen functor. Moreover, if $\bfA$ is equipped with a compatible simplicial structure, then the induced simplicial structure on $\Mod_{A}(\bfA)$ endows $\Mod_{A}(\bfA)$ with the structure of a simplicial model category.
\end{proposition}

\begin{proof}
The proof is similar to that of Proposition \ref{comboline}. We first observe that $\Mod_{A}(\bfA)$ is presentable (Corollary \ref{underwhear}). Let
$T: \Mod_{A}(\bfA) \rightarrow \bfA$ be the forgetful functor. Then $T$ admits a left adjoint
given by the formula $F(X) = A \otimes X$, and a right adjoint given by the formula $G(X) = ^{A}\!X$. 
Since $\bfA$ is combinatorial, there exists a (small) collection of morphisms $I = \{ i_{\alpha}: C \rightarrow C' \}$ which generates the class of cofibrations in $\bfA$, and a (small) collection of morphisms $J = \{ j_{\alpha}: D \rightarrow D' \}$ which generates the class of trivial cofibrations in $\bfA$. Let $\overline{F(I)}$ be the saturated class of morphisms in $\Mod_{A}(\bfA)$ generated by $\{ F(i): i \in I \}$, and let $\overline{F(J)}$ be defined similarly. Unwinding the definitions, we see that a morphism in $\Mod_{A}(\bfA)$ is a trivial fibration if and only if it has the right lifting property with respect to $F(i)$, for every $i \in I$. Invoking the small object argument, we deduce that every morphism $f: M \rightarrow N$ in $\Mod_{A}(\bfA)$ admits a factorization
$M \stackrel{f'}{\rightarrow} N' \stackrel{f''}{\rightarrow} N$ where $f' \in \overline{F(I)}$ and
$f''$ is a trivial fibration. Similarly, we can find an analogous factorization where
$f' \in \overline{F(J)}$ and $f''$ is a fibration. Using standard arguments, we are reduced to the problem of showing that each morphism belonging to $\overline{F(J)}$ is a weak equivalence in
$\Mod_{A}(\bfA)$. Let $S$ be the collection of all morphisms $f: M \rightarrow N$ in
$\Mod_{A}(\bfA)$ such that $T(f)$ is a trivial cofibration in $\bfA$. We wish to prove that
$\overline{F(J)} \subseteq S$. Since $T$ preserves colimits, we conclude that $S$ is saturated; it will therefore suffice to show that for each $j \in J$, $F(j) \in S$. In other words, we must show that
if $j: X \rightarrow Y$ is a trivial cofibration in $\bfA$, then the induced map $A \otimes X \rightarrow A \otimes Y$ is again a trivial cofibration. This follows immediately from the definition of a monoidal model category, in view of our assumption that $A$ is a cofibrant object of $\bfA$. This completes the proof that $\Mod_{A}(\bfA)$ is a model category.

The forgetful functor $T: \Mod_{A}(\bfA) \rightarrow \bfA$ is a right Quillen functor by construction.
To see that $T$ is also a left Quillen functor, it suffices to show that the right adjoint
$G: \bfA \rightarrow \Mod_{A}(\bfA)$ preserves fibrations and trivial fibrations. In view of the definition of fibrations and trivial fibrations in $\Mod_{A}(\bfA)$, this is equivalent to the assertion that the composition $T \circ G: \bfA \rightarrow \bfA$ preserves fibrations and trivial fibrations. This follows immediately from the definition of a monoidal model category, since $A \in \bfA$ is cofibrant. 

Now suppose that $\bfA$ is equipped with a compatible simplicial structure. We claim that
$\Mod_{A}(\bfA)$ inherits the structure of a simplicial model category. For this, we suppose
that $f: M \rightarrow N$ is a fibration in $\Mod_{A}(\bfA)$ and that
$g: X \rightarrow Y$ is a cofibration of simplicial sets. We wish to show that the induced map
$M^{Y} \rightarrow M^{X} \times_{ N^{X} } N^{Y}$ is a fibration in $\Mod_{A}(\bfA)$, which is trivial if either $f$ or $g$ is trivial. This follows immediately from the analogous statement in the simplicial model category $\bfA$.
\end{proof}

\begin{remark}
Proposition \ref{marchingduck} admits the following generalization: suppose that $\bfA$ is a combinatorial model category, and that $F: \bfA \rightarrow \bfA$ is a left Quillen functor which is equipped with the structure of a monad. Then the category of algebras over $F$ inherits a model structure, where the fibrations and weak equivalences are defined at the level of the underlying objects of $\bfA$. The proof is exactly as given above. 
\end{remark}

\begin{remark}
If $\bfA$ is a symmetric monoidal model category satisfying the monoid axiom, then the hypothesis that the algebra $A$ be cofibrant is superfluous; see \cite{monmod}. 
\end{remark}

\begin{theorem}\label{duckblur}
Let $\bfA$ be a combinatorial simplicial model category equipped with a compatible monoidal structure, let 
$\theta: \Nerve(\Alg(\bfA)^{\degree}) \rightarrow \Alg( \Nerve(\bfA^{\degree}) )$
be as in Theorem \ref{beckify}, and let $A \in \Alg(\bfA)^{\degree}$. Then the canonical map 
$\Nerve( \Mod_{A}(\bfA)^{\degree} ) \rightarrow \Mod_{ \theta(A)}( \Nerve(\bfA^{\degree}) )$
is an equivalence of $\infty$-categories. 
\end{theorem}

\begin{proof}
Consider the diagram
$$ \xymatrix{ \Nerve( \Mod_A(\bfA)^{\degree}) \ar[rr] \ar[dr]^{G} & & \Mod_{\theta(A)}( \Nerve( \bfA^{\degree}) )\ar[dl]^{G'} \\
& \Nerve(\bfA^{\degree}). & }$$
We will show that this diagram satisfies the hypotheses of Corollary \ref{littlerbeck}:
\begin{itemize}
\item[$(a)$] The $\infty$-categories $\Nerve( \Mod_A(\bfA)^{\degree} )$ and
$\Mod_{\theta(A)}( \Nerve(\bfA^{\degree}))$ admit geometric realizations of simplicial objects. In fact, both of these $\infty$-categories are presentable. For $\Nerve( \Alg(\bfA)^{\degree})$, this follows from Proposition \toposref{notthereyet} and \ref{marchingduck}. For $\Mod_{\theta(A)}( \Nerve(\bfA)^{\degree})$, we first observe that $\Nerve(\bfA)^{\degree}$ is presentable (Proposition \toposref{notthereyet}) and that the functor $A \otimes \bigdot$ preserves colimits, and then apply Corollary \ref{underwhear}.

\item[$(b)$] The functors $G$ and $G'$ admit left adjoints $F$ and $F'$. The existence of a left adjoint to $G$ follows from the fact that $G$ is determined by a right Quillen functor. The existence of a left adjoint to $G'$ follows from Proposition \ref{pretara}.

\item[$(c)$] The functor $G'$ is conservative and preserves geometric realizations of simplicial objects. This follows from Corollaries \ref{thetacart} and \ref{gloop}.

\item[$(d)$] The functor $G$ is conservative and preserves geometric realizations of simplicial objects. The first assertion is immediate from the definition of the weak equivalences in
$\Mod_A(\bfA)$, and the second follows from the fact that $G$ is also a left Quillen functor.

\item[$(e)$] The natural map $G' \circ F' \rightarrow G \circ F$ is an equivalence. 
Unwinding the definitions, we are reduced to proving that if $N \in \bfA^{\degree}$, then 
the natural map $N \rightarrow A \otimes N$ induces an equivalence $F'(N) \simeq A \otimes N$. This is a special case of Proposition \ref{pretara}.
\end{itemize}
\end{proof}

\subsection{The $\infty$-Category of Tensored $\infty$-Categories}\label{hugree}

According to Remark \ref{otherlander}, there are at least three possible ways to define a monoidal $\infty$-category:
\begin{itemize}
\item[$(A)$] As a coCartesian fibration $\calC^{\otimes} \rightarrow \Nerve(\cDelta)^{op}$; this is the point of view we have chosen to adopt.
\item[$(B)$] As a monoid object in the $\infty$-category $\Cat_{\infty}$.
\item[$(C)$] As an algebra object of $\Cat_{\infty}$, where we view $\Cat_{\infty}$ as endowed with the Cartesian monoidal structure.
\end{itemize}

The purpose of this section is to obtain a similar picture for the theory of $\infty$-categories tensored over $\calC$. Our first step is to establish an analogue of Proposition \ref{ungbat}, which will allow us to describe module objects in an $\infty$-category endowed with a Cartesian monoidal structure.

\begin{definition}\index{left monoid action}\index{monoid action}
Let $\calC$ be an $\infty$-category. A {\it left monoid action in $\calC$} is a functor
$F: \Delta^1 \times \Nerve(\cDelta)^{op} \rightarrow \calC$ with the following properties:
\begin{itemize}
\item[$(1)$] The restriction $F | \{1\} \times \Nerve(\cDelta)^{op} \rightarrow \calC$ is a monoid object of $\calC$ (Definition \ref{monob}). 
\item[$(2)$] For $n \geq 0$, the maps
$F(0, \{n\}) \leftarrow F(0, [n]) \rightarrow F(1, [n])$
exhibit $F(0,[n])$ as a product $F(0,\{n\}) \times F(1, [n])$ in $\calC$.
\end{itemize}
We let $\Mon^{L}(\calC)$ denote the full subcategory of $\Fun( \Delta^1 \times \Nerve(\cDelta)^{op}, \calC)$ spanned by the left monoid actions.\index{ZZZMonLC@$\Mon^{L}(\calC)$}
\end{definition}

We now have the following analogue of Proposition \ref{ungbat}:

\begin{proposition}\label{ungbat2}
Let $\calC$ be an $\infty$-category which admits finite products, let
$p: \calC^{\times} \rightarrow \Nerve(\cDelta)^{op}$ be the associated $($Cartesian$)$ monoidal $\infty$-category $($see Proposition \ref{commonton}$)$, and let 
$\calC^{\times,L} \rightarrow \calC^{\times}$ exhibit $\calC$ as left tensored over itself $($Example \ref{sumai}$)$. Then composition with the
Cartesian structure $\calC^{\times} \rightarrow \calC$ induces an equivalence of $\infty$-categories
$\Mod( \calC ) \rightarrow \Mon^{L}(\calC).$ 
\end{proposition}

\begin{proof}
Let $t: \cDelta \rightarrow \cDelta$ be the functor $[n] \mapsto [n] \star [0]$, and $\alpha: \id \rightarrow t$ the canonical transformation, which we may identify with a functor $[1] \times \cDelta \rightarrow \cDelta$. Let $\calI = \cDelta^{\times} \times_{\cDelta} ( [1] \times \cDelta )$.
Objects of $\calI$ can be identified with triples $(a, [n], i \leq j)$, where 
$0 \leq a \leq 1$ and $0 \leq i \leq j \leq n+a$. Let $\calI_0$ be the full subcategory (equivalent to $[1] \times \cDelta$) spanned by objects of the form $(a, [n], 0 \leq n+a)$. 

Unwinding the definitions, we see that $\Mod(\calM)$ is equivalent to the full subcategory of
$\Fun( \Nerve(\calI)^{op}, \calC)$ spanned by those functors $F: \Nerve(\calI)^{op} \rightarrow \calC$ which have the following properties:
\begin{itemize}
\item[$(i)$] For $0 \leq a \leq 1$ and $0 \leq i \leq j \leq n+a$, the natural maps
$F(a, [n], i \leq j ) \rightarrow F(a, [n], k \leq k+1)$
exhibit $F(a, [n], i \leq j)$ as a product
$\prod_{i \leq k < j} F(a, [n], k \leq k+1)$ in the $\infty$-category $\calC$.
\item[$(ii)$] For every {\em convex} map $\alpha: [m] \rightarrow [n]$, and every $0 \leq i \leq j \leq m$, the associated map
$F( 0, [n], \alpha(i) \leq \alpha(j) ) \rightarrow F( 0, [m], i \leq j)$ is an equivalence.
\item[$(iii)$] For every $0 \leq i \leq j \leq n$, the canonical map
$F(1, [n], i \leq j) \rightarrow F(0, [n], i \leq j)$ is an equivalence in $\calC$.
\item[$(iv)$] For every morphism $\alpha: [m] \rightarrow [n]$ such that $\alpha(m) = n$, the
map $$F(1, [n], n \leq n+1) \rightarrow F(1, [m], m \leq m+1)$$ is an equivalence.
\end{itemize}

On the other hand, $\Mon^{L}(\calC)$ can be identified with the $\infty$-category of functors
$F_0: \Nerve(\calI_0)^{op} \rightarrow \calC$ with the following properties:
\begin{itemize}
\item[$(i')$] For each $n \geq 0$, the natural map
$F_0(0, [n], 0 \leq n) \rightarrow \prod_{0 \leq i < n} F_0(0, \{i,i+1\}, i \leq i+1)$
is an equivalence.
\item[$(ii')$] For $n \geq 0$, the inclusion $[0] \simeq \{n\} \subseteq [n]$ induces
an equivalence
$$F_0(1, [n], 0 \leq n+1) \simeq F_0(0, [n], 0 \leq n) \times F_0(1, [0], 0 \leq 1).$$
\end{itemize}

The map $\Mod(\calM) \rightarrow \Mon^{L}(\calC)$ is given by restriction from $\calI$ to $\calI_0$.
In virtue of Proposition \toposref{lklk}, it will suffice to verify the following:
\begin{itemize}
\item[$(a)$] Every functor $F_0: \Nerve(\calI_0)^{op} \rightarrow \calC$ which satisfies
$(i')$ and $(ii')$ admits a right Kan extension to $\Nerve(\calI)$.
\item[$(b)$] A functor $F: \Nerve(\calI)^{op} \rightarrow \calC$ satisfies conditions
$(i)$ through $(iv)$ if and only if the restriction $F_0 = F | \Nerve(\calI_0)^{op}$ satisfies
$(i')$ and $(ii')$, and $F$ is a right Kan extension of $F_0$.
\end{itemize}

Assertion $(a)$ is obvious, since for every object $X = (a, [n], i \leq j)$, the
category $\calI_0 \times_{\calI} \calI_{X/}$ has a final object. Moreover, we see that a functor
$F: \Nerve(\calI)^{op} \rightarrow \calC$ is a right Kan extension of $F| \Nerve(\calI_0)^{op}$ if and only if the following condition is satisfied:

\begin{itemize}
\item[$(iii')$] For every $0 \leq a \leq 1$ and every  $0 \leq i \leq j \leq n$, the natural map
$$F(a, [n], i \leq j) \rightarrow F(0, \{ i, i+1, \ldots j \}, i \leq j)$$ is an equivalence.
\item[$(iv')$] For $0 \leq i \leq n+1$, the map
$F(1, [n], i \leq n+1) \rightarrow F(1, \{i, \ldots, n\}, i \leq n+1)$
is an equivalence.
\end{itemize}

It remains only to show that conditions $(i)$ through $(iv)$ are equivalent to conditions $(i')$ through $(iv')$. This is straightforward and left to the reader (we observe that each of the relevant conditions depends only on the underlying functor $\calI^{op} \rightarrow \h{\calC}$). 
\end{proof}

Our next goal is to relate monoid actions in $\Cat_{\infty}$ with tensored $\infty$-categories.
For this, we will need to introduce a bit of notation.

\begin{notation}
Let $\bfA = (\mSet)_{/ \Nerve(\cDelta)^{op} }$ denote the category of {\em marked} simplicial sets
$(X, \calE)$ equipped with a map $X \rightarrow \Nerve(\cDelta)^{op}$. We will regard
$\bfA$ as endowed with the model structure {\em opposite} to the marked model structure defined in \S \toposref{twuf}. The fibrant objects of $\bfA$ are of the form $X^{\natural}$, where
$p: X \rightarrow \Nerve(\cDelta)^{op}$ is a coCartesian fibration and $X^{\natural} = (X, \calE)$, where $\calE$ is the collection of $p$-coCartesian edges of $X$.

We let $\bfA^{[1]}$ denote the category of functors from $[1]$ to $\bfA$, endowed with the {\em injective} model structure (see \S \toposref{quasilimit3}). Then $\bfA^{[1]}$ inherits from $\bfA$ the structure of a simplicial model category. We let $\widetilde{\CatMod}$ denote the underlying
$\infty$-category $\Nerve( (\bfA^{[1]})^{\degree})$, obtained by taking the simplicial nerve of the category of fibrant-cofibrant objects of $\bfA^{[1]}$. 
\end{notation}

Note that every object of $\bfA^{[1]}$ is cofibrant. Consequently, we may identify
$( \bfA^{[1]})^{\degree}$ with the category of fibrant objects of $\bfA^{[1]}$. These, in turn, can be described as {\it fibrations} between fibrant objects of $\bfA$. The following result describes the structure of such a fibration:

\begin{proposition}\label{umpertein}
Suppose given a commutative triangle of $\infty$-categories
$$ \xymatrix{ X \ar[rr]^{f} \ar[dr]^{p} & & Y \ar[dl]^{q} \\
& S & }$$
where $p$ and $q$ are coCartesian fibrations, and $f$ carries $p$-coCartesian morphisms in $X$ to $q$-coCartesian morphisms in $Y$. The following conditions are equivalent:
\begin{itemize}
\item[$(1)$] The map $f$ is a categorical fibration of simplicial sets.
\item[$(2)$] The map $f$ induces a fibration $X^{\natural} \rightarrow Y^{\natural}$ in the
category $(\mSet)_{/S}$, endowed with the {\em opposite} of the model structure defined in
\S \toposref{twuf} $($here $X^{\natural}$ and $Y^{\natural}$ denote the marked simplicial sets obtained by marking all of the coCartesian edges of $X$ and $Y$, respectively$)$.
\end{itemize}
\end{proposition}

\begin{proof}
According to Theorem \toposref{bigdiag}, the forgetful functor
$(\mSet)_{/S} \rightarrow (\sSet)_{/S}$ is a right Quillen functor, where
$(\sSet)_{/S}$ is endowed with the Joyal model structure. Thus $(2) \Rightarrow (1)$.
Conversely, suppose that $f$ is a categorical fibration. 
We must show that every lifting problem
$$ \xymatrix{ A \ar@{^{(}->}[d]^{i} \ar[r]^{f} & X^{\natural} \ar[d]^{p} \\
B \ar[r]^{g} \ar@{-->}[ur]^{G} & Y^{\natural} }$$
in $(\mSet)_{/S}$ has a solution, provided that $i$ is a trivial cofibration in
$(\mSet)_{/S}$. We first invoke Proposition \toposref{markedfibrant} to deduce that $X^{\natural}$ and $Y^{\natural}$ are fibrant objects of $(\mSet)_{/S}$. Consequently, $f$ extends to a map
$F: B \rightarrow X^{\natural}$. The maps $p \circ F, g: B \rightarrow Y^{\natural}$ agree on
$A$. Since $Y^{\natural}$ is fibrant, we conclude that there is a homotopy
from $p \circ F$ to $g$. Thus, there exists a contractible Kan complex $K$ containing
a pair of vertices $v$ and $v'$, and a map $h: B \times K^{\flat} \rightarrow Y^{\natural}$ such
that $h| B \times \{v\} = p \circ F$ and $h | B \times \{v'\} = g$, and $h | A \times K^{\flat}$
coincides with the composition
$A \times K^{\flat} \rightarrow A \stackrel{f}{\rightarrow} Y^{\natural}.$
The inclusion
$$ (B \times \{v\} ) \coprod_{A \times \{v\} } (A \times K^{\flat}) \subseteq B \times K^{\flat}$$
is a categorical equivalence (on the level of the underlying simplicial sets). Invoking assumption $(1)$, we conclude that the homotopy $h$ lifts to a map of simplicial sets $H: B \times K \rightarrow X$ such that $H | A \times K$ coincides with the composition
$ A \times K \rightarrow A \stackrel{f}{\rightarrow} X$
and $H | B \times \{v\} = F$. Since the collection of $p$-coCartesian edges of $X$ is stable under equivalence, we deduce that $H$ underlies a map of marked simplicial sets
$B \times K^{\flat} \rightarrow X^{\natural}$. The restriction $H| B \times \{v'\}$ is therefore a map
$B \rightarrow X^{\natural}$ with the desired properties.
\end{proof}

Suppose given a diagram
$\calM^{\otimes} \stackrel{p}{\rightarrow} \calC^{\otimes} \rightarrow \Nerve(\cDelta)^{op}$
which exhibits an $\infty$-category $\calM$ as left-tensored over a monoidal $\infty$-category $\calC$. According to Proposition \ref{umpertein}, we can view $p$ as a fibration between fibrant objects of $\bfA$; that is, as an object of the $\infty$-category $\widetilde{\CatMod}$. Of course, not every object of $\widetilde{\CatMod}$ is of this form.

\begin{definition}\index{ZZZCatmod@$\CatMod$}
Let $\CatMod$ denote the full subcategory of $\widetilde{\CatMod}$ corresponding to those diagrams
$\calM^{\otimes} \rightarrow \calC^{\otimes} \rightarrow \Nerve(\cDelta)^{op}$
which exhibit $\calC^{\otimes}$ as a monoidal $\infty$-category and $\calM^{\otimes}$ as an $\infty$-category left-tensored over $\calC^{\otimes}$.
\end{definition}

 \begin{corollary}\label{spg}
 There are canonical equivalences of $\infty$-categories
 $$ \CatMod \simeq \Mon^{L}(\Cat_{\infty}) \simeq \Mod( \Cat_{\infty}),$$
where $\Cat_{\infty}$ is endowed with the Cartesian monoidal structure and regarded as left-tensored over itself.
\end{corollary}

\begin{proof}
According to Proposition \toposref{gumby4}, there is a (canonical) equivalence of
$\infty$-categories
$\widetilde{\CatMod} \simeq \Fun( \Delta^1, \Nerve(\bfA^{\degree})).$ Theorem \toposref{straightthm} allows us to identify $\Nerve(\bfA^{\degree})$ with 
$\Fun( \Nerve(\cDelta)^{op}, \Cat_{\infty})$. Consequently, we obtain an equivalence
$\widetilde{\CatMod} \simeq \Fun( \Delta^1 \times \Nerve(\cDelta)^{op}, \Cat_{\infty} ).$
Unwinding the definitions, we see that under this equivalence, $\CatMod$ corresponds to the full subcategory $\Mon^{L}(\Cat_{\infty}) \subseteq \Fun( \Delta^1 \times \Nerve(\cDelta)^{op}, \Cat_{\infty})$. This proves the existence of the first equivalence. The existence of the second follows from Proposition \ref{ungbat2}.
\end{proof}

\begin{corollary}\label{specialtime}
Let $p: \calC^{\otimes} \rightarrow \Nerve(\cDelta)^{op}$ be a monoidal $\infty$-category, let
$q: \calC^{\otimes, L} \rightarrow \calC^{\otimes}$ exhibit $\calC$ as left-tensored over itself,
let $q': \calM^{\otimes} \rightarrow \calC^{\otimes}$ be an arbitrary $\infty$-category left-tensored over $\calC$, and let
$\calE \subseteq \bHom_{\calC^{\otimes}}( \calC^{\otimes, L}, \calM^{\otimes} )$
be the full subcategory spanned by those functors which carry $q \circ p$-coCartesian edges to $q' \circ p$-coCartesian edges. Then evaluation on the unit object $1_{\calC}$ induces an equivalence
of $\infty$-categories
$\theta: \calE \rightarrow \calM = \calM^{\otimes}_{[0]}.$
\end{corollary}

\begin{proof}
It will suffice to prove that, for every simplicial set $K$, the map $\theta$ induces a bijection from
the set of equivalence classes of objects of $\Fun(K,\calE)$ to the set of equivalence classes of objects of $\Fun(K,\calM)$. Replacing $\calM^{\otimes}$ by the fiber product
$\calC^{\otimes} \times_{ \Fun(K, \calC^{\otimes}) } \Fun(K, \calM^{\otimes}),$
we may reduce to the case where $K = \Delta^0$. In this case, it suffices to show that
$\theta$ induces an equivalence from the largest Kan complex contained in $\calE$ to the largest Kan complex contained in $\calM$. The latter Kan complex is equivalent to the mapping space
$\bHom_{ \Cat_{\infty}}( \ast, \calM)$. Using Corollary \ref{spg}, we can identify the former with
$\bHom_{ \Mod_{\calC^{\otimes}}( \Cat_{\infty} )}( \calC^{\otimes,L}, \calM^{\otimes} ).$ The desired result is now a special case of Corollary \ref{ingurtime} (applied to the monoidal $\infty$-category $\Cat_{\infty}$, endowed with the Cartesian monoidal structure).
\end{proof}

\subsection{Endomorphism Algebras}\label{cupper}

Let $\calC$ be a monoidal category containing an object $X$, and suppose that there
exists an {\em endomorphism object $\End(X)$}: that is, an object of $\calC$ equipped with evaluation map\index{endomorphism object} $e: \End(X) \otimes X \rightarrow X$
which induces a bijection $\Hom_{\calC}( Y, \End(X) ) \rightarrow \Hom_{\calC}( Y \otimes X, X)$
for every $Y \in \calC$ (such an object always exists if the monoidal category $\calC$ is {\em closed}). Then $\End(X)$ is automatically equipped with the structure of an algebra object of $\calC$: the unit map $1_{\calC} \rightarrow \End(X)$ classifies the identity map $\id_{X}: X \rightarrow X$, while the multiplication $\End(X) \otimes \End(X) \rightarrow \End(X)$
classifies the composition
$$ \End(X) \otimes (\End(X) \otimes X) \rightarrow \End(X) \otimes X \rightarrow X.$$

Our goal in this section is to obtain an $\infty$-categorical generalization of the above construction.
With an eye toward later applications, we will treat a slightly more general problem. Suppose that
$\calC$ is a monoidal $\infty$-category, $\calM$ an $\infty$-category which is left-tensored over $\calC$, and $M$ is an object of $\calM$. Let us say that an object $C \in \calC$ {\it acts on $M$} if we are given a map $C \otimes M \rightarrow M$. We would like to extract an object $\End(M) \in \calC$ which is {\em universal} among objects which act on $M$, and show that $\End(M)$ is an algebra object of $\calC$ (and further, that $M$ can be promoted to an object of
$\Mod_{\End(M)}(\calM)$). 

Our first goal is to give a careful formulation of the universal property desired of $\End(M)$.
We would like to have, for every algebra object $A \in \Alg(\calC)$, a homotopy equivalence of
$\bHom_{ \Alg(\calC)}( A, \End(M) )$ with a suitable classifying space for actions of $A$ on $M$. The natural candidate for this latter space is the fiber product
$\{ A \} \times_{ \Alg(\calC) } \Mod(\calM) \times_{\calM} \{ M \},$
which can be viewed as a fiber of the projection map $\theta: \Mod(\calM) \times_{\calM} \{M\} \rightarrow \Alg(\calC)$. We will see below that $\theta$ is a right fibration (Corollary \ref{hughp}).
Our problem is therefore to find an object $\End(M) \in \Alg(\calC)$ which
{\em represents} the right fibration $\theta$, in the sense that we have a 
an equivalence $\Mod(\calM) \times_{ \calM} \{ M \} \simeq \Alg(\calC)_{/ \End(M) }$ of
right fibrations over $\Alg(\calC)$ (see \S \toposref{quasilimit7}).\index{ZZZEndM@$\End(M)$}

The next step is to realize the $\infty$-category $\Mod( \calM) \times_{ \calM } \{ M \}$ as (equivalent to) the $\infty$-category of {\em algebra objects} in a suitable monoidal $\infty$-category
$\calC[M]$. Roughly speaking, we will think of objects of $\calC[M]$ as pairs 
$(C, \eta)$, where $C \in \calC$ and $\eta: C \otimes M \rightarrow M$ is a morphism in $\calM$.
The monoidal structure on $\calC[M]$ may be described informally by the formula
$(C, \eta) \otimes (C', \eta') = (C \otimes C', \eta''),$
where $\eta''$ denotes the composition
$$ C \otimes C' \otimes M \stackrel{ \id \otimes \eta'}{\rightarrow} C \otimes M
\stackrel{\eta}{\rightarrow} M.$$
The desired object $\End(M)$ can be viewed as a {\em final} object of $\calC[M]$. Provided that
this final object exists, it automatically has the structure of an algebra object of $\calC[M]$ (Corollary \ref{firebaugh}). The image of $\End(M)$ under the (monoidal) forgetful functor
$\calC[M] \rightarrow \calC$ will therefore inherit the structure of an algebra object of $\calC$.

We are now ready to begin with a detailed definition of the monoidal $\infty$-category $\calC[M]$.

\begin{notation}\label{skinnytrip}
We define a category $\calJ$ as follows:
\begin{itemize}
\item[$(a)$] An object of $\calJ$ is either pair $( [n], i \leq j)$, where $[n] \in \cDelta$ and
$0 \leq i \leq j \leq n$, or a pair $([n], \ast)$, where $[n] \in \cDelta$ and $\ast$ is a fixed symbol.
\item[$(b)$] Morphisms in $\calJ$ are given as follows:
$$ \Hom_{\calJ}( ( [m], i \leq j ), ( [n], i' \leq j' ) ) = \{ \alpha \in \Hom_{\cDelta}( [m], [n]):
i' \leq \alpha(i) \leq \alpha(j) \leq j' \} $$
$$ \Hom_{ \calJ}( ([m], \ast), ([n], i \leq j) ) = \{ ( \alpha, k): \alpha \in \Hom_{\cDelta}( [m], [n]), i \leq k \leq j \}$$ 
$$ \Hom_{ \calJ}( ([m], i \leq j ), ( [n], \ast) ) = \emptyset \quad \Hom_{ \calJ}( ([m], \ast), ([n], \ast) ) = \Hom_{\cDelta}([m], [n]).$$
\end{itemize}

Let $\cDelta'$ denote a new copy $\infty$-category $\cDelta$ (to avoid confusion below), and define functors
$$ \psi: \calJ \rightarrow \cDelta \quad \psi': \calJ \rightarrow \cDelta'$$
by the formulas
$$ \psi( [n], i \leq j ) = \psi( [n], \ast) = [n] \quad \psi'( [n], i \leq j) = \{ i, i+1, \ldots, j \} \quad \psi'( [n], \ast) = [0]. $$
We will identify $\cDelta$ with the full subcateory of $\calJ$ spanned by the objects $([n], \ast)$ (note that the remainder of $\calJ$ can be identified with the category $\cDelta^{\times}$ introduced in Notation \ref{hugegrin}). Similarly, we will identify $\cDelta'$ with the full subcategory
of $\calJ$ spanned by the objects $( [n], 0 \leq n)$. 

Let $\calM^{\otimes} \stackrel{q}{\rightarrow} \calC^{\otimes} \stackrel{p}{\rightarrow}
\Nerve(\cDelta')^{op}$ exhibit $\calM = \calM^{\otimes}_{[0]}$ as left-tensored over the monoidal
$\infty$-category $\calC = \calC^{\otimes}_{[1]}$, and let $M$ be an object of $\calM$. We
define a simplicial set $\widetilde{\calC[M]}^{\otimes}$ equipped with a map
$\widetilde{ \calC[M]}^{\otimes} \rightarrow \Nerve(\cDelta)^{op}$ as follows.
Let $K$ be an arbitrary simplicial set equipped with a map $e: K \rightarrow \Nerve(\cDelta)^{op}$. 
Then $\Hom_{ \Nerve(\cDelta)^{op} }( K, \widetilde{\calC[M]}^{\otimes})$ is in bijection with the set of 
commutative diagrams
$$ \xymatrix{ K \times_{ \Nerve(\cDelta)^{op} } \Nerve(\cDelta)^{op} \ar[r] \ar@{^{(}->}[d] & \{M\} \ar@{^{(}->}[d] \\
K \times_{ \Nerve(\cDelta)^{op} } \Nerve(\calJ)^{op} \ar[r] \ar[d] & \calM^{\otimes} \ar[d]^{p \circ q} \\
\Nerve( \cDelta')^{op} \ar@{=}[r] & \Nerve(\cDelta')^{op}. }$$

For each $[n] \in \cDelta$, we let $\calJ_{ [n] }$ denote the fiber
$\calJ \times_{ \cDelta } \{ [n] \}$. An object of $\calJ_{[n]}$ can be identified with either
a pair of integers $i$ and $j$ such that $0 \leq i \leq j \leq n$, or the symbol $\ast$.
A vertex of $\widetilde{\calC[M]}^{\otimes}$ can be identified with the following data:
\begin{itemize}
\item[$(i)$] An object $[n] \in \cDelta$.
\item[$(ii)$] A functor $f: \Nerve(\calJ_{[n]})^{op} \rightarrow \calM^{\otimes}$, covering
the map $\Nerve(\calJ_{[n]})^{op} \rightarrow \Nerve(\cDelta')^{op}$ induced by the functor $\psi'$.
\end{itemize}

We let $\calC[M]^{\otimes}$ denote the full simplicial subset of $\widetilde{\calC[M]}^{\otimes}$
spanned by those objects which classify functors $f: \Nerve( \calJ_{[n]})^{op} \rightarrow \calM^{\otimes}$ satisfying the following additional conditions:
\begin{itemize}
\item[$(1)$] For every morphism $\alpha$ in $\calJ_{[n]}$, the edge
$(q \circ f)(\alpha) \in \Hom(\Delta^1, \calC^{\otimes})$ is $p$-coCartesian.
\item[$(2)$] Let $\alpha: ( [n], \ast) \rightarrow ([n], i \leq j)$ be the morphism in
$\calJ_{[n]}$ corresponding to the element $j \in \{ i, \ldots, j\}$. Then
$f(\alpha)$ is $(p \circ q)$-coCartesian.
\end{itemize}\index{ZZZCM@$\calC[M]$}\index{ZZZCMotimes@$\calC[M]^{\otimes}$}
We let $\calC[M]$ denote the fiber $\calC[M]^{\otimes}_{[1]}$. 
\end{notation}

Our first main result is:

\begin{proposition}\label{umptump}
Let $\calC$ be a monoidal $\infty$-category, $\calM$ an $\infty$-category which is left-tensored over $\calC$, and $M \in \calM$ an object. Then:
\begin{itemize}
\item[$(1)$] The map $\calC[M]^{\otimes} \rightarrow \Nerve(\cDelta)^{op}$ constructed in
Notation \ref{skinnytrip} is a monoidal $\infty$-category.
\item[$(2)$] Restriction to $\cDelta' \subseteq \calJ$ induces a monoidal functor
$\calC[M]^{\otimes} \rightarrow \calC^{\otimes}$.
\end{itemize}
\end{proposition}

We will give the proof at the end of this section. First, we would like to have a better understanding of the category $\calC[M]$. An object of
$\calC[M]$ can be identified with a diagram
$$ \xymatrix{ M_0 \ar[dr]^{\sim} & N \ar@<.5ex>[d] \ar@<-.5ex>[d] \ar[l] \ar[r] & M_1 \ar[dl]^{\sim} \\
& M & }$$
in $\calM^{\otimes}$. 
Here $N \in \calM^{\otimes}_{[1]} \simeq \calC \times \calM$.
We may therefore identify $N$ with a pair of objects $C \in \calC$, $M' \in \calM$. 
The triangle on the right determines equivalences $M' \simeq M_1 \simeq M$ in
$\calM$, while the triangle on the left determines an equivalence $M_0 \simeq M$ and
a map $C \otimes M' \rightarrow M_0$. Consequently,
every object of $\calC[M]$ determines a morphism $\alpha: C \otimes M \rightarrow M$ in $\calM$, which is well-defined up to homotopy. In fact, we have the following more precise statement:

\begin{proposition}\label{ugher}
Let $\calC$ be a monoidal $\infty$-category, $\calM$ an $\infty$-category which is left-tensored over $\calC$, and $M \in \calM$ an object. Then:
\begin{itemize}
\item[$(1)$] Consider the functor $\calC \rightarrow \calM$ given by $C \mapsto C \otimes M$.
There exists an equivalence of $\infty$-categories $e: \calC \times_{ \calM} \calM^{/M} \rightarrow \calC[M]$ such that the composition of $e$ with the forgetful functor $\calC[M] \rightarrow \calC$
is equivalent to the projection $\calC \times_{ \calM} \calM^{/M} \rightarrow \calC$.

\item[$(2)$] The forgetful functor $\calC[M] \rightarrow \calC$ is a right fibration.

\item[$(3)$] Let $N \in \calC[M]$ be an object which classifies a morphism
$\alpha: C \otimes M \rightarrow M$ (see the discussion preceding Lemma \ref{staplegun}). Then
$N$ is a final object of $\calC[M]$ if and only if $\alpha$ exhibits $C$ as a morphism object
$\Mor_{\calM}(M,M)$ (Definition \ref{supiner}).
\end{itemize}
\end{proposition}

Again, we will defer the proof until the end of this section.

\begin{remark}\label{girlytime}
Proposition \ref{ugher} implies that the fiber
$\calC[M] \times_{ \calC} \{ 1_{\calC} \}$ can be identified with the Kan complex
$\calM^{/M} \times_{\calM} \{M\} \simeq \bHom_{\calM}(M,M)$. The monoidal structure
on $\calC[M]$ induces a coherently associative multiplication on $\bHom_{\calM}(M,M)$, which simply encodes the composition in the $\infty$-category $\calM$. In particular, if
$\widetilde{1}_{\calC}$ is an object of $\calC[M] \times_{\calC} \{1_{\calC} \}$ which classifies
an {\em equivalence} from $M$ to itself, then $\widetilde{1}_{\calC}$ is an {\em invertible object} of $\calC[M]$ (see Example \ref{satline}).
\end{remark}

\begin{corollary}\label{sunn}
Let $\calC$ be a monoidal $\infty$-category, $\calM$ an $\infty$-category which is left-tensored over $\calC$, and $M \in \calM$ an object. Then:
\begin{itemize}
\item[$(1)$] For each $n \geq 0$, the forgetful functor
$f: \calC[M]^{\otimes} \rightarrow \calC^{\otimes}$ induces a right fibration
$f_{[n]}: \calC[M]^{\otimes}_{[n]} \rightarrow \calC^{\otimes}_{[n]}$.
\item[$(2)$] Consider the diagram
$$ \xymatrix{ \calC[M]^{\otimes} \ar[rr]^{f} \ar[dr]^{q} & & \calC^{\otimes} \ar[dl]^{p} \\
& \Nerve(\cDelta)^{op}. & }$$
A morphism $\alpha$ in $\calC[M]^{\otimes}$ is $q$-coCartesian if and only if
$f(\alpha)$ is $p$-coCartesian.
\item[$(3)$] Let $A$ be a section of $q$. Then $A \in \Alg( \calC[M])$ if and only if
$f \circ A \in \Alg(\calC)$. Similarly, an object $A_0 \in \bHom_{ \Nerve(\cDelta)^{op} }( \Nerve(\cDelta^{\nounit})^{op}, \calC[M]^{\otimes})$ belongs to $\Alg^{\nounit}(\calC[M])$ if and only if
$f \circ A_0 \in \Alg^{\nounit}(\calC)$.
\item[$(4)$] Composition with $f$ induces isomorphisms of simplicial sets
$$\Alg(\calC[M]) \rightarrow \Alg(\calC) \times_{ \bHom_{ \Nerve(\cDelta)^{op}}( \Nerve(\cDelta)^{op}, \calC^{\otimes} )} \bHom_{ \Nerve(\cDelta)^{op} }( \Nerve(\cDelta)^{op} , \calC[M]^{\otimes}).$$
$$\Alg^{\nounit}(\calC[M]) \rightarrow \Alg^{\nounit}(\calC) \times_{ \bHom_{ \Nerve(\cDelta)^{op}}( \Nerve(\cDelta^{\nounit})^{op}, \calC^{\otimes} )} \bHom_{ \Nerve(\cDelta)^{op} }( \Nerve(\cDelta^{\nounit})^{op} , \calC[M]^{\otimes}).$$
\end{itemize}
\end{corollary}

\begin{proof}
Each of the functors $f_{[n]}$ is equivalent to the $n$th power of the forgetful functor
$\calC[M] \rightarrow \calC$. Consequently, $(1)$ follows immediately from Proposition \ref{ugher} and Lemma \ref{staplegun}. In particular, we deduce that a morphism $\alpha$ in
$\calC[M]^{\otimes}_{[n]}$ is an equivalence if and only if $f_{[n]}(\alpha)$ is an equivalence.
Since $f$ preserves coCartesian edges (Proposition \ref{umptump}), assertion $(2)$ follows.
The implications $(2) \Rightarrow (3)$ and $(3) \Rightarrow (4)$ are obvious.
\end{proof}

In the situation of Proposition \ref{umptump}, we obtain an induced map
$\Alg( \calC[M] ) \rightarrow \Alg(\calC)$. We may therefore think of
an algebra object of $\calC[M]$ as an algebra object of $\calC$ equipped with some kind of additional structure. The following result makes this idea precise:

\begin{proposition}\label{poofer}
Let $\calC$ be a monoidal $\infty$-category, $\calM$ an $\infty$-category which is left-tensored over $\calC$, and $M \in \calM$ an object. Then composition with the functor
$\psi'$ of Notation \ref{skinnytrip} induces categorical equivalences
$$ \Mod( \calM) \times_{\calM} \{ M \} \rightarrow \Alg( \calC[M] ) \quad  \Mod^{\nounit}(\calM) \times_{ \calM} \{M\} \rightarrow \Alg^{\nounit}( \calC[M] ).$$
\end{proposition}

In other words, we can identify algebra objects of $\calC[M]$ with algebra objects
$A \in \Alg(\calC)$, together with an action of $A$ on the fixed object $M \in \calM$.
The proof will be given at the end of this section. 

\begin{proposition}\label{lave}
Suppose given a diagram
$$ \xymatrix{ \calC^{\otimes} \ar[rr]^{f} \ar[dr]^{q} & & \calD^{\otimes} \ar[dl]^{p} \\
& \Nerve(\cDelta)^{op} & }$$
which exhibits $\calC = \calC^{\otimes}_{[1]}$ and $\calD= \calD^{\otimes}_{[1]}$ as monoidal $\infty$-categories. Suppose further that $f$ is a categorical fibration, a monoidal functor, and that
$f$ induces a right fibration $\calC \rightarrow \calD$. Then composition with $f$ induces right fibrations
$$ \Alg(\calC) \rightarrow \Alg(\calD) \quad \Alg^{\nounit}(\calC) \rightarrow \Alg^{\nounit}(\calD).$$
\end{proposition}

\begin{corollary}\label{hughp}
Let $\calC$ be a monoidal $\infty$-category, $\calM$ an $\infty$-category which is left-tensored over $\calC$, and $M \in \calM$ an object. Then:
\begin{itemize}
\item[$(1)$] The forgetful functors
$$\Alg( \calC[M] ) \rightarrow \Alg(\calC) \quad  \Alg^{\nounit}(\calC[M]) \rightarrow \Alg^{\nounit}(\calC)$$
are right fibrations of simplicial sets.
\item[$(2)$] The forgetful functors
$$ \Mod( \calM ) \times_{ \calM } \{ M \} \rightarrow \Alg(\calC) \quad  \Mod^{\nounit}( \calM ) \times_{ \calM } \{ M \} \rightarrow \Alg^{\nounit}( \calC) $$
are right fibrations of simplicial sets. 
\end{itemize}
\end{corollary}

\begin{proof}
Assertion $(1)$ follows from Propositions \ref{lave}, \ref{ugher}, and \ref{umptump}. Assertion
$(2)$ follows from $(1)$, Proposition \ref{poofer}, and Lemma \ref{staplegun}.
\end{proof}

\begin{corollary}\label{usebilly}\index{endomorphism object!universal property of}
Let $\calC$ be a monoidal $\infty$-category, $\calM$ an $\infty$-category which is left-tensored over $\calC$, and let $\overline{M} \in \Mod(\calM)$ a (left) module object having images
$M \in \calM$ and $A \in \Alg(\calC)$. Suppose that the multiplication map
$A \otimes M \rightarrow M$
exhibits $A$ as equivalent (in $\calC$) to a morphism object $\Mor_{\calM}(M,M)$. Then,
for every algebra object $B \in \Alg(\calC)$, we have a canonical isomorphism
$\bHom_{ \Alg(\calC)}(B,A) \simeq \Mod_{B}(\calM) \times_{ \calM } \{ M \}$
in the homotopy category $\calH$ of spaces.
\end{corollary}

\begin{proof}
Consider the diagram
$\Alg(\calC)_{/A} \leftarrow ( \Mod(\calM) \times_{ \calM } \{M \})_{/ \overline{M} }
\rightarrow \Mod(\calM) \times_{\calM} \{M \}.$
Propositions \ref{poofer} and \ref{ugher} imply that $\overline{M}$ is
a final object of $\Mod(\calM) \times_{\calM} \{ M\}$, so that the right map is a trivial Kan fibration. Since the map $\Mod(\calM) \times_{\calM} \{ M\} \rightarrow \Alg(\calC)$ is a right fibration (Corollary \ref{hughp}), the left map is also a trivial Kan fibration. 
Passing to fibers over the object $B \in \Alg(\calC)$, we obtain the desired result.
\end{proof}

\begin{remark}\label{usebully}
In the situation of Corollary \ref{usebilly}, let $A_0 \in \Alg^{\nounit}(\calC)$ denote the image of $A$ under the forgetful functor. Then the same argument shows that, for every nonunital algebra
object $B_0 \in \Alg^{\nounit}(\calC)$, we have a canonical isomorphism
$\bHom_{ \Alg^{\nounit}(\calC)}(B_0, A_0) \simeq \Mod^{\nounit}_{B_0}( \calM) \times_{ \calM } \{M \}$ in the homotopy category $\calH$ of spaces.
\end{remark}

The proof of Proposition \ref{lave} is based on the following easy lemma:

\begin{lemma}\label{lavalam}
Consider a diagram of simplicial sets
$$ \xymatrix{ X \ar[rr]^{f} \ar[dr]^{q} & & Y \ar[dl]^{p} \\
& S. & }$$
Suppose that:
\begin{itemize}
\item[$(i)$] The maps $p$ and $q$ are coCartesian fibrations.
\item[$(ii)$] The map $f$ carries $q$-coCartesian edges of $X$ to $p$-coCartesian edges of $Y$.
\item[$(iii)$] The map $f$ is a categorical fibration.
\item[$(iv)$] For every vertex $s$ of $S$, the induced map $X_{s} \rightarrow Y_{s}$ is a right fibration.
\end{itemize}
Then composition with $f$ induces a right fibration
$ \bHom_{S}( S, X) \rightarrow \bHom_{S}(S,Y)$. 
\end{lemma}

\begin{proof}
According to Corollary \toposref{helper}, it will suffice to show that the map
$$ \phi: \Fun( \Delta^1, \bHom_{S}(S,X) ) \rightarrow
\Fun( \Delta^1, \bHom_{S}(S,Y) ) \times_{ \Fun( \{1\}, \bHom_{S}(S,Y) )}
\Fun( \{1\}, \bHom_{S}(S,X) )$$
is a trivial Kan fibration. Set
$$ X' = \Fun( \Delta^1, X) \times_{ \Fun( \Delta^1, S) } S, \quad Y' = \Fun( \Delta^1, Y) \times_{ \Fun( \Delta^1, S) } S, \quad Z = Y' \times_{ \Fun( \{1\}, Y) } \Fun( \{1\}, X),$$
so that $f$ determines a map $\overline{\phi}: X' \rightarrow Z$, and $\phi$ can be identified with
the induced map $\bHom_{S}(S,X') \rightarrow \bHom_{S}(S,Z)$. It will therefore suffice to show that $\overline{\phi}$ is a trivial Kan fibration. Since $\overline{\phi}$ is clearly a categorical fibration, it will suffice to show that $\overline{\phi}$ is a categorical equivalence.

We have a commutative diagram
$$ \xymatrix{ X' \ar[dr]^{q'} \ar[r]^{\overline{\phi}} & Z \ar[d]^{r} \ar[r]^{\psi} & Y' \ar[dl]^{p'} \\
& S. & }$$
Proposition \toposref{doog} and $(i)$ implies that $p'$ and $q'$ are coCartesian fibrations, and
$(ii)$ implies that $\psi \circ \overline{\phi}$ carries $q'$-coCartesian edges to $p'$-coCartesian edges. Combining $(ii)$, $(iii)$, and Proposition \ref{umpertein}, we deduce that $r$ is also a coCartesian fibration, and that $\overline{\phi}$ carries $q'$-coCartesian edges to $r$-coCartesian edges. According to Proposition \toposref{apple1}, the map $\overline{\phi}$ is a categorial equivalence if and only if it induces a categorical equivalence $\overline{\phi}_{s}: X'_{s} \rightarrow Z_{s}$, for each vertex $s \in S$. We now observe that $(iv)$ and Corollary \toposref{helper} imply that $\overline{\phi}_{s}$ is a trivial Kan fibration.
\end{proof}

\begin{proof}[Proof of Proposition \ref{lave}]
We will give the proof in the unital case; the assertion for nonunital algebras follows using the same argument. According to Corollary \ref{sunn}, we have a pullback diagram
$$ \xymatrix{ \Alg(\calC) \ar[r]^{\theta} \ar[d] & \Alg(\calD) \ar[d] \\
\bHom_{ \Nerve(\cDelta)^{op} }( \Nerve(\cDelta)^{op}, \calC^{\otimes})
\ar[r]^{\theta'} & \bHom_{ \Nerve(\cDelta)^{op} }( \Nerve(\cDelta)^{op}, \calD^{\otimes}). }$$
Corollary \ref{sunn} and Lemma \ref{lavalam} imply that $\theta'$ is a right fibration.
It follows that $\theta$ is also a right fibration.
\end{proof}

We now return to the proofs of Propositions \ref{umptump}, \ref{ugher}, and \ref{poofer}. We first need to establish some technical preliminaries.

\begin{lemma}\label{starduck}
Suppose given a pullback square of $\infty$-categories
$$ \xymatrix{ X' \ar[r]^{q} \ar[d]^{p'} & X \ar[d]^{p} \\
S' \ar[r]^{q'} & S, }$$
where $p$ is a coCartesian fibration.
Let $f: K^{\triangleleft} \rightarrow X'$ be an arbitrary diagram. Then
$f$ is a $p'$-limit diagram if and only if $q \circ f$ is a $p$-limit diagram.
\end{lemma}

\begin{proof}
Replacing the above diagram by
$$ \xymatrix{ X'_{/f} \ar[r] \ar[d] & X_{/qf} \ar[d] \\
S'_{/p'f} \ar[r] & S_{/pqf} }$$
(and invoking Proposition \toposref{verylonger}), we can reduce to the case where
$K = \emptyset$. In this case, we can identify $f$ with an object $x' \in X'$. 
Let $s' = p'(x')$ and $s = q'(s')$. Corollary \toposref{superduck} implies that
$x'$ is $p'$-final if and only if $x'$ is a final object of the fiber $X'_{s'}$. Similarly, 
$x = q(x')$ is $p$-final if and only if $x$ is a final object of the fiber $X_{s}$. The desired result now follows from the observation that $X'_{s'} \simeq X_{s}$.
\end{proof}

\begin{lemma}\label{starfull}
Let $\calM^{\otimes} \stackrel{q}{\rightarrow} \calC^{\otimes} \stackrel{p}{\rightarrow} \Nerve(\cDelta)^{op}$ exhibit $\calM = \calM^{\otimes}_{[0]}$ as left-tensored over
$\calC = \calC^{\otimes}_{[1]}$. Let $n \geq 2$, and suppose we are given a diagram $\sigma_0$:
$$ \xymatrix{ X^{L} \ar[dr]^{\alpha} & & X^{R} \ar[dl]^{\beta} \\
& X^{LR} & }$$
in $\calM^{\otimes}$, lifting the diagram $\tau_0$:
$$ \xymatrix{ \{ 0, \ldots, n-1\} & & \{ 1, \ldots, n\} \\
& \{ 1, \ldots, n-1\} \ar[ur] \ar[ul] & }$$
in $\cDelta$. Suppose furthermore that $\alpha$ is a $(p \circ q)$-coCartesian morphism in
$\calM^{\otimes}$, and that $q(\beta)$ is a $p$-coCartesian morphism in $\calC^{\otimes}$.
Then:
\begin{itemize}
\item[$(1)$] The diagram $\tau$:
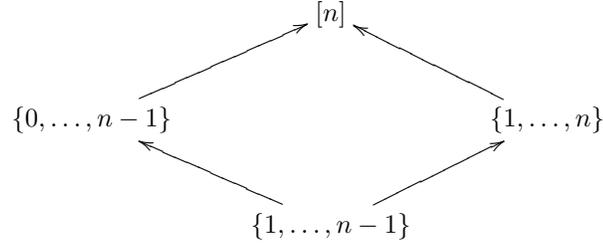
$$ \xymatrix{ & [n] & \\
\{ 0, \ldots, n-1\} \ar[ur] & & \{ 1, \ldots, n\} \ar[ul] \\
& \{ 1, \ldots, n-1\} \ar[ur] \ar[ul] & }$$
in $\cDelta$ can be lifted to a $(p \circ q)$-limit diagram $\sigma$:
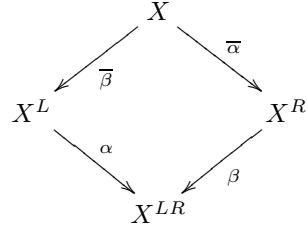
$$ \xymatrix{ & X \ar[dl]^{\overline{\beta}} \ar[dr]^{\overline{\alpha}} & \\
X^{L} \ar[dr]^{\alpha} & & X^{R} \ar[dl]^{\beta} \\
& X^{LR} & }$$
in $\calM^{\otimes}$. 

\item[$(2)$] Consider an arbitrary diagram $\sigma$ as in $(1)$, which is compatible with
both $\sigma_0$ and $\tau$. Then $\sigma$ is a $(p \circ q)$-limit diagram if and only if $\overline{\alpha}$ is a $(p \circ q)$-coCartesian morphism in $\calM^{\otimes}$ and $q(\overline{\beta})$ is $p$-coCartesian morphism in $\calC^{\otimes}$.
\end{itemize}

\end{lemma}

\begin{proof}
We first prove $(1)$. 
Let us identify $\tau$ with a map $\Delta^1 \times \Delta^1 \rightarrow \Nerve(\cDelta)^{op}$, and form a pullback diagram
$$ \xymatrix{ \calN \ar[d]^{r} \ar[r] & \calM^{\otimes} \ar[d]^{p \circ q} \\
\Delta^1 \times \Delta^1 \ar[r]^{\tau} \ar[r] & \Nerve(\cDelta)^{op}. }$$
We will regard $\calN$ as a simplicial subset of $\calM^{\otimes}$ containing the diagram
$\sigma_0$. In view of Lemma \ref{starduck} and the fact that $(p \circ q)$ is a coCartesian fibration, it will suffice to prove that $\sigma_0$ can be extended to an $r$-limit diagram $\sigma: \Delta^1 \times \Delta^1 \rightarrow \calN$.

Unwinding the definitions, we see that the coCartesian fibration $r: \calN \rightarrow \Delta^1 \times \Delta^1$ is classified by the following diagram of $\infty$-categories:
$$ \xymatrix{ & \calC^{\otimes}_{ \{0, 1\} } \times \calC^{\otimes}_{ \{ 1, \ldots, n-1\} } \times \calM^{\otimes}_{ \{ n-1, n\} }
\ar[dl]^{T} \ar[dr] & \\
\calC^{\otimes}_{ \{0,1\} } \times \calC^{\otimes}_{ \{1, \ldots, n-1\} } \times \calM^{\otimes}_{ \{n-1\} }\ar[dr] & & 
\calC^{\otimes}_{ \{1, \ldots, n-1\} } \times \calM^{\otimes}_{ \{n-1, n\} } \ar[dl]^{T} \\
& \calC^{\otimes}_{ \{1, \ldots, n-1\} } \times \calM^{\otimes}_{ \{n-1\} }. & }$$
Here $$T: \calM^{\otimes}_{ \{n-1, n\} } \simeq \calC \times \calM \stackrel{\otimes}{\rightarrow} \calM \simeq \calM^{\otimes}_{ \{n-1\} }$$ denotes the left action of $\calC$ on $\calM$.
Let $\calN'$ denote the relative nerve of this diagram, so that we have an equivalence
$f: \calN' \rightarrow \calN$ of coCartesian fibrations over $\Delta^1 \times \Delta^1$, and
$r': \calN' \rightarrow \Delta^1 \times \Delta^1$ the projection.

The diagram $\sigma_0$ is equivalent to a composition $f \circ \sigma'_0$. It will suffice to show that the analogues of $(1)$ and $(2)$ hold for $\widetilde{\sigma}_0$:
\begin{itemize}
\item[$(1')$] The map $\sigma'_0$ can be extended to an $r'$-limit diagram $\sigma': \Delta^1 \times \Delta^1 \rightarrow \calN'$ (which is simultaneously a section of $r'$).
\item[$(2')$] Let $\sigma'$ be an arbitrary section of $r'$ which extends $\sigma'_0$, and denote the diagram $f \circ \sigma'$ as follows:
$$ \xymatrix{ & X' \ar[dl]^{\overline{\beta}'} \ar[dr]^{\overline{\alpha}'} & \\
{X'}^{L} \ar[dr]^{\alpha'} & & {X'}^{R} \ar[dl]^{\beta} \\
& {X'}^{LR}. & }$$
Then $\sigma'$ is an $r'$-limit diagram if and only if $\overline{\alpha}'$ is a $(p \circ q)$-coCartesian morphism in $\calM^{\otimes}$ and $q(\overline{\beta}')$ is $p$-coCartesian morphism in $\calC^{\otimes}$.
\end{itemize}

The diagram $\sigma'_0$ determines objects
$$ (C, D, N) \in \calC^{\otimes}_{ \{0,1\} } \times \calC^{\otimes}_{ \{1, \ldots, n-1\} } \times \calM^{\otimes}_{ \{n-1\} }$$
$$ (D', M) \in \calC^{\otimes}_{ \{1, \ldots, n-1\} } \times \calM^{\otimes}_{ \{n-1, n\} } \quad (D'', N') \in \calC^{\otimes}_{ \{1, \ldots, n-1\} }\times \calM^{\otimes}_{ \{n-1\} }.$$
and morphisms
$$ D \stackrel{\gamma_0}{\rightarrow} D'' \stackrel{\gamma_1}{\leftarrow} D', \quad N \stackrel{\delta_0}{\rightarrow} N' \stackrel{\delta_1}{\leftarrow} T(M)$$
Since $\alpha$ is $(p \circ q)$-coCartesian, we conclude that $\gamma_0$ and $\delta_0$ are
equivalences in $\calC^{\otimes}_{ \{1, \ldots, n-1\} }$. Similarly, since $q(\beta)$ is $p$-coCartesian, we deduce that $\delta_0$ is an equivalence in $\calM^{\otimes}_{ \{n-1\} }$. Modifying $\sigma'_0$ if necessary, we may assume that
$D = D' = D''$, $N= N'$, and the maps $\gamma_0$, $\gamma_1$ and $\delta_0$ are all identities.
In this case, we can take $\sigma'$ to be the diagram
$$ \xymatrix{ & (C, D, M) \ar[dl]^{\id \times \id \times \delta_1} \ar[dr] & \\
(C, D, N) \ar[dr] & & (D, M) \ar[dl]^{\id \times \delta_1} & \\
& (D,N). & }$$
An easy calculation shows that $\sigma'$ is an $r'$-limit diagram. This proves $(1)$. Moreover, 
$\sigma'$ satisfies the criterion of $(2')$. Since $r'$-limit diagrams extending $\sigma'_0$ are uniquely determined up to equivalence, we deduce the ``only if'' direction of $(2')$. 

To prove the ``if'' direction of $(2')$, let us suppose given an arbitrary section $\overline{\sigma}'$ of $r'$ which extends $\sigma'_0$, depicted below:
$$ \xymatrix{ & (C_0, D_0, M_0) \ar[dl]^{f_0 \times f_1 \times f_2} \ar[dr]^{g_1 \times g_2} & \\
(C, D, N) \ar[dr] & & (D, M) \ar[dl]^{\id \times \delta_1} & \\
& (D,N). & }$$
Since $\sigma'$ is an $r'$-limit diagram, there exists a natural transformation
$s: \overline{\sigma}' \rightarrow \sigma'$ which is the identity, except possibly on the initial object.
Suppose that $\overline{\sigma}'$ satisfies the criterion of $(2')$. Then the maps
$f_0$, $f_1$, and $g_2$ are equivalences. Using the two-out-of-three property, we see that
$s$ induces an equivalence $(C_0, D_0, M_0) \rightarrow (C,D,M)$, so that $s$ is itself an equivalence. Since $\overline{\sigma}'$ is equivalent to $\sigma'$, it is also an $r'$-limit diagram, as we wished to prove.
\end{proof}

\begin{proof}[Proof of Proposition \ref{umptump}]
Let $\calM^{\otimes} \stackrel{q}{\rightarrow} \calC^{\otimes} \stackrel{p}{\rightarrow} \Nerve(\cDelta)^{op}$ exhibit $\calM = \calM^{\otimes}_{[0]}$ as left-tensored over
$\calC = \calC^{\otimes}_{[1]}$. The functor $\psi: \calJ \rightarrow \cDelta$ is an op-fibration of categories, so the induced map
$\Nerve(\calJ)^{op} \rightarrow \Nerve(\cDelta)^{op}$ is a Cartesian fibration. We now define
a simplicial set $X$ equipped with a map $f: X \rightarrow \Nerve(\cDelta)^{op}$, characterized by the following universal property: for every map of simplicial sets $K \rightarrow \Nerve(\cDelta)^{op}$, we have a canonical bijection
$\Hom_{\Nerve(\cDelta)^{op}}(K, X) \simeq \Hom_{ \Nerve(\cDelta')^{op} }( K \times_{ \Nerve(\cDelta)^{op}} \Nerve(\calJ)^{op}, \calM^{\otimes} ).$
Invoking Corollary \toposref{skinnysalad}, we deduce:
\begin{itemize}
\item[$(i)$] The map $f$ is a coCartesian fibration.
\item[$(ii)$] Let $\overline{\alpha}$ be a morphism in $X$, covering a map
$\alpha: [m] \rightarrow [n]$ in $\cDelta$. Then $\overline{\alpha}$ is $f$-coCartesian if and only if the following conditions are satisfied:
\begin{itemize}
\item[$(a)$] For every $0 \leq i \leq j \leq m$, the induced map $\overline{\alpha}( [n], \alpha(i) \leq \alpha(j) ) \rightarrow \overline{\alpha}( [m], i \leq j)$ is a $(q \circ p)$-coCartesian morphism in $\calM^{\otimes}$.
\item[$(b)$] The map $\overline{\alpha}( [n], \ast) \rightarrow \overline{\alpha}( [m], \ast)$ is a $(q \circ p)$-coCartesian morphism in $\calM^{\otimes}$. 
\end{itemize}
\end{itemize}

Restriction to the full subcategory $\cDelta \subseteq \calJ$ (spanned by the objects of the form
$( [n], \ast)$) gives rise to a commutative diagram
$$ \xymatrix{ X \ar[dr]^{f} \ar[rr]^{s} & & \calM \times \Nerve(\cDelta)^{op} \ar[dl]^{f'} \\
& \Nerve(\cDelta)^{op}. & }$$
The map $f'$ is obviously a coCartesian fibration, and $(ii)$ implies that $s$ carries $f$-coCartesian edges to $f'$-coCartesian edges. Moreover, we have a canonical isomorphism
$ \widetilde{ \calC[M]}^{\otimes} \simeq X \times_{ \calM \times \Nerve(\cDelta)^{op} } ( \{ M \} \times \Nerve(\cDelta)^{op} ).$
Combining Proposition \ref{umpertein} with Proposition \toposref{markedfibrant}, we deduce:
\begin{itemize}
\item[$(i')$] The projection $g: \widetilde{\calC[M]}^{\otimes} \rightarrow \Nerve(\cDelta)^{op}$ is a coCartesian fibration.
\item[$(ii')$] Let $\overline{\alpha}$ be a morphism in $\widetilde{\calC[M]}^{\otimes}$, covering a map $\alpha: [m] \rightarrow [n]$ in $\cDelta$. Then $\overline{\alpha}$ is $g$-coCartesian if and only if, for every $0 \leq i \leq j \leq m$, the induced map
$\overline{\alpha}( [n], \alpha(i) \leq \alpha(j) ) \rightarrow \overline{\alpha}( [m], i \leq j)$ is a $(q \circ p)$-coCartesian morphism in $\calM^{\otimes}$.
\end{itemize}

It follows easily from $(ii')$ that the map $g$ restricts to a coCartesian fibration
$g_0: \calC[M]^{\otimes} \rightarrow \Nerve(\cDelta)^{op}$, and that a morphism in
$\calC[M]^{\otimes}$ is $g_0$-coCartesian if and only if it is $g$-coCartesian.

We now prove $(1)$. We must show that, for each $n \geq 0$, the coCartesian fibration
$g_0$ induces an equivalence
$\calC[M]^{\otimes}_{[n]} \simeq \calC[M]^{\otimes}_{ \{0,1\} } \times \ldots \times \calM^{\otimes}_{ \{n-1, n\} }.$
For $n=0$ this is easy: the $\infty$-category $\calC[M]^{\otimes}_{[0]}$ is isomorphic to the full subcategory of $\calM^{/M}$ spanned by the final objects, and therefore a contractible Kan complex (Proposition \toposref{initunique}). For each object $[m] \in \cDelta$, let $\calJ_{[m]}$ denote the fiber product $\calJ \times_{ \cDelta} \{ [m] \}$. Let $P$ be the set of pairs of integers $i, j \in [n]$ satisfying $i \leq j$.
We endow $P$ with the following partial ordering:
$$ ( i \leq j) \leq (i' \leq j')$$ if and only if the interval
$\{i, \ldots, j\}$ is contained in the interval $\{ i', \ldots, j' \}$; in other words, if and only if
$0 \leq i' \leq i \leq j \leq j' \leq n$. We regard $P$ as equipped with a functor
$P \rightarrow \cDelta'$, given by $(i \leq j) \mapsto \{i, i+1, \ldots, j\} \simeq [j-i]$.
For every downward-closed subset $Q \subseteq P$, let $\calJ^{Q}_{[n]}$ denote the
full subcategory of $\calJ_{[n]}$ spanned by the objects $( [n], \ast)$ and
$\{ ( [n], i \leq j) \}_{ (i \leq j) \in Q  }$, and let
$\calE(Q)$ denote the full subcategory of 
$$ \bHom_{ \Nerve( \cDelta')^{op}}( \Nerve(\calJ_{[n]}^{P_0})^{op}, \calM^{\otimes})$$
spanned by those functors $F: \Nerve(\calJ_{[n]}^{0})^{op} \rightarrow \calM^{\otimes}$
which satisfy the following conditions:
\begin{itemize}
\item[$(a)$] Let $(i \leq j) \in Q$ and let
$0 \leq i \leq i' \leq j \leq n$. Then the induced map
$F( [n], i \leq j) \rightarrow F([n], i' \leq j)$ is a $(p \circ q)$-coCartesian morphism in $\calM^{\otimes}$.

\item[$(b)$] Let $(i \leq j) \in Q$ and let
$0 \leq i \leq i' \leq j' \leq j \leq n$. Then $q$ carries the induced map
$F( [n], i \leq j) \rightarrow F([n], i' \leq j')$ to a $p$-coCartesian morphism in $\calC^{\otimes}$.

\item[$(c)$] Let $(i \leq j) \in Q$ and let
$\alpha: ( [n], \ast) \rightarrow ([n], i \leq j)$ be the morphism in $\calJ_{[n]}$ classified by the element $j \in \{i, i+1, \ldots, j\}$. Then $F(\alpha)$ is a $(p \circ q)$-coCartesian morphism in $\calM^{\otimes}$.
\end{itemize}

It is convenient to reformulate the conditions $(a)$, $(b)$, and $(c)$. Using
a transitivity argument, it is easy to see that $(a)$ is equivalent to the following apparently weaker condition:

\begin{itemize}
\item[$(a')$] Let $(i < j) \in Q$. Then the induced map
$F( [n], i \leq j) \rightarrow F([n], i+1 \leq j)$ is a $(p \circ q)$-coCartesian morphism in $\calM^{\otimes}$.
\end{itemize}

Assuming $(a')$, another transitivity argument allows us to reformulate $(b)$ as follows:

\begin{itemize}
\item[$(b')$] Let $(i < j) \in Q$. Then $q$ carries the induced map
$F( [n], i \leq j) \rightarrow F([n], i \leq j-1)$ to a $p$-coCartesian morphism in $\calC^{\otimes}$.
\end{itemize}

Finally, assuming that $(a')$ is satisfied, another transitivity argument gives the following reformulation of $(c)$:

\begin{itemize}
\item[$(c')$] Let $(i \leq i) \in Q$. Then
the induced map $F( [n], i \leq i) \rightarrow F([n], \ast)$ is $(p \circ q)$-coCartesian.
\end{itemize}

Evaluation at the object $([n],\ast)$ induces a map $\calE(Q) \rightarrow \calM$; let
$\calE_0(Q)$ denote the fiber $\calE(Q) \times_{ \calM } \{M\}$. We observe that
$\calC[M]^{\otimes}_{[n]}$ is canonically isomorphic to $\calE_0(P)$.
Let $P_{\leq 1} \subseteq P$ be the collection of all pairs $(i \leq j) \in P$ such that $j \leq i+1$.
The simplicial set $\calE_0( P_{\leq 1} )$ is canonically isomorphic to the fiber product
$$ \calC[M]^{\otimes}_{ \{0,1\} } \times_{ \calC[M]^{\otimes}_{\{1\} } }
\ldots \times_{ \calC[M]^{\otimes}_{ \{n-1\} }} \calC[M]^{\otimes}_{ \{n-1, n\} },$$
and is therefore (since this fiber product is also a homotopy fiber product) equivalent
to the product $\prod_{0 \leq i < n} \calC[M]^{\otimes}_{ \{i, i+1\} }$. To complete the proof of $(1)$, it will suffice to show that the restriction map $r_0: \calE_0(P) \rightarrow \calE_0( P_{\leq 1})$ is a trivial Kan fibration. 

We have a Cartesian rectangle of simplicial sets
$$ \xymatrix{ \calE_0(P) \ar[r]^{r_0} \ar[d] & \calE_0(P_{\leq 1}) \ar[r] \ar[d] & \{M\} \ar[d] \\
\calE(P) \ar[r]^{r} & \calE(P_{\leq 1}) \ar[r] & \calM. }$$
It will therefore suffice to show that the map $r$ is a trivial Kan fibration.
We will prove the following more general statement: for every 
$P_{\leq 1} \subseteq Q \subseteq Q' \subseteq P$, the restriction map
$\calE(Q') \rightarrow \calE(Q)$ is a trivial Kan fibration. Using a transitivity argument, we
can reduce the case where $Q'$ is obtained from $Q$ by adjoining a single element
$(i \leq j) \in P-Q$. Since $P_{\leq 1} \subseteq Q$, we have $j \geq i+2$.

In view of Proposition \toposref{lklk}, it will suffice to prove the following, for every
$F_0 \in \calE(Q)$:

\begin{itemize}
\item[$(I)$] There exists a functor
$F \in \bHom_{ \Nerve(\cDelta')^{op} }( \Nerve(\calJ^{Q'}_{[n]})^{op}, \calM^{\otimes} )$
which is a $(p \circ q)$-right Kan extension of $F_0$.

\item[$(II)$] Let $F \in \bHom_{ \Nerve(\cDelta')^{op} }( \Nerve(\calJ^{Q'}_{[n]})^{op}, \calM^{\otimes})$
be an arbitrary extension of $F_0$. Then $F$ is a $(p \circ q)$-right Kan extension of
$F_0$ if and only if $F \in \calE(Q')$.
\end{itemize}

Let $\calI = \calJ^{Q}_{[n]} \times_{ \calJ_{[n]}} (\calJ_{[n]})_{/ ( [n], i \leq j)}$, 
and let $\calI_0$ denote the full subcategory of $\calI$ spanned by the objects
$$ X^{L}: ( [n], i \leq j-1) \rightarrow ( [n], i \leq j), \quad X^{R}: ( [n], i+1 \leq j) \rightarrow ([n], i \leq j)$$ 
$$X^{LR}: ( [n], i+1 \leq j-1) \rightarrow ([n], i \leq j).$$
Let $F'_0$ denote the composition
$ \Nerve(\calI)^{op} \rightarrow \Nerve(\calJ^{Q}_{[n]})^{op} \stackrel{F_0}{\rightarrow}
\calM^{\otimes}.$ Using Lemma \toposref{kan2} (and the equivalence between
conditions $(a)$, $(b)$, and $(c)$ with their analogues $(a')$, $(b')$, and $(c')$)
we are reduced to proving:

\begin{itemize}
\item[$(I')$] There exists a $(p \circ q)$-limit diagram $F'$, rendering following diagram commutative:
$$ \xymatrix{ \Nerve(\calI)^{op} \ar[r]^{F'_0} \ar@{^{(}->}[d] & \calM^{\otimes} \ar[d]^{p \circ q} \\
( \Nerve(\calI)^{op} )^{\triangleleft} \ar[r] \ar@{-->}[ur]^{F'} & \Nerve(\cDelta')^{op}.}$$
\item[$(II')$] Let $F'$ be an arbitrary map which renders the above diagram commutative. Then
$F'$ is a $(p \circ q)$-limit diagram if and only if $F'( \{X^R\}^{\triangleleft} )$ is
$(p \circ q)$-coCartesian morphism in $\calM^{\otimes}$, 
and $( q \circ F')( \{ X^L \}^{\triangleleft} )$ is $p$-coCartesian morphism in 
$\calC^{\otimes}$.
\end{itemize}

We observe that the inclusion $\calI_0 \subseteq \calI$ has a left adjoint, so the
induced map $\Nerve(\calI_0) \rightarrow \Nerve(\calI)$ is cofinal. Consequently, it suffices to prove the analogues of $(I')$ and $(II')$ obtained by replacing $\calI$ by $\calI_0$. In this case, the desired result is an immediate conseqence of Lemma \ref{starfull}. This completes the proof of $(1)$.

We now prove $(2)$. We wish to show that the restriction map
$\calC[M]^{\otimes} \rightarrow \calC^{\otimes}$ carries $g_0$-coCartesian edges to
$p$-coCartesian edges. Let $\overline{\alpha}$ be a $g_0$-coCartesian edge of
$\calC[M]^{\otimes}$, covering a map $\alpha: [m] \rightarrow [n]$ in $\cDelta$. Let
$\beta: \overline{\alpha}( [n], 0 \leq n) \rightarrow \overline{\alpha}( [m], 0 \leq m)$
be the induced map in $\calM^{\otimes}$. We wish to show that $q(\beta)$ is a $p$-coCartesian morphism in $\calC$. We observe that $q(\beta)$ factors as a composition
$$ q( \overline{\alpha}( [n], 0 \leq n)) \stackrel{q(\beta')}{\rightarrow} q(\overline{\alpha}( [n], \alpha(0) \leq \alpha(m) )) \stackrel{q(\beta'')}{\rightarrow} q(\overline{\alpha}( [m], 0 \leq m)).$$
Condition $(ii')$ implies that $\beta''$ is $(p \circ q)$-coCartesian, so that
$q(\beta'')$ is $p$-coCartesian. Moreover, since the domain of $\overline{\alpha}$ belongs
to $\calC[M]^{\otimes}$, the map $q(\beta')$ is $p$-coCartesian. It follows that
$q(\beta)$ is $p$-coCartesian, as desired.
\end{proof}

We now turn to the proof of Proposition \ref{ugher}. Once again, we will need a lemma.

\begin{lemma}\label{staplegun}
Suppose given a diagram of simplicial sets
$$ \xymatrix{ X' \ar[r] \ar[d]^{p'} & X \ar[d]^{p} \\
S' \ar[r] & S }$$
where $p$ is a right fibration, $p'$ a categorical fibration, and the horizontal arrows are categorical equivalences. Then $p'$ is also a right fibration.
\end{lemma}

\begin{proof}
Replacing $X$ by $X \times_{S'} S$ (and invoking Proposition \toposref{basechangefunky}), we can reduce to the case where $S = S'$ and the bottom horizontal map is the identity. We wish to show that every lifting problem of the form
$$ \xymatrix{ A \ar@{^{(}->}[d]^{i} \ar[r] & X' \ar[d]^{p'} \\
B \ar[r] \ar@{-->}[ur] & S' }$$
admits a solution, provided that $i$ is right anodyne. Applying Proposition \toposref{simpex} (in the model category $(\sSet)_{/S'}$, we can reduce to the problem of solving the associated mapping problem
$$ \xymatrix{ A \ar@{^{(}->}[d]^{i} \ar[r] & X \ar[d]^{p} \\
B \ar[r] \ar@{-->}[ur] & S, }$$
which is possibly in virtue of our assumption that $p$ is a right fibration.
\end{proof}

\begin{proof}[Proof of Proposition \ref{ugher}]
We first prove $(1)$. 
Let $\calN$ denote the fiber product $\calM^{\otimes} \times_{ \Nerve(\cDelta')^{op} } \Nerve(\calJ_{[1]})^{op}$, and $p: \calN \rightarrow \Nerve(\calJ_{[1]})^{op}$ the projection. Then $p$ is classified by a functor $F$ from $\Nerve(\calJ_{[1]})^{op}$ to the $\infty$-category $\Cat_{\infty}$. Unwinding the definitions, we see that $F$ is equivalent to the functor described by the diagram of $\infty$-categories
$$ \xymatrix{ \calM \ar@{=}[dr] & \calC \times \calM \ar@<.5ex>[d] \ar@<-.5ex>[d] \ar[l]^{\otimes} \ar[r]^{\pi} & \calM \ar@{=}[dl] \\
& \calM. & }$$
Here $\otimes: \calC \times \calM \rightarrow \calM$ denotes the tensor product functor, and
$\pi: \calC \times \calM \rightarrow \calM$ denotes the projection onto the second factor.
This diagram is described by a functor $F': \calJ_{[1]}^{op} \rightarrow \sSet$. Let 
$\calN'$ denote the relative nerve $\Nerve_{F'}( \calJ_{[1]}^{op} )$, so that
we have an equivalence $\calN \rightarrow \calN'$ of coCartesian fibrations over
$\Nerve( \calJ_{[1]})^{op}$. This equivalence induces an equivalence of $\infty$-categories
$ \calC[M]' \rightarrow \calC[M],$
where $\calC[M]'$ denotes the fiber product
$$ \{ M \} \times_{ \Fun( \{2\}, \calM) } \Fun'( \Delta^2, \calM) \times_{ \Fun( \{0\}, \calM ) }
\Fun( \{0\}, \calC \times \calM) \times_{ \Fun( \{0\}, \calM ) }
\Fun''( \Delta^2, \calM) \times_{ \Fun( \{2\}, \calM) } \{ M \};$$
here $\Fun'( \Delta^2, \calM)$ denotes the full subcategory of
$\Fun(\Delta^2, \calM)$ spanned by those diagrams
$$ \xymatrix{ & N' \ar[dr]^{\gamma} & \\
N \ar[ur]^{\beta} \ar[rr] & & N'' }$$
where $\gamma$ is an equivalence, and $\Fun''(\Delta^2, \calM)$ the full subcategory spanned by those diagrams where $\beta$ and $\gamma$ are both equivalences. Since $\Delta^2$ is weakly contractible, the diagonal map $\calM \rightarrow \Fun''( \Delta^2, \calM)$ is a categorical equivalence. An easy argument shows that this diagonal map induces a categorical equivalence
$\calC[M]'' \rightarrow \calC'[M]$, where $\calC[M]''$ denotes the fiber product
$$ \{ M \} \times_{ \Fun( \{2\}, \calM) } \Fun'( \Delta^2, \calM) \times_{ \Fun( \{0\}, \calM ) }
\Fun( \{0\}, \calC);$$
here $\calC$ maps to $\calM$ via the functor $C \mapsto C \otimes M$. 
We now observe that evaluation along the long edge of $\Delta^2$ induces a trivial Kan fibcation
$\Fun'( \Delta^2, \calM) \rightarrow \Fun( \Delta^{ \{0,2\} }, \calM)$.
We therefore obtain a trivial Kan fibration
$\calC[M]'' \rightarrow \calC \times_{\calM} \calM^{/M}.$
Let $s$ denote a section to this trivial fibration (for example, the section given by composition with a degeneracy map $\Delta^2 \rightarrow \Delta^1$), and define $e$ to be the composition
$$\calC \times_{\calM} \calM^{/M} \stackrel{s}{\rightarrow} \calC[M]''
\subseteq \calC[M]' \rightarrow \calC[M].$$
It is easy to see that $e$ satisfies the requirements of $(1)$.
Assertion $(2)$ follows from $(1)$ and Lemma \ref{staplegun}, since the projection $\calM^{/M} \rightarrow \calM$ is a right fibration.

We now prove $(3)$. Note that the objects of $\calC[M]''$ can be identified with pairs
$(C, \alpha)$, where $C \in \calC$ and $\alpha: C \otimes M \rightarrow M$ is a morphism in $\calM$. In view of $(1)$, it will suffice to show that an object $(C, \alpha) \in \calC[M]''$ is final if and only if
$\alpha$ exhibits $C$ as a morphism object $\Mor_{\calM}(M,M)$. Fix another object
In view of $(1)$, it will suffice to show that an object
$(C',\alpha') \in \calC[M]''$. Since the projection
$\calC[M]'' \rightarrow \calC$ is a right fibration, the induced map $q: \calC[M]''_{/ (C, \alpha)} \rightarrow \calC_{/C}$ is a trivial Kan fibration. We therefore obtain a homotopy fiber sequence
$$ \bHom_{ \calC[M]''}( (C', \alpha'), (C, \alpha) ) \rightarrow
\bHom_{ \calC}( C', C) \stackrel{\gamma}{\rightarrow} \calC[M]'' \times_{ \calC} \{ C' \}.$$ 
It follows that $(C, \alpha)$ is an initial object of $\calC[M]''$ if and only if $\gamma$ is a homotopy equivalence, for every choice of object $C' \in \calC$. We now observe that the codomain of $\gamma$ is canonically identified with $\bHom_{\calM}( C' \otimes M, M)$, so that
the $\gamma$ is an equivalence for {\em every} $C' \in \calC$ if and only if
$\alpha$ exhibits $C$ as a morphism object $\Mor_{\calM}(M,M)$. This proves $(3)$.
\end{proof}

\begin{proof}[Proof of Proposition \ref{poofer}]
We will give the proof for unital algebras; the nonunital case can be established using the same argument. Let $\calM^{\otimes} \stackrel{q}{\rightarrow} \calC^{\otimes} \stackrel{p}{\rightarrow}
\Nerve(\cDelta')^{op}$ exhibit $\calM = \calM^{\otimes}_{[0]}$ as left-tensored over the monoidal
$\infty$-category $\calC = \calC^{\otimes}_{[1]}$.
In what follows, we will abuse notation by identifying the object $M \in \calM$ with the
constant functor $\Nerve(\calI)^{op} \rightarrow \calM$ taking the value $M$, for a variety of
categories $\calI$.

Let $r: \calC[M]^{\otimes} \rightarrow \Nerve(\cDelta)^{op}$ denote the projection map. Using
the definition of $\calC[M]^{\otimes}$ and the description of the $r$-coCartesian edges given in the proof of Proposition \ref{umptump}, we see that the $\infty$-category
$\Alg(\calC[M])$ can be identified with the full subcategory of
$\bHom_{ \Nerve(\cDelta')^{op} }( \Nerve(\calJ)^{op}, \calM^{\otimes}) \times_{ \Fun( \Nerve(\cDelta)^{op}, \calM) } \{M\}$ spanned by those functors $f$ with the following properties:
\begin{itemize}
\item[$(a)$] Given a collection of integers $0 \leq i \leq i' \leq j' \leq j \leq n$, let $\alpha$ denote the induced map $f( [n], i \leq j) \rightarrow f( [n], i' \leq j')$. THen $q(\alpha)$ is $p$-coCartesian.
\item[$(b)$] Given a collection of integers $0 \leq i \leq j \leq n$, let
$\beta: ([n], \ast) \rightarrow ([n], i \leq j)$ be the map in $\calJ$ corresponding to the element
$j \in \{i, \ldots, j\}$. Then $f(\beta)$ is $(p \circ q)$-coCartesian. 
\item[$(c)$] Given a {\em convex} map $\gamma: [m] \rightarrow [n]$ in $\cDelta$ and a pair
of integers $0 \leq i \leq j \leq m$, the induced map
$f( [n], \alpha(i) \leq \alpha(j) ) \rightarrow f([m], i \leq j)$ is an equivalence in $\calM^{\otimes}$.
\end{itemize}

We let $X$ denote the full subcategory of
$\bHom_{ \Nerve(\cDelta')^{op} }( \Nerve(\calJ)^{op}, \calM^{\otimes})$ spanned by those functors $f$ which satisfy. We observe that if $f \in X$, then $f$ satisfies the following additional condition:

\begin{itemize}
\item[$(d_0)$] The functor $f$ carries every morphism $\alpha: ( [m], \ast) \rightarrow ([n], \ast)$ to an equivalence in $\calM$. 
\end{itemize}

To prove this, we first use the two-out-of-three property to reduce to the case where $m=0$. 
The map $\alpha$ is then classified by an element $i \in [n]$. We have a commutative diagram
in $\calJ^{op}$:
$$ \xymatrix{ ( [n], i \leq i) \ar[r]^{\beta} \ar[d]^{\gamma} & ( [0], 0 \leq 0) \ar[d]^{\gamma'} \\
( [n], \ast) \ar[r]^{\alpha} & ( [0], \ast ). }$$
Assumptions $(b)$ and $(c)$ guarantees that $f$ carries $\beta$, $\gamma$, and
$\gamma'$ to equivalences in $\calM$, so that $f$ carries $\beta'$ to an equivalence in $\calM$ as well.

Let $\calJ_{+}$ denote the category obtained from $\calJ$ by adjoining a new element
$( [-1], \ast)$, with
$$ \Hom_{\calJ_{+}}( ( [-1], \ast) , ( [n], \ast) ) = \{ \ast \} \quad \Hom_{\calJ_{+}}( ( [-1], \ast ) , ( [n], i \leq j) ) = \{ i, \ldots, j \}$$
$$ \Hom_{\calJ_{+}}(u, ( [-1], \ast) ) = \begin{cases} \{ \ast \} & \text{if } u = ( [-1], \ast) \\
\emptyset & \text{otherwise.} \end{cases}$$ 
The functor $\psi': \calJ \rightarrow \cDelta'$ extends canonically to a functor
$\psi'_{+}: \calJ_{+} \rightarrow \cDelta'$, with $\psi'_{+}( [-1], \ast) = [0]$.
Let $\widetilde{X}$ denote the full subcategory of 
$\bHom_{ \Nerve(\cDelta')^{op} }( \Nerve(\calJ_{+})^{op}, \calM^{\otimes})$ spanned by those functors $f$ which satisfy $(a)$, $(b)$, $(c)$, and the following stronger form of
$(d_0)$:

\begin{itemize}
\item[$(d)$] The functor $f$ carries every morphism $\alpha: ( [m], \ast) \rightarrow ([n], \ast)$ in $\calJ_{+}$ to an equivalence in $\calM$. 
\end{itemize}
 
Recall that we have identified $\cDelta$ with the full subcategory of
$\calJ$ spanned by the objects $([n], \ast)$. We observe that the inclusion
$\cDelta \simeq \cDelta \times_{ \calJ_{+} } ( \calJ_{+})_{ ([-1], \ast)/}
\subseteq \calJ \times_{ \calJ_{+} } (\calJ_{+})_{ ([-1], \ast)/}$ has a right adjoint, given by the formula
$$ ( \alpha: ( [-1], \ast) \rightarrow ( [n], \bigdot) ) \mapsto ( [n], \ast). $$
Combining this observation with $(d)$, Lemma \toposref{kan2}, Corollary \toposref{silt}, and
the observation that $\Nerve(\cDelta')^{op}$ is weakly contractible, we deduce:

\begin{itemize}
\item[$(i)$] Every functor $f_0 \in X$ admits a $(p \circ q)$-left Kan extension
$f \in \bHom_{ \Nerve(\cDelta')^{op} }( \Nerve(\calJ_{+})^{op}, \calM^{\otimes})$.
\item[$(ii)$] Let $f \in \bHom_{ \Nerve(\cDelta')^{op} }( \Nerve(\calJ_{+})^{op}, \calM^{\otimes})$ be arbitrary. Then $f \in \widetilde{X}$ if and only if $f_0 = f | \Nerve(\calJ)^{op}$ belongs to $X$ and $f$ is a $(p \circ q)$-left Kan extension of $f_0$.
\end{itemize}

Recall that we have identified $\cDelta'$ with the full subcategory of $\calJ$ spanned by the objects
$\{ ( [n], 0 \leq n) \}_{ n \geq 0}$. Similarly, we will category $\cDelta^{\times}$ of
Notation \ref{hugegrin} with the full subcategory of $\calJ$ spanned by the objects
$\{ ( [n], i \leq j) \}_{0 \leq i \leq j \leq n}$. Finally, let $\calJ'$ denote the full subcategory
of $\calJ_{+}$ spanned by $\cDelta^{\times}$ together with the object $( [-1], \ast)$, so that we have a sequence of inclusions
$$ \cDelta' \subseteq \cDelta^{\times} \subseteq \calJ' \subseteq \calJ_{+}.$$

We make the following observations:
\begin{itemize}
\item[$(1)$] For every object $( [n], i \leq j) \in \cDelta'$, the category
$\cDelta' \times_{ \cDelta^{\times} } ( \cDelta^{\times})_{/ ([n], i \leq j) }$ has
a final object, given by the map $( \{ i, \ldots, j\}, i \leq j) \rightarrow ( [n], i \leq j)$.

\item[$(2)$] The category $\cDelta^{\times} \times_{ \calJ'} \calJ'_{ ([-1], \ast)/}$ has an initial object, given by the map $( [-1], \ast) \rightarrow ( [0], 0 \leq 0)$. 

\item[$(3)$] For every object $([n], \ast) \in \calJ_{+}$, the category
$\calJ' \times_{ \calJ_{+} } \calJ_{ / ([n], \ast) }$ has a unique (final) object, given by
$( [-1], \ast) \rightarrow ([n], \ast)$. 
\end{itemize}

Applying Lemma \toposref{kan2}, we deduce:
\begin{itemize}
\item[$(1')$] Every functor $f_0 \in \bHom_{ \Nerve(\cDelta')^{op} }( \Nerve(\cDelta')^{op}, \calM^{\otimes})$
admits a $(p \circ q)$-right Kan extension 
$$f \in \bHom_{ \Nerve(\cDelta')^{op} }( \Nerve(\cDelta^{\times})^{op}, \calM^{\otimes}).$$ Moreover, 
an arbitrary functor $f \in \bHom_{ \Nerve(\cDelta')^{op} }( \Nerve(\cDelta^{\times})^{op}, \calM^{\otimes})$ is a $(p \circ q)$-right Kan extension of the restriction $f | \Nerve(\cDelta')^{op}$ if and only if the following condition is satisfied:
\begin{itemize}
\item[$(a')$] For every object $([n], i \leq j) \in \cDelta^{\times}$, the canonical map
$f( [n], i \leq j) \rightarrow f( \{i, \ldots, j\}, i \leq j)$ is an equivalence in $\calM^{\otimes}$.
\end{itemize}
\item[$(2')$] Every functor $f_0 \in \bHom_{ \Nerve(\cDelta')^{op} }( \Nerve(\cDelta^{\times})^{op}, \calM^{\otimes})$ admits a $(p \circ q)$-left Kan extension 
$$f \in \bHom_{ \Nerve(\cDelta')^{op} }( \Nerve(\calJ')^{op}, \calM^{\otimes}).$$ Moreover, 
an arbitrary functor $f \in \bHom_{ \Nerve(\cDelta')^{op} }( \Nerve(\calJ')^{op}, \calM^{\otimes})$ is a $(p \circ q)$-left Kan extension of the restriction $f | \Nerve(\cDelta^{\times})^{op}$ if and only if the following condition is satisfied:
\begin{itemize}
\item[$(b')$] The canonical map $f( [0], 0 \leq 0) \rightarrow f( [-1], \ast)$ is an equivalence in $\calM$.
\end{itemize}

\item[$(3')$] Every functor $f_0 \in \bHom_{ \Nerve(\cDelta')^{op} }( \Nerve(\calJ')^{op}, \calM^{\otimes})$ admits a $(p \circ q)$-right Kan extension 
$$f \in \bHom_{ \Nerve(\cDelta')^{op} }( \Nerve(\calJ_{+})^{op}, \calM^{\otimes}).$$ Moreover, 
an arbitrary functor $f \in \bHom_{ \Nerve(\cDelta')^{op} }( \Nerve(\calJ_{+})^{op}, \calM^{\otimes})$ is a $(p \circ q)$-right Kan extension of the restriction $f | \Nerve(\calJ')^{op}$ if and only if the following condition is satisfied:
\begin{itemize}
\item[$(c')$] For each $n \geq 0$, the canonical map $f( [n], \ast) \rightarrow f( [-1], \ast)$ is
an equivalence in $\calM$.
\end{itemize}
\end{itemize}

Let $Y$ denote the full subcategory of $\bHom_{ \Nerve(\cDelta')^{op} }( \Nerve(\calJ_{+})^{op}, \calM^{\otimes})$ spanned by those functors $f$ which satisfy $(a')$, $(b')$, and $(c')$, and the following additional condition:
\begin{itemize}
\item[$(d')$] The restriction $f | \Nerve(\cDelta')^{op}$ belongs to $\Mod(\calM)$.
\end{itemize}

Applying Proposition \toposref{lklk} repeatedly, we deduce that the restriction map $Y \rightarrow \Mod(\calM)$ is a trivial Kan fibration. Our next goal is to show that $Y = \widetilde{X}$. In other words, if we fix a functor $f \in \bHom_{ \Nerve(\cDelta')^{op} }( \Nerve(\calJ_{+})^{op}, \calM^{\otimes}),$ then $f$ satisfies $(a)$, $(b)$, $(c)$, and $(d)$ if and only if $f$ satisfies $(a')$, $(b')$, $(c')$ and $(d')$.

Suppose first that $f$ satisfies $(a)$, $(b)$, $(c)$, and $(d)$. Condition $(a')$ follows immediately from $(c)$, and conditions $(b')$ and $(c')$ follow from $(d)$. To verify $(d')$, we first observe that
$(a)$ and $(c)$ imply that $q \circ f|\Nerve(\cDelta')^{op}$ is an algebra object of $\calC$. It remains only to show that if $\alpha: ( [m], 0 \leq m) \rightarrow ( [n], 0 \leq n)$ satisfies $\alpha(m) = n$, then the induced map $f( [n], 0 \leq n) \rightarrow f( [m], 0 \leq m)$ is $(p \circ q)$-coCartesian. Using the equivalence $\calM^{\otimes}_{[m]} \simeq \calC^{\otimes}_{[m]} \times \calM$ and the fact that
$(q \circ f) | \Nerve(\cDelta')^{op} \in \Alg(\calC)$, we may reduce to the case $m=0$. In this case, the desired result follows from $(b)$. 

Now let us suppose that $f$ satisfies $(a')$, $(b')$, $(c')$ and $(d')$. Condition $(c)$ follows immediately from $(a')$ Condition $(d)$ follows from $(c')$ using a two-out-of-three argument. 
To prove $(b)$, we observe that the map $\beta: ( [n], \ast) \rightarrow ( [n], i \leq j)$ corresponding
to $j \in \{ i, \ldots, j\}$ fits into a commutative diagram
$$ \xymatrix{ ( [-1], \ast) \ar[d]^{\beta_0} \ar[r]^{\beta_1} & ( [0], 0 \leq 0) \ar[r]^{\beta_2} & ( \{ i, \ldots, j\}, i \leq j ) \ar[d]^{\beta_3} \\
( [n], \ast) \ar[rr]^{\beta} & & ( [n], i \leq j ). }$$
Condition $(c')$ implies that $f(\beta_0)$ and $f(\beta_1)$ are equivalences, and condition
$(a')$ guarantees that $f(\beta_3)$ is an equivalence. It follows that 
$f(\beta)$ is $(p \circ q)$-coCartesian if and only if $f(\beta_2)$ is $(p \circ q)$-coCartesian, which follows from $(d')$. The verification of $(a)$ is similar: let $0 \leq i \leq i' \leq j' \leq i \leq n$, and let
$\alpha: ( [n], i' \leq j') \rightarrow ([n], i \leq j)$ be the induced morphism in $\calJ$. Then
$\alpha$ fits into a commutative diagram
$$ \xymatrix{ ( \{i', \ldots, j' \}, i' \leq j') \ar[r]^{\alpha_1} \ar[d]^{\alpha_0} & ( \{i, \ldots, j\}, i \leq j) \ar[d]^{\alpha_2} \\
( [n], i' \leq j') \ar[r]^{\alpha} & ( [n], i \leq j). }$$
Condition $(a')$ guarantees that $f(\alpha_0)$ and $f(\alpha_2)$ are equivalences, so that
$(q \circ f)(\alpha)$ is $p$-coCartesian if and only if $(q \circ f)(\alpha_2)$ is $p$-coCartesian. Since
$\alpha_2$ is a convex morphism, this follows from assumption $(d')$.

We have a diagram of trivial Kan fibrations
$ X \leftarrow \widetilde{X} = Y \rightarrow \Mod(\calM). $
For every full subcategory $\calI \subseteq \calJ_{+}$, let $\calE(\calI)$ denote the
full subcategory of $\Fun( \Nerve(\calI)^{op}, \calM)$ spanned by those functors which
carry each morphism in $\calI$ to an equivalence in $\calM$. Let $\calJ''$ be the full subcategory of $\calJ_{+}$ spanned by the objects of $\cDelta_{+}$ and the object $([0], 0 \leq 0)$, and
$\calJ''_0$ the full subcategory spanned by the single object $( [0], 0 \leq 0)$.

Consider the diagram
$$ \Alg( \calC[M]) \simeq X \times_{ \calE(\cDelta)} \{M\} 
\stackrel{\phi_0}{\leftarrow} \widetilde{X} \times_{ \calE(\cDelta)} \{M\}
\stackrel{\phi_1}{\leftarrow} \widetilde{X} \times_{ \calE(\calJ'')} \{M\}
\stackrel{\phi_2}{\rightarrow} \widetilde{X} \times_{ \calE(\calJ''_0)} \{M\} 
\stackrel{\phi_3}{\rightarrow} \Mod(\calM) \times_{ \calM } \{M \}.$$
The above arguments imply that $\phi_0$ and $\phi_3$ are trivial Kan fibrations.
We now observe that the categories
$\calJ''$, $\cDelta$, and $\calJ''_0$ all have weakly contractible nerves (for
$\calJ''$, this follows from the observation that $([0], 0 \leq 0)$ is a final object). Consequently, the restriction maps $\calE(\cDelta) \leftarrow \calE(\calJ'') \rightarrow \calE( \calJ''_0) \simeq \calM$
are trivial Kan fibrations. Using this, we deduce that $\phi_1$ and $\phi_2$ are categorical equivalences.

The functor $\Mod(\calM) \times_{\calM} \{M\} \rightarrow \Alg( \calC[M] )$ induced by composition with $\psi'$ factors as a composition
$$ \Mod(\calM) \times_{ \calM} \{M\} \stackrel{\phi}{\rightarrow} \widetilde{X}
\times_{ \calE(\calJ'')} \{M\} 
\stackrel{\phi_1}{\rightarrow} \widetilde{X} \times_{ \calE(\cDelta)} \{M\}
\stackrel{\phi_0}{\rightarrow} \Alg(\calC[M]).$$ 
Consequently, it will suffice to show that $\phi$ is a categorical equivalence. 
We now observe that $\phi$ is a section of $\phi_3 \circ \phi_2$, hence a categorical equivalence by the two-out-of-three property.
\end{proof}

\subsection{Application: Existence of Units}\label{giddug}

Recall that a {\it nonunital ring} is an abelian group $(A,+)$ equipped with an associative multiplication
$A \otimes A \rightarrow A.$
Every associative ring determines a nonunital ring, simply by forgetting the multiplicative identity element. On the other hand, if $A$ is an associative ring, then the ring structure on $A$ is uniquely determined by underlying nonunital ring of $A$. In other words, if $A$ is a nonunital ring which admits a multiplicative identity $1$, then $1$ is uniquely determined. The proof is simple: if $1$ and $1'$ are both identities for the multiplication on $A$, then $ 1 = 1 1' = 1'. $
Our goal in this section is to prove an $\infty$-categorical analogue of this result.
More precisely, we will show that if $A$ is a nounital algebra object of a monoidal $\infty$-category $\calC$ which admits a quasi-unit, then $A$ can be extended to an algebra object of $\calC$ in an essentially unique way. 

In ordinary category theory, this is a tautology. However, in the $\infty$-categorical setting the result is not quite as obvious; the unit in an algebra object $A$ of $\calC$ is required to satisfy a hierarchy of coherence conditions with respect to the multiplication on $A$, while the definition of a quasi-unit involves only the homotopy category $\h{\calC}$. Nevertheless, we have the following result:

\begin{theorem}\label{uniqueunit}\index{quasi-unit!uniqueness of}
Let $\calC$ be a monoidal $\infty$-category. Then the restriction map
$\Alg(\calC) \rightarrow \Alg^{\nounit}(\calC)$ induces a trivial Kan fibration
$\theta: \Alg(\calC) \rightarrow \Alg^{\qunit}(\calC)$.
\end{theorem}

Before giving the proof, let us sketch the main idea. Suppose that $A$ is a nonunital ring, and we
wish to promote $A$ to an associative ring. Let $M = A$, regarded as a (nonunital) right module over itself. Left multiplication induces a homomorphism of nonunital algebras $\phi: A \rightarrow \Hom_{A}(M,M)$. If $A$ admits a left unit $1$, then $A$ is freely generated by $1$ as a right $A$-module, so that evaluation at $1$ induces an isomorphism $\Hom_{A}(M, M) \simeq M$. Under this isomorphism, $\phi$ corresponds to the map $a \mapsto a1$. If the element $1 \in A$ is also a right unit, then $\phi$ is an isomorphism. On the other hand, $\End_{A}(M, M)$ is manifestly an associative ring. To translate this sketch into the setting of higher category theory, we will need the following lemma, which will be proven at the end of this section:

\begin{lemma}\label{unique2}
Let $\calC$ be a monoidal $\infty$-category, and let $A \in \Alg^{\qunit}(\calC)$. There
exists an $\infty$-category $\calM$ which is left-tensored over $\calC$ and an object
$M \in \Mod^{\qunit}_{A}(\calM)$ which exhibits $A$ as a morphism object $\Mor_{\calM}(M,M)$.
\end{lemma}

\begin{proof}[Proof of Theorem \ref{uniqueunit}]
The map $\theta$ is evidently a categorical fibration. It will therefore suffice to show that
$\theta$ is a categorical equivalence. We first show that $\theta$ is essentially surjective. Let $A_0$ be a quasi-unital algebra object of $\calC$. According to Lemma \ref{unique2}, we can find an $\infty$-category $\calM$ which is left-tensored over $\calC$ and a quasi-unital module 
$$\overline{M}_0 \in \Mod^{\qunit}_{A_0}(\calM) \subseteq \Mod^{\nounit}(\calM)$$
which exhibits $A_0$ as a morphism object $\Mor_{\calM}(M, M)$; here we let $M$ denote the image of $\overline{M}_0$ under the forgetful functor $\Mod^{\qunit}(\calM) \rightarrow \calM$.
We have a commutative diagram
$$ \xymatrix{ \Mod(\calM) \times_{ \calM } \{ M \} \ar[r] \ar[d] & \Alg( \calC[M] ) \ar[d]^{\theta'} \ar[r] & \Alg(\calC) \ar[d]^{\theta} \\
\Mod^{\nounit}(\calM) \times_{ \calM} \{ M \} \ar[r] & \Alg^{\nounit}( \calC[M] ) \ar[r] & \Alg^{\nounit}(\calC). }$$
Let $\overline{A}_0$ be the image of $\overline{M}_0$ in
$\Alg^{\nounit}( \calC[M])$. To prove that $A_0$ belongs to the essential image of
$\theta$, it will suffice to show that $\overline{A}_0$ belongs to the essential image of
$\theta'$. 

Proposition \ref{ugher} and Remark \ref{firebaugh2} imply that
$\overline{A}_0$ is a final object of $\Alg^{\nounit}( \calC[M])$. On the other hand,
Proposition \ref{ugher} an Corollary \ref{firebaugh} imply that
$\Alg(\calC[M])$ admits a final object $\overline{A}$, and Remark \ref{firebaugh2} implies that
$\theta'( \overline{A} )$ is a final object of $\Alg^{\nounit}( \calC[M])$. It follows that
$\theta'(\overline{A}) \simeq \overline{A}_0$, so that $\overline{A}_0$ belongs to the essential image of $\theta'$ as desired.

We now prove that $\theta$ is fully faithful. Fix objects $A,B \in \Alg(\calC)$, and set
$A_0 = \theta(A) \in \Alg^{\qunit}(\calC), B_0 = \theta(B) \in \Alg^{\qunit}(\calC)$.
We wish to prove that the map
$\bHom_{\Alg(\calC)}( B,A) \rightarrow \bHom_{ \Alg^{\qunit}(\calC)}(B_0, A_0)$ is a homotopy equivalence. Use Lemma \ref{unique2} to choose $\overline{M}_0 \in \Mod_{A_0}^{\qunit}( \calM)$ as above. Consider the diagram
$$ \xymatrix{ \Mod_{A}(\calM) \ar[d] & \Mod_{A}(\calM) \times_{\calM} \{M \} \ar[l] \ar[r] \ar[d] & \{A\} \times_{\Alg(\calC) } \Alg(\calC[M]) \ar[d]^{\psi} \\
\Mod^{\qunit}_{A_0}(\calM) & \Mod^{\qunit}_{A_0}(\calM) \times_{\calM} \{ M \} \ar[l] \ar[r] & 
\{A_0 \} \times_{\Alg^{\qunit}(\calC)} \Alg^{\qunit}( \calC[M]). }$$
The left square is a pullback, and the left vertical map is a trivial Kan fibration (Proposition \ref{uniquemo}). The horizontal maps on the right are both categorical equivalences (Proposition \ref{poofer}). Using the two-out-of-three property, we deduce that $\psi$ is a categorical equivalence. Since $\psi$ is also a categorical fibration, it is a trivial Kan fibration; we may therefore choose
$\overline{M} \in \Mod_{A}(\calM)$ lifting $\overline{M}_0$.

According to Corollary \ref{usebilly} and Remark \ref{usebully}, we have canonical homotopy equivalences
$$ \phi: \bHom_{ \Alg(\calC)}(B,A) \simeq \Mod_{B}(\calM) \times_{ \calM } \{M \}$$
$$ \phi^{\nounit}: \bHom_{ \Alg^{\nounit}(\calC)}( B_0, A_0) \simeq \Mod^{\nounit}_{B_0}(\calM) \times_{\calM} \{M \}.$$ 
Let $f: B_0 \rightarrow A_0$ be a map of nonunital algebras, and let
$\overline{N}_0$ be the corresponding object of $\Mod^{\nounit}_{B_0}(\calM) \times_{ \calM } \{M\}$. Then $f$ is quasi-unital if and only if the composition
$u: 1_{\calC} \rightarrow B \stackrel{f}{\rightarrow} A$
is homotopic to the unit of $A$; here we abuse notation by identifying $A$ and $B$ with their images in $\calC$. Using the equivalence $A \simeq \Mor_{\calM}(M,M)$, we can identify 
$u \in \pi_0 \bHom_{\calC}( 1_{\calC}, \Mor_{\calM}(M,M) )$ with a point
$u' \in \pi_0 \bHom_{\calM}( M, M)$; then $f$ is quasi-unital if and only if $u'$ is homotopic to the identity. The map $u'$ can be identified with the composition
$M \simeq 1_{\calC} \otimes M \rightarrow B \otimes M \rightarrow M,$
so that $f$ is quasi-unital if and only if $\overline{N}_0$ is quasi-unital. It follows that
$\phi^{\nounit}$ restricts to a homotopy equivalence
$$ \phi^{\qunit}: \bHom_{ \Alg^{\qunit}(\calC)}( B_0, A_0) \simeq \Mod^{\qunit}_{B_0}(\calM) \times_{\calM} \{ M \}.$$

We wish to prove that $\theta$ induces a homotopy equivalence
$\bHom_{ \Alg(\calC)}(B, A) \rightarrow \bHom_{ \Alg^{\qunit}(\calC)}( B_0, A_0 )$. In view of the above identifications, it will suffice to show that the restriction map
$$g: \Mod_{B}(\calM) \times_{\calM} \{ M \} \rightarrow \Mod_{B_0}^{\qunit}(\calM) \times_{ \calM } \{M\}$$
is a homotopy equivalence. Proposition \ref{uniquemo} implies that $g$
is a trivial Kan fibration.
\end{proof}

\begin{definition}\label{spreck2}\index{ZZZModqunitM@$\Mod^{\qunit}(\calM)$}
Let $\calC$ be a monoidal $\infty$-category, let $\calM$ be an $\infty$-category which is left-tensored over $\calC$. We let $\Mod^{\qunit}(\calM)$ denote the full subcategory of
$\Mod^{\nounit}(\calM) \times_{ \Alg^{\nounit}(\calC) } \Alg^{\qunit}(\calC)$
spanned by those nonunital $A$-modules $M$ satisfying the following condition:
\begin{itemize}
\item[$(\ast)$] Let $u: 1_{\calC} \rightarrow A$ be a quasi-unit for $A$. Then the composition
$$ M \simeq 1_{\calC} \otimes M \stackrel{u}{\rightarrow} A \otimes M \rightarrow M$$
is an equivalence.
\end{itemize}
\end{definition}

\begin{remark}
Let $\calC$ be a monoidal $\infty$-category, let $\calM$ be an $\infty$-category which is left-tensored over $\calC$, let $A$ be an algebra object of $\calC$ and let $A_0$ denote the 
underlying nonunital algebra. Then the fiber
$$\Mod^{\qunit}(\calM) \times_{ \Alg^{\qunit}(\calC) } \{A_0 \} \subseteq \Mod^{\nounit}_{A}(\calM)$$ 
coincides with the full subcategory $\Mod^{\qunit}_{A}(\calM) \subseteq \Mod^{\nounit}_{A}(\calM)$
introduced in Definition \ref{spreck2}.
\end{remark}

\begin{corollary}\label{uniquemod}
Let $\calC$ be a monoidal $\infty$-category, let $\calM$ be an $\infty$-category which is left-tensored over $\calC$. Then:
\begin{itemize}
\item[$(1)$] The restriction functor
$ \theta: \Mod(\calM) \rightarrow \Mod^{\qunit}(\calM) \times_{ \Alg^{\qunit}( \calC) } \Alg(\calC)$
is a trivial Kan fibration.
\item[$(2)$] The restriction functor $\Mod(\calM) \rightarrow \Mod^{\qunit}(\calM)$ is a trivial Kan fibration.
\end{itemize}
\end{corollary}

\begin{proof}
It is clear that $(2)$ follows from $(1)$ and Theorem \ref{uniqueunit}. It will therefore suffice to prove $(1)$.
Consider the diagram
$$ \xymatrix{ \Mod^{\qunit}(\calM) \ar[r] \ar[d]^{p} & \Mod^{\nounit}(\calM) \ar[d]^{p'} \\
\Alg^{\qunit}(\calC) \ar[r] & \Alg^{\nounit}(\calC). }$$
According to Remark \ref{thetacart2}, the map $p'$ is a Cartesian fibration; moreover, a morphism
$f: M \rightarrow M'$ in $\Mod^{\nounit}(\calM)$ is $p'$-Cartesian if and only if it induces an equivalence in $\calM$. Suppose that $f$ is the image of a morphism $\overline{f}: \overline{M} \rightarrow \overline{M}'$ in $\Mod^{\nounit}(\calM) \times_{ \Alg^{\qunit}(\calC) } \Alg^{\nounit}(\calC)$. If $f$ is $p'$-Cartesian, we deduce easily that $\overline{M} \in \Mod^{\qunit}(\calM)$ if and only if $\overline{M}' \in \Mod^{\qunit}(\calM)$. It follows that $p$ is also a Cartesian fibration, and that a morphism $\overline{f}: \overline{M} \rightarrow \overline{M}'$ in $\Mod^{\qunit}(\calM)$ is 
$p$-Cartesian if and only if it induces an equivalence in $\calM$.

Consider next the diagram
$$ \xymatrix{ \Mod(\calM) \ar[d]^{q} \ar[r]^-{\theta'} & \Mod^{\qunit}(\calM)
\times_{ \Alg^{\qunit}(\calC)} \Alg(\calC) \ar[d]^{q'} \\
\Alg(\calC) \ar@{=}[r] & \Alg(\calC). }$$
The map $q'$ is a pullback of $p$. It follows that $q'$ is a Cartesian fibration, and that a morphism of $\Mod^{\qunit}(\calM) \times_{ \Alg^{\qunit}(\calC)} \Alg(\calC)$ is $q'$-Cartesian if and only if its image in $\calM$ is an equivalence. Corollary \ref{thetacart} implies that $q$ is also a Cartesian fibration, and that $\theta'$ carries $q$-Cartesian edges to $q'$-Cartesian edges. According to Corollary \toposref{usefir}, $\theta$ is a categorical equivalence provided that it induces
a categorical equivalence of fibers $\Mod_{A}(\calM) \rightarrow \Mod^{\qunit}_{A}(\calM)$ over every object $A \in \Alg(\calC)$. The desired result now follows immediately from Proposition \ref{uniquemo}.
\end{proof}

\begin{remark}
Theorem \ref{uniqueunit} can be considered as a special case of Corollary \ref{uniquemod}.
To see this, let $q: \calC^{\otimes} \rightarrow \Nerve(\cDelta)^{op}$ be a monoidal $\infty$-category. The identity map $\calC^{\otimes} \rightarrow \calC^{\otimes}$ exhibits
$\calM= \calC^{\otimes}_{[0]}$ as left-tensored over $\calC = \calC^{\otimes}_{[1]}$. We have a commutative diagram
$$ \xymatrix{ \Mod(\calM) \ar[d] \ar[r] & \Mod^{\qunit}(\calM) \ar[r] \ar[d] & \Mod^{\nounit}(\calM) 
\ar[d] \\
\Alg(\calC) \ar[r] & \Alg^{\qunit}(\calC) \ar[r] & \Alg^{\nounit}(\calC) }$$
where the vertical maps are isomorphisms.
\end{remark}

We now return to the proof of Lemma \ref{unique2}. The idea is to take $\calM$ to be the
$\infty$-category of nonunital {\em right} $A$-modules in $\calC$, and $M \in \calM$ to be $A$ itself, regarded as a right $A$-module. 




\begin{proof}[Proof of Lemma \ref{unique2}]
Let $\cDelta^{\nounit}_{+}$ be the category obtained by adjoining an initial object
$[-1]$ to $\cDelta^{\nounit}$. Let $\psi: \cDelta \times \cDelta^{\nounit}_{+} \rightarrow \cDelta$ denote the {\em join} functor
$$ \psi( [m], [n] ) = [m] \star [n] \simeq [m+n+1],$$
and let $\psi_0 = \psi | \cDelta \times \cDelta^{\nounit}$. We observe that
there is a canonical natural transformation of functors $\alpha: \pi_{2} \rightarrow \psi_0$, where
$\pi_2: \cDelta \times \cDelta^{\nounit} \rightarrow \cDelta$ denotes projection onto the second factor (followed by the inclusion $\cDelta^{\nounit} \subseteq \cDelta$). Set
$$ T = (\Nerve(\cDelta_{+}^{\nounit})^{op} \times \{0\} )
\coprod_{ \Nerve(\cDelta^{\nounit})^{op} \times \{0\} } ( \Nerve(\cDelta^{\nounit})^{op} \times \Delta^1)$$
so that $\alpha$ determines a map of simplicial sets
$h: \Nerve(\cDelta)^{op} \times T \rightarrow \Nerve(\cDelta)^{op}.$

Let $p: \calC^{\otimes} \rightarrow \Nerve(\cDelta)^{op}$ exhibit $\calC = \calC^{\otimes}_{[1]}$ as a monoidal $\infty$-category.
We define a simplicial set $\widetilde{\calM}^{\otimes}$ equipped with a map
$\widetilde{\calM}^{\otimes} \rightarrow \Nerve(\cDelta)^{op}$ by the following universal property:
for every simplicial set $K$ equipped with a map $K \rightarrow \Nerve(\cDelta)^{op}$, we have
a canonical bijection of $\Hom_{ \Nerve(\cDelta)^{op} }( K, \widetilde{\calM}^{\otimes} )$ with
the set of maps $f: K \times T \rightarrow \calC^{\otimes}$ for which the diagram
$$ \xymatrix{ 
K \times ( \Nerve( \cDelta^{\nounit})^{op} \times \{1\} ) \ar[r] \ar[d] & \Nerve(\cDelta^{\nounit})^{op} \ar[d]^{A} \\ 
K \times T \ar[d] \ar[r]^{f} & \calC^{\otimes} \ar[d]^{p} \\
\Nerve(\cDelta)^{op} \times T \ar[r]^{h} & \Nerve(\cDelta)^{op}}$$
is commutative. 

A vertex of $\widetilde{\calM}^{\otimes}$ classifies a map of simplicial sets
$f: T \rightarrow \calC^{\otimes}$. Let $\calM^{\otimes}$ be the full simplicial subset
of $\widetilde{\calM}^{\otimes}$ spanned by those vertices corresponding to functors $f$ which satisfy the following additional conditions:
\begin{itemize}
\item[$(i)$] For every object $[n] \in \cDelta^{\nounit}$, the restriction 
$f| \{ [n] \} \times \Delta^1$ is a $p$-coCartesian edge of $\calC^{\otimes}$.
\item[$(ii)$] For every convex morphism $\alpha: [m] \rightarrow [n]$ in
$\cDelta^{\nounit}$ such that $\alpha(0) = 0$, the morphism $f(\alpha \times \{0\})$ is $p$-coCartesian.
\item[$(iii)$] For every object $[n] \in \cDelta^{\nounit}$, if $\alpha: [-1] \rightarrow [n]$ denotes the
unique morphism in $\cDelta^{\nounit}_{+}$, then $f( \alpha \times \{0\} )$ is $q$-coCartesian.
\end{itemize}

It is not difficult to see that the projection $\calM^{\otimes} \rightarrow \Nerve(\cDelta)^{op}$ is a coCartesian fibration. Restriction to the object $\{ [-1] \} \times \{ 0 \}$ determines a map
$q: \calM^{\otimes} \rightarrow \calC^{\otimes}$, which exhibits $\calM = \calM^{\otimes}_{[0]}$ as left-tensored over $\calC$. Restriction to the simplicial subset
$\Nerve(\cDelta^{\nounit})^{op} \times \Delta^1 \subseteq T$ determines a trivial Kan fibration
$\calM \rightarrow \Mod^{\nounit, R}_{A}( \calC)$. Finally, we observe that the composition
$\Nerve(\cDelta^{\nounit})^{op} \times T
\stackrel{h}{\rightarrow} \Nerve( \cDelta^{\nounit})^{op}
\stackrel{A}{\rightarrow} \calC^{\otimes} $
determines a quasiunital left $A$-module $M \in \Mod_{A}^{\nounit}(\calM)$. 

In what follows, we will abuse notation by identifying the nonunital algebra $A$ and the
$A$-module $M$ with their images in $\calC$ and $\calM$, respectively. 
To complete the proof, it will suffice to show that the (left) action of $A$ on $M$
$\theta: A \otimes M \rightarrow M$
exhibits $A$ as a morphism object $\Mor_{\calM}(M,M)$. In other words, we must show that for every object $C \in \calC$, composition with $\theta$ induces a homotopy equivalence
$\bHom_{\calC}(C, A) \rightarrow \bHom_{\calM}( C \otimes M, M).$

We must show that for every Kan complex $K$, $\theta$ induces a bijection
$$ [ K, \bHom_{\calC}(C,A) ] \rightarrow [K, \bHom_{\calM}( C \otimes M, M)];$$
here $[K,X]$ denotes the set of maps from $K$ to $X$ in the homotopy category $\calH$ of spaces.
Replacing $\calC$ by $\Fun(K,\calC)$ and $A$ by the nonunital algebra
$A' \in \Alg^{\nounit}( \Fun(K, \calC)) \simeq \Fun(K, \Alg^{\nounit}(\calC) )$ corresponding
to the constant map $K \rightarrow \{A\} \subseteq \Alg^{\nounit}(\calC)$, we can reduce to the case where $K = \emptyset$. In other words, it will suffice to show composition with $\theta$ induces a bijection
$q: \pi_0 \bHom_{\calC}(C,A) \rightarrow \pi_0 \bHom_{\calM}(C \otimes M, M).$

Our next step is to construct an inverse to $q$. Let $u: 1_{\calC} \rightarrow A$ be a quasi-unit.
Let $\phi: C \otimes M \rightarrow M$ be an arbitrary morphism in $\calM$. Then
$\phi$ determines a map $C \otimes A \rightarrow A$ in $\calC$. Let $q'(\phi)$ denote the composition
$$ C \simeq C \otimes 1_{\calC} \stackrel{u}{\rightarrow} C \otimes A \rightarrow A.$$
We may view $q'$ as a map of sets from
$\pi_0 \bHom_{\calM}(C \otimes M, M)$ to
$\pi_0 \bHom_{\calC}(C,A).$
The composition $q' \circ q: \pi_0 \bHom_{\calC}(C,A) \rightarrow \pi_0 \bHom_{\calC}(C,A)$ is induced by the map $r_u: A \rightarrow A$ given by right multiplication by $u$.
Since $u$ is a {\em right} unit of $A$, we deduce that $q' \circ q$ is the identity. In particular, $q$ is injective.

To complete the proof, it will suffice to show that $q$ is surjective. For this, we use the results of \S \ref{cupper}. According to Proposition \ref{ugher}, if $C$ is an object of $\calC$, then giving a map
$\phi: C \otimes M \rightarrow M$ is equivalent to lifting $C$ to an object $\widetilde{C} \in \calC[M]$. In particular, the left action of $A$ on $M$ gives rise to a canonical element
$\widetilde{A} \in \calC[M]$. Then $\phi$ belongs to the image of $q$ if and only if
there exists a map $\widetilde{C} \rightarrow \widetilde{A}$ in $\calC[M]$. Consequently, the surjectivity of $q$ is equivalent to the following assertion:

\begin{itemize}
\item[$(\ast)$] For every object $\widetilde{C} \in \calC[M]$, there exists a morphism
$\widetilde{C} \rightarrow \widetilde{A}$ in $\calC[M]$.
\end{itemize}

Proposition \ref{ugher} asserts that the projection $\calC[M] \rightarrow \calC$ is a right fibration. Consequently, the quasi-unit $u: 1_{\calC} \rightarrow A$ can be lifted to a map
$\widetilde{u}: \widetilde{1}_{\calC} \rightarrow \widetilde{\calA}$ in $\calC[M]$. Since $u$ is a {\em left} unit of $A$, the object $\widetilde{1}_{\calC}$ classifies a map
$v: 1_{\calC} \otimes M \rightarrow M$
which determines an equivalence in $\calC$. Using Remark \ref{thetacart2}, we conclude that
$v$ is an equivalence in $\calM$. It follows from Remark \ref{girlytime} that $\widetilde{1}_{\calC}$ is an invertible object of $\calC[M]$, so that the functor $\widetilde{\calC} \mapsto \widetilde{\calC} \otimes \widetilde{1}_{\calC}$ is an equivalence from $\calC[M]$ to itself. Consequently, condition $(\ast)$ is equivalent to:

\begin{itemize}
\item[$(\ast')$] For every object $\widetilde{C} \in \calC[M]$, there exists a morphism
$\widetilde{C} \otimes \widetilde{1}_{\calC} \rightarrow \widetilde{A}$ in $\calC[M]$.
\end{itemize}

In view of the existence of $\widetilde{u}: \widetilde{1}_{\calC} \rightarrow \widetilde{A}$, it will suffice to prove the following slightly stronger assertion:

\begin{itemize}
\item[$(\ast'')$] For every object $\widetilde{C} \in \calC[M]$, there exists a morphism
$\widetilde{C} \otimes \widetilde{A} \rightarrow \widetilde{A}$ in $\calC[M]$.
\end{itemize}

Applying Proposition \ref{ugher} again, we see that $(\ast'')$ is equivalent to the following assertion: for every map $\phi: C \otimes M \rightarrow M$, there exists a commutative diagram
$$ \xymatrix{ C \otimes A \otimes M \ar[d] \ar[r] & C \otimes M \ar[d]^{\phi} \\
A \otimes M \ar[r] & M; }$$
in $\calM$, where the horizontal arrows are given by the canonical action of $A$ on $M$.
This is a straightforward consequence of our construction of $\calM$ and $M$; we leave the details to the reader.

\end{proof}

\section{Monads and the Barr-Beck Theorem}\label{barri}

Suppose given a pair of adjoint functors
$ \Adjoint{F}{\calC}{\calD}{G}$
between ordinary categories. Then:
\begin{itemize}\index{monad}
\item[$(A)$] The composition $T = G \circ F$ has the structure of a {\it monad} on $\calC$; that is,
an algebra object of the category $\Fun(\calC, \calC)$ of endofunctors of $\calC$. Here the unit map 
$\id_{\calC} \rightarrow T$ given by the unit of the adjunction between $F$ and $G$, and the product is given by the composition
$$ T \circ T = G \circ (F \circ G) \circ F \rightarrow G \circ \id_{\calD} \circ F = T$$
where the second map is given by a compatible counit $v$ for the adjunction between $F$ and $G$.

\item[$(B)$] For every object $D \in \calD$, the object $G(D)$ has the structure of a module over the monad $T$, given by the map $T G(D) = ((G \circ F) \circ G)(D) = (G \circ (F \circ G))(D)
\stackrel{v}{\rightarrow} G(D)$. This construction determines a functor $\theta$ from $\calD$ to the category of $T$-modules in $\calC$.

\item[$(C)$] In many cases, the functor $\theta$ is an equivalence of categories. The {\it Barr-Beck theorem} provides necessary and sufficient conditions on the functor $G$ to guarantee that this is the case. We refer the reader to \cite{maclane} for a detailed statement (or to Theorem \ref{barbeq} for our $\infty$-categorical version, which subsumes the classical statement).\index{Barr-Beck theorem}
\end{itemize}

Our goal in this section is to obtain $\infty$-categorical generalizations of assertions $(A)$ through $(C)$. Our first step is to define the notion of a {\it monad} on an $\infty$-category $\calC$. We will do so by introducing a monoidal structure on the $\infty$-category $\Fun(\calC, \calC)$, which is determined by the (strictly) associative composition product on $\Fun(\calC, \calC)$. The details of this construction are outlined in \S \ref{moadsec}, using a relative version of the nerve construction which is of some independent interest. We can then define a {\it monad} on $\calC$ to be an algebra object of $\Fun(\calC, \calC)$.

The next step is to establish an analogue of assertion $(A)$.
Suppose given a pair of adjoint functors $\Adjoint{F}{\calC}{\calD}{G}$ between $\infty$-categories. The composition $T = G \circ F$ is equipped with natural transformations $\id \rightarrow T$, $T \circ T \rightarrow T$, given by composition with unit and counit maps for the adjunction between $\calC$ and $\calD$. Moreover, the classical proof of the associative law can be carried through unchanged in the homotopy category, so that $T$ can be viewed as an algebra object of the monoidal category $\h{ \Fun(\calC, \calC)}$. However, endowing $T$ with the structure of a monad on $\calC$ is more difficult. To carry it out, we need to construct an elaboration of the monoidal structure on $\Fun(\calC, \calC)$, which involves not only $\calC$ but also $\calD$. 
We will outline this construction in \S \ref{moadthird}.
Using this construction, we will define an {\it adjunction datum}: roughly speaking, this is a structure relating the $\infty$-categories $\calC$ and $\calD$ that simultaneously determines a pair of adjoint functors $\Adjoint{F}{\calC}{\calD}{G}$ and a monad $T$ on $\calC$, where $T$ is equivalent to
$G \circ F$. Our main result, Theorem \ref{partA}, asserts that every functor $F: \calC \rightarrow \calD$ which admits a right adjoint can be extended, in an essentially unique way, to an adjunction datum. Surprisingly, this turns out to be quite a bit harder to prove than the Barr-Beck theorem itself; we will therefore postpone the proof until \S \ref{digtwo}.

Once we have extracted the relevant monad $T \in \Alg( \Fun(\calC, \calC) )$, we can consider the $\infty$-category of $T$-modules $\Mod_{T}(\calC)$. In \S \ref{datpoint}, we will show that
there is a canonical functor $\calD \rightarrow \Mod_{T}(\calC)$, which is well-defined up to contractible ambiguity. More precisely, we will construct a diagram
$$ \calD \stackrel{\phi}{\leftarrow} \overline{\AdjDiag}_{U}(\calC, \calD) \stackrel{\psi}{\rightarrow}
\Mod_{T}(\calC),$$
and prove that $\phi$ is a trivial Kan fibration (Corollary \ref{hooku}). This is our $\infty$-categorical analogue of $(B)$.

In \S \ref{monoid7}, we will formulate and prove an $\infty$-categorical version of the Barr-Beck theorem (Theorem \ref{barbeq}), which gives a necessary and sufficient condition for the functor $\psi$ appearing above to be an equivalence of $\infty$-categories. The same argument will be used to prove a very useful assertion concerning the existence of simplicial resolutions (Proposition \ref{littlebeck}), whose formulation does not require the theory of monads.

\subsection{$\infty$-Categories of Endofunctors}\label{moadsec}

Let $\calM$ be an $\infty$-category. In \S \ref{hugr2.0}, we introduced the definition of a monoidal $\infty$-category $\calC$ equipped with a left action of $\calC$ on $\calM$. In this section, we will construct the universal example of such a left action, given by $\calC = \Fun(\calM, \calM)$. It is clear that $\calC$ is a monoid object in the ordinary category of simplicial sets, and that $\calM$ is equipped with a (left) action of $\calC$ in the classical sense. Our first goal is to convert this data into a left action of a monoidal $\infty$-category, in the sense of Definition \ref{ulult}. 

We begin with a few generalities. Let $\calI$ be an ordinary category, and let $f: \Nerve(\calI) \rightarrow \Cat_{\infty}^{op}$ be a diagram. We can associate to $f$ a Cartesian fibration $X \rightarrow \Nerve(\calI)$. Let us briefly review the construction (see \S \toposref{universalfib}). We first identify $f$ with a simplicial functor $F: \sCoNerve[ \Nerve(\calI) ]^{op} \rightarrow \mSet$, where $\mSet$ denotes the category of marked simplicial sets (\S \toposref{twuf}). The functor $F$ is a weakly fibrant object of
$(\mSet)^{\sCoNerve[ \Nerve(\calI) ]^{op}}$, so that after applying the unstraightening functor
$\Un^{+}_{\Nerve(\calI)}$ we obtain a fibrant object of $(\mSet)_{/\Nerve(\calI)}$, which we can identify with the desired Cartesian fibration $p: X \rightarrow \Nerve(\calI)$. For many purposes, this construction is unnecessarily complicated. For example, the fiber of $p$ over an object $I \in \calI$ is equivalent to $f(I)$, but not isomorphic to $f(I)$. In the case where $f$ arises as the nerve of a functor $\calI \rightarrow \Set_{\Delta}^{op}$, there is an equivalent construction which is quite a bit simpler. We will present this construction in a dual form (since our primary interest is in {\em coCartesian} fibrations).

\begin{definition}\label{sulke}\index{relative nerve}\index{nerve!relative}\index{ZZZNervefI@$\Nerve_{f}(\calI)$}
Let $\calI$ be a category, and let $f: \calI \rightarrow \sSet$ be a functor. We define a new simplicial
set $\Nerve_{f}(\calI)$, the {\it nerve of $\calI$ relative to $f$}, as follows. For every finite linearly ordered set $J$, a map $\Delta^J \rightarrow \Nerve_{f}(\calI)$ consists of the following data:
\begin{itemize}
\item[$(1)$] A functor $\sigma$ from $J$ to $\calI$. 
\item[$(2)$] For every nonempty subset $J' \subseteq J$ having a maximal element $j'$, a
map $\tau(J'): \Delta^{J'} \rightarrow f( \sigma(j'))$. 
\item[$(3)$] For nonempty subsets $J'' \subseteq J' \subseteq J$, with maximal elements $j'' \in J''$, $j' \in J'$, the diagram $$ \xymatrix{ \Delta^{J''} \ar[r]^{\tau(J'')} \ar@{^{(}->}[d] & f( \sigma(j'') ) \ar[d] \\
\Delta^{J'} \ar[r]^{\tau(J')} & f( \sigma(j') ) }$$
commutes.
\end{itemize}
\end{definition}

\begin{remark}
The simplicial set $\Nerve_{f}(\calI)$ of Definition \ref{sulke} depends functorially on $f$. When
$f$ takes the constant value $\Delta^0$, there is a canonical isomorphism
$\Nerve_{f}(\calI) \simeq \Nerve(\calI)$. In particular, for {\em any} value of $f$, there is a canonical
map $\Nerve_{f}(\calI) \rightarrow \Nerve(\calI)$; the fiber of this map over an object $I \in \calI$ can be identified with the simpicial set $f(I)$. 
\end{remark}

\begin{remark}
Let $\calI$ be denote the linearly ordered set $[n]$, regarded as a category, and let
$f: \calI \rightarrow \sSet$ correspond to a composable sequence of morphisms
$$ \phi: X_0 \rightarrow \ldots \rightarrow X_n.$$
Then $\Nerve_{f}(\calI)$ is closely related to the mapping simplex $M^{op}(\phi)$ (see \S \toposref{funkystructure}). More specifically, there is a canonical map $\Nerve_{f}(\calI) \rightarrow M^{op}(\phi)$ compatible with the projection to $\Delta^{n}$, which induces an isomorphism on each fiber.
\end{remark}

\begin{lemma}\label{sulken2}
Let $\calI$ be a category and
let $\alpha: f \rightarrow f'$ be a natural transformation of functors
$f,f': \calI \rightarrow \sSet$.
\begin{itemize}
\item[$(1)$] Suppose that, for each $I \in \calI$, the map $\alpha(I): f(I) \rightarrow f'(I)$ is
an inner fibration of simplicial sets. Then the induced map $\Nerve_{f}(\calI) \rightarrow
\Nerve_{f'}(\calI)$ is an inner fibration.
\item[$(2)$] Suppose that, for each $I \in \calI$, the simplicial set $f(I)$ is an $\infty$-category. Then $\Nerve_{f}(\calI)$ is an $\infty$-category.
\item[$(3)$] Suppose that, for each $I \in \calI$, the map $\alpha(I): f(I) \rightarrow f'(I)$ is
a categorical fibration of $\infty$-categories. Then the induced map
$\Nerve_{f}(\calI) \rightarrow \Nerve_{f'}(\calI)$ is a categorical fibration of $\infty$-categories.
\end{itemize}
\end{lemma}

\begin{proof}
Consider a commutative diagram
$$ \xymatrix{ \Lambda^n_i \ar[r] \ar@{^{(}->}[d] & \Nerve_{f}(\calI) \ar[d]^{p} \\
\Delta^n \ar@{-->}[ur] \ar[r] & \Nerve_{f'}(\calI), }$$
and let $I$ be the image of $\{n\} \subseteq \Delta^n$ under the
bottom map. If $0 \leq i < n$, then the lifting problem depicted in the diagram above is equivalent to the existence of a dotted arrow in an associated diagram
$$ \xymatrix{ \Lambda^n_i \ar[r]^{g} \ar@{^{(}->}[d] & f(I) \ar[d]^{\alpha(I)} \\
\Delta^n \ar@{-->}[ur] \ar[r] & f'(I).}$$
If $\alpha(I)$ is an inner fibration and $0 < i < n$, then we conclude that this lifting problem admits a solution. This proves $(1)$. 

To prove $(2)$, we apply $(1)$ in the special case where $f'$ is the constant functor taking the value $\Delta^0$. It follows that $\Nerve_{f}(\calI) \rightarrow \Nerve(\calI)$ is an inner fibration, so that $\Nerve_{f}(\calI)$ is an $\infty$-category.

We now prove $(3)$. According to Corollary \toposref{gottaput}, an inner fibration
of $\calC \rightarrow \calD$ of $\infty$-categories is a categorical fibration if and only if the following condition is satisfied:
\begin{itemize}
\item[$(\ast)$] For every equivalence $e: D \rightarrow D'$ in $\calD$, and every
object $C \in \calC$ lifting $D$, there exists an equivalence $\overline{e}: C \rightarrow C'$
in $\calC$ lifting $e$.
\end{itemize}

We can identify equivalences in $\Nerve_{f'}(\calI)$ with triples
$(g: I \rightarrow I', X, e: X' \rightarrow Y)$ where $g$ is an isomorphism in $\calI$, $X$ is an object
of $f'(I)$, $X'$ is the image of $X$ in $f'(I')$, and $e: X' \rightarrow Y$ is an equivalence in
$f'(I')$. Given a lifting $\overline{X}$ of $X$ to $f(I)$, we can apply the assumption that
$\alpha(I')$ is a categorical fibration (and Corollary \toposref{gottaput}) to lift $e$ to an equivalence $\overline{e}: \overline{X}' \rightarrow \overline{Y}$ in $f(I')$. This produces the desired equivalence
$(g: I \rightarrow I', \overline{X}, \overline{e}: \overline{X}' \rightarrow \overline{Y})$ in
$\Nerve_{f}(\calI)$.
\end{proof}

\begin{proposition}\label{sulken}
Let $\calI$ be a category, and let $f: \calI \rightarrow \sSet$ be a functor such that
$f(I)$ is an $\infty$-category for each $I \in \calI$. Then:
\begin{itemize}
\item[$(1)$] The projection $p: \Nerve_{f}(\calI) \rightarrow \Nerve(\calI)$ is a coCartesian fibration of simplicial sets. 

\item[$(2)$] Let $e$ be an edge of $\Nerve_{f}(\calI)$, covering a morphism $I \rightarrow I'$
in $\calI$. Then $e$ is $p$-coCartesian if and only if the corresponding edge of $f(I')$ is an equivalence.

\item[$(3)$] The coCartesian fibration $p$ is associated to the functor $\Nerve(f): \Nerve(\calI) \rightarrow \Cat_{\infty}$ $($see \S \toposref{universalfib}$)$. 
\end{itemize}
\end{proposition}

\begin{proof}
Lemma \ref{sulken2} implies that $p$ is an inner fibration.
Consider a commutative diagram
$$ \xymatrix{ \Lambda^n_0 \ar[r] \ar@{^{(}->}[d] & \Nerve_{f}(\calI) \ar[d]^{p} \\
\Delta^n \ar@{-->}[ur] \ar[r] & \Nerve(\calI), }$$
and let $I$ be the image of $\{n\} \subseteq \Delta^n$ under the
bottom map. Then the lifting problem depicted in the diagram above is equivalent to the existence of a dotted arrow in an associated diagram
$$ \xymatrix{ \Lambda^n_0 \ar[r]^{g} \ar@{^{(}->}[d] & f(I) \\
\Delta^n. \ar@{-->}[ur] & }$$
If $n > 1$, an extension exists provided
that $g$ carries the initial edge of $\Lambda^n_0$ to an equivalence in $f(I)$. This proves
the ``if'' direction of $(2)$.

We next observe that for every morphism $h: I \rightarrow I'$ in $\calI$ and every object $x \in f(I)$, there exists morphism $\overline{h}: x \rightarrow x'$ in $\Nerve_{f}(\calI)$ which lifts
$h$ and classifies an equivalence in $f(I')$; in fact, we can choose $\overline{h}$ to correspond to the identity map from $h_{!}(x) \in f(I')$ to itself. The above argument shows that 
$\overline{h}$ is $p$-coCartesian. This completes the proof of $(1)$. The ``only if'' direction of $(2)$ now follows from the fact that $p$-coCartesian lifts of morphisms in $\Nerve(\calI)$ are unique up to equivalence.

To prove $(3)$, we use the formalism of marked simplicial sets (see \S \toposref{twuf}). 
Let $f^{op}: \Nerve(\calI) \rightarrow \mSet$ denote the functor given by the formula 
$f^{op}(I) = ( f(I)^{op} )^{\natural}$, and let $Z$ be the result of applying the unstraightening functor 
$\Un^{+}_{ \Nerve(\calI)^{op}}$ to $f^{op}$. Then $Z$ is a fibrant object of $(\mSet)_{/\Nerve(\calI)^{op}}$, and is therefore of the form $(X^{op})^{\natural}$, where
$q: X \rightarrow \Nerve(\calI)$ is a coCartesian fibration. Unwinding the definitions, we see
that there is a commutative diagram
$$ \xymatrix{ \Nerve_{f}(\calI) \ar[dr]^{p} \ar[rr]^{r} & & X \ar[dl]^{q} \\
& \Nerve(\calI) & }$$
where $r$ carries $p$-coCartesian edges to $q$-coCartesian edges. It follows from 
Theorem \toposref{straightthm} (applied over a point) that $r$ induces an equivalence
of $\infty$-categories after passing to the fiber over each object $I \in \calI$. 
The desired result now follows from Corollary \toposref{usefir}.
\end{proof}

\begin{notation}
Let $\calM$ be a simplicial set. We define functors 
$E, \overline{E}: \cDelta^{op} \rightarrow \sSet$ as follows:
\begin{itemize}
\item[$(1)$]
If $n \geq 0$ and $K$ is a simplicial set, then elements 
$\Hom_{\sSet}( K, E([n]) )$ can be identified with collections of maps
$\{ \sigma_{ij} \in \Hom_{K}( K \times \calM, K \times \calM) \}_{0 \leq i \leq j \leq n}$ such that each
$\sigma_{ii}$ is the identity, and $\sigma_{ij} \circ \sigma_{jk} = \sigma_{ik}$ for $0 \leq i \leq j \leq k \leq n$.

\item[$(2)$] If $n \geq 0$ and $K$ is a simplicial set, then elements 
$\Hom_{\sSet}( K, \overline{E}([n]) )$ can be identified with collections of maps
$\{ \sigma_{ij} \in \Hom_{K}( K \times \calM, K \times \calM) \}_{0 \leq i \leq j \leq n}$,
$\{ \tau_{i} \in \Hom_{K}( K, K \times \calM) \}_{0 \leq i \leq n}$ such that each
$\sigma_{ii}$ is the identity, $\sigma_{ij} \circ \sigma_{jk} = \sigma_{ik}$ for $0 \leq i \leq j \leq k \leq n$, and $\tau_i = \sigma_{ij} \circ \tau_{j}$ for $0 \leq i \leq j \leq n$. 

\item[$(3)$]  For each morphism $\alpha: [m] \rightarrow [n]$ in $\cDelta$,
the associated maps $$E([n]) \rightarrow E([m]), \quad \overline{E}([n]) \rightarrow \overline{E}([m])$$
are given by composition with $\alpha$.
\end{itemize}

We set $\End^{\otimes}(\calM) = \Nerve_{E}(\cDelta^{op})$ and
$\overline{\End}^{\otimes}(\calM) = \Nerve_{ \overline{E}}( \cDelta^{op} )$
(see Definition \ref{sulke}). We observe that there is a natural transformation of functors
$\overline{E} \rightarrow E$, given by forgetting the morphisms $\tau_{i}$; this natural transformations induces a map of simplicial sets $\overline{\End}^{\otimes}(\calM) \rightarrow \End^{\otimes}(\calM)$. 
\end{notation}\index{ZZZEndotimescalM@$\End^{\otimes}(\calM)$}\index{ZZZEndcalM@$\End(\calM)$}\index{ZZZbarEndotimescalM@$\overline{\End}^{\otimes}(\calM)$}

\begin{proposition}\label{obvus}
Let $\calM$ be an $\infty$-category. Then:
\begin{itemize}
\item[$(1)$] The map $p: \End^{\otimes}(\calM) \rightarrow \Nerve(\cDelta)^{op}$ determines a
monoidal structure on the $\infty$-category $\Fun(\calM, \calM) \simeq \End^{\otimes}_{[1]}(\calM)$. 

\item[$(2)$] The map $q: \overline{\End}^{\otimes}(\calM) \rightarrow \End^{\otimes}(\calM)$
exhibits $\calM \simeq \overline{\End}^{\otimes}_{[0]}(\calM)$ as left-tensored over
$\Fun(\calM, \calM)$. 
\end{itemize}
\end{proposition}

\begin{proof}
Proposition \ref{sulken} implies that $p$ and $p \circ q$ are coCartesian fibrations, while
Lemma \ref{sulken2} impies that $q$ is a categorical fibration. We observe
that the fiber of $p$ over a vertex $[n] \in \Nerve(\cDelta)^{op}$ is isomorphic to $\Fun(\calM, \calM)^{n}$, while the fiber of $q$ over $[n]$ is isomorphic to $\Fun(\calM, \calM)^{n} \times \calM$. It is easy to see that these identifications are compatible with the associated functors, so that $p$ and $q$ satisfies the hypotheses of Definitions \ref{mainef} and \ref{ulult}, respectively.
\end{proof}

We will refer to the monoidal structure on $\Fun(\calM, \calM)$ supplied by Proposition \ref{obvus} as the {\it composition monoidal structure}.\index{composition monoidal structure}

\begin{definition}\index{monad}
Let $\calM$ be an $\infty$-category. A {\it monad} on $\calM$ is an algebra object of $\Fun(\calM, \calM)$ (with respect to the composition monoidal structure). If $T$ is a monad on $\calM$, we let
$\Mod_{T}(\calM)$ denote the associated $\infty$-category of $T$-modules in $\calM$.
\end{definition}

\begin{remark}
More informally, a monad on an $\infty$-category $\calM$ consists of an endofunctor
$T: \calM \rightarrow \calM$ equipped with maps $1 \rightarrow T$ and $T \circ T \rightarrow T$
which satisfy the usual unit and associativity conditions, up to coherent homotopy. A {\it $T$-module} is then an object $M \in \calM$ equipped with a structure map $T(M) \rightarrow M$ which
is compatible with the algebra structure on $T$, again up to coherent homotopy.
\end{remark}

\begin{remark}
Let $\calM$ be an $\infty$-category, which we regard as an object of $\Cat_{\infty}$, and
regard $\Cat_{\infty}$ as endowed with the Cartesian monoidal structure.
According to Corollary \ref{spg}, giving a monoidal $\infty$-category $\calC$ over which
$\calM$ is left-tensored is equivalent to lifting $\calM$ to a (left) module object of
$\Cat_{\infty}$. In view of Proposition \ref{poofer}, this is equivalent to producing an algebra object
of the monoidal $\infty$-category $\Cat_{\infty}[\calM]$. In particular, $\overline{\End}^{\otimes}(\calM)$ determines an algebra object of $A \in \Alg( \Cat_{\infty}[\calM])$.

For every $\infty$-category $\calD$, the action of $\End(\calM)$ on $\calM$ determines a homotopy equivalence (even an isomorphism of simplicial sets)
$\bHom_{ \Cat_{\infty} }( \calD, \End(\calM) ) \rightarrow \bHom_{ \Cat_{\infty} }( \calD \times \calM, \calM).$
Combining Proposition \ref{ugher} with Corollary \ref{firebaugh}, we deduce that
$A$ is a final object of $\Alg( \Cat_{\infty}[\calM] )$. In other words, 
$\End^{\otimes}(\calM)$ is {\em universal} among monoidal $\infty$-categories which act on $\calM$.
\end{remark}

\subsection{Adjunction Data}\label{moadthird}

Suppose given a pair of adjoint functors
$ \Adjoint{F}{\calC}{\calD}{G}$
between $\infty$-categories $\calC$ and $\calD$. Our objective in this section
is to associate to the pair $(F,G)$ a monad $T \in \Alg( \Fun(\calC,\calC))$, given informally by the formula $T= G \circ F$. 

The basic obstacle we need to overcome is that an adjunction is {\em overdetermined} by the pair of functors $F$ and $G$. Specifying only the functor $F: \calC \rightarrow \calD$ determines a right adjoint $G$ to $F$ up to contractible ambiguity, provided that $G$ exists. If we also specify the functor $G$, we should really include additional data which identifies $G$ with an adjoint to $F$.
Such data is provided by {\em either} a unit
$ u: \id_{\calC} \rightarrow G \circ F$
or a counit
$ v: F \circ  G \rightarrow \id_{\calD}$
for the adjunction. However, if we specify {\em both} a unit and a counit, then the adjunction is again overdetermined. In classical category theory, one imposes a condition on a unit and counit: they are said to be {\em compatible} if the composite transformations
$$ \alpha:F \stackrel{u}{\rightarrow} F \circ (G \circ F) = (F \circ G) \circ F \stackrel{v}{\rightarrow} F$$
$$ \beta: G \stackrel{u}{\rightarrow} (G \circ F) \circ G = G \circ (F \circ G) \stackrel{v}{\rightarrow} G$$
coincide with the identity. In the higher categorical setting, we should instead require the existence of homotopies $h: \alpha \simeq \id_{F}$, $h': \beta \simeq \id_{G}$. Once again, we obtain the correct theory if we specify {\em either} $h$ or $h'$, but specifying {\em both} will overdetermine the adjunction. This leads us to formulate further compatibilities between $h$ and $h'$, and so forth.
There are two different strategies for dealing with the situation:

\begin{itemize}
\item[$(1)$] Specify a minimal amount of data: for example, the single functor $F$.
\item[$(2)$] Specify {\em all} of the relevant data: the functors $F$ and $G$, the unit $u$ and counit $v$, the homotopies $h$ and $h'$, and all of their higher-dimensional relatives.
\end{itemize}

For most applications, approach $(1)$ is entirely sufficient. However, if we wish to extract
a monad from the adjoint pair $(F,G)$, then we are forced to adopt approach $(2)$.
Our primary goal in this section is to make precise sense of $(2)$, by introducing the notion of
an {\it adjunction datum} between the $\infty$-categories $\calC$ and $\calD$. Our main result Theorem \ref{partA}, which asserts that $(1)$ and $(2)$ are actually equivalent to one another: in other words, if $F: \calC \rightarrow \calD$ is a functor which admits a right adjoint, then $F$ can be promoted (in an essentially unique way) to an adjunction datum. The proof of this result is very technical (probably the most difficult result in this paper) and will be given in \S \ref{digtwo}.

We begin with some generalities. Recall that the classical theory of monoidal categories can be regarded as a special case of the theory of (weak) {\it $2$-categories}. More precisely, giving a monoidal category $(\calC, \otimes)$ is essentially the same thing as giving a $2$-category
$\calD$ having only a single, fixed object $X \in \calD$. The correspondence assigns to
$X \in \calD$ the category $\Hom_{\calD}( X,X)$, with monoidal structure given by composition in
$\calD$. We can obtain a generalization of the theory of monoidal categories by allowing
$\calD$ to have many objects. We now describe the analogous generalization of the theory of monoidal $\infty$-categories.

\begin{definition}\index{ZZZcDeltaS@$\cDelta_{S}$}
Let $S$ be a set of symbols. We define a category $\cDelta_{S}$ as follows.
An object of $\cDelta_{S}$ consists of a pair $([n], c)$, where $[n] \in \cDelta$ and
$c: [n] \rightarrow S$ is an arbitrary map. Given a pair of objects $( [m], c), ([n], c') \in \cDelta_{S}$, we set
$$ \Hom_{ \cDelta_{S} } ( ( [m], c) , ([n], c') ) = \{ \alpha \in \Hom_{\cDelta}([m], [n])
: ( \forall 0 \leq i \leq m) [ c(i) = (c' \circ \alpha)(i) ] \}.$$
\end{definition}

We will think of an object of $\cDelta_{S}$ as a (nonempty) finite sequence
$[ c(0) , c(1), \ldots, c(n) ]$ of elements of $S$; in other words, as a nonempty {\it word}
in the alphabet $S$.

One can define an {\it $S$-colored monoidal $\infty$-category} to be a coCartesian
fibration $\calC^{\otimes} \rightarrow \Nerve( \cDelta_{S})^{op}$, satisfying an appropriate analogue of the condition $(\ast)$ of Definition \ref{mainef}. We can think of an $S$-colored monoidal $\infty$-category as encoding the structure of an $(\infty,2)$-category whose objects are in bijection with the set $S$. We will refrain from going into the details, since we are primarily interested in only a single example.

\begin{notation}\label{defEE}
Let $\calC$ and $\calD$ be simplicial sets, and let
$S = \{ \calC, \calD \}$. 
We define a functor
$E: \cDelta_{S}^{op} \rightarrow \sSet$ as follows:
\begin{itemize}
\item[$(1)$] Let $n \geq 0$ be a finite nonempty linearly ordered set, and let $c: [n] \rightarrow \{ \calC, \calD \}$ be a map, and let $K \in \sSet$ be arbitrary. Then elements of 
$\Hom_{\sSet}(K, E([n],c) )$ can be identified with collections of maps
$ \{ \sigma_{ij} \in \Hom_{K}( K \times c(j), K \times c(i) ) \}_{0 \leq i \leq j \leq n}$ such that
$\sigma_{ii} = \id$ for $i \in J$, and $\sigma_{ij} \circ \sigma_{jk} = \sigma_{ik}$ for $i \leq j \leq k$.
\item[$(2)$] Given a map $f: ([m], c) \rightarrow ([n], c')$ in $\cDelta_{S}$, the associated map $E([n], c') \rightarrow E([m], c)$ is given by pullback along $f$.
\end{itemize}
\end{notation}

\begin{remark}
In Notation \ref{defEE}, and all that follows, we will implicitly assume that
$\calC \neq \calD$, so that the set $S = \{ \calC, \calD \}$ contains exactly two symbols. (This is no loss of generality, since we can always achieve this state of affairs by replacing $\calC$ or $\calD$ by an isomorphic simplicial set.) This is purely a matter of notational convenience.
\end{remark}

\begin{definition}\label{strad0}\index{ZZZEndotimesCD@$\End^{\otimes}(\calC, \calD)$}
Let $\calC$ and $\calD$ be $\infty$-categories, and let
$S = \{ \calC, \calD\}$. We let $\End^{\otimes}(\calC,\calD)$ denote the relative nerve $\Nerve_{E}( \cDelta_{S}^{op} )$.
\end{definition}

\begin{remark}
The projection map $q: \End^{\otimes}(\calC, \calD) \rightarrow \Nerve(\cDelta_{S})^{op}$ is an example of an $S$-colored monoidal $\infty$-category. One can think of this monoidal $\infty$-category as follows: the collection of all $\infty$-categories really constitutes an $(\infty,2)$-category, since we can associate to every pair
of $\infty$-categories $\calE$ and $\calE'$ an $\infty$-category $\Fun(\calE, \calE')$. Then $q$ encodes the full subcategory spanned by the pair of $\infty$-categories $\calC$ and $\calD$.
\end{remark}

\begin{definition}\label{strad2}
Let $\calC$ and $\calD$ be $\infty$-categories, and let $S = \{ \calC, \calD\}$. 
We will say that a morphism $\alpha: ([m],c) \rightarrow ( [n], c')$ in $\cDelta_{S}$ is {\it $\calC$-convex} if, whenever $j \in [n]$ satisfies $c'(j) = \calC$ and
$\alpha(i) \leq j \leq \alpha(i')$ for suitably chosen $i,i' \in [m]$, there exists a unique
$j_0 \in [m]$ such that $\alpha(j_0) = j$.
We will say that a section $U$ of the projection $q: \End^{\otimes}(\calC, \calD) \rightarrow
\Nerve(\cDelta_{S})^{op}$ is an {\it adjunction datum} if it carries
every $\calC$-convex morphism in $\cDelta_{S}$ to a $q$-coCartesian morphism of
$\End^{\otimes}(\calC, \calD)$. 
We let $\AdjDiag(\calC, \calD)$ denote the full subcategory of $\bHom_{ \Nerve(\cDelta_S)^{op}}( \Nerve(\cDelta_{S})^{op}, \End^{\otimes}(\calC,\calD))$ spanned by the adjunction data.\index{adjunction datum}\index{ZZZAdjDiag@$\AdjDiag(\calC,\calD)$}\index{convex!$\calC$}\index{$\calC$-convex}
\end{definition}

\begin{remark}\label{mondet}
We can identify $\cDelta \simeq \cDelta_{ \{\calC\} }$ with the full subcategory of $\cDelta_{S}$ spanned by the objects $([n], c)$, where $c: [n] \rightarrow \{ \calC, \calD \}$ takes the constant value $\calC$. Under this identification, a morphism in $\cDelta$ is $\calC$-convex if and only if it is convex in the sense of Definition \ref{defcondef}.

The fiber product 
$\End^{\otimes}(\calC, \calD) \times_{\Nerve( \cDelta_{S} )^{op}} \Nerve(\cDelta)^{op}$ is canonically isomorphic to the monoidal $\infty$-category $\End^{\otimes}(\calC)$ defined in \S \ref{moadsec}. Restriction determines a functor $\AdjDiag(\calC, \calD) \rightarrow \Alg( \End(\calC) )$, so that every adjunction datum determines a monad on $\calC$.

We can use the same argument to produce a restriction map
$\AdjDiag(\calC, \calD) \rightarrow \Alg( \End(\calD) ).$ However, this map is less interesting: the monad on $\calD$ induced by an adjunction datum is always equivalent
to the initial object of $\Alg( \End(\calD) )$ (since every morphism of
$\cDelta_{ \{\calD \} }$ is $\calC$-convex in $\cDelta_{S}$.
\end{remark}

We now proceed to unwind the details of Definition \ref{strad2}.

\begin{remark}\label{humberton}
Let $\calC$ and $\calD$ be $\infty$-categories, let $S= \{\calC, \calD\}$, and let $U \in \AdjDiag(\calC, \calD)$. Let us agree to denote an object of $\cDelta_{S}$ by a finite (nonempty) string of elements of $S$. Then:
\begin{itemize}

\item[$(i)$] Evaluation of $U$ at the object $[\calD,\calC] \in \cDelta_{S}$ determines a functor
$F: \calC \rightarrow \calD$.

\item[$(ii)$] Evaluation of $U$ at the object $[\calC,\calD] \in \cDelta_{S}$ determines another functor
$G: \calD \rightarrow \calC$.

\item[$(iii)$] Evaluation of $U$ at the object $[\calD,\calD] \in \cDelta_{S}$ determines a functor 
$i: \calD \rightarrow \calD$. The morphism $[\calD,\calD] \rightarrow [\calD]$ determines
a natural transformation $\id_{\calD} \rightarrow i$ in $\Fun(\calD, \calD)$. Since the morphism
$[ \calD, \calD ] \rightarrow \calD$ is $\calC$-convex, this natural transformation is an equivalence, so that $i$ is (canonically) equivalent to the identity functor on $\calD$.

\item[$(iv)$] Evaluation of $U$ at the object $[\calC,\calC] \in \cDelta_{S}$ determines a functor
$T: \calC \rightarrow \calC$. The morphism $[\calC,\calC] \rightarrow [\calC]$ determines
a natural transformation $\id_{\calC} \rightarrow T$. This transformation is generally not an equivalence. Rather, it is the unit map of a monad structure on $T$ (see Remark \ref{mondet}).

\item[$(v)$] Evaluation of $U$ at the object $[\calD,\calC,\calD] \in \cDelta_{S}$ determines a diagram
$$ \xymatrix{ \calD \ar[rr]^{F' \circ G'} \ar[dr]^{G'} & & \calD \\
& \calC \ar[ur]^{F'}. & }$$
Moreover, applying $s$ to the inclusions of $[\calD, \calC]$, $[\calD, \calC]$, and
$[\calD, \calD]$ into $[\calD, \calC, \calD]$ determines natural transformations
$\alpha: F' \rightarrow F$, $\beta: G' \rightarrow G$, and $F' \circ G' \rightarrow i$. Since the first two of these inclusions are $\calC$-convex, we deduce that $\alpha$ and $\beta$ are equivalences.
Combining these three transformations with $(iii)$, we obtain a natural transformation $v: F \circ G \rightarrow \id_{\calD}$, well-defined up to homotopy.

\item[$(vi)$] Evaluation of $U$ at the object $[\calC,\calD,\calC] \in \cDelta_{S}$ determines a diagram
$$ \xymatrix{ \calC \ar[rr]^{G'' \circ F''} \ar[dr]^{F''} & & \calC \\
& \calD \ar[ur]^{G''}. & }$$
Moreover, applying $s$ to the inclusions of $[\calC,\calD]$, $[\calC,\calC]$ and $[\calD,\calC]$
into $[\calC,\calD,\calC]$ yields natural transformations $F'' \rightarrow F$, 
$G'' \rightarrow G$, and $G'' \circ F'' \rightarrow T$. All three of these inclusions are $\calC$-convex, so the induced natural transformations are all equivalences.
In particular, the transformation $\id_{\calC} \rightarrow T$ of $(iv)$ determines a natural transformation
$u: \id_{\calC} \rightarrow G \circ F$, which is well-defined up to homotopy.
\end{itemize}
\end{remark}

\begin{lemma}\label{hungertown}
Let $\calC$ and $\calD$ be $\infty$-categories, and suppose given an adjunction datum $U \in \AdjDiag(\calC, \calD)$. Let $F: \calC \rightarrow \calD$, $G: \calD \rightarrow \calC$,
$u: \id_{\calC} \rightarrow G \circ F$, and $v: F \circ G \rightarrow \id_{\calD}$ be defined as in 
Remark \ref{humberton}. Then $u$ is the unit of an adjunction between $F$ and $G$, and $v$ is the counit of an adjunction between $F$ and $G$. In particular, $F$ and $G$ are adjoint to one another.
\end{lemma}

\begin{proof}
Let $S= \{\calC, \calD \}$. We will show that $u$ and $v$ determine an adjunction between the underlying functors between the homotopy categories $\h{\calC}$ and $\h{\calD}$. For this, it suffices to show that the
compositions
$$ F \stackrel{u}{\rightarrow} F \circ (G \circ F) = (F \circ G) \circ F \stackrel{v}{\rightarrow} F$$
$$ G \stackrel{u}{\rightarrow} (G \circ F) \circ G = G \circ (F \circ G) \stackrel{v}{\rightarrow} G$$
are homotopic to the identity.

To prove the first claim, we consider the commutative diagram
$$ 
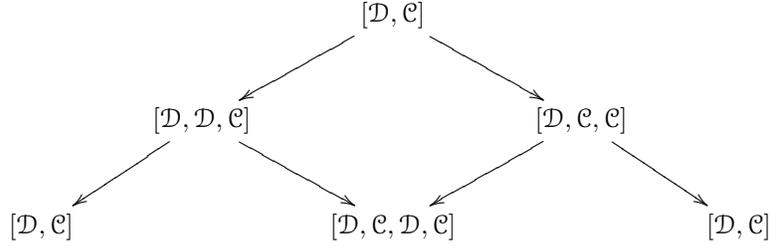
\xymatrix{ & & [ \calD, \calC] \ar[dr] \ar[dl] & & \\
& [\calD, \calD, \calC] \ar[dl] \ar[dr] & & [\calD, \calC, \calC] \ar[dl] \ar[dr] & \\
[\calD, \calC] & & [\calD, \calC, \calD,\calC] & & [\calD, \calC] }$$
in $\cDelta_{S}$, where each of the corresponding maps of linearly ordered sets preserves both the initial and final objects. Applying the functor $U$ and evaluating at
those initial and final objects, we obtain a diagram in the homotopy category $\h{\Fun(\calC, \calD)}$, which is equivalent to
$$ 
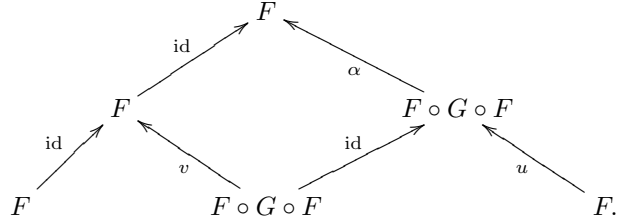
\xymatrix{ & & F & & \\
& F \ar[ur]^{\id} & & F \circ G \circ F \ar[ul]^{\alpha} & \\
F \ar[ur]^{\id} & & F \circ G \circ F \ar[ul]^{v} \ar[ur]^{\id} & & F \ar[ul]^{u}. }$$
The desired result now follows from the observation that the composition of $\alpha$
with $u$ is (homotopic to) the identity.

The proof of the second claim is similar, but we consider instead the diagram
$$ 
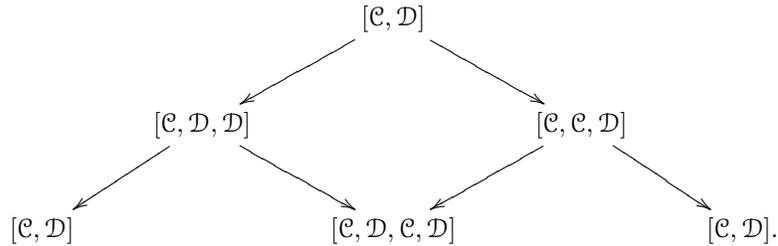
\xymatrix{ & & [\calC,\calD] \ar[dr] \ar[dl] & & \\
& [\calC, \calD, \calD] \ar[dl] \ar[dr] & & [\calC, \calC, \calD] \ar[dl] \ar[dr] & \\
[\calC, \calD] & & [\calC, \calD, \calC,\calD] & & [\calC, \calD]. }$$
Applying $s$ and evaluating at the initial and final objects, we obtain a diagram in
$\h{ \Fun(\calD, \calC)}$ which is equivalent to
$$ 
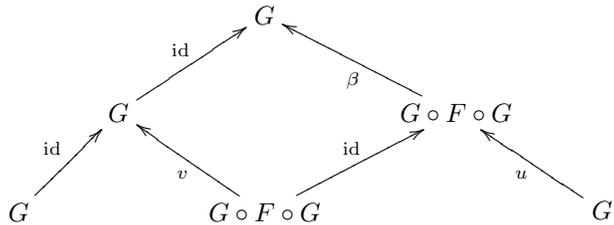
\xymatrix{ & & G & & \\
& G \ar[ur]^{\id} & & G \circ F \circ G \ar[ul]^{\beta} & \\
G \ar[ur]^{\id} & & G \circ F \circ G \ar[ul]^{v} \ar[ur]^{\id} & & G \ar[ul]^{u}, }$$
and the desired result follows from the observation that the composition of $u$ with $\beta$
is the identity on $G$.
\end{proof}

In the situation of Lemma \ref{hungertown}, either the unit $u$ or the counit $v$ are sufficient to identify $F$ with the left adjoint of $G$. An adjunction datum $U$ therefore supplies {\em two} such identifications. However, these two identifications are compatible with one another, up to coherent homotopy: this compatibility is encoded in the values $U$ on more complicated
objects of $\cDelta_{S}$. The result is that, up to a contractible space of choices, the entire diagram $U$ is determined by the functor $F$:

\begin{theorem}\label{partA}
Let $\calC$ and $\calD$ be $\infty$-categories, let $S = \{\calC, \calD\}$, and let
$\theta: \AdjDiag(\calC, \calD) \rightarrow \Fun(\calC, \calD)$
be given by evaluation at $[\calD, \calC] \in \cDelta_{S}$. Then $\theta$ induces
a trivial Kan fibration $\AdjDiag(\calC, \calD) \rightarrow \Fun'(\calC, \calD)$, where
$\Fun'(\calC, \calD)$ denotes the subcategory of $\Fun(\calC, \calD)$ whose objects are functors
from $\calC$ to $\calD$ which admit right adjoints, and whose morphisms are equivalences of functors. In particular, $\AdjDiag(\calC, \calD)$ is a Kan complex.
\end{theorem}

In other words, if a functor $F: \calC \rightarrow \calD$ admits a right adjoint $G$, then $F$ can be extended (in an essentially unique way) to an adjunction datum $U \in \AdjDiag(\calC, \calD)$. Invoking Remark \ref{mondet}, we see that $F$ determines a monad on $\calC$, whose underlying functor from $\calC$ to $\calC$ can be identified with $G \circ F$ (Remark \ref{humberton}). We may therefore regard Theorem \ref{partA} as an $\infty$-categorical analogue of assertion $(A)$ (see the introduction to \S \ref{barri}). The proof will be given in \S \ref{digtwo}.

\subsection{Pointed Adjunction Data}\label{datpoint}

Our goal in this section is to establish an $\infty$-categorical analogue of assertion
$(B)$ of \S \ref{barri}. Let $\calC$ and $\calD$ be a pair of $\infty$-categories, fixed throughout this section, and let $S = \{ \calC, \calD \}$. For each adjunction datum $U \in \AdjDiag(\calC, \calD)$, we will construct a diagram
$ \calD \leftarrow \overline{\AdjDiag}_{U}(\calC, \calD) \rightarrow \Mod_{T}(\calC),$
where $T$ is the monad on $\calC$ determined by $U$. Moreover, we will show that
the left map is a trivial Kan fibration (Corollary \ref{hooku}), so that $U$ determines a map
$\calD \rightarrow \Mod_{T}(\calC)$, which is well-defined up to a contractible space of choices.
 
Our first step is to introduce an elaboration of
the simplicial set $\End^{\otimes}(\calC, \calD)$.

\begin{notation}
We define a functor $\overline{E}: \cDelta_{S}^{op} \rightarrow \sSet$ as follows:
\begin{itemize}
\item[$(1)$] Let $n \geq 0$, let $c: [n] \rightarrow \{ \calC, \calD \}$ be a map, and let $K$ be a simplicial set. Then elements of $\Hom_{\sSet}(K, \overline{E}([n],c) )$ can be identified with collections of maps
$$ \{ \sigma_{ij} \in \Hom_{K}( K \times c(j), K \times c(i) ) \}_{0 \leq i \leq j \leq n}, \{ \tau_{i} \in \Hom_{K}( K, K \times c(i) \}_{0 \leq i \leq n}$$ such that $\sigma_{ii} = \id$ for $0 \leq i \leq n$, $\sigma_{ij} \circ \sigma{jk} = \sigma_{ik}$ for $0 \leq i \leq j \leq k \leq n$, and $\sigma_{ij} \circ \tau_{j} = \tau_{i}$ for $0 \leq i \leq j \leq n$.

\item[$(2)$] Given a map $f: ([m], c) \rightarrow ([n], c')$ in $\cDelta_{S}$, the associated map $\overline{E}([n],c') \rightarrow \overline{E}([m],c)$ is given by pullback along the induced map $[m] \rightarrow [n]$.
\end{itemize}

We let $\overline{\End}^{\otimes}(\calC, \calD) = \Nerve_{\overline{E}}(\cDelta_{S}^{op})$ (see Definition \ref{sulke}).\index{ZZZbarEndCD@$\overline{\End}^{\otimes}(\calC,\calD)$}
\end{notation}

\begin{remark}
Let $\calC$ and $\calD$ are $\infty$-categories. Proposition \ref{sulken} implies that the projection $$\overline{\End}^{\otimes}(\calC, \calD) \rightarrow \Nerve( \cDelta_{S})^{op}$$ is a coCartesian fibration of simplicial sets. 
We observe that there is a canonical projection $\overline{\End}^{\otimes}(\calC, \calD) \rightarrow
\End^{\otimes}(\calC, \calD)$, given by forgetting the maps $\tau_{i}$.
Moreover, if we identify $\cDelta$ with the full subcategory of $\cDelta_{S}$ spanned by those functions which take the constant value $\calC$, then the fiber product
$$\overline{\End}^{\otimes}(\calC, \calD) \times_{\Nerve( \cDelta_{S} )^{op} }
\Nerve(\cDelta)^{op}$$
is canonically isomorphic to the simplicial set $\overline{\End}^{\otimes}(\calC)$ defined in \S \ref{moadsec}.
\end{remark}

\begin{definition}\label{termus}
Consider the maps
$ \overline{\End}^{\otimes}(\calC, \calD) \stackrel{p}{\rightarrow}
\End^{\otimes}(\calC, \calD) \stackrel{q}{\rightarrow} \Nerve( \cDelta_{S})^{op}.$\index{adjunction datum!pointed}\index{pointed adjunction datum}
\index{ZZZbarAdjDiag@$\overline{\AdjDiag}(\calC, \calD)$}
A {\it pointed adjunction datum} is a map
$\overline{U}: \Nerve( \cDelta_{S} )^{op} \rightarrow  \overline{\End}^{\otimes}(\calC, \calD)$ with the following properties:
\begin{itemize}
\item[$(1)$] The map $\overline{U}$ is a section of $q \circ p$; that is, $q \circ p \circ \overline{U} = \id$.
\item[$(2)$] The map $p \circ \overline{U}: \Nerve( \cDelta_{S})^{op} \rightarrow \End^{\otimes}(\calC, \calD)$
is an adjunction datum.
\item[$(3)$] Suppose that $\alpha: ([m],c) \rightarrow ([n],c')$ is a morphism in $\cDelta_{S}$ such that $\alpha(m) = n$. Then $\overline{U}(\alpha)$ is a locally $p$-coCartesian morphism in $\overline{\End}^{\otimes}(\calC, \calD)$.
\item[$(4)$] Let $\alpha$ be the unique morphism from
$[\calC]$ to $[\calC, \calD]$ in $\cDelta_{S}$. Then $\overline{U}(\alpha)$ is a locally $p$-coCartesian morphism in $\overline{\End}^{\otimes}(\calC, \calD)$.
\end{itemize}
We let $\overline{\AdjDiag}(\calC, \calD)$ denote the full subcategory of 
$\bHom_{\Nerve(\cDelta_{S})^{op}}(\Nerve(\cDelta_{S})^{op}, 
\overline{\End}^{\otimes}(\calC, \calD))$ spanned by the pointed adjunction data.
\end{definition}

\begin{remark}\label{supinor}
Restriction to the subcategory $\Nerve( \cDelta_{ \{\calC\} })^{op} \subseteq \Nerve(\cDelta_{S})^{op}$ determines a functor $\overline{\AdjDiag}(\calC,\calD) \rightarrow \Mod(\calC),$
where we regard $\calC$ as left-tensored over $\End^{\otimes}(\calC)$.
Composition with the projection $\overline{\End}^{\otimes}(\calC, \calD) \rightarrow
\End^{\otimes}(\calC, \calD)$ determines another functor
$\overline{\AdjDiag}(\calC, \calD) \rightarrow \AdjDiag(\calC, \calD).$
Finally, evaluation on the object $[\calD] \in \cDelta_{S}$ determines a functor $\overline{\AdjDiag}(\calC, \calD) \rightarrow \calD$.
\end{remark}

Our main result can be stated as follows:

\begin{proposition}\label{stumpus}
Let $\calC$ and $\calD$ be $\infty$-categories, and let
$\theta: \overline{\AdjDiag}(\calC, \calD) \rightarrow \AdjDiag(\calC, \calD) \times \calD$
be the product of the maps described in Remark \ref{supinor}. Then $\theta$ is a trivial fibration of simplicial sets.
\end{proposition}

\begin{corollary}\label{hooku}
Let $\calC$ and $\calD$ be $\infty$-categories, let $s \in \AdjDiag(\calC, \calD)$, and let
$\overline{\AdjDiag}_s(\calC,\calD)$ denote the fiber product
$\overline{\AdjDiag}(\calC,\calD) \times_{ \AdjDiag(\calC, \calD) } \{s\}$. Then the evaluation map
$\overline{\AdjDiag}_{s}(\calC, \calD) \rightarrow \calD$ is a trivial fibration of simplicial sets.
\end{corollary}

The proof will require the following somewhat technical lemma:

\begin{lemma}\label{mikest}
Consider the projection map
$p: \overline{\End}^{\otimes}(\calC, \calD) \rightarrow \End^{\otimes}(\calC, \calD)$. Then:
\begin{itemize}
\item[$(1)$] The map $p$ is a locally coCartesian fibration.
\item[$(2)$] Suppose that $\alpha: ([m],c) \rightarrow ([n],c')$ is a morphism in
$\cDelta_{S}$ which induces a {\em convex} morphism $[m] \rightarrow [n]$, and let
$f$ be a morphism in $\overline{\End}^{\otimes}(\calC, \calD)$ which projects to $\alpha$. Then $f$ is $p$-coCartesian if and only if $f$ is locally $p$-coCartesian.
\item[$(3)$] Suppose that $\alpha: ([m],c) \rightarrow ([n],c')$ is a morphism in
$\cDelta_{S}$ such that $\alpha(m) = n$, and let
$f$ be a morphism in $\overline{\End}^{\otimes}(\calC, \calD)$ which projects to $\alpha$. The following conditions are equivalent:
\begin{itemize}
\item[$(i)$] The morphism $f$ is locally $p$-coCartesian.
\item[$(ii)$] The morphism $f$ is locally $p$-Cartesian.
\item[$(iii)$] The morphism $f$ is $p$-Cartesian.
\end{itemize}
\end{itemize}
\end{lemma}

\begin{proof}
Assertion $(1)$ follows from Proposition \toposref{fibertest}. The map $p$ is generally not a coCartesian fibration, because the class of locally $p$-coCartesian morphisms is not stable under composition. Nevertheless, suppose given a $2$-simplex
$$ \xymatrix{ & Y \ar[dr]^{e''} & \\
X \ar[ur]^{e'} \ar[rr]^{e} & & Z }$$
in $\overline{\End}^{\otimes}(\calC, \calD)$ where the morphisms $e'$ and $e''$ are locally $p$-coCartesian, covering a diagram
$$ \xymatrix{ & ([m], c') \ar[dl]^{\alpha'} & \\
([n],c'') & & ([k], c) \ar[ul]^{\alpha''} \ar[ll]^{\alpha} }$$
in the category $\cDelta_{S}$. Using the description of the class of locally $p$-coCartesian morphisms provided by Proposition \toposref{fibertest}, we conclude that $e$ is
locally $p$-coCartesian provided that the following condition is satisfied:
\begin{itemize}
\item[$(\ast)$] The map $\alpha$ determines a bijection
$$ \{ i \in [m]: i > \alpha''(k)  \} \rightarrow \{ j \in [n]: \alpha(k) < j \leq \alpha'(m)  \}.$$
(This implies, in particular, that if $\alpha''(k) < i \leq m$, then $\alpha(k) < \alpha(i)$.)
\end{itemize}

We note that condition $(\ast)$ is always satisfied if the map $\alpha'$ induces a {\em convex} map
$[m] \hookrightarrow [n]$. Assertion $(2)$ now follows from Lemma \toposref{charloccart}. Similarly, condition $(\ast)$ is automatically satisfied if $\alpha''(k) = m$. Using Lemma \toposref{gruft}, we deduce the equivalence $(ii) \Leftrightarrow (iii)$ of $(3)$. The equivalence $(i) \Leftrightarrow (ii)$ follows from the observation that if $f: X \rightarrow Y$ is as in $(3)$, then the associated functor
$\overline{\End}^{\otimes}(\calC, \calD)_{X} \rightarrow \overline{\End}^{\otimes}(\calC, \calD)_{Y}$ is an equivalence of $\infty$-categories.
\end{proof}

\begin{remark}\label{recaster}
In the situation of Definition \ref{termus}, condition $(4)$ is can be replaced by the following apparently stronger condition:
\begin{itemize}
\item[$(4')$] Let $n > 0$, let $c: [n] \rightarrow \{ \calC, \calD \}$ satisfy 
$c(n) = \calD$, let $c_0 = c | \{0, \ldots, n-1\}$, and let $\alpha: ( [n-1], c_0) \rightarrow
( [n],c)$ be the inclusion. Then $s(\alpha)$ is a $p$-coCartesian morphism in
$\overline{\End}^{\otimes}(\calC, \calD)$.  
\end{itemize}
It is clear that $(4') \Rightarrow (4)$. Conversely, suppose that $\overline{U}$ satisfies the conditions of Definition \ref{termus}, and let $\alpha$ be as in the statement of $(4')$. We wish to prove that
$\overline{U}(\alpha)$ is $p$-coCartesian. Lemma \ref{mikest} shows that
$\overline{U}(\alpha)$ is $p$-coCartesian if and only if $\overline{U}(\alpha)$ is locally $p$-coCartesian. Using
assumption $(2)$, we deduce that $\overline{U}(\alpha)$ is locally $p$-coCartesian if and only if
$s(\alpha)$ is locally $p$-Cartesian.

If $c(n-1) = \calD$, then we have a commutative diagram
$$ \xymatrix{ & ([n-1], c_0) \ar[dr]^{\alpha} & \\
( [n], c) \ar[ur]^{\beta} \ar[rr]^{\id} & & ([n], c) }$$
in $\cDelta_{S}$.  In view of Proposition \toposref{protohermes}, it will suffice to show that $\overline{U}(\id)$ and $\overline{U}(\beta)$ are $p$-Cartesian. For $\overline{U}(\id)$ this is obvious; for
$\overline{U}(\beta)$ we combine Lemma \ref{mikest} with assumption $(3)$.

Suppose instead that $c(n-1) = \calC$. Let $c' = c | \{n-1, n\}$ and let $c'_0 = c | \{n-1\}$. We have a commutative diagram
$$ \xymatrix{ ( \{n-1\}, c'_0) \ar[r]^{\alpha'} \ar[d]^{\beta} & ( \{n-1, n\}, c' ) \ar[d]^{\beta'} \\
( [n-1], c_0) \ar[r]^{\alpha} & ( [n], c) }$$
in the category $\cDelta_{S}$. We wish to prove that $\overline{U}(\alpha)$ is locally
$p$-Cartesian. Using $(2)$ and Lemma \ref{mikest}, we deduce that
$\overline{U}(\beta)$ is $p$-Cartesian. According to Lemma \toposref{charloccart}, it will suffice to show
that $\overline{U}( \alpha \circ \beta)$ is locally $p$-Cartesian. Using $(2)$, we reduce to showing that $\overline{U}(\alpha \circ \beta)$ is locally $p$-coCartesian. Conditions $(3)$ and $(4)$ imply that
$\overline{U}(\alpha')$ and $\overline{U}(\beta')$ are locally $p$-coCartesian. Invoking Lemma \ref{mikest}, we deduce that $\overline{U}(\beta')$ is $p$-coCartesian. Lemma \toposref{charloccart} now implies that
$\overline{U}( \beta' \circ \alpha') \simeq \overline{U}( \alpha \circ \beta)$ is locally $p$-coCartesian, as desired.
\end{remark}

\begin{remark}\label{recaster1}
In the situation of Definition \ref{termus}, condition $(3)$ is equivalent to the following
apparently weaker condition: 
\begin{itemize}
\item[$(3')$] For every object $([n],c) \in \cDelta_{S}$, let $c_0 = c | \{n\}$ and let
$\alpha: ( \{n\}, c_0) \rightarrow ( [n], c)$ denote the inclusion. Then $\overline{U}(\alpha)$ is locally $p$-coCartesian.
\end{itemize}
It is clear that $(3) \Rightarrow (3')$. Conversely, suppose that $(3')$ is satisfied, and consider an arbitrary morphism $\beta: ( [m], c) \rightarrow ([n], c')$ such that $\beta(m) = n$. Let 
$c_0 = c | \{m\}$, and form a commutative diagram
$$ \xymatrix{ & ( [m], c) \ar[dr]^{\beta} & \\
( \{m\}, c_0 ) \ar[ur]^{\alpha} \ar[rr] & & ( [n], c'). }$$
We wish to show that $\overline{U}(\beta)$ is locally $p$-coCartesian. In view of Lemma \ref{mikest}, it will suffice to show $\overline{U}(\beta)$ is $p$-Cartesian. Using Lemma \toposref{charloccart}, it will suffice to show that $\overline{U}(\alpha)$ and $\overline{U}( \beta \circ \alpha)$ are $p$-Cartesian. This follows immediately from $(3')$ and Lemma \ref{mikest}.
\end{remark}

\begin{remark}\label{recaster2}
Suppose we are in the situation of Definition \ref{termus}, and that
$\overline{U}: \Nerve( \cDelta_{S} )^{op} \rightarrow  \overline{\End}^{\otimes}(\calC, \calD)$ satisfies
conditions $(1)$, $(2)$, and $(4')$ (see Remark \ref{recaster}). Then $(3)$ is equivalent to the following apparently weaker condition:
\begin{itemize}
\item[$(3'')$] Let $([n],c)$ be an object of $\cDelta_{S}$ such that $c(n) = \calD$, let $c_0 = c | \{n\}$ and let $\alpha: ( \{n\}, c_0) \rightarrow ( [n], c)$ be the inclusion. Then $\overline{U}(\alpha)$ is locally $p$-coCartesian. 
\end{itemize}
It is clear that $(3) \Rightarrow (3'')$. To prove the converse, it will suffice to show that
$(3'')$ implies condition $(3')$ of Remark \ref{recaster1}. Suppose given 
an object $([n], c) \in \cDelta_{S}$, let $c_0 = c | \{n\}$ and let
$\alpha: ( \{n\}, c_0) \rightarrow ([n],c)$ be the inclusion. We wish to prove that
$\overline{U}(\alpha)$ is locally $p$-coCartesian. If $c(n) = \calC$ this follows immediately from $(3'')$.
Otherwise, define $c': [n+1] \rightarrow \{ \calC, \calD\}$ by the formula
$$ c'(i) = \begin{cases} c(i) & \text{if } 0 \leq i \leq n \\
\calD & \text{if } i = n+1, \end{cases}$$
and let $c'_0 = c' | \{n,n+1\}$. 
We have a commutative diagram
$$ \xymatrix{ ( \{n\}, c_0) \ar[r]^{\alpha} \ar[d]^{\beta'} & ([n],c) \ar[d]^{\beta} \\
( \{n, n+1\}, c'_0) \ar[r]^{\alpha'} & ([n+1], c') }$$
in the category $\cDelta_{S}$. Condition $(4')$ implies that $\overline{U}(\beta)$ is $p$-coCartesian. Lemma \toposref{charloccart} now asserts that $\overline{U}(\alpha)$ is locally $p$-coCartesian if and only if
$\overline{U}( \beta \circ \alpha)$ is locally $p$-coCartesian. For this, it will suffice to show that
$\overline{U}(\alpha')$ and $\overline{U}(\beta')$ are $p$-coCartesian. For $\overline{U}(\beta')$, this follows from
$(4')$. To show that $\overline{U}(\alpha')$ is $p$-coCartesian, we first invoke Lemma \ref{mikest} to reduce
to the problem of showing that $\overline{U}(\alpha')$ is locally $p$-coCartesian, then $(2)$ to reduce to the problem of showing that $\overline{U}(\alpha)$ is locally $p$-Cartesian. We now set $c'_1 = c' | \{n+1\}$ and consider the commutative diagram
$$ \xymatrix{ & ( \{n, n+1\}, c'_0) \ar[dr]^{\alpha'} & \\
( \{n+1\}, c'_1 ) \ar[ur]^{\gamma} ) \ar[rr] & & ( [n+1], c'). }$$
In view of Lemma \toposref{charloccart}, it will suffice to show that
$\overline{U}(\gamma)$ and $\overline{U}(\alpha' \circ \gamma)$ are $p$-Cartesian. This follows immediately from $(3'')$ and Lemma \ref{mikest}. 
\end{remark}

\begin{proof}[Proof of Proposition \ref{stumpus}]
Lemma \ref{sulken2} implies that the projection $\overline{\End}^{\otimes}(\calC, \calD) \rightarrow \End^{\otimes}(\calC,\calD)$ is a categorical fibration. It follows easily that $\theta$ is a categorical fibration. It will therefore suffice to show that $\theta$ is an equivalence of $\infty$-categories.

Theorem \ref{partA} implies that $\AdjDiag(\calC, \calD)$ is a Kan complex. Since
the projection $\theta_0: \overline{ \AdjDiag}(\calC, \calD) \rightarrow \AdjDiag(\calC, \calD)$ is a categorical fibration, Proposition \toposref{groob} implies that $\theta_0$ is a coCartesian fibration. 
Invoking Corollary \toposref{usefir}, we deduce that $\theta$ is a categorical equivalence if and only if, for every adjunction datum $U \in \AdjDiag(\calC, \calD)$, the induced map
$\overline{\AdjDiag}_{U}(\calC, \calD) \rightarrow \calD$ is a categorical equivalence. In other words, we are reduced to proving a slightly weaker version of Corollary \ref{hooku}.

Form a pullback diagram
$$ \xymatrix{ \calN \ar[r] \ar[d]^{p} & \overline{\End}^{\otimes}(\calC, \calD) \ar[d] \\
\Nerve( \cDelta_{S})^{op} \ar[r]^{U} & \End^{\otimes}(\calC, \calD). }$$
Using Lemma \ref{mikest}, we deduce:
\begin{itemize}
\item[$(a1)$] The map $p$ is a locally coCartesian fibration.
\item[$(a2)$] Suppose that $\alpha: ([m],c) \rightarrow ([n],c')$ is a morphism in
$\cDelta_{S}$ which induces a {\em convex} morphism $[m] \rightarrow [n]$, and let
$f$ be a morphism in $\calN$ which projects to $\alpha$. Then $f$ is $p$-coCartesian if and only if $f$ is locally $p$-coCartesian.
\end{itemize}

Let $\calI$ be the subcategory of $[1] \times \cDelta_{S}^{op}$ defined as follows:
\begin{itemize}
\item An object $(i,[n],c)$ of $[1] \times \cDelta_{S}^{op}$ belongs to $\calI$
if and only if either $i = 1$ or $c(n) = \calD$.
\item Given a pair of objects $(i, [m], c), (j, [n], c') \in \calI$, a morphism
$(i, [m], c) \rightarrow (j, [n], c')$ in $[1] \times \cDelta_{S}^{op}$ belongs to
$\calI$ if and only if either $j=1$, or the induced map $\alpha: [n] \rightarrow [m]$ satisfies
$\alpha(n) = m$.
\end{itemize}

Let $\calI_0, \calI_1 \subseteq \calI$ be the preimages of the objects $0,1 \in [1]$, respectively, so we have canonical identifications $\calI_0 \subseteq \cDelta_{S}^{op} \simeq \calI_1$.
Let 
$$\calE = \bHom_{ \Nerve( \cDelta_{S})^{op}}( \Nerve(\calI), \calN)$$
$$\calE_0 = \bHom_{ \Nerve( \cDelta_{S})^{op}}( \Nerve(\calI_0), \calN) \quad \calE_1 = \bHom_{ \Nerve( \cDelta_{S})^{op}}( \Nerve(\calI_1), \calN).$$
Then we have restriction maps $ \calE_0 \leftarrow \calE \rightarrow \calE_1.$
Moreover, $\calE_1$ can be identified with the $\infty$-category of sections of $p$.
Consequently, Remarks \ref{recaster} and \ref{recaster2} allow us to identify $\overline{\AdjDiag}_{U}(\calC, \calD)$ with the full subcategory of
$\calE_1$ spanned by those sections $u: \Nerve(\cDelta_{S})^{op} \rightarrow \calN$
which have the following properties:
\begin{itemize}
\item[$(b1)$] Let $( [n], c) \in \cDelta_{S}$ be such that
$c(n) = \calD$, let $c_0 = c | \{n\}$, and let $\alpha: ( \{n\}, c_0) \rightarrow ([n], c)$ be the inclusion. Then $u(\alpha)$ is a locally $p$-coCartesian morphism of $\calN$.

\item[$(b2)$] Let $n > 0$, let $c: [n] \rightarrow \{ \calC, \calD \}$ satisfy 
$c(n) = \calD$, let $c_0 = c | \{0, \ldots, n-1\}$, and let $\alpha: ( [n-1], c_0) \rightarrow
( [n],c)$ be the inclusion. Then $u(\alpha)$ is a $p$-coCartesian morphism of
$\calN$.  
\end{itemize}

Note that, for every object $(1, [n], c) \in \calI_1$, the $\infty$-category
$\Nerve(\calI_0)_{/ (1, [n],c)}$ has a final object, given by $(0, [n+1], c')$, where $c'$ is defined by the formula
$$ c'(i) = \begin{cases} c(i) & \text{if } i \leq n \\
\calD & \text{if } i = n+1. \end{cases}$$
Likewise, for every object $(0, [n], c) \in \calI_0$, the $\infty$-category
$\Nerve(\calI_1)_{(0,[n],c)/}$ has an initial object, given by $(1, [n], c)$. Combining
this observation, $(a1)$, $(a2)$, and Proposition \toposref{lklk}, we deduce the following:

\begin{itemize}
\item[$(c1)$] Every functor $f_0 \in \calE_0$ admits a $p$-left Kan extension
$f: \Nerve(\calI) \rightarrow \calN$. Moreover, an arbitrary extension
$f': \Nerve(\calI) \rightarrow \calN$ of $f_0$ is a $p$-left Kan extension of $f_0$
if and only if it is compatible with the projection to $\Nerve( \cDelta_{S})^{op}$ and
carries each morphism $\alpha: (0, [n+1],c') \rightarrow (1, [n], c)$ to a locally $p$-coCartesian edge of $\calN$, where $c'$ is defined as above. Moreover, the projection
$\calE^0 \rightarrow \calE_0$ is a trivial Kan fibration, where $\calE^0 \subseteq \calE$ denotes the full subcategory spanned by those functors $f$ which are $p$-left Kan extensions of
$f | \Nerve(\calI_0)$.

\item[$(c2)$] Every functor $f_1 \in \calE_1$ admits a $p$-right Kan extension
$f: \Nerve(\calI) \rightarrow \calN$. Moreover, an arbitrary extension
$f': \Nerve(\calI) \rightarrow \calN$ of $f_1$ is a $p$-right Kan extension of $f_1$
if and only if it is compatible with the projection to $\Nerve( \cDelta_{S})^{op}$ and
carries each morphism $\beta: (0, [n],c) \rightarrow (1, [n], c)$ to an equivalence in $\calN$. Moreover, the projection
$\calE^1 \rightarrow \calE_1$ is a trivial Kan fibration, where $\calE^1 \subseteq \calE$ denotes the full subcategory spanned by those functors $f$ which are $p$-right Kan extensions of
$f | \Nerve(\calI_1)$.
\end{itemize}

Applying Proposition \toposref{leftkanadj}, we obtain a pair of adjoint functors $\Adjoint{F}{\calE_0}{\calE_1}{G},$ where:
\begin{itemize}
\item The functor $F$ is given by a composition
$ \calE_0 \stackrel{\psi}{\rightarrow} \calE^0 \subseteq \calE \rightarrow \calE_1$
where $\psi$ is a section of the trivial Kan fibration $\calE^0 \rightarrow \calE_0$ (a left Kan extension functor) and the last map is given by restriction.
\item The functor $G$ is given by a composition
$ \calE_1 \stackrel{\psi'}{\rightarrow} \calE^1 \subseteq \calE \rightarrow \calE_0$
where $\psi'$ is a section of the trivial Kan fibration $\calE^1 \rightarrow \calE_1$ (a right Kan extension functor) and the last map is given by restriction. If we choose $\psi'$ to be given by composition with the retraction of $\calI$ onto $\calI_1$, then we can identify $G$ with the
restriction map $\calE_1 \rightarrow \calE_0$.
\end{itemize}

Unwinding the definitions, and using the fact that $U([\calD, \calD])$ is equivalent to the identity map from $\calD$ to itself (see Remark \ref{humberton}), we deduce that the unit 
$\id_{\calE_0} \rightarrow G \circ F$ is an equivalence. Consequently, the functor $F$ is fully faithful. The essential image of $F$ can be identified with the full subcategory
$\calE'_1 \subseteq \calE_1$ spanned by those sections of $p$ which satisfy condition $(b2)$.

Form a pullback diagram
$$ \xymatrix{ \calN_0 \ar[r] \ar[d]^{p_0} & \calN \ar[d]^{p} \\
\Nerve(\calI_0) \ar@{^{(}->}[r] & \Nerve( \cDelta_{S})^{op}. }$$
If $u \in \calE_0$, then $F(u) \in \calE'_1$ satisfies $(b1)$ if and only if $u$ satisfies the 
obvious analogue of $(b1)$. In view of Lemma \ref{mikest}, this can be reformulated as follows:

\begin{itemize}
\item[$(b1')$] Let $( [n], c) \in \calI_0$, let $c_0 = c | \{n\}$, and let $\alpha: ( \{n\}, c_0) \rightarrow ([n], c)$ be the inclusion. Then $u(\alpha)$ is a $p_0$-Cartesian morphism of $\calN_0$.
\end{itemize}

Let $\calE''_{0} \subseteq \calE_0$ be the full subcategory spanned by those functors $u$ which satisfy condition $(b1')$. Let $\calI'_0 \subseteq \calI_0$ be the full subcategory spanned by the object $[\calD] \in \cDelta_{S}$.
We observe that condition $(b1')$ is simply a reformulation of the requirement that
$u$ be a $p$-right Kan extension of $u | \Nerve(\calI'_0)^{op}$. Moreover, Lemma \ref{mikest} implies that 
every partial section $u_0: \Nerve(\calI'_0)^{op} \rightarrow \calN_0$ admits a $p$-right Kan extension. It follows from Proposition \toposref{lklk} that the restriction map
$\calE''_0 \rightarrow \bHom_{ \Nerve(\calI_0)^{op} }(  \Nerve(\calI'_0)^{op}, \calN_0)
\simeq \calD$ is a trivial Kan fibration, hence a categorical equivalence as desired.
\end{proof}

\subsection{The Barr-Beck Theorem}\label{monoid7}

Let $F: \calC \rightarrow \calD$ be a functor between $\infty$-categories which admits a right adjoint $G$. In \S \ref{moadthird}, we saw that $F$ determines an adjunction datum $U \in \AdjDiag(\calC, \calD)$, which restricts to give a monad $T \in \Alg( \Fun(\calC, \calC))$. Moreover, the functor
$G$ factors canonically (up to homotopy) as a composition
$ \calD \stackrel{G'}{\rightarrow} \Mod_{T}(\calC) \rightarrow \calC.$
Our goal in this section is to establish a criterion which can be used to test whether $G'$ is an equivalence of $\infty$-categories. First, we need to introduce a bit of notation.

\begin{notation}
The category $\cDelta_{-\infty}$ is defined as follows: the objects of $\cDelta_{\infty}$ are
integers $n \geq -1$, and $\Hom_{ \cDelta_{-\infty} }( [m], [n])$ is the set of monotone maps
$[m] \cup \{ -\infty\} \rightarrow [n] \cup \{ -\infty \}$ which preserve the base point $-\infty$ (which
is regarded as a least element of both $[m] \cup \{-\infty \}$ and $[n] \cup \{ -\infty \}$). 
We have inclusions of subcategories
$ \cDelta \subseteq \cDelta_{+} \subseteq \cDelta_{-\infty},$
where the latter identifies $\cDelta_{+}$ with the subcategory of $\cDelta_{-\infty}$ having the same objects, and a map $f: [m] \cup \{ -\infty\} \rightarrow [n] \cup \{-\infty\}$ belongs to 
$\cDelta_{+}$ if and only if $f^{-1}( -\infty) = \{ -\infty \}$. 
\end{notation}

\begin{definition}\label{umpa}\index{split simplicial object}\index{simplicial object!split}\index{simplicial object!$G$-split}
Let $\calC$ be an $\infty$-category. We will say that an augmented simplicial object
$U: \Nerve(\cDelta_{+})^{op} \rightarrow \calC$ is {\it split} if $U$ extends to a functor
$\Nerve( \cDelta_{-\infty})^{op} \rightarrow \calC$. We will say that a simplicial object
$U: \Nerve(\cDelta)^{op} \rightarrow \calC$ is {\it split} if it extends to a split augmented simplicial object. Given a functor $G: \calD \rightarrow \calC$, we will say that an (augmented) simplicial object $U$ of $\calD$ is {\it $G$-split} if $G \circ U$ is split, when regarded as a simplicial object of $\calC$.
\end{definition}

\begin{remark}
Using the terminology of Definition \ref{umpa}, we can formulate Lemma \toposref{aclock} as follows: every split augmented simplicial object is a colimit diagram. (In \cite{topoi}, we used a slightly different notation, employing a category $\cDelta_{\infty}$ in place of the
category $\cDelta_{-\infty}$ defined above. However, these two categories are equivalent to one another, via the functor which reverses order.)
\end{remark}

\begin{remark}\label{bsplit}
Let $F: \calC \rightarrow \calD$ be a functor between $\infty$-categories, and let
$V$ be a split simplicial object of $\calC$. Then $F \circ V$ is a split simplicial object of $\calD$.
Lemma \toposref{aclock} implies that $V$ admits a colimit in $\calC$, and that $F$ preserves that colimit.
\end{remark}

\begin{theorem}[$\infty$-Categorical Barr-Beck Theorem]\label{barbeq}\index{Barr-Beck theorem}
Let $\calC$ and $\calD$ be $\infty$-categories, and suppose given an adjunction datum $U \in \AdjDiag(\calC, \calD)$. Let $G: \calD \rightarrow \calC$ be the induced functor, $T$ the induced monad on $\calC$, and let $\psi: \overline{\AdjDiag}_{U}(\calC,\calD) \rightarrow \Mod_{T}(\calC)$ be the forgetful functor $($here $\overline{\AdjDiag}_{U}(\calC, \calD)$ is canonically equivalent to $\calD$, in view of Corollary \ref{hooku}$)$. Then $\psi$ is an equivalence of categories if and only if the following conditions are satisfied:
\begin{itemize}
\item[$(1)$] The functor $G: \calD \rightarrow \calC$ is {\em conservative}; that is,
a morphism $f: D \rightarrow D'$ in $\calD$ is an equivalence if and only if
$G(f)$ is an equivalence in $\calC$.
\item[$(2)$] Let $V$ be a simplicial object of $\calD$ which is $G$-split. Then
$V$ admits a colimit in $\calD$, and that colimit is preserved by $G$.
\end{itemize}
\end{theorem}

We will give the proof at the end of this section.

\begin{remark}\label{sonica}
Hypotheses $(1)$ and $(2)$ of Theorem \ref{barbeq} can be rephrased as the following single condition:
\begin{itemize}
\item[$(\ast)$] Let $V: \Nerve(\cDelta)^{op} \rightarrow \calD$ be a simplicial object of $\calD$ which is $G$-split. Then $V$ has a colimit in $\calD$. Moreover, an arbitrary extension
$\overline{V}: \Nerve(\cDelta)^{op} \rightarrow \calD$ is a colimit diagram if and only if $G \circ \overline{V}$ is a colimit diagram.
\end{itemize}
It is clear that $(1)$ and $(2)$ imply $(\ast)$, and that $(\ast)$ implies $(2)$. To prove that
$(\ast)$ implies $(1)$, let us consider an arbitrary morphism $f: D' \rightarrow D$ in $\calD$.
The map $f$ determines an augmented simplicial object $\overline{V}$ of $\calD$, with $$\overline{V}([n]) = \begin{cases} D' & \text{if } n \geq 0 \\
D & \text{if } n=-1. \end{cases}$$
The underlying simplicial object $V = \overline{V} | \Nerve(\cDelta)^{op}$ is constant, 
and therefore $G$-split. Since $\Nerve(\cDelta)$ is contractible, $\overline{V}$ is a colimit diagram
if and only if $f$ is an equivalence, and $G \circ \overline{V}$ is a colimit diagram if and only if
$G(f)$ is an equivalence. If $(\ast)$ is satisfied, then these conditions are equivalent, so that
$G$ is conservative as desired.
\end{remark}

\begin{remark}
Let $\calC$ be an $\infty$-category, and regard $\calC$ as left-tensored over $\End(\calC)$.
We have a commutative diagram
$$ \xymatrix{ \Mod(\calC) \ar[rr]^{\theta} \ar[dr]^{p} & & \calC \times \Alg(\End(\calC)) \ar[dl]^{p'} \\
& \Alg(\End(\calC)), & }$$
where $\theta$ is determined by the $p$ and the forgetful functor $\Mod(\calC) \rightarrow \calC$.
The functor $p'$ is obviously a Cartesian fibration, and Corollary \ref{thetacart} implies that
$p$ is a Cartesian fibration as well. Proposition \ref{eggur} implies that $\theta$ carries $p$-Cartesian edges to $p'$-Cartesian edges. Consequently, $\theta$ classifies a natural
transformation of functors $F \rightarrow F'$, where $F: \Alg(\End(\calC))^{op} \rightarrow \Cat_{\infty}$ is the functor classified by $p$ (so that $F(A) \simeq \Mod_{T}(\calC)$ for
every monad $T \in \Alg(\End(\calC))$) and $F'$ is the constant functor taking the value $\calC \in \Cat_{\infty}$. We may identify this transformation with a functor $\alpha: \Alg(\End(\calC))^{op} \rightarrow \Cat_{\infty}^{/ \calC}$. We can interpret Theorem \ref{barbeq} as describing the essential image of the functor $\alpha$: namely, a functor $G: \calD \rightarrow \calC$ belongs to the essential image of $\alpha$ if and only if $G$ admits a left adjoint and satisfies conditions
$(1)$ and $(2)$ of Theorem \ref{barbeq}. In a future paper, we will show that the functor $\alpha$ is fully faithful. In other words, we may identify monads on $\calC$ with (certain) $\infty$-categories lying over $\calC$.
\end{remark}

In practice, it is often easier to apply the following slightly weaker version of the Barr-Beck theorem:

\begin{corollary}\label{usualbb}
Let $G: \calD \rightarrow \calC$ and $\psi: \calD \rightarrow \Mod_{T}(\calC)$ be as in Theorem \ref{barbeq}. Suppose that:
\begin{itemize}
\item[$(1)$] The functor $G$ is conservative.
\item[$(2)$] The $\infty$-category $\calD$ admits geometric realizations of simplicial objects.
\item[$(3)$] The functor $G$ preserves geometric realizations of simplicial objects.
\end{itemize}
Then the functor $\psi$ is an equivalence.
\end{corollary}

The Barr-Beck theorem is an extremely useful tool. However, for many applications (such as our construction of coproducts of algebras in \S \ref{limalg}) it will suffice to know the following less precise statement (which makes no mention of the theory of monads):

\begin{proposition}\label{littlebeck}
Let $\Adjoint{F}{\calC}{\calD}{G}$ be a pair of adjoint functors between $\infty$-categories which
satisfies conditions $(1)$ and $(2)$ of Theorem \ref{barbeq}. For every object $D \in \calD$, there exists a simplicial object $D_{\bigdot}: \Nerve(\cDelta)^{op} \rightarrow \calD$ having colimit $D$, such that each $D_{n}$ lies in the essential image of $F$.
\end{proposition}

\begin{proof}
Using Theorem \ref{partA} and Proposition \ref{stumpus}, we may assume that
the functors $F$ and $G$ are determined by an adjunction datum
$U \in \AdjDiag( \calC, \calD)$, and that $D \in \calD$ is the image of $[ \calD ] \in \cDelta_{S}$ under a pointed adjunction datum $\overline{U}: \Nerve( \cDelta_{S})^{op} 
\rightarrow \overline{\End}^{\otimes}(\calC, \calD)$ which lies over $U$.

We define a functor $\phi: \cDelta_{+} \rightarrow \cDelta_{S}$ by the formula
$$ [n] \mapsto ([1] \star [n], c),$$
where $c$ is given by the equation
$$ c(i) = \begin{cases} \calC & \text{if } i = 0 \\
\calD & \text{if } i = 1 \\
\calC & \text{if } 1 < i \leq n+2. \end{cases}$$
Composing with $\overline{s}$ and evaluating at $1 \in [n+2] \simeq [1] \star [n]$, we obtain
an augmented simplicial object $D_{\bigdot}: \Nerve( \cDelta_{+})^{op} \rightarrow \calD$.
It follows immediately from Definition \ref{termus} that $D_{-1} \simeq D$. To complete the proof, it will suffice to show that $D_{\bigdot}$ is a colimit diagram. In view of Remark \ref{sonica}, it will suffice to show that $D_{\bigdot}$ is $G$-split.  

Let $C_{\bigdot}: \Nerve( \cDelta_{+})^{op} \rightarrow \calC$ be the functor obtained from $\psi$ by composing with $\overline{U}$ and evaluating at $0 \in [n+2]$. Then $C_{\bigdot}$ is the image
of $D_{\bigdot}$ under a functor $G_{\bigdot}: \Nerve(\cDelta_{+})^{op} \rightarrow
\Fun( \calD, \calC)$ determined by the augmentation datum $s$. Invoking the definition of an augmentation datum, we conclude that $G_{\bigdot}$ is equivalent to the constant augmented simplicial object taking the value $G \in \Fun(\calD, \calC)$. It will therefore suffice to show that
the augmented simplicial object $C_{\bigdot}$ is split.

Let $\phi': \cDelta_{+} \rightarrow \cDelta_{S}$ be defined by the formula
$[n] \mapsto ( [0] \star [n], c'),$
where $c'$ takes the constant value $\calC$. Let $C'_{\bigdot}$ be the augmented simplicial
object of $\calC$ obtained by composing $\phi'$ with $\overline{s}$ and evaluating
at $0 \in [n+1] \simeq [0] \star [n]$. There is an evident natural transformation of functors
$\phi' \rightarrow \phi$, which induces a map of simplicial objects
$C_{\bigdot} \rightarrow C'_{\bigdot}$. Invoking the definition of a pointed adjunction datum, we deduce that this transformation is an equivalence. It will therefore suffice to prove that $C'_{\bigdot}$ is split. For this, it suffices to show that $\phi'$ can be extended to a functor
$\cDelta_{-\infty} \rightarrow \cDelta_{S}$, which is evident.
\end{proof}

\begin{remark}
The proof of Proposition \ref{littlebeck} actually establishes a stronger result: namely, 
the simplicial resolution of an object $D \in \calD$ can be chosen functorially in $D$.
Informally, we can describe this functorial choice as follows. Let $V = F \circ G$. A dual to Theorem \ref{partA} implies that $V$ underlies a {\it comonad} on the $\infty$-category $\calD$. Consequently, given $D \in \calD$ we can construct an augmented simplicial object $E_{\bigdot}$, with $D_{n} \simeq V^{n+1} D$. This augmented simplicial object is $G$-split, and therefore a colimit diagram. Consequently, $D \simeq D_{-1}$ can be obtained as the colimit of
the underlying simplicial object $D_{\bigdot}$, and for $n \geq 0$ we have
$D_{n} \simeq F( (G \circ V^n) D)$, so that $D_{n}$ belongs to the essential image of $F$.
\end{remark}

We conclude with few applications of Proposition \ref{littlebeck}:

\begin{corollary}\label{littlerbeck}
Suppose given a commutative diagram of $\infty$-categories
$$ \xymatrix{ \calC \ar[rr]^{U} \ar[dr]^{G} & & \calC' \ar[dl]^{G'} \\
& \calD. & }$$
Suppose furthermore that:
\begin{itemize}
\item[$(a)$] The $\infty$-categories $\calC$ and $\calC'$ admit geometric realizations
for simplicial objects.
\item[$(b)$] The functors $G$ and $G'$ admit left adjoints, which we will denote by $F$ and $F'$, respectively.
\item[$(c)$] The functor $G'$ is conservative and preserves geometric realizations of simplicial objects. 
\end{itemize}
Then:
\begin{itemize}
\item[$(1)$] The functor $U$ admits a left adjoint $T$.
\item[$(2)$] The functor $U$ is an equivalence if and only if the following additional conditions are satisfied:
\begin{itemize}
\item[$(d)$] The functor $G$ is conservative and preserves geometric realizations of simplicial objects.
\item[$(e)$] The unit of the adjunction $\Adjoint{T}{\calC'}{\calC}{U}$ induces an equivalence of functors 
$$ G' \circ F' \rightarrow G' \circ U \circ T \circ F' \simeq G \circ F.$$
\end{itemize}
\end{itemize}
\end{corollary}

\begin{proof}
We first prove $(1)$. Let $j': \calC' \rightarrow \Fun(\calC', \SSet)^{op}$ denote the Yoneda embedding, and consider the composition
$$\overline{j}': \calC' \stackrel{j'}{\rightarrow} \Fun(\calC', \SSet)^{op} \stackrel{ \circ U}{\rightarrow}
\Fun(\calC, \SSet)^{op}.$$ 
Let $\calC'_0 \subseteq \calC$ be the full subcategory spanned by those objects $C' \in \calC'$ for which
$\overline{j}'(C')$ belongs to the essential image of the Yoneda embedding $j: \calC \rightarrow
\Fun(\calC, \SSet)^{op}$. We wish to prove that $\calC'_0 = \calC'$.

Proposition \toposref{yonedaprop} and assumption $(a)$ guarantee that the essential image of $j$ is stable under the formation of geometric realizations of simplicial objects. Moreover, Proposition \toposref{yonedaprop} also implies that $\overline{j}'$ preserves geometric realizations of simplicial objects. It follows that $\calC'_0$ is stable under the formation of geometric realizations in $\calC'$.
In view of Proposition \ref{littlebeck}, it will suffice to show that $\calC'_0$ contains the essential image of $F'$. This is clear, since $\overline{j}'( F'(D) ) \simeq j( F(D))$. This completes the proof of $(1)$.

We now prove $(2)$. If $U$ is an equivalence, then $(e)$ is obvious (since the unit itself is an equivalence) and $(d)$ follows from $(c)$. To prove the converse, choose compatible unit and counit transformations 
$$u: \id \rightarrow U \circ T \quad v: T \circ U \rightarrow \id.$$
We wish to prove that $u$ and $v$ are equivalences. Since $G$ is conservative, the functor $U$ is also conservative. Consequently, the transformation $v$ is an equivalence if and only if
$U(v): U \circ T \circ U \rightarrow U$ is an equivalence. Since $U(v)$ has a section
determined by $u$, it will suffice to prove that $u$ is an equivalence.

Let $\calC'_1$ be the full subcategory of $\calC'$ spanned by those objects
$C$ for which the map $u(C): C \rightarrow (U \circ T)(C)$ is an equivalence. 
Since $G = G' \circ U$ and $G'$ both preserve geometric realizations and $G'$ is conservative, we conclude that $U$ preserves geometric realizations. The functor $T$ is a left adjoint, and therefore preserves geometric realizations. It follows that $U \circ T$ preserves geometric realizations, so 
that $\calC'_{1}$ is stable under geometric realizations in $\calC'$. In view of Propsition \ref{littlebeck}, it will suffice to show that $\calC'_{1}$ contains the essential image of $F'$.
In other words, we must show that $u$ induces an equivalence $F' \rightarrow U \circ T \circ F'$. 
Using our assumption that $G'$ is conservative, we reduce immediately to $(e)$.
\end{proof}

Our next result makes use of the notation and terminology of \S \stableref{stable11}.

\begin{corollary}\label{progen}
Suppose given a pair of adjoint functors between $\infty$-categories
$\Adjoint{F}{\calC}{\calD}{G}$.
Assume that:
\begin{itemize}
\item[$(i)$] The $\infty$-category $\calD$ admits filtered colimits and geometric realizations, and 
$G$ preserves filtered colimits and geometric realizations.
\item[$(ii)$] The $\infty$-category $\calC$ is projectively generated (Definition \stableref{defpro}).
\item[$(iii)$] The functor $G$ is conservative.
\end{itemize}
Then:
\begin{itemize}
\item[$(1)$] The $\infty$-category $\calD$ is projectively generated.
\item[$(2)$] An object $D \in \calD$ is compact and projective if and only if there exists
a compact projective object $C \in \calC$ such that $D$ is a retract of $F(C)$.
\item[$(3)$] The functor $G$ preserves all sifted colimits.
\end{itemize}
\end{corollary}

\begin{proof}
Let $\calC^{0}$ denote the full subcategory of $\calC$ spanned by the compact projective objects.
Let $\calD^{0}$ denote the essential image of $F | \calC^{0}$. Using
assumption $(i)$, we deduce that $\calD^{0}$ consists of compact projective objects of $\calD$.
Without loss of generality, we may assume that the $\infty$-category $\calD$ is minimal, so that
$\calD^{0}$ is small. Moreover, since $\calC^{0}$ is stable under finite coproducts in $\calC$, and
$F$ preserves finite coproducts, we conclude that $\calD^{0}$ admits finite coproducts (which are also finite coproducts in $\calD$). 

Let $\calD' = \calP_{\Sigma}( \calD^{0} )$ (see Definition \stableref{vardef}). Using
Proposition \stableref{surottt}, we deduce that the inclusion $\calD^{0} \subseteq \calD$
is homotopic to a composition
$$ \calD^{0} \stackrel{j}{\rightarrow} \calD' \stackrel{f}{\rightarrow} \calD,$$
where the functor $f$ preserves filtered colimits and geometric realizations. Combining Proposition \stableref{smearof} with Proposition \ref{littlebeck} and assumption $(ii)$, we conclude that $f$ is an equivalence of $\infty$-categories. This proves $(1)$. Moreover, the proof shows that
$\calD^{0}$ is spanned by a set of compact projective generators for $\calD$ (Definition \stableref{defpro}), so that assertion $(2)$ follows from Proposition \stableref{protus}. Assertion $(3)$ now follows from 
Proposition \stableref{surottt}.
\end{proof}

We now return to the proof of Theorem \ref{barbeq}. One direction is essentially contained in the following Lemma:

\begin{lemma}\label{hasbeen}
Let $\calC$ be a monoidal $\infty$-category, $\calM$ an $\infty$-category which is left-tensored over $\calC$, $A$ an algebra object, and $\theta: \Mod_{A}(\calM) \rightarrow \calM$ the forgetful functor. Then:
\begin{itemize}
\item[$(1)$] Every $\theta$-split simplicial object of $\Mod_{A}(\calM)$ admits a colimit in $\calM$.
\item[$(2)$] The functor $\theta$ preserves colimits of $\theta$-split simplicial objects.
\end{itemize}
\end{lemma}

\begin{proof}
Note that a simplicial object of $\Mod_{A}(\calM)$ can be viewed as a {\em bisimplicial object}
of $\calM^{\otimes}$. In order to avoid confusion, we let $K$ denote the simplicial set
$\Nerve(\cDelta)^{op}$ when we wish to emphasize the role of $\Nerve(\cDelta)^{op}$ as indexing simplicial objects, while we use the notation $\Nerve(\cDelta)^{op}$ to denote the base of the coCartesian fibration $\calM^{\otimes} \rightarrow \Nerve(\cDelta)^{op}$.

Form a pullback diagram
$$ \xymatrix{ \calN \ar[d]^{p} \ar[r] & \calM^{\otimes} \ar[d] \\
\Nerve(\cDelta)^{op} \ar[r]^{A} & \calC^{\otimes}, }$$
so that $p$ is a locally coCartesian fibration (Lemma \ref{excel}). 
Let $u: K \rightarrow \Mod_{A}(\calM)$ be a $\theta$-split simplicial object, corresponding to a map
$V: K \times \Nerve(\cDelta)^{op} \rightarrow \calN$. We observe that every fiber of $p$ can be identified with
$\calM$, and each of the induced maps
$V_{[n]}: K \times \{ [n] \} \rightarrow \calN_{[n]}$
can be identified with the composition of $V$ with the forgetful functor $\theta$.
It follows that each of the simplicial objects $V_{[n]}$ splits. Using Lemma \toposref{aclock}, we deduce that each $V_{[n]}$ admits a colimit 
$\overline{V}_{[n]}: K^{\triangleright} \rightarrow \calN_{[n]}$, and that these colimits are preserved by each of the associated functors $\calN_{[n]} \rightarrow \calN_{[m]}$ (Remark \ref{bsplit}).
Applying Proposition \toposref{relcolfibtest}, we conclude that each $\overline{V}_{[n]}$ is a $p$-colimit diagram. We now invoke Lemma \ref{surtybove} to deduce the existence of a map
$\overline{V}: K^{\triangleright} \times \Nerve(\cDelta)^{op} \rightarrow \calN$
which induces a $p$-colimit diagram on each fiber $K^{\triangleright} \times \{ [n] \}$. 
Moreover, Lemma \ref{surtybove} implies that $\overline{V}$ defines a colimit diagram
$\overline{v}: K^{\triangleright} \rightarrow \bHom_{ \Nerve(\cDelta)^{op}}( \Nerve(\cDelta)^{op}, \calN)$. Since $v = \overline{v}|K$ factors through the full subcategory 
$\Mod_{A}(\calM) \subseteq \bHom_{ \Nerve(\cDelta)^{op}}( \Nerve(\cDelta)^{op}, \calN)$, it is easy to see that $\overline{v}$ factors through $\Mod_{A}(\calM)$. It follows that
$\overline{v}$ is a colimit of $v$ in $\Mod_{A}(\calM)$. This proves $(1)$, and $(2)$ follows from the observation that $\theta \circ \overline{v} = \overline{V} | K^{\triangleright} \times \{ [0] \}$ is a colimit diagram in $\calM$ by construction.
\end{proof}

\begin{proof}[Proof of Theorem \ref{barbeq}]
The ``only if'' direction follows from Corollary \ref{goop} and Lemma \ref{hasbeen}. For the converse, let us suppose that $(1)$ and $(2)$ are satisfied. Set $S = \{ \calC, \calD \}$, and
form a pullback diagram
$$ \xymatrix{ \calN \ar[r] \ar[d]^{p} & \overline{\End}^{\otimes}(\calC, \calD) \ar[d] \\
\Nerve(\cDelta_{S})^{op} \ar[r]^{s} & \End^{\otimes}(\calC, \calD). }$$
We define full subcategories
$\cDelta_{\calC} \subseteq \cDelta_{S}^{0} \subseteq \cDelta_{S}$
as follows:
\begin{itemize}
\item An object $([n],c) \in \cDelta_{S}$ belongs to $\cDelta_{S}^{0}$ if and only if
$c(n) = \calC$.
\item An object $([n],c) \in \cDelta_{S}$ belongs to $\cDelta_{\calC}$ if and only if
$c$ is a constant map taking the value $\calC$.
\end{itemize}
Define $\infty$-categories
$$ \calE = \bHom_{ \Nerve( \cDelta_{S})^{op} }( \Nerve(\cDelta_{S})^{op}, \calN),$$
$$ \calE' = \bHom_{ \Nerve( \cDelta_{S})^{op} }( \Nerve(\cDelta_{S}^{0})^{op}, \calN), \quad \calE'' = \bHom_{ \Nerve( \cDelta_{S})^{op} }( \Nerve(\cDelta_{\calC})^{op}, \calN),$$
so that we have restriction maps
$$ \calE \rightarrow \calE' \rightarrow \calE''.$$
We observe that $\Mod_{T}(\calC)$ can be identified with a full subcategory of $\calE''$, and
$\overline{\AdjDiag}_{U}(\calC, \calD)$ can be identified with a full subcategory of $\calE$.
Let $\calE'_0 \subseteq \calE'$ be the full subcategory spanned by those functors
$f: \Nerve( \cDelta^{0}_{\calC, \calD} )^{op} \rightarrow \calN$ such that $f$ is $p$-right Kan extension of $f_0 = f | \Nerve( \cDelta_{\calC})^{op}$, and $f_0 \in \Mod_{T}(\calC)$. 
Let $\calE_0 \subseteq \calE$ be the full subcategory spanned by those functors
$g: \Nerve( \cDelta_{S})^{op} \rightarrow \calN$ such that $g$ is a $p$-left Kan extension of $g_0 = g | \Nerve(\cDelta_{S}^0)^{op}$, and $g_0 \in \calE'_0$. We will prove:
\begin{itemize}
\item[$(i)$] Every $f_0: \Nerve( \cDelta_{\calC})^{op} \rightarrow \calN$ which belongs to
$\Mod_{T}(\calC)$ admits a $p$-right Kan extension to $\Nerve( \cDelta_{S}^{0})^{op}$.
\item[$(ii)$] Every $g_0: \Nerve( \cDelta_{S}^{0})^{op} \rightarrow \calN$ whcih belongs to $\calE'_0$ admits a $p$-right Kan extension to $\Nerve( \cDelta_{S})^{op}$.
\end{itemize}
It will follow from Proposition \toposref{lklk} that the restriction maps
$\calE_0 \rightarrow \calE'_0 \rightarrow \Mod_{T}(\calC)$
are trivial Kan fibrations. We will then complete the proof by showing:
\begin{itemize}
\item[$(iii)$] The full subcategory $\calE_0 \subseteq \calE$ coincides with
$\overline{\AdjDiag}_{U}(\calC, \calD)$.
\end{itemize}

Observe that, for every object
$([n],c) \in \cDelta_{S}^{0}$, the $\infty$-category
$ \Nerve(\cDelta_{\calC}^{op})_{([n],c)/}$ has an initial object, given by a map
$\alpha: ([k],c') \rightarrow ([n],c)$ which induces a bijection
$[k] \rightarrow \{ i \in [n]: c(i) = \calC \}$. Since $([n],c) \in \cDelta_{S}^{0}$, the
map $\alpha$ satisfies $\alpha(k) = n$. Combining Lemmas \ref{mikest} and \toposref{kan2}, we deduce the following slightly stronger version of $(i)$: 

\begin{itemize}
\item[$(i')$] Every functor $f_0 \in \calE''$ admits a $p$-right Kan extension
$f: \Nerve(\cDelta_{S}^{0})^{op} \rightarrow \calN$. Moreover, an arbitrary functor
$f \in \calE'$ is a $p$-right Kan extension of $f | \Nerve( \cDelta_{\calC})^{op}$ if and only if
for every map $\alpha: ([k],c') \rightarrow ([n],c)$ in $\cDelta_{S}^{0}$ which induces a bijection $[k] \rightarrow \{ i \in [n]: c(i) = \calC \}$, the morphism $f(\alpha)$ is locally $p$-coCartesian.
\end{itemize}

Combining $(i')$, Proposition \toposref{lklk}, and the definition of $\Mod_{T}(\calC)$, we deduce:

\begin{itemize}
\item[$(i'')$] The restriction map $\calE'_0 \rightarrow \Mod_{T}(\calC)$ is a trivial Kan fibration.
Moreover, an arbitrary functor $f \in \calE'$ belongs to $\calE'_0$ if and only if, for every morphism
$\alpha: ([k],c') \rightarrow ([n],c)$ in $\cDelta_{S}$ such that
$\alpha(k) = n$, the morphism $f(\alpha)$ is locally $p$-coCartesian.
\end{itemize}

The proof of $(ii)$ is a bit more subtle. Fix an functor $g_0: \Nerve( \cDelta_{S}^{0})^{op} \rightarrow \calN$ which belongs to $\calE'_0$. We wish to prove that $g_0$ admits a $p$-left Kan extension to $\Nerve(\cDelta_{S})^{op}$. According to Lemma \toposref{kan2}, it will suffice to show that, for every object $([n], c) \in \cDelta_{S}$, the induced diagram
$ (\Nerve( \cDelta^{0}_{\calC, \calD} )^{op})_{/ ([n],c)}
\rightarrow \Nerve(\cDelta_{S}^{0})^{op} \stackrel{g_0}{\rightarrow}
\calN$ admits a $p$-colimit. If $([n],c) \in \cDelta_{S}^{0}$, then there is nothing to prove.
Assume instead that $c(n) = \calD$. We define a functor $\psi: \cDelta \rightarrow
(\cDelta_{S}^{0})_{([n],c)/}$ by the formula
$$ \psi( [k] ) = ( [n+k+1], c' )$$
$$c'(i) =  \begin{cases} c(i) & \text{if } 0 \leq i \leq n \\
\calC & \text{if } n+1 \leq i \leq n+k+1. \end{cases}$$
We claim that the induced map
$\Nerve(\cDelta)^{op} \rightarrow (\Nerve( \cDelta_{S}^{0})^{op})_{ /([n],c)}$
is cofinal. According to Theorem \toposref{hollowtt}, it will suffice to show that for every object
$ \alpha: ( [n], c) \rightarrow ([m], c')$ in $(\cDelta_{S}^{0})_{([n],c)/}$, the category
$$\cDelta' = \cDelta \times_{ (\cDelta_{S}^{0})_{([n],c)/} }
((\cDelta_{S}^{0})_{([n],c)/})_{ /([m],c')}$$
has weakly contractible nerve. Let $J = \{ i \in [m]: (\alpha(n) \leq i) \wedge (c'(i) = \calC) \}$.
Then $J$ contains $m$, and is therefore nonempty. Unwinding the definitions, we see that
$\cDelta'$ can be identified with $\cDelta_{/J}$, which has a final object.

Using Proposition \toposref{relexist}, we are reduced to proving that the diagram
$$ \xymatrix{ \Nerve(\cDelta^{op}) \ar@{^{(}->}[d] \ar[r]^{g_0 \circ \psi} & \calN \ar[d]^{p} \\
\Nerve(\cDelta^{op})^{\triangleright} \ar@{-->}[ur] \ar[r]^{\overline{h}} & \Nerve(\cDelta_{S})^{op} }$$
admits an extension as indicated, which is a $p$-colimit diagram. Let
$h = \overline{h} | \Nerve(\cDelta)^{op}$. We observe that $\overline{h}$ determines a 
natural transformation $H: \Nerve(\cDelta)^{op} \times \Delta^1 \rightarrow
\Nerve(\cDelta_{S})^{op}$ from $h$ to the constant map taking the value
$([n],c)$. Using Lemma \ref{mikest}, we can lift
$H$ to a $p$-coCartesian natural transformation
$\widetilde{H}: \Nerve(\cDelta)^{op} \times \Delta^1 \rightarrow \calN$
from $g_0 \circ \psi$ to a simplicial object $V: \Nerve(\cDelta)^{op} \rightarrow \calN_{([n],c)}$.
Using Proposition \toposref{chocolatelast}, we are reduced to proving that $V$ can be extended to a $p$-colimit diagram.

We now observe that the fiber $\calN_{([n],c)}$ can be identiied with $c(n) = \calD$. We will prove that $V$ is a $G$-split simplicial object of $\calD$. In view of assumption $(2)$, this will imply that
$V$ can be extended to a colimit diagram $\overline{V}: \Nerve(\cDelta^{op})^{\triangleright} \rightarrow \calN_{ ([n],c)} \simeq \calD$. We need to know a little bit more: namely, that
$\overline{U}$ is a $p$-colimit diagram in $\calN$. In view of Proposition \toposref{relcolfibtest}, it will suffice to show that for every morphism $\alpha: ([m], c') \rightarrow ([n], c)$ in $\cDelta_{S}$, the composition
$\Nerve(\cDelta^{op})^{\triangleright}
\stackrel{ \overline{V}}{\rightarrow} \calN_{([n],c)} 
\stackrel{\alpha_{!}}{\rightarrow} \calN_{([m],c')}$
is a colimit diagram. We observe that $\calN_{([m],c')}$ can be identified with the $\infty$-category $c'(m) \in \{ \calC, \calD\}$, and that $\alpha_{!}$ is equivalent to some iterated composition of the adjoint functors $F$ and $G$. If $\alpha_{!}$ is not an equivalence, then it is equivalent to a composition $B \circ G$, for some functor $B: \calC \rightarrow c'(m)$. Invoking assumption $(2)$, we deduce that
$G \circ \overline{V}$ is colimit diagram in $\calC$. Moreover, since $G \circ V$ is split, the
augmented simplicial object $G \circ \overline{V}$ is likewise split. It follows that the image
of $G \circ \overline{V}$ under any functor remains a split augmented simplicial object, and therefore a colimit diagram (Lemma \toposref{aclock}). 

It remains to show that $V$ is a $G$-split simplicial object of $\calD$. Our first step us to show that this assertion is independent of the original object $([n], c) \in \cDelta_{S}$. To prove this, consider the map $\beta: ([0], c') \rightarrow ([n], c)$, where $\beta(0) = n$ and
$c'(0) = \calD$. We define maps
$$\psi': \cDelta \rightarrow (\cDelta_{S}^{0})_{([0],c')/},$$ 
$$ H' : \Nerve(\cDelta)^{op} \times \Delta^1 \rightarrow \Nerve( \cDelta_{S})^{op}, \quad \widetilde{H}' : \Nerve(\cDelta)^{op} \times \Delta^1 \rightarrow \calN $$
as above, so that $\widetilde{H}'$ is a $p$-coCartesian transformation from
$g_0 \circ \psi'$ to a simplicial object $V': \Nerve(\cDelta) \rightarrow \calN_{([0],c')} \simeq \calD$.
We will show that $V$ and $V'$ are equivalent (as simplicial objects of $\calD$), so that
$V$ is $G$-split if and only if $V'$ is $G$-split. 

Let us be a bit more precise. Let $\overline{ ([n],c) }$ denote the constant functor
$\Nerve(\cDelta)^{op} \rightarrow \Nerve(\cDelta_{S})^{op}$ taking the value
$([n],c)$, and define $\overline{ ([0], c')}$ likewise. The morphism $\beta$ determines a 
natural transformation $\overline{\beta}: \overline{ ([n],c) } \rightarrow
\overline{ ([0],c') }$ in the $\infty$-category of simplicial objects
$\Fun( \Nerve(\cDelta)^{op}, \Nerve(\cDelta_{S})^{op} )$. Let
$\overline{p}: \Fun( \Nerve(\cDelta)^{op}, \calN) \rightarrow
\Fun( \Nerve(\cDelta)^{op}, \Nerve(\cDelta_{S})^{op} )$ be given by composition with $p$. Using Lemma \ref{mikest}, we can choose an $\overline{p}$-coCartesian transformation
$\overline{B}: V \rightarrow V''$ in $\Fun( \Nerve(\cDelta)^{op}, \calN)$. We will prove that
$V'$ and $V''$ are equivalent (as simplicial objects of $\calN_{([0],c')}$. 

The functors $\psi$ and $\psi'$ determine simplicial objects
$\Nerve(\cDelta)^{op} \rightarrow \Nerve( \cDelta_{S})^{op}$, which we will denote by
$\phi$ and $\phi'$, respectively. The map $\beta$ determines a natural transformation
$\widetilde{\beta}: \phi \rightarrow \phi'$. We have a commutative diagram
$$ \xymatrix{ \phi \ar[r]^{H} \ar[d]^{\widetilde{\beta}} & 
\overline{ ([n],c) } \ar[d]^{\overline{\beta}} \\
\phi' \ar[r]^{H'} & \overline{ ([0], c)}, }$$
which we can describe by a map $\Delta^1 \times \Delta^1 \rightarrow \Fun( \Nerve(\cDelta)^{op}, \Nerve( \cDelta_{S})^{op} )$. Using
Lemma \ref{mikest}, we deduce that $\overline{p}$ restricts to a coCartesian fibration
$$ \Fun( \Nerve(\cDelta)^{op}, \calN) \times_{ \Fun( \Nerve(\cDelta)^{op}, \Nerve(\cDelta_{S})^{op})} (\Delta^1 \times \Delta^1) \rightarrow \Delta^1 \times \Delta^1,$$
so that we obtain an associated diagram of $\infty$-categories and functors
$$ \xymatrix{ \Fun( \Nerve(\cDelta)^{op}, \calN)_{\phi} \ar[r]^{H_{!}} \ar[d]^{\widetilde{\beta}_!} &
\Fun( \Nerve(\cDelta)^{op}, \calN_{([n], c)}) \ar[d]^{ \overline{\beta}_!} \\
\Fun( \Nerve(\cDelta)^{op}, \calN)_{\phi'} \ar[r]^{H'_!} & \Fun( \Nerve(\cDelta)^{op}, \calN_{([0],c')})
}$$
which commutes up to (canonical) homotopy. Since $V$ is defined to be the image under
$H_!$ of $g_0 \circ \psi \in \Fun( \Nerve(\cDelta)^{op}, \calN)_{\phi}$ and
$V'$ is defined to be the image under $H'_{!}$ of $g_0 \circ \psi' \in \Fun( \Nerve(\cDelta)^{op}, \calN)_{\phi'}$, it will suffice to prove that $g_0 \circ \psi'$ is equivalent
to $\widetilde{\beta}_!( g_0 \circ \psi)$. The desired equivalence is provided by the
natural transformation
$g_0 \circ \widetilde{\beta}: \Nerve(\cDelta)^{op} \times \Delta^1 \rightarrow \calN,$
which is $p$-coCartesian in view of Lemma \ref{mikest}, $(i'')$, and our assumption that $g_0 \in \calE'_0$.

We now return to the task of proving that $V$ is $G$-split. Applying the above argument twice, we can reduce to the case where $n > 0$ and $c(n-1) = \calC$. Let $c_0 = c | [n-1]$, and let
$\gamma: ([n-1], c_0) \rightarrow ([n], c)$ be the associated map in $\cDelta_{S}$. 
Proceeding as above, we let $\overline{ ([n-1, c_0)}$ denote the constant functor
$\Nerve(\cDelta)^{op} \rightarrow \Nerve( \cDelta_{S})^{op}$ taking the value
$([n-1], c_0)$, $\overline{\gamma}: \overline{ ([n],c)} \rightarrow
\overline{ ([n-1], c_0)}$ the associated natural transformation in
$\Fun( \Nerve(\cDelta)^{op}, \Nerve(\cDelta_{S})^{op} ),$ and
$\overline{\gamma}_! : \Fun( \Nerve(\cDelta)^{op}, \calN_{([n],c)}
) \rightarrow \Fun( \Nerve(\cDelta)^{op}, \calN_{([n-1], c_0)}$
the functor associated to the locally coCartesian fibration $\overline{p}$. Then
$\overline{\gamma}_!$ can be identified with the functor
$ \Fun( \Nerve(\cDelta)^{op}, \calD) \rightarrow \Fun( \Nerve(\cDelta)^{op}, \calC)$
given by composition with $G$. It will therefore suffice to show that
$\overline{\gamma}_!(V)$ is a split simplicial object of
$\calN_{([n-1],c_0)}$. 

We have a commutative diagram
$$ \xymatrix{ & \overline{([n],c)} \ar[dr]^{ \overline{\gamma} } & \\
\phi \ar[ur]^{H} \ar[rr]^{\delta} & & \overline{([n-1], c_0}) }$$
described by a $2$-simplex $\Delta^2 \rightarrow \Fun( \Nerve(\cDelta)^{op}, \Nerve(\cDelta_{S})^{op})$. Applying Lemma \ref{mikest}, we deduce that $\overline{p}$ restricts to a coCartesian fibration
$ \Fun( \Nerve(\cDelta)^{op}, \calN) \times_{ \Fun( \Nerve(\cDelta)^{op}, 
\Nerve(\cDelta_{S})^{op} )} \Delta^2 \rightarrow \Delta^2,$
so that we have an equivalence $\delta_! \simeq \overline{\gamma}_! \circ H_!$ of functors
from $\Fun( \Nerve(\cDelta)^{op}, \calN)_{\phi} \rightarrow \Fun( \Nerve(\cDelta)^{op}, \calN_{([n-1], c_0)})$. It
will therefore suffice to show that $\delta_{!}( g_0 \circ \psi)$ is split.

Define $\theta: \cDelta_{-\infty} \rightarrow (\cDelta_{S})_{([n-1],c_0)/}$ as follows:
\begin{itemize}
\item The functor $\theta$ carries $[k]$ to $([n+k], c'_0)$, where
$$ c'_0(i) = \begin{cases} c_0(i) & \text{if } 0 \leq i < n \\
\calC & \text{if } n \leq i \leq n+k. \end{cases}$$
\item Given a map $\alpha: \{ - \infty \} \cup [k] \rightarrow \{ - \infty \} \cup [k']$
in $\cDelta_{-\infty}$, the associated map $\theta(\alpha)$ is given by the formula
$$ \theta(\alpha)(i) = \begin{cases} i & \text{if } i < n \\
n-1 & \text{if } i \geq n \text{ and } \alpha(i)= - \infty \\
n+\alpha(i-n) & \text{if } i \geq n \text{ and } \alpha(i) \neq - \infty. \end{cases}$$
\end{itemize}
Let $\overline{\theta}: \Nerve(\cDelta_{-\infty})^{op} \rightarrow
\Nerve( \cDelta_{S} )^{op}$ be the map determined by $\theta$, and let
$\epsilon$ be the induced natural transformation from $\overline{\theta}$ to the constant map
$\Nerve(\cDelta_{-\infty})^{op} \rightarrow \Nerve( \cDelta_{S})^{op}$ taking the value
$([n-1], c_0)$. Using Lemma \ref{mikest}, we can choose a $p$-coCartesian transformation
$\widetilde{\epsilon}: (g_0 \circ \overline{\theta}) \rightarrow W$, where $W$ is a map
$\Nerve(\cDelta_{-\infty})^{op} \rightarrow \calN_{ ([n-1],c_0)}$. We now observe that
$W$ is an extension of $\delta_!( g_0 \circ \psi)$, so that $\delta_1( g_0 \circ \psi)$ is split.
This completes the proof of $(ii)$.

We now prove $(iii)$. Suppose first that $g \in \calE$, and let $g_0 = g| \Nerve( \cDelta^{0}_{\calC, \calD})^{op}$. We wish to show that $g \in \calE_0$ if and only if $g \in \overline{\AdjDiag}_{U}(\calC, \calD)$. In view of $(i'')$, either of these conditions implies that $g_0 \in \calE'_0$. We will therefore assume that $g_0 \in \calE'_0$. By definition, $g \in \calE_0$ if and only if $g$ is a $p$-left Kan extension of $g_0$. Unwinding the definitions, this condition holds if and only if
$g$ is a $p$-left Kan extension of $g_0$ at $([n],c)$, for every object $([n],c) \in \cDelta_{S}$. This condition is automatic if $c(n) = \calC$. Suppose instead that $c(n) = \calD$. Let $\psi: \cDelta \rightarrow ( \cDelta^{0}_{\calC, \calD})_{([n],c)/}$ be defined as above. We observe that $\psi$ extends to a map $\overline{\psi}: \cDelta_{+} \rightarrow
(\cDelta_{S})_{([n],c)/}$, which preserves initial objects. Since $\psi$ induces a cofinal map $\Nerve(\cDelta)^{op} \rightarrow (\Nerve( \cDelta^{0}_{\calC, \calD})^{op}$, we conclude
that $g$ is a $p$-left Kan extension of $g_0$ at $([n],c)$ if and only if 
$g \circ \overline{\psi}$ is a $p$-colimit diagram. Using Lemma \ref{mikest}, we can construct a $p$-coCartesian transformation from $g \circ \overline{\psi}$ to an augmented simplicial object
$\overline{V}: \Nerve(\cDelta^{op})^{\triangleright} \rightarrow \calN_{([n],c)}$. In view of Proposition \toposref{chocolatelast}, $g$ is a $p$-left Kan extension of $g_0$ at $([n],c)$ if and only if 
$\overline{V}$ is a $p$-colimit diagram. The arguments given above show that this condition is satisfied if and only if $\overline{V}$ is a colimit diagram in $\calN_{([n],c)}$. Moreover, this condition is independent of the choice of $([n],c)$. In particular, we may suppose that
$n=1$ and $c$ is given by the formula
$$ c(i) = \begin{cases} \calC & \text{if } i = 0 \\
\calD & \text{if } i = 1. \end{cases} $$

Let $c_0$ denote the restriction of $c$ to $[0] \subseteq [1]$.
Using Lemma \ref{mikest}, we can choose a $p$-coCartesian transformation
$ \Nerve(\cDelta^{op})^{\triangleright} \times \Delta^1 \rightarrow \calN$
from $\overline{X}$ to an augmented simplicial object $\overline{X}: \Nerve(\cDelta^{op})^{\triangleright} \rightarrow \calN_{([0],c_0)}$. If we view $\overline{V}$ as an augmented simplicial object of $\calD$, then $\overline{X}$ is equivalent to the augmented simplicial object of $\calC$ obtained by composing $\overline{V}$ with the functor $G$. Since the underlying simplicial object of $\overline{V}$ is $G$-split, Remark \ref{sonica} implies that $\overline{V}$ is a colimit diagram
if and only $\overline{X}$ is a colimit diagram. Let $W: \Nerve(\cDelta_{-\infty})^{op} \rightarrow
\calN_{([0],c_0)}$ be the map constructed in the argument above, and let
$W_0 = W | \Nerve( \cDelta^{op})^{\triangleright}$ be the underlying augmented simplicial object.
Then $g$ determines a natural transformation of augmented simplicial objects
$\alpha: \overline{X} \rightarrow W_0$, which induces an equivalence $\alpha_{[k]}: \overline{X}( [k] ) \simeq W_0( [k] )$ for all $k \geq 0$. Lemma \toposref{aclock} implies that $W_0$ is a colimit diagram in $\calN_{([0],c_0)}$. Consequently, $\overline{X}$ is a colimit diagram if and only if
$\alpha_{[-1]}: \overline{X}([-1]) \rightarrow W_0([-1])$ is an equivalence. In other words, $g$ is a $p$-left Kan extension of $g_0$ if and only if $\alpha_{[-1]}$ is an equivalence. Unwinding the definitions, we observe that this latter condition is equivalent to the requirement that $g$ carry
the inclusion $([0], c_0) \rightarrow ([1], c)$ to a locally $p$-coCartesian morphism in $\calN$.
We now observe that, in view $(i'')$ and our assumption that $g_0 \in \calE'_0$, this last condition is equivalent to the requirement $g \in \overline{\AdjDiag}_{U}(\calC, \calD)$.
\end{proof}

\subsection{Existence of Adjunction Data}\label{digtwo}

Our goal in this section is to give the proof of Theorem \ref{partA}. Fix a pair of $\infty$-categories $\calC$ and $\calD$, and let $S = \{ \calC, \calD \}$. We wish to prove that every functor
$F: \calC \rightarrow \calD$ which admits a right adjoint can be promoted to an adjunction datum $U$. The proof is somewhat technical, and can be safely skipped by the reader; the ideas introduced here will not be needed elsewhere.

The passage from $F$ to $U$ will proceed in stages. To describe these stages in more detail, we need to introduce some terminology.

\begin{definition}\index{adjunction!nonunital}\index{nonunital!adjunction}
Let $\calC$ and $\calD$ be $\infty$-categories. A {\it nonunital adjunction} is
a pair of functors $F: \calC \rightarrow \calD$, $G: \calD \rightarrow \calC$, and 
a pair of maps $F \circ G \stackrel{v}{\rightarrow} U \stackrel{\alpha}{\leftarrow} \id_{\calD}$
in the $\infty$-category $\Fun(\calD, \calD)$, where $\alpha$ is an equivalence. We let
$\NUAdj(\calC, \calD)$ denote the full subcategory of
$$ (\Fun( \calC, \calD) \times \Fun( \calD, \calC) ) \times_{ \End(\calD) } 
\Fun( \Lambda^2_2, \End(\calD) ) \times_{ \End(\calD) } \{ \id_{\calD} \}$$
spanned by the nonunital adjunctions. We will refer to $\NUAdj(\calC, \calD)$ as the
{\it $\infty$-category of nonunital adjunctions}.
We let $\UAdj( \calC, \calD)$ denote the full subcategory of $\NUAdj(\calC, \calD)$ spanned
by those nonunital adjunctions
$$ F \circ G \stackrel{v}{\rightarrow} U \stackrel{\sim}{\leftarrow} \id_{\calD}$$
for which $v$ is the counit of an adjunction between $F$ and $G$.
\end{definition}

\begin{notation}
Let the subcategory $\cDelta_{S}^{\nounit} \subseteq \cDelta_{S}$ be defined as follows:
\begin{itemize}
\item Every object of $\cDelta_{S}$ belongs to $\cDelta_{S}^{\nounit}$.
\item A morphism $\alpha: ([m], c) \rightarrow ([n], c')$ of $\cDelta_{S}$ belongs to
$\cDelta_{S}^{\nounit}$ if and only if, whenever $\alpha(i) = \alpha(j)$ and
$c(i) = c(j) = \calC$, we have $i = j$.
\end{itemize}
\end{notation}

\begin{definition}\index{ZZZAdjDiagnounit@$\AdjDiag^{\nounit}(\calC, \calD)$}
Let $q: \End^{\otimes}(\calC, \calD) \rightarrow \Nerve( \cDelta_{S})^{op}$ denote the projection.\index{adjunction datum!nonunital}\index{nonunital!adjunction datum}
A {\it nonunital adjunction datum} is a functor $F \in \bHom_{ \Nerve( \cDelta_{S})^{op}}
( \Nerve( \cDelta^{\nounit}_{S})^{op}, \End^{\otimes}(\calC, \calD) )$ which carries
every $\calC$-convex morphism in $\cDelta_{S}^{\nounit}$ to a $q$-coCartesian
morphism in $\End^{\otimes}(\calC, \calD)$ (note that every $\calC$-convex morphism in
$\cDelta_{S}$ is a morphism of $\cDelta_{S}^{\nounit}$).

We let $\AdjDiag^{\nounit}(\calC, \calD)$ denote the full subcategory of
$\bHom_{ \Nerve( \cDelta_{S})^{op}}
( \Nerve( \cDelta^{\nounit}_{S})^{op}, \End^{\otimes}(\calC, \calD) )$ spanned by the nonunital adjunction diagrams.
\end{definition}

The inclusion $\cDelta^{\nounit}_{S} \subseteq \cDelta_{S}$ induces a restriction functor
$\AdjDiag(\calC, \calD) \rightarrow \AdjDiag^{\nounit}(\calC, \calD).$
Similarly, evaluation on the diagram $ [\calD] \leftarrow [\calD, \calD] \rightarrow [\calD, \calC, \calD] $ in $\cDelta^{\nounit}_{S}$ determines a functor
$r: \AdjDiag^{\nounit}(\calC, \calD) \rightarrow \NUAdj(\calC, \calD)$.
Lemma \ref{hungertown} implies that the composition
$ \AdjDiag(\calC, \calD) \rightarrow \AdjDiag^{\nounit}(\calC, \calD)
\rightarrow \NUAdj(\calC, \calD)$
factors through the full subcategory $\UAdj(\calC, \calD)$. Consequently, we obtain
a commutative diagram
$$ \xymatrix{ \AdjDiag(\calC, \calD) \ar[r] \ar[d] & \AdjDiag^{\nounit}(\calC, \calD) \ar[d] \\
\UAdj(\calC, \calD) \ar@{^{(}->}[r] & \NUAdj(\calC, \calD). }$$

We will break the proof of Theorem \ref{partA} into three steps:

\begin{proposition}\label{partA1}
Let $\calC$ and $\calD$ be $\infty$-categories. Then the projection map $\NUAdj(\calC, \calD) \rightarrow \Fun(\calC, \calD)$
induces a categorical equivalence
$\UAdj(\calC) \rightarrow \Fun'(\calC, \calD),$ where
$\Fun'(\calC, \calD)$ is the subcategory of $\Fun(\calC, \calD)$ whose objects are functors which admit right adjoints and whose morphisms are natural equivalences.
\end{proposition}

\begin{proposition}\label{partA2}
Let $\calC$ and $\calD$ be $\infty$-categories. The restriction map
$r: \AdjDiag^{\nounit}(\calC, \calD) \rightarrow \NUAdj(\calC, \calD)$ is a trivial Kan fibration.
\end{proposition}

\begin{proposition}\label{partA3}
Let $\calC$ and $\calD$ be $\infty$-categories. Then the restriction map
$\AdjDiag(\calC, \calD) \rightarrow \AdjDiag^{\nounit}(\calC, \calD)$
induces an equivalence of $\infty$-categories
$$ \AdjDiag( \calC, \calD) \rightarrow \AdjDiag^{\nounit}(\calC, \calD)
\times_{ \NUAdj(\calC, \calD)} \UAdj(\calC, \calD).$$
\end{proposition}

Let us grant these assertions for the moment.

\begin{proof}[Proof of Theorem \ref{partA}]
We observe that if $F \in \AdjDiag(\calC, \calD)$, the induced object of
$\NUAdj(\calC, \calD)$ belongs to $\UAdj(\calC, \calD)$ (Lemma \ref{hungertown}).
Consequently, the forgetful functor $\theta: \AdjDiag(\calC, \calD) \rightarrow \Fun(\calC, \calD)$ factors as a composition of categorical equivalences
$$ \AdjDiag(\calC, \calD) \stackrel{ \theta_0}{\rightarrow}
\AdjDiag^{\nounit}(\calC, \calD) \times_{ \NUAdj(\calC, \calD) }
\UAdj(\calC, \calD) \stackrel{\theta_1}{\rightarrow} \UAdj(\calC, \calD) \stackrel{\theta_2}{\rightarrow} \Fun'(\calC, \calD),$$ followed by the inclusion of $\Fun'(\calC, \calD)$ into $\Fun(\calC, \calD)$. 
It follows that $\theta$ is a categorical equivalence. Since $\theta$ is evidently a categorical fibration, it is a trivial Kan fibration.
\end{proof}

The proofs of Propositions \ref{partA1}, \ref{partA2}, and \ref{partA3} are very different from one another. Proposition \ref{partA1} is the easiest: it merely expresses the idea that
when a functor $F: \calC \rightarrow \calD$ admits a right adjoint $G$, then $G$ is uniquely determined. Proposition \ref{partA2} is more difficult. Roughly speaking, we need to show that
every nonunital adjunction $F_0$ can be promoted to a nonunital adjunction datum
$F: \Nerve( \cDelta^{\nounit}_{S})^{op} \rightarrow \End^{\otimes}(\calC, \calD)$. Our strategy is to choose an appropriate filtration of $\cDelta^{\nounit}_{S}$, and construct
$F$ from $F_0$ by forming a series of relative Kan extensions. 

Proposition \ref{partA3} involves rather different ideas. It asserts that a nonunital adjunction diagram
$F_0$ can be extended to an adjunction diagram $F$ if and only if the desired extension exists at the level of the homotopy category; moreover, in this case, $F$ is essentially unique. Our strategy is to reduce this assertion to an analogous result for nonunital algebras, and then invoke the 
results of \S \ref{giddug}.

We begin with the proof of Proposition \ref{partA1}. The first step is to reduce the size of 
the $\infty$-category $\NUAdj(\calC, \calD)$.

\begin{notation}
Let $\UAdj^{-}(\calC, \calD)$ denote the full subcategory of
$$ ( \Fun(\calC, \calD) \times \Fun(\calD, \calC) ) \times_{ \End(\calD) }
\End(\calD)^{ \id_{\calD} / }$$
spanned by those transformations $(F \circ G) \rightarrow \id_{\calD}$
which are counits for adjunctions between $F$ and $G$.
We will regard $\UAdj^{-}(\calC, \calD)$ as a simplicial subset of
$\UAdj(\calC, \calD)$, whose vertices correspond to diagrams
$F \circ G \rightarrow U \stackrel{\alpha}{\leftarrow} \id_{\calD}$
for which $\alpha$ is a degenerate morphism in $\End(\calD)$.
\end{notation}

\begin{lemma}\label{humblestar}
The inclusion $i: \UAdj^{-}(\calC, \calD) \subseteq \UAdj(\calC, \calD)$
is an equivalence of $\infty$-categories.
\end{lemma}

\begin{proof}
We have a homotopy Cartesian square
$$ \xymatrix{ \UAdj^{-}( \calC, \calD) \ar[r]^{i} \ar[d] & \UAdj(\calC, \calD) \ar[d] \\
\Delta^0 \ar[r] & \calE, }$$
where $\calE$ is the $\infty$-category of initial objects of $\End(\calD)^{\id_{\calD}/}$. Since
$\calE$ is a contractible Kan complex, we conclude that $i$ is a categorical equivalence.
\end{proof}

\begin{lemma}\label{stubby}
The $\infty$-category $\UAdj^{-}(\calC,\calD)$ is a Kan complex.
\end{lemma}

\begin{proof}
We note that a morphism $f$ in $\UAdj^{-}(\calC, \calD)$ can be identified with a pair of natural transformations $\alpha: F \rightarrow F'$, $\beta: G \rightarrow G'$, and a commutative square
$$ \xymatrix{ F \circ G \ar[r]^{v} \ar[d]^{\alpha \circ \beta} & \id_{\calD} \ar[d]^{=} \\
F' \circ G' \ar[r]^{v'} & \id_{\calD} }$$
in the $\infty$-category $\End(\calD)$.
We wish to prove that $\alpha$ and $\beta$ are equivalences. For this, we observe that
$\alpha$ induces a natural transformation of adjoint functors $\alpha': G' \rightarrow G$, and
$\beta$ induces a natural transformation of adjoint functors $\beta': F' \rightarrow F$. Using the commutativity of the above square, we deduce that $\alpha'$ is a homotopy inverse to $\beta$ and $\beta'$ is a homotopy inverse to $\alpha$.
\end{proof}

\begin{lemma}\label{stardust}
The projection map
$$ \theta: \NUAdj(\calC, \calD) \rightarrow \Fun(\calC, \calD)$$
induces a categorical equivalence $\theta_0: \UAdj^{-}(\calC) \rightarrow \Fun'(\calC, \calD)$, where
$\Fun'(\calC, \calD)$ is the subcategory of $\Fun(\calC, \calD)$ whose objects are functors which admit right adjoints and whose morphisms are natural equivalences.
\end{lemma}

\begin{proof}
It is clear that the essential image of $\theta_0$ consists of functors which admit right adjoints, and
Lemma \ref{stubby} implies that the image of every morphism in $\UAdj^{-}(\calC, \calD)$ is
an equivalence in $\Fun(\calC, \calD)$. Consequently, $\theta_0$ factors through
$\Fun'(\calC, \calD) \subseteq \Fun(\calC, \calD)$. It is easy to see that $\theta_0$ is a categorical fibration. Since the domain and codomain of $\theta_0$ are Kan complexes (Lemma \ref{stubby}), we deduce that $\theta_0$ is a Kan fibration (Propositions \toposref{groob}, \toposref{goey}, and Lemma \toposref{toothie2}). It will therefore suffice to show that every fiber of $\theta_0$ is contractible. We now observe that, if $F \in \Fun'(\calC, \calD)$, then the 
$\theta_0^{-1} \{ F \}$ can be identified with the $\infty$-category of initial objects
of the fiber product $\Fun( \calD, \calC) \times_{ \End(\calD) } \End(\calD)^{ \id_{\calD}/ }.$
Since this $\infty$-category is nonempty, it is a contractible Kan complex (Proposition \toposref{initunique}).
\end{proof}

\begin{proof}[Proof of Proposition \ref{partA1}]
Lemmas \ref{stardust} and \ref{humblestar} imply that the restriction map
$\UAdj(\calC, \calD) \rightarrow \Fun(\calC, \calD)$ factors as a composition
$\UAdj(\calC, \calD) \stackrel{r_0}{\rightarrow} \Fun'(\calC, \calD) \subseteq \Fun(\calC, \calD).$
Since $r_0$ is a categorical fibration, it will suffice to show that $r_0$ is a categorical equivalence. This follows immediately from Lemmas \ref{stardust} and \ref{humblestar}.
\end{proof}

We now turn to the proof of Proposition \ref{partA2}. In what follows, we let $q$ denote the projection map
$\End^{\otimes}(\calC,\calD) \rightarrow \Nerve(\cDelta_{S})^{op}.$

\begin{notation}\label{ablebear}
We define full subcategories
$ \calI_0 \subseteq \calI_1 \subseteq \calI_2 \subseteq \calI_3 \subseteq \cDelta_{S}^{\nounit}$ as follows:
\begin{itemize}
\item[$(J_0)$] An object $([n], c) \in \cDelta_{S}^{\nounit}$ belongs to
$\calI_0$ if and only if $\{ i \in [n]: c(i) = \calC \}$ is empty.
\item[$(J_1)$] An object $([n], c) \in \cDelta_{S}^{\nounit}$ belongs to
$\calI_1$ if and only if either $( [n], c) \in \calI_0$, or
$([n], c) = [ \calD, \calC, \calD ]$.
\item[$(J_2)$] An object $([n],c) \in \cDelta_{S}^{\nounit}$ belongs to 
$\calI_2$ if and only if $\{ i \in [n]: c(i) = \calC \}$ has at most one element.
\item[$(J_3)$] An object $([n], c) \in \cDelta_{S}^{\nounit}$ belongs to
$\calI_3$ if and only if $\{ i \in [n]: c(i) = \calC \}$ does not contain any pair of consecutive
integers $\{i, i+1\}$. 
\end{itemize}
For every full subcategory $\calI \subseteq \cDelta_{S}^{\nounit}$, we let
$\AdjDiag^{\nounit}_{\calI}(\calC, \calD) \subseteq
\bHom_{ \Nerve( \cDelta_{S})^{op} }( \Nerve(\calI)^{op}, \End^{\otimes}(\calC, \calD) )$
denote the full subcategory spanned by those functors $F$ which carry every $\calC$-convex morphism in
$\calI$ to a $q$-coCartesian morphism in $\End^{\otimes}(\calC, \calD)$.
\end{notation}

Proposition \ref{partA2} is a consequence of the following sequence of lemmas:

\begin{lemma}\label{partA20}
The restriction map $r_0: \AdjDiag^{\nounit}_{\calI_1}(\calC, \calD) \rightarrow \NUAdj(\calC, \calD)$
is a trivial Kan fibration.
\end{lemma}

\begin{lemma}\label{partA21}
Let $F_0 \in \AdjDiag^{\nounit}_{\calI_1}( \calC, \calD)$. Then:
\begin{itemize}
\item[$(1)$] There exists a functor $F \in \bHom_{ \Nerve(\cDelta_{S})^{op}}( \Nerve(\calI_2)^{op}, \End^{\otimes}(\calC, \calD) )$ which is a $q$-left Kan extension of $F_0$.
\item[$(2)$] Let $F \in \bHom_{ \Nerve(\cDelta_{S})^{op}}( \Nerve(\calI_2)^{op}, \End^{\otimes}(\calC, \calD) )$ be an arbitrary extension of $F_0$. Then $F$ is a $q$-left Kan extension of $F_0$ if and only if $F \in \AdjDiag^{\nounit}_{\calI_2}(\calC, \calD)$. 
\end{itemize}
\end{lemma}

\begin{lemma}\label{partA22}
Let $F_0 \in \AdjDiag^{\nounit}_{\calI_2}( \calC, \calD)$. Then:
\begin{itemize}
\item[$(1)$] There exists a functor $F \in \bHom_{ \Nerve(\cDelta_{S})^{op}}( \Nerve(\calI_3)^{op}, \End^{\otimes}(\calC, \calD) )$ which is a $q$-right Kan extension of $F_0$.
\item[$(2)$] Let $F \in \bHom_{ \Nerve(\cDelta_{S})^{op}}( \Nerve(\calI_3)^{op}, \End^{\otimes}(\calC, \calD) )$ be an arbitrary extension of $F_0$. Then $F$ is a $q$-right Kan extension of $F_0$ if and only if $F \in \AdjDiag^{\nounit}_{\calI_3}(\calC, \calD)$. 
\end{itemize}
\end{lemma}

\begin{lemma}\label{partA23}
Let $F_0 \in \AdjDiag^{\nounit}_{\calI_3}( \calC, \calD)$. Then:
\begin{itemize}
\item[$(1)$] There exists a functor $F \in \bHom_{ \Nerve(\cDelta_{S})^{op}}( \Nerve(\cDelta_{S}^{\nounit})^{op}, \End^{\otimes}(\calC, \calD) )$ which is a $q$-left Kan extension of $F_0$.
\item[$(2)$] Let $F \in \bHom_{ \Nerve(\cDelta_{S})^{op}}( \Nerve(\cDelta_{S}^{\nounit})^{op}, \End^{\otimes}(\calC, \calD) )$ be an arbitrary extension of $F_0$. Then $F$ is a $q$-left Kan extension of $F_0$ if and only if $F \in \AdjDiag^{\nounit}(\calC, \calD)$. 
\end{itemize}
\end{lemma}

\begin{proof}[Proof of Proposition \ref{partA2}]
The restriction map $r: \AdjDiag^{\nounit}(\calC, \calD) \rightarrow \NUAdj(\calC, \calD)$
factors as a composition
\begin{eqnarray*}
\AdjDiag^{\nounit}(\calC, \calD) & \stackrel{r_3}{\rightarrow} & \AdjDiag^{\nounit}_{\calI_3}(\calC, \calD) \\
& \stackrel{r_2}{\rightarrow} & \AdjDiag^{\nounit}_{\calI_2}(\calC, \calD) \\
& \stackrel{r_1}{\rightarrow} & \AdjDiag^{\nounit}_{\calI_1}(\calC, \calD) \\
& \stackrel{r_0}{\rightarrow} & \NUAdj(\calC, \calD).
\end{eqnarray*}
Lemma \ref{partA20} implies that $r_0$ is a trivial Kan fibration. Combining Lemmas
\ref{partA21}, \ref{partA22}, and \ref{partA23} with Proposition \toposref{lklk}, we deduce
that $r_1$, $r_2$, and $r_3$ are trivial Kan fibrations. It follows that 
$r = r_0 \circ r_1 \circ r_2 \circ r_3$ is also a trivial Kan fibration.
\end{proof}

\begin{proof}[Proof of Lemma \ref{partA20}]
We have a coCartesian fibration $p: \Nerve(\calI_1)^{op} \rightarrow \Delta^1$, where
the inverse image of $\{1\} \subseteq \Delta^1$ is $\Nerve(\calI_0)^{op}$. Let
$K \subseteq \Nerve(\calI_1)^{op}$ be the union of $\Nerve(\calI_0)$ and
the $p$-coCartesian edge joining $[ \calD, \calC, \calD]$ to $[\calD, \calD]$. 
Using Proposition \toposref{simplexplay}, we deduce that the inclusion
$K \subseteq \Nerve(\calI_1)^{op}$ is a categorical equivalence. Consequently, 
the restriction map
$$ \bHom_{ \Nerve(\cDelta_{S})^{op} }( \Nerve(\calI_1)^{op}, \End^{\otimes}(\calC, \calD) )
\rightarrow \bHom_{ \Nerve(\cDelta_{S})^{op}}( K, \End^{\otimes}(\calC, \calD) )$$
is a trivial Kan fibration. We also observe that a functor $F \in \bHom_{ \Nerve(\cDelta_{S})^{op} }( \Nerve(\calI_1)^{op}, \End^{\otimes}(\calC, \calD) )$ belongs to
$\AdjDiag^{\nounit}_{\calI_1}(\calC, \calD)$ if and only if $F | \Nerve(\calI_0)^{op}$ belongs to
$\AdjDiag^{\nounit}_{\calI_0}(\calC, \calD)$. It will therefore suffice to show that the restriction map
$$r'_0: \bHom_{ \Nerve(\cDelta_{S})^{op}}( K, \End^{\otimes}(\calC, \calD) )
\times_{ \bHom_{ \Nerve(\cDelta_{S})^{op} }( \Nerve(\calI_0)^{op}, \End^{\otimes}(\calC, \calD) ) } \AdjDiag^{\nounit}_{\calI_0}(\calC, \calD) \rightarrow \NUAdj(\calC, \calD)$$
is a trivial Kan fibration. We observe that $r'_0$ is a pullback of the restriction map
$r''_0: \AdjDiag^{\nounit}_{\calI_0}(\calC, \calD) \rightarrow \calE,$
where $\calE$ is the $\infty$-category of initial objects in $\End(\calD)^{ \id_{\calD}/ }$ and
$r''_0$ is given by evaluation on the morphism $[\calD] \rightarrow [\calD, \calD]$ in
$\Nerve(\cDelta_{S}^{\nounit})^{op}$. 

We wish to show that $r''_0$ is a trivial Kan fibration. It is easy to see that $r''_0$ is a categorical fibration. It will therefore suffice to show that $r''_0$ is a categorical equivalence. For this, it suffices to show that the source and target of $r''_0$ are contractible Kan complexes. We observe that
Proposition \ref{gurgle} allows us to identify $\AdjDiag^{\nounit}_{\calI_0}(\calC, \calD)$ with the
$\infty$-category of initial objects in $\Alg( \End(\calD) )$. The desired contractibility assertions now follow immediately from Proposition \toposref{initunique}.
\end{proof}

\begin{proof}[Proof of Lemma \ref{partA21}]
We observe that for every object in $( [n],c) \in \calI_2$, the category
$\calI_{1} \times_{ \calI_2 } (\calI_{2})_{([n],c)/}$ has an initial object. If
$( [n], c) \in \calI_0$, this initial object is given by the identity map from
$([n],c)$ to itself; otherwise, it is given by the unique map
$( [n], c) \rightarrow [\calD, \calC, \calD]$. Since $q$ is a coCartesian fibration, Lemma \toposref{kan2} immediately implies the following analogues of $(1)$ and $(2)$:
\begin{itemize}
\item[$(1')$] For every functor $F_0 \in \bHom_{ \Nerve(\cDelta_{S})^{op}}( \Nerve(\calI_1)^{op}, \End^{\otimes}(\calC, \calD) ),$ there exists a functor $$F \in \bHom_{ \Nerve(\cDelta_{S})^{op}}( \Nerve(\calI_2)^{op}, \End^{\otimes}(\calC, \calD) )$$ which is a $q$-left Kan extension of $F_0$.
\item[$(2')$] Let $F \in \bHom_{ \Nerve(\cDelta_{S})^{op}}( \Nerve(\calI_2)^{op}, \End^{\otimes}(\calC, \calD) )$. Then $F$ is a $q$-left Kan extension of $F_0 = F | \Nerve(\calI_1)^{op}$
if and only if it satisfies the following condition:
\begin{itemize}
\item[$(\ast)$] For every object $( [n], c) \in \calI_2$ which does not belong to $\calI_0$, $F$ carries the canonical map $\alpha: ([n],c) \rightarrow [\calD, \calC, \calD]$ in $\calI_2$ to a $q$-coCartesian morphism of $\End^{\otimes}(\calC, \calD)$.
\end{itemize}
\end{itemize}
It is clear that $(1')$ implies $(1)$. To prove that $(2')$ implies $(2)$, we must show that
a functor $$F \in \bHom_{ \Nerve(\cDelta_{S})^{op}}( \Nerve(\calI_2)^{op}, \End^{\otimes}(\calC, \calD) )$$ belongs to $\AdjDiag^{\nounit}_{\calI_2}(\calC, \calD)$ if and only if
$F$ satisfies $(\ast)$ and $F | \Nerve(\calI_1)$ belongs to $\AdjDiag^{\nounit}_{\calI_1}( \calC, \calD)$. The ``only if'' direction is clear, since each of the morphisms $\alpha$ which appears in the statement of $(\ast)$ is $\calC$-convex.

Conversely, suppose that $F_0 = F | \Nerve(\calI_1)$ belongs to $\AdjDiag^{\nounit}_{\calI_1}(\calC, \calD)$ and that $F$ satisfies $(\ast)$. Let $\alpha: ([m], c) \rightarrow ([n], c')$ be a $\calC$-convex morphism in $\calI_2$. We must show that $F(\alpha)$ is $q$-coCartesian. If $([n],c') \in \calI_0$, this follows from our assumption that $F_0$ satisfies $(\ast)$. If $([n], c') \notin \calI_0$, then we have a commutative diagram
$$ \xymatrix{ & ( [n], c') \ar[dr]^{\beta} & \\
([m], c) \ar[ur]^{\alpha} \ar[rr]^{\gamma} & & [ \calD, \calC, \calD ]. }$$
According to Proposition \toposref{protohermes}, it will suffice to show that
$F(\beta)$ and $F(\gamma)$ are $q$-coCartesian. For $\beta$, this follows from
condition $(\ast)$. If
$([m],c) \notin \calI_0$, then condition $(\ast)$ also guarantees that $F(\gamma)$ is $q$-coCartesian. If $([m], c) \in \calI_0$, then $\gamma$ factors as a composition
$ ( [m],c) \stackrel{\gamma'}{\rightarrow} [ \calD] \stackrel{\gamma''}{\rightarrow} [ \calD, \calC, \calD]. $
Then $F(\gamma')$ is $q$-coCartesian (since $\gamma'$ is a $\calC$-convex morphism in
$\calI_0$) and $F(\gamma'')$ is $q$-coCartesian (since the $\infty$-category $\End^{\otimes}(\calC, \calD)_{[\calD]}$ is a Kan complex). Applying Proposition \toposref{protohermes}, we deduce that $F(\gamma)$ is $q$-coCartesian, as desired.
\end{proof}

\begin{proof}[Proof of Lemma \ref{partA22}]
Fix an object $([n], c) \in \calI_3$, and let $\calJ = \calI_2 \times_{ \calI_3} (\calI_3)_{/([n],c)}$. Our first goal is to prove:
\begin{itemize}
\item[$(a)$] Consider the commutative diagram
$$ \xymatrix{ \Nerve(\calJ^{op}) \ar[r]^{f} \ar@{^{(}->}[d] & \End^{\otimes}(\calC, \calD) \ar[d]^{q} \\
\Nerve(\calJ^{op})^{\triangleleft} \ar[r] \ar@{-->}[ur]^{\overline{f}} & \Nerve(\cDelta_{S})^{op} }$$
where $f$ denotes the composition 
$$\Nerve(\calJ^{op}) \rightarrow \Nerve(\calI_2^{op}) \stackrel{F_0}{\rightarrow} \End^{\otimes}(\calC, \calD).$$ Then there exists a map $\overline{f}$, as indicated, which is a $q$-limit diagram.
\item[$(b)$] An arbitrary extension $\overline{f}$ rendering the above diagram commutative is a
$q$-limit diagram if and only if the following condition is satisfied:
\begin{itemize}
\item[$(\star)$] For every {\em convex} morphism $\alpha: ([m], c') \rightarrow ([n],c)$ in
$\calJ$, the image of $\{ \alpha\}^{\triangleleft} \subseteq \Nerve(\calJ^{op})^{\triangleleft}$ under
$\overline{f}$ is a $q$-coCartesian morphism of $\End^{\otimes}(\calC, \calD)$.
\end{itemize}
\end{itemize}

Let $\calJ_0$ be the full subcategory of $\calJ$ spanned by those objects
$\alpha: ([m],c') \rightarrow ([n],c)$ such that the induced map $[m] \rightarrow [n]$ is {\em injective}. 
We can identify the category $\calJ_0$ with the partially ordered sets of nonempty subsets $S \subseteq [n]$ having the property that $\{ j \in S: c(j) = \calC \}$ has at most one element.
Consider the following further properties which $S$ may or may not satisfy:

\begin{itemize}
\item[$(P_1)$] The set $S$ contains $\{ j \in [n]: c(j) = \calD \}$. 
\item[$(P_2)$] For some $0 < i < n$ such that $c(i) = \calC$, we have $S = \{i-1, i+1\}$. 
\item[$(P_3)$] For some $0 < i < n$ such that $c(i) = \calC$, we have $S = \{ i-1, i, i+1 \}$.
\item[$(P_4)$] The set $S$ consists of a single element.
\item[$(P_5)$] For some $0 \leq i < n$, we have $S = \{ i, i+1 \}$.
\end{itemize}

For $1 \leq i \leq 5$, let $\calJ_i$ denote the full subcategory of $\calJ_0$ spanned by those objects 
which satisfy $(P_j)$ for {\em some} $j \geq i$. Consider the following assertions:

\begin{itemize}
\item[$(a_i)$] Consider the commutative diagram
$$ \xymatrix{ \Nerve(\calJ_i^{op}) \ar[r]^{f_i} \ar@{^{(}->}[d] & \End^{\otimes}(\calC, \calD) \ar[d]^{q} \\
\Nerve(\calJ^{op}_i)^{\triangleleft} \ar[r] \ar@{-->}[ur]^{\overline{f}_i} & \Nerve(\cDelta_{S})^{op} }$$
where $f_i$ denotes the composition 
$\Nerve(\calJ_i^{op}) \subseteq \Nerve(\calJ^{op}) \rightarrow \Nerve(\calI_2^{op}) \stackrel{F_0}{\rightarrow} \End^{\otimes}(\calC, \calD).$ Then there exists a map $\overline{f}_i$, as indicated, which is a $q$-limit diagram.
\item[$(b_i)$] An arbitrary extension $\overline{f}_i$ rendering the above diagram commutative is a
$q$-limit diagram if and only if the following condition is satisfied:
\begin{itemize}
\item[$(\star'_i)$] For every morphism $\alpha: ([m], c') \rightarrow ([n],c)$ in
$\calJ_0$ satisfying $P_5$, the image of $\{ \alpha\}^{\triangleleft} \subseteq \Nerve(\calJ_0^{op})^{\triangleleft}$ under $\overline{f}_i$ is a $q$-coCartesian morphism of $\End^{\otimes}(\calC, \calD)$.
\end{itemize}
\end{itemize}

We note that, since $F_0 \in \AdjDiag^{\nounit}_{\calI_2}$, if $\overline{f}$ is as in $(b)$ and
$0 \leq i \leq 5$, then $\overline{f}$ satisfies $(\star)$ if and only if $\overline{f} | \Nerve(\calJ_i^{op})$ satisfies $(\star'_5)$. 

The inclusions $\calJ_1 \subseteq \calJ$,
$\calJ_2 \subseteq \calJ_3$, and $\calJ_4 \subseteq \calJ_5$ admit left adjoints, so the induced inclusions $\Nerve(\calJ_1) \subseteq \calJ$, $\Nerve(\calJ_2) \subseteq \Nerve(\calJ_3)$, and
$\Nerve(\calJ_4) \subseteq \Nerve(\calJ_5)$ are cofinal. It follows that
$(a) \Leftrightarrow (a_1)$, $(a_2) \Leftrightarrow (a_3)$, and $(a_4) \Leftrightarrow (a_5)$.
Using the assumption that $F_0 \in \AdjDiag^{\nounit}_{\calI_2}$, we deduce that
$f_1$ is a $q$-right Kan extension of $f_2$ and that $f_3$ is a $q$-right Kan extension of $f_4$. 
Applying Lemma \toposref{kan0}, we deduce that $(a_1) \Leftrightarrow (a_2)$, and
$(a_3) \Leftrightarrow (a_4)$. Composing these equivalences, we see that $(a)$ is equivalent to
$(a_5)$. The same argument shows that $(b)$ is equivalent to $(b_5)$. We now observe that
$(a_5)$ and $(b_5)$ are obvious, since $q$ is a coCartesian fibration which exhibits
$\End^{\otimes}(\calC, \calD)_{ ([n],c)}$ as equivalent to the product
$$\prod_{0 \leq i < n} \End^{\otimes}(\calC, \calD)_{ ( \{i, i+1\}, c | \{i, i+1\} )}.$$

Assertion $(1)$ now follows from $(a)$ and Lemma \toposref{kan0}. Moreover, 
$(b)$ implies the following analogue of $(2)$:
\begin{itemize}
\item[$(2')$] Let $F \in \bHom_{ \Nerve(\cDelta_{S})^{op}}( \Nerve(\calI_3)^{op}, \End^{\otimes}(\calC, \calD) )$ be an arbitrary extension of $F_0$. Then $F$ is a $q$-right Kan extension of $F_0$ if and only if $F$ satisfies the following condition:
\begin{itemize}
\item[$(\ast)$] For every {\em convex} morphism $\alpha: ([m], c) \rightarrow ( [n], c')$ in
$\calI_3$, where $( [m], c) \in \calI_2$, the image $F(\alpha)$ is a $q$-coCartesian
morphism in $\End^{\otimes}(\calC, \calD)$.
\end{itemize}
\end{itemize}
To complete the proof, it will suffice to show that if $F \in \bHom_{ \Nerve(\cDelta_{S})^{op}}( \Nerve(\calI_3)^{op}, \End^{\otimes}(\calC, \calD) )$ is an extension of
$F_0$, then $F$ satisfies $(\ast)$ if and only if $F$ carries every $\calC$-convex morphism
in $\calI_3$ to a $q$-coCartesian morphism in $\End^{\otimes}(\calC, \calD)$. The ``if'' direction is obvious, since every convex morphism in $\calI_3$ is $\calC$-convex. For the converse,
let us suppose that $F$ satisfies $(\ast)$ and that
$\alpha: ([m], c) \rightarrow ([n], c')$ is a $\calC$-convex morphism in
$\calI_3$. We wish to prove that $F(\alpha)$ is $q$-coCartesian. Using the product decomposition of $\End^{\otimes}(\calC, \calD)_{ ([m],c)}$, we can reduce to the case where $m = 1$, so that
$( [m], c) \in \calI_2$. We now observe that $\alpha$ admits a unique factorization as a composition
$ ([m], c) \stackrel{\alpha'}{\rightarrow} ( [n_0], c'_0) \stackrel{\alpha''}{\rightarrow}
( [n], c),$
where $\alpha'$ preserves endpoints and $\alpha''$ is convex. Assumption
$(\ast)$ guarantees that $F(\alpha'')$ is $q$-coCartesian. Moreover, using the assumption that
$\alpha$ is $\calC$-convex, we deduce that $([n_0], c'_0) \in \calI_2$ and
the morphism $\alpha'$ is $\calC$-convex. Since $F_0 \in \AdjDiag^{\nounit}_{\calI_2}( \calC, \calD)$, we conclude that $F(\alpha') = F_0(\alpha')$ is $q$-coCartesian. Applying Proposition
\toposref{protohermes}, we deduce that $F(\alpha)$ is $q$-coCartesian, as desired.
\end{proof}

The proof of Lemma \ref{partA23} will require the following preliminary:

\begin{lemma}\label{hornhold}
Let $p: X \rightarrow S$ be a coCartesian fibration of $\infty$-categories, let $K$ be a weakly contractible simplicial set, and suppose given a diagram
$$ \xymatrix{ K \ar@{^{(}->}[d] \ar[r]^{f_0} & X \ar[d]^{p} \\
K^{\triangleright} \ar[r]^{g} \ar@{-->}[ur]^{f} & S. }$$
Suppose that $f_0$ carries each edge of $K$ to a $p$-coCartesian edge of $X$. Then:
\begin{itemize}
\item[$(1)$] There exists a map $f$ as indicated in the diagram, which is a $p$-colimit diagram.
\item[$(2)$] Let $f$ be an arbitrary morphism which renders the above diagram commutative. The following conditions are equivalent:
\begin{itemize}
\item[$(i)$] The map $f$ is a $p$-colimit diagram.
\item[$(ii)$] For every vertex $k$ of $K$, the image of $\{k\}^{\triangleright}$ under $f$ is
a $p$-coCartesian morphism of $X$.
\item[$(iii)$] There exists a vertex $k$ of $K$ such that the image of $\{ k \}^{\triangleright}$ under $f$ is a $p$-coCartesian morphism of $X$.
\end{itemize}
\end{itemize}
\end{lemma}

\begin{proof}
Using Proposition \toposref{chocolatelast}, we can reduce to the case where
the map $g$ is constant at some object $s$ of $S$, so that $f_0$ factors through the fiber
$X_{s}$. Corollary \toposref{silt} guarantees the existence of an extension $f$ satisfying
$(ii)$. Moreover, Corollary \toposref{silt} also ensures that for every edge
$s \rightarrow s'$ in $S$, the image of $f$ under the associated functor $X_{s} \rightarrow X_{s'}$ is a colimit diagram in $X_{s'}$. Proposition \toposref{relcolfibtest} implies that $f$ is a $p$-colimit diagram. This proves $(1)$, as well as the implication $(ii) \rightarrow (i)$ of $(2)$. The implication
$(i) \rightarrow (ii)$ follows from the uniqueness of relative colimit diagrams, and the equivalence $(ii) \Leftrightarrow (iii)$ follows from the connectedness of $K$.
\end{proof}

\begin{proof}[Proof of Lemma \ref{partA23}]
Fix an object $([m], c) \in \cDelta_{S}^{\nounit}$, and let
$\calJ = \calI_3 \times_{ \cDelta_{S} } ( \cDelta_{S} )_{ ([m],c) }$.
Let $\calJ_0$ denote the full subcategory of $\calJ$ spanned by those morphisms
$\alpha: ([m],c) \rightarrow ([n], c')$ with the following property:
\begin{itemize}
\item[$(\star)$] The map $\alpha$ preserves endpoints, and for $0 \leq i < m$, either $\alpha(i+1) = \alpha(i) + 1$, or
$c(i) = c(i+1) = \calC$, $\alpha(i+1) > \alpha(i) + 1$, and 
$c'(j) = \calD$ for $\alpha(i) < j < \alpha(i+1)$. 
\end{itemize}
The inclusion $\calJ_0 \subseteq \calJ$ admits a right adjoint, so that
the inclusion $\Nerve(\calJ_0)^{op} \subseteq \Nerve(\calJ)^{op}$ is cofinal. 
We observe that $\calJ_0$ is equivalent to a product of copies of the category
$\cDelta$, where the product is indexed by integers $i$ such that
$0 \leq i < i+1 \leq m$ and $c(i) = c(i+1) = \calC$. In particular, the simplicial set
$\Nerve(\calJ_0)$ is weakly contractible.

We observe that, since $F_0 \in \AdjDiag^{\nounit}_{\calI_3}(\calC, \calD)$, the composition
$$\Nerve(\calJ_0)^{op} \subseteq \Nerve(\calJ)^{op} \rightarrow
\Nerve(\calI_3)^{op} \stackrel{F_0}{\rightarrow} \End^{\otimes}(\calC, \calD)^{op}$$ carries each morphism in $\calJ_0$ to $q$-coCartesian edge of $\End^{\otimes}(\calC, \calD)$.
Combining Lemma \ref{hornhold} with Lemma \toposref{kan2}, we deduce $(1)$ and the following analogue of $(2)$:
\begin{itemize}
\item[$(2')$] Let $F \in \bHom_{ \Nerve(\cDelta_{S})^{op}}( \Nerve(\cDelta_{S}^{\nounit})^{op}, \End^{\otimes}(\calC, \calD) )$ be an arbitrary extension of $F_0$. Then $F$ is a $q$-left Kan extension of $F_0$ if and only if $F$ satisfies the following condition:
\begin{itemize}
\item[$(\ast)$] For every morphism $\alpha: ([m],c) \rightarrow
([n],c')$ in $\cDelta_{S}^{\nounit}$ which satisfies $\star$, the image
$F(\alpha)$ is $q$-coCartesian.
\end{itemize}
\end{itemize}
To complete the proof, it will suffice to show that if
$F \in \bHom_{ \Nerve(\cDelta_{S})^{op}}( \Nerve(\cDelta_{S}^{\nounit})^{op}, \End^{\otimes}(\calC, \calD) )$ is an extension of $F_0$, then $F$ satisfies $(\ast)$ if and only if $F$
carries every $\calC$-convex morphism of $\cDelta_{S}^{\nounit}$ to a $q$-coCartesian morphism in $\End^{\otimes}(\calC, \calD)$. The ``if'' direction is obvious, since every morphism
which satisfies $(\star)$ is $\calC$-convex.

For the converse, let us suppose that $F$ satisfies $(\ast)$ and that
$\alpha: ([m], c) \rightarrow ([n],c')$ is a $\calC$-convex morphism of $\cDelta_{S}^{\nounit}$. We wish to prove that $F(\alpha)$ is $q$-coCartesian. First, choose a morphism
$\beta: ( [n], c') \rightarrow ( [p], c'')$ satisfying $(\star)$. Condition $(\ast)$ guarantees that
$F(\beta)$ is $q$-coCartesian. In view of Proposition \toposref{protohermes}, it will suffice to show that $F( \beta \circ \alpha)$ is $q$-coCartesian. In other words, we may replace
$([n],c')$ by $([p],c'')$ and thereby reduce to the case where $([n], c') \in \calI_3$. 
We now observe that $\alpha$ factors as a composition
$$ ( [m], c) \stackrel{\alpha'}{\rightarrow} ( [n_0], c'_0 ) \stackrel{\alpha''}{\rightarrow} ( [n], c),$$
where $\alpha'$ satisfies $(\star)$ and $\alpha''$ is $\calC$-convex.
Condition $(\ast)$ implies that
$F(\alpha')$ is $q$-coCartesian, and our assumption that $F_0 \in \AdjDiag^{\nounit}_{\calI_3}(\calC, \calD)$ guarantees that $F(\alpha'')$ is $q$-coCartesian. Invoking Proposition \toposref{protohermes}, we deduce that $F(\alpha)$ is $q$-coCartesian, as desired.
\end{proof}

We now give the proof of Proposition \ref{partA3}. We begin with a rough outline.
Every nonunital adjunction datum $F_0$ determines a {\em nonunital monad} on
$\calC$: that is, a nonunital algebra $T_0 \in \Alg^{\nounit}( \End(\calC) )$. If $F_0$ extends to an adjunction, then $T$ extends to an algebra object in $\End(\calC)$. We would like to prove a converse to this result: that is, the only obstruction to the existence of a adjunction datum extending
$F_0$ is the existence of a monad extending $T_0$. Actually, this is not quite true: we will need to know not only that $T$ exists, but that various (nonunital) $T_0$-modules can be extended to
(unital) $T$-modules. Nevertheless, the arguments below will allow us to reduce to the problems regarding the existence and uniqueness of units which were treated in \S \ref{digunit} and \S \ref{digtwo}.

We begin by introducing a bit of notation.

\begin{notation}
Let $\calJJ$ denote the full subcategory of $\Fun( [1], \cDelta_{S} )$ spanned by those
morphisms $\alpha: ( [m], c_0) \rightarrow ([n], c)$ which induce a {\em bijection}
$ \{ i \in [m]: c_0(i) = \calD \} \rightarrow \{ j \in [n]: c(j) = \calD \}.$
We define a subcategory $\calJJ^{\nounit} \subseteq \calJJ$ as follows:
\begin{itemize}
\item Every object of $\calJJ$ belongs to $\calJJ^{\nounit}$.
\item Given a pair of objects $\alpha, \alpha' \in \calJJ$, a morphism
$$ \xymatrix{ ( [m], c_0) \ar[r]^{\alpha} \ar[d]^{\beta} & ( [n], c) \ar[d]^{\gamma} \\
( [m'], c'_0) \ar[r]^{\alpha'} & ( [n'], c') }$$
in $\calJJ$ belongs to $\calJJ^{\nounit}$ if and only if the morphism $\beta$ induces an {\em injection} $\{ i \in [m_0]: c'_0(i) = \calC \} \rightarrow
\{ j \in [m_1], c'_1(j) = \calC \}$.
\end{itemize}

We will identify $\cDelta_{S}$ with the full subcategory of
$\calJJ$ spanned by constant functors $[1] \rightarrow \cDelta_{S}$. 
We observe that the intersection $\calJJ^{\nounit} \cap \cDelta_{S}$ can be identified with the full subcategory $\cDelta^{\nounit}_{S} \subseteq \cDelta_{S}$. 
The inclusions $\cDelta_{S} \subseteq \calJJ$ and
$\cDelta^{\nounit}_{S} \subseteq \calJJ^{\nounit}$
admit right adjoints $\psi$ and $\psi^{\nounit}$, both given by the formula
$$ ( \alpha: ([m],c_0) \rightarrow ( [n], c) ) \mapsto \id_{ ( [m], c_0 )}.$$

We let $\widetilde{\AdjDiag}(\calC, \calD)$ denote the full subcategory of
$ \bHom_{ \Nerve(\cDelta_{S})^{op} }(
\Nerve(\calJJ)^{op}, \End^{\otimes}(\calC, \calD) )$
spanned by those functors $F: \Nerve(\calJJ)^{op} \rightarrow \End^{\otimes}(\calC, \calD)$
such that $q \circ F = \Nerve(\psi)^{op}$, and whenever $\gamma: \alpha \rightarrow \alpha'$ is a morphism
in $\calJJ$ such that $\psi(\gamma): \psi(\alpha) \rightarrow \psi(\alpha')$ is $\calC$-convex, the morphism $F(\gamma)$ is a $q$-coCartesian. Similarly, we let
$\widetilde{\AdjDiag}^{\nounit}(\calC, \calD)$ denote the full subcategory of
$ \bHom_{ \Nerve(\cDelta_{S})^{op} }(
\Nerve(\calJJ^{\nounit})^{op}, \End^{\otimes}(\calC, \calD) )$
spanned by those functors $F_0: \Nerve(\calJJ^{\nounit})^{op} \rightarrow \End^{\otimes}(\calC, \calD)$ such that $q \circ F_0 = \Nerve(\psi^{\nounit})^{op}$, and whenever $\gamma: \alpha \rightarrow \alpha'$ is a morphism
in $\calJJ^{\nounit}$ such that $\psi(\gamma): \psi(\alpha) \rightarrow \psi(\alpha')$ is $\calC$-convex, the morphism $F_0(\gamma)$ is a $q$-coCartesian.
\end{notation}

\begin{lemma}\label{hyperstork}
The restriction maps
$$ \widetilde{\AdjDiag}(\calC, \calD) \rightarrow \AdjDiag(\calC, \calD) \quad \widetilde{\AdjDiag}^{\nounit}(\calC, \calD) \rightarrow \AdjDiag^{\nounit}(\calC, \calD)$$
are trivial Kan fibrations.
\end{lemma}

\begin{proof}
We will give the proof in the unital case; the nonunital case uses exactly the same argument.
According to Proposition \toposref{lklk}, it will suffice to show:
\begin{itemize}
\item[$(a)$] A functor 
$F \in \bHom_{ \Nerve(\cDelta_{S})^{op} }(
\Nerve(\calJJ)^{op}, \End^{\otimes}(\calC, \calD) )$ belongs to
$\widetilde{\AdjDiag}(\calC, \calD)$ if and only if
$F_0 = F | \Nerve( \cDelta_{S})^{op}$ belongs to 
$\AdjDiag(\calC, \calD)$, and $F$ is a $q$-right Kan extension of $F_0$.
\item[$(b)$] Every functor $F_0 \in \AdjDiag(\calC, \calD)$ admits a $q$-right Kan extension
$$F \in \bHom_{ \Nerve(\cDelta_{S})^{op} }(
\Nerve(\calJJ)^{op}, \End^{\otimes}(\calC, \calD) ).$$
\end{itemize}

We observe that a functor $F \in \bHom_{ \Nerve(\cDelta_{S})^{op} }(
\Nerve(\calJJ)^{op}, \End^{\otimes}(\calC, \calD) )$ is a $q$-right Kan extension of
$F_0 = F | \Nerve( \cDelta_{S})^{op}$ at an object $\alpha \in \calJ$ if and only if
$F$ carries the counit map $v_{\alpha}: \psi(\alpha) \rightarrow \alpha$ to an equivalence
in $\End^{\otimes}(\calC, \calD)$. Assertion $(b)$ now follows from Lemma \toposref{kan2}, and the
``only if'' direction of $(a)$ is obvious from the definitions. For the converse, let us suppose that
$F$ carries each $v_{\alpha}$ to an equivalence in $\End^{\otimes}(\calC, \calD)$, and that
$F_0 \in \AdjDiag(\calC, \calD)$. We wish to prove that $F \in \widetilde{\AdjDiag}(\calC, \calD)$. 
Let $\gamma: \alpha \rightarrow \alpha'$ be a morphism in $\calJJ$ such that $\psi(\gamma)$ is
$\calC$-convex. We wish to prove that $F(\gamma)$ is $q$-coCartesian. Since $F$ is a $q$-right Kan extension of $F_0$, the morphisms $F(\gamma)$ and $F( \psi(\gamma) )$ are equivalent. Moreover, $F( \psi(\gamma) ) = F_0( \psi(\gamma) ) $ is $q$-coCartesian in virtue of our assumption that $F_0 \in \AdjDiag(\calC, \calD)$.
\end{proof}

The advantage of the $\infty$-categories $\widetilde{\AdjDiag}(\calC, \calD)$ and
$\widetilde{\AdjDiag}^{\nounit}(\calC, \calD)$ is that they are assembled out of ``local'' information. 
To be more precise, we need to introduce a bit more notation.

\begin{notation}
Let $\cDelta'_{S}$ denote an isomorphic copy of the category $\cDelta_{S}$ (but with a new name, to avoid confusion). Let $\pi: \calJJ \rightarrow \cDelta'_{S}$ be defined by the formula
$$(\alpha: ([m], c) \rightarrow ([n], c')) \mapsto ( [n], c').$$
We define a pair of simplicial sets $\calM$, $\calM^{\nounit}$ equipped with maps
$\calM \rightarrow \calM^{\nounit} \rightarrow \Nerve( \cDelta'_{S})^{op}.$
These simplicial sets are characterized by the property that, for every simplicial set $K$
equipped with a map $K \rightarrow \Nerve( \cDelta'_{S})^{op}$, there are canonical bijections
$$ \Hom_{ \Nerve(\cDelta'_{S})^{op} }( K, \calM )
\simeq \Hom_{ \Nerve(\cDelta_{S})^{op} }( K \times_{ \Nerve(\cDelta'_{S})^{op}}
\Nerve(\calJJ)^{op}, \End^{\otimes}(\calC, \calD) )$$
$$ \Hom_{ \Nerve(\cDelta'_{S})^{op} }( K, \calM^{\nounit} )
\simeq \Hom_{ \Nerve(\cDelta_{S})^{op} }( K \times_{ \Nerve(\cDelta'_{S})^{op}}
\Nerve(\calJJ^{\nounit})^{op}, \End^{\otimes}(\calC, \calD) ).$$
\end{notation}

It follows immediately from the definitions that we can identify 
$\widetilde{\AdjDiag}^{\nounit}(\calC, \calD)$ with a full subcategory of the
$\infty$-category $\bHom_{ \Nerve(\cDelta'_{S})^{op} }
( \Nerve(\cDelta^{\nounit}_{S})^{op}, \calM)$ of sections of $p^{\nounit}$. Similarly,
we can identify $\widetilde{\AdjDiag}(\calC, \calD)$ with the full subcategory
$$\bHom_{ \Nerve(\cDelta'_{S})^{op} }
( \Nerve(\cDelta'_{S})^{op}, \calM)
\times_{ \bHom_{ \Nerve(\cDelta'_{S})^{op} }
( \Nerve(\cDelta^{\nounit}_{S})^{op}, \calM) } \widetilde{\AdjDiag}^{\nounit}(\calC, \calD)$$ of $\bHom_{ \Nerve(\cDelta'_{S})^{op} }
( \Nerve(\cDelta'_{S})^{op}, \calM)$.

Proposition \ref{partA3} is be a consequence of the following result:

\begin{lemma}\label{sturk}
There exists a subcategory $\calM^{\nounit}_0 \subseteq \calM^{\nounit}$
with the following properties:
\begin{itemize}
\item[$(1)$] Let $\calM_0 = \calM^{\nounit}_0 \times_{ \calM^{\nounit} } \calM$. Then
the projection $\calM_0 \rightarrow \calM^{\nounit}_0$ is a trivial Kan fibration.

\item[$(2)$] Let $F \in \widetilde{\AdjDiag}^{\nounit}(\calC, \calD)$ be such that the image of
$F$ in $\NUAdj(\calC, \calD)$ belongs to $\UAdj(\calC, \calD)$. Then the associated
section $\Nerve(\cDelta'_{S} )^{op} \rightarrow \calM^{\nounit}$ of $p^{\nounit}$ factors through $\calM^{\nounit}_0$.
\end{itemize}
\end{lemma}

Let us assume Lemma \ref{sturk} for the moment.

\begin{proof}[Proof of Proposition \ref{partA3}]
In view of Lemma \ref{hyperstork}, it will suffice to show that 
the restriction map $r: \widetilde{\AdjDiag}(\calC, \calD) \rightarrow \widetilde{\AdjDiag}^{\nounit}(\calC, \calD)$ induces an equivalence of $\infty$-categories
$$ r_0: \widetilde{\AdjDiag}(\calC, \calD) \rightarrow \widetilde{\AdjDiag}^{\nounit}(\calC, \calD) \times_{ \NUAdj(\calC, \calD) } \UAdj(\calC, \calD).$$
Lemma \ref{hungertown} guarantees that $r$ factors through 
$\widetilde{\AdjDiag}^{\nounit}(\calC, \calD) \times_{ \NUAdj(\calC, \calD) } \UAdj(\calC, \calD)$, so that $r_0$ is well-defined. Now suppose that $\calM^{\nounit}_0$ and $\calM_0$ are as in
Lemma \ref{sturk}. Then $r_0$ is a pullback of the map
$$\bHom_{ \Nerve(\cDelta'_{S})^{op} } ( \Nerve( \cDelta'_{S} ),
\calM_0 ) \rightarrow \bHom_{ \Nerve(\cDelta'_{S})^{op} } ( \Nerve( \cDelta'_{S} ), \calM^{\nounit}_0 ).$$
Since the projection $\calM_0 \rightarrow \calM^{\nounit}_0$ is a trivial Kan fibration, we conclude that $r_0$ is a trivial Kan fibration.
\end{proof}

In order to produce the full subcategory $\calM^{\nounit}_0 \subseteq \calM^{\nounit}$ whose existence is asserted by Lemma \ref{sturk}, we will need to understand the $\infty$-categories
$\calM$ and $\calM_0$ a bit better. For this, we need a bit more notation.

\begin{notation}
Let $\cDelta_{S, +}$ be the category obtained from $\cDelta_{S}$ by adjoining a new initial object; we will think of this initial object as corresponding to a pair $(\emptyset, c)$, where $c: \emptyset \rightarrow \{ \calC, \calD \}$ denotes the empty word.

Fix an object $\sigma = ( [n], c) \in \cDelta_{S,+}$. 
We define functors
$\shift^{L}_{\sigma}, \shift^{R}_{\sigma}: \cDelta \rightarrow \cDelta_{S}$ as follows:
$$ \shift^{L}_{\sigma}( [m] ) = ( [m] \star [n], c' ) \quad \shift^{R}_{\sigma}( [m] ) = ( [n] \star [m], c''),$$
where 
$$c'(i) = \begin{cases} \calC & \text{if } 0 \leq i \leq m \\
c(j) & \text{if } i = n+1+j > n \end{cases} \quad
c''(i) = \begin{cases} c(i) & \text{if } 0 \leq i \leq n \\
\calC & \text{if } i>n. \end{cases}$$
We will let $\Fun^{\otimes}( \sigma, \calC)$ denote the fiber product
$\Nerve(\cDelta)^{op} \times_{ \Nerve(\cDelta_{S})^{op} }
\End^{\otimes}(\calC, \calD),$
where $\Nerve(\cDelta)^{op}$ maps to $\Nerve(\cDelta_{S})^{op}$ via
$\shift^L_{\sigma}$, and $\Fun^{\otimes}(\calC, \sigma)$ the fiber product
$ \Nerve(\cDelta)^{op} \times_{ \Nerve(\cDelta_{S})^{op} }
\End^{\otimes}(\calC, \calD)$
where $\Nerve(\cDelta)^{op}$ maps to $\Nerve(\cDelta_{S})^{op}$ via
$\shift^{R}_{\sigma}$. 

We observe that there is a commutative diagram
$$ \xymatrix{ \Fun^{\otimes}( \calC, \sigma) \ar[r] \ar[dr] & \End^{\otimes}(\calC) \ar[d] & \Fun^{\otimes}(\sigma, \calC) \ar[dl] \ar[l] \\
& \Nerve(\cDelta)^{op}. & }$$
This diagram exhibits $\Fun(\sigma, \calC) = \Fun^{\otimes}(\sigma, \calC)_{[0]}$ as left-tensored over the monoidal $\infty$-category $\End(\calC)$ and $\Fun(\calC, \sigma) = \Fun^{\otimes}(\calC, \sigma)_{[0]}$ as right-tensored over the monoidal $\infty$-category $\End(\calC)$. 
Note also that if $\sigma = ( \emptyset, c)$, then the horizontal arrows in the above diagram are isomorphisms.
\end{notation}

\begin{notation}\label{psied}
Fix an object $( [n], c) \in \cDelta'_{S}$, and let $\Sigma = \{ i \in [n]: c(i) = \calC \}$. 
We observe that the fiber $\pi^{-1} \{ ( [n], c) \}$ can be identified with a full subcategory 
$\calD$ of the product
$\cDelta^{\Sigma}_{+} = \prod_{i \in T} \cDelta_{+};$
this subcategory either coincides with $\cDelta^{\Sigma}_{+}$ (if $\Sigma \neq [n]$) or
with the subcategory obtained by deleting the initial object of $\cDelta^{S}_{+}$ (if
$\Sigma = [n]$). 

Let $E^{( [n], c)}: \calD \rightarrow \sSet$ denote the composition
$\calD \simeq \pi^{-1} \{ ( [n], c) \} \subseteq \calJJ
\stackrel{\psi}{\rightarrow} \cDelta_{S} \stackrel{E}{\rightarrow} \sSet,$
where the functor $E$ is defined as in Notation \ref{defEE}. Let
$E^{([n],c)}_{+}: \cDelta^{\Sigma}_{+} \rightarrow \sSet$ be given by the formula
$$E^{([n],c)}_{+}( X )  = 
\begin{cases} E^{([n],c)}(X) & \text{if } X \in \calD \\
\Delta^0 & \text{otherwise.} \end{cases}$$

Let $\cDelta_{+}^{\nounit} \subseteq \cDelta_{+}$ be the subcategory consisting of all objects of
$\cDelta_{+}$ and {\em injective} maps between them. For every subset $\Sigma_0 \subseteq \Sigma$, we set
$$T(\Sigma_0) = (\prod_{ i \in \Sigma-\Sigma_0} \Nerve( \cDelta^{\nounit}_{+})^{op} ) \times_{ \prod_{ j \in \Sigma_0} }
\Nerve( \cDelta_{+})^{op}.$$
$$ V^{\Sigma_0}( [n], c) = \bHom_{ \Nerve( \cDelta_{+}^{\Sigma})^{op} }(
T(\Sigma_0), \Nerve_{ E^{ ([n], c)}_{+} }( (\cDelta^{\Sigma}_{+})^{op}) ).$$
Let $V^{+}( [n], c) = V^{\Sigma}( [n], c)$ and $V^{-}( [n], c) = V^{\emptyset}( [n], c)$. We observe
that $V^{+}$ and $V^{-}$ can be viewed as functors from
$(\cDelta'_{S})^{op}$ to the category of simplicial sets, and that we have
canonical isomorphisms
$$ \calM \simeq \Nerve_{V^{+}}( (\cDelta'_{S})^{op} ) \quad \calM^{\nounit} \simeq \Nerve_{V^{-}}( ( \cDelta'_{S})^{op} ).$$
Moreover, the restriction map $\calM \rightarrow \calM^{\nounit}$ is induced by a natural transformation $V^{+} \rightarrow V^{-}$. More generally we have canonical maps
$ V^{\Sigma_0}( [n], c) \rightarrow V^{\Sigma_1}( [n], c)$ for every $\Sigma_1 \subseteq \Sigma$, given by restriction of functors.
\end{notation}

\begin{notation}
Let $( [n], c)$ and $\Sigma$ be as in Notation \ref{psied}. For every morphism
$\alpha: ( [m], c_0) \rightarrow ( [n], c)$ in $\calJJ$, we set
$$ [m]^{+} = \{ j \in [m]: \alpha(j) > i \}  \quad  [m]^{-} = \{ j \in [m]: \alpha(j) < i \} $$
$$ \sigma^{L}_{\alpha} = ( [m]^{+}, c_0 | [m]^{+} \} \in \cDelta_{S, +} \quad
 \sigma^{R}_{\alpha} = ( [m]^{-}, c_0 | [m]^{-} \} \in \cDelta_{S, +}. $$

Let $\pi: (\cDelta^{\nounit}_{+})^{\Sigma} \rightarrow ( \cDelta^{\nounit}_{+})^{ \Sigma - \{ i\} }$
be the restriction map, and
fix $X \in ( \cDelta^{\nounit}_{+})^{\Sigma - \{ i \}}$. We observe that, if
$\overline{X} \in \pi^{-1} \{ X \}$ corresponds
to an object $\alpha \in \calJJ_{ ( [n], c)}$, then
$\sigma^{L}_{\alpha}$, $\sigma^{R}_{\alpha} \in \cDelta_{S, +}$ are independent of the choice of $\overline{X}$. We will denote these words by $\sigma^{L}_X$ and $\sigma^{R}_{X}$, respectively. Restriction to $\pi^{-1} \{ X\}$ determines functors
$$ l_{i,X}: V^{-}( [n], c) \rightarrow \bHom_{ \Nerve(\cDelta)^{op} }
( \Nerve( \cDelta^{\nounit})^{op}, \Fun^{\otimes}( \sigma^{L}_X, \calC))$$
$$ r_{i,X}: V^{-}( [n], c) \rightarrow \bHom_{ \Nerve(\cDelta)^{op} }
( \Nerve( \cDelta^{\nounit})^{op}, \Fun^{\otimes}( \calC, \sigma^{R}_X)).$$

We will say that an object $v \in V^{-}( [n], c)$ is {\it $i$-good} if, for each $X \in ( \cDelta^{\nounit}_{+})^{\Sigma - \{ i\} }$, the functors $l_{i,X}$ and $r_{i,X}$ carry $v$ to quasi-unital module objects
of $\Fun^{\otimes}( \sigma^{L}_X, \calC)$ and $\Fun^{\otimes}( \calC, \sigma^{R}_{X})$, respectively In this case, both quasi-unital modules have the same underlying quasi-unital algebra
$A_{i,X}: \Nerve( \cDelta^{\nounit})^{op} \rightarrow \End^{\otimes}(\calC)$.
We will say that a morphism $v \rightarrow v'$ in $V^{-}( [n], c)$ is {\it $i$-good} if $v$ and $v'$ are
$i$-good and the induced map of nonunital algebras
$ A_{i,X} \rightarrow A'_{i,X}$
is quasi-unital, for every $X \in ( \cDelta^{\nounit}_{+})^{ \Sigma- \{i \}}$. 

Let us say that an object $v \in V^{-}( [n], c)$ is {\it good} if it is $i$-good for each $i \in \Sigma$, and
we will say that a morphism $v \rightarrow v'$ in $V^{-}( [n], c)$ is {\it good} if it is $i$-good for each
$i \in \Sigma$. More generally, for each $\Sigma_0 \subseteq \Sigma$, we will say that an object or morphism in
$V^{\Sigma_0}( [n], c)$ is {\it good} if its image in $V^{-}( [n], c)$ is good. The collection of good objects and morphsims in $V^{\Sigma_0}( [n], c)$ determines a subcategory $V^{\Sigma_0}_0( [n], c) \subseteq
V^{\Sigma_0}( [n], c)$. In particular, we have subcategories (generally not full)
$$ V^{+}_0( [n], c) \subseteq V^{+}( [n], c) \quad  V^{-}_0( [n], c) \subseteq V^{-}( [n], c).$$
It is not difficult to see that this construction determines a pair of functors
$V^{+}_0, V^{-}_0: ( \cDelta'_{S})^{op} \rightarrow \sSet.$
We now define
$$ \calM_0 \simeq \Nerve_{V^{+}_0}( (\cDelta'_{S})^{op} ), \quad \calM^{\nounit}_0 \simeq \Nerve_{V^{-}_0}( ( \cDelta'_{S})^{op} ).$$
\end{notation}

We are now almost ready to give the proof of Lemma \ref{sturk}. We require only one more preliminary:

\begin{lemma}\label{stabbus}
Let $A \subseteq B$ be a right anodyne inclusion of simplicial sets, $K$ an arbitrary simplicial set, and $p: X \rightarrow Y$ a categorical fibration. Then the map
$$ \Fun( B^{\triangleright} \times K, X)
\rightarrow \Fun( B^{\triangleright} \times K, Y)
\times_{ \Fun( (B \times K) \coprod_{ A \times K } (A^{\triangleright} \times K), Y)}
\Fun( (B \times K) \coprod_{ A \times K } (A^{\triangleright} \times K), X)$$
is a trivial Kan fibration.
\end{lemma}

\begin{proof}
Replacing $p$ by the induced map $X^K \rightarrow Y^K$, we can reduce to the case
where $K = \emptyset$. It now suffices to show that the inclusion
$i: B \coprod_{A} A^{\triangleright} \rightarrow B^{\triangleright}$ is a trivial cofibration (with respect to the Joyal model structure). This is equivalent to the assertion that $i$ has the left lifting property with respect to every categorical fibration. In fact, $i$ has the left lifting property with respect to every inner fibration (see Proposition \toposref{sharpen2}).  
\end{proof}

\begin{remark}\label{podus}
According to Proposition \toposref{cofbasic}, an inclusion of simplicial sets is right anodyne if and only if it is cofinal. In particular, Lemma \ref{stabbus} applies in the special case
where $B = \Nerve(\cDelta)^{op}$ and $A = \Nerve( \cDelta^{\nounit})^{op}$ (Lemma \toposref{bball4}). 
\end{remark}

\begin{proof}[Proof of Lemma \ref{sturk}]
We evidently have a Cartesian diagram
$$ \xymatrix{ \calM_0 \ar[r] \ar[d]^{f} & \calM \ar[d] \\
\calM^{\nounit}_0 \ar[r] & \calM^{\nounit}. }$$
It remains only to show that the subcategory $\calM^{\nounit}_0 \subseteq \calM^{\nounit}$ satisfies conditions $(1)$ and $(2)$.

We first prove $(1)$. We wish to show that the projection $f: \calM_0 \rightarrow \calM^{\nounit}_0$ is a trivial Kan fibration. Since $f$ is a categorical fibration (Lemma \ref{sulken2}), it will suffice to show that $f$ is a categorical equivalence. According to Proposition \toposref{apple1}, it will suffice to show that for every object $( [n], c) \in \cDelta'_{S}$, the induced map
$f_{([n], c)}: V^{+}_0( [n], c) \rightarrow V^{-}_{0}( [n], c)$
is a categorical equivalence.

Let $\Sigma = \{ i \in [n]: c(i) = \calC \}$. We will prove that, for each $\Sigma_1 \subseteq \Sigma_0 \subseteq \Sigma$, the restriction map $V^{\Sigma_0}_0( [n], c) \rightarrow V^{\Sigma_1}_0([n],c)$ is a trivial Kan fibration. Arguing inductively, it suffices to treat the case where $\Sigma_0 = \Sigma_1 \cup \{i\}$. Let 
$T = \cDelta_{+}^{ \Sigma_1} \times ( \cDelta^{\nounit}_{+})^{\Sigma - \Sigma_0}$, and set
$$Y =   \bHom_{ \Nerve( \cDelta_{+}^{\Sigma})^{op} }(
\Nerve(T)^{op} \times \Nerve(\cDelta)^{op} , \Nerve_{ E^{ ([n], c)}_{+} }( (\cDelta^{S}_{+})^{op}) )$$
$$Y^{\nounit} =  \bHom_{ \Nerve( \cDelta_{+}^{\Sigma})^{op} }(
\Nerve(T)^{op} \times \Nerve(\cDelta^{\nounit})^{op} , \Nerve_{ E^{ ([n], c)}_{+} }( (\cDelta^{\Sigma}_{+})^{op}) ).$$
Let $Y^{\nounit}_0$ be a subcategory of $Y^{\nounit}$ which contains the image of
the restriction map $V^{S_1}_0( [n], c) \rightarrow Y^{\nounit}$, and set $Y_0 = Y \times_{ Y^{\nounit} } Y_0$. We have a commutative diagram of simplicial sets
$$ \xymatrix{ V^{\Sigma_0}_0( [n], c) \ar[r]^-{g} \ar[dr] & V^{\Sigma_1}_0([n],c) \times_{ Y^{\nounit} } Y \ar[r] \ar[d]^{g''} & Y_0 \ar[r] \ar[d]^{g'} & Y \ar[d] \\
& V^{\Sigma_0}_0( [n], c) \ar[r] & Y^{\nounit}_0 \ar[r] & Y^{\nounit}. }$$
Using Lemma \ref{stabbus} (and Remark \ref{podus}), we deduce that $g$ is a trivial Kan fibration. Since each square in the diagram is a pullback, it will suffice to show that $g'$ is also a trivial Kan fibration, provided that $Y_0^{\nounit}$ is appropriately chosen.

We now analyze the map $g'$. For each object $t$ of $T$, consider the composition
$$ E_{t}: \cDelta^{op} \simeq (\cDelta^{ \{i \} })^{op} \times \{t\} \subseteq
(\cDelta_{+}^{S})^{op} \stackrel{ E^{ ( [n], c)}_{+} }{\rightarrow} \sSet.$$
Set 
$$H(t) = \bHom_{ \Nerve(\cDelta)^{op} }( \Nerve(\cDelta)^{op}, \Nerve_{ E_{t} }( \cDelta^{op} )), \quad H^{\nounit}(t) = \bHom_{ \Nerve(\cDelta)^{op} }( \Nerve(\cDelta^{\nounit})^{op}, \Nerve_{ E_{t} }( \cDelta^{op} )).$$
We observe that, as in Notation \ref{psied}, the object $t \in T$ determines a pair of words
$ \sigma^{L}_{t}, \sigma^{R}_t \in \cDelta_{\calC, \calD, +}$ and functors
$$ l_{t}: H^{\nounit}(t) \rightarrow \bHom_{ \Nerve(\cDelta)^{op} }
( \Nerve( \cDelta^{\nounit})^{op}, \Fun^{\otimes}( \sigma^{L}_t, \calC)$$
$$ r_{t}: H^{\nounit}(t) \rightarrow \bHom_{ \Nerve(\cDelta)^{op} }
( \Nerve( \cDelta^{\nounit})^{op}, \Fun^{\otimes}( \calC, \sigma^{R}_t).$$
Let $H^{\nounit}_0(t)$ denote the intersection
of $l_{t}^{-1} \Mod^{\qunit}( \Fun( \sigma^{L}_t, \calC))$ with
$r_{t}^{-1} \Mod^{\qunit}( \Fun( \calC, \sigma^{R}_{t} ) )$ (as subcategories of
$H^{\nounit}(t)$, and let
$H_0(t) = H(t) \times_{ H^{\nounit}(t)} H^{\nounit}_0(t)$. 
This construction determines functors
$$ H, H^{\nounit}, H_0, H_0^{\nounit}: T^{op} \rightarrow \sSet.$$
Set $Z = \Nerve_{H}( T^{op} )$, and define $Z_0$, $Z^{\nounit}$, and $Z^{\nounit}_0$ similarly.
We have canonical isomorphisms
$$ Y \simeq \bHom_{ \Nerve(T)^{op} }( \Nerve(T)^{op}, Z) \quad Y^{\nounit} \simeq \bHom_{ \Nerve(T)^{op} }( \Nerve(T)^{op}, Z^{\nounit})$$
Set $Y_0 = \bHom_{ \Nerve(T)^{op} }( \Nerve(T)^{op}, Z_0 )$ and
$Y^{\nounit}_0 = \bHom_{ \Nerve(T)^{op} }( \Nerve(T)^{op}, Z^{\nounit}_0)$. It is clear that we have
a Cartesian diagram
$$ \xymatrix{ Y_0 \ar[d]^{g'} \ar[r] & Y \ar[d] \\
Y^{\nounit}_0 \ar[r] & Y^{\nounit}, }$$
and that the restriction map $V^{S_1}_0( [n], c) \rightarrow Y^{\nounit}$ factors through
$Y^{\nounit}_0$. It remains only to show that $g'$ is a trivial fibration. For this, it will
suffice to show that the map $h: Z_0 \rightarrow Z^{\nounit}_{0}$ is a trivial Kan fibration.
The map $h$ is a categorical fibration (Lemma \ref{sulken2}), so it will suffice to show that $h$ is a categorical equivalence. Invoking Proposition \toposref{apple1}, we are reduced to proving that
for each $t \in T$, the map $h_{t}: H_0(t) \rightarrow H_0^{\nounit}(t)$ is an equivalence of $\infty$-categories. We observe that $h_{t}$ is induced by a map between the pullback squares
$$ \xymatrix{ H_0(t) \ar[r] \ar[d] & \Mod( \Fun( \sigma^L_{t}, \calC) ) \ar[d] & H_0^{\nounit}(t) \ar[r] \ar[d] & \Mod^{\qunit}( \Fun( \sigma^{L}_{t}, \calC)) \ar[d] \\
\Mod( \Fun( \calC, \sigma^{R}_{t}) ) \ar[r] & \Alg( \End(\calC) ) & \Mod^{\qunit}( \Fun( \calC, \sigma^{R}_{t} )) \ar[r] & \Alg^{\qunit}( \End(\calC) ). }$$
The desired result now follows from Theorem \ref{uniqueunit} and (two applications of)
Corollary \ref{uniquemod}. This completes the proof of $(1)$.

We now prove $(2)$. Let $s \in \widetilde{\AdjDiag}^{\nounit}(\calC, \calD)$, and let
$F \circ G \stackrel{v}{\rightarrow} U \simeq \id_{\calD} $
denote the associated object in $\NUAdj(\calC, \calD)$. Suppose that $v$ is the counit of an adjunction between $F$ and $G$. 
We wish to prove that the associated section of the projection $\calM \rightarrow \Nerve( \cDelta'_{S})^{op}$ factors through $\calM_0$. We observe that $s$ determines a
nonunital algebra $T \in \Alg^{\nounit}( \End(\calC) )$, whose underlying object of
$\End(\calC)$ can be identified with $G \circ F$. Unwinding the definitions, we see that
the desired result is equivalent to the following assertions:
\begin{itemize}
\item[$(i)$] The nonunital algebra $T \in \Alg^{\nounit}( \End(\calC ))$ admits a quasi-unit $u$.
\item[$(ii)$] The left actions of $T$ on $G$ and $G \circ F$ are quasi-unital.
\item[$(iii)$] The right actions of $T$ on $F$ and $F \circ G$ are quasi-unital.
\end{itemize}
A quasi-unit $u$ for $T$ can be identified with a morphism $u: \id_{\calC} \rightarrow G \circ F$. 
Moreover, it is easy to see that $u$ satisfies $(i)$, $(ii)$, and $(iii)$ if and only if the compositions
$$ F \simeq F \circ \id_{\calC} \stackrel{u}{\rightarrow} F \circ (G \circ F) = (F \circ G) \circ F
\stackrel{v}{\rightarrow} \id_{\calD} \circ F = F$$
$$ G \simeq \id_{\calC} \circ G \stackrel{u}{\rightarrow} (G \circ F) \circ G = G \circ (F \circ G) \stackrel{v}{\rightarrow} G \circ \id_{\calD} = G$$
are homotopic to the identity. In other words, $(2)$ is equivalent to the assertion that there exists
a map $u: \id_{\calC} \rightarrow G \circ F$ which is a compatible unit for the adjunction determined by the counit $v$. This follows immediately from our assumption that the image of $s$ lies
in $\UAdj(\calC, \calD) \subseteq \NUAdj(\calC, \calD)$.
\end{proof}

\section{The Monoidal Structure on Stable Homotopy Theory}\label{monoid6.0}

Let $\Spectra$ denote the $\infty$-category of {\em spectra} (as defined in \S \stableref{stable8}). The homotopy category $\h{\Spectra}$ can be identified with the classical {\em stable homotopy category}. Given a pair of spectra $X, Y \in \h{\Spectra}$, one can define new spectrum called the {\it smash product} of $X$ and $Y$. The smash product operation determines a monoidal structure on $\h{\Spectra}$. In this section, we will show that this monoidal structure is determined by a monoidal structure which exists on the $\infty$-category $\Spectra$ itself. There are at least three ways to see this.

\begin{itemize}
\item[$(S1)$] Choose a simplicial model category $\bfA$ equipped with a compatible monoidal structure, whose underlying $\infty$-category is equivalent to $\Spectra$. For example, we can take $\bfA$ to be the category of {\it symmetric spectra} (see \cite{symmetricspectra}). According to  Proposition \ref{hurgoven}, the underlying $\infty$-category $\Nerve(\bfA^{\degree}) \simeq \Spectra$ is endowed with a monoidal structure. The advantage of this perspective is that it permits us to easily compare the algebras and modules considered in this paper with more classical approaches to the theory of structured ring spectra. For example, Theorem \ref{beckify} implies that $\Alg(\Spectra)$ is (equivalent to) the $\infty$-category underlying the model category of algebras in symmetric spectra (that is, strictly associative monoids in $\bfA$); see Example \ref{bulwork}.

The main disadvantage of this approach is that it seems to require auxiliary data (namely, a strictly associative model for the smash product functor), which could be supplied in many different ways. From a conceptual point of view, the existence of such a model ought to be irrelevant: the very purpose of higher category theory is to provide a formalism which allows us to avoid assumptions like strict associativity. 

\item[$(S2)$] Let $\LFun(\Spectra, \Spectra)$ denote the full subcategory of $\Fun(\Spectra, \Spectra)$ spanned by those functors from $\Spectra$ to $\Spectra$ which preserve small colimits. Corollary \stableref{choccrok} asserts that evaluation on the sphere spectrum yields an equivalence of $\infty$-categories $\LFun(\Spectra, \Spectra) \rightarrow \Spectra$. On the other hand, since $\LFun(\Spectra, \Spectra)$ is stable under composition in $\Fun(\Spectra, \Spectra)$, the composition monoidal structure on $\Fun(\Spectra, \Spectra)$ induces a monoidal structure on
$\LFun(\Spectra, \Spectra)$. This definition has the virtue of being very concrete (the smash product operation is simply given by composition of functors), and it allows us to identify the algebra objects of $\LFun(\Spectra, \Spectra)$: they are precisely the colimit-preserving monads on the $\infty$-category $\Spectra$ (for an application of this last observation, see Theorem \ref{schwedeshipley}). The disadvantage of this definition is that it is very ``associative'' in nature, and therefore does not generalize easily to show that $\Spectra$ is a {\em symmetric} monoidal $\infty$-category (as we will see in \cite{symmetric}). 

\item[$(S3)$] Let $\LPress$ be the $\infty$-category whose objects are presentable $\infty$-categories and whose morphisms are colimit-preserving functors (see \S \toposref{colpres}), and let
$\LPressStab$ be the full subcategory of $\LPress$ spanned by those presentable $\infty$-categories which are stable. As we will explain below, $\LPressStab$ admits a monoidal structure. Moreover, the algebra objects of $\LPressStab$ can be identified with monoidal $\infty$-categories $\calC$ which are stable, presentable, and have the property that the bifunctor $\otimes: \calC \times \calC \rightarrow \calC$ is colimit preserving separately in each variable. We will show that
$\Spectra$ is the unit object of $\LPressStab$ (with respect to its tensor structure). It follows from
Proposition \ref{gurgle} that $\Spectra$ can be endowed with the structure of an algebra object of $\LPressStab$, and in fact is {\em initial} among such algebra objects. This establishes not only the existence of a monoidal structure on $\Spectra$, but also a universal property which can be used to prove uniqueness (Corollary \ref{surcoi}). 
\end{itemize}

We will follow approach $(S3)$. Our first step is to construct a monoidal
structure on the $\infty$-category $\LPress$ of presentable $\infty$-categories. This construction will be carried out in \S \ref{jurmit}. Roughly speaking, given a pair of presentable $\infty$-categories
$\calC$ and $\calD$, the tensor product $\calC \otimes \calD$ is universal among presentable $\infty$-categories which receive a bifunctor $\calC \times \calD \rightarrow \calC \otimes \calD$
which preserves (small) colimits separately in each variable.
In \S \ref{hummingburg}, we will see that this monoidal structure induces (via the constructions of \S \ref{locol}) another monoidal structure on the $\infty$-category of {\em stable} presentable $\infty$-categories, which we can use to carry out the reasoning outlined in the above discussion.

We will define an {\it $A_{\infty}$-ring} to be an algebra object of the $\infty$-category $\Spectra$.
In \S \ref{monoid9}, we will study some of the basic formal properties of the $\infty$-category
of $A_{\infty}$-rings. For example, we will show that an $A_{\infty}$-ring $A$ satisfying
$\pi_{i} A \simeq 0$ for all $i \neq 0$ is essentially the same thing as an associative ring, in the sense of classical algebra (Proposition \ref{umberilt}).\index{$A_{\infty}$-ring}

If $A$ is an $A_{\infty}$-ring, then we have an associated $\infty$-category
$\Mod^{L}_{A} = \Mod_{A}(\Spectra)$ of {\it left $A$-module spectra}, or simply {\it left $A$-modules}. In \S \ref{monoid11}, we will study the $\infty$-category $\Mod^{L}_{A}$ in some detail. In particular, we will show that $\Mod^{L}_{A}$ is a {\em stable} $\infty$-category, so that the homotopy category
$\h{\Mod^{L}_{A}}$ is triangulated (Proposition \ref{stabult}). In the case where $A$ corresponds to 
an ordinary associative ring, $\h{\Mod^{L}_{A}}$ can be identified with the usual {\em derived category} of that associative ring (see Proposition \ref{derivdisc} and Remark \ref{ultra}).

There is also an (entirely dual) theory of {\em right modules} over $A$. Given a right $A$-module $M$ and every left $A$-module $N$, one can construct a {\it relative tensor product} $M \otimes_{A} N \in \Spectra$. We will describe this construction (in a very general setting) in \S \ref{balpair}.
In \S \ref{monoid13}, we will use the theory of relative tensor products to introduce a theory of
{\em flat} modules. Roughly speaking, a left $A$-module $N$ is flat if the tensor product functor
$M \mapsto M \otimes_{A} N$ has good exactness properties (at least when $A$ is connective). Our main result is an analogue of Lazard's theorem: any flat $A$-module can be obtained as a filtered colimit of free $A$-modules (Theorem \ref{lazard}). To prove this, we will introduce a spectral sequence which can be used (in favorable cases) to compute the homotopy groups $\pi_{n} (M \otimes_{A} N)$ in terms of the homotopy groups of $M$, $A$, and $N$.

A theorem of Schwede and Shipley asserts that a stable $\infty$-category $\calC$ can be realized as the $\infty$-category of left modules over an $A_{\infty}$-ring $A$ if and only if $\calC$ is presentable and compactly generated (we will give a proof of this result in \S \ref{monoid11}, using the $\infty$-categorical Barr-Beck theorem). In particular, the $\infty$-category $\Mod_{A}^{L}$ is compactly generated. The compact objects of $\Mod_{A}^{L}$ are called {\it perfect} $A$-modules. We will study the class of perfect $A$-modules in \S \ref{finprop}, together with the somewhat larger class of {\it almost perfect} $A$-modules which arise frequently in applications.

\subsection{Tensor Products of Presentable $\infty$-Categories}\label{jurmit}

Our goal in this section is to show that the $\infty$-category $\LPress$ of presentable $\infty$-categories admits a monoidal structure. This monoidal structure can be described informally as follows: given a pair of presentable $\infty$-categories $\calC$ and $\calD$, the tensor product
$\calC \otimes \calD$ is the recipient of a universal bifunctor $\calC \times \calD \rightarrow \calC \otimes \calD$ which is ``bilinear''; that is, which preserves colimits separately in each variable.

We begin by introducing a bit of notation.

\begin{notation}
Let $\widehat{\Cat}_{\infty}$ denote the $\infty$-category of (not necessarily small) $\infty$-categories. Then $\widehat{\Cat}_{\infty}$ admits finite products. Consequently, there
exists a monoidal $\infty$-category $\widehat{\Cat}_{\infty}^{\otimes}$ endowed with a Cartesian structure $\widehat{\Cat}_{\infty}^{\otimes} \rightarrow \widehat{\Cat}_{\infty}$ which induces
an equivalence $( \widehat{\Cat}_{\infty}^{\otimes})_{[1]} \simeq \widehat{\Cat}_{\infty}$.

We can construct an explicit model for $\widehat{\Cat}_{\infty}^{\otimes}$ as follows. 
Let $\bfA$ be the category of (not necessarily small) marked simplicial sets, as in Example \ref{exalcun}. Then $\bfA$ is endowed with the structure of a monoidal model category, with monoidal structure given by Cartesian product. We now let $\widehat{\Cat}_{\infty}^{\otimes} = \Nerve( \bfA^{\otimes, \degree})$, as defined in the statement of Proposition \ref{hurgoven}. More concretely: 
\begin{itemize}
\item[$(i)$] The objects of $\widehat{\Cat}_{\infty}^{\otimes}$ are finite sequences
$[ X_1, \ldots, X_n ]$, where each $X_i$ is an $\infty$-category.
\item[$(ii)$] Given a pair of objects $[X_1, \ldots, X_n], [Y_1, \ldots, Y_m] \in \widehat{\Cat}_{\infty}^{\otimes}$, a morphism from $[X_1, \ldots, X_n]$ to $[Y_1, \ldots, Y_m]$ consists of an
order-preserving map $f: [m] \rightarrow [n]$ and a collection of functors 
$\eta_{i}: X_{f(i-1) +1} \times \ldots \times X_{f(i)} \rightarrow Y_{i}$. 
\end{itemize}

We define a subcategory $\MonLPres \subseteq \widehat{\Cat}^{\otimes}_{\infty}$
as follows:

\begin{itemize}\index{ZZZMonLPres@$\MonLPres$}
\item[$(iii)$] An object $[X_1, \ldots, X_n] \in \widehat{\Cat}^{\otimes}_{\infty}$ belongs to
$\MonLPres$ if and only if each $X_i$ is a presentable $\infty$-category.
\item[$(iv)$] Let $[X_1, \ldots, X_n], [Y_1, \ldots, Y_m] \in \MonLPres$, and let
$F: [X_1, \ldots, X_n] \rightarrow [Y_1, \ldots, Y_m]$ be a morphism in $\widehat{\Cat}_{\infty}^{\otimes}$ covering a map $f: [m] \rightarrow [n]$ in $\cDelta$. Then $F$ belongs to
$\MonLPres$ if and only if each of the induced functors
$\eta_{i}: X_{f(i-1) +1} \times \ldots \times X_{f(i)} \rightarrow Y_{i}$
is colimit-preserving in each variable.
\end{itemize}
\end{notation}

Our next goal is to show that $\MonLPres$ is a monoidal $\infty$-category. First, we need a lemma. 

\begin{notation}\index{ZZZRFun@$\RFun(\calC, \calD)$}\index{ZZZRFunn@$\LFun(\calC, \calD)$}
Let $\calC$ and $\calD$ be $\infty$-categories which admit small limits. Then we let $\RFun(\calC, \calD)$ denote the full subcategory of $\Fun(\calC, \calD)$ spanned by those functors which preserve small limits. Similarly, if $\calC$ and $\calD$ admit small colimits, we let $\LFun(\calC, \calD)$ denote the full subcategory of $\Fun(\calC, \calD)$ which preserve small colimits.
\end{notation}

\begin{warning}
We used the same notation in \S \toposref{afunc5} for a slightly different purpose: there,
$\LFun( \calC, \calD)$ denoted the $\infty$-category of functors from $\calC$ to $\calD$ which
were left adjoints, and $\RFun(\calC, \calD)$ the $\infty$-category of functors which were right adjoints. In what follows below, we will work with {\em presentable} $\infty$-categories, so that the two notations are {\em almost} consistent with one another in view of the adjoint functor theorem(Corollary \toposref{adjointfunctor}).
\end{warning}

\begin{lemma}\label{studymax}
Let $\calC$ and $\calD$ be presentable $\infty$-categories. Then $\RFun( \calC^{op}, \calD)$ is a presentable $\infty$-category.
\end{lemma}

\begin{proof}
Using Theorem \toposref{pretop} and the results of \S \toposref{invloc}, we can choose a small $\infty$-category $\calC'$, a small collection $S$ of morphisms in $\calP(\calC')$, and an equivalence $\calC \simeq S^{-1} \calP(\calC')$. Then
$$\RFun( \calP(\calC')^{op}, \calD) \simeq \LFun( \calP(\calC'), \calD^{op})^{op} \simeq
\Fun( \calC', \calD^{op})^{op} \simeq \Fun( {\calC'}^{op}, \calD)$$
is presentable by Proposition \toposref{presexp}, where the second equivalence is given by composition with the Yoneda embedding (Theorem \toposref{charpresheaf}). For each morphism
$\alpha \in S$, let $\calE(\alpha)$ denote the full subcategory of $\RFun( \calP(\calC')^{op}, \calD)$ spanned by those functors which carry $\alpha$ to an equivalence in $\calD$. Then
$\RFun(\calC^{op}, \calD)$ is equivalent to the intersection $\bigcap_{\alpha \in S} \calE(\alpha)$. 
In view of Lemma \toposref{stur3}, it will suffice to show that each $\calE(\alpha)$ is a localization
of $\RFun( \calP(\calC')^{op}, \calD)$. We now observe that $\calE(\alpha)$ is given by a pullback diagram
$$ \xymatrix{ \calE(\alpha) \ar[r] \ar[d] & \RFun( \calP(\calC'')^{op}, \calD) \ar[d] \\
\calE \ar@{^{(}->}[r] & \Fun( \Delta^1, \calD), }$$
where $\calE$ denotes the full subcategory of $\Fun(\Delta^1, \calD)$ spanned by the equivalences. According to Lemma \toposref{stur2}, it will suffice to show that $\calE$ is an accessible localization of $\Fun(\Delta^1, \calD)$, which is clear.
\end{proof}

\begin{lemma}\label{jackus}\index{tensor product!of presentable $\infty$-categories}
Let $\calC_1, \ldots, \calC_{n}$ be a finite collection of presentable $\infty$-categories. Then
there exists a presentable $\infty$-category $\calC_1 \otimes \ldots \otimes \calC_{n}$ and a functor
$f: \calC_1 \times \ldots \times \calC_{n} \rightarrow \calC_{1} \otimes \ldots \otimes \calC_{n}$
with the following properties:
\begin{itemize}
\item[$(1)$] The functor $f$ preserves colimits separately in each variable.
\item[$(2)$] For every presentable $\infty$-category $\calD$, composition with $f$ induces
an equivalence from $\LFun( \calC_1 \otimes \ldots \otimes \calC_n, \calD)$ onto the full subcategory of $\Fun( \calC_1 \times \ldots \times \calC_n, \calD)$ spanned by those functors which preserve colimits separately in each variable. 
\end{itemize}
\end{lemma}

\begin{proof}
The proof goes by induction on $n$. If $n = 0$, we take $\calC_1 \otimes \ldots \otimes \calC_n = \SSet$, choose the functor $f$ to classify a final object of $\SSet$, and apply Theorem \toposref{charpresheaf}. If $n=1$ there is nothing to prove. Assume therefore that $n > 1$.

Let $\calD$ be an arbitrary presentable $\infty$-category, and let
$\Fun'(\calC_1 \times \ldots \times \calC_n, \calD)$ be the full subcategory of $\Fun( \calC_1 \times \ldots \times \calC_n, \calD)$ which preserve colimits separately in each variable. Then we have a canonical isomorphism
$\Fun'( \calC_1 \times \ldots \calC_n, \calD) \simeq
\Fun'( \calC_1\times \ldots \calC_{n-1}, \LFun( \calC_n, \calD)).$
Here the $\infty$-category $\LFun(\calC_n, \calD)$ is presentable (Proposition \toposref{intmap}), so that the inductive hypothesis allows us to identify this $\infty$-category with
$\LFun( \calC_1 \otimes \ldots \otimes \calC_{n-1}, \LFun( \calC_n, \calD)).$
If $n > 2$, then we can conclude by identifing the last $\infty$-category with
$\LFun( (\calC_1 \otimes \ldots \otimes \calC_{n-1}) \otimes \calC_{n}, \calD).$

Suppose instead that $n = 2$. Using Corollary \toposref{adjointfunctor} and Proposition \toposref{adjobs}, we can identify $\LFun( \calC_2, \calD)$ with the full subcategory of $\RFun( \calD, \calC_2)^{op}$ spanned by those functors which are accessible. Consequently, we
get a fully faithful embedding
\begin{eqnarray*}
\Fun'( \calC_1 \times \calC_2, \calD) & \rightarrow &
\LFun( \calC_1, \RFun( \calD, \calC_2)^{op}) \\
& \simeq & \LFun( \calC_1, \LFun( \calD^{op}, \calC_2^{op})) \\
& \simeq & \LFun( \calD^{op}, \LFun( \calC_1, \calC_2^{op})) \\
& \simeq & \RFun( \calD, \RFun( \calC_1^{op}, \calC_2) )^{op}
\end{eqnarray*}
whose essential image consists of the collection of accessible functors from
$\calD$ to $\RFun( \calC_1^{op}, \calC_2)$. We now apply Lemma \ref{studymax} to conclude that $\RFun( \calC_1^{op}, \calC_2)$, so that (using Corollary \toposref{adjointfunctor} and Proposition \toposref{adjobs} again) $\Fun'( \calC_1 \times \calC_2, \calD)$ can be identified with
$\LFun( \RFun( \calC_1^{op}, \calC_2), \calD)$. We now conclude by defining $\calC_1 \otimes \calC_2$ to be $\RFun( \calC_1^{op}, \calC_2)$. 
\end{proof}

\begin{remark}\label{huger}
Let $\calC_1, \ldots, \calC_n$ be a finite sequence of presentable $\infty$-categories.
It follows from the proof of Lemma \ref{jackus} that, for $n > 0$, the tensor product can be canonically identified with the iterated functor $\infty$-category
$ \RFun( \calC^{op}_1, \RFun( \calC^{op}_2, \ldots \RFun( \calC^{op}_{n-1}, \calC_n) \ldots )).$
Combining this observation
with Theorem \toposref{surbus}, we conclude that the bifunctor
$\otimes: \LPress \times \LPress \rightarrow \LPress$ preserves colimits separately in each variable
(remember that {\em colimits} in $\LPress$ can also be computed as {\em limits} in $\RPres$, which are computed by forming limits in $\Cat_{\infty}$ by Theorem \toposref{surbus}). Alternatively, one can observe that $\LPress$ is actually a {\em closed} monoidal category,
with internal mapping objects given by $\LFun(\calC, \calD)$ (see Proposition \toposref{intmap}).
\end{remark}

\begin{proposition}\label{appbilg}\index{monoidal $\infty$-category!of presentable $\infty$-categories}
The natural map $q: \MonLPres \rightarrow \Nerve(\cDelta)^{op}$ determines a monoidal structure on the $\infty$-category $\LPress \simeq (\MonLPres)_{[1]}$. Moreover, the inclusion functor
$\MonLPres \subseteq \widehat{\Cat}^{\otimes}_{\infty}$ is lax monoidal.
\end{proposition}

\begin{proof}
We will show that $q$ is a coCartesian fibration; the remaining condition that $q$ induces equivalences $(\MonLPres)_{[n]} \simeq ( \LPress)^{n}$ will then follow by inspection. According to 
Proposition \toposref{gotta}, it will suffice to verify the following conditions:
\begin{itemize}
\item[$(1)$] Let $\alpha: [m] \rightarrow [n]$ be a morphism in $\cDelta$, and let
$X_{[n]}$ be an object of $\MonLPres$ lying over $[n]$. Then there exists a locally $q$-coCartesian
morphism $X_{[n]} \rightarrow X_{[m]}$ in $\MonLPres$ covering $\alpha$.
\item[$(2)$] Let $\alpha: [m] \rightarrow [n]$ and $\beta: [n] \rightarrow [p]$ be morphisms
in $\cDelta$, covered by locally $q$-coCartesian morphisms
$$ X_{[p]} \rightarrow X_{[n]} \rightarrow X_{[m]} $$
in $\MonLPres$. Then the composite morphism $X_{[p]} \rightarrow X_{[m]}$ is locally $q$-coCartesian.
\end{itemize}

Using the product structure on the fibers of $q$, we can reduce to the case where $m=1$ and the morphisms $\alpha: [m] \rightarrow [n]$ and $\beta: [n] \rightarrow [p]$ preserve the endpoints of intervals. Unwinding the definitions, we obtain the following reformulation of conditions $(1)$ and $(2)$:

\begin{itemize}
\item[$(1')$] Let $\calC_1, \ldots, \calC_{n}$ be a finite collection of presentable $\infty$-categories. There exists a presentable $\infty$-category $\calC_1 \otimes \ldots \otimes \calC_{n}$ and a functor
$f: \calC_1 \times \ldots \times \calC_{n} \rightarrow \calC_{1} \otimes \ldots \otimes \calC_{n},$
which preserves colimits separately in each variable and which is universal in the following sense:
\begin{itemize}
\item[$(\ast)$] Let $\calD$ be a presentable $\infty$-category, let 
$\Fun'( \calC_1 \times \ldots \times \calC_n, \calD) \subseteq \Fun( \calC_1 \times \ldots \times \calC_n, \calD)$ be the full subcategory spanned by those functors which preserve colimits separately in each variable. Then the functor
$ \LFun( \calC_1 \otimes \ldots \otimes \calC_n, \calD) \rightarrow \Fun'( \calC_1 \times \ldots \times \calC_n, \calD)$
given by composition with $f$ induces a homotopy equivalence between maximal Kan complexes contained in $\LFun( \calC_1 \otimes \ldots \otimes \calC_n, \calD)$ and
$\Fun'( \calC_1 \times \ldots \times \calC_n, \calD)$.
\end{itemize}
\item[$(2')$] Given an endpoint-preserving map $\alpha: [n] \rightarrow [p]$ and a collection of presentable $\infty$-categories $\calC_1, \ldots, \calC_{p}$. Then the canonical map
$$ \calC_1 \otimes \ldots \otimes \calC_{p} \rightarrow ( \calC_{ \alpha(0)+1 } \otimes \ldots \otimes \calC_{\alpha(1)} ) \otimes \ldots \otimes ( \calC_{ \alpha(n-1)+1} \otimes \ldots \otimes \calC_{\alpha(n)} )$$
is an equivalence of (presentable) $\infty$-categories.
\end{itemize}
Assertion $(1')$ follows immediately from Lemma \ref{jackus}, and $(2')$ follows from Remark \ref{huger}.
\end{proof}

\begin{remark}\label{otherlandar}
In view of Remark \ref{otherlander}, we can identify $\Alg( \widehat{\Cat}_{\infty})$
with the $\infty$-category of monoidal $\infty$-categories. The (lax monoidal) inclusion
$ \MonLPres \subseteq \widehat{\Cat}_{\infty}^{\otimes}$
allows us to identify $\Alg( \LPress)$ with a subcategory of $\Alg( \widehat{\Cat}_{\infty})$. Unwinding the definitions, we see that a monoidal $\infty$-category $\calC$ belongs to $\CAlg(\LPress)$ if and only if $\calC$ is presentable, and the tensor product operation $\otimes: \calC \times \calC \rightarrow \calC$ preserves (small) colimits separately in each variable.
\end{remark}

\begin{example}
Let $\calX$ and $\calY$ be $\infty$-topoi. Then $\calX \otimes \calY$ is an $\infty$-topos, and can be identified with the (Cartesian) product of $\calX$ and $\calY$ in the $\infty$-category of $\infty$-topoi. For a proof of a slightly weaker assertion, we refer the reader to Theorem \toposref{conprod}. The general statement can be proved using the same argument.
\end{example}

\begin{example}\label{intin}
The proof of Proposition \ref{appbilg} shows that the $\infty$-category $\SSet$ is the unit object of $\LPress$. In particular, for every presentable $\infty$-category $\calC$ we have canonical equivalences $ \calC \simeq \calC \otimes \SSet \simeq \RFun(\calC^{op}, \SSet).$
The essential surjectivity of the composition is a restatement of the representability criterion of Proposition \toposref{representable}.
\end{example}

\begin{example}\label{presreower}
Recall that, for every $\infty$-category $\calC$, the $\infty$-category $\calC_{\ast}$ of {\it pointed objects} of $\calC$ is defined to be the full subcategory of $\Fun( \Delta^1, \calC)$ spanned by those functors $F: \Delta^1 \rightarrow \calC$ for which $F(0)$ is a final object of $\calC$.
The canonical isomorphism of simplicial sets
$ \RFun( \calC^{op}, \calD_{\ast}) \simeq \RFun( \calC^{op}, \calD)_{\ast}$
induces an equivalence $\calC \otimes \calD_{\ast} \simeq (\calC \otimes \calD)_{\ast}$ for every pair of presentable $\infty$-categories $\calC, \calD \in \LPress$. In particular, we have a canonical equivalence $\calC \otimes \SSet_{\ast} \simeq \calC_{\ast}$.
\end{example}

\begin{example}\label{srower}
Let $\SSet_{\infty}$ denote the $\infty$-category of spectra. Then $\SSet_{\infty}$ can be identified with a homotopy limit of the tower
$$ \ldots \stackrel{\Omega}{\rightarrow} \SSet_{\ast} \stackrel{\Omega}{\rightarrow} \SSet_{\ast}.$$
Consequently, for every presentable $\infty$-category $\calC$, we have equivalences
$$ \calC \otimes \SSet_{\infty} \simeq \RFun( \calC^{op}, \SSet_{\infty} )
\simeq \holim \{ \RFun( \calC^{op}, \SSet_{\ast}) \}
\simeq \holim \{ \calC_{\ast} \} \simeq \Stab(\calC),$$
where $\Stab(\calC)$ denotes the {\it stabilization} of $\calC$ defined in \S \stableref{stable9}.
\end{example}

\subsection{The Smash Product Monoidal Structure}\label{hummingburg}

In this section, we will apply the results of \S \ref{jurmit} to construct a monoidal structure on the $\infty$-category $\Spectra$ of spectra. Our first step is to construct a monoidal structure on
the $\infty$-category of {\em stable} presentable monoidal $\infty$-categories.

\begin{notation}\index{ZZZMonStab@$\MonStab$}
Let $\MonStab \subseteq \MonLPres$ denote the full subcategory spanned by the objects
$[ X_1, \ldots, X_n]$ where each $X_{i}$ is a {\em stable} presentable $\infty$-category.
\end{notation}

\begin{remark}\label{suty}
Let $\LPressStab$ denote the full subcategory of $\LPress$ spanned by the stable $\infty$-categories.
In view of Corollary \stableref{mapprop}, the inclusion $\LPressStab \subseteq \LPress$ admits a left
adjoint $L: \LPress \rightarrow \LPressStab$, given by the formula
$$ \calC \mapsto \Stab(\calC).$$
In view of Example \ref{srower}, we can identify $L$ with the functor
$\calC \mapsto \calC \otimes \SSet_{\infty}$. It follows easily that $L$
is compatible with the monoidal structure on $\LPress$. (We could equally well make the same argument with ``stable'' replaced by ``pointed'', and Example \ref{presreower} in place of
Example \ref{srower}).
\end{remark}

\begin{proposition}\label{tupid}\index{monoidal $\infty$-category!of presentable stable $\infty$-categories}

\begin{itemize}
\item[$(1)$] The projection $\MonStab \rightarrow \Nerve(\cDelta)^{op}$ determines a monoidal structure on the $\infty$-category 
$\LPressStab \simeq ( \MonStab )_{[1]}$.
\item[$(2)$] The inclusion functor $i:\MonStab \subseteq \MonLPres$ is a lax monoidal functor.
\item[$(3)$] The localization functor $\Stab: \LPress \rightarrow \LPressStab$ extends to a monoidal
functor from $\MonLPres$ to $\MonStab$.
\item[$(4)$] The unit object of $\LPressStab$ is the $\infty$-category $\Spectra$ of spectra.
\end{itemize}
\end{proposition}

\begin{proof}
Assertions $(1)$ through $(3)$ follow from Proposition \ref{localjerk} and Remark \ref{suty}. Assertion $(4)$ follows from $(3)$, Example \ref{intin}, and the observation that
$\Spectra \simeq \Stab(\SSet)$.
\end{proof}

\begin{remark}\label{otherlandur}
Applying the reasoning of Remark \ref{otherlandar}, we deduce that
$\Alg( \LPressStab)$ can be identified with a full subcategory of $\Alg( \LPressStab)$. A
monoidal $\infty$-category $\calC$ belongs to this full subcategory if and only if
$\calC$ is stable, presentable, and the bifunctor $\otimes: \calC \times \calC \rightarrow \calC$ preserves colimits separately in each variable.
\end{remark}

\begin{remark}\label{tupi}
The inclusion functor $i: \MonStab \subseteq \MonLPres$ is {\em almost} a monoidal functor.
Unwinding the definitions, we see that $i$ is monoidal if and only if, for every
finite collection $[ \calC_1, \ldots, \calC_n]$ of presentable stable $\infty$-categories, the
tensor product 
$\calC_1 \otimes \ldots \otimes \calC_n$ is again stable. This is true provided that $n > 0$; in fact, Remark \ref{huger} implies that $\calC_1 \otimes \ldots \otimes \calC_n$ is stable provided that
$\calC_j$ is stable for some $j \in \{1, \ldots, n\}$. However, in the degenerate case $n=0$, the
tensor product is equivalent to $\SSet$, which is not stable. In other words, $i$ fails to be a monoidal functor only because it does not preserve the unit object.
\end{remark}

\begin{corollary}\label{surcoi}
Let $\widehat{\Cat}_{\infty}^{\Mon}$
denote the $\infty$-category of $($not necessarily small$)$ monoidal $\infty$-categories. Let
$\widehat{\Cat}_{\infty}^{\Stabb, \Mon}$ denote the subcategory whose objects are
required to be stable presentable monoidal $\infty$-categories such that the bifunctor
$\otimes$ preserves colimits in separately in each variable, and whose morphisms
are given by colimit-preserving monoidal functors. Then:

\begin{itemize}
\item[$(1)$] The $\infty$-category $\widehat{\Cat}_{\infty}^{\Stabb, \Mon}$ has an initial object $\calC^{\otimes}$.

\item[$(2)$] The underlying $\infty$-category $\calC^{\otimes}_{[1]}$ is equivalent to the $\infty$-category of spectra.

\item[$(3)$] Let $\calD^{\otimes}$ be an arbitrary monoidal $\infty$-category. Suppose that:
\begin{itemize}
\item[$(i)$] The underlying $\infty$-category $\calD = \calD^{\otimes}_{[1]}$ is stable and presentable.
\item[$(ii)$] The functor $\Spectra \rightarrow \calD$ determined by the unit object of $\calD$ $($see Corollary \stableref{choccrok}$)$ is an equivalence of $\infty$-categories.
\item[$(iii)$] The bifunctor $\otimes: \calD \times \calD \rightarrow \calD$ preserves small colimits in each variable.
\end{itemize}
Then there exists an equivalence of monoidal $\infty$-categories $\calC^{\otimes} \rightarrow \calD^{\otimes}$. Moreover, the collection of such equivalences is parametrized by a contractible Kan complex.
\end{itemize}
\end{corollary}

\begin{proof}
In view of Remark \ref{otherlandur}, we can identify $\widehat{\Cat}_{\infty}^{\Stabb, \Mon}$ with the $\infty$-category of algebra objects of $\LPressStab$. Proposition \ref{gurgle} implies that $\widehat{\Cat}_{\infty}^{\Stabb, \Mon}$ has an initial object $\calC^{\otimes}$. This proves $(1)$. Moreover, Proposition \ref{gurgle} also asserts that the underlying $\infty$-category
$\calC = \calC^{\otimes}_{[1]}$ is equivalent to the unit object of $(\MonStab)_{[1]}$, which
is the $\infty$-category of spectra (Proposition \ref{tupid}); this proves $(2)$. 

Suppose that $\calD^{\otimes}$ is a monoidal $\infty$-category satisfying $(i)$, $(ii)$, and $(iii)$. Since $\calC^{\otimes}$ is an initial object of $\widehat{\Cat}_{\infty}^{\Stabb, \Mon}$, there exists a monoidal functor $f: \calC^{\otimes} \rightarrow \calD^{\otimes}$, unique up to a contractible space of choices. We claim that $f$ is an equivalence of monoidal $\infty$-categories. According to Remark \ref{eggsal}, it will suffice to show that $f$ induces an equivalence $f_{[1]}: \calC \rightarrow \calD = \calD^{\otimes}_{[1]}$ of ordinary underlying categories. Corollary \stableref{choccrok} implies that $f_{[1]}$ is determined, up to equivalence, by the image of the sphere spectrum
$\Sphere \in \calC$. Since $\Sphere$ is the unit object of $\calC$, $f_{[1]}$ carries $\Sphere$ to the unit object of $\calD$. Condition $(ii)$ now implies that $f_{[1]}$ is an equivalence, as desired.
\end{proof}

It follows from Corollary \ref{surcoi} that the $\infty$-category $\Spectra$ admits an essentially unique monoidal structure, which may be characterized by the following properties:

\begin{itemize}
\item[$(i)$] The bifunctor $\otimes: \Spectra \times \Spectra \rightarrow \Spectra$ preserves small colimits separately in each variable.
\item[$(ii)$] The unit object of $\Spectra$ is the sphere spectrum $\Sphere$.
\end{itemize}

We will refer to this monoidal structure on $\Spectra$ as the {\it smash product monoidal structure}.
From the uniqueness, we conclude that the smash product monoidal structure on $\Spectra$ is equivalent to any monoidal structure obtained via the constructions described in $(S1)$ or $(S2)$ above.\index{smash product!monoidal structure}\index{monoidal $\infty$-category!of spectra}

\begin{proposition}\label{discur}
The forgetful functor
$\Mod_{\Spectra}( \LPress ) \rightarrow \LPress$ is fully faithful, and induces an equivalence
$\Mod_{\Spectra}( \LPress ) \rightarrow \LPressStab$.
\end{proposition}

In other words, every presentable stable $\infty$-category is tensored over $\Spectra$ in an essentially unique way.

\begin{proof}
Corollaries \ref{puterry} and \ref{surcoi} imply that the forgetful functor
$\Mod_{ \Spectra} (\LPressStab) \rightarrow \LPressStab$ is an equivalence of $\infty$-categories.
To complete the proof, we must show that the inclusion $\Mod_{\Spectra}( \LPressStab)
\subseteq \Mod_{\Spectra}(\LPress)$ is also an equivalence of $\infty$-categories. In other words, we must show that if a presentable $\infty$-category $\calC$ is tensored over $\Spectra$ (in a colimit-preserving fashion), then $\calC$ is stable.

We first show that $\calC$ is pointed. Let $\otimes: \Spectra \times \calC \rightarrow \calC$ be the action of $\Spectra$ on $\calC$. Let $\Sphere \in \Spectra$ denote the sphere spectrum, let $0 \in \Spectra$ be a zero object, and let $\alpha: \Sphere \rightarrow 0$ be a morphism. Let
$1_{\calC} \in \calC$ be a final object. Then we obtain an induced map
$1_{\calC} \simeq \Sphere \otimes 1_{\calC} \rightarrow 0 \otimes 1_{\calC}.$
Since the functor $\bigdot \otimes C$ preserves colimits, the tensor product $0 \otimes 1_{\calC}$ is an initial object of $\calC$. We now apply Remark \stableref{tustpoint} to conclude that $\calC$ is pointed.

We have a pushout diagram
$$ \xymatrix{ \Sphere \ar[r] \ar[d] & 0 \ar[d] \\
0 \ar[r] & \Sphere[1]. }$$
We therefore obtain an associated pushout diagram
$$ \xymatrix{ \Sphere \otimes \bigdot \ar[r] \ar[d] & 0 \otimes \bigdot \ar[d] \\
0 \otimes \bigdot \ar[r] & \Sphere[1] \otimes \bigdot }$$
in the $\infty$-category $\Fun(\calC, \calC)$. It follows that tensor product with $\Sphere[1]$ can be identified with the suspension functor $\Sigma: \calC \rightarrow \calC$. Consequently, the suspension $\Sigma: \calC \rightarrow \calC$ is an equivalence: a homotopy inverse is given by tensor product with $\Sphere[-1] \in \Spectra$. Corollary \stableref{charstut} implies that $\calC$ is stable, as desired.
\end{proof}

\begin{remark}\label{smallstuck}
Let $\calC$ be a small stable $\infty$-category. Then the $\infty$-category
$\Ind(\calC)$ is a presentable stable $\infty$-category, so we may regard
$\Ind(\calC)$ as tensored over $\Spectra$ in an essentially unique way such that the tensor product functor
$$ \otimes: \Spectra \times \Ind(\calC) \rightarrow \Ind(\calC)$$
preserves colimits separately in each variable. Let $\calC' \subseteq \Ind(\calC)$ denote the essential image of the Yoneda embedding $j: \calC \rightarrow \Ind(\calC)$, so that $j$ induces an equivalence of $\calC$ with $\calC'$. Let $\FinSpectra$ denote the full subcategory of $\Spectra$ spanned by
the {\em finite} spectra (see \S \stableref{stable8}). We observe that $\FinSpectra$ is generated under finite colimits by the collection of $n$-spheres $\{ \Sphere[n] \}_{n \in \Z}$. Moreover, the operation
$$\Sphere[n] \otimes \bigdot: \Ind(\calC) \rightarrow \Ind(\calC)$$
can be identified with the shift functor $C \mapsto C[n]$. Since $\calC'$ is a stable subcategory
of $\Ind(\calC)$ (Proposition \stableref{kappstable}), we conclude that the tensor product $\otimes$
induces a functor
$$ \FinSpectra \otimes \calC' \rightarrow \calC'.$$
We observe that the full subcategory $\FinSpectra \subseteq \Spectra$ is closed under tensor products and contains the unit object $\Sphere \in \Spectra$, and therefore inherits a monoidal structure
(Proposition \ref{yukyik}). Using the same argument, one can show that $\calC'$
inherits the structure of an $\infty$-category tensored over $\FinSpectra$. Identifying $\calC$ with
$\calC'$ via the Yoneda embedding $j$, we conclude that every (small) stable $\infty$-category
$\calC$ can be regarded as tensored over the $\infty$-category of finite spectra, in such a way that the tensor product
$$ \otimes: \FinSpectra \times \calC \rightarrow \calC$$
preserves finite colimits in each variable. With a bit more effort, one can prove an analogue of Proposition \ref{discur} in this context. 
\end{remark}

\begin{remark}\label{bigstuck}
Proposition \ref{discur} and Remark \ref{smallstuck} have analogues in the setting of pointed $\infty$-categories:
\begin{itemize}
\item[$(1)$] The $\infty$-category $\SSet_{\ast}$ of pointed spaces admits an (essentially unique)
monoidal structure such that the tensor product $\otimes: \SSet_{\ast} \times \SSet_{\ast} \rightarrow \SSet_{\ast}$ preserves colimits separately in each variable.
\item[$(2)$] The forgetful functor $\Mod_{ \SSet_{\ast}}( \LPress) \rightarrow \LPress$ is fully faithful, and its essential image consists of pointed presentable $\infty$-categories.
\item[$(3)$] Let $\calC$ be a pointed $\infty$-category which admits finite colimits, and let
$\FinSpace \subseteq \SSet_{\ast}$ denote the full subcategory of $\SSet_{\ast}$ spanned by the finite pointed spaces (see \S \stableref{stable8}). Then $\FinSpace$ is a monoidal subcategory
of $\SSet_{\ast}$, and $\calC$ admits the structure of an $\infty$-category tensored over $\FinSpace$ in such a way that the tensor product
$$ \otimes: \FinSpace \times \calC \rightarrow \calC$$
preserves finite colimits in each variable.
\end{itemize}
\end{remark}

We conclude with the following counterpart to Corollary \ref{surcoi}:

\begin{proposition}\label{dadd}
Let $\calC^{\otimes} \rightarrow \Nerve(\cDelta)^{op}$ be a monoidal $\infty$-category which satisfies conditions $(i)$ through $(iii)$ of Corollary \ref{surcoi} $($so that $\calC^{\otimes}$ is equivalent, as a monoidal $\infty$-category, to the smash product monoidal structure on $\Spectra${}$)$. Let $q: \calM^{\otimes} \rightarrow \calC^{\otimes}$ be an $\infty$-category left-tensored over
$\calC^{\otimes}$. Then $q$ is equivalent $($as an object of $\widehat{\CatMod}${}$)$ to the canonical action of $\calC^{\otimes}$ on itself, if and only if the following conditions are satisfied:
\begin{itemize}
\item[$(1)$] The underlying $\infty$-category $\calM = \calM^{\otimes}_{[0]}$ is equivalent to the $\infty$-category $\Spectra$.
\item[$(2)$] The tensor product functor $\otimes: \calC \times \calM \rightarrow \calM$ is colimit-preserving in each variable.
\end{itemize}
\end{proposition}

\begin{proof}
The necessity of $(1)$ and $(2)$ is clear. For the converse, we invoke the equivalence
$ \widehat{ \CatMod} \simeq \Mod( \widehat{\Cat}_{\infty} )$
of Corollary \ref{spg}. In virtue of assumption $(2)$, we can identify both
$q: \calM^{\otimes} \rightarrow \calC^{\otimes}$ and the canonical action of $\calC^{\otimes}$ on itself with objects of $\Mod_{ \calC^{\otimes}}( \LPress )$. Combining Corollary \ref{surcoi} with Proposition \ref{discur}, we deduce that the forgetful functor $\Mod_{\calC^{\otimes}}( \LPress) \rightarrow \LPressStab$ is an equivalence. It therefore suffices to prove that the underlying $\infty$-categories $\calC$ and $\calM$ are equivalent, which follows from $(1)$.
\end{proof}

\begin{warning}
Given a pair of objects $X,Y \in \h{\Spectra}$, the smash product of $X$ and $Y$ is usually denoted by $X \wedge Y$. We will depart from this convention by writing instead $X \otimes Y$ for the smash product. 
\end{warning}

\subsection{Associative Ring Spectra}\label{monoid9}

In this section, we will introduce the theory of {\it $A_{\infty}$-ring spectra}, or, as we will call them, {\it $A_{\infty}$-rings}. Roughly speaking, an $A_{\infty}$-ring is to an ordinary associative ring as a spectrum is to an abelian group, or as a homotopy type is to a set. In particular, if we restrict our attention to discrete $A_{\infty}$-rings, then we recover classical ring theory (Proposition \ref{umberilt}). 

A large portion of classical (noncommutative) algebra can be generalized to the setting of $A_{\infty}$-rings; we will see some examples in the next few sections. For the time being, we will concern ourselves only with the definition and basic formal properties.

\begin{definition}\index{$A_{\infty}$-ring}
An {\it $A_{\infty}$-ring} is an algebra object of the $\infty$-category of $\Spectra$ (endowed with its smash product monoidal structure). We let $\AInfty$ denote the $\infty$-category
$\Alg(\Spectra)$ of $A_{\infty}$-rings.\index{ZZZAInfty@$\AInfty$}
\end{definition}

Let $R$ be an $A_{\infty}$-ring. We will generally not distinguish notationally between $R$ and its underlying spectrum. In particular, for each $n \in \Z$, we let $\pi_n R$ denote the $n$th homotopy group of the underlying spectrum. We observe that $\pi_n R$ can be identified with the
set $\pi_0 \bHom_{\Spectra}( \Sphere[n], R)$, where $\Sphere$ denotes the sphere spectrum.
Since $\Sphere$ is the identity for the smash product, there is a canonical equivalence
$\Sphere \otimes \Sphere \simeq \Sphere$; using the fact that $\otimes$ is exact in each variable, we deduce the existence of equivalences $\Sphere[n] \otimes \Sphere[m] \simeq \Sphere[n+m]$ for all $n,m \in \Z$. The composition
$$ \bHom_{\Spectra}( \Sphere[n], R) \times \bHom_{\Spectra}( \Sphere[m], R)
\rightarrow \bHom_{\Spectra}( \Sphere[n] \otimes \Sphere[m], R \otimes R) \rightarrow
\bHom_{\Spectra}( \Sphere[n+m], R)$$
determines a map of abelian groups $\pi_n R \otimes \pi_m R \rightarrow \pi_{n+m} R$.
It is not difficult to see that these maps endow 
$\pi_{\ast} R = \bigoplus_{n} \pi_n R$
with the structure of a graded associative ring, which depends functorially on $R$.
In particular, $\pi_0 R$ is an ordinary associative ring, and each $\pi_{n} R$ has the structure of a $\pi_0 R$-bimodule.

Let $\Omega^{\infty}: \Spectra \rightarrow \SSet$ be the $0$th space functor.
If $R$ is an $A_{\infty}$-ring, we will refer to $X=\Omega^{\infty} R$ as the {\it underlying space} of $R$. The underlying space $X$ is equipped with an addition $X \times X \rightarrow X$ (determined by the fact that it is the $0$th space of a spectrum) and a multiplication $X \times X \rightarrow X$ (determined by the map $R \otimes R \rightarrow R$); these maps endow $X$ with the structure of a ring object in the homotopy category $\calH$ of spaces. However, the structure of an $A_{\infty}$-ring is much richer: not only do the ring axioms on $X$ hold up to homotopy, they hold up to {\em coherent} homotopy. 

The functor $\Omega^{\infty}: \AInfty \rightarrow \SSet$ is not conservative: a map of $A_{\infty}$-rings $f: A \rightarrow B$ which induces a homotopy equivalence of underlying spaces need not be an equivalence in $\AInfty$. We observe that $f$ is an equivalence of $A_{\infty}$-rings if and only if it is an equivalence of spectra; that is, if and only if $\pi_n(f): \pi_n A \rightarrow \pi_n B$ is an isomorphism of abelian groups for all $n \in \Z$. However, $\Omega^{\infty}(f)$ is a homotopy equivalence of spaces provided only that $\pi_n(f)$ is an isomorphism for $n \geq 0$; this is generally a weaker condition.

\begin{example}\label{hyperex}
Let $\calC$ be a stable $\infty$-category, and let $X \in \calC$ be an object.
Then it is possible to extract from $\calC$ an $A_{\infty}$-ring spectrum
$\End_{\calC}(X)$ with the property that $\pi_{n} \End_{\calC}(X) \simeq \Ext^{-n}_{\calC}(X,X)$ for all $n \in \Z$, and the ring structure on $\pi_{\ast} \End_{\calC}(X)$ is given by composition
in the triangulated category $\h{\calC}$. We will describe the argument in the case where $\calC$ is presentable. According to Proposition \ref{discur}, the $\infty$-category $\calC$ is naturally left-tensored over $\Spectra$. Proposition \ref{enterich} implies that $\calC$ is also enriched over
$\Spectra$, so that there exists a morphism object $\Mor_{\calC}(X,X)$. 
Proposition \ref{ugher} implies that $\End_{\calC}(X) = \Mor_{\calC}(X,X)$ can be lifted to a final object of the monoidal $\infty$-category $\Spectra[X]$, and therefore inherits the structure of an algebra object. The identification of the homotopy groups of $\End_{\calC}(X)$ follows from
the homotopy equivalence
$\bHom_{ \Spectra}( \Sphere[n], \End_{\calC}(X) ) \simeq \bHom_{\calC}( \Sphere[n] \otimes X, X).$
\end{example}

\begin{remark}
Combining the uniqueness assertion of Corollary \ref{surcoi} with Theorem \ref{beckify}, we conclude that $\AInfty$ is equivalent to the underlying $\infty$-category of strictly associative monoids in any sufficiently nice monoidal model category of spectra (see Example \ref{bulwork}). With minor modifications, the same argument can be applied to the model of spectra described in \cite{EKMM} (though Theorem \ref{beckify} does not quite apply in its present form). 
\end{remark}

\begin{remark}
Combining the uniqueness assertion of Corollary \ref{surcoi} with the construction $(S2)$ of \S \ref{monoid6.0}, we conclude that $\AInfty$ can be identified with the $\infty$-category of
colimit-preserving monads on $\Spectra$.
\end{remark}

Recall that a spectrum $X$ is said to be {\it connective} if $\pi_{n} X \simeq 0$ for $n < 0$.\index{connective}

\begin{lemma}\label{smashstab}
The t-structure on the $\infty$-category $\Spectra$ determined by the class of connective objects is compatible with the smash product monoidal structure $($in the sense of Definition \ref{tuppa}$)$. In other words, the full subcategory $\connSpectra \subseteq \Spectra$ spanned by the connective objects is closed under smash products and contains the unit object. Consequently, the monoidal structure on $\Spectra$ determines a monoidal structure on $\connSpectra$. 
\end{lemma}

\begin{proof}
It follows from the results of \S \stableref{stable16} that $\connSpectra$ is the smallest full subcategory of $\Spectra$ which contains the sphere spectrum $\Sphere \in \Spectra$ and is stable under colimits and extensions. 
Let $\calC$ be the full subcategory of $\Spectra$ spanned by those spectra $X$ such that, 
for all $Y \in \connSpectra$, $X \otimes Y$ is connective. We wish to prove that $\connSpectra \subseteq \calC$. Since the smash product preserves colimits separately in each variable, we conclude that $\calC$ is closed under colimits and extensions in $\Spectra$. It will therefore suffice to prove that $\Sphere \in \calC$. This is clear, since $\Sphere$ is the unit object of $\Spectra$.
\end{proof}

We will say that an $A_{\infty}$-ring $R$ is {\it connective} if its underlying spectrum is connective.\index{connective!$A_{\infty}$-ring}\index{$A_{\infty}$-ring!connective}
We let $\AInfty^{\conn}$ denote the full subcategory of $\AInfty$ spanned by the connective objects. 
Equivalently, we may view $\AInfty^{\conn}$ as the $\infty$-category $\Alg( \connSpectra )$ of algebra objects in connective spectra.\index{ZZZAInftyconn@$\AInfty^{\conn}$}

When restricted to connective $A_{\infty}$-rings, the functor $\Omega^{\infty}$ detects equivalences: if $f: A \rightarrow B$ is a morphism in $\AInfty^{\conn}$ such that $\Omega^{\infty}(f)$ is an equivalence, then $f$ is an equivalence. We observe that
the functor $\Omega^{\infty}: \AInfty^{\conn} \rightarrow \SSet$ is a composition of a pair of functors
$ \Alg( \connSpectra ) \rightarrow \connSpectra \rightarrow \SSet,$
both of which preserve geometric realizations (Corollary \ref{filtfem} and \stableref{denkmal}) and admit left adjoints. It follows from Theorem \ref{barbeq} that $\AInfty^{\conn}$ can be identified with the $\infty$-category of modules over a suitable monad on $\SSet$. In other words, we can view connective $A_{\infty}$-rings as spaces equipped with some additional structures. Roughly speaking, these additional structures consist of an addition and multiplication which satisfy the ring axioms, up to coherent homotopy.\index{ZZZconnSpectra@$\connSpectra$}


\begin{definition}\label{cancon}\index{connective!cover}
Let $R$ be an $A_{\infty}$-ring. A {\it connective cover} of $R$ is a morphism
$\phi: R' \rightarrow R$ of $A_{\infty}$-rings with the following properties:
\begin{itemize}
\item[$(1)$] The $A_{\infty}$-ring $R'$ is connective.
\item[$(2)$] For every connective $A_{\infty}$-ring $R''$, composition with $\phi$ induces a homotopy equivalence
$$ \bHom_{ \AInfty}(R'', R') \rightarrow \bHom_{\AInfty}(R'',R). $$
\end{itemize}
\end{definition}

\begin{remark}
In the situation of Definition \ref{cancon}, we will generally abuse terminology and simply refer to $R'$ as a connective cover of $R$, in the case where the map $\phi$ is implicitly understood.
\end{remark}

\begin{proposition}\label{canconex}
\begin{itemize}
\item[$(1)$] Every $A_{\infty}$-ring $R$ admits a connective cover.
\item[$(2)$] An arbitrary map $\phi: R' \rightarrow R$ of $A_{\infty}$-rings is a connective cover of $R$ if and only if $R'$ is connective, and the induced map $\pi_{n} R' \rightarrow \pi_{n} R$ is an isomorphism for $n \geq 0$.
\item[$(3)$] The inclusion $\AInfty^{\conn} \subseteq \AInfty$ admits a right adjoint $G$, which carries each $A_{\infty}$-ring $R$ to a connective cover $R'$ of $R$.
\end{itemize}
\end{proposition}

\begin{proof}
Combine Lemma \ref{smashstab}, Proposition \ref{yukyik} and Remark \ref{yuka}.
\end{proof}

Recall that an object $X$ of an $\infty$-category $\calC$ is said to be {\it $n$-truncated} if the
mapping spaces $\bHom_{\calC}(Y,X)$ are $n$-truncated, for every $Y \in \calC$ (see \S \toposref{truncintro}). Proposition \ref{algprec} implies that the $\infty$-categories $\AInfty$ and $\AInfty^{\conn}$ are both presentable, so that we have a good theory of truncation functors.

\begin{proposition}\label{equitrunc}
Let $R$ be a connective $A_{\infty}$-ring and let $n$ be a nonnegative integer. The following conditions are equivalent:
\begin{itemize}
\item[$(1)$] As an object of $\AInfty^{\conn}$, $R$ is $n$-truncated. 
\item[$(2)$] As an object of $\connSpectra$, $R$ is $n$-truncated.
\item[$(3)$] The space $\Omega^{\infty}(R)$ is $n$-truncated.
\item[$(4)$] For all $m > n$, the homotopy group $\pi_{m} R$ is trivial.
\end{itemize}
\end{proposition}

\begin{proof}
The equivalence $(3) \Leftrightarrow (4)$ is easy (Remark \toposref{humpter}), and 
the equivalence $(2) \Leftrightarrow (3)$ was explained in Warning \stableref{umper}. 
The implication $(1) \Rightarrow (2)$ follows from Proposition \toposref{eaa}, since
the forgetful functor $\AInfty^{\conn} \rightarrow \connSpectra$ preserves small limits (Corollary \ref{slimycomp}). 

We now prove that $(2) \Rightarrow (1)$. Assume that $R$ is $n$-truncated as a spectrum.
Let $T: ( \AInfty^{\conn} )^{op} \rightarrow \SSet$ be the functor represented by $R$. Let
$\calC \subseteq \AInfty^{\conn}$ be the full subcategory of $\AInfty^{\conn}$ spanned by those objects $B$ such that $T(B)$ is $n$-truncated. We wish to prove that $\calC = \AInfty^{\conn}$. Since $T$ preserves limits (Proposition \toposref{yonedaprop}) and the class of $n$-truncated spaces is stable under limits (Proposition \toposref{altum}), we conclude that
$\calC$ is stable under small colimits in $\AInfty^{\conn}$. Let $F$ be a left adjoint to the forgetful functor $\AInfty^{\conn} \rightarrow \connSpectra$ (Theorem \ref{hutmunn}). Proposition \ref{littlebeck} implies that $\calC$ is generated under colimits by the essential image of $F$. Consequently, it will suffice to show that $F(M) \in \calC$, for every $M \in \connSpectra$. Equivalently, we must show that the space
$\bHom_{ \AInfty^{\conn} }( F(M), R) \simeq \bHom_{ \connSpectra}( M, R)$
is $n$-truncated, which follows from $(2)$.
\end{proof}

\begin{remark}
An $A_{\infty}$-ring $R$ is $n$-truncated as an object of $\AInfty$ if and only if it is equivalent to zero (since the $\infty$-category $\Spectra$ has no nontrivial $n$-truncated objects).
\end{remark}

Let $\tau_{\leq n}: \connSpectra \rightarrow \connSpectra$ be the truncation functor on connective spectra, and let $\tau^{\Alg}_{\leq n}: \AInfty^{\conn} \rightarrow \AInfty^{\conn}$ be
the truncation functor on connective $A_{\infty}$-rings. Since the forgetful functor
$\theta: \AInfty^{\conn} \rightarrow \connSpectra$ preserves $n$-truncated objects, there is a canonical natural transformation
$\alpha: \tau_{\leq n} \circ \theta \rightarrow \theta \circ \tau_{\leq n}^{\Alg}.$
Our next goal is to show that $\alpha$ is an equivalence.

\begin{proposition}\label{superjump}
Let $\connSpectra$ be the $\infty$-category of connective spectra, endowed with
the smash product monoidal structure, and let $(\connSpectra)_{\leq n}$ be the $\infty$-category of $n$-truncated objects of $\connSpectra$. Then:
\begin{itemize}
\item[$(1)$] The localization functor $\tau_{\leq n}: \connSpectra \rightarrow \connSpectra$ is compatible with the smash product monoidal structure, in the sense of Definition \ref{compatmon}.
\item[$(2)$] The smash product monoidal structure on $\Spectra_{\geq 0}$ induces a monoidal structure on $(\connSpectra)_{\leq n}$ and an identification
$\Alg( (\connSpectra)_{\leq n} ) \simeq \tau^{\Alg}_{\leq n} \AInfty^{\conn}$.
\end{itemize}
\end{proposition}

\begin{proof}
Assertion $(1)$ is a special case of Proposition \ref{jumperr}. Assertion $(2)$ follows from $(1)$ and Proposition \ref{equitrunc}.
\end{proof}

In other words, if $R$ is a connective $A_{\infty}$-ring, then the ring structure on $R$ determines
a ring structure on $\tau_{\leq n} R$ for all $n \geq 0$. 

\begin{definition}\index{discrete!$A_{\infty}$-ring}\index{$A_{\infty}$-ring!discrete}
\index{ZZZAInftydisc@$\AInfty^{\disc}$}
We will say that an $A_{\infty}$-ring is {\it discrete} if it is connective and $0$-truncated.
We let $\AInfty^{\disc}$ denote the full subcategory of $\AInfty$ spanned by the discrete objects.
\end{definition}

Since the mapping spaces in $\AInfty^{\disc}$ are $0$-truncated, it follows that
$\AInfty^{\disc}$ is equivalent to the nerve of an ordinary category. We conclude this section by identifying the relevant category.

\begin{proposition}\label{umberilt}
The functor $R \mapsto \pi_0 R$ determines an equivalence from
$\AInfty^{\disc}$ to the $($nerve of the$)$ category of associative rings.
\end{proposition}

\begin{proof}
According to Proposition \ref{superjump}, we can identify $\AInfty^{\disc}$ with the
$\infty$-category of algebra objects of the heart $\Spectralheart$, which inherits a monoidal structure from $\Spectra$ in view of Proposition \ref{jumperr} and Lemma \ref{smashstab}. Proposition \stableref{specster} allows us to identify $\Spectralheart$ with (the nerve of) the category of abelian groups. Moreover, the induced monoidal structure on $\Spectralheart$ has $\pi_0 \Sphere \simeq \Z$ as unit object, and the tensor product functor $\otimes$ preserves colimits separately in each variable. It follows that this monoidal structure coincides (up to canonical equivalence) with the usual monoidal structure on $\Spectralheart$, given by tensor products of abelian groups. Consequently, we may identify $\Alg( \Spectralheart )$ with the (nerve of the) category of
associative rings.
\end{proof}


\subsection{Modules over $A_{\infty}$-Rings}\label{monoid11}

Let $R$ be an $A_{\infty}$-ring. In this section, we will describe the associated theory of {\it $R$-module spectra}, or simply {\it $R$-modules}. This can be regarded as a generalization of homological algebra: if $R$ is an ordinary ring (regarded as a discrete $A_{\infty}$-ring via Proposition \ref{umberilt}), then the homotopy category of $R$-module spectra coincides with the classical {\em derived category} of $R$ (Proposition \ref{derivdisc}); in particular, the theory of $R$-module spectra is a generalization of the usual theory of $R$-modules. 

For any $A_{\infty}$-ring $R$, the $\infty$-category $\Mod_{R}$ of $R$-modules is stable (Proposition \ref{stabult}), so its homotopy category $\h{ \Mod_{R} }$ is triangulated. 
According to a result of Schwede and Shipley (Theorem \ref{schwedeshipley}), a large variety of stable $\infty$-categories have the form $\Mod_{R}$, for appropriately chosen $R$. In particular, $\h{ \Mod_{R} }$ is generally not equivalent to the derived category of an abelian category.

\begin{definition}
Let $R$ be an $A_{\infty}$-ring. We let $\Mod_{R}$ denote the $\infty$-category
$\Mod_{R}(\Spectra)$; we will refer to $\Mod_{R}$ as {\it the $\infty$-category of $($left$)$ $R$-modules}.
\end{definition}

We will generally not distinguish between an $R$-module $M$ and the underlying spectrum.
In particular, the homotopy groups $\pi_{n} M$ are defined to be the homotopy groups of the underlying spectrum. The action map $R \otimes M \rightarrow M$ induces
bilinear maps $\pi_{n} R \times \pi_{m} M \rightarrow \pi_{n+m} M$, which endow $\pi_{\ast} M$ with the structure of a graded left module over $\pi_{\ast} R$. We will say that $M$ is {\it connective} if
its underlying spectrum is connective; that is, if $\pi_{n} M \simeq 0$ for $n < 0$.

\begin{remark}
Roughly speaking, if we think of $R$ as a space equipped with the structure of an associative ring up to coherent homotopy, then an $R$-module can be thought of as another space which has an addition and an action of $R$, up to coherent homotopy in the same sense. This intuition is really only appropriate in the case where $R$ and $M$ are connective, since the homotopy groups in negative degree have no ready interpretation in terms of underlying spaces.
\end{remark}

Our first goal is to prove that the $\infty$-category of modules over an $A_{\infty}$-ring is stable. This is a consequence of the following more general assertion:

\begin{proposition}\label{stabult}
Let $\calC$ be an $\infty$-category equipped with a monoidal structure and a left action on an $\infty$-category $\calM$. Assume that $\calM$ is a stable $\infty$-category, and let
$R \in \Alg(\calC)$ be such that the functor $M \mapsto R \otimes M$ is exact.
Then $\Mod_{R}(\calM)$ is a stable $\infty$-category. Moreover, if $\calN$ is an arbitrary stable $\infty$-category, then a functor $\calN \rightarrow \Mod_{R}(\calM)$ is exact if and only if the composite functor 
$\calN \rightarrow \Mod_{R}(\calM) \rightarrow \calM$ is exact. In particular, the forgetful functor $\Mod_{R}(\calM) \rightarrow \calM$ is exact.
\end{proposition}

\begin{proof}
This is a special case of Proposition \ref{poststorkus}.
\end{proof}

If $R$ is a connective $A_{\infty}$-ring, the homotopy groups of an $R$-module $M$ can be interpreted in terms of an appropriate t-structure on $\Mod_{R}$.

\begin{notation}\index{ZZZModRgeq0@$\Mod_{R}^{\geq 0}$}\index{ZZZMoRleq0@$\Mod_{R}^{\leq 0}$}
If $R$ is an $A_{\infty}$-ring, we let $\Mod_{R}^{\geq 0}$ be the full subcategory of
$\Mod_{R}$ spanned by those $R$-modules $M$ for which $\pi_{n} M \simeq 0$ for $n < 0$, and $\Mod_{R}^{\leq 0}$ the full subcategory of $\Mod_{R}$ spanned by those $R$-modules $M$ for which $\pi_{n} M \simeq 0$ for $n > 0$.
\end{notation}

\begin{notation}
Let $R$ be an $A_{\infty}$-ring, and let $M$ and $N$ be {\em left} $R$-modules. We let
$\Ext^i_{R}(M,N)$ denote the abelian group $\pi_{0} \bHom_{ \Mod_{R} }(M, N[i])$.
\end{notation}

\begin{proposition}\label{tmod}
Let $R$ be a connective $A_{\infty}$-ring. Then:
\begin{itemize}
\item[$(1)$] The full subcategory $\Mod_{R}^{\geq 0} \subseteq \Mod_{R}$ is the smallest
full subcategory which contains $R$ $($regarded as an $R$-module in the natural way; see Example \ref{algitself}$)$ and is stable under small colimits.

\item[$(2)$] The subcategories $\Mod_{R}^{\geq 0}, \Mod_{R}^{\leq 0}$ determine
an accessible t-structure on $\Mod_{R}$ $($see \S \stableref{stable16}$)$.

\item[$(3)$] The t-structure described in $(2)$ is both left and right complete, and the
functor $\pi_0$ determines an equivalence of the heart $\Mod^{\heartsuit}_{R}$ with the $($nerve of the$)$ ordinary category of $\pi_0 R$-modules.

\item[$(4)$] The subcategories $\Mod_{R}^{\geq 0}, \Mod_{R}^{\leq 0} \subseteq \Mod_{R}$ are stable under small products and small filtered colimits.
\end{itemize}
\end{proposition}

\begin{proof}
According to Proposition \stableref{kura}, there exists an accessible t-structure $(\Mod_{R}', \Mod_{R}'')$ with the following properties:
\begin{itemize}
\item[$(a)$] An object $M \in \Mod_{R}$ belongs to $\Mod_{R}''$ if and only if
$\Ext^{i}_{R}(R, M) \simeq 0$ for $i < 0$.
\item[$(b)$] The $\infty$-category $\Mod_{R}'$ is the smallest full subcategory of $\Mod_{R}$ which contains the object $R$ and is stable under extensions and small colimits. 
\end{itemize}
Corollary \ref{tara} implies that $R$ (regarded as an object of $\Mod_{R}$) corepresents the composition $ \Mod_{R} \rightarrow \Spectra \stackrel{\Omega^{\infty}}{\rightarrow} \SSet.$
It follows that $\Mod_{R}'' = \Mod_{R}^{\leq 0}$.
Because the forgetful functor $\Mod_{R} \rightarrow \Spectra$ preserves small colimits (Corollary \ref{gloop}), we conclude that $\Mod_{R}^{\geq 0}$ is stable under extensions and small colimits.
Since $R$ is connective, $R \in \Mod_{R}^{\geq 0}$, so that $\Mod_{R}' \subseteq \Mod_{R}^{\geq 0}$. Let $\calC$ be the smallest full subcategory of $\Mod_{R}$ which contains $R$ and is stable under small colimits, so that $\calC \subseteq \Mod'_{R}$. We will complete the proof of $(1)$ and $(2)$ by showing that $\calC = \Mod_{R}^{\geq 0}$. 

Let $M \in \Mod_{R}^{\geq 0}$. We will construct a diagram
$$ M(0) \rightarrow M(1) \rightarrow M(2) \rightarrow \ldots $$
in $(\Mod_{R})_{/M}$ with the following properties:
\begin{itemize}
\item[$(i)$] Let $i \geq 0$, and let $K(i)$ be a kernel of the map $M(i) \rightarrow M$. 
Then $\pi_j K(i) \simeq 0$ for $j < i$. 
\item[$(ii)$] The $R$-module $M(0)$ is a coproduct of copies of $R$.
\item[$(iii)$] For $i \geq 0$, there is a pushout diagram
$$ \xymatrix{ F[i] \ar[r] \ar[d] & 0 \ar[d] \\
M(i) \ar[r] & M(i+1), }$$
where $F$ is a coproduct of copies of $R$.
\end{itemize}

We begin by choosing $M(0)$ to be any coproduct of copies of $R$ equipped with a map
$M(0) \rightarrow M$ which induces a surjection $\pi_0 M(0) \rightarrow \pi_0 M$; for example, we can take $M(0)$ to be a coproduct of copies of $R$ indexed by $\pi_0 M$.
Let us now suppose that the map $f: M(i) \rightarrow M$ has been constructed, with $K(i) = \ker(f)$ such that $\pi_j K(i) \simeq 0$ for $j < i$. We now choose $F$ to be a coproduct of copies of $R$
and a map $g: F[i] \rightarrow K(i)$ which induces a surjection $\pi_0 F \rightarrow \pi_i K(i)$. 
Let $h$ denote the composite map $F[i] \rightarrow K(i) \rightarrow M(i)$, and let
$M(i+1) = \coker(h)$. The canonical nullhomotopy of $K(i) \rightarrow M(i) \rightarrow M$
induces a factorization
$$ M(i) \rightarrow M(i+1) \stackrel{f'}{\rightarrow} M$$
of $f$. We observe that there is a canonical equivalence $\ker(f') \simeq \coker(g)$, so that
$\pi_j \ker(f') \simeq 0$ for $j \leq i$.

Let $M(\infty)$ be the colimit of the sequence $\{ M(i) \}$, and let $K$ be the kernel of the canonical map $M(\infty) \rightarrow M$. Then $K$ can be identified with a colimit of the sequence $\{ K(i) \}_{i \geq 0}$. Since the formation of homotopy groups is preserves filtered colimits, we conclude that
$\pi_{j} K \simeq \colim \pi_{j} K(i) \simeq 0$. Thus $M(\infty) \simeq M$, so that $M \in \calC$ as desired.

Assertion $(4)$ follows from the corresponding result for $\Spectra$, since the forgetful functor
$\Mod_{R} \rightarrow \Spectra$ preserves all limits and colimits (Corollaries \ref{goop} and \ref{gloop}). Since $\Mod_{R} \rightarrow \Spectra$ is a conservative functor, an $R$-module $M$ is zero if and only if $\pi_{n} M$ is zero for all $n \in \Z$. It follows from Proposition \stableref{cosparrow} that $\Mod_{R}$ is both right and left complete.

Let $F$ be the functor from $\Mod^{\geq 0}_{R}$ to the (nerve of the) ordinary category of
left $\pi_0 R$-modules, given by $M \mapsto \pi_0 M$. It is easy to see that $F$ preserves colimits, and that the restriction of $F$ to $\Mod^{\heartsuit}_{R}$ is an exact functor. We wish to prove that
$F_0 = F| \Mod^{\heartsuit}_{R}$ is an equivalence.
We first show that the restriction of $F_0$ is fully faithful. Fix $N \in \Mod^{\heartsuit}_{R}$, and let $\calD$ be the full subcategory of $\Mod^{\geq 0}_{R}$ spanned by those objects $M$ for which the map
$\pi_0 \bHom_{ \Mod_{R} }( M, N) \rightarrow \Hom( F( \tau_{\leq 0} M), F(N) )$
is bijective, where the right hand side indicates the group of $\pi_0 R$-module homomorphims. 
It is easy to see that $\calD$ is stable under colimits and contains $R$. The first part of the proof shows that $\calD = \Mod_{R}^{\geq 0}$. In particular, $F_0$ is fully faithful. 

It remains to show that $F_0$ is essentially surjective. Since $F_0$ is fully faithful and exact, the essential image of $F_0$ is closed under the formation of cokernels. It will therefore suffice to show that every free left $\pi_0 R$-module belongs to the essential image of $F_0$. Since $F_0$ preserves coproducts, it will suffice to show that $\pi_0 R$ itself belongs to the essential image of $F_0$. We now conclude by observing that $F_0( \tau_{\leq 0} R) \simeq \pi_0 R$.
\end{proof}

Let $R$ be a connective $A_{\infty}$-ring, let $\calA$ be the abelian category
of modules over the (ordinary) ring $\pi_0 R$. Then $\calA$ has enough projective objects, so
we can consider the derived $\infty$-category $\calD^{\dplus}(\calA)$ as described in \S \stableref{stable10}. 
Part $(3)$ of Proposition \ref{tmod} determines an equivalence $\Nerve(\calA) \simeq
\Mod_{R}^{\heartsuit}$. Applying Corollary \stableref{truceborn}, we deduce the existence of an (essentially unique) right t-exact functor $\theta: \calD^{\dplus}( \calA) \rightarrow \Mod_{R}$. 

\begin{proposition}\label{derivdisc}\index{discrete!$A_{\infty}$-ring}
\index{$A_{\infty}$-ring!discrete}\index{derived category}
Let $R$ be a connective $A_{\infty}$-ring, and let 
$\theta: \calD^{\dplus}(\calA) \rightarrow \Mod_{R}$ be as described above.
The following conditions are equivalent:
\begin{itemize}
\item[$(1)$] The $A_{\infty}$-ring $R$ is discrete.
\item[$(2)$] The functor $\theta$ is fully faithful, and induces an equivalence
of $\calD^{\dplus}(\calA)$ with the $\infty$-category of right bounded objects of $\Mod_{R}$.
\end{itemize}
\end{proposition}

\begin{proof}
Let $P \in \calA$ be projective object corresponding to the free $R$-module on one generator.
Then, for $M \in \calD^{\dplus}(\calA)$, we have a canonical isomorphism
$ \Ext^0_{ \calD^{\dplus}(\calA) }( P, M ) \simeq \pi_0 M.$
If $(2)$ is satisfied, then we deduce the existence of a canonical isomorphisms
$$\Ext^0_{R}( \theta(P), M) \simeq \pi_0 M \simeq \Ext^0_{R}(R,M)$$
for $M \in \Mod_{R}^{\geq 0}$. Thus $\theta(P)$ and $R$ are isomorphic in the homotopy category $\h{\Mod_{R}}$. Since $\theta(P)$ is discrete, we conclude that $R$ is discrete, which proves $(1)$.

For the converse, let us suppose that $R$ is discrete. Let us regard (the nerve of) $\calA$ as a full subcategory of both $\calD^{\dplus}(\calA)$ and $\Mod_{R}$. For $M,N \in \calA$, let
$\Ext^{i}_{\calA}(M,N)$ denote the abelian group $\pi_0 \bHom_{ \calD^{\dplus}(\calA)}(M,N)$ (in other words, $\Ext^i_{\calA}(M,N)$ is the classical Yoneda $\Ext$-group computed in the abelian category $\calA$). We claim that 
canonical map $\Ext^{i}_{\calA}(M,N) \rightarrow \Ext^i_{R}(M, N)$
is an isomorphism. For $i < 0$, both sides vanish. The proof in general goes by induction on $i$, the case $i = 0$ being trivial. For $i > 0$, we choose an exact sequence
$$ 0 \rightarrow K \rightarrow P \rightarrow M \rightarrow 0$$ in $\calA$, where $P$ is a free $\pi_0 R$-module. We have a commutative diagram of abelian groups with exact rows
$$ \xymatrix{ 
\Ext^{i-1}_{A}(P, N) \ar[r] \ar[d]^{\psi_1} & \Ext^{i-1}_{\calA}(K,N) \ar[r] \ar[d]^{\psi_2} &
\Ext^{i}_{\calA}( M, N) \ar[r] \ar[d]^{\psi_3} & \Ext^i_{\calA}( P, N) \ar[d] \\
\Ext^{i-1}_{R}(P, N) \ar[r] & \Ext^{i-1}_{R}(K,N) \ar[r] &
\Ext^{i}_{R}( M, N) \ar[r] & \Ext^i_{R}( P, N). }$$
We wish to show that $\psi_3$ is an isomorphism. Since $\psi_1$ and $\psi_2$ are bijective by the inductive hypothesis, it will suffice to show that
$\Ext^i_{\calA}(P,N) \simeq 0 \simeq \Ext^i_{R}(P,N).$
The first equivalence follows from the fact that $P$ is a projective object of $\calA$. For the second,
we observe that as an object of $\Mod_{R}$, $P$ coincides with a coproduct of copies of $R$ (in virtue of assumption $(1)$). Consequently, $\Ext^i_{R}(P, N)$ can be identified with a product of copies of $\pi_{-i} N$, which vanishes since $i > 0$ and $N \in \Mod_{R}^{\geq 0}$.

Now suppose that $M \in \calA$, and consider the full subcategory $\calC \subseteq \calD^{\dplus}(\calA)$ spanned by those objects $N$ for which the canonical map
$ \Ext^{i}_{\calD^{\dplus}(\calA)}(M,N) \rightarrow \Ext^i_{R}( \theta(M), \theta(N) )$
is an isomorphism for all $i \in \Z$. Applying the five lemma to the relevant long exact sequences, we conclude that $\calC$ is stable under extensions in $\calD^{\dplus}(\calA)$. The above argument shows that $\calC$ contains the heart of $\calD^{\dplus}(\calA)$; it therefore contains the full subcategory $\calD^{b}(\calA)$ of bounded object of $\calD^{\dplus}(\calA)$.

Now let $\calC' \subseteq \calD^{\dplus}(\calA)$ spanned by those objects $M$ having the property that for {\em every} $N \in \calD^{b}(\calA)$, the canonical map 
$ \Ext^{i}_{\calD^{\dplus}(\calA)}(M,N) \rightarrow \Ext^i_{R}( \theta(M), \theta(N) )$
is an isomorphism for $i \in \Z$. Repeating the above argument, we conclude that
$\calD^{b}(\calA) \subseteq \calC'$. In particular, the restriction $\theta | \calD^{b}(\calA)$ is fully faithful.

We claim that the essential image of $\theta | \calD^{b}(\calA)$ consists of precisely the bounded objects of $\Mod_{R}$. Let $M \in \Mod_{R}$ be a bounded object. We wish to prove that
$M$ belongs to the essential image of $\theta$. Without loss of generality, we may suppose that
$M \in \Mod_{R}^{\geq 0}$. Since $M$ is bounded, we have also $M \in \Mod_{R}^{\leq n}$ for some $n \geq 0$. We now work by induction on $n$. If $n= 0$, then $M$ belongs to the heart of 
$\Mod_{R}$ and the result is obvious. If $n > 0$, then we have a distinguished triangle
$$ \tau_{\geq n} M \rightarrow M \rightarrow \tau_{\leq n-1} M \rightarrow \tau_{\geq n} M[1].$$
Since $\theta$ is exact and fully faithful, it will suffice to show that $\tau_{\geq n} M[-n]$ and
$\tau_{\leq n-1} M$ belong to the essential image of $\theta$, which follows from the inductive hypothesis.

The preceding argument shows that $\theta$ induces an equivalence
$ \calD^{b}(\calA) \rightarrow \Mod_{R}^{b}$ between the full subcategories of bounded objects. 
We now conclude by observing that both $\calD^{\dplus}(\calA)$ and $\Mod_{R}$ are left-complete.
\end{proof}

\begin{remark}\label{ultra}
Let $R$ and $\calA$ be as in Proposition \ref{derivdisc}. The abelian category $\calA$
also has enough injective objects. The proof of Proposition \ref{derivdisc} can be repeated, without essential change, to obtain an identification of $\calD^{\dminus}(\calA)$ with the $\infty$-category
of left bounded objects of $\Mod_{R}$. In other words, we may view $\Mod_{R}$ as a candidate for an unbounded derived $\infty$-category of $R$-modules. Using a completeness argument, it is not difficult to show that this coincides with the usual unbounded derived category of $R$-modules; we will not pursue this point further.
\end{remark}

We conclude this section by addressing the following question: given an $\infty$-category $\calC$, under what conditions does there exist an $A_{\infty}$-ring $R$ and an equivalence
$\calC \rightarrow \Mod_{R}$? Of course, there might be several candidates for $R$. For example, if $R$ is a discrete ring, then the module categories (in the ordinary or derived sense) of
$R$ and $M_2(R)$ are equivalent, where $M_2(R)$ denotes the ring of $2 \times 2$ matrices with coefficients in $R$. To eliminate this ambiguity, we should specify an object $C \in \calC$
which is the hypothetical preimage of $R \in \Mod_{R}$ under the functor $\theta$. In this case, we have the following theorem of Schwede and Shipley:

\begin{theorem}\label{schwedeshipley}[Schwede-Shipley \cite{schwedeshipley}]\index{Schwede-Shipley theorem}
Let $\calC$ be a stable $\infty$-category. Then $\calC$ is equivalent to $\Mod_{R}$, for some $A_{\infty}$-ring $R$, if and only if there exists an object $C \in \calC$ satisfying the following conditions:

\begin{itemize}
\item[$(1)$] The $\infty$-category $\calC$ is presentable.
\item[$(2)$] The object $C \in \calC$ is compact.
\item[$(3)$] The object $C$ generates $\calC$ in the following sense: if $D \in \calC$ has
the property that $\Ext^{n}_{\calC}(C,D)$ for all $n \in \Z$, then $D \simeq 0$.
\end{itemize}
\end{theorem}

\begin{proof}
Suppose first that $\calC \simeq \Mod_{R}$. Then $(1)$ follows from Corollary \ref{underwhear}. To prove $(2)$ and $(3)$, we take $C$ to be $R$ itself, regarded as an $R$-module as explained in Example \ref{algitself}. Corollary \ref{tara} implies that $R$ corepresents the composite functor
$ \Mod_{R} \rightarrow \Spectra \stackrel{\Omega^{\infty}}{\rightarrow} \SSet,$
which preserves filtered colimits in virtue of Corollary \ref{gloop}. This proves $(2)$. If $\Ext^{i}_{\Mod_{R}}( R, D) \simeq \pi_{-i} D$ vanishes for all $i \in \Z$, so that $D \simeq 0$. This proves $(3)$.

Conversely, suppose that $\calC$ satisfies conditions $(1)$, $(2)$, and $(3)$ for an appropriately chosen object $C \in \calC$. Applying Corollary \stableref{choccrok}, we deduce that there is a colimit-preserving functor $F: \Spectra \rightarrow \calC$, determined up to equivalence by the requirement that $F( \Sphere) \simeq C$, where $\Sphere$ denotes the sphere spectrum. Applying the adjoint functor theorem (Corollary \toposref{adjointfunctor}), we deduce that $F$ admits a right adjoint $G$.
We now apply Theoream \ref{partA} to extend the functor $F$ to a strong adjunction diagram
$s \in \AdjDiag(\Spectra, \calC)$. By restrction, $s$ determines a monad $T \in \Alg( \End(\Spectra))$. 

Let $\End_0(\Spectra)$ be the full subcategory of $\End(\Spectra)$ spanned by those functors which preserve colimits. Then $\End_0(\Spectra)$ is closed under composition. We let $\End^{\otimes}_0( \Spectra)$ denote the corresponding full (monoidal) subcategory of $\End^{\otimes}(\Spectra)$, and let 
$$\overline{\End}_{0}^{\otimes}(\Spectra) = \overline{\End}^{\otimes}(\Spectra)
\times_{ \End^{\otimes}(\Spectra) } \End_0^{\otimes}(\Spectra).$$
The composition
$ g_n: \calC \stackrel{G}{\rightarrow} \Spectra \stackrel{\Omega^{\infty-n}_{\ast}}{\rightarrow} 
\SSet_{\ast} \rightarrow \SSet$
can be identified with the composition of a shift functor from $\calC$ to itself (an equivalence of $\infty$-categories) and the functor corepresented by $C \in \calC$. Since $C$ is assumed to be compact, we conclude that $g_n$ preserves filtered colimits. Since the forgetful functor
$\SSet_{\ast} \simeq \SSet_{\ast/} \rightarrow \SSet$ detects filtered colimits, we conclude that
$\Omega^{\infty-n} \circ G$ preserves filtered colimits. Since $\Spectra$ is defined as the homotopy inverse limit of a tower
$$ \ldots \stackrel{\Omega}{\rightarrow} \SSet_{\ast} \stackrel{\Omega}{\rightarrow} \SSet_{\ast}$$
of continuous functors, we conclude that $G$ itself preserves filtered colimits. Since $G$ is exact, it preserves all colimits. It follows that $T$ can be identified with an algebra object of $\End_0(\Spectra)$.

Invoking Remark \ref{supinor} and Corollary \ref{hooku}, we obtain a functor
$\theta: \calC \rightarrow \Mod_{T}( \Spectra)$. Assumption $(3)$ implies that $G$ is conservative. Since $G$ preserves geometric realizations, Corollary \ref{usualbb} implies that $\theta$ is an equivalence. Proposition \ref{dadd} implies that the tensored $\infty$-category
$$ \overline{\End}_0^{\otimes}(\Spectra) \rightarrow \End_{0}^{\otimes}(\Spectra) \rightarrow
\Nerve(\cDelta)^{op}$$
is equivalent (in $\CatMod$) to the standard left action of $\Spectra$ on itself, where $\Spectra$ is endowed with the smash product monoidal structure. It follows that we can identify
$T$ with an $A_{\infty}$-ring $R$, and that there is an equivalence
$\Mod_{T}(\Spectra) \simeq \Mod_{R}(\Spectra)$. Composing this equivalene with $\theta$, we obtain the desired result.
\end{proof}

\begin{remark}
Let $\calC$ be a presentable stable $\infty$-category containing an object $C$.
The $A_{\infty}$-ring $R$ appearing in the proof of Theorem \ref{schwedeshipley} can be identified with the $A_{\infty}$-ring $\End_{\calC}(C)$ constructed in Example \ref{hyperex}.
\end{remark}

\subsection{Balanced Pairings and the Bar Construction}\label{balpair}

Let $A$ be an associative ring. If $M$ is a right $A$-module and $N$ is a left $R$-module spectrum, then we can define an abelian group $M \otimes_{A} N$. The functor
$(M,N) \mapsto M \otimes_{A} N$
preserves colimits in separately in each variable. Moreover, if $F$ is a colimit-preserving functor $F$ from left (right) $A$-modules to abelian groups, then $F(A)$ inherits a right (left) $A$-module structure, and there is a canonical isomorphism $F(N) \simeq F(A) \otimes_{A} N$ ($F(M) \simeq M \otimes_{A} F(A)$). Our goal in this section is to obtain an analogous picture, where we allow
$A$ to be an $A_{\infty}$-ring and $M$ and $N$ to be module spectra over $A$.

We will begin with a much more general situation: a monoidal $\infty$-category $\calC$, an $\infty$-category $\calM$ which is right-tensored over $\calC$, an $\infty$-category $\calN$ which is left-tensored over $\calC$, and a pairing $\calM \times \calN \rightarrow \calD$
which is {\it balanced} in a suitable sense. Given this data, we will define a relative
tensor product functor
$$ \Mod_{A}(\calM) \times \Mod_{A}(\calN) \rightarrow \calD$$
for every algebra object $A \in \Alg(\calC)$ (see Definition \ref{barcon}).
We will proceed to establish some of the basic properties of the relative tensor product; 
In particular, we will prove a ``base change'' or ``push-pull'' formula, which asserts that if 
$f: A \rightarrow B$ is a map of algebras, $M$ is a (right) $B$-module and $N$ is a (left) $A$-module, then there is a canonical equivalence
$ M \otimes_{A} N \simeq M \otimes_{B} (B \otimes_{A} N),$
where $B \otimes_{A} N$ is the $B$-module induced from $N$ (see Proposition \ref{ussr} for a more precise formulation).

The principal example of interest to us are when $\calM$ or $\calN$ (or both) coincide with the original monoidal $\infty$-category $\calC$. In this case, there are canonical examples of
balanced pairings
$$ \calM \times \calC \rightarrow \calM, \quad \calC \times \calN \rightarrow \calN;$$
see Definition \ref{cannapair} (in \S \ref{monoid13}, we will specialize further to the case where
$\calC$ is the $\infty$-category of spectra, with the smash product monoidal structure).

To begin, let us suppose that $\calC$ is a monoidal $\infty$-category, and that $A$ is an algebra object of $\calC$. We will let $\Mod_{A}^{L}(\calC)$ denote the $\infty$-category of {\em left} $A$-modules in $\calC$, as defined in \S \ref{hugr}. As explained in Remark \ref{hugrr}, we have an entirely dual theory of {\em right} modules; we let $\Mod_{A}^{R}(\calC)$ denote the $\infty$-category of right modules over $A$. We wish to define a functor
$\otimes_{A}: \Mod_{A}^{R}(\calC) \times \Mod_{A}^{L}(\calC) \rightarrow \calC$.

Fix objects $M \in \Mod_{A}^{R}(\calC)$, $N \in \Mod_{A}^{L}(\calC)$. We will abuse notation by identifying $M$ and $N$ with the underlying object of the $\infty$-category $\calC$. Using the monoidal structure on $\calC$, we can form the tensor product $M \otimes N$. This should be regarded as a ``first approximation'' to $M \otimes_{A} N$. We will have a natural map $\alpha: M \otimes N \rightarrow M \otimes_{A} N$. We do not expect $\alpha$ to be an equivalence; rather, $M \otimes_{A} N$ should be obtained from $M \otimes N$ by forcing
the two natural actions of $A$ to coincide. To be more precise, consider the diagram
$$\xymatrix{ M \otimes A \otimes N \ar@<.4ex>[r]^{f} \ar@<-.4ex>[r]_{g} & M \otimes N}.$$
where $f$ is given by the right action of $A$ on $M$, and $g$ is given by the left action of $A$ on $N$. We should expect $\alpha \circ f \simeq \alpha \circ g$. In other words, we expect that $\alpha$
should factor through the coequalizer of the maps $f$ and $g$. In classical category theory, we can {\em define} $M \otimes_{A} N$ to be this coequalizer. However, in the setting of higher category theory, we are not yet done. Namely, we should not expect an equality $\alpha \circ f = \alpha \circ g$; rather, we expect a homotopy $h: \alpha \circ f \rightarrow \alpha \circ g$. Consider now the diagram
$$\xymatrix{ M \otimes A \otimes A \otimes N \ar@<.4ex>[r]^-{f'} \ar@<-.4ex>[r]_-{g'} & M \otimes N}.$$
where, as before, $f'$ is obtained using the action of $A$ on $M$ (twice) and $g'$ is obtained using the action of $A$ on $N$ (twice). The homotopy $h$ can be used to obtain two {\em different} homotopies between $\alpha \circ f'$ and $\alpha \circ g'$: first by applying $h$ to each copy of $A$ individually, and second by applying $h$ together with the observation that $f'$ and $g'$ factor through the map $M \otimes A \otimes A \otimes N \rightarrow M \otimes A \rightarrow N$
determined by the algebra structure on $A$. In order to obtain the correct relative tensor product
$M \otimes_{A} N$, we need to guarantee that these two homotopies coincide, up to a higher homotopy. Of course, this leads to still higher obstructions which need to be taken under consideration.

To efficiently organize all of the relevant data, it is useful to consider a simplicial object
$\Baar_{A}(M,N)_{\bigdot}$, which may be informally described as follows:
\begin{itemize}\index{ZZZBaar@$\Baar_{A}(M,N)_{\bigdot}$}
\item[$(1)$] For each $n \geq 0$, the object $\Baar_{A}(M,N)_{n} \in \calC$ is given by
the tensor product $M \otimes A \otimes \ldots \otimes A \otimes N$, where $n$ factors of $A$ appear. 
\item[$(2)$] If $i = 0$, the face map
$d_i : \Baar_{A}(M,N)_{n} \rightarrow \Baar_{A}(M,N)_{n-1}$ is given by the right action of $A$ on $M$. If $i = n$, $d_i$ is given by the left action of $A$ on $N$. If $0 < i < n$, then $d_i$ is given by the algebra structure on $A$, applied to the $i$ and $(i+1)$st factors.
\item[$(3)$] For $0 \leq i \leq n$, the degeneracy map
$s_i: \Baar_{A}(M,N)_{n} \rightarrow \Baar_{A}(M,N)_{n+1}$ is given by the composition
$$ M \otimes A \otimes \ldots \otimes A \otimes N \simeq
M \otimes A \otimes \ldots \otimes 1 \otimes \ldots \otimes A \otimes N
\rightarrow M \otimes A \otimes \ldots \otimes A \otimes \ldots \otimes A \otimes N$$
where the second map is given by the unit map of the $(i+1)$st factor of $A$.
\end{itemize}

We can then define $M \otimes_{A} N$ to be the geometric realization of the simplicial object $\Baar_{A}(M,N)_{\bigdot}$. The functor
$(M,N) \mapsto | \Baar_{A}(M,N)_{\bigdot} |$ is called the ({\it two-sided}) {\it bar construction}. We now recast the preceding discussion in more formal terms.\index{bar construction}

\begin{definition}\label{barcon}\index{balanced pairing}
Let $p: \calC^{\otimes} \rightarrow \Nerve(\cDelta)^{op}$ be a monoidal $\infty$-category. Let
$q: \calM^{\otimes} \rightarrow \calC^{\otimes}$ be an $\infty$-category which is right-tensored over $\calC^{\otimes}$, and let $q': \calN^{\otimes} \rightarrow \calC^{\otimes}$ be an $\infty$-category which is left-tensored over $\calC$. A {\it balanced pairing} is a functor
$ F: \calM^{\otimes} \times_{ \calC^{\otimes} } \calN^{\otimes} \rightarrow \calD$
with the following property: let $(\alpha, \beta)$ be a morphism in 
$\calM^{\otimes} \times_{ \calC^{\otimes} } \calN^{\otimes}$ such that
$\alpha$ is a $(p \circ q)$-coCartesian morphism in $\calM^{\otimes}$ and $\beta$ is 
a $(p \circ q')$-coCartesian morphism in $\calN^{\otimes}$. Then $F(\alpha,\beta)$ is an equivalence in $\calD$.

Under these conditions, we let $\Baar_{\bigdot}$ denote the composition
$$ \Mod^{R}(\calM) \times_{ \Alg(\calC) } \Mod^{L}(\calN)
\subseteq \Fun( \Nerve( \cDelta)^{op}, \calM^{\otimes} \times_{\calC^{\otimes}} \calN^{\otimes})
\stackrel{F \circ}{\rightarrow} \Fun( \Nerve(\cDelta)^{op}, \calD).$$
We will refer to $\Baar_{\bigdot}$ as the {\it two-sided bar construction}. If $A \in \Alg(\calC)$, $M \in \Mod^{R}_{A}(\calM)$ and $N \in \Mod^{L}_{A}(\calN)$, we will denote the image of $(M,N)$ under
$\Baar_{\bigdot}$ by $\Baar_{A}(M,N)_{\bigdot}$. 

Suppose now that $\calD$ admits geometric realizations for simplicial objects. We then
define the {\it relative tensor product functor} to be the composition
$$\Mod^{R}(\calM) \times_{ \Alg(\calC) } \Mod^{L}(\calN) \stackrel{\Baar_{\bigdot}}{\rightarrow}
\Fun( \Nerve(\cDelta)^{op}, \calD) \stackrel{||}{\rightarrow} \calD,$$
where the second arrow is a geometric realization functor (that is, a left adjoint to the diagonal embedding $\calD \rightarrow \Fun( \Nerve(\cDelta)^{op}, \calD)$).
\end{definition}\index{relative tensor product}

\begin{remark}
In the situation of Definition \ref{barcon}, we can identify objects of
$\Mod^R(\calM) \times_{ \Alg(\calC) } \Mod^L(\calN)$ with triples
$(M,A,N)$, where $A$ is an algebra object of $\calC$, $M$ is a right $A$-module, and
$N$ is a left $A$-module. In this case, we will denote the image of $(M,A,N)$ under the relative tensor product functor by $M \otimes_{A} N$.
\end{remark}

\begin{remark}
Let $\calM^{\otimes} \stackrel{q}{\rightarrow} \calC^{\otimes} \stackrel{q'}{\leftarrow} \calN^{\otimes}$
be as in Definition \ref{barcon}, and let $F: \calM^{\otimes} \times_{ \calC^{\otimes} } \calN^{\otimes} \rightarrow \calD$ be a balanced pairing. Let $\calM = \calM^{\otimes}_{[0]}$, $\calN = \calN^{\otimes}_{[0]}$. Then the inclusion $\calM^{\otimes}_{[0]} \times_{ \calC^{\otimes}_{[0]}} \calN^{\otimes}_{[0]} \subseteq \calM \times \calN$ is an equivalence of $\infty$-categories, so that $F$ the restriction of $F$ to the inverse image of $[0] \in \cDelta$ determines
a functor $F_0: \calM \times \calN \rightarrow \calD$, well-defined up to equivalence.  We will denote this functor by $F_0(M,N) = \langle M, N \rangle.$\index{ZZZbrackets@$\langle M,N \rangle$}
Likewise, the restriction of $F$ to the inverse image of $[1] \in \cDelta$ determines a functor 
$F_1: \calM \times \calC \times \calN \rightarrow \calD$
Moreover, since $F$ carries coCartesian morphisms to equivalences in $\calD$, we have natural equivalences $$F_0( M \otimes C, N) \simeq F_1(M,C,N) \simeq F_0( M, C \otimes N),$$ which we can view as a compatibility condition on $\langle, \rangle$:
$$ \gamma: \langle M \otimes C, N \rangle \simeq \langle M, C \otimes N \rangle.$$
The restriction of $F$ to the inverse image of the remainder of $\Nerve(\cDelta)^{op}$ can be viewed as expressing the compatibility of $\gamma$ with the associative tensor product on the $\infty$-category $\calC$.
\end{remark}

Let $\calC$ be a monoidal $\infty$-category, let $\calM$ be an $\infty$-category which is right-tensored over $\calC$. Our next goal is to show that the prescription
$ \langle M, C \rangle = M \otimes C$ can be extended to a balanced pairing
$\calM^{\otimes} \times_{ \calC^{\otimes} } \calC^{\otimes, L} \rightarrow \calM$.

\begin{proposition}\label{umat}
Let $\calC$ be a monoidal $\infty$-category with unit object $1_{\calC}$, let $\calM$ be an $\infty$-category which is right-tensored over $\calC$, let $\calD$ be an arbitrary $\infty$-category, and let
$$ \Fun'( \calM^{\otimes} \times_{\calC^{\otimes}} \calC^{\otimes, L}, \calD) \subseteq
\Fun( \calM^{\otimes} \times_{\calC^{\otimes}} \calC^{\otimes,L}, \calD)$$
denote the full subcategory spanned by the balanced pairings. Then the functor
$$ F \mapsto \langle 1_{\calC}, \bigdot \rangle $$
induces an equivalence of $\infty$-categories
$ \Fun'( \calM^{\otimes} \times_{\calC^{\otimes}} \calC^{\otimes, L}, \calD)
\rightarrow \Fun(\calM, \calD).$
\end{proposition}

\begin{proof}
We will use the theory of marked simplicial sets described in \S \toposref{twuf}. However,
we make a slight departure from the notation employed there: given a {\em coCartesian} fibration of simplicial sets $p: X \rightarrow S$, we let $X^{\natural}$ denote the marked simplicial set
$(X, \calE)$, where $\calE$ is the set of $p$-coCartesian edges of $X$. This notation is potentially ambiguous, since it depends not only on $X$ but also on the map $p$. In practice, we will have
either $S \simeq \Nerve( \cDelta)^{op}$ or $S \simeq \Delta^0$ (in the latter case, $X$ is an
$\infty$-category and $\calE$ is the collection of all equivalences in $X$).

Let us introduce a bit of notation. 
Let $q: \calM^{\otimes} \rightarrow \calC^{\otimes} \rightarrow \Nerve(\cDelta)^{op}$ exhibit $\calM = \calM^{\otimes}_{[0]}$ as right-tensored over $\calC = \calC^{\otimes}_{[1]}$. 
Let $\calN$ denote the fiber product $\calM^{\otimes} \times_{\calC^{\otimes}} \calC^{\otimes, L}$.
Let $\alpha: \Delta^1 \times \Nerve(\cDelta)^{op} \rightarrow \Nerve(\cDelta)^{op}$ be defined as in Example \ref{sumai}. We define a simplicial set $\calM^{\otimes,L}$ equipped with a map
$\calM^{\otimes,L} \rightarrow \Nerve( \cDelta)^{op}$ via the formula
$$ \Hom_{ \Nerve(\cDelta)^{op}}(K, \calM^{\otimes,L})
= \Hom'_{ \Nerve( \cDelta)^{op}}(\Delta^1 \times K, \calM^{\otimes})$$
where $\Delta^1 \times K$ maps to $\Nerve( \cDelta)^{op}$ via the composition
$\Delta^1 \times K \rightarrow \Delta^1 \times \Nerve(\cDelta)^{op} \stackrel{\alpha}{\rightarrow} \Nerve(\cDelta)^{op},$
and $$\Hom'_{ \Nerve(\cDelta)^{op}}( \Delta^1 \times K, \calM^{\otimes}) \subseteq
\Hom_{\Nerve(\cDelta)^{op}}( \Delta^1 \times K, \calM^{\otimes})$$ denotes the subset consisting of those maps
$f: \Delta^1 \times K \rightarrow \calM^{\otimes}$ which carry each edge $\Delta^1 \times \{k\}$ to a 
$q$-coCartesian edge of $\calM^{\otimes}$.

Let $\Nerve(\cDelta')^{op}$ denote the $\infty$-category $\Nerve(\cDelta)^{op}$, regarded as an object of $(\sSet)_{/ \Nerve(\cDelta)^{op}}$ via the map $\alpha_0: \alpha | \{0\} \times \Nerve(\cDelta)^{op}$, and let $\calM^{\otimes}_{L}$ denote the fiber product $\calM^{\otimes} \times_{ \Nerve(\cDelta)^{op}} \Nerve(\cDelta')^{op}$. Restriction determines a functor $\calM^{\otimes,L} \rightarrow \calM^{\otimes}_{L}$. Regard $\calM^{\otimes}_{L}$ as a simplicial set
over $\Nerve(\cDelta)^{op}$ via the composition
$\calM^{\otimes}_{L} \rightarrow \Nerve(\cDelta')^{op} \simeq \Nerve(\cDelta)^{op}$.

Let $\calN$ denote the fiber product $\calM^{\otimes} \times_{\calC^{\otimes}} \calC^{\otimes, L}$.
We observe that $\Fun'( \calN, \calD)$ and $\Fun(\calM, \calD)$ can be identified with
$\bHom^{\flat}( \calN^{\natural}, \calD^{\natural})$ and $\bHom^{\flat}( \calM^{\natural}, \calD^{\natural})$, respectively. In view of Proposition \toposref{markdefeq}, it will suffice to show that
the map
$\calM \simeq \{1_{\calC} \} \times \calM \subseteq \calC \times \calM
\simeq \calN_{[0]} \subseteq \calN$
induces an equivalence of marked simplicial sets $\calM^{\natural} \subseteq \calN^{\natural}$.
Corollary \toposref{usefir} implies that the restriction maps
$\calN \leftarrow \calM^{\otimes,L} \rightarrow \calM^{\otimes}_{L}.$
are trivial Kan fibrations. Consequently, it will suffice to prove that the composition
$\psi: \calM \simeq \{1_{\calC} \} \times \calM \subseteq \calC \times \calM
\simeq (\calM^{\otimes}_{L})_{[0]} \subseteq \calM^{\otimes}_{L}$
induces a marked equivalence $\calM^{\natural} \rightarrow ( \calM^{\otimes}_{L} )^{\natural}$. 

Let $F: \Nerve(\cDelta)^{op} \rightarrow \Cat_{\infty}$ classify the coCartesian fibration $q$.
According to Proposition \toposref{charcatcolimit}, the marked simplicial set
$(\calM^{\otimes}_{L})^{\natural}$ can be identified with the colimit of diagram
$\Nerve( \cDelta')^{op} \rightarrow \Nerve(\cDelta)^{op} \stackrel{F}{\rightarrow} \Cat_{\infty}.$
Lemma \toposref{aclock} implies that this colimit is equivalent to $F( [0] ) \simeq \calM$. 
We now observe that the composition of this equivalence with the map $\psi$ is given
by the functor $M \mapsto 1_{\calC} \otimes M$, which is equivalent to $\id_{\calM}$ and therefore an equivalence.
\end{proof}

\begin{definition}\label{cannapair}\index{balanced pairing!canonical}\index{canonical balanced pairing}
Let $\calC$ be a monoidal $\infty$-category with unit object $1_{\calC}$, let $\calM$ be an $\infty$-category which is right-tensored over $\calC$. A {\it canonical balanced pairing}
is a balanced pairing
$\calM^{\otimes} \times_{\calC^{\otimes}} \calC^{\otimes, L} \rightarrow \calM$
whose image under the equivalence described in Proposition \ref{umat} is equivalent
to the identity functor $\id_{\calM} \in \Fun( \calM, \calM)$.
\end{definition}

In other words, a canonical balanced pairing is a balanced pairing for which
the induced functor $\langle 1_{\calC}, \bigdot \rangle: \calM \rightarrow \calM$ is equivalent to the identity.

\begin{remark}
In the situation of Proposition \ref{umat}, composition with a canonical balanced pairing
$$ \calM^{\otimes} \times_{\calC^{\otimes}} \calC^{\otimes, L} \rightarrow \calM$$ induces a homotopy inverse to the equivalence
$ \Fun'( \calM^{\otimes} \times_{\calC^{\otimes}} \calC^{\otimes, L}, \calD)
\rightarrow \Fun(\calM, \calD).$
\end{remark}

\begin{remark}
If $\calM$ is instead right-tensored over $\calC$, then we have a dual notion of
a canonical balanced pairing 
$\calC^{\otimes, R} \times_{\calC^{\otimes}} \calM^{\otimes} \rightarrow \calM$. We note that
this leads to two different definitions of a canonical balanced pairing
$\calC^{\otimes, R} \times_{\calC^{\otimes}} \calC^{\otimes, L} \rightarrow \calC.$
However, these two notions coincide, since the functors
$$C \mapsto \langle 1_{\calC}, C \rangle, \quad  C \mapsto \langle C, 1_{\calC} \rangle $$
are equivalent to one another.
\end{remark}

Proposition \ref{umat} is really a special case of the following more general result
(applied to a canonical pairing $\Cat_{\infty}^{\times, R} \times_{ \Cat_{\infty} } \Cat_{\infty}^{\times, L}
\rightarrow \Cat_{\infty}$):

\begin{proposition}\label{usss}
Let $\calC$ be a monoidal $\infty$-category, $\calM$ an $\infty$-category right-tensored over 
$\calC$, $A$ an algebra object of $\calC$, $\calD$ an $\infty$-category which admits geometric realizations, $F: \calM^{\otimes} \times_{ \calC^{\otimes} } \calC^{\otimes,L} \rightarrow \calD$ a balanced pairing, and 
$$ \otimes_{A}: \Mod^{R}_{A}(\calM) \times \Mod^{L}_{A}(\calC) \rightarrow \calD$$
the relative tensor product functor.
Then the canonical map
$$ \phi: \langle M, 1_{\calC} \rangle \rightarrow \langle M, A \rangle
\simeq \Baar_A(M,A)_1 \rightarrow | \Baar_{A}(M,A)_{\bigdot} | \simeq M \otimes_{A} A$$
determines an equivalence between $\bigdot \otimes_{A} A$ and $\langle \bigdot, 1_{\calC} \rangle$ $($regarded as functors from $\Mod_{A}(\calM)$ to $\calM${}$)$.
\end{proposition}

\begin{proof}
Let $\calN$, $\calM^{\otimes,L}$, and $\calM^{\otimes}_{L}$ be defined as in the proof of Proposition \ref{umat}, and let
$F'$ denote the composition
$\calM^{\otimes}_{L} \stackrel{s}{\rightarrow} \calM^{\otimes,L} \rightarrow
\calN \stackrel{F'}{\rightarrow} \calD.$
Define subcategories $\cDelta'' \subseteq \cDelta' \subseteq \cDelta$ as follows:
\begin{itemize}
\item[$(i)$] Every object of $\cDelta$ is an object of $\cDelta'$; an object
$[n] \in \cDelta$ belongs to $\cDelta''$ if and only if $n \geq 0$.
\item[$(ii)$] A morphism $\alpha: [m] \rightarrow [n]$ of $\cDelta$ belongs to
$\cDelta'$ if and only if $\alpha(m) = n$. If $m, n > 0$, then $\alpha$ belongs to
$\cDelta''$ if and only if $\alpha^{-1} \{n\} = \{m\}$. 
\end{itemize}
We have an equivalence of categories $\cDelta \simeq \cDelta''$, given by
$[n] \mapsto [n] \star [0]$.

Observe that $\Nerve(\cDelta')^{op}$ contains $[0]$ as a final object (in fact, as a zero object). Let $\calM^{\otimes}_{\infty}$ denote the fiber product
$\Nerve(\cDelta')^{op} \times_{ \Nerve(\cDelta)^{op} } \calM^{\otimes},$ and let $q_{\infty}: \calM^{\otimes}_{\infty} \rightarrow \Nerve(\cDelta')^{op}$
be the projection. Lemma \toposref{aclock} and Proposition \toposref{charcatcolimit} imply that the inclusion
$(\calM^{\otimes}_{L})^{\natural} \subseteq (\calM^{\otimes}_{\infty})^{\natural}$
is an equivalence of marked simplicial sets. It follows that $F'$ extends to a map
$\overline{F}': \calM^{\otimes}_{\infty} \rightarrow \calD,$
where $\overline{F}'$ carries $q_{\infty}$-coCartesian morphisms to equivalences in $\calD$.

The simplicial object $\Baar_{A}(M,A)_{\bigdot}$ can be identified with the composition
$$ \Nerve(\cDelta)^{op} \simeq \Nerve(\cDelta'')^{op} \stackrel{M}{\rightarrow} \calM^{\otimes}_{L} \stackrel{F'}{\rightarrow} \calD.$$
It follows that $\Baar_{A}(M,A)_{\bigdot}$ extends to a functor
$U: \Nerve(\cDelta_{\infty})^{op} \stackrel{M}{\rightarrow}
\calM^{\otimes}_{\infty} \stackrel{ \overline{F}'}{\rightarrow} \calD.$
Lemma \toposref{aclock} implies that $U$ determines an equivalence
$$ M \otimes_{A} A \simeq | \Baar_{A}(M,A)_{\bigdot} | \simeq
U([0]) \simeq \overline{F}'( M([0])).$$
We note that there is a $q_{\infty}$-coCartesian edge joining $(M([0]), 1_{\calC})
\in \calM \times \calC \simeq \calM^{\otimes}_{[1]}$ to $M([0])$ in $\calM^{\otimes}_{\infty}$.
Since $\overline{F}'$ carries $q_{\infty}$-coCartesian edges to equivalences in $\calD$, we conclude that the map $\phi$ is an equivalence, as desired.
\end{proof}

\begin{corollary}\label{sharkon}
Let $\calM^{\otimes} \stackrel{q}{\rightarrow} \calC^{\otimes} \stackrel{q'}{\leftarrow} \calN^{\otimes}$
be as in Definition \ref{barcon}, let $F: \calM^{\otimes} \times_{ \calC^{\otimes} } \calN^{\otimes} \rightarrow \calD$ be a balanced pairing. Suppose that $\calD$ admits geometric realizations.
Let $A \in \Alg(\calC)$, $M \in \Mod_{A}(\calM)$, $N_0 \in \calN$, and let $N \in \Mod_{A}(\calN)$ be a free $A$-module generated by $N_0$. Then the composite map
$$ \langle M, N_0 \rangle \rightarrow \langle M, N \rangle
= \Baar_{A}(M,N)_{\bigdot} \rightarrow | \Baar_{A}(M,N)_{\bigdot} | \simeq M \otimes_{A} N$$
is an equivalence in $\calD$.
\end{corollary}

\begin{proof}
Using Corollary \ref{specialtime}, we can choose a map of $\infty$-categories
$F: \calC^{\otimes, L} \rightarrow \calN^{\otimes}$ of $\infty$-categories left-tensored over
$\calC^{\otimes}$, such that $F$ carries the unit object of $\calC$ to $N_0 \in \calN$.
In view of Proposition \ref{pretara}, we may identify $N$ with the image under $F$
of $A$, regarded as a left $A$-module over itself as in Example \ref{algitself}.
We may therefore reduce to the case where $\calN^{\otimes} = \calC^{\otimes,L}$, and
$N_0 = 1_{\calC} \in \calC$. The desired result now follows from Proposition \ref{usss}.
\end{proof}



We now establish some formal properties enjoyed by the two-sided bar construction and the relative tensor product.

\begin{proposition}\label{colimp}
Let $\calM^{\otimes} \stackrel{q}{\rightarrow} \calC^{\otimes} \stackrel{q'}{\leftarrow} \calN^{\otimes}$
be as in Definition \ref{barcon}, let $F: \calM^{\otimes} \times_{ \calC^{\otimes} } \calN^{\otimes} \rightarrow \calD$ be a balanced pairing. Let $K$ be a simplicial set and let $A \in \Alg(\calC)$. 
Suppose that:
\begin{itemize}
\item[$(i)$] The $\infty$-categories $\calD$ and $\calN = \calN^{\otimes}_{[0]}$ admit $K$-indexed colimits.
\item[$(ii)$] The $\infty$-category $\calD$ admits geometric realizations for simplicial objects.
\item[$(iii)$] For each $M \in \calM$, the pairing $\langle M, \bigdot \rangle:
\calN \rightarrow \calD$ preserves $K$-indexed colimits.
\end{itemize}
Then for each $M \in \Mod^{R}_{A}(\calM)$, the relative tensor product functor $M \otimes_{A} \bigdot: \Mod^{L}_{A}(\calN) \rightarrow \calD$ preserves $K$-indexed colimits.
\end{proposition}

\begin{proof}
In view of Lemma \toposref{limitscommute}, it will suffice to show that for every $n \geq 0$, the functor $$N \mapsto \Baar_{A}(M, N)_{\bigdot}$$ preserves $K$-indexed colimits. We now observe that
there is canonical equivalence $\Baar_{A}(M,N)_{\bigdot} \simeq \langle M \otimes A^{\otimes n}, N \rangle$ of functors $\calN \rightarrow \calD$. The desired result now follows from $(iii)$.
\end{proof}

Let $\calM^{\otimes} \stackrel{q}{\rightarrow} \calC^{\otimes} \stackrel{q'}{\leftarrow} \calN^{\otimes}$ be as in Definition \ref{barcon}, and let $F: \calM^{\otimes} \times_{ \calC^{\otimes} } \calN^{\otimes} \rightarrow \calD$ be a balanced pairing. Let $A$ be an algebra object of $\calC$, and let
$$ \theta: \Mod^{R}_{A}(\calM) \rightarrow \calM, \quad \theta': \Mod^{L}_{A}(\calN) \rightarrow \calN$$ denote the forgetful functors. We observe that, for every
$M \in \Mod^{R}_{A}(\calM)$ and $N \in \Mod^{L}_{A}(\calN)$, we have a canonical map
$$\langle \theta(M), \theta'(N) \rangle \simeq \Baar_{A}(M,N)_{\bigdot} \rightarrow | \Baar_{A}(M,N)_{\bigdot} |
\simeq M \otimes_{A} N.$$

\begin{proposition}\label{pairweese}
Let $\calM^{\otimes} \stackrel{q}{\rightarrow} \calC^{\otimes} \stackrel{q'}{\leftarrow} \calN^{\otimes}$
be as in Definition \ref{barcon}, let $F: \calM^{\otimes} \times_{ \calC^{\otimes} } \calN^{\otimes} \rightarrow \calD$ be a balanced pairing. Suppose that $\calD$ admits geometric realizations for simplicial objects, and let $A$ be an initial object of $\Alg(\calC)$ $($see Proposition \ref{gurgle}$)$. Then for every $M \in \Mod_{A}^{R}(\calM)$, $N \in \Mod_{A}^{L}(\calN)$, the natural map
$\phi: \langle \theta(M), \theta'(N) \rangle \rightarrow M \otimes_{A} N$ is an equivalence in $\calD$.
\end{proposition}

\begin{proof}
The relative tensor product $M \otimes_{A} N$ is defined as the geometric realization of the simplicial object $\Baar_{A}(M,N)_{\bigdot}: \Nerve(\cDelta)^{op} \rightarrow \calD$. 
Since $\Nerve(\cDelta)^{op}$ is weakly contractible, Corollary \toposref{silt} implies that
$\phi$ is an equivalence provided that $\Baar_{A}(M,N)_{\bigdot}$ carries each morphism in
$\cDelta$ to an equivalence in $\calD$. Unwinding the definitions, we see that each of these morphisms can be identified with the map
$$ \langle \theta(M), A^{\otimes m} \otimes \theta(N) \rangle \rightarrow \langle \theta(N), A^{\otimes n} \otimes \theta(N) \rangle$$
determined by a map of linearly ordered sets $[m] \rightarrow [n]$. The desired result now follows easily from the fact that $A$ is a unit object of $\calC$ (Proposition \ref{gurgle}).
\end{proof}

\begin{lemma}\label{peacestick}
Let $\calC$ be an $\infty$-category equipped with a monoidal structure, and let
$\calM$ be an $\infty$-category which is left-tensored over $\calC$. Suppose that $\calM$ admits geometric realizations, and that for each $C \in \calC$ the functor $C \otimes \bigdot$ preserves geometric realizations. Let $f: A \rightarrow A'$ be a morphism in $\Alg(\calC)$. Then the associated functor $\theta: \Mod_{A'}(\calM) \rightarrow \Mod_{A}(\calM)$ $($see Corollary \ref{thetacart}$)$ admits a left adjoint.
\end{lemma}

\begin{proof}
Apply Corollary \ref{littlerbeck} to the homotopy commutative diagram
$$ \xymatrix{ \Mod_{A'}(\calM) \ar[rr] \ar[dr] & & \Mod_{A}(\calM) \ar[dl] \\
& \calM, & }$$
together with Proposition \ref{pretara} and Corollary \ref{gloop}.
\end{proof}



\begin{proposition}[Push-Pull Formula]\label{ussr}\index{push-pull formula}
Let $$\calM^{\otimes} \stackrel{q}{\rightarrow} \calC^{\otimes} \stackrel{q'}{\leftarrow} \calN^{\otimes}$$
be as in Definition \ref{barcon}, let $F: \calM^{\otimes} \times_{ \calC^{\otimes} } \calN^{\otimes} \rightarrow \calD$ be a balanced pairing. Suppose that:
\begin{itemize}
\item[$(i)$] The $\infty$-categories $\calD$ and $\calN$ admit geometric realizations.
\item[$(ii)$] For every $C \in \calC$, the functor $C \otimes \bigdot: \calN \rightarrow \calN$
preserves geometric realizations.
\item[$(iii)$] For every $M \in \calM$, the pairing
$\langle M, \bigdot \rangle: \calN \rightarrow \calD$ preserves geometric realizations.
\end{itemize}
Let $\Mod^{R}(\calM) \stackrel{r}{\rightarrow} \Alg(\calC) \stackrel{r'}{\leftarrow} \Mod^{L}(\calN)$
be the induced functors, and let 
$$\alpha = (\alpha_{R}, \alpha', \alpha_{L}): (M,A,N) \rightarrow (M',A',N')$$ be a morphism
in $\Mod^{R}(\calM) \times_{ \Alg(\calC) } \Mod^{L}(\calN)$. Suppose that
$\alpha_{R}$ is an $r$-Cartesian morphism of $\Mod^{R}(\calM)$ and
$\alpha_{L}$ is an $r'$-coCartesian morphism of $\Mod^{L}(\calN)$. Then $\alpha$ induces an equivalence $M \otimes_{A} N \rightarrow M' \otimes_{A'} N'$ in $\calD$.
\end{proposition}

\begin{proof}
Let us first regard $\alpha_{R}: M \rightarrow M'$ and $\alpha': A \rightarrow A'$ as fixed.
Lemma \ref{peacestick} implies that $r'$ is a coCartesian fibration, so that
for each $N \in \Mod_{A}(\calM)$ there exists an (essentially unique) $r'$-coCartesian morphism
$\alpha_{L}: N \rightarrow \alpha'_{!} N$ lifting $\alpha'$. Let
$\Mod'_{A}(\calN) \subseteq \Mod_{A}(\calN)$ denote the full subcategory spanned by those $A$-modules $N$ for which the induced map $M \otimes_{A} N \rightarrow M' \otimes_{A'} (\alpha'_{!} N)$ is an equivalence. Proposition \ref{colimp} implies that $\Mod'_{A}(\calN)$ is closed under geometric realizations in $\Mod_{A}(\calN)$ (which admits geometric realizations in view of
Corollary \ref{gloop}). Proposition \ref{littlebeck} implies that $\Mod_{A}(\calN)$ is generated, under geometric realizations, by the essential image of the functor
$\beta_{!}: \Mod_{A_0}(\calN) \rightarrow \Mod_{A}(\calN)$, where $A_0$ denotes the initial object of $\Alg(\calC)$. Consequently, we may assume that $N \simeq \beta_! N_0$. Let $M_0$ denote the image of $M$ under the forgetful functor $\Mod_{A}(\calM) \rightarrow \Mod_{A_0}(\calM)$.
We have a commutative triangle
$$ \xymatrix{ & M \otimes_{A} N \ar[dr]^{\psi} & \\
M_0 \otimes_{A_0} N_0 \ar[rr]^{\psi''} \ar[ur]^{\psi'} & & M' \otimes_{A'} N' }$$
in $\calD$. Consequently, to prove that $\psi$ is an equivalence, it will suffice to prove that
$\psi'$ and $\psi''$ are equivalences. In other words, we may reduce to the case where $A$ is an initial object of $\Alg(\calC)$. The desired result now follows immediately from Proposition \ref{pairweese} and Corollary \ref{sharkon}.
\end{proof}

\begin{corollary}\label{tara}
Let $\calC$, $\calM$, and $f: A \rightarrow A'$ be as in Lemma \ref{peacestick}, and let $F: \Mod_{A}(\calM) \rightarrow \Mod_{A'}(\calM)$ be the left adjoint to the forgetful functor
$\Mod_{A'}(\calM) \rightarrow \Mod_{A}(\calM)$, and let $\theta: \Mod_{A'}(\calM) \rightarrow \calM$ be the forgetful functor. Then the composition
$\Mod_{A}(\calM) \stackrel{F}{\rightarrow} \Mod_{A'}(\calM) \stackrel{\theta}{\rightarrow} \calM$
can be identified with the relative tensor product
$M \mapsto A' \otimes_{A} M$
determined by a canonical balanced pairing $\calC^{\otimes, R} \times_{\calC^{\otimes}} \calM^{\otimes} \rightarrow \calM$.
\end{corollary}

\begin{proof}
According to Proposition \ref{usss}, the forgetful functor
$\theta$ can be identified with the relative tensor product $A' \otimes_{A'} \bigdot$.
Proposition \ref{ussr} implies that $A' \otimes_{A'} F(M)$ is (canonically) equivalent to
$A' \otimes_{A} M$. Here we identify $A'$ with its image in $\Mod^{R}_{A}(\calC)$ under the forgetful functor $\Mod^{R}_{A'}(\calC) \rightarrow \Mod^{R}_{A}(\calC)$. 
\end{proof}

\begin{remark}
Let $\calC$, $\calM$, and $f: A \rightarrow A'$ satisfy the hypotheses of Lemma \ref{peacestick}. We will sometimes abuse notation by indicating the associated functor $\Mod^{L}_{A}(\calM) \rightarrow \Mod^{L}_{A'}(\calM)$ by $M \mapsto A' \otimes_{A} M$. This notation is partially justified by Corollary \ref{tara}.
\end{remark}

\subsection{Flat Modules}\label{monoid13}

Let $A$ be a connective $A_{\infty}$-ring. In this section, we will show that there is a good theory of flat and projective $A$-modules, which reduces to the classical theory in the case where $A$ is discrete.

\begin{definition}\index{module!free}\index{free!module over a ring spectrum}\label{mofree}
Let $A$ be an $A_{\infty}$-ring. We will say that a left $A$-module $M$ is {\it free} if
is equivalent to a coproduct of copies of $A$ (where we view $A$ as a left module over itself, as in Example \ref{algitself}). We will say that a free module $M$ is {\it finitely generated} if it can be
written as a finite coproduct of copies of $A$.\index{module!finitely generated}

Suppose that $A$ is connective. We will say that a map $f: M \rightarrow N$ of connective (left) $A$-modules is {\it surjective} if it induces a surjection of $\pi_0 A$-modules $\pi_0 M \rightarrow \pi_0 N$.\index{surjective!map of module spectra}
\end{definition}

\begin{warning}
The terminologies introduced in Definitions \ref{mofree} and \ref{squirch} are not quite compatible with one another. Let $M$ be a (left) module over an $A_{\infty}$-ring $A$. Then $M$ is free
(in the sense of Definition \ref{mofree}) if and only if it is freely generated (in the sense of
Definition \ref{squirch}) by a coproduct of copies of the sphere spectrum. If $M_0$ is a general spectrum, then the left $A$-module $A \otimes M_0$ freely generated by $M_0$ is typically not free in the sense of Definition \ref{mofree}. 
\end{warning}

\begin{remark}
Using the long exact sequence of homotopy groups associated to an exact triangle, we conclude that a map $f: M \rightarrow N$ of connective modules over a connective $A_{\infty}$-ring is surjective if and only if $\ker(f)$ is also connective.
\end{remark}

\begin{definition}\index{projective!module}\index{module!projective}\label{projmoj}
Let $A$ be a connective $A_{\infty}$-ring. We will say that a left $A$-module $P$ is {\it projective} if it
a projective object of the $\infty$-category $( \Mod_{A}^{L})_{\geq 0}$ of connective left $A$-modules, in the sense of Definition \stableref{humber}.
\end{definition}

\begin{remark}
The terminology of Definition \ref{projmoj} is potentially ambiguous: a projective
left $A$-module is typically not projective as an object of $\Mod_{A}^{L}$. However, there is little risk of confusion, since the $\infty$-category $\Mod^{L}_{A}$ has no nonzero projective objects.
\end{remark}

\begin{proposition}\label{charp}
Let $A$ be a connective $A_{\infty}$-ring, and let $P$ be a connective left $A$-module. The following conditions are equivalent:
\begin{itemize}
\item[$(1)$] The left $A$-module $P$ is projective.
\item[$(2)$] There exists a free $A$-module $M$ such that $P$ is a retract of $M$.
\end{itemize}
\end{proposition}

\begin{proof}
Suppose first that $P$ is projective.
Choose a map of left $A$-modules $p: M \rightarrow P$, where $M$ is free and the induced map
$\pi_0 M \rightarrow \pi_0 P$ is surjective (for example, we can take $M$ to be a direct sum of copies of $A$ indexed by the set $\pi_0 P$). Invoking Proposition \stableref{charprojjj}, we deduce that $p$ admits a section (up to homotopy), so that $P$ is a retract of $M$. This proves $(2)$. To prove the converse, we observe that the collection of projective left $A$-modules is stable under retracts. It will therefore suffice to show that every free left $A$-module is projective. This follows immediately from the characterization given in Proposition \stableref{charprojjj}.
\end{proof}

\begin{remark}\label{charpp}
It follows from the proof of Proposition \ref{charp} that, if $\pi_0 P$ is a finitely generated left module
over $\pi_0 A$, then we can choose $M$ to be a finitely generated free $A$-module.
\end{remark}

\begin{corollary}\label{chadwick}
Let $A$ be a connective $A_{\infty}$-ring. Then the $\infty$-category
$(\Mod_{A}^{L})_{\geq 0}$ is projectively generated (Definition \stableref{defpro}). 
Moreover, the following conditions on a connective left $A$-module $P$ are equivalent:
\begin{itemize}
\item[$(1)$] The $A$-module $P$ is projective, and $\pi_0 P$ is finitely generated as a $\pi_0 A$-module.
\item[$(2)$] The $A$-module $P$ is a compact projective object of $( \Mod_{A}^{L})_{\geq 0}$.
\item[$(3)$] There exists a finitely generated free $A$-module $M$ such that $P$ is a retract of $M$.
\end{itemize}
\end{corollary}

\begin{proof}
The first assertion and the equivalence $(2) \Leftrightarrow (3)$ follow 
by applying Corollary \ref{progen} to the composition
$$ (\Mod^{L}_{A})_{\geq 0} \rightarrow \connSpectra \stackrel{\Omega^{\infty}}{\rightarrow} \SSet,$$
and invoking Example \stableref{swine}. The equivalence $(1) \Leftrightarrow (3)$ follows from Remark \ref{charpp}.
\end{proof}

Recall that, if $M$ is a left module over an ordinary ring $A$, we say that $M$ is {\it flat} if the
functor $N \mapsto N \otimes_{A} M$ is exact.

\begin{definition}\label{defflat}\index{flat!module spectrum}\index{module object!flat}
Let $M$ be a left module over an $A_{\infty}$-ring $A$. We will say that $M$ is {\it flat}
if the following conditions are satisfied:
\begin{itemize}
\item[$(1)$] The homotopy group $\pi_0 M$ is flat as a left module over $\pi_0 A$, in the classical sense. 
\item[$(2)$] For each $n \in \Z$, the natural map $\pi_n A \otimes_{ \pi_0 A} \pi_{0} M \rightarrow \pi_{n} M$ is an isomorphism of abelian groups.
\end{itemize}
\end{definition}

\begin{remark}
Let $A$ be a connective $A_{\infty}$-ring. Then every flat left $A$-module is also connective.
\end{remark}

\begin{remark}
Let $A$ be a discrete $A_{\infty}$-ring. A left $A$-module $M$ is flat if and only if $M$ is discrete, and $\pi_0 M$ is flat over $\pi_0 A$ (in the sense of classical ring theory). In other words, Definition \ref{defflat} is compatible with the classical definition of flatness, if we identify discrete $A_{\infty}$-rings and modules with the underlying classical objects.
\end{remark}

In order to use Definition \ref{defflat} effectively, we need to understand the relationship between
the relative tensor products constructed in \S \ref{balpair} and the tensor products in classical noncommutative algebra. 

\begin{notation}\index{ZZZTor@$\Tor_{p}{A}(M,N)_{q}$}
Let $A$ be an ordinary associative ring, let $M$ be a right $A$-module and $N$ a left $A$-module.
Homological algebra associates to the triple $(M,A,N)$ a collection of groups
$\Tor_{p}^{A}(M,N)$, given by the left derived functors of the (classical) tensor product
$\otimes_{A}$ (see \cite{weibel}). We observe that if $A$ is a graded ring and the modules $M$ and $N$ are compatibly graded, then $\Tor_{p}^{A}(M,N)$ inherits a grading (as can be seen by
performing homological algebra in the setting of {\em graded} $A$-modules). In this case, we will denote the $q$th graded piece of $\Tor_{p}^{A}(M,N)$ by $\Tor_{p}^{A}(M,N)_{q}$. 
\end{notation}

\begin{proposition}\label{siiwe}\index{spectral sequence}
Let $A$ be an $A_{\infty}$-ring, let $M$ be a right $A$-module, and let $N$ be a right $A$-module.
We regard $\pi_{\ast} M$ and $\pi_{\ast} N$ as graded modules over the graded ring $\pi_{\ast} A$.
There exists a convergent spectral sequence with $E_2$-page
$$ E_2^{p,q} = \Tor_{p}^{\pi_{\ast} A}( \pi_{\ast} M, \pi_{\ast} N)_{q}
\Rightarrow \pi_{p+q} ( M \otimes_{A} N ).$$
\end{proposition}

\begin{proof}
We will say that a left $A$-module $P$ is {\it quasi-free} if $P$ can be obtained
as a direct sum of $A$-modules of the form $A[n]$, where $n \in \Z$. We will say that
a map of left $A$-modules $P \rightarrow N$ is {\it quasi-surjective} if the induced map
$\pi_{\ast} P \rightarrow \pi_{\ast} N$ is a surjection of graded $\pi_{\ast} A$-modules. 
We observe that for every left $A$-module $N$, there exists a quasi-surjection
$P \rightarrow N$, where $P$ is quasi-free. For example, we can take $P$ to
be a direct sum $\oplus_{\eta} A[n]$, where $\eta$ ranges over $\pi_{n} N$.

We now construct a filtered object
$$ \ldots \rightarrow Q(-1) \rightarrow Q(0) \rightarrow Q(1) \rightarrow \ldots $$
in the $\infty$-category $(\Mod^{L}_{A})_{/N}$. We begin by setting $Q(i) = 0$ for
$i < 0$. Suppose that $i \geq 0$, and that the map $f(i-1): Q(i-1) \rightarrow N$ has already been constructed. Choose a quasi-surjection $P(i) \rightarrow \coker f(i-1)$, where $P(i)$ is quasi-free, and form a pullback square
$$ \xymatrix{ Q(i) \ar[r] \ar[d]^{f(i)} & P(i) \ar[d] \\
N  \ar[r] & \coker f(i-1) }$$
in the $\infty$-category $( \Mod_{A}^{L} )_{Q(i-1)/}$ (and therefore also in $\Mod_{A}^{L}$). 
We have a map of distinguished triangles
$$ \xymatrix{ Q(i-1) \ar[r] \ar@{=}[d] & Q(i) \ar[r] \ar[d]^{f(i)} & P(i) \ar[d] \ar[r] & P(i)[1] \ar[d] \\
Q(i-1) \ar[r]^{f(i-1)} & N \ar[r] & \coker f(i-1) \ar[r] & Q(i-1)[1] . }$$

We have a quasi-surjective map $\delta(0): P(0) \simeq Q(0) \rightarrow N$. For
each $i > 0$, let $\delta(i)$ denote the composite map
$$ P(i+1) \rightarrow Q(i)[1] \rightarrow P(i)[1].$$
We claim that the induced sequence of
left $\pi_{\ast} A$-modules
$$ \ldots \rightarrow \pi_{\ast + 1} P(1) \rightarrow \pi_{\ast} P(0) \rightarrow \pi_{\ast} N \rightarrow 0$$ is exact. The exactness at $\pi_{\ast} N$ is obvious. We will prove the exactness
at $\pi_{\ast+ i} P(i)$ for $i > 0$; the proof for $i = 0$ requires only simple changes of notation and is left to the reader. We first observe that the composition
$\delta(i) \circ \delta(i+1)[-1]$ factors through $Q(i) \rightarrow P(i) \rightarrow Q(i-1)[1]$, which is nullhomotopic. Now suppose that $\eta \in \pi_{n} P(i)$ lies in the kernel of the
map $\pi_{n} P(i) \rightarrow \pi_{n-1} P(i-1)$. Let $\eta_0 \in \pi_{n-1} Q(i-1)$ denote the
image of $\eta$. Using the exactness of the sequence
$$ \pi_{n-1} Q(i-2) \rightarrow \pi_{n-1} Q(i-1) \rightarrow \pi_{n-1} P(i),$$
we conclude that $\eta_0$ is the image of some $\eta_1 \in \pi_{n-1} Q(i-2)$, such
that the image of $\eta_1$ in $\pi_{n-1} N$ is zero. Using the exactness of the
sequence $$ \pi_{n} \coker f(i-2) \rightarrow \pi_{n-1} Q(i-2) \rightarrow \pi_{n-1} N$$
and the assumption that the map $P(i-1) \rightarrow \coker f(i-2)$ is quasi-surjective,
we conclude that $\eta_1$ is itself the image of some $\eta_2 \in \pi_n P(i-1)$.
Since the composition $P(i-1) \rightarrow Q(i-2)[1] \rightarrow Q(i-1)[1]$ is null, we conclude
that $\eta_0 = 0$. The exactness of the sequence
$$ \pi_{n} Q(i-1) \rightarrow \pi_{n} Q(i) \rightarrow \pi_{n} P(i) \rightarrow \pi_{n-1} Q(i-1)$$
shows that $\eta$ is the image of an element $\widetilde{\eta} \in \pi_{n} Q(i)$. Moreover,
we are free to modify $\widetilde{\eta}$ by adding the image of any element
in $\pi_{n} Q(i-1)$. Since the map $\pi_{n} Q(i-1) \rightarrow \pi_{n} N$ is surjective,
we may assume without loss of generality that the image of $\widetilde{\eta}$
in $\pi_{n} N$ is zero. Using the exactness of the sequence
$$ \pi_{n+1} \coker f(i) \rightarrow \pi_{n} Q(i) \rightarrow \pi_{n} N$$
and the assumption that the map $P(i+1) \rightarrow \coker f(i)$ is quasi-surjective,
we conclude that $\widetilde{\eta}$ is the image of an element of $\pi_{n+1} P(i+1)$, as desired.

Let $Q(\infty)$ be a colimit of the sequence
$$ \ldots Q(-1) \rightarrow Q(0) \rightarrow Q(1) \rightarrow \ldots $$
in the $\infty$-category of left $A$-modules. Applying Proposition \stableref{conseq}
to the above sequence (regarded as a filtered object of the $\infty$-category of spectra), we conclude that there is a spectral sequence $\{ E^{p,q}_{r}, d_r \}$ which converges to $\pi_{\ast} Q(\infty)$. The exactness statement established above shows that this spectral sequence
degenerates at $r = 2$, with
$$ E^{p,q}_{2} = \begin{cases} 0 & \text{if } p \neq 0 \\
\pi_{q} N & \text{if } p = 0. \end{cases}$$
It follows that the natural map $Q(\infty) \rightarrow N$ is an equivalence, so that we can
identify $N$ with a colimit of the sequence $Q$. 

Proposition \ref{colimp} implies that $M \otimes_{A} N$ can be identified with a colimit of the filtered spectrum
$$ \ldots \rightarrow M \otimes_{A} Q(0) \rightarrow M \otimes_{A} Q(1) \rightarrow \ldots $$
Applying Proposition \stableref{conseq} again, we deduce the existence of another spectral sequence
$$E^{p,q}_{1} = \pi_{p+q} ( M \otimes_{A} P(p) ) \Rightarrow \pi_{p+q} ( M \otimes_{A} N).$$
To complete the proof, it suffices to observe that since each $P(p)$ is quasi-free, the natural map $\pi_{\ast} M \otimes_{ \pi_{\ast} A} \pi_{\ast} P(p) \rightarrow \pi_{\ast} (M \otimes_{A} P(p) )$
is an isomorphism. It follows that the $E_2^{p,q}$ term of the above spectral sequence
can be identified with the homologies of the complex
$$ \ldots \rightarrow \pi_{\ast} M \otimes_{\pi_{\ast} A} \pi_{\ast+1} P(1)
\rightarrow \pi_{\ast} M \otimes_{ \pi_{\ast} A} \pi_{\ast} P(0),$$
which coincide with the $\Tor$-groups
$\Tor_{\ast}^{\pi_{\ast} A}( \pi_{\ast} M, \pi_{\ast} N)$.
\end{proof}

\begin{remark}
In the situation of Proposition \ref{siiwe}, the $E_1$-page of the spectral sequence depends
on the choice of the filtered object
$$ \ldots Q(-1) \rightarrow Q(0) \rightarrow Q(1) \rightarrow \ldots $$
However, one can show that the $E_2$-page is independent of this choice, and is functorially determined by the triple $(M,A,N)$.
\end{remark}

\begin{corollary}\label{presiwe}
Let $A$ be an $A_{\infty}$-ring, $M$ a right $A$-module, and $N$ a left $A$-module.
Suppose that $A$, $M$, and $N$ are discrete. Then there exists a canonical isomorphism
$$ \pi_{n}(M \otimes_{A} N) \simeq \Tor_{n}^{\pi_0 A}( \pi_0 M, \pi_0 N).$$
\end{corollary}

\begin{corollary}\label{siwe}
Let $A$ be an $A_{\infty}$-ring, $M$ a right $A$-module, and $N$ a left $A$-module.
Suppose that $N$ is flat. For each $n \in \Z$, the canonical homomorphism
$$ \theta: \pi_n M \otimes_{ \pi_0 A} \pi_0 N \rightarrow \pi_{n} (M \otimes_{A} N)$$
is an isomorphisms of abelian groups.
\end{corollary}

\begin{proof}
If $N$ is flat, then $\Tor_{p}^{\pi_{\ast} A}( \pi_{\ast} M, \pi_{\ast} N)$ vanishes
for $p > 0$, and is isomorphic to $\pi_{\ast} M \otimes_{ \pi_0 A} \pi_{0} N$ for $p = 0$. 
It follows that the spectral sequence of Proposition \ref{siiwe} degenerates at the $E_2$-page and yields the desired result.
\end{proof}

\begin{corollary}\label{huty}
Let $A$ be a connective $A_{\infty}$-ring, $M$ a connective right $A$-module, and $N$ a connective left $A$-module. Then:
\begin{itemize}
\item[$(1)$] The relative tensor product $M \otimes_{A} N$ is connective.
\item[$(2)$] There is a canonical isomorphism
$\pi_0( M \otimes_{A} N) \simeq \pi_0 M \otimes_{ \pi_0 A} \pi_0 N$
in the category of abelian groups.
\end{itemize}
\end{corollary}

\begin{proof}
This follows from the spectral sequence of Proposition \ref{siiwe}, since
$E_{2}^{p,q}$ vanishes for $p < 0$ or $q < 0$, while
$E_{2}^{0,0} \simeq \pi_0 M \otimes_{ \pi_0 A} \pi_0 N$.
\end{proof}

Our next goal is to prove an analogue of Lazard's theorem, which characterizes
the class of flat modules over a {\em connective} $A_{\infty}$-ring.

\begin{lemma}\label{sudsy}
Let $A$ be an $A_{\infty}$-ring.
\begin{itemize}
\item[$(1)$] The $A_{\infty}$-ring $A$ is flat, when viewed as a left module over itself.
\item[$(2)$] The collection of flat $A$-modules is stable under coproducts, retracts, and filtered colimits.
\item[$(3)$] Every free $A$-module is flat. If $A$ is connective, then every projective
$A$-module is flat. 
\end{itemize}
\end{lemma}

\begin{proof}
Assertions $(1)$ and $(2)$ are obvious, and $(3)$ follows from Proposition \stableref{charprojjj}.
\end{proof}

\begin{theorem}[Lazard's Theorem]\label{lazard}\index{Lazard's Theorem}
Let $A$ be a connective $A_{\infty}$-ring, and let $N$ be a connective left $A$-module.
The following conditions are equivalent:
\begin{itemize}

\item[$(1)$] The left $A$-module $N$ can be obtained as a filtered colimit of finitely generated
free modules.

\item[$(2)$] The left $A$-module $N$ can be obtained as a filtered colimit of
projective left $A$-modules.

\item[$(3)$] The left $A$-module $N$ is flat.

\item[$(4)$] The functor $M \mapsto M \otimes_{R} N$ is left t-exact; in other words,
it carries $(\Mod_{A}^{R})_{\leq 0}$ into $(\Spectra)_{\leq 0}$. 

\item[$(5)$] If $M$ is a discrete right $A$-module, then $M \otimes_{R} N$ is discrete.
\end{itemize}
\end{theorem}

\begin{proof}
The implications $(1) \Rightarrow (2) \Rightarrow (3)$ are obvious. 
The implication $(3) \Rightarrow (4)$ follows from Corollary \ref{siwe}, and
$(4) \Rightarrow (5)$ follows from Corollary \ref{huty}. 

We next show that $(5) \Rightarrow (3)$.
Suppose that $(5)$ is satisfied. The functor $M \mapsto M \otimes_{R} N$ is exact, and carries the heart of $\Mod_{A}^{R}$ into itself. It therefore induces an exact functor from
the abelian category of right $\pi_0 A$-modules to itself. According to Corollary \ref{huty}, this functor is given by (classical) tensor product with the left $\pi_0 A$-module $\pi_0 N$. From the exactness we conclude that $\pi_0 N$ is a flat $\pi_0 A$-module. 

We now prove that the natural map
$\phi: \pi_n A \otimes_{ \pi_0 A} \pi_0 N \rightarrow \pi_{n} N$ is an isomorphism
for all $n \in \Z$. For $n < 0$, this follows from the assumption that both $N$ and $A$ are connective. If $n=0$ there is nothing to prove. We may therefore assume that $n > 0$, and we work by induction on $n$. Let $M$ be the discrete right $A$-module corresponding to
$\pi_0 A$. The inductive hypothesis implies that $\Tor_{p}^{ \pi_{\ast} A}( \pi_{\ast} M, \pi_{\ast} N)_{q}$ vanishes unless $p=q=0$ or $q \geq n$. These $\Tor$-groups can be identified
with the $E_2$-terms of the spectral sequence of Proposition \ref{siiwe}, which computes
the homotopy groups of the discrete spectrum $M \otimes_{A} N$. Consequently, we have
$E^{0,0}_{\infty} = E^{0,0}_{2} \simeq \pi_0 N$. A simple calculation shows that
$E^{0,n}_{2} \simeq \coker(\phi)$, and that if $\phi$ is surjective then
$E^{1,n}_{2} \simeq \ker(\phi)$. To complete the proof, it suffices to prove that
$E^{i,n}_{2} \simeq \ast$ for $0 \leq i \leq 1$. To see this, we observe that the
vanishing of the groups $E^{i-r,n+r-1}_{2}$ and $E^{i+r,n-r+1}_{2}$ for $r \geq 2$
implies that $E^{i,n}_{2} \simeq E^{i,n}_{\infty}$, and the latter is a subquotient
of $\pi_{i+n}( M \otimes_{A} N)$, which vanishes in view of assumption $(5)$.

To complete the proof, it will suffice to show that $(3)$ implies
$(1)$. Let $\calC$ be the full subcategory of $\Mod_{A}$ spanned by
a set of representatives for all free, finitely generated $A$-modules. Then $\calC$ is small,
and consists of compact projective objects of $(\Mod_{A})_{\geq 0}$, which generate
$(\Mod_{A})_{\geq 0}$ under colimits. It follows (see \S \stableref{stable11}) that
the inclusion $\calC \subseteq \Mod_{A}$ induces an equivalence
$\calP_{\Sigma}(\calC) \rightarrow (\Mod_{A})_{\geq 0}$. Applying Lemma \toposref{longwait1}, we conclude that the identity functor from $(\Mod_{A})_{\geq 0}$ is a left Kan extension of
its restriction to $\calC$. It follows that for every connective left $A$-module $N$, the canonical
diagram $\calC_{/N} = \calC \times_{ \Mod_{A} } (\Mod_{A})_{/N} \rightarrow \Mod_{A}$ has $N$ as a colimit. To complete the proof, it will suffice to show that $\calC_{/N}$ is filtered
provided that $N$ is flat.

According to Proposition \toposref{stook}, it will suffice to verify the following conditions:
\begin{itemize}
\item[$(i)$] For every finite collection of objects $\{ X_i \}$ of $\calC_{/N}$, there
exists an object $X \in \calC_{/N}$ together with morphisms $X_i \rightarrow X$.

\item[$(ii)$] For every pair $X, Y \in \calC_{/N}$, every nonnegative integer $n \geq 0$, and every map $S^n \rightarrow \bHom_{\calC_{/N}}(X,Y)$ in the homotopy category $\calH$, there exists a morphism $Y \rightarrow Z$ in $\calC_{/N}$ such that the induced map $S^n \rightarrow \bHom_{\calC}(X,Z)$ is nullhomotopic.
\end{itemize}

Assertion $(i)$ follows immediately from the stability of $\calC$ under finite coproducts.
We now prove $(ii)$. Suppose given a pair of maps $f: X \rightarrow N$
and $g: Y \rightarrow N$, respectively. We have a homotopy fiber sequence
$$ \bHom_{\calC_{/N}}(f,g) \rightarrow \bHom_{\calC}(X,Y) \rightarrow
\bHom_{\Mod_{A}}( X, N ). $$
Since $\Mod_{A}$ is stable, $\bHom_{\calC_{/N}}(f,g)$ is a torsor over $\bHom_{\Mod_{A}}( X, \ker(g) )$. It follows that any map
$S^n \rightarrow \bHom_{\Mod_{A}}(X, \ker(g))$ determines a homotopy class
$\eta \in \Ext^{-n}_{\Mod_A}( X, \ker(g) )$. We wish to prove that
there exists a commutative diagram
$$ \xymatrix{ Y \ar[rr] \ar[dr]^{g} & & Z \ar[dl]^{h} \\
& N & }$$
such that the image of $\eta$ in $\Ext^{-n}_{\Mod_{A}}(X, \ker(h))$ vanishes.
Arguing iteratively, we can reduce to the case where $X \simeq A$, so that
$\eta$ can be identified with an element of $\pi_{n} \ker(g)$. 

Let $\eta' \in \pi_{n} Y$ denote the image of $\eta$. Our first step is to choose a diagram
as above with the property that the image of $\eta'$ in $\pi_{n} Z$ vanishes. 
We observe that $\eta'$ lies in the kernel of the natural map
$\pi_{n}(g): \pi_{n} Y \rightarrow \pi_{n} N \simeq 
\Tor_0^{\pi_0 A}(\pi_{n} A,\pi_0 N)$. The classical version of Lazard's theorem (see \cite{lazard}) implies that $\pi_0 N$ is isomorphic to a filtered colimit of free left $\pi_0 A$-modules. It follows that
there exists a commutative diagram
$$ \xymatrix{ & P \ar[dr]^{\overline{h}} & \\
\pi_0 Y \ar[rr]^{\pi_0 g} \ar[ur]^{\overline{k}}  & & \pi_0 N }$$
of left $\pi_0 A$-modules, where $P$ is a finitely generated free module, and the image of
$\eta'$ in $\pi_{n} A \otimes_{\pi_0 A} P$ vanishes. 
Using the freeness of $P$, we can realize $\overline{h}$ as $\pi_0 h$, where $h: Z \rightarrow N$
is a morphism of left $A$-modules, and $Z$ is a finitely generated free module with
$\pi_0 Z \simeq P$. Similarly, we can realize $\overline{k}$ as $\pi_0 k$, where
$k: Y \rightarrow Z$ is a morphism of left $A$-modules. Using the freeness of $Y$ again, we conclude that the diagram 
$$ \xymatrix{ & Z \ar[dr]^{h} & \\
Y \ar[rr]^{g} \ar[ur]^{k} & & N }$$
commutes in the homotopy category $\h \Mod_{A}$, and can therefore be lifted to a commutative triangle in $\Mod_{A}$. By construction, the image of $\eta'$ in $\pi_{n} Z$ vanishes.
Replacing $Y$ by $Z$, we may reduce to the case where $\eta' = 0$.

We now invoke the exactness of the sequence sequence
$$ \pi_{n+1} N \rightarrow \pi_{n} \ker(g) \rightarrow \pi_{n} Y$$
to conclude that $\eta$ is the image of a class $\eta'' \in \pi_{n+1} N
\simeq \pi_{n+1} A \otimes_{\pi_0 A} \pi_0 N$. Invoking Lazard's theorem once more, we
deduce the existence of a commutative diagram of left $\pi_0 A$-modules  
$$ \xymatrix{ & Q \ar[dr] & \\
\pi_0 Y \ar[rr]^{\pi_0 g} \ar[ur] & & \pi_0 N }$$
where $Q$ is a finitely generated free module, and $\eta''$ is the image of some element of
$\pi_{n+1} A \otimes_{ \pi_0 A} N$. Arguing as before, we may assume that the preceding diagram is induced by a commutative triangle of left $A$-modules
$$ \xymatrix{ & Z \ar[dr] & \\
Y \ar[rr]^{g} \ar[ur] & & N. }$$
Replacing $Y$ by $Z$, we may assume that $\eta''$ lies in the image of the map
$\pi_{n+1} Y \rightarrow \pi_{n+1} N$. The exactness of the sequence
$$ \pi_{n+1} Y \rightarrow \pi_{n+1} N \rightarrow \pi_{n} \ker(g)$$
now implies that $\eta = 0$, as desired. This completes the proof of the implication $(3) \Rightarrow (1)$.
\end{proof}

We now study the behavior of flatness under base change.

\begin{proposition}\label{urwise}
Let $f: A \rightarrow B$ be a map of $A_{\infty}$-rings, let
$G: \Mod_{B} \rightarrow \Mod_{A}$ be the forgetful functor, and let
$F: \Mod_{A} \rightarrow \Mod_{B}$ be a left adjoint to $G$
$($given by $M \mapsto B \otimes_{A} M$, in view of Corollary \ref{tara}$)$.
Then:
\begin{itemize}
\item[$(1)$] The functor $F$ carries free $($projective, flat$)$ $A$-modules to
free $($projective, flat$)$ $B$-modules.
\item[$(2)$] Suppose that $B$ is free $($projective, flat$)$ as a left $A$-module $($that is,
$G(B)$ is free, projective, or flat$)$. Then $G$ carries free $($projective, flat$)$ $B$-modules to
free $($projective, flat$)$ $A$-modules.
\item[$(3)$] Suppose that, for every $n \geq 0$, the map $f$ induces an isomorphism
$\pi_{n} A \rightarrow \pi_{n} B$. Then $F$ induces an equivalence of categories
$\Mod_{A}^{\flat} \rightarrow \Mod_{B}^{\flat}$; here $\Mod_{A}^{\flat}$ denotes the full subcategory of $\Mod_{A}$ spanned by the flat left $A$-modules, and $\Mod_{B}^{\flat}$ is defined likewise.
\end{itemize}
\end{proposition}

\begin{proof}
Assertions $(1)$ and $(2)$ are obvious. To prove $(3)$, we first choose $A'$ to be a connective cover of $A$ (see Proposition \ref{canconex}). We have a homotopy commutative triangle
of $\infty$-categories
$$ \xymatrix{ & \Mod_{A}^{\flat} \ar[dr] & \\
\Mod_{A'}^{\flat} \ar[ur] \ar[rr] & & \Mod_{B}^{\flat}. }$$
It therefore suffices to prove the analogous assertion for the morphisms $A' \rightarrow A$
and $A' \rightarrow B$. In other words, we may reduce to the case where $A$ is connective.

Since $A$ is connective, the $\infty$-category $\Mod_{A}$ admits a t-structure. Let
$F'$ denote the composite functor
$$ (\Mod_{A})_{\geq 0} \subseteq \Mod_{A} \stackrel{F}{\rightarrow} \Mod_{B}.$$
Then $F'$ has a right adjoint, given by the composition $\tau_{\geq 0} \circ G$. 
Assertion $(1)$ implies that $F'$ preserves flatness, and a simple calculation of
homotopy groups shows that $G'$ preserves flatness as well. Consequently, $F'$ and $G'$
induce adjoint functors
$$ \Adjoint{F''}{\Mod_{A}^{\flat}}{\Mod_{B}^{\flat}}{G''}.$$
It now suffices to show that the unit and counit of the adjunction are equivalences.
In other words, we must show:
\begin{itemize}
\item[$(i)$]
For every flat left $A$-module $M$, the unit map
$ M \rightarrow \tau_{\geq 0} B \otimes_{A} M$
is an equivalence. For this, it suffices to show that
$\pi_i M \simeq \pi_{i} \tau_{\geq 0} B \otimes_{A} M$ is an isomorphism for $i \in \Z$. If
$i < 0$, then both groups vanish, so there is nothing to prove. If $i \geq 0$, then we must show that $\pi_i M \simeq \pi_i (B \otimes_{A} M)$, which follows immediately from Corollary \ref{siwe} and the assumption that $\pi_i A \simeq \pi_i A$.
\item[$(ii)$] For every flat left $B$-module $N$, the counit map
$B \otimes_{A} \tau_{\geq 0} G(N) \rightarrow N$ is an equivalence. In other words,
we must show that for each $j \in \Z$, the map
$ \pi_{j} (B \otimes_A \tau_{\geq 0} G(N) ) \rightarrow \pi_{j} N$ is an isomorphism of abelian groups.
Since $G(N)$ is flat over $A$, Corollary \ref{siwe} implies that the left side is given
by $\pi_{j} B \otimes_{ \pi_0 A} \pi_0 N$. The desired result now follows immediately from our assumption that $N$ is flat.
\end{itemize}
\end{proof}

In general, if $M$ is a flat left
$A$-module, then ``global" properties of $M$ as an $A$-module
are often controlled by ``local" properties of $\pi_0 M$, viewed as a module
over the ordinary ring $\pi_0 A$. Our next pair of results illustrates this principle.

\begin{lemma}\label{sarefla}
Let $A$ be an $A_{\infty}$-ring, and let $f: M \rightarrow N$ be a map of flat left $A$-modules.
Then $f$ is an equivalence if and only if it induces an isomorphism $\pi_0 M \rightarrow \pi_0 N$.
\end{lemma}

\begin{proof}
This follows immediately from the definition of flatness.
\end{proof}

\begin{proposition}\label{redline}
Let $A$ be a connective $A_{\infty}$-ring. A flat left
$A$-module $M$ is projective if and only if $\pi_0 M$ is a
projective module over $\pi_0 A$.
\end{proposition}

\begin{proof}
Suppose first that $\pi_0 M$ is freely generated by elements $\{ \eta_i \}_{i \in I}$. Let
$P = \coprod_{i \in I} A$, and let $f: P \rightarrow M$ be a map represented by 
$\{ \eta_i \}_{i \in I}$. By construction, $f$ induces an isomorphism $\pi_0 P \rightarrow \pi_0 M$. 
Lemma \ref{sarefla} implies that $f$ is an equivalence, so that $M$ is free.

In the general case, there exists a free $\pi_0 A$-module $F_0$ and a direct sum decomposition
$F_0 \simeq N_0 \oplus \pi_0 M$. Replacing $F_0$ by $\oplus_{ n \geq 0} F_0$ if necessary, we may assume that $N_0$ is itself free. The projection map $F_0 \rightarrow \pi_0 M$ is induced by
a map $g: F \rightarrow M$ of left $A$-modules, where $F$ is free. Then $g$ induces a
surjection $\pi_0 F \simeq F_0 \rightarrow \pi_0 M$. Using the flatness of $M$, we conclude that the maps $$\pi_i F \simeq \pi_i A \otimes_{ \pi_0 A} \pi_{0} F
\rightarrow \pi_i A \otimes_{ \pi_0 A} \pi_0 M \simeq \pi_i M$$
are also surjective. Let $N$ be a kernel of $g$, so that we have a commutative diagram
$$ \xymatrix{ 0 \ar[r] & \pi_i A \otimes_{ \pi_0 A} \pi_0 N
\ar[r] \ar[d]^{\phi'} & \pi_{i} A \otimes_{ \pi_0 A} \pi_{0} F \ar[r] \ar[d]^{\phi} & \pi_i A \otimes_{ \pi_0 A} \pi_0 M \ar[r] \ar[d]^{\phi''} & 0 \\
 0 \ar[r] & \pi_i N \ar[r] & \pi_i F \ar[r] & \pi_i M \ar[r] & 0. }$$
 Using the flatness of $F$ and $M$, we deduce that the upper row is short exact, and that the maps
 $\phi$ and $\phi''$ are isomorphisms. The snake lemma implies that $\phi$ is an isomorphism; moreover, $\pi_0 N$ is isomorphic to the kernel of a surjection between flat $\pi_0 A$-modules, and is therefore itself flat. It follows that $N$ is a flat $A$-module. Since $\pi_0 N$ is free, the first part of the proof shows that $N$ is itself free.
  
Let $p: N \rightarrow F$ denote the natural map. Since $\pi_0 M$ is projective, the inclusion
$\pi_0 N \subseteq \pi_0 F$ is split. Since $F$ is free, we can lift this splitting to a morphism $q: F \rightarrow N$. Then $q \circ p: N \rightarrow N$ induces the identity map from
$\pi_0 N$ to itself. Since $N$ is free, we conclude that there is a homotopy $q \circ p \simeq \id_{N}$. It follows that $N$ is a direct summand of $F$ in the homotopy category $\h{\calC}$. 
Consequently, $M \simeq \coker(p)$ can be identified with the complementary summand, and is therefore projective.
\end{proof}

\subsection{Finiteness Properties of Rings and Modules}\label{finprop}

In this section, we will discuss some finiteness conditions on the $\infty$-categories of modules over an $A_{\infty}$-ring $R$. We begin by introducing the definition of a {\em perfect} $R$-module. Roughly speaking, an $R$-module $M$ is perfect if it can be obtained by gluing together
finitely many (possibly shifted) copies of $R$, or is a retract of such an $R$-module. Alternatively, we can describe the class of perfect $R$-modules as the compact objects of the $\infty$-category $\Mod_{R}$ (Proposition \ref{fergus}).

In general, the condition that an $R$-module be perfect is very strong. For example, if $R$ is a discrete commutative ring and $M$ is a finitely generated (discrete) module over $R$, then $M$ need not be perfect when viewed as an object of the stable $\infty$-category $\Mod_{R}$. In general, this is true if and only if $R$ is a regular Noetherian ring. To remedy the situation, we will introduce the weaker notion of an {\it almost perfect} $R$-module. This notion has a much closer relationship with finiteness conditions in the classical theory of rings and modules. In particular, we will show that a left $R$-module $M$ is almost perfect if and only if
each $\pi_i M$ is a finitely presented left module over $\pi_0 R$, and $\pi_{i} M \simeq 0$ for $i \ll 0$ (Proposition \ref{bbh}), at least when $R$ itself satisfies a suitable finiteness condition (that is, when $R$ is {\it left coherent} in the sense of Definition \ref{umbul}).

\begin{definition}\label{perfet}\index{perfect module}\index{module object!perfect}
Let $A$ be an $A_{\infty}$-ring. We let $\Mod_{A}^{L, \perf}$ denote the
smallest stable subcategory of $\Mod_{A}^{L}$ which contains $A$ (regarded as a left module over itself) and is closed under retracts. Similarly, we let $\Mod_{A}^{R, \perf}$ denote the smallest stable subcategory of
$\Mod_{A}^{R}$ which contains $A$ and is closed under retracts.
We will say that a left (right) $A$-module $M$ is {\it perfect} if it belongs to
$\Mod_{A}^{L, \perf}$ ($\Mod_{A}^{R, \perf}$).\index{ZZZModAperf@$\Mod_{A}^{\perf}$}
\end{definition}

\begin{proposition}\label{fergus}
Let $A$ be an $A_{\infty}$-ring. Then:
\begin{itemize}
\item[$(1)$] The $\infty$-category $\Mod_{A}^{L}$ is compactly generated.
\item[$(2)$] An object of $\Mod_{A}^{L}$ is compact if and only if it is perfect.
\end{itemize}
\end{proposition}

\begin{proof}
The compact objects of $\Mod_{A}$ form a stable subcategory which is closed under the formation of retracts. Moreover, Proposition \ref{pretara} implies $A \in \Mod^{L}_{A}$ corepresents the composition $\Mod_{A} \rightarrow \Spectra \stackrel{\Omega^{\infty}}{\rightarrow} \SSet,$
where both maps preserve filtered colimits. It follows that $A$ is compact as an $A$-module. This proves that every perfect object of $\Mod_{A}^L$ is compact.

According to Proposition \toposref{uterr}, the inclusion $f: \Mod_{A}^{L, \perf} \subseteq \Mod_{A}^{L}$ induces a fully faithful functor $F: \Ind( \Mod_{A}^{L, \perf} ) \rightarrow \Mod_{A}^{L}$. To complete the proof, it will suffice to show that $F$ is essentially surjective.
Since $f$ is right exact, $F$ preserves all colimits (Proposition \toposref{sumatch}), so the essential image of $F$ is stable under colimits. If $F$ is not essentially surjective, then Proposition \stableref{kura} implies that there exists a nonzero $N \in \Mod_{A}^{L}$ such that $\bHom_{ \Mod_{A}^{L}}(N', N) \simeq \ast$ for all $N'$ belonging to the essential image of $F$. In particular, taking
$N' = A[n]$, we conclude that $\pi_{n} N \simeq \ast$. It follows that $N$ is a zero object
of $\Mod_{A}^{L}$, contrary to our assumption. 
\end{proof}

\begin{proposition}\label{funker}
Let $A$ be an $A_{\infty}$-ring. The relative tensor product functor
$$ \otimes_{A}: \Mod^{R}_A \times \Mod^{L}_{A} \rightarrow \Spectra$$
induces fully faithful embeddings
$$ \theta: \Mod^{R}_{A} \rightarrow \Fun( \Mod^{L}_{A}, \Spectra) \quad  \theta': \Mod^{L}_{A} \rightarrow \Fun( \Mod^{R}_{A}, \Spectra).$$
A functor $f: \Mod^{L}_{A} \rightarrow \Spectra$ $(${}$f'': \Mod^{R}_{A} \rightarrow \Spectra${}$)$ belongs to the essential image of $\theta$ $(${}$\theta'${}$)$ if and only if $f$ $(${}$f'${}$)$ preserves small colimits.
\end{proposition}

\begin{proof}
Let $\calC$ be the full subcategory of $\Fun( \Mod^{L}_{A}, \Spectra)$ spanned by those functors which preserve colimits. Proposition \toposref{intmap} implies that $\calC$ is presentable, and
Proposition \ref{colimp} implies that $\theta$ factors through $\calC$.
We will show that $\theta$ induces an equivalence $G: \Mod^{R}_{A} \rightarrow \calC$; the analogous assertion for left modules follows by the same argument.

Proposition \ref{fergus} implies that $\Mod^{L}_{A}$ is equivalent to the $\infty$-category $\Ind( \Mod^{L,\perf}_{A} )$. It follows from Propositions \toposref{uterrr}, \toposref{sumatch} and \stableref{funrose} that $\calC$ is equivalent to the $\infty$-category $\Fun^{\Ex}( \Mod^{L, \perf}_{A}, \Spectra)$ of exact functors from $\Mod^{L, \perf}_{A}$ to spectra. In particular, for every perfect left $A$-module $N$, evaluation on $\calN$ induces a functor $\calC \rightarrow \Spectra$ that preserves all small limits and colimits. 

Let $G': \calC \rightarrow \Spectra$ be given by evaluation on $A$, regarded as a (perfect) left module over itself. If $\alpha: f \rightarrow f'$ is a natural transformation of functors from $\Mod^{L, \perf}_{A}$ to $\Spectra$, then the full subcategory of $\h{\Mod^{L, \perf}_{A}}$ spanned by objects $C$ such that $\alpha(C): f(C) \rightarrow f'(C)$ is an equivalence is a triangulated subcategory
of $\h{ \Mod^{L,\perf}_{A}}$ which is stable under retracts. If follows that $G'$ is {\em conservative}: if $G'(\alpha)$ is an equivalence, then $\alpha$ is an equivalence. Since $G'$ preserves small limits, we deduce also that $G'$ detects small limits: if $\overline{p}: K^{\triangleleft} \rightarrow \calC$ is such that $G' \circ \overline{p}$ is a limit diagram in $\Spectra$, then $\overline{p}$ is a limit diagram in $\calC$. Similarly, $G'$ detects small colimits.

The composite functor $G' \circ G: \Mod^{R}_{A} \rightarrow \Spectra$ can be identified with the forgetful functor, in view of Proposition \ref{pretara}. It follows that $G' \circ G$ preserves all limits and colimits (Corollaries \ref{goop} and \ref{gloop}). Since $G'$ detects small limits and colimits, we deduce that $G$ preserves small limits and colimits. Corollary \toposref{adjointfunctor} implies that $G$ and $G'$ admit left adjoints, which we will denote by $F$ and $F'$.

Choose unit and counit transformations 
$$u: \id \rightarrow G \circ F, \quad v: F \circ G \rightarrow \id.$$
We wish to prove that $u$ and $v$ are equivalences. Since $G' \circ G$ detects equivalences,
the functor $G$ detects equivalences, so that $v$ is an equivalence if and only if
$G(v): G \circ F \circ G \rightarrow G$ is an equivalence. Since $G(v)$ has a section
determined by $u$, it will suffice to prove that $u$ is an equivalence.

Let $\calC'$ be the full subcategory of $\calC$ spanned by those objects
$f \in \calC$ for which the map $u(f): f \rightarrow (G \circ F)(f)$ is an equivalence. Since $G$ and $F$ preserve small colimits, we deduce that $\calC'$ is stable under shifts and colimits in $\calC$.
Proposition \ref{littlebeck} implies that $\calC$ is generated, under geometric realizations, by the essential image of $F'$. Let $\Spectraprime \subseteq \Spectra$ be the inverse image of $\calC'$ under $F'$. Since $F'$ is a colimit-preserving functor, we deduce that $\Spectraprime \subseteq \Spectra$ is closed under shifts and colimits. Since $\Spectra$ is generated under colimits by the objects $\Sphere[n]$, where $n \in \Z$ and $\Sphere \in \Spectra$ denotes the sphere spectrum, it will suffice to show that $\Sphere \in \Spectraprime$. 

Proposition \ref{pretara} allows us to identify $(F \circ F')(\Sphere)$ with $\Sphere \otimes A \simeq A$, regarded as a left module over itself. Proposition \ref{usss} implies that $(G \circ F)(F'(\Sphere))$ can be identified with the forgetful functor $f_0: \Mod_{A}^{L} \rightarrow \Spectra$. We are reduced to proving that the unit map
$\Sphere \rightarrow A \simeq f_0(A)$ induces an equivalence
$F'(\Sphere) \simeq f_0$ in $\calC$. 

Applying Proposition \stableref{urtusk}, we can identify $\calC$ with the full subcategory $\calD \subseteq \Fun( \Mod^{L, \perf}_{A}, \SSet) = \calP( \Mod^{L, \perf, op}_{A})$ spanned by those functors which preserve finite limits. Under this equivalence, $f_0$ corresponds to the composition
$ \Mod^{L, \perf}_{A} \subseteq \Mod^{L}_{A} \rightarrow \Spectra \stackrel{\Omega^{\infty}}{\rightarrow} \SSet, $
while $F'(\Sphere)$ corresponds to the image of $\ast \in \SSet$ under the composition
$\SSet \stackrel{\Sigma^{\infty}}{\rightarrow} \Spectra \stackrel{F'}{\rightarrow} \calC \simeq \calD$ which is the left adjoint to the the functor $\calD \rightarrow \SSet$ given by evaluation
at $A \in \Mod^{L, \perf}_{A}$. To complete the proof, it will suffice to show that the unit
$1 \in \pi_0 A$ exhibits the composite functor
$ \Mod^{L}_{A} \rightarrow \Spectra \stackrel{\Omega^{\infty}}{\rightarrow} \SSet$
as corepresented by $A \in \Mod^{L}_{A}$. In other words, we must show that for every
$M \in \Mod^{L}_{A}$, the canonical map
$\bHom_{\Mod^{L}_{A}}(A, M) \rightarrow \Omega^{\infty} M$
is a homotopy equivalence. Using Proposition \ref{pretara}, we can reduce to the case where
$A$ is the unit object of $\Alg(\Spectra)$. In view of Corollary \ref{puterry}, we are reduced to proving that if $M \in \Spectra$, then the canonical map
$ \bHom_{\Spectra}( \Sphere, M) \rightarrow \Omega^{\infty} M$
is an equivalence, which is clear (since $\Sphere \simeq \Sigma^{\infty}(\ast)$ by definition).
\end{proof}

\begin{remark}
Proposition \ref{funker} admits the following generalization. Suppose that $A$ and $B$ are $A_{\infty}$-rings. Then the $\infty$-category of colimit preserving functors from $\Mod^{L}_{A}$ to $\Mod^{L}_{B}$ is equivalent to an $\infty$-category of {\it $A$-$B$-bimodules}. A precise formulation (and proof) of this statement require some elaboration on the ideas presented in this section. For a proof in the language of model categories, we refer the reader to \cite{schwedeshipley}. 
\end{remark}

\begin{proposition}\label{ferfet}
Let $A$ be an $A_{\infty}$-ring, and let $M$ be a left $A$-module. The following conditions are equivalent:

\begin{itemize}
\item[$(1)$] The left $A$-module $M$ is perfect.
\item[$(2)$] The left $A$-module $M$ is a compact object of
$\Mod_{A}^{L}$.
\item[$(3)$] There exists a right $A$-module $M^{\vee}$ such that
the composition
$\Mod_{A}^{L} \stackrel{ M^{\vee} \otimes_{A} \bigdot}{\rightarrow}
\Spectra \stackrel{\Omega^{\infty}}{\rightarrow} \SSet$
is equivalent to the functor corepresented by $M$.
\end{itemize}
Moreover, if these conditions are satisfied, then the object
$M^{\vee} \in \Mod_A^{R}$ is also perfect.
\end{proposition}

\begin{proof}
The equivalence $(1) \Leftrightarrow (2)$ is a consequence of Proposition \ref{fergus}. 
Let $\calC$ denote the full subcategory of $\Fun( \Mod_{A}^{L}, \Spectra)$ spanned by
those functors which are continuous and exact. Proposition \ref{funker} yields an equivalence of $\infty$-categories $\Mod_{A}^{R} \rightarrow \calC$. 
According to Proposition \stableref{urtusk}, composition with $\Omega^{\infty}$ induces a fully faithful embedding $\calC \rightarrow \Fun( \Mod_{A}^{L}, \SSet)$, whose essential image
$\calC'$ consists precisely of those functors which are continuous and left exact. The functor
co-represented by $M$ is automatically left-exact, and is continuous if and only if $M$ is compact. This proves that $(2) \Leftrightarrow (3)$.

Let $j: (\Mod_{A}^{L})^{op} \rightarrow \Fun( \Mod_{A}^{L}, \SSet)$ denote the dual Yoneda embedding, so that $j$ restricts to a map $j': (\Mod_{A}^{L, \perf})^{op} \rightarrow \calC'$.
Composing $j'$ with a homotopy inverse to the equivalence
$ \Mod_{A}^{R} \rightarrow \calC \rightarrow \calC',$
we obtain a ``dualization'' map $(\Mod_{A}^{L, \perf})^{op} \rightarrow \Mod_{A}^{R}$, which we will indicate by $M \mapsto M^{\vee}$. Let
$\calD \subseteq \Mod_{A}^{L, \perf}$ be the full subcategory spanned by those objects
$M$ such that $M^{\vee}$ is perfect. We wish to show that $\calD = \Mod_{A}^{L, \perf}$. The functor $M \mapsto M^{\vee}$ is exact, and $\Mod_{A}^{R, \perf}$ is a stable subcategory
of $\Mod_{A}^{R}$ which is closed under retracts. It follows that $\calD$ is a stable subcategory of $\Mod_{A}^{L}$ which is closed under retracts. It now suffices to observe that $A^{\vee} \simeq A$, so that $A \in \calD$.
\end{proof}

\begin{corollary}\label{whumph}
Let $A$ be a connective $A_{\infty}$-ring, and let $M$ be a perfect left $A$-module.
Then:
\begin{itemize}
\item[$(1)$] The homotopy groups $\pi_{m} M$ vanish for $m \ll 0$.
\item[$(2)$] Suppose that $\pi_{m} M \simeq \ast$ for all $m < k$. Then
$\pi_{k} M$ is a finitely presented module over $\pi_0 A$.
\end{itemize}
\end{corollary}

\begin{proof}
We have an equivalence $M \simeq \varinjlim( \tau_{\geq -n} M)$. Since
$M$ is compact (Proposition \ref{fergus}), we conclude that 
$M$ is a retract of some $\tau_{\geq -n} M$, so that $\pi_{m} M \simeq \ast$ for
$m < -n$. This proves $(1)$. To prove $(2)$, we observe that the inclusion
of the heart of $\Mod_{A}$ into $(\Mod_{A})_{\geq 0}$ preserves filtered colimits, so
the right adjoint $\tau_{\leq 0}: (\Mod_A)_{\geq 0} \rightarrow \heart{\Mod_{A}}$
preserves compact objects. It now suffices to observe that the compact objects in
the ordinary category of $\pi_0 A$-modules are precisely the finitely presented modules.
\end{proof}

In the situation of Proposition \ref{ferfet}, we will refer to $M^{\vee}$ as a {\it dual} of $M$. We next prove that this notion of duality determines an antiequivalence between the $\infty$-categories $\Mod_{A}^{L, \perf}$ and $\Mod_{A}^{R, \perf}$. 

\begin{lemma}\label{dudername}
Let $\calC$ and $\calD$ be $\infty$-categories, and let
$F: \calC \times \calD \rightarrow \SSet$ be a bifunctor. The following
conditions are equivalent:
\begin{itemize}
\item[$(1)$] The induced map $f: \calC \rightarrow \Fun(\calD, \SSet) = \calP(\calD^{op})$ is fully faithful, and the essential image of $f$ coincides with the essential image of the Yoneda embedding $\calD^{op} \rightarrow \calP(\calD^{op})$. 
\item[$(2)$] The induced map $f': \calD \rightarrow \Fun(\calC, \SSet) = \calP(\calC^{op})$ is fully faithful, and the essential image of $f$ coincides with the essential image of the Yoneda embedding $\calC^{op} \rightarrow \calP(\calC^{op})$. 
\end{itemize}
\end{lemma}

\begin{proof}
Let
$$ G_{\calC}: \calC \times \calC^{op} \rightarrow \SSet, \quad G_{\calD}: \calD^{op} \times \calD \rightarrow \SSet$$
be the maps used in the definition of the Yoneda embeddings (see \S \toposref{presheaf1}). Then $(1)$ is equivalent to the existence of an equivalence $\alpha: \calD^{op} \rightarrow \calC$
such that the composition
$\calD^{op} \times \calD \rightarrow \calC \times \calD \stackrel{F}{\rightarrow} \SSet$
is homotopic to $G_{\calD}$. If $\beta$ is a homotopy inverse to $\alpha$, then the
composition
$\calC \times \calC^{op} \stackrel{\beta^{op}}{\rightarrow} \calC \times \calD
\stackrel{F}{\rightarrow} \SSet$
is homotopic to $G_{\calC}$, which proves $(2)$. The converse follows by the same argument.
\end{proof}

We will say that a functor $F: \calC \times \calD \rightarrow \SSet$ is a {\it perfect pairing} if it satisfies the hypotheses of Lemma \ref{dudername}. In this case, $F$ determines an equivalence between $\calC$ and $\calD^{op}$, well-defined up to homotopy. The proof of Proposition \ref{ferfet} yields the following:\index{perfect pairing}

\begin{proposition}\label{munk}
Let $A$ be an $A_{\infty}$-ring. Then the bifunctor
$$ \Mod^{R,\perf}_{A} \times \Mod^{L, \perf}_{A}
\stackrel{\otimes_{A}}{\rightarrow} \Spectra \stackrel{\Omega^{\infty}}{\rightarrow} \SSet $$
is a perfect pairing.
\end{proposition}

\begin{remark}
Let $A$ be an $A_{\infty}$-ring. It follows from Proposition \ref{munk} that the $\infty$-category of right $A$-modules is formally determined by the $\infty$-category of left $A$-modules. Namely, $\Mod_{A}^{R}$ is equivalent to the $\infty$-category of $\Ind$-objects of $(\Mod_{A}^{L, \perf})^{op}$, where
$\Mod_{A}^{L, \perf}$ is the $\infty$-category of compact objects of $\Mod_{A}^{L}$.
\end{remark}

Let $A$ be an $A_{\infty}$-ring, and let $M$ be a left $A$-module.
Roughly speaking, $M$ is perfect if it can be obtained from the zero module through a finite process of attaching copies of $A$ and forming retracts. This is a very strong condition which is often violated in practice. For example, suppose that $A$ is a commutative Noetherian ring (viewed as a discrete $A_{\infty}$-ring), and that $M$ is a discrete $A$-module, such that
$\pi_0 M$ is finitely generated over $A$ in the sense of classical commutative algebra. In this case, we can choose a resolution
$$ \ldots \rightarrow P_2 \rightarrow P_1 \rightarrow P_0 \rightarrow \pi_0 M$$
where each $P_i$ is a free $A$-module of finite rank. However, we cannot usually guarantee that $P_n \simeq 0$ for $n \gg 0$; in general this is possible only when $A$ is regular (\cite{eisenbud}). The module $M$ is not generally perfect as an $A$-module spectrum (see Example \ref{tukus} below). Nevertheless, the fact that we can choose each $P_i$ to be of finite rank imposes a strong finiteness condition on $M$, which we now formulate.

\begin{definition}\index{compact!to order $n$}
Let $\calC$ be a compactly generated $\infty$-category. We will say that an object $C \in \calC$ is {\it compact to order $n$} if
$\tau_{\leq n}(C)$ is a compact object of $\tau_{\leq n} \calC$ (see \S \toposref{truncintro} for an explanation of this notation). We will say that
$C \in \calC$ is {\it almost compact} if it is compact to order $n$ for all $n \geq 0$.
\end{definition}

\begin{lemma}
Let $\calC$ be a compactly generated $\infty$-category, let $C$ be an object of $\calC$, and let $m \leq n$ be integers. If $C$ is compact to order $n$, then $C$ is compact to order $m$. If $C$ is compact, then it is compact to order $n$ for all integers $n$.
\end{lemma}

\begin{proof}
Corollary \toposref{hunterygreen} implies that the $\infty$-categories $\tau_{\leq m} \calC$ and
$\tau_{\leq n} \calC$ are compactly generated and stable under filtered colimits in $\calC$, so that
$\tau_{\leq m} \calC$ is stable under filtered colimits in $\tau_{\leq n} \calC$. The desired conclusion now follows from Proposition \toposref{comppress}.
\end{proof}

\begin{definition}\index{almost perfect module}\index{module object!almost perfect}
Let $A$ be a connective $A_{\infty}$-ring. We will say that a left $A$-module $M$ is {\it almost perfect} if there exists an integer $k$ such that for all $n \geq 0$, $M \in (\Mod_A)_{\geq k}$ and
is compact to order $n$ as an object of $(\Mod_{A})_{\geq k}$. 
We let $\Mod_{A}^{\aperf}$ denote the full subcategory of $\Mod_{A}$ spanned by the almost perfect left $A$-modules.\index{ZZZModAaperf@$\Mod_{A}^{\aperf}$}
\end{definition}

\begin{proposition}\label{almor}
Let $A$ be a connective $A_{\infty}$-ring. Then:
\begin{itemize}
\item[$(1)$] The full subcategory $\Mod_{A}^{\aperf} \subseteq
\Mod_{A}^{L}$ is closed under 
translations and finite colimits, and is therefore a stable subcategory of $\Mod_{A}$.
\item[$(2)$] The full subcategory $\Mod_{A}^{\aperf} \subseteq
\Mod_{A}$ is closed under the formation of retracts.
\item[$(3)$] Every perfect left $A$-module is almost perfect.
\item[$(4)$] The full subcategory $( \Mod_{A}^{\aperf})_{\geq 0} \subseteq \Mod_{A}$ is closed under the formation of geometric realizations of simplicial objects.
\item[$(5)$] Let $M$ be a left $A$-module which is connective and almost perfect.
Then $M$ can be obtained as the geometric realization of a simplicial left $A$-module
$P_{\bigdot}$ such that each $P_{n}$ is a free $A$-module of finite rank.
\end{itemize}
\end{proposition}

\begin{proof}
Assertions $(1)$ and $(2)$ are obvious, and $(3)$ follows from Proposition \ref{fergus}.
To prove $(4)$, it suffices to show that the collection of compact objects of
$\tau_{\leq n} (\Mod_{A})_{\geq 0}$ is closed under geometric realizations, which follows
from Lemma \stableref{simpenough}.

We now prove $(5)$. In view of Theorem \stableref{dkan} and Remark \stableref{saltine}, it will suffice to show that $M$ can be obtained as the colimit of a sequence
$$ D(0) \stackrel{f_1}{\rightarrow} D(1) \stackrel{f_2}{\rightarrow} D(2) \rightarrow \ldots$$
where each $\coker(f_n)[-n]$ is a free $A$-module of finite rank; here we agree by convention that
$f_0$ denotes the zero map $0 \rightarrow D(0)$. The construction goes by induction.
Suppose that the diagram
$$ D(0) \rightarrow \ldots \rightarrow D(n) \stackrel{g}{\rightarrow} M$$
has already been constructed, and that $N= \ker(g)$ is $(n-1)$-connected. 
Part $(1)$ implies that $N$ is almost perfect, so that the bottom homotopy group
$\pi_n N$ is a compact object in the category of left $\pi_0 A$-modules. It follows that
there exists a map $\beta: Q[n] \rightarrow N$, where $Q$ is a free left $A$-module of finite rank, and $\beta$ induces a surjection $\pi_0 Q \rightarrow \pi_n N$. We now define
$D(n+1) = \coker(\beta)$, and construct a diagram
$$ D(0) \rightarrow \ldots \rightarrow D(n) \rightarrow D(n+1) \stackrel{g'}{\rightarrow} M.$$
Using the octahedral axiom, we obtain a distinguished triangle
$$ Q[n] \rightarrow \ker(g) \rightarrow \ker(g') \rightarrow Q[n+1],$$
and the associated long exact sequence of homotopy groups proves that $\ker(g')$ is $n$-connected. 

In particular, we conclude that for fixed $m$, the maps $\pi_{m} D(n) \rightarrow \pi_m M$ are isomorphisms for $n \gg 0$, so that the natural map $\varinjlim D(n) \rightarrow M$ is
an equivalence of left $A$-modules, as desired.
\end{proof}

Using Proposition \ref{almor}, we can give the following characterization of the $\infty$-category of (connective) almost perfect modules over a connective $A_{\infty}$-ring.

\begin{corollary}
Let $A$ be a connective $A_{\infty}$-ring, let $\calC$ denote the full subcategory
of $\Mod_{A}$ spanned by those connective left $A$-modules which are connective and almost perfect, and let $\calC^{0} \subseteq \calC$ denote the full subcategory of $\calC$ spanned by the objects $\{ A^n \}_{n \geq 0}$. Let $\calD$ be an arbitrary $\infty$-category which admits geometric realizations for simplicial objects, and let $\Fun_{\sigma}(\calC, \calD)$ be the full subcategory
of $\Fun(\calC, \calD)$ spanned by those functors which preserve geometric realizations of simplicial objects. Then the restriction functor $\Fun_{\sigma}(\calC, \calD) \rightarrow \Fun( \calC^{0}, \calD)$ is an equivalence of $\infty$-categories.
\end{corollary}

\begin{proof}
Let $\calC'$ be the smallest full subcategory of $\calP(\calC^{0})$ which contains the essential image of the Yoneda embedding and is stable under geometric realizations of simplicial objects, and let $j: \calC^0 \rightarrow \calC'$ be the Yoneda embedding.
Using Remark \toposref{poweryoga}, we conclude that composition with $j$ induces an
equivalence $\Fun_{\sigma}(\calC', \calD) \rightarrow \Fun(\calC^{0}, \calD)$ for any
$\infty$-category $\calD$ which admits geometric realizations of simplicial objects. In particular,
the inclusion $\calC^{0} \subseteq \calC$ extends (up to homotopy) to a functor
$F: \calC' \rightarrow \calC$ which commutes with geometric realizations. To complete the proof, it will suffice to show that $F$ is an equivalence of $\infty$-categories.
Using the fact that
each $A^{n}$ is a projective object of $(\Mod_{A})_{\geq 0}$, we deduce that $F$ is fully faithful. Part $(5)$ of Proposition \ref{almor} implies that $F$ is essentially surjective.
\end{proof}

For a general connective $A_{\infty}$-ring $A$, the t-structure on
$\Mod_{A}$ does not restrict to a t-structure on the full subcategory
$\Mod_{A}^{\aperf}$. Roughly speaking, one would expect the
heart of $\Mod_{A}^{\aperf}$ to be equivalent to the ordinary category of finitely presented $\pi_0 A$-modules. In general, this is not an abelian category. We can correct this defect by introducing an appropriate hypothesis on $A$. We begin by recalling a definition from classical algebra.

\begin{definition}\label{umbul}\index{coherent!ring}\index{left coherent!ring}
An associative ring $R$ is {\it left coherent} if every finitely generated left ideal of $R$ is finitely presented (as a left $R$-module).
\end{definition}

\begin{example}
Every left Noetherian ring is left coherent.
\end{example}

For completeness, we include a proof of the following classical result:

\begin{lemma}\label{swirly}
Let $R$ be a left coherent ring. Then:
\begin{itemize}
\item[$(1)$] Every finitely generated submodule of $R^n$ is finitely presented.
\item[$(2)$] Every finitely generated submodule of a finitely presented left $R$-module is finitely presented.
\item[$(3)$] If $f: M \rightarrow N$ is a map of finitely presented left $R$-modules, then
$\ker(f)$ and $\coker(f)$ are finitely presented.
\end{itemize}
\end{lemma}

\begin{proof}
We first make the following elementary observations, which do not require the assumption that $R$ is left coherent:
\begin{itemize}
\item[$(a)$] Suppose $f: M \rightarrow N$ is an epimorphism of left $R$-modules. If
$M$ is finitely generated and $N$ is finitely presented, then $\ker(f)$ is finitely generated.
\item[$(b)$] Let $0 \rightarrow M' \rightarrow M \rightarrow M'' \rightarrow 0$ be a short exact sequence of left $R$-modules. If $M'$ and $M''$ are finitely presented, then $M$ is finitely presented.
\end{itemize}

We now prove $(1)$ using induction on $n$. When $n=0$ there is nothing to prove. Suppose that $n > 0$ and that $M \subseteq R^n$ is finitely generated. Form a diagram 
$$ \xymatrix{ 0 \ar[r] \ar[d] & M' \ar[r] \ar[d] & M \ar[r] \ar[d] & M'' \ar[r] \ar[d] & 0 \ar[d] \\
0 \ar[r] & R^{n-1} \ar[r] & R^n \ar[r] & R \ar[r] & 0 }$$
where the vertical maps are monomorphisms. Then $M''$ can be identified with a finitely generated left ideal of $R$. Since $R$ is left coherent, we conclude that $M''$ is finitely presented.
Using $(a)$, we deduce that $M'$ is itself finitely generated. The inductive hypothesis now implies that $M'$ is finitely presented, so that we can use $(b)$ to conclude that $M$ is finitely presented.

We next prove $(2)$. Suppose that $f: M \rightarrow N$ is a monomorphism, where
$N$ is finitely presented and $M$ is finitely generated. Choose an epimorphism
$g: R^n \rightarrow N$, and form a pullback diagram
$$ \xymatrix{ K \ar[r] \ar[d] & R^n \ar[d] \\
M \ar[r] & N. }$$
Then $K$ can be identified with the kernel of the induced map $R^n \rightarrow N/M$, and 
is therefore finitely generated. Part $(1)$ implies that $K$ is finitely presented. The induced map
$K \rightarrow M$ is an epimorphism, whose kernel is isomorphic to $\ker(g)$ and is therefore finitely generated. It follows that $M$ is finitely presented, as desired.

We now prove $(3)$. It is clear that $\coker(f)$ is finitely presented (this does not require the left coherence of $R$). We next show that $\ker(f)$ is finitely presented. The image of
$f$ is a finitely generated submodule of $N$, and therefore finitely presented by $(2)$. 
Consequently, we may replace $N$ by $\im(f)$, and thereby reduce to the case where $f$
is an epimorphism. We now apply $(a)$ to deduce that $\ker(f)$ is finitely generated.
Invoking $(2)$ again, we conclude that $\ker(f)$ is finitely presented as desired.
\end{proof}

\begin{definition}\index{left coherent!$A_{\infty}$-ring}\index{coherent!$A_{\infty}$-ring}\label{pokus}
Let $A$ be an $A_{\infty}$-ring. We will say that $A$ is {\it left coherent} if
the following conditions are satisfied:
\begin{itemize}
\item[$(1)$] The $A_{\infty}$-ring $A$ is connective.
\item[$(2)$] The associative ring $\pi_0 A$ is left coherent (in the sense of Definition \ref{umbul}).
\item[$(3)$] For each $n \geq 0$, the homotopy group
$\pi_n A$ is finitely presented as a left module over $\pi_0 A$.
\end{itemize}
\end{definition}

\begin{proposition}\label{bbh}
Let $A$ be an $A_{\infty}$-ring and $M$ an $A$-module. 
Suppose that $A$ is left coherent. Then $M$ is almost perfect if and only if the following conditions are satisfied:
\begin{itemize}
\item[$(i)$] For $m \ll 0$, $\pi_m M = 0$.

\item[$(ii)$] For all $m \in \Z$, $\pi_m M$ is finitely presented as a $\pi_0 A$-module.
\end{itemize}
\end{proposition}

\begin{proof}
Without loss of generality, we may assume that $M$ is connective. Suppose first that $M$ is almost perfect. We will prove by induction on $n$ that $\pi_n M$ is finitely presented as a $\pi_0 A$-module. If $n=0$, this simply reduces to the observation that the compact objects of the ordinary category of left $\pi_0 A$-modules are precisely the finitely presented $\pi_0 A$-modules. In particular, we can choose a finitely generated free $A$-module $P$ and a map $\alpha: P \rightarrow M$ which induces a surjection $\pi_0 P \rightarrow \pi_0 M$. Since $A$ is coherent, the homotopy groups $\pi_m P$ are finitely presented $\pi_0 A$-modules. Let $K = \ker(\alpha)$. Then $K$ is connective by construction, and almost perfect by Proposition \ref{almor}. The inductive hypothesis implies that $\pi_{i} K$ is finitely presented for $0 \leq i < n$. 

We have a short exact sequence $$0 \rightarrow \coker( \pi_{n} K \rightarrow \pi_{n} P) \rightarrow \pi_{n} M \rightarrow \ker( \pi_{n-1} K \rightarrow \pi_{n-1} P) \rightarrow 0.$$
Using Lemma \ref{swirly}, we deduce that the outer terms are finitely generated, so that $\pi_{n} M$ is finitely generated. Applying the same reasoning, we conclude that $\pi_{n} K$ is finitely generated, so that $\coker( \pi_{n} K \rightarrow \pi_{n} P)$ is finitely presented. Using the exact sequence again, we conclude that $\pi_n M$ is finitely presented.

Now suppose that the connective left $A$-module $M$ satisfies condition $(ii)$. We will prove that $M$ can be obtained as the geometric realization of a simplicial left $A$-module $P_{\bigdot}$ such that each $P_n$ is a free $A$-module of finite rank. As in the proof of Proposition \ref{almor}, 
it will suffice to show that $M$ is the colimit of a sequence
$$ D(0) \stackrel{f_1}{\rightarrow} D(1) \stackrel{f_2}{\rightarrow} D(2) \rightarrow \ldots$$
where each $\coker(f_n)[-n]$ is a free $A$-module of finite rank. Supposing that the partial sequence $$ D(0) \rightarrow \ldots \rightarrow D(n) \stackrel{g}{\rightarrow} M$$
has been constructed, with the property that $\ker(g)$ is $(n-1)$-connected. If $\pi_n \ker(g)$ is finitely generated as a $\pi_0 A$-module, then we can proceed as in the proof of Proposition \ref{almor}. To verify this, we observe that $D(n)$ is almost perfect and therefore satisfies
$(ii)$ (by the first part of the proof). We now use the exact sequence
$$ 0 \rightarrow \coker( \pi_{n+1} D(n) \rightarrow \pi_{n+1} M) \rightarrow \pi_{n} \ker(g) \rightarrow \ker( \pi_{n} D(n) \rightarrow \pi_{n} M) \rightarrow 0 $$
to conclude that $\pi_{n} \ker(g)$ is finitely presented.
\end{proof}

\begin{proposition}\label{itu}
Let $A$ be a connective $A_{\infty}$-ring. The following conditions
are equivalent:
\begin{itemize}
\item[$(1)$] The $A_{\infty}$-ring $A$ is left coherent.

\item[$(2)$] For every left $A$-module $M$, if $M$ is almost perfect, then $\tau_{\geq 0} M$ is almost perfect.

\item[$(3)$] The full subcategories
$\Mod_{A}^{\aperf} \cap \Mod_{A}^{\geq 0}$ and
$\Mod_{A}^{\aperf} \cap \Mod_{A}^{\leq 0}$ determine a t-structure on $\Mod_{A}^{\aperf}$.
\end{itemize}
\end{proposition}

\begin{proof}
The implication $(1) \Rightarrow (2)$ follows from the description
of almost perfect modules given in Proposition \ref{bbh}. The equivalence
$(2) \Leftrightarrow (3)$ is obvious. We will show that $(3) \Rightarrow (1)$.

Suppose that $(3)$ is satisfied. We note that the
first non-vanishing homotopy group of any almost perfect
$A$-module is a finitely presented module over $\pi_0 A$. 
Applying $(2)$ to the module $A[-n]$, we deduce that
$\pi_n A$ is a finitely presented $\pi_0 A$-module. To complete
the proof, it suffices to show that $\pi_0 A$ is left coherent.

Let $R= \pi_0 A$, and regard $R$ as a discrete left $A$-module. Using condition $(3)$, we deduce that $R$ is almost perfect. Let $I \subseteq R$ be a finitely generated left ideal. Then $I$ is the image (in the classical sense) of a map $f: R^{n} \rightarrow R$. Then $\ker(f)$ (taken in the $\infty$-category $\Mod_{A}$) is almost perfect. We have a short exact sequence
$$ 0 \rightarrow \pi_0 \ker(f) \rightarrow R^n \rightarrow I \rightarrow 0.$$
Since $\ker(f)$ is almost perfect, condition $(3)$ implies that $\pi_0 \ker(f)$ is finitely generated, so that $I$ is finitely presented as a $\pi_0 A$-module. This completes the proof of $(1)$.
\end{proof}

\begin{remark}
Let $A$ be a left coherent $A_{\infty}$-ring, and regard $\Mod_{A}^{\aperf}$ as endowed with the t-structure described in Proposition \ref{itu}. Then 
$\Mod_{A}^{\aperf}$ is right bounded and left complete, and the functor
$M \mapsto \pi_0 M$ determines an equivalence from the
heart of $\Mod_{A}^{\aperf}$ to the (nerve of the) category of
finitely presented left modules over $\pi_0 A$. 
\end{remark}

We conclude with a few remarks about the interaction between finiteness and flatness conditions on modules.

\begin{proposition}\label{projj}
Let $A$ be a connective $A_{\infty}$-ring, and let $M$ be a connective left $A$-module. The following conditions are equivalent:
\begin{itemize}
\item[$(1)$] The left $A$-module $M$ is a retract of a finitely generated free $A$-module.
\item[$(2)$] The left $A$-module $M$ is flat and almost perfect.
\end{itemize}
\end{proposition}

\begin{proof}
The implication $(1) \Rightarrow (2)$ is obvious. Conversely, suppose that $M$ is flat and almost perfect. Then $\pi_0 M$ is a left module over $\pi_0 A$ which is finitely presented and flat, and therefore projective. Using Proposition \ref{redline}, we deduce that $M$ is projective.
Choose a map $f: P \rightarrow M$, where $P$ is a free module of finite rank and the induced map $\pi_0 P \rightarrow \pi_0 M$ is surjective. Since $M$ is projective, the map $f$ splits, so that
$M$ is a summand of $P$ in $\h{\calC}$. This proves $(1)$.
\end{proof}

In particular, if an almost perfect left $A$-module $M$ is flat, then $M$ is perfect. We conclude with a mild generalization of this statement, where the flatness hypothesis is relaxed.

\begin{definition}\label{toramp}\index{$\Tor$-amplitude}
Let $A$ be a connective $A_{\infty}$-ring. We will say that a left
$A$-module $M$ has {\it $\Tor$-amplitude $\leq n$} if, for every
discrete right $A$-module $N$, the homotopy groups
$\pi_i( N \otimes_{A} M)$ vanish for $i > n$. We will say that
$M$ is {\it of finite $\Tor$-amplitude} if it has $\Tor$-amplitude $\leq n$ for some integer $n$.
\end{definition}

\begin{remark}\label{lazardapp}
In view of Theorem \ref{lazard}, a connective left $A$-module $M$ has $\Tor$-amplitude $\leq 0$ if and only if $M$ is flat.
\end{remark}

\begin{proposition}\label{lastone}
Let $A$ be a connective $A_{\infty}$-ring.
\begin{itemize}
\item[$(1)$] If $M$ is a left $A$-module of $\Tor$-amplitude $\leq n$, then
$M[k]$ has $\Tor$-amplitude $\leq n+k$.

\item[$(2)$] Let 
$$M' \rightarrow M \rightarrow M'' \rightarrow M'[1]$$ be a distinguished triangle of left $A$-modules. If $M'$ and $M''$ have $\Tor$-amplitude $\leq n$, then so does $M$.

\item[$(3)$] Let $M$ be a left $A$-module of $\Tor$-amplitude $\leq n$. Then any retract of $M$ has $\Tor$-amplitude $\leq n$.

\item[$(4)$] Let $M$ be an almost perfect left module over $A$. Then $M$ is perfect if and only if $M$ has finite $\Tor$-amplitude.
\end{itemize}
\end{proposition}

\begin{proof}
The first three assertions follow immediately from the exactness of the functor
$N \mapsto N \otimes_{A} M$. It follows that the collection left $A$-modules of finite
$\Tor$-amplitude is stable under retracts and finite colimits, and contains the module $A[n]$ for every integer $n$. This proves the ``only if'' direction of $(4)$. For the converse, let us suppose that
$M$ is an almost perfect of finite $\Tor$-amplitude. We wish to show that $M$ is perfect.
We first apply $(1)$ to reduce to the case where
$M$ is connective. The proof now goes by induction on the $\Tor$-amplitude $n$ of $M$. If
$n = 0$, then $M$ is flat and we may conclude by applying Proposition \ref{projj}. We may therefore assume $n > 0$.

Since $M$ is almost perfect, there exists a free left $A$-module $P$ of finite rank and
a distinguished triangle $$ M' \rightarrow P \stackrel{f}{\rightarrow} M \rightarrow M'[1]$$
where $f$ is surjective. To prove that $M$ is perfect, it will suffice to show that $P$ and $M'$ are perfect. It is clear that $P$ is perfect, and it follows from Proposition \ref{almor} that $M'$ is almost perfect. Moreover, since $f$ is surjective, $K$ is connective. We will show that $M'$ is of $\Tor$-amplitude $\leq n-1$; the inductive hypothesis will then imply that $M$ is perfect, and the proof will be complete. 

Let $N$ be a discrete right $A$-module. We wish to prove that $\pi_{k}( N \otimes_{A} K) \simeq 0$ for $k \geq n$. Since the functor $N \otimes_{A} \bigdot$ is exact, 
we obtain for each an exact sequence of homotopy groups
$$ \pi_{k+1}( N \otimes_{A} M) \rightarrow \pi_{k}( N \otimes_{A} M') \rightarrow
\pi_{k}( N \otimes_{A} P).$$
The left entry vanishes in virtue of our assumption that $M$ has $\Tor$-amplitude $\leq n$.
We now complete the proof by observing that $\pi_{k}( N \otimes_{A} P)$ is a finite
direct sum of copies of $\pi_{k} N$, and therefore vanishes because $k \geq n > 0$ and
$N$ is discrete.
\end{proof}

\begin{remark}\label{tukkus}
Let $A$ be a connective $A_{\infty}$-ring, and let $\calC$ be the smallest stable subcategory
of $\Mod_{A}$ which contains all finitely generated projective modules. Then
$\calC = \Mod_{A}^{\perf}$. The inclusion $\calC \subseteq \Mod_{A}^{\perf}$ is obvious.
To prove the converse, we must show that every object $M \in \Mod_{A}^{\perf}$ belongs to $\calC$.
Invoking Corollary \ref{whumph}, we may reduce to the case where $M$ is connective. We then work by induction on the (necessarily finite) $\Tor$-amplitude of $M$. If $M$ is $\Tor$-amplitude zero, then $M$ is flat and the desired result follows from Proposition \ref{projj}. In the general case, we choose a finitely generated free $A$-module $P$ and a map $f: P \rightarrow M$ which
induces a surjection $\pi_0 P \rightarrow \pi_0 M$ (which is possible in view of Corollary \ref{whumph}). As in the proof of Proposition \ref{lastone}, we may conclude that
$\ker(f)$ is a connective perfect module of smaller $\Tor$-amplitude than $M$, so that
the $\ker(f) \in \calC$ by the inductive hypothesis. Since $P \in \calC$ and $\calC$ is stable under the formation of cokernels, we conclude that $M \in \calC$ as desired.
\end{remark}

\begin{example}\label{tukus}
Let $A$ be a discrete $A_{\infty}$-ring, let $\Mod_{A}^{\dplus} \subseteq \Mod_{A}$ be the full subcategory consisting of bounded-above objects, and let let $M \in \Mod^{\dplus}_{A}$. 
Using an inverse to the functor $\theta$ of Proposition \ref{derivdisc}, we can identify any 
$M \in \Mod_{A}^{\dplus}$ with a (bounded above) complex $P_{\bigdot}$ of projective left $\pi_0 A$-modules. It follows from Remark \ref{tukkus} that $M$ is perfect if and only if $P_{\bigdot}$ can be chosen to have only finitely many terms, each of which is finitely generated over $\pi_0 A$. 
\end{example}




\begin{thebibliography}{99}

\bibitem{giraud} Artin, M. {\it Th\'{e}orie des topos et cohomologie
\'{e}tale des sch\'{e}mas.} SGA 4. Lecture Notes in Mathematics
269, Springer-Verlag, Berlin and New York, 1972.

\bibitem{artinmazur} Artin, M. and B. Mazur. {\it \'{E}tale Homotopy.} Lecture Notes in Mathematics 100, Springer-Verlag, Berlin and New York, 1969.

\bibitem{virtual} Behrend, K. and B. Fantechi. {\it The intrinsic
normal cone.} Inventiones Mathematicae 128 (1997) no. 1, 45-88.

\bibitem{BBD} Beilinson, A. , Bernstein, J. and P. Deligne.
{\it Faisceaux pervers.} Asterisuqe 100, Volume 1, 1982.

\bibitem{eilenbergsteenrod} Eilenberg, S. and N.E. Steenrod. {\it Axiomatic approach to homology theory.} Proc. Nat. Acad. Sci. U.S.A. 31, 1945, 117-120.

\bibitem{eisenbud} Eisenbud, D. {\it Commutative algebra.} Springer-Verlag, New York, 1995. 

\bibitem{EKMM} Elmendorf, A.D., Kriz, I. , Mandell, M.A., and J.P.
May. {\it Rings, modules and algebras in stable homotopy theory.}
Mathematical Surveys and Monographs 47, American Mathematical
Society, 1997.

\bibitem{bezout} Fulton, W. {\it Algebraic curves.} W.A.
Benjamin, Inc., New York, 1969.

\bibitem{goerssjardine} Goerss, P. and J.F. Jardine. {\it Simplicial Homotopy Theory.} Progress in Mathematics, Birkhauser, Boston, 1999.

\bibitem{goodwillie} Goodwillie, T. {\it Calculus III: Taylor Series.} Geometry and Topology, Volume 7 (2003) 645-711.

\bibitem{hatcher} Hatcher, A. {\it Algebraic Topology}. Cambridge University Press, 2002.

\bibitem{bordism} Hook, E.C. {\it Equivariant cobordism and
duality.} Transactions of the American Mathematical Society 178
(1973) 241-258.

\bibitem{stablemodel} Hovey, M. {\it Model Categories.}
Mathematical Surveys and Monographs 63, AMS, Providence, RI, 1999.

\bibitem{symmetricspectra} Hovey, M., Shipley, B. and J. Smith. {\it Symmetric spectra.} Journal of the American Mathematical Society 13, 2000, no. 1, 149-208.

\bibitem{illusie} Illusie, L. {\it Complexe cotangent et d\'{e}formations
I}. Lecture Notes in Mathematics 239, Springer-Verlag, 1971.

\bibitem{illusie2} Illusie, L. {\it Complexe cotangent et d\'{e}formations
II}. Lecture Notes in Mathematics 283, Springer-Verlag, 1972.

\bibitem{kerz} Kerz, M. {\it The complex of words and Nakaoka stability.} Homology, Homotopy and Applications, volume 7(1), 2005, pp. 77-85.

\bibitem{categoricalring} Laplaza, M. {\it Coherence for
distributivity.} Coherence in categories, 29-65. Lecture Notes in
Mathematics 281, Springer-Verlag, 1972.

\bibitem{lazard} Lazard, Daniel. {\it Sur les modules plats.} C.R.
Acad. Sci. Paris 258, 1964, 6313-6316.

\bibitem{topoi} Lurie, J. {\it Higher Topos Theory.} Available for download at http://www.math.harvard.edu/~lurie/ .

\bibitem{DAGStable} Lurie, J. {\it Derived Algebraic Geometry I: Stable $\infty$-Categories.} Available for download.

\bibitem{monoidal} Lurie, J. {\it Derived Algebraic Geometry II: Noncommutative Algebra.} Available for download.

\bibitem{symmetric} Lurie, J. {\it Derived Algebraic Geometry III: Commutative Algebra.} Available for download.

\bibitem{deformation} Lurie, J. {\it Derived Algebraic Geometry IV: Deformation Theory.} In preparation.

\bibitem{structured} Lurie, J. {\it Derived Algebraic Geometry V: Structured Spaces.} In preparation.

\bibitem{maclane} MacLane, S. {\it Categories for the Working Mathematician.} Second edition. Graduate Txts in Mathematics, 5. Springer-Verlag, New York, 1998.

\bibitem{gabriel} Mitchell, B. {\it A quick proof of the
Gabriel-Popesco theorem.} Journal of Pure and Applied Algebra 20
(1981), 313-315.

\bibitem{neeman} Neeman, A. {\it Triangulated categories.} Annals
of Mathematics Studies, 148. Princeton University Press, 2001.

\bibitem{homotopicalalgebra} Quillen, D. {\it Homotopical Algebra.} Lectures Notes in Mathematics 43, SpringerÐVerlag, Berlin, 1967. 

\bibitem{homotopyvarieties} Rosicki, J. {\it On Homotopy Varieties.} Available as math.CT/0509655 .

\bibitem{schwede} Schwede, S. {\it Spectra in model categories and applications to the algebraic
cotangent complex.} Journal of Pure and Applied Algebra 120 (1997), pp.
77-104.

\bibitem{monmod} Schwede, S. and B. Shipley. {\it Algebras and Modules in Monoidal Model Categories.} Proceedings of the London Mathematical Society (80) 2000, 491-511.

\bibitem{schwedeshipley} Schwede, S. and B. Shipley. {\it Stable model categories are categories of modules.} Topology 42, 2003, no. 1, 103-153.

\bibitem{intersection} Serre, Jean-Pierre. {\it Local algebra.}
Springer-Verlag, 2000.

\bibitem{shipley} Shipley, B. {\it A Convenient Model Category for Commutative Ring Spectra.} Homotopy theory: relations with algebraic geometry, group cohomology, and algebraic $K$-theory. Contemp. Math. volume 346 pp. 473-483, American Mathematical Society, Providence, RI, 2004. 

\bibitem{srinivas} Srinivas, V. {\it Algebraic K-Theory.} Birkhauser, Boston, 1993.

\bibitem{toenK} To\"{e}n, B. and G. Vezzosi. {\it A remark on K-theory and S-categories.} Topology 43, No. 4 (2004), 765-791

\bibitem{weibel} Weibel, C. {\it An Introduction to Homological Algebra.} Cambridge University Press, 1995.

\end{thebibliography}
\end{document}